\documentclass{article}
\usepackage{amssymb}
\usepackage{pifont}
\usepackage{amsmath,amsthm,amssymb,amsfonts}
\usepackage{comment}
\usepackage{amsfonts}
\usepackage{amsmath, amsthm}
\usepackage{amssymb, amsxtra}
\usepackage{mathtools}
\usepackage{graphicx}

\usepackage{amscd}
\usepackage{multicol}
\usepackage{epstopdf}
\usepackage{color}

\usepackage[numbers,sort&compress]{natbib}
\usepackage{epsf}
\usepackage{float}
\usepackage{extarrows}
\usepackage[a4paper, text={.8\paperwidth,.83\paperheight}, ratio=1:1]{geometry}
\usepackage[format=hang,font=small,textfont=it]{caption}
\usepackage{authblk}

\usepackage[english]{babel}
\usepackage{latexsym}
\usepackage{nicefrac}

\usepackage{pgf, tikz, pgfplots}
\pgfplotsset{compat=1.15}
\usepackage{mathrsfs}
\usetikzlibrary{arrows}
\usetikzlibrary{calc}
\usepackage[tight]{subfigure}
\usepackage{hyperref}

\renewcommand{\labelenumi}{$($\arabic{enumi}$)$}
\renewcommand{\labelenumii}{(\Alph{enumi})}

\usepackage{color}
\definecolor{greenbean}{RGB}{199,237,204}

\usepackage[T1]{fontenc}%
\usepackage[utf8]{inputenc}%
\usepackage{lmodern}%
\usepackage{textcomp}%
\usepackage{lastpage}%
\usepackage{tikz}%
\usepackage{float}
\usepackage{amsfonts, amsmath}
\usepackage{hyperref}
\usepackage{environ}
\usetikzlibrary{arrows}

\makeatletter
\newsavebox{\measure@tikzpicture}
\NewEnviron{scaletikzpicturetowidth}[1]{%
	\def\tikz@width{#1}%
	\def\tikzscale{1}\begin{lrbox}{\measure@tikzpicture}%
		\BODY
	\end{lrbox}%
	\pgfmathparse{#1/\wd\measure@tikzpicture}%
	\edef\tikzscale{\pgfmathresult}%
	\BODY
}
\makeatother
\newtheorem{thm}{Theorem}[section]
\newtheorem{lemma}[thm]{Lemma}
\newtheorem{prop}[thm]{Proposition}
\newtheorem{cor}[thm]{Corollary}
\newtheorem{eg}{Example}[section]

\theoremstyle{definition}
\newtheorem{Def}[thm]{Definition}
\newtheorem{Construction}[thm]{Construction}
\newtheorem{remark}{Remark}

\def\CC{\mathbb{C}}
\def\cpt{\CC \mathbb{P}^2}
\newcommand{\pcpt}[1]{{\pi_1(\cpt - #1, *)}}

\DeclarePairedDelimiterX\setc[2]{\{}{\}}{\,#1 \;\delimsize\vert\; #2\,}

\newcommand{\notion}[1]{\textit{#1}}
\newcommand{\Bigp}[1]{\Big(#1\Big)}
\renewcommand{\baselinestretch}{1.5}

\newcommand{\ug}[1]{\Gamma_#1}
\newcommand{\Ggal}{\widetilde{G}}

\DeclareMathOperator{\sing}{Sing}
\DeclareMathOperator{\pr}{pr}

\makeatletter
\newcommand{\uGammaSq}[2]{\begin{equation}\label{#1}
		\ug{#2}^2\uGammaSqChecknextarg}
	\newcommand{\uGammaSqChecknextarg}{\@ifnextchar\bgroup{\uGammaSqGobblenextarg}{ = e,
\end{equation} }}
\newcommand{\uGammaSqGobblenextarg}[1]{ = \ug{#1}^2\@ifnextchar\bgroup{\uGammaSqGobblenextarg}{  = e,
\end{equation}}}
\makeatother

\newcommand\Corref[1]{{Corollary~\ref{#1}}} 
\newcommand\CAref[1]{{Chapter~\ref{#1}}}    
\newcommand\Rref[1]{{Remark~\ref{#1}}}      
\newcommand\Sref[1]{{Section~\ref{#1}}}     
\newcommand\ssref[1]{{Subsection~\ref{#1}}} 
\newcommand\Pref[1]{{Proposition~\ref{#1}}} 
\newcommand\Tref[1]{{Theorem~\ref{#1}}}     
\newcommand\Dref[1]{{Definition~\ref{#1}}}  
\newcommand\DGref[1]{{Diagram~\ref{#1}}}    
\newcommand\Lref[1]{{Lemma~\ref{#1}}}       
\newcommand\Eqref[1]{Equation~(\ref{#1})}   
\newcommand\Iref[1]{{Item~(\ref{#1})}}      
\newcommand\Fref[1]{{Figure~\ref{#1}}}      
\newcommand\eq[1]{{(\ref{#1})}}             

\newcommand{\ubegineq}[1]{\begin{equation}\label{#1}}
\newcommand{\uendeq}{\end{equation}}
\newcommand{\trip}[2] {{\langle \ug{#1}, \ug{#2} \rangle}}
\newcommand{\simptrip}[2] {{\langle {#1}, {#2} \rangle}}
\newcommand{\tripFirstPrime}[2] {{\langle \ug{#1}', \ug{#2} \rangle}}
\newcommand{\tripSecondPrime}[2] {{\langle \ug{#1}, \ug{#2}' \rangle}}
\newcommand{\tripBothPrime}[2] {{\langle \ug{#1}', \ug{#2}' \rangle}}
\newcommand{\comm}[2] {{[\ug{{#1}}, \ug{{#2}}]}}
\newcommand{\commFirstPrime}[2] {{[\ug{#1}', \ug{#2}]}}
\newcommand{\commSecondPrime}[2] {{[\ug{#1}, \ug{#2}']}}
\newcommand{\simpcomm}[2] {{[{#1}, {#2}]}}

\makeatletter
\newcommand{\uProjRel}[2]{
\begin{equation}\label{#1}
	\ug{#2'}\ug{#2}\uProjRelChecknextarg}
\newcommand{\uProjRelChecknextarg}{\@ifnextchar\bgroup{\uProjRelGobblenextarg}{ = e
\end{equation}}}
\newcommand{\uProjRelGobblenextarg}[1]{\ug{#1'}\ug{#1}\@ifnextchar\bgroup{\uProjRelGobblenextarg}{  = e. \end{equation}}}
\makeatother

\def\Z{\mathbb{Z}}
\newcommand{\mbb}[1]{\mathbb{#1}}
\newcommand{\todo}[1]{\textcolor{red}{[TODO: #1]}}

\newcommand{\onlyinsubfile}[1]{#1}
\newcommand{\notinsubfile}[1]{}

\begin{document}

\title{Zariski pairs of conic-line arrangements of degrees 7 and 8 via fundamental groups\footnotetext{Email addresses: Meirav Amram (corresponding author): meiravt@sce.ac.il; Robert Shwartz: robertsh@ariel.ac.il; \\ Uriel Sinichkin: sinichkin@mail.tau.ac.il; Sheng-Li Tan: sltan@math.ecnu.edu.cn;  Hiro-O Tokunaga: tokunaga@tmu.ac.jp.\\
	2020 Mathematics Subject Classification. 14H30, 14H50, 14Q05, 20F55. \\{\bf Key words}: Zariski pairs, conic-line arrangements, fundamental group, Coxeter groups.
}}
\author[1]{Meirav Amram}
\author[2]{Robert Shwartz}
\author[3]{Uriel Sinichkin}
\author[4]{Sheng-Li Tan}
\author[5]{Hiro-O Tokunaga}

\affil[1]{\small{Department of Mathematics, Shamoon College of Engineering, Ashdod, Israel; ORCID  ID: 0000-0003-4912-4672}}
\affil[2]{\small{Department of Mathematics, Ariel University, Ariel, Israel}}
\affil[3]{\small{School of Mathematical Sciences, Tel Aviv University, Tel Aviv, Israel}}
\affil[4]{\small{School of Mathematical Sciences, Shanghai Key Laboratory of PMMP, East China Normal University, Shanghai 200241, P. R. China}}
\affil[5]{\small{Department of Mathematical Sciences, Tokyo Metropolitan University, Japan}}

\maketitle

\abstract{
We find a new Zariski pair with non-isomorphic fundamental groups that consists of degree $ 8 $ conic-line arrangements. Each arrangement has three conics and two lines.
We use the Zariski-van Kampen Theorem and some known Coxeter groups to determine the fundamental groups.

Two examples of degree $7$ Zariski pairs that were introduced in 2014 by the last named author, are given as well. They consist of a pair of conic-line arrangements with three conics in each (and thus, each has a single line) and a pair with two conics in each (and thus, each has three lines).
We were able to provide alternative proof of the fact those are indeed Zariski pairs by our methods.
}

\section{Introduction}

The study of Zariski pairs, originated by Zariski \cite{zariski29}, has its roots in the question of which plane curves can be realized as the ramification locus of a branched covering.
This issue is still wide-open and has given rise to various problems in algebraic geometry, such as the Chisini conjecture.
The study of Zariski pairs problems related to finding and classifying them, has led many researchers to study this phenomenon.
Before we mention the contributions of these researchers, we provide the definition of the Zariski pair as presented in \cite{Bartolo_tokunaga_survey}.
\begin{Def}\label{def:zariski_pairs}
A pair of reduced plane curves $ \mathcal{B}_1,\mathcal{B}_2 \subseteq \mathbb{CP}^2 $ is called a \emph{Zariski pair} if:
\begin{enumerate}
	\item There exist tubular neighborhoods $ T(\mathcal{B}_1) $ and  $ T(\mathcal{B}_2) $ and a homeomorphism
	$ h: T(\mathcal{B}_1)\to T(\mathcal{B}_2) $ with $ h(\mathcal{B}_1)=\mathcal{B}_2 $.
	
	\item There is no homeomorphism $ f:\mathbb{CP}^2\to \mathbb{CP}^2 $ with $ f(\mathcal{B}_1)=\mathcal{B}_2 $.
\end{enumerate}
\end{Def}

\begin{remark}
	In some texts, the definition of a Zariski pair requires the stronger assumption that $\cpt - \mathcal{B}_1$ and $\cpt - \mathcal{B}_2$ are not homeomorphic.
\end{remark}

The first example was presented by Zariski himself and consists of 2 sextics with 6 cusps each, where the cusps lie on a conic for one of the curves and are not in such a special position for the other \cite{zariski29}.
Since then, almost a century has passed and various techniques for finding Zariski pairs have been developed, especially by Artal-Bartolo, Cogolludo-Agustin, and Tokunaga in   \cite{Bartolo_tokunaga_survey}, and by Artal-Bartolo in \cite{Bartolo94}. Artal-Bartolo pointed that two curves with isomorphic combinatorics and different topologies are named a Zariski pair.
Shimada extended Artal-Bartolo's method of Kummer coverings from \cite{Bartolo94}, to construct two infinite series of Zariski pairs of growing degrees in \cite{shimada}.
In \cite[Appendix A]{shustin_appendix}, Shustin used patchworking construction to construct plane curves of degree $ \nu (\nu-1) $ (for $ 3\le \nu \le 10 $), which cannot appear as branch loci of degree $ \nu $ surfaces, yielding examples of Zariski pairs of degrees $ d\in \{ 6, 12, 20, 30, 42, 56, 72, 90 \} $. 
Bannai and Tokunaga \cite{bannai_tokunaga2019zariski}, Artal-Bartolo, Bannai, Shirane, and Tokunaga \cite{bartolo_tokunaga2020torsion, bartolo_tokunaga2020torsion2}, and Takahashi and Tokunaga \cite{takahashi_tokunaga2020explicit} recently used the arithmetic of elliptic curves to construct examples of Zariski pairs and Zariski tuples. 
Oka directly compared the fundamental groups $ \pcpt{\mathcal{B}_1} $ and $ \pcpt{\mathcal{B}_2} $ of the complements of $\mathcal{B}_1$ and $\mathcal{B}_2$ in \cite{Oka1999FlexCA, Oka2002}, which is also Zariski's original point of view on the problem. Artal-Bartolo, Cogolludo-Agustin, and Martin-Morales \cite{artal1} found examples of Zariski pairs in weighted projective planes, distinguished by the Alexander polynomial, which is a very efficient way to study Zariski pairs. Other invariants were studied to deal with pairs of curves and with the complexity of groups. Artal-Bartolo and Dimca \cite{Dimca} studied properties of fundamental groups and Alexander polynomials of plane curves,   describing some topological properties of curves with an abelian fundamental group. Artal-Bartolo, Cogolludo-Agustin, and Carmona \cite{Artal2} proved that braid monodromy of an affine plane curve determines the topology of a related projective plane curve (and the monodromy can be applied to the Zariski-van Kampen Theorem to get a fundamental group). Artal-Bartolo, Carmona, Cogolludo-Agustin, Luengo, and Melle \cite{Artal3} presented invariants such as the Alexander polynomial and characteristic varieties,  that distinguish many pairs. They also reviewed the Zariski-van Kampen method for finding a presentation of the fundamental group of the complement of an algebraic curve in $\mathbb{P}^2$ and stressed that it is a sensitive invariant for distinguishing many pairs of curves, which therefore has a prior importance in the classification of curves. In general, one of the best sources for studying fundamental groups related to geometric objects or considered as geometric objects by themselves, is the book written by Drutu-Kapovich \cite{Kap}.

We focus on the last technique mentioned, namely the comparison of the groups $\pcpt{\mathcal{B}_1}$ and $\pcpt{\mathcal{B}_2}$.
This technique is very strong in the sense that many of the known examples of Zariski pairs, for which the fundamental groups \footnote{Here and henceforth we use the term "fundamental group of $\mathcal{B}$" to denote $ \pcpt{\mathcal{B}} $. } were computed, have distinct fundamental groups. Note that the technique poses the challenge of solving the isomorphism problem for a pair of finitely presented groups, a famously undecidable problem.
We use this technique in the current work to construct a new example of Zariski pair, as well as reproving the correctness of some known Zariski pairs.
The study of Coxeter groups is a naturally appearing algebra in such an endeavor, a method that we hope can pave the way for a more systematic investigation of Zariski pairs. The technique consists of finding braids related to singularities in the conic-line arrangements, then by the Zariski-van Kampen Theorem \cite{vanKampen33}, we can derive presentations of those curves by means of generators and relations. We then define the Coxeter quotients of the fundamental groups and by simplifications of relations and using Coxeter groups, we are able to determine the fundamental groups, thereby showing that they are non-isomorphic. We note that there are Zariski pairs with isomorphic fundamental groups, see \cite{AC}, \cite{new}, and \cite{Ben} as examples.

In \cite{DEG_isotopy}, Degtyarev proves that the rigid isotopy class of curves of degree at most 5
is determined by the combinatorial data. In particular, his work implies that, up to
degree 5, there are no Zariski pairs. He also classifies in \cite{DEG} the isotopy classes of the
curves of degree 6 with simple singularities.

Most of the known examples of Zariski pairs, including the original pair by Zariski, are pairs of irreducible curves.
One can try to find Zariski pairs among curves at the other extreme - in the unions of curves of low degree.
Arrangements of lines were studied quite extensively.
The first Zariski pairs of line arrangements were found by Rybnikov in \cite{Rybnikov_line_arrangement}.
Conversely, it was shown by Nazir and Yoshinaga in \cite{Nazir_Yoshinaga_line_arrangements} that there are no Zariski pairs of arrangements of 8 lines, and by Ye in \cite{Ye2013ClassificationOM} that no such pairs occur for arrangements of 9 lines.
In this scenario Zariski pairs gain additional importance - it was shown by Randell in \cite{Randell} that Zariski pairs of line arrangements cannot lie in the same component of the moduli space, thereby giving us a powerful tool to study the moduli space of line arrangements.

The case of conic-line arrangements (unions of conics and lines) is much less studied.
As far as we know, there are no results for the moduli space of conic-line arrangements that are  similar to Randell's above-mentioned theorem.
As for Zariski pairs, there are few examples - Namba and Tsuchihashi, in \cite{Namba2004OnTF},  used the fundamental group of the complement to construct an example of a Zariski pair of a degree 8 conic-line arrangement; Tokunaga used the theory of dihedral covers to construct a  few such examples of degree 7 in \cite{tokunaga2014}. 
The fundamental groups of conic-line arrangements were studied also in \cite{AGT_order6_conic_line_arr} for arrangements of degree up to 6, but no Zariski pairs of degree 6 were found in the cases considered there. 
Therefore, a question worth exploring in future work can be for arrangements of degree 7 and higher than that. 
A similar  question was also asked lately in \cite{ben2}. 
The work of Friedman and Garber \cite{Friedmann_Garber_conicline_fundamental}, where they generalized some known results regarding line arrangements to conic-line arrangements, seems to be  relevant to our investigation.
Although we did not use the Friedman and Garber results in the current work, we are certain it could be useful for any further investigation.

This work is significant in the exploration of the moduli space of conic-line arrangements of degrees 7 and 8, as it introduces new examples of Zariski pairs with non-isomorphic fundamental groups. Constructing such examples in low degrees is notably difficult. By overcoming this challenge, our research contributes to a deeper understanding of the topological and algebraic properties of these arrangements. The results not only assist in the classification of Zariski pairs but also pave the way for further investigations into the moduli space of conic-line arrangements.

The fundamental groups appearing in the paper are complicated, and therefore special quotients of them, obtained by squares of generators, are used. Such usage works excellently in the classification of algebraic surfaces where special coverings are involved, but is not known or perhaps not recognized in the classification of algebraic curves. Therefore, the use of such quotients is unique and interesting, and it is possible to explicitly determine what these groups are.

The structure of the current paper is as follows:
In Section \ref{sec:new_pair} we present a pair of conic-line arrangements of degree 8, then in Section \ref{presentation} we calculate the fundamental groups of their complements in $\mathbb{CP}^2$. In Section \ref{COX} we recall some Coxeter groups theory and then define square quotients of the groups, which assist us in Section \ref{sec:existing} to determine that we have a Zariski pair. Moreover, we give examples of two pairs of conic-line arrangements that were studied by Tokunaga in 2014, and show that the method used in the paper enables us to provide additional proof of the fact they form Zariski pairs.

\paragraph{Acknowledgments:}
Thanks are given to Enrique Artal-Bartolo and Jose Ignacio Cogolludo Agustin for their helpful comments and fruitful conversations.
This research was supported by ISF-NSFC joint research program (grant No. 2452/17).
We thank two anonymous referees for their great generosity, their mathematical suggestions, and also for the suggestion for a much richer, clear and orderly structure of the paper, which was implemented.

\section{Description of the new Zariski pair}\label{sec:new_pair}

In this section we give some details about Zariski pairs, their combinatorics, and related definitions.
Then we present the new Zariski pair we have found, which consists of two conic-line arrangements of degree 8.

In accordance with \cite{Bartolo_tokunaga_survey}, the first condition of Definition \ref{def:zariski_pairs} can be easily checked by inspecting the combinatorics of the curves. 
To simplify the notation, we will restrict to the case of conic-line arrangements, which gives the following definition.

\begin{Def}\label{def:combinatorial_type}
The \emph{combinatorial type} of a conic-line arrangement $ \mathcal{B} $ is the tuple
$$ \left( \text{Irr}(\mathcal{B}), \deg, \sing(\mathcal{B}), \Sigma_{\text{top}}, \sigma_{\text{top}}, \{ \mathcal{B}(P) \}_{P\in \sing(\mathcal{B})}, \{ \beta_P \}_{P\in \sing(\mathcal{B})}, \{ i_P \}_{P\in \sing(\mathcal{B})} \right).  $$
In the above expression: 
\begin{itemize}
	\item
	$ \text{Irr}(\mathcal{B}) $ is the set of irreducible components of $ \mathcal{B} $, $ \deg :\text{Irr}(\mathcal{B})\to \mathbb{N} $ assigns to each component its degree.
	\item
	$ \sing(\mathcal{B}) $ is the set of singular points of $ \mathcal{B} $, $ \Sigma_{\text{top}} $ is the set of topological types of singularities occurring in $ \mathcal{B} $, and $ \sigma_{\text{top}}:\sing(\mathcal{B})\to \Sigma_{\text{top}} $ assigns to each singular point its topological type.
	\item
	$ \mathcal{B}(P) $ is the set of local branches at a point $ P\in \sing(\mathcal{B}) $ and $ \beta_P:\mathcal{B}(P)\to \text{Irr}(\mathcal{B}) $ assigns to each local branch its irreducible component of $ \mathcal{B} $.
	\item
	$ i_P:\binom{\mathcal{B}(P)}{2}\to \mathbb{N} $ assigns to each pair of local branches at $ P\in \sing(\mathcal{B}) $ their intersection multiplicity.
\end{itemize}

Two conic-line arrangements $ \mathcal{B}_1, \mathcal{B}_2 $ have \emph{the same combinatorial type} (or simply \emph{the same combinatorics}) if there exist bijections $ \varphi_{\text{Irr}}:\text{Irr}(\mathcal{B}_1) \to \text{Irr}(\mathcal{B}_2) $, $ \varphi_{\sing}: \sing(\mathcal{B}_1) \to \sing(\mathcal{B}_2) $, and $ \varphi_P:\mathcal{B}_1(P) \to \mathcal{B}_2(\varphi_{\sing}(P)) $ for every $ P\in \sing(\mathcal{B}_1) $, such that $ \varphi_{\text{Irr}} $ preserves degree, $ \varphi_{\sing} $ preserves topological type, $ \beta_{\varphi_{\sing}(P)}\circ \varphi_P = \varphi_{\text{Irr}}\circ \beta_{P} $ for all $ P\in \sing(\mathcal{B}_1) $, and $ i_{\varphi_{\sing}(P)}\left(\varphi_P(C_1), \varphi_P(C_2)\right) = i_{P}(C_1, C_2) $ for all $ P\in \sing(\mathcal{B}_1) $ and $ C_1, C_2\in \mathcal{B}_1(P) $.
\end{Def}

If two conic-line arrangements have the same combinatorics, they satisfy the first condition of Definition \ref{def:zariski_pairs}, see \cite[Remark 3]{Bartolo_tokunaga_survey}; so this is a relatively simple condition to check.

The more intricate task is checking the second condition of Definition \ref{def:zariski_pairs}.
For that we note that any homeomorphism $ (\cpt, \mathcal{B}_1)\to (\cpt, \mathcal{B}_2) $ gives rise to any isomorphism $ \pcpt{\mathcal{B}_1}\to \pcpt{\mathcal{B}_2} $, so if those groups are not isomorphic, the pair $ (\mathcal{B}_1, \mathcal{B}_2) $ is a Zariski pair.

\begin{figure}[H]
	\centering
	\begin{tikzpicture} [scale=0.35, style = very thick]
		\begin{scope}[xshift=-330]
			\draw (0,0) ellipse (7 and 4);
			\draw (6.7, 2.7) node {\footnotesize $C_2$};
			
			\draw (0,0) ellipse (4 and 7);
			\draw (2.7, 7) node {\footnotesize $ C_3 $};
			
			\draw (0,0) ellipse (4 and 4);
			\draw (-2, 2) node {\footnotesize $ C_1 $};
			
			\draw (1, 9.1) -- (-9.1, -1);
			\draw (0.5, 9.5) node {\footnotesize $ L_1 $};
			
			\draw (1, -9.1) -- (-9.1, 1);
			\draw (1.1, -8.5) node {\footnotesize $ L_2 $};
			
			\draw (-4.7, 0) node {\footnotesize $Q_1$};
			\draw (4.7, 0) node {\footnotesize $Q_2$};
			\draw (0, 4.7) node {\footnotesize $Q_3$};
			\draw (0, -4.7) node {\footnotesize $Q_4$};
			
			\draw (0, -10) node {$\mathcal{B}_1$};
		\end{scope}
		
		
		\begin{scope}[xshift=330]
			\draw (0,0) ellipse (7 and 4);
			\draw (6.7, 2.7) node {\footnotesize $C_2$};
			
			\draw (0,0) ellipse (4 and 7);
			\draw (2.7, 7) node {\footnotesize $ C_3 $};
			
			\draw (0,0) ellipse (4 and 4);
			\draw (-2, 2) node {\footnotesize $ C_1 $};
			
			\draw (1, 9.1) -- (-9.1, -1);
			\draw (0.5, 9.5) node {\footnotesize $ L_1 $};
			\draw (-1, -9.1) -- (9.1, 1);
			\draw (9.5, 0.5) node {\footnotesize $ L_3 $};
			
			\draw (-4.7, 0) node {\footnotesize $Q_1$};
			\draw (4.7, 0) node {\footnotesize $Q_2$};
			\draw (0, 4.7) node {\footnotesize $Q_3$};
			\draw (0, -4.7) node {\footnotesize $Q_4$};
			
			\draw (0, -10) node {$\mathcal{B}_2$};
		\end{scope}
	\end{tikzpicture}

\caption{The arrangements $ \mathcal{B}_1,\mathcal{B}_2 $.}
\label{fig_new_pair}
\end{figure}
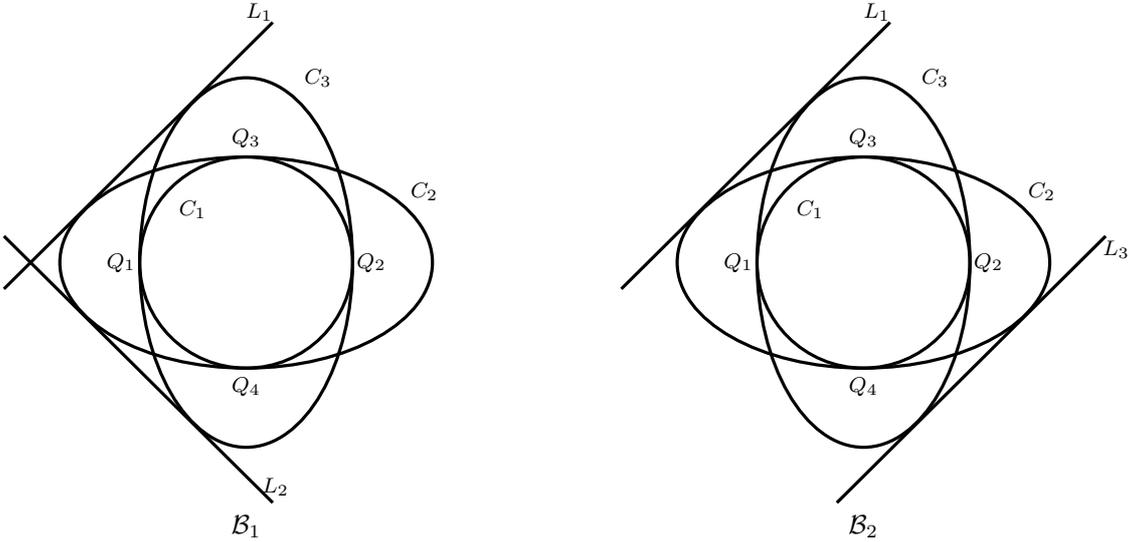

The main result of this paper is the introduction of a new Zariski pair. The Zariski pair is described in Construction \ref{construcetion_tokunaga_arrangement5_construction}. For the reader's convenience, we present the real picture in Figure \ref{fig_new_pair}. However, note that we are interested in the complex plane, so Figure \ref{fig_new_pair} should be used with care.  Construction \ref{construcetion_tokunaga_arrangement5_construction} presents two curves that are conic-line arrangements of degree 8, we call them $\mathcal{B}_1$ and $\mathcal{B}_2$. Each curve has three conics and two lines., i.e., $$\mathcal{B}_1=C_1+C_2+C_3+L_1+L_2 \ \ \  \text{and} \ \ \ \mathcal{B}_2 = C_1+C_2+C_3+L_1+L_3.$$

\begin{Construction}\label{construcetion_tokunaga_arrangement5_construction}
Let $ C_i \left(i=1,2,3\right)$ be smooth conics and let $ L_j \left(j=1,2,3\right) $ be lines as follows (see Figure \ref{fig_new_pair}):
\begin{itemize}
\item $ C_2 $ and $ C_3 $ meet transversely.
\item $ C_1 $ is tangent to both $ C_3 $ and $ C_2 $ such that the intersection multiplicities at intersection points are all equal to 2. We set $ C_1\cap C_2=\{ Q_3, Q_4 \} $ and $ C_1\cap C_3 = \{Q_1, Q_2 \} $.
\item $ L_1 $, $ L_2 $,  and $ L_3 $ are tangent to both $ C_2 $ and $ C_3 $ and intersect $ C_1 $ transversely. Because $ L_1 $, $ L_2 $, and $ L_3 $ are all distinct lines, each pair of lines meets transversely at one point.
	The intersection point of $ L_1 $ and $ L_3 $ in Figure \ref{fig_new_pair} is not depicted because  it is contained in the line at infinity.
\end{itemize}
\end{Construction}

The equations of conics $C_i$ and lines $L_i$ ($i=1,2,3$) in $ \mathcal{B}_1$ and $\mathcal{B}_2 $ are given below with affine coordinates and projective coordinates as well:

\begin{tabular}{c c c}
  conic/line \ & \ affine coordinates \ & \ projective coordinates \\
  $C_1$ \ & \ $x^2+y^2-1=0$ & $X^2+Y^2-Z^2=0$ \\
  $C_2$ \ & \ $\frac{x^2}{3}+y^2-1=0$  \ & \ $\frac{X^2}{3}+Y^2-Z^2=0$ \\
  $C_3$ \ & \ $x^2+\frac{y^2}{3}-1=0$  \ & \ $X^2+\frac{Y^2}{3}-Z^2=0$ \\
  $L_1$ \ & \ $y-x-2=0$ \ & \ $Y-X-2Z=0$ \\
  $L_2$ \ & \ $y+x+2=0$ \ & \ $Y+X+2Z=0$  \\
  $L_3$ \ & \ $y-x+2=0$ \ & \ $Y-X+2Z=0$  \\
\end{tabular}

\section{Presentations of the fundamental groups}\label{presentation}

In this section we give a short background on the Zariski-van Kampen Theorem \cite{vanKampen33} and an explanation of how to derive the fundamental group of a complement of a curve in $\mathbb{CP}^2$. In Lemmas \ref{lemma-fundamental1} and \ref{lemma-fundamental2}, we then compute the fundamental groups $\pcpt{\mathcal{B}_1}$ and $\pcpt{\mathcal{B}_2}$ of the complements of $\mathcal{B}_1$ and $\mathcal{B}_2$ in $\mathbb{CP}^2$.

To compute the fundamental group $ \pcpt{\mathcal{B}} $ for a curve $ \mathcal{B} $ we use the Zariski-van Kampen algorithm as described in \cite{vanKampen33} which, for the sake of completeness, we describe here.
We begin by computing the affine fundamental group $ \pi_1(\mbb{C}^2-\mathcal{B}, *) $; that is, we pick generic line $ L\subseteq \cpt $ and choose coordinates such that $ L $ is the line at infinity.
We then consider the projection $ \pr: \mbb{C}^2 \to\mbb{C}^1 $ given by $ (x,y)\mapsto x $.
The genericity conditions ensure that no tangent line to $ \mathcal{B} $ at a singular point can be parallel to the $ y $-axis.
Let $ q_1,\dots, q_N\in \mathbb{CP}^1 $ be the branch locus of $ \pr|_\mathcal{B} $; that is, the images of the singular points of $ \mathcal{B} $ and the images of points of $ \mathcal{B} $ where the tangent to $ \mathcal{B}$ is parallel to the $ y $-axis (the latter ones are called \emph{branch points}).
Pick a base point $ y_0\in \mathbb{CP}^1-\{q_1,\dots,q_N\} $ and a base point $ x_0\in \pr^{-1}(y_0)-\mathcal{B} $ in its fiber.
One can be convinced that any loop in $ \pi_1(\CC^2 - \mathcal{B}, x_0 ) $ is equivalent to a loop whose image under $ \pr $ avoids points $ q_1,\dots, q_N $.
Covering $ \CC^1-\{q_1,\dots,q_N\} $ by simply-connected open neighborhoods of $ y_0 $ and using the Zariski-van Kampen Theorem, we see that any loop in $ \pi_1(\CC^2-\mathcal{B}, x_0) $ is equivalent to a loop  entirely contained in the fiber $ \pr^{-1}(y_0) - \mathcal{B} $.
So, the fundamental group $ \pi_1(\CC^2-\mathcal{B}, x_0) $ is a quotient of $ \pi_1(\pr^{-1}(y_0)-\mathcal{B}, x_0) $,  which is isomorphic to the free group on $ \deg \mathcal{B}$ generators.

To find the relations that define $ \pcpt{\mathcal{B}} $, we consider the monodromy action of $ \pi_1(\CC^1 - \{ q_1, \dots, q_N \}, y_0) $ on $ \pi_1(\pr^{-1}(y_0)-\mathcal{B}, x_0) $.
For that, for every element $ [\gamma]\in \pi_1(\CC^1 - \{q_1,\dots, q_N\}, y_0) $, we choose a loop $ \widetilde{\gamma} $ in $ \CC^2 - \mathcal{B} $ with $pr(\widetilde{\gamma})=\gamma$ (the precise choise is not important, we chose a loop $\widetilde{\gamma}\subset \{ (x,a+ib)\in \CC^2 \mid b<<0\}$).
Then, for every $[\Gamma]\in \pi_1(\pr^{-1}(y_0)-\mathcal{B}, x_0) $, the loop $\widetilde{\gamma} $ gives rise to a homotopy in $ \CC^2 - \mathcal{B} $ between $[\Gamma]$ and $[\widetilde{\gamma}\Gamma\widetilde{\gamma}^{-1}]$ and so we get a relation in $ \pi_1(\CC^2-\mathcal{B}, x_0) $.
In fact, those are all the relations in $ \pi_1(\CC^2-\mathcal{B}, x_0) $.
To get a representation of $ \pcpt{\mathcal{B}} $, one can use the Zariski-van Kampen Theorem again, which gives one additional relation, called a \emph{projective relation}, which corresponds to the fact that a loop around all the points in $ \mathcal{B}\cap \text{pr}^{-1}(y_0) $ is null-homotopic in $ \mathbb{CP}^2-\mathcal{B} $.

Obviously, it is enough to consider the relations arising from a generating set of the group $ \pi_1(\mathbb{C}^1-\{q_1,\dots,q_N\}, y_0) $.
In all our cases we are able to pick coordinates such that $ q_1,\dots,q_N $ are all real. Those coordinates are not those that are shown in Figure \ref{fig_new_pair}, for more details see \cite[Appendix]{arxiv_version}. We then pick a real coordinate $ y_0 $ to be larger than $ \max\{ q_1, \dots, q_N \} $, and also a generating set $ [\gamma_1],\dots,[\gamma_N] $ of $ \pi_1(\mathbb{C}^1-\{q_1,\dots,q_N\}, y_0) $, such that $ \gamma_i $ is a loop that goes in the upper half plane to $ q_i $, then performs a counter-clockwise twist around $ q_i $, and finally returns to $ y_0 $ in the upper half plane.
The calculation of the monodromy action is then separated into a local calculation around $ q_i $ and a conjugation corresponding to replacing the basepoint from some point that lies close to $ q_i $, to point $ y_0 $.
The data of the local relations and conjugations corresponding to all the singular points $ q_i $ can be represented in a monodromy table. See \cite{AGT_order6_conic_line_arr} for more details regarding this process.

The generators of $ \pi_1(\pr^{-1}(y_0)-\mathcal{B}, x_0) $ we use are all represented by loops that enclose a unique point of $ \mathcal{B}\cap\pr^{-1}(y_0) $.
Thus, every irreducible component $ X\in Irr(\mathcal{B}) $ corresponds to a $ \deg X $ generator.
For a line $ L\in Irr(\mathcal{B}) $ we denote the unique generator by $ \Gamma_L $.
For a conic $ C\in Irr(\mathcal{B}) $ we have two corresponding generators: $ \Gamma_{C} $ and $ \Gamma_{C'} $.

The local relation around a branch point of a component $ X $ gives the relation $ \Gamma_X=\Gamma_{X'} $; around a node between components $ X $ and $ Y $ it gives $ \bigl[\Gamma_X, \Gamma_Y \bigr]=e $ (as usual, square brackets denote the commutator); around a tangency point between components $ X $ and $ Y $ we have the relation $  \Gamma_X\Gamma_Y\Gamma_X\Gamma_Y\Gamma_X^{-1}\Gamma_Y^{-1}\Gamma_X^{-1}\Gamma_Y^{-1} = e $, and we denote it as $ \bigl\{\Gamma_X,\Gamma_Y\bigr\}=e $.

\bigskip

In the following lemmas we give the fundamental groups $ \pcpt{\mathcal{B}_1} $ and $ \pcpt{\mathcal{B}_2} $.
These presentations were obtained by implementing the Zariski-van Kampen algorithm as described above, the full computations appear in \cite[Appendix A,B]{arxiv_version}.
The publicly available package {\tt Sirocco} \cite{sirocco} in {\tt Sagemath} can also be used to verify the results.

\begin{lemma}\label{lemma-fundamental1}
Group $ \pcpt{\mathcal{B}_1} $ is generated by   $\Gamma_{C_1}, \Gamma_{C_2}, \Gamma_{C_3} , \Gamma_{L_2}$ and has the following relations:
\begin{equation}\label{yes1}
\begin{split}
\bigl\{\Gamma_{C_1} , \Gamma_{C_2}  \bigr\}=e,
\end{split}
\end{equation}
\begin{equation}\label{}
\begin{split}
\bigl\{\Gamma_{C_1} , \Gamma_{C_3}  \bigr\}=e,
\end{split}
\end{equation}
\begin{equation}\label{}
\begin{split}
\bigl\{\Gamma_{C_2} , \Gamma_{L_2}  \bigr\}=e,
\end{split}
\end{equation}
\begin{equation}\label{}
\begin{split}
\bigl\{\Gamma_{C_1}\Gamma_{C_3}\Gamma_{C_1}^{-1} , \Gamma_{L_2}  \bigr\}=e,
\end{split}
\end{equation}
\begin{equation}\label{}
\begin{split}
\bigl[\Gamma_{C_2} , \Gamma_{C_3}  \bigr]=e,
\end{split}
\end{equation}
\begin{equation}\label{}
\begin{split}
\bigl[\Gamma_{C_2}^{-1}\Gamma_{C_1}\Gamma_{C_2} , \Gamma_{L_2}  \bigr]=e,
\end{split}
\end{equation}
\begin{equation}\label{}
\begin{split}
\bigl[\Gamma_{C_3}^{-1}\Gamma_{C_1}\Gamma_{C_3} , \Gamma_{L_2}  \bigr]=e,
\end{split}
\end{equation}
\begin{equation}\label{}
\begin{split}
\bigl[\Gamma_{C_1}\Gamma_{C_3}\Gamma_{C_2}\Gamma_{L_2}\Gamma_{C_1}\Gamma_{C_3}\Gamma_{C_2} , \Gamma_{L_2}  \bigr]=e,
\end{split}
\end{equation}
\begin{equation}\label{}
\begin{split}
\bigl[\Gamma_{L_2}\Gamma_{C_1}\Gamma_{C_3}\Gamma_{C_1}^{-1}\Gamma_{L_2}^{-1} , \Gamma_{C_2}  \bigr]=e,
\end{split}
\end{equation}
\begin{equation}\label{}
\begin{split}
\bigl[\Gamma_{L_2}\Gamma_{C_1}\Gamma_{C_2}\Gamma_{C_3}\Gamma_{C_1}\Gamma_{C_3}\Gamma_{L_2}\Gamma_{C_1}
\Gamma_{C_3}\Gamma_{C_1}^{-1} , \Gamma_{C_2}  \bigr]=e,
\end{split}
\end{equation}
\begin{equation}\label{yes11}
\begin{split}
\bigl[\Gamma_{L_2}\Gamma_{C_1}\Gamma_{C_3}\Gamma_{C_2}\Gamma_{C_1}\Gamma_{C_2}\Gamma_{L_2}\Gamma_{C_1}
\Gamma_{C_2}\Gamma_{C_1}^{-1} , \Gamma_{C_3}  \bigr]=e,
\end{split}
\end{equation}
\begin{equation}\label{rel1}
\begin{split} 
& \Gamma_{L_2}^{-1}(\Gamma_{C_1}\Gamma_{C_3}\Gamma_{C_2})^{-1}\Gamma_{L_2}^{-1}
(\Gamma_{C_1}\Gamma_{C_3}\Gamma_{C_2})^{-1}\Gamma_{C_1}(\Gamma_{C_1}\Gamma_{C_3}\Gamma_{C_2})\Gamma_{L_2}
(\Gamma_{C_1}\Gamma_{C_3}\Gamma_{C_2})\Gamma_{L_2}= \\&=\Gamma_{C_1}\Gamma_{C_3}^{-1}\Gamma_{C_2}\Gamma_{C_1}\Gamma_{C_2}^{-1}\Gamma_{C_3}\Gamma_{C_1}^{-1},
\end{split}
\end{equation}
\begin{equation}\label{rel2}
\begin{split}
\Gamma_{C_1}\Gamma_{C_2}\Gamma_{C_1}^{-1}  =\Gamma_{L_2}^{-1}(\Gamma_{C_1}\Gamma_{C_3}\Gamma_{C_2})^{-1}\Gamma_{L_2}^{-1}\Gamma_{C_2}\Gamma_{L_2}
(\Gamma_{C_1}\Gamma_{C_3}\Gamma_{C_2})\Gamma_{L_2},
\end{split}
\end{equation}
\begin{equation}\label{rel3}
\begin{split}
\Gamma_{C_1}\Gamma_{C_3}\Gamma_{C_1}^{-1}  = \Gamma_{L_2}^{-1}(\Gamma_{C_1}\Gamma_{C_3}\Gamma_{C_2})^{-1}\Gamma_{L_2}^{-1}\Gamma_{C_3}\Gamma_{L_2}
(\Gamma_{C_1}\Gamma_{C_3}\Gamma_{C_2})\Gamma_{L_2}.
\end{split}
\end{equation}
\end{lemma}

\begin{lemma}\label{lemma-fundamental2}
Group $ \pcpt{\mathcal{B}_2} $ is generated by   $\Gamma_{C_1}, \Gamma_{C_2}, \Gamma_{C_3} , \Gamma_{L_3}$ and has the following relations:
\begin{equation}\label{eq:fund1_repr1}
\begin{split}
\Gamma_{C_1}\Gamma_{C_3}\Gamma_{C_2}\Gamma_{C_1}\Gamma_{C_2}^{-1}\Gamma_{C_3}^{-1}
\Gamma_{C_1}^{-1}=\Gamma_{L_3} \Gamma_{C_1}  \Gamma_{C_3}^2 \Gamma_{C_1} \Gamma_{C_3}^{-2} \Gamma_{C_1}^{-1} \Gamma_{L_3}^{-1},
\end{split}
\end{equation}
\begin{equation}\label{eq:fund1_repr2}
\begin{split}
\bigl\{\Gamma_{C_1} , \Gamma_{C_3}  \bigr\}=e,
\end{split}
\end{equation}
\begin{equation}\label{eq:fund1_repr3}
\begin{split}
\bigl\{\Gamma_{C_1} , \Gamma_{C_2}  \bigr\}=e,
\end{split}
\end{equation}
\begin{equation}\label{eq:fund1_repr4}
\begin{split}
\bigl\{\Gamma_{C_2} , \Gamma_{L_3}  \bigr\}=e,
\end{split}
\end{equation}
\begin{equation}\label{eq:fund1_repr5}
\begin{split}
\bigl\{\Gamma_{C_3} , \Gamma_{L_3}  \bigr\}=e,
\end{split}
\end{equation}
\begin{equation}\label{eq:fund1_repr6}
\begin{split}
\bigl[\Gamma_{C_2} , \Gamma_{C_3}  \bigr]=e,
\end{split}
\end{equation}
\begin{equation}\label{eq:fund1_repr7}
\begin{split}
\bigl[\Gamma_{C_2} , \Gamma_{C_1}\Gamma_{C_3} \Gamma_{C_1}^{-1} \bigr]=e,
\end{split}
\end{equation}
\begin{equation}\label{eq:fund1_repr8}
\begin{split}
\bigl[\Gamma_{C_1} , \Gamma_{L_3}  \bigr]=e,
\end{split}
\end{equation}
\begin{equation}\label{eq:fund1_repr9}
\begin{split}
\bigl[\Gamma_{C_2}^{-1}\Gamma_{C_1}\Gamma_{C_2}, \Gamma_{L_3}  \bigr]=e,
\end{split}
\end{equation}
\begin{equation}\label{eq:fund1_repr10}
\begin{split}
\bigl[\Gamma_{C_3}^{-1}\Gamma_{C_1}\Gamma_{C_3}, \Gamma_{L_3}  \bigr]=e,
\end{split}
\end{equation}
\begin{equation}\label{eq:fund1_repr11}
\begin{split}
\bigl[\Gamma_{L_3}^{-1}\Gamma_{C_3}\Gamma_{L_3}, \Gamma_{C_1}\Gamma_{C_2} \Gamma_{C_1}^{-1} \bigr]=e,
\end{split}
\end{equation}
\begin{equation}\label{eq:fund1_repr12}
\begin{split}
\bigl[\Gamma_{L_3}^{-1}\Gamma_{C_3}\Gamma_{L_3}, \Gamma_{C_2}  \bigr]=e,
\end{split}
\end{equation}
\begin{equation}\label{eq:fund1_repr13}
\begin{split}
\bigl[\Gamma_{C_3}\Gamma_{C_1}\Gamma_{C_3}^{-1}, \Gamma_{C_2} \Gamma_{C_1}\Gamma_{C_2} \Gamma_{C_3}\Gamma_{L_3}\Gamma_{C_1}\Gamma_{C_3} \bigr]=e,
\end{split}
\end{equation}
\begin{equation}\label{eq:fund1_repr14}
\begin{split}
\bigl[\Gamma_{L_3},\Gamma_{C_3}^{-1}\Gamma_{C_2} \Gamma_{C_1} \Gamma_{C_3}\Gamma_{C_2} \Gamma_{L_3}\Gamma_{C_1}\Gamma_{C_3}\Gamma_{C_3}\bigr]=e.
\end{split}
\end{equation}
\end{lemma}

\section{Square quotients and Coxeter groups}\label{COX}

To show that $ \mathcal{B}_1,\mathcal{B}_2 $ form a Zariski pair, we use the groups $\pcpt{\mathcal{B}_1} $ and $\pcpt{\mathcal{B}_2} $, which we have computed in Lemmas \ref{lemma-fundamental1} and \ref{lemma-fundamental2}.
We succeed in this task by showing that there is no isomorphism between $\pcpt{\mathcal{B}_1} $ and $\pcpt{\mathcal{B}_2} $ that agrees with a bijection $ \varphi_{\text{Irr}} $, as in Definition \ref{def:combinatorial_type}.
If $ X $ is an irreducible component of $ \mathcal{B}_1 $ and $ \Upsilon_X $ is any loop around it, then the image of $\Upsilon_X$ must be a loop around $ \varphi_{\text{Irr}}(X) $; meaning that it can be conjugate to either the chosen generator $ \Gamma_{\varphi_{\text{Irr}}(X)} $ or its inverse $ \Gamma_{\varphi_{\text{Irr}}(X)}^{-1} $.

We define the \emph{square quotients}
$$G_1 := \pcpt{\mathcal{B}_1} / \langle \Gamma_X^2 \; | \; X\in Irr( \mathcal{B}_1) \rangle \ \ \  \text{and} \ \ \ G_2 := \pcpt{\mathcal{B}_2} / \langle \Gamma_X^2 \; | \; X \in Irr( \mathcal{B}_2) \rangle.$$
In addition to relieving us from the need to distinguish between $ \Gamma_{\varphi_{\text{Irr}}(X)} $ and $ \Gamma_{\varphi_{\text{Irr}}(X)}^{-1} $ when using the above argument, this allows us to utilize the rich theory of Coxeter groups (for the definition of Coxeter groups and an introduction to the subject, see, for example, \cite{Bjorner2005}).

The Coxeter group that arises in the current work is $\widetilde{S}_2^C$ which is a Coxeter group of type $ \widetilde{C}_2 $. We use it to prove, via $G_1$ and $G_2$,  that $\mathcal{B}_1$ and $\mathcal{B}_2$ form a Zariski pair.
The elements of $\widetilde{S}_2^C$ are $5$-tuples $(\sigma; \epsilon_1, n_1; \epsilon_2, n_2)$ where $\sigma\in S_2$ is a permutation, $\epsilon_1, \epsilon_2\in \{\pm 1\}$ are signs, and $n_1,n_2\in \mathbb{Z}$ are integers.
The composition law is 
$$ (\sigma; \epsilon_1, n_1; \epsilon_2, n_2)\cdot (\sigma'; \epsilon_1', n_1'; \epsilon_2', n_2') = (\sigma\circ\sigma'; \epsilon_1'\cdot \epsilon_{\sigma'(1)}, n_1'+\epsilon_1'\cdot n_{\sigma'(1)}; \epsilon_2'\cdot \epsilon_{\sigma'(2)}, n_2'+\epsilon_2'\cdot n_{\sigma'(2)}).$$
We denote the elements $(\sigma; 1, n_1; 1, n_2)$, $(\sigma; -1, n_1; 1, n_2)$, $(\sigma; 1, n_1; -1, n_2)$, and $(\sigma; -1, n_1; -1, n_2)$) by 
$$  \begin{pmatrix}
1 & 2 \\
\sigma(1)^{ n_1} & \sigma(2)^{n_2}
\end{pmatrix},  \begin{pmatrix}
1 & 2 \\
\overline{\sigma(1)}^{n_1} & \sigma(2)^{n_2}
\end{pmatrix} , \begin{pmatrix}
1 & 2 \\
\sigma(1)^{n_1} & \overline{\sigma(2)}^{n_2}
\end{pmatrix}, \text{ and } \begin{pmatrix}
1 & 2 \\
\overline{\sigma(1)}^{n_1} & \overline{\sigma(2)}^{n_2}
\end{pmatrix}. $$
This group is generated by the Coxeter system
$$ \alpha:= \begin{pmatrix} 1 & 2 \\ \overline{1} & 2 \end{pmatrix} \; ; \; \beta := \begin{pmatrix} 1 & 2 \\ 2 & 1 \end{pmatrix} \; ; \; \gamma := \begin{pmatrix} 1 & 2 \\ 1 & \overline{2}^1 \end{pmatrix}. $$
The relations between those generators (as one can easily be convinced) are
\begin{equation}\label{eq:c2_tilde_def_sq}
	\alpha^2 = \beta^2 = \gamma^2 = e,
\end{equation}
\begin{equation}\label{eq:c2_tilde_def_comm}
	\bigl[\alpha, \gamma \bigr] = e,
\end{equation}
\begin{equation}\label{eq:c2_tilde_def_quad}
	\bigl\{\alpha, \beta \bigr\} = \bigl\{ \beta, \gamma \bigr\} = e.
\end{equation}
For a detailed description of Coxeter groups, see \cite{Bjorner2005}, where Coxeter groups of type $ \widetilde{C}_2 $ are examined in Subsection 8.4.

Another group that appears later in the proof is $D_8$, the dihedral group of order 8, which is also a Coxeter group (of type $B_2$) and has the following presentation:
\begin{equation}\label{eq:representation_dihedral}
	D_8 = \bigl\langle s_1,s_2 \ | \ s_1^2=s_2^2=e, (s_1s_2)^4=e \bigr\rangle.
\end{equation}

\begin{remark}\label{rem:d8_subgroup_c2_tilde}
	Note that $D_8$ is isomorphic to the subgroup $\{(\sigma; \epsilon_1, 0; \epsilon_2, 0) \mid \sigma\in S_2, \epsilon_1, \epsilon_2\in \{\pm 1\}\}$ of the Coxeter group $\widetilde{S}_2^C$ by the map 
	\begin{equation}\label{eq:d8_subgroup_c2_tilde}
		 s_1\mapsto \begin{pmatrix} 1 & 2 \\ 2 & 1 \end{pmatrix} \; ; \; s_2\mapsto \begin{pmatrix} 1 & 2 \\ \overline{1} & 2 \end{pmatrix}. 
	\end{equation}

	This allows us to construct a natural homomorphism $D_8\to S_2$ by taking the composition with the natural homomorphism $\widetilde{S}_2^C\to S_2$.
	The permutation of two elements corresponding to an element of $D_8$ can also be described as the permutation of the diagonals of the square induced by the usual action of $D_8$ on this square.
\end{remark}
\bigskip

In the following propositions, we provide the presentations for $G_1$ and $G_2$.

\begin{lemma}\label{lemma:diheral_representation}
	The group generated by $\Gamma_{C_1}$, $\Gamma_{C_2}$, and $\Gamma_{C_3}$ subject to the relations 
	\begin{align*}
		\Gamma_{C_1}^2 = \Gamma_{C_2}^2 = \Gamma_{C_3}^2 = e,\\
		\bigl\{\Gamma_{C_1}, \Gamma_{C_2}\bigr\} = \bigl\{\Gamma_{C_1}, \Gamma_{C_3}\bigr\} = e,\\
		\bigl[\Gamma_{C_2}, \Gamma_{C_3}\bigr] = \bigl[\Gamma_{C_1}, \Gamma_{C_2}\Gamma_{C_3}\bigr]=e,
	\end{align*}
	
	is isomorphic to $D_8 \times \mathbb{Z}_2$, where $D_8$ is the dihedral group of order $8$.
\end{lemma}

\begin{proof}
	Denote the group in question by $G$.
	By the relations \eqref{eq:c2_tilde_def_sq}-\eqref{eq:c2_tilde_def_quad}, the map from the Coxeter group of type $\tilde{C}_2$ to $G$ given by 
	\begin{equation*}
		\begin{pmatrix} 1 & 2 \\ 2 & 1 \end{pmatrix} \mapsto \Gamma_{C_1} \; ; \; \begin{pmatrix} 1 & 2 \\ \overline{1} & 2 \end{pmatrix} \mapsto \Gamma_{C_2} \; ; \; \begin{pmatrix} 1 & 2 \\ 1 & \overline{2}^1 \end{pmatrix} \mapsto \Gamma_{C_3},
	\end{equation*}
	is well defined, and obviously it is surjective.
	This means that $G$ is the quotient of the Coxeter group of type $\tilde{C}_2$ by a preimage of the relation $\bigl[\Gamma_{C_1}, \Gamma_{C_2}\Gamma_{C_3}\bigr]$ which can be computed to be 
	\begin{equation*}
		\bigl[\Gamma_{C_1}, \Gamma_{C_2}\Gamma_{C_3}\bigr] = \left[ \begin{pmatrix}1 & 2 \\ 2 & 1 \end{pmatrix}, \begin{pmatrix}1 & 2 \\ \bar{1} & \bar{2}^1\end{pmatrix} \right] = \begin{pmatrix}
			1 & 2 \\ 1^1 & 2^{-1}
		\end{pmatrix}.
	\end{equation*}
	Conjugating this relation by $\Gamma_{C_2}$ we see that $\begin{pmatrix*}1 & 2 \\ 1^1 & 2^1\end{pmatrix*}$ vanishes in $G$ as well.
	
	We can conclude from this that the image of $(\sigma; \epsilon_1, n_1; \epsilon_2, n_2)$ in $G$ is characterized by the triple $(\sigma, \epsilon_1, \epsilon_2)$ and the parity of the sum $n_1+n_2$.
	We note that the triple $(\sigma, \epsilon_1, \epsilon_2)$ is exactly an element of $D_8$, for example, considering the representation \eqref{eq:representation_dihedral} with the generators $\Gamma_{C_1}$ and $\Gamma_{C_2}$ as $s_1$ and $s_2$.
	Moreover, it is a straightforward verification that the product of two elements of $G$ can be computed by taking the product of the underlying signed permutations and summing the parities of the exponents, so we get the desired direct product.
\end{proof}

\begin{lemma}\label{lemma:dihedral_semidirect_representation}
	The group generated by $ \alpha $, $ \beta$, and $\gamma$ subject to the relations
	\begin{equation}\label{eq:lemma_4_2_statement1}
		\alpha^2 = \beta^2 = \gamma^2 = e,
	\end{equation}
	\begin{equation}\label{eq:lemma_4_2_statement2}
		\bigl[\alpha, \gamma\bigr] = e,
	\end{equation}
	\begin{equation}\label{eq:lemma_4_2_statement3}
		\bigl\{\alpha, \beta \bigr\} = \bigl\{\beta, \gamma\bigr\}  = e,
	\end{equation}
	\begin{equation}\label{eq:lemma_4_2_statement4}
		\bigl[\alpha, \beta\gamma\beta \bigr] = e,
	\end{equation}
	is isomorphic to the semi-direct product $ D_8 \ltimes \mathbb{Z}_2^2$, where $D_8$ is the dihedral group of order $8$ and the action on $\mathbb{Z}_2^2$ is given by permuting the coordinates according to the underlying permutation as in Remark \ref{rem:d8_subgroup_c2_tilde}.
\end{lemma}

\begin{proof}
	Let $G$ be the group in question.
	By relations \eqref{eq:c2_tilde_def_sq}-\eqref{eq:c2_tilde_def_quad} the map from the Coxeter group of type $\tilde{C}_2$ to $G$ given by 
	\begin{equation}\label{eq:dihedral_semi_map_from_coxeter}
		\begin{pmatrix} 1 & 2 \\ \overline{1} & 2 \end{pmatrix} \mapsto \alpha \; ; \; \begin{pmatrix} 1 & 2 \\ 2 & 1 \end{pmatrix} \mapsto \beta \; ; \; \begin{pmatrix} 1 & 2 \\ 1 & \overline{2}^1 \end{pmatrix} \mapsto \gamma,
	\end{equation}
	is well defined and surjective.
	The preimage of \eqref{eq:lemma_4_2_statement4} under this map is $\begin{pmatrix}1 & 2\\ 1 & 2^{-2}\end{pmatrix}$.
	We therefore conclude that an element of $G$ is defined by the signed permutation mapped to it by the map described by \eqref{eq:dihedral_semi_map_from_coxeter}, but the exponents only defined modulo $2$.
	Using this we can describe the short exact sequence 
	\begin{equation*}
		1 \to \mathbb{Z}_2^2 \to G \to D_8 \to 1
	\end{equation*}
	where the first nontrivial map is 
	\begin{equation*}
		\begin{pmatrix} x \\ y \end{pmatrix} \mapsto \begin{pmatrix}1 & 2 \\ 1^x & 2^y\end{pmatrix},
	\end{equation*}
	and the second map is (here we use the representation \eqref{eq:representation_dihedral} of $D_8$)
	\begin{equation*}
		\alpha\mapsto s_1 \; ; \; \beta\mapsto s_2 \; ; \; \gamma \mapsto s_2s_1s_2.
	\end{equation*}
	In fact, the latter map has a section $s_1\mapsto \alpha, \ s_2\mapsto \beta$, so $G\cong D_8\ltimes \mathbb{Z}_2^2$.
\end{proof}

\section{Zariski pairs}\label{sec:existing}

In this section, we state and prove the main result of the paper, i.e. that $\mathcal{B}_1$ and $\mathcal{B}_2$ form a Zariski pair (Theorem \ref{thm:tokunaga_pair3}), using the resulting square quotients $G_1$ and $G_2$. 
We will show that there is no isomorphism between $\pcpt{\mathcal{B}_1} $ and $\pcpt{\mathcal{B}_2} $ that agrees with a bijection $ \varphi_{\text{Irr}} $, as in Definition \ref{def:combinatorial_type}.
We determine both groups as extensions of Coxeter groups of different types, and use this fact to distinguish between them.
Thereafter, we show that the same techniques applied in the paper also work for Tokunaga's examples, which we call $\mathcal{B}_3, \mathcal{B}_4, \mathcal{B}_5$, and $\mathcal{B}_6$.

\begin{prop}\label{prop:G1_representation}
	The quotient $ G_1 := \pcpt{\mathcal{B}_1} / \langle \Gamma_X^2 \; | \; X\in Irr( \mathcal{B}_1) \rangle $ is isomorphic to $D_8 \ltimes \tilde{S}^C_2 $, where $D_8$ is the dihedral group of order $8$, $\tilde{S}^C_2$ is a Coxeter group of type $\tilde{C}_2$, and the action of the semidirect product is given by conjugating the element of $\widetilde{S}_2^C$ by the permutation in $S_2$ corresponding to the element of $D_8$ as in Remark \ref{rem:d8_subgroup_c2_tilde}.
\end{prop}

\begin{proof}
	By a straightforward simplification we can deduce that $G_1$ is generated by  $\Gamma_{C_1}, \Gamma_{C_2}, \Gamma_{C_3}$, and $ \Gamma_{L_2} $, subject to the following relations:
	\begin{equation}\label{}
		\Gamma_{C_1}^2 = \Gamma_{C_2}^2 = \Gamma_{C_3}^2 = \Gamma_{L_2}^2 =e,
	\end{equation}
	\begin{equation}\label{B5-new1}
		\begin{split}
			\bigl\{\Gamma_{C_1} , \Gamma_{C_2}  \bigr\}=e,
		\end{split}
	\end{equation}
	\begin{equation}\label{}
		\begin{split}
			\bigl\{\Gamma_{C_1} , \Gamma_{C_3}  \bigr\}=e,
		\end{split}
	\end{equation}
	\begin{equation}\label{B5-new2}
		\begin{split}
			\bigl\{\Gamma_{C_2} , \Gamma_{L_2}  \bigr\}=e,
		\end{split}
	\end{equation}
	\begin{equation}\label{}
		\begin{split}
			\bigl\{\Gamma_{C_3} , \Gamma_{L_2}  \bigr\}=e,
		\end{split}
	\end{equation}
	\begin{equation}\label{}
		\begin{split}
			\bigl\{\Gamma_{C_1}\Gamma_{C_3}\Gamma_{C_1} , \Gamma_{L_2}  \bigr\}=e,
		\end{split}
	\end{equation}
	\begin{equation}\label{}
		\begin{split}
			\bigl[\Gamma_{C_2} , \Gamma_{C_3}  \bigr]=e,
		\end{split}
	\end{equation}
	\begin{equation}\label{B5-new6}
		\begin{split}
			\bigl[\Gamma_{C_1},\Gamma_{C_2}\Gamma_{C_3}  \bigr]=e,
		\end{split}
	\end{equation}
	\begin{equation}\label{B5-new3}
		\begin{split}
			\bigl[\Gamma_{C_2}\Gamma_{C_1}\Gamma_{C_2} , \Gamma_{L_2}  \bigr]=e,
		\end{split}
	\end{equation}
	\begin{equation}\label{}
		\begin{split}
			\bigl[\Gamma_{C_3}\Gamma_{C_1}\Gamma_{C_3} , \Gamma_{L_2}  \bigr]=e,
		\end{split}
	\end{equation}
	\begin{equation}\label{B5-new5}
		\begin{split}
			\bigl[\Gamma_{C_1}\Gamma_{C_3}\Gamma_{C_1} , \Gamma_{L_2}\Gamma_{C_2} \Gamma_{L_2}  \bigr]=e,
		\end{split}
	\end{equation}
	\begin{equation}\label{B5-new4}
		\begin{split}
			\bigl[\Gamma_{C_1}\Gamma_{C_3}\Gamma_{C_2}\Gamma_{L_2}\Gamma_{C_1}\Gamma_{C_3}\Gamma_{C_2} , \Gamma_{L_2}  \bigr]=e.
		\end{split}
	\end{equation}
	
	Setting $\alpha:=\Gamma_{C_1}\Gamma_{C_2}\Gamma_{C_1}\Gamma_{C_3}$ and $\Gamma_{C_1'}:=\Gamma_{C_2}\Gamma_{C_1}\Gamma_{C_2}$ we can see, using Lemma \ref{lemma:diheral_representation}, that $G_1$ is generated by $\Gamma_{C_1'}$, $\Gamma_{C_2}$, $\Gamma_{L_2}$, and $\alpha$ subject to the relations
	\begin{equation}\label{eq:G_1_rep_simpl1}
		\Gamma_{C_1'}^2 = \Gamma_{C_2}^2 = \Gamma_{L_2}^2=\alpha^2=e,
	\end{equation}
	\begin{equation}\label{eq:G_1_rep_simpl2}
		\bigl\{ \Gamma_{C_2}, \Gamma_{L_2}\bigr\} = e,
	\end{equation}
	\begin{equation}\label{eq:G_1_rep_simpl3}
		\bigl[\alpha, \Gamma_{C_1'}\bigr] = \bigl[\alpha, \Gamma_{C_2}\bigr] = e,
	\end{equation}
	\begin{equation}\label{eq:G_1_rep_simpl4}
		\bigl[\alpha, \Gamma_{L_2}\Gamma_{C_2}\Gamma_{L_2}\bigr] = e,
	\end{equation}
	\begin{equation}\label{eq:G_1_rep_simpl5}
		\bigl\{ \alpha, \Gamma_{L_2} \bigr\} = e,
	\end{equation}
	\begin{equation}\label{eq:G_1_rep_simpl6}
		\bigl[ \Gamma_{C_1'}, \Gamma_{L_2} \bigr] = e,
	\end{equation}
	\begin{equation}\label{eq:G_1_rep_simpl7}
		\bigl\{ \Gamma_{C_1'}, \Gamma_{C_2} \bigr\} = e.
	\end{equation}
	
	We first consider the group $H$ generated by $\Gamma_{C_2}$, $\Gamma_{L_2}$, and $\alpha$ subject to relations \eqref{eq:G_1_rep_simpl1}-\eqref{eq:G_1_rep_simpl5}.
	By Lemma \ref{lemma:dihedral_semidirect_representation}, $H\cong D_8\ltimes \mathbb{Z}_2^2$, so we can represent any element of $H$ by an element of $D_8$ times a vector $\left(x, y\right)\in \mathbb{Z}_2^2$, where $\alpha=s_1$, $\Gamma_{L_2}=s_2$, and $\Gamma_{C_2}=s_2s_1s_2\left(0, 1\right)$, the element $s_1$ commutes with elements of the form $\left(x, y\right)$, and $\left(x, y\right)s_2 = s_2\left(y, x\right)$ holds ($s_1$ and $s_2$ are the elements of $D_8$ as in the representation \eqref{eq:representation_dihedral}).
	
	Substituting this into \eqref{eq:G_1_rep_simpl3}, \eqref{eq:G_1_rep_simpl6}, and \eqref{eq:G_1_rep_simpl7}, we get that $\Gamma_{C_1'}$ commutes with any elements of $D_8\subseteq H$, and $\bigl\{\Gamma_{C_1'}, (1,0) \bigr\}=\bigl\{\Gamma_{C_1'}, (0,1)\bigr\} = e$.
	This means that $\Gamma_{C_1'}$, $(1,0)$, and $(0,1)$ satisfy the relations \eqref{eq:c2_tilde_def_sq}-\eqref{eq:c2_tilde_def_quad} that conclude the proof.
\end{proof}

\begin{prop}\label{prop:G2_representation}
	The quotient $ G_2 := \pcpt{\mathcal{B}_2} / \langle \Gamma_X^2 \; | \; X \in Irr( \mathcal{B}_2) \rangle $ is isomorphic to $D_8 \ltimes \left( D_8\times \mathbb{Z}_2\right) $, where $D_8$ is the dihedral group of order $8$, the action of the semidirect product is given by 
	$$ g \centerdot (h, x) = (h, x + sgn(g) + sgn(h)),$$
	and $sgn:D_8\to \mathbb{Z}_2$ is the function that assigns to each element of $D_8$ the sign of the corresponding permutation (as in Remark \ref{rem:d8_subgroup_c2_tilde}). 
\end{prop}

\begin{proof}
	Considering relation \eqref{eq:fund1_repr1} together with \eqref{eq:fund1_repr8} and $\Gamma_{C_1}^2=\Gamma_{C_3}^2=e$, we get that $\bigl[\Gamma_{C_1}, \Gamma_{C_2}\Gamma_{C_3}\bigr]=e$ holds in $G_2$.
	By Lemma \ref{lemma:diheral_representation}, this together with \eqref{eq:fund1_repr2}, \eqref{eq:fund1_repr3}, and \eqref{eq:fund1_repr6}, indicate that the group generated by $\Gamma_{C_1}$, $\Gamma_{C_2}$, and $\Gamma_{C_3}$ in $G_2$ is a quotient of $D_8\times \mathbb{Z}_2$.
	Denote $\alpha:=\Gamma_{C_3}\Gamma_{C_1}\Gamma_{C_2}\Gamma_{C_1}$.
	The proof of Lemma \ref{lemma:diheral_representation} shows that $\alpha$ is of order $2$ and commutes with both $\Gamma_{C_1}$ and $\Gamma_{C_2}$.
	Plugging $\Gamma_{C_3}=\alpha\Gamma_{C_1}\Gamma_{C_2}\Gamma_{C_1}$ in \eqref{eq:fund1_repr1}-\eqref{eq:fund1_repr14} and setting $\Gamma_X^2=e$ for all $X\in Irr(\mathcal{B}_2)$ we get that $G_2$ is generated by $\Gamma_{C_1}$, $\Gamma_{C_2}$, $\Gamma_{L_3}$, and $\alpha$ subject to the following relations:
	\begin{equation}\label{eq:G2_rep_rel1}
		\Gamma_{C_1}^2=\Gamma_{C_2}^2=\Gamma_{L_3}^2=\alpha^2=e,
	\end{equation}
	\begin{equation}\label{eq:G2_rep_rel2}
		\bigl\{\Gamma_{C_1}, \Gamma_{C_2}\bigr\} = \bigl\{\Gamma_{C_2}, \Gamma_{L_3}\bigr\} = \bigl[\Gamma_{C_1}, \Gamma_{L_3}\bigr] = e,
	\end{equation}
	\begin{equation}\label{eq:G2_rep_rel3}
		\bigl[\Gamma_{C_2}\Gamma_{C_1}\Gamma_{C_2}, \Gamma_{L_3}\bigr] = e,
	\end{equation}
	\begin{equation}\label{eq:G2_rep_rel4}
		\bigl[\alpha, \Gamma_{C_1}\bigr] = \bigl[\alpha, \Gamma_{C_2}\bigr] = e,
	\end{equation}
	\begin{equation}\label{eq:G2_rep_rel5}
		\bigl[\alpha, \Gamma_{L_3}\Gamma_{C_2}\Gamma_{L_3}\bigr] = e,
	\end{equation}
	\begin{equation}\label{eq:G2_rep_rel6}
		\bigl\{\alpha\Gamma_{C_2}, \Gamma_{L_3}\bigr\} = e.
	\end{equation}
	
	We first consider the group $H$ generated by $\Gamma_{C_1}$, $\Gamma_{C_2}$, and $\Gamma_{L_3}$ subject to relations \eqref{eq:G2_rep_rel1}-\eqref{eq:G2_rep_rel3}.
	By Lemma \ref{lemma:dihedral_semidirect_representation}, $H\cong D_8\ltimes \mathbb{Z}_2^2$, so we can represent any element of $H$ by an element of $D_8$ times a vector $\left(x, y\right)\in \mathbb{Z}_2^2$, where $\Gamma_{C_1}=s_1$, $\Gamma_{C_2}=s_2$, and $\Gamma_{L_3}=s_2s_1s_2\left(0, 1\right)$, the element $s_1$ commutes with elements of the form $\left(x, y\right)$, and $\left(x, y\right)s_2 = s_2\left(y, x\right)$ holds.
	
	We now turn our attention to relations \eqref{eq:G2_rep_rel4}-\eqref{eq:G2_rep_rel6}.
	Substituting the representation described above for $\Gamma_{C_1}$, $\Gamma_{C_2}$, and $\Gamma_{L_3}$ we get that 
	\begin{equation*}
		\bigl[\alpha, \begin{pmatrix}1\\ 1\end{pmatrix}\bigr] = \bigl\{ \alpha, \begin{pmatrix}1 \\ 0\end{pmatrix} \bigr\} =e.
	\end{equation*}
	This means that the group generated by $\alpha$ and $\mathbb{Z}_2^2$ inside $G_2$ is isomorphic to $D_8\times \mathbb{Z}_2$.
	Combining this with the facts that $H\cong D_8\ltimes \mathbb{Z}_2^2$, and that $\alpha$ commutes with both $\Gamma_{C_1}$ and $\Gamma_{C_2}$, we get the desired result.
\end{proof}

We now state a theorem that describes our method.
For a curve $\mathcal{C}$, an irreducible component $X\in Irr( \mathcal{C})$, and an element $\Gamma\in \pcpt{\mathcal{C}}$, we will say that \emph{$\Gamma$ goes around $X$} if the image of $\Gamma$ in $\pcpt{\overline{\left( \mathcal{C} - X\right)}}$ is null homotopic.
\begin{thm}\label{thm:coxeter_non_isomorphic}
	Let $\mathcal{C}_1, \mathcal{C}_2 \subseteq \cpt$ be two curves with the same combinatorial type.
	For every irreducible component $X\in Irr( \mathcal{C}_1)$ pick an element $\Gamma_X\in \pcpt{\mathcal{C}_1}$ that goes around $X$, and similarly pick an element $\Gamma_Y\in \pcpt{\mathcal{C}_2}$ that goes around $Y$ for every irreducible component $Y\in Irr( \mathcal{C}_2)$.
	If the groups ${\pcpt{\mathcal{C}_1}}/{\left\langle \Gamma_X^2 | X\in Irr( \mathcal{C}_1) \right\rangle}$ and ${\pcpt{\mathcal{C}_2}}/{\left\langle \Gamma_Y^2 | Y\in Irr(\mathcal{C}_2) \right\rangle}$ are not isomorphic, then $\mathcal{C}_1$ and $\mathcal{C}_2$ form a Zariski pair.
\end{thm}

\begin{proof}
	
	Assume, towards a contradiction, that there exists a homeomorphism $\varphi : \left(\cpt, \mathcal{C}_1\right)\to \left(\cpt, \mathcal{C}_2\right)$.
	Then $\varphi $ sends a component of $\mathcal{C}_1$ to a component of $\mathcal{C}_2$, so the induced map on fundamental groups $\varphi_*:\pcpt{\mathcal{C}_1} \to \pcpt{\mathcal{C}_2}$ sends $\Gamma_{X}\in\pcpt{\mathcal{C}_1}$ to either a conjugation of some $\Gamma_{Y}\in\pcpt{\mathcal{C}_2}$ or to a conjugation of some $\Gamma_{Y}^{-1}\in\pcpt{\mathcal{C}_2}$.
	In particular, the image under $\varphi_*$ of the normal closure of the group $\left\langle \Gamma_X^2 \; \mid \; X\in Irr(\mathcal{C}_1) \right\rangle $ is the normal closure of $\left\langle \Gamma_Y^2 \; \mid \; Y\in Irr( \mathcal{C}_2) \right\rangle $, meaning that we get a well defined isomorphism ${\pcpt{\mathcal{C}_1}}/{\left\langle \Gamma_X^2 | X\in Irr(\mathcal{C}_1) \right\rangle} \xrightarrow{\sim} {\pcpt{\mathcal{C}_2}}/{\left\langle \Gamma_Y^2 | Y\in Irr(\mathcal{C}_2) \right\rangle}$, in contradiction to the assumption that they are not isomorphic.
\end{proof}

\begin{cor}\label{thm:tokunaga_pair3}
Curves $ \mathcal{B}_1 $ and $ \mathcal{B}_2 $ form a Zariski pair.
\end{cor}

\begin{proof}
	This is immediate by Theorem \ref{thm:coxeter_non_isomorphic} and the fact that Proposition \ref{prop:G1_representation} implies that $G_1$ is infinite (as all affine Coxeter groups are infinite), while Proposition \ref{prop:G2_representation} shows $G_2$ is finite.
	
\end{proof}
\begin{remark}
	The arrangements $\mathcal{B}_1$ and $\mathcal{B}_2$ can be considered as degenerations of conic arrangements considered by Namba and Tsuchihashi \cite{Namba2004OnTF}. We can see that
    eight tacnodes given by $(C_2 + C_3)\cap (C_1 + L_1 +L_2)$ are not
    on any conic, while those given by $(C_2 + C_3)\cap (C_1 + L_1 + L_3)$ are a conic (two lines). Hence likewise \cite[Lemma 4.20]{Bartolo_tokunaga_survey}, we see that for each $n \ge 3$ there exists a dihedral 
    cover $\pi_n : X_n \to \mathbb{P}^2$ such that
    (i) the Galois group is the dihedral group of order $2n$,
    (ii) the branch locus is $\mathcal{B}_2$, and (iii) the ramification index along $C_2+ C_3$ (resp. $C_1 + L_1 +L_3$)
     is $2$ (resp. $n$).  On the other hand, by the results in
     \cite{tokunaga2014}, e.g., Theorems 3.2 and 3.3, we see that
     for any odd prime $p$, there exist no dihedral covers as above for $\mathcal{B}_1$.
    
\end{remark}

\bigskip

In \cite{tokunaga2014}, Tokunaga studied dihedral covers naturally arising from
elliptic surfaces called elliptic dihedral covers. He  applied
the results to study the existence/non-existence of dihedral
covers branched along conic-line arrangements, and gave 
 two examples of Zariski pairs of conic-line arrangements of degree 7.
Here, we recall these two pairs of arrangements and denote them as $\mathcal{B}_3, \mathcal{B}_4$ and  $\mathcal{B}_5, \mathcal{B}_6$. In the following two examples we show that our method can prove the fact that those pairs of arrangements are Zariski pairs, see Theorems \ref{Zar1} and \ref{Zar2}. We leave calculations of the groups to the reader, they are long and we do not feel it is appropriate to plant them into the paper. However, full calculations of the groups appear in \cite{arxiv_version}.

\begin{eg}
Let  $\mathcal{B}_3$ and $\mathcal{B}_4$ be two conic-line arrangements of degree 7, see Figure \ref{fig_tokunaga2}.
\end{eg}

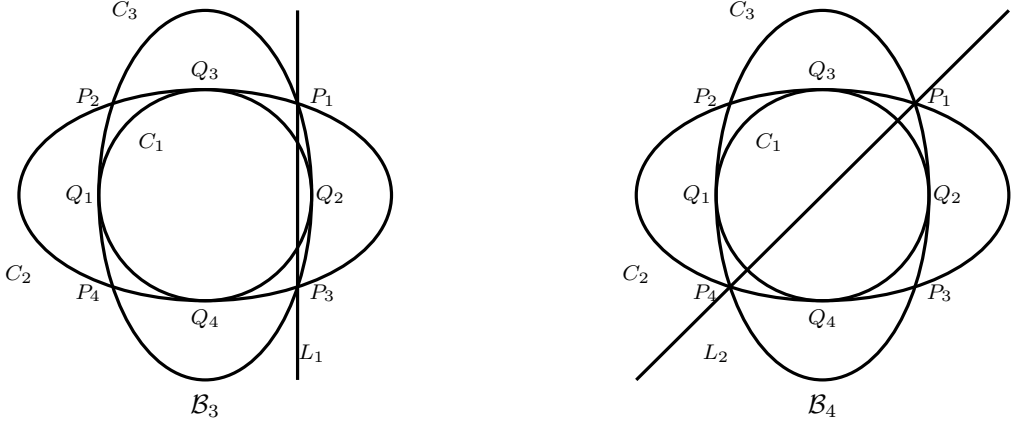
\begin{figure}[H]
\centering
\begin{tikzpicture} [scale=0.35, style = very thick]
	\begin{scope}[xshift=-330]
		\draw (0,0) ellipse (7 and 4);
		\draw (-7, -3) node {\footnotesize $C_2$};
		
		\draw (0,0) ellipse (4 and 7);
		\draw (-3, 7) node {\footnotesize $ C_3 $};
		
		\draw (0,0) ellipse (4 and 4);
		\draw (-2, 2) node {\footnotesize $ C_1 $};
		
		\draw (3.47, -7) -- (3.47, 7); 
		\draw (4, -6) node {\footnotesize $ L_1 $};
		
		\draw (-4.7, 0) node {\footnotesize $Q_1$};
		\draw (4.7, 0) node {\footnotesize $Q_2$};
		\draw (0, 4.7) node {\footnotesize $Q_3$};
		\draw (0, -4.7) node {\footnotesize $Q_4$};
		
		\draw (4.4, 3.7) node {\footnotesize $ P_1 $};
		\draw (-4.4, 3.7) node {\footnotesize $ P_2 $};
		\draw (4.4, -3.7) node {\footnotesize $ P_3 $};
		\draw (-4.4, -3.7) node {\footnotesize $ P_4 $};
		\draw (0, -8) node {$\mathcal{B}_3$};
	\end{scope}
	
	\begin{scope}[xshift=330]
		\draw (0,0) ellipse (7 and 4);
		\draw (-7, -3) node {\footnotesize $C_2$};
		
		\draw (0,0) ellipse (4 and 7);
		\draw (-3, 7) node {\footnotesize $ C_3 $};
		
		\draw (0,0) ellipse (4 and 4);
		\draw (-2, 2) node {\footnotesize $ C_1 $};
		
		\draw (-7, -7) -- (7, 7); 
		\draw (-4, -6) node {\footnotesize $ L_2 $};
		
		\draw (-4.7, 0) node {\footnotesize $Q_1$};
		\draw (4.7, 0) node {\footnotesize $Q_2$};
		\draw (0, 4.7) node {\footnotesize $Q_3$};
		\draw (0, -4.7) node {\footnotesize $Q_4$};
		
		\draw (4.4, 3.7) node {\footnotesize $ P_1 $};
		\draw (-4.4, 3.7) node {\footnotesize $ P_2 $};
		\draw (4.4, -3.7) node {\footnotesize $ P_3 $};
		\draw (-4.4, -3.7) node {\footnotesize $ P_4 $};
		
		\draw (0, -8) node {$\mathcal{B}_4$};
	\end{scope}
\end{tikzpicture}

\caption{The arrangements $ \mathcal{B}_3,\mathcal{B}_4 $.}
\label{fig_tokunaga2}
\end{figure}

Here is an explanation how to construct $ \mathcal{B}_3$ and $\mathcal{B}_4 $.
\begin{Construction}\label{construcetion_tokunaga_arrangement2_construction}
Let $ C_i \left(i=1,2,3\right)$ be smooth conics and let $ L_j \left(j=1,2\right) $ be lines as follows:
\begin{itemize}
	\item $ C_2 $ and $ C_3 $ meet transversely. We put $ C_2\cap C_3=\{ P_1, P_2, P_3, P_4 \} $.
	\item $ C_1 $ is tangent to both $ C_3 $ and $ C_2 $ such that the intersection multiplicities at intersection points are all equal to 2. We set $ C_1\cap C_2=\{ Q_3, Q_4 \} $ and $ C_1\cap C_3 = \{Q_1, Q_2 \} $.
	\item $ L_1 $ passes through $ P_1 $ and $ P_3 $.
	\item $ L_2 $ passes through $ P_1 $ and $ P_4 $.
	\item $ L_1, L_2 $ meet $ C_1 $ transversely.
\end{itemize}
Denote $ \mathcal{B}_3=C_1\cup C_2\cup C_3\cup L_1 $ and $ \mathcal{B}_4=C_1\cup C_2\cup C_3 \cup L_2 $.
\end{Construction}

Using the Zariski-van Kampen Theorem, many simplifications of relations, and known structures of Coxeter groups, we can determine the square quotients $G_3$ and $G_4$.
\begin{prop}
Group $ G_3 := \pcpt{\mathcal{B}_3}/\langle \Gamma_X^2 \; | \; X\in Irr(\mathcal{B}_3) \rangle  $ is isomorphic to $ D_8 \ltimes \mbb{Z}_2^2$ where the action of the semidirect product is the same as in Lemma \ref{lemma:dihedral_semidirect_representation}, and group $G_4 := \pcpt{\mathcal{B}_4}/\langle \Gamma_X^2 \; | \; X \in Irr(\mathcal{B}_4) \rangle $ is a Coxeter group of type $\widetilde{C}_2$.
\end{prop}
And the proposition gives us the following conclusion.
\begin{thm}\label{Zar1}
$ \mathcal{B}_3 $ and $ \mathcal{B}_4 $ form a Zariski pair.
\end{thm}

 \bigskip

\begin{eg}
Let $ \mathcal{B}_5$ and $\mathcal{B}_6 $ be a pair of conic-line arrangements of degree 7, see Figure \ref{fig_tokunaga1}.
\end{eg}

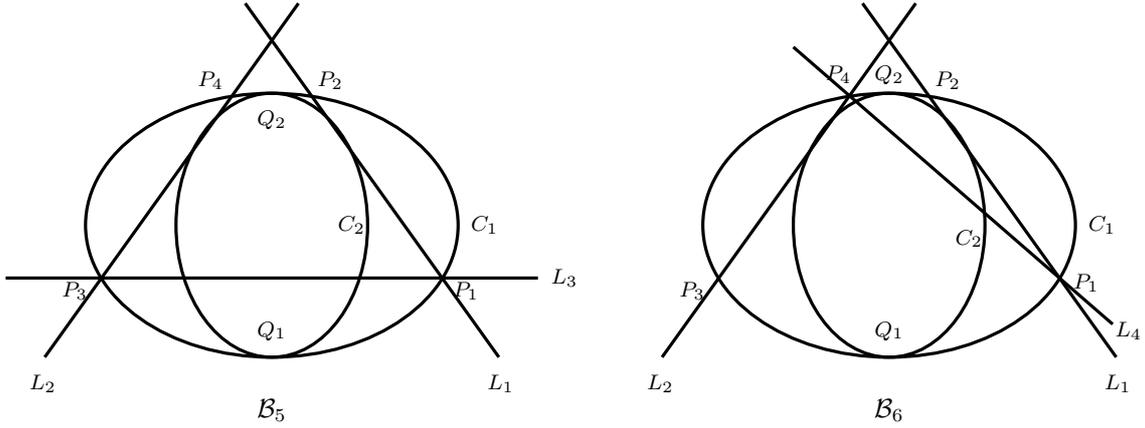
\begin{figure}[H]
\centering
\begin{tikzpicture} [scale=0.35, style = very thick]
	\begin{scope}[xshift=-330]
		\draw (0,0) ellipse (7 and 5);
		\draw (8,0) node {\footnotesize $C_1$};
		\draw (0,0) ellipse (3.6 and 5);
		\draw (3,0) node {\footnotesize $C_2$};
		\draw (-10, -2) -- (10, -2);
		\draw (11, -2) node {\footnotesize $ L_3 $};
		\draw (-1,8.406) -- (8.533, -5);
		\draw (8.6, -6) node {\footnotesize $ L_1 $};
		\draw (1,8.406) -- (-8.533, -5);
		\draw (-8.6, -6) node {\footnotesize $ L_2 $};
		
		\draw (0, -4) node {\footnotesize $ Q_1 $};
		\draw (0, 4) node {\footnotesize $ Q_2 $};
		
		\draw (7.3, -2.5) node {\footnotesize $ P_1 $};
		\draw (-7.4, -2.5) node {\footnotesize $ P_3 $};
		\draw (2.2, 5.5) node {\footnotesize $ P_2 $};
		\draw (-2.3, 5.5) node {\footnotesize $ P_4 $};
		
		\draw (0, -7) node {$\mathcal{B}_5$};
	\end{scope}
	
	\begin{scope}[xshift=330]
		\draw (0,0) ellipse (7 and 5);
		\draw (8,0) node {\footnotesize $C_1$};
		\draw (0,0) ellipse (3.6 and 5);
		\draw (3,-0.5) node {\footnotesize $C_2$};
		\draw (-1,8.406) -- (8.533, -5);
		\draw (8.6, -6) node {\footnotesize $ L_1 $};
		\draw (1,8.406) -- (-8.533, -5);
		\draw (-8.6, -6) node {\footnotesize $ L_2 $};
		\draw (8.4, -3.75) -- (-3.6, 6.75);
		\draw (9, -4) node {\footnotesize $ L_4 $};
		
		\draw (0, -4) node {\footnotesize $ Q_1 $};
		\draw (0, 5.7) node {\footnotesize $ Q_2 $};
		
		\draw (7.4, -2.2) node {\footnotesize $ P_1 $};
		\draw (-7.4, -2.5) node {\footnotesize $ P_3 $};
		\draw (2.2, 5.5) node {\footnotesize $ P_2 $};
		\draw (-1.9, 5.7) node {\footnotesize $ P_4 $};
		
		\draw (0, -7) node {$\mathcal{B}_6$};
	\end{scope}
	
\end{tikzpicture}
\caption{The arrangements $ \mathcal{B}_5,\mathcal{B}_6 $.}
\label{fig_tokunaga1}
\end{figure}

We construct the arrangements in the following way.
\begin{Construction}\label{construcetion_tokunaga_arrangement1_construction}
Let $ C_i \left(i=1,2\right)$ be smooth conics and let $ L_j \left(j=1,\dots,4\right) $ be lines as follows:
\begin{itemize}
	\item $ C_1 $ and $ C_2 $ are tangent at $ 2 $ distinct points $ Q_1 $ and $ Q_2 $.
	\item $ C_2 $ is tangent to $ L_1 $ and to $ L_2 $.
	\item $ L_1 $ and $ L_2 $ meet $ C_1 $ transversely. We put $ C_1\cap L_1=\left\{P_1, P_2\right\}, \; C_1\cap L_2=\left\{P_3,P_4\right\} $.
	\item $ L_3 $ is the line connecting $ P_1 $ and $ P_3 $, $ L_4 $ is the line connecting $ P_1 $ and $ P_4 $.
	\item Both $ L_3 $ and $ L_4 $ meet $ C_2 $ transversely.
\end{itemize}
Denote $ \mathcal{B}_5=C_1\cup C_2\cup L_1\cup L_2\cup L_3 $ and $ \mathcal{B}_6=C_1\cup C_2\cup L_1\cup L_2 \cup L_4 $.
\end{Construction}

In the following proposition we write down the presentations that we derive for the square quotients $G_5$ and $G_6$.
\begin{prop}
The quotient $ G_5 := \pcpt{\mathcal{B}_5} / \langle \Gamma_X^2 \; | \; X\in \mathcal{B}_5 \rangle $ is generated by $ \Gamma_{C_1}, \Gamma_{C_2}, \Gamma_{L_1}, \Gamma_{L_2} $, and $ \Gamma_{L_3} $, subject to the following relations:
\begin{equation}\label{arrangement1_final_squares}
	\Gamma_{C_1}^2 = \Gamma_{C_2}^2 = \Gamma_{L_1}^2 = \Gamma_{L_2}^2 = \Gamma_{L_3}^2 = e,
\end{equation}
\begin{equation}\label{arrangement1_final1}
	\bigl[\Gamma_{L_3},\Gamma_{L_1}\bigr] = \bigl[\Gamma_{L_3}, \Gamma_{L_2}\bigr] = \bigl[\Gamma_{L_3}, \Gamma_{C_1}\bigr] = \bigl[\Gamma_{L_3}, \Gamma_{C_2}\bigr] = \bigl[\Gamma_{L_1}, \Gamma_{L_2}\bigr] = \bigl[\Gamma_{C_1}, \Gamma_{L_1}\bigr] = \bigl[\Gamma_{C_1}, \Gamma_{L_2}\bigr] = e,
\end{equation}
\begin{equation}\label{arrangement1_final2}
	\bigl[\Gamma_{L_1}, \Gamma_{C_2}\Gamma_{C_1}\Gamma_{C_2}\bigr] = \bigl[\Gamma_{L_2}, \Gamma_{C_2}\Gamma_{C_1}\Gamma_{C_2}\bigr] = \bigl[\Gamma_{C_2}, \Gamma_{L_1}\Gamma_{L_2}\bigr] = e,
\end{equation}
\begin{equation}\label{arrangement1_final3}
	\bigl\{\Gamma_{C_1}, \Gamma_{C_2}\bigr\} = \bigl\{ \Gamma_{C_2}, \Gamma_{L_1} \bigr\} = \bigl\{ \Gamma_{C_2}, \Gamma_{L_2}\bigr\} = e,
\end{equation}
\begin{equation}\label{arrangement1_final4}
	\Gamma_{C_1}\Gamma_{C_2}\Gamma_{C_1}\Gamma_{L_1}\Gamma_{C_2}\Gamma_{L_3}\Gamma_{L_2} = e.
\end{equation}
In particular, $ \Gamma_{L_3} $ is in the center of $ G_5 $.

Group $ G_6 := \pcpt{\mathcal{B}_6} / \langle \Gamma_{X}^2 \; | \; X\in Irr(\mathcal{B}_6)  \rangle $ is generated by $ \Gamma_{C_1},\Gamma_{C_2}, \Gamma_{L_1} $, and $ \Gamma_{L_2} $, subject to the following relations:
	\begin{equation}\label{arrangement2_second_final_1}
		\Gamma_{C_1}^2 = \Gamma_{C_2}^2 = \Gamma_{L_1}^2 = \Gamma_{L_2}^2 = e,
	\end{equation}
	\begin{equation}\label{arrangement2_second_final_2}
		\bigl[\Gamma_{L_1}, \Gamma_{L_2} \bigr] = \bigl[\Gamma_{C_1}, \Gamma_{L_1}\bigr] = \bigl[\Gamma_{C_1}, \Gamma_{L_2}\bigr] = \bigl[\Gamma_{L_1}\Gamma_{L_2}, \Gamma_{C_2}\bigr] = e,
	\end{equation}
	\begin{equation}\label{arrangement2_second_final_3}
		\bigl\{ \Gamma_{C_1}, \Gamma_{C_2} \bigr\} = \bigl\{ \Gamma_{C_2}, \Gamma_{L_1} \bigr\} = \bigl\{ \Gamma_{C_2}, \Gamma_{L_2} \bigr\} = e,
	\end{equation}
	\begin{equation}\label{arrangement2_second_final_4}
		\bigl[\Gamma_{C_2}, \Gamma_{L_1}\Gamma_{C_1}\Gamma_{C_2}\Gamma_{C_1}\Gamma_{L_1}\bigr] = e.
	\end{equation}
Group $ G_6 $ is isomorphic to the direct product $ G_6'\times \mbb{Z}_2 $, where $ G_6' $ is generated by $ \Gamma_{C_1},\Gamma_{C_2} $, and $ \Gamma_{L_1} $, subject to the following relations:
	\begin{equation}\label{arrangement2_without_L2_1}
		\Gamma_{C_1}^2 = \Gamma_{C_2}^2 = \Gamma_{L_1}^2 = e,
	\end{equation}
	\begin{equation}\label{arrangement2_without_L2_2}
		\bigl[\Gamma_{C_1}, \Gamma_{L_1}\bigr] = e,
	\end{equation}
	\begin{equation}\label{arrangement2_without_L2_3}
		\bigl\{ \Gamma_{C_1}, \Gamma_{C_2} \bigr\} = \bigl\{ \Gamma_{C_2}, \Gamma_{L_1} \bigr\} = e,
	\end{equation}
	\begin{equation}\label{arrangement2_without_L2_4}
		\bigl[\Gamma_{C_2}, \Gamma_{L_1}\Gamma_{C_1}\Gamma_{C_2}\Gamma_{C_1}\Gamma_{L_1}\bigr] = e.
	\end{equation}		
The element $ \Gamma_{L_4} = \Gamma_{C_2}\Gamma_{C_1}\Gamma_{L_1}\Gamma_{C_2}\Gamma_{C_1}\Gamma_{L_2} \in G_6 $ is not contained in the center of $ G_6 $.
\end{prop}

And the conclusion is:

\begin{thm}\label{Zar2}
$\mathcal{B}_5$ and $\mathcal{B}_6$ form a Zariski pair.
\end{thm}
\begin{remark}{The fact that $(\mathcal{B}_3, \mathcal{B}_4)$ and 
$(\mathcal{B}_5, \mathcal{B}_6)$ are Zariski pairs implies that
it is impossible to deform $\mathcal{B}_3$ (resp. $\mathcal{B}_5$)
 to $\mathcal{B}_4$ (resp. $\mathcal{B}_6$) with keeping the combinatorics. In \cite{bannai-tokunaga-yorisaki}, the realization
 spaces of these two combinatorics were studied. Here, the
 realization space means the set  of curves with fixed combinatorics.
 It is shown that there exist exactly two connected components for
 each case.
}    
\end{remark}

\newpage
\appendix

\section{Explanation about the calculations}\label{sec:appendix_explanation}
In the following appendices we compute a presentation of $ \pcpt{\mathcal{B}_i} $ for $ i=1,\dots,6 $, using the van Kampen algorithm.
The calculation is separated into a local study around a singularity $ q_i $ and a conjugation that can be written as a composition of conjugations (i.e.,  diffeomorphisms) that correspond to each $ q_j $ with $ j<i $.
This data is represented in a \emph{monodromy table}, where each row corresponds to one of the points $ q_i $ and consists of the number $ i $, a verbal description of the arrangement in the neighborhood of $ q_i $, the relation arising from a loop that goes around $ q_i $, and the diffeomorphism.
Instead of the relation, we write only the \emph{skeleton}, that is, we only denote the positions of the points in the fiber that are affected by the monodromy action.
The exact monodromy action can either be seen from the description of the singularity (one twist for branch points, double twists for nodes, and fourfold twists for tangencies), or, in the case of the complicated singular points in $ \mathcal{B}_1 $ and $ \mathcal{B}_2 $, the exact monodromy action is computed in \cite{AGT_order6_conic_line_arr}.
We use the following notation for diffeomorphisms - $ \Delta\langle a,a+1,\dots,a+k \rangle $ denotes a half-twist of the points $ a,\dots,a+k $ around a common point.
The diffeomorphism $ \Delta^{1/2}\langle a \rangle $ changes the points $ a $ and $ a+1 $ from real to imaginary or vice versa, where the exact action is indicated by a subscript.

To compute the monodromy action of a loop in $ \pi_1(\mbb{C}^1 - \{q_1,\dots,q_N\}) $, we take the local computation around $ q_i $ and conjugate it first, by the diffeomorphism corresponding to $ q_{i-1} $, then by the one corresponding to $ q_{i-2} $, and so on, until $ q_1 $.
The intermediate computations are represented by a braid connecting the affected points in the fiber.
The real points in the fiber are drawn as black circles and ordered by their appearance on the real axis.
The imaginary points in the fiber are depicted by red squares and ordered (from left to right) by the increasing value of the imaginary part.
The place in the sequence of points where the red squares (i.e., the imaginary points) appear, can be thought of as the origin (the intersection of the real and imaginary axis).

In each of the following appendices, we give figures of the curves $\mathcal{B}_i$ for $i=1,\dots,6$, the monodromy tables, relations in $\pcpt{\mathcal{B}_i} $, simplifications and proofs of Zariski pairs with non-isomorphic fundamental groups.

\newpage

\section{New Zariski pair with non-isomorphic fundamental groups}

In this section we provide the monodromy tables and computations of braids and relations for $ \pi_1(\mathbb{CP}^2\setminus \mathcal{B}_1) $ and $ \pi_1(\mathbb{CP}^2\setminus \mathcal{B}_2) $.
It also includes the coordinates we used to perform the calculations.

\subsection{Affine coordinates}
We pick affine coordinates $ x= \frac{2.4150635094611 T + 1.68301270189222 X}{0.707106781186548 T + 0.707106781186548 X + 2 Z} $ and $ y = \frac{- 0.816987298107781 T + 1.9150635094611 X}{0.707106781186548 T + 0.707106781186548 X + 2 Z} $, so $ T = 0.319177251576849 x - 0.280502116982037 y $, $ X = 0.13616454968463 x + 0.402510584910183 y $, and $ Z = - 0.160987637714845 x - 0.0431365075170868 y + 0.5 $.
The equations become:
$$ C_{1} = V\left( 0.0944978830179 x^{2} - 0.0833333333333 x y + 0.160987637714 x + 0.23883545031 y^{2} + 0.0431365075170 y - 0.25 \right) $$
$$ C_{2} = V\left( 0.400120236790 x^{2} - 0.620512701892 x y + 0.160987637714 x + 0.474879763209 y^{2} + 0.0431365075170 y - 0.25 \right) $$
$$ C_{3} = V\left( 0.150120236790 x^{2} + 0.245512701892 x y + 0.160987637714 x + 0.724879763209 y^{2} + 0.0431365075170 y - 0.25 \right) $$
$$ L_{1} = V\left( 0.00302305115846 x - 0.731240783452 y + 0.559016994374 \right) $$
$$ L_{2} = V\left( 0.363002352625 x - 0.634784620332 y - 0.559016994374 \right) $$
$$ L_{3} = V\left( 0.275352150527729 x + 0.0737803863680778 y + 0.559016994374 \right) $$

\subsection{Computation of $ \pcpt{\mathcal{B}_1} $}
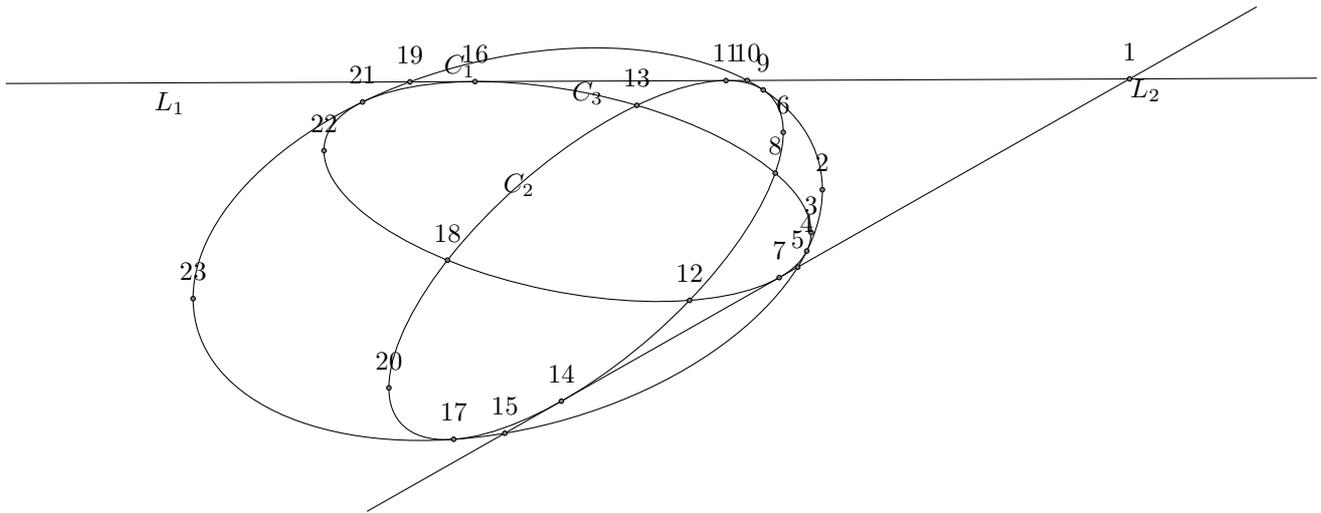
\begin{figure}[H]
\begin{center}
\begin{scaletikzpicturetowidth}{\textwidth}

\end{scaletikzpicturetowidth}
\end{center}
\end{figure}

\subsection{Monodromy table}

\begin{tabular}{| c | c | c | c |}
\hline \hline
Vertex number& Vertex description& Skeleton& Diffeomorphism\\ \hline \hline
1& Node between $L_2$ and $L_1$& $ \langle 1 - 2 \rangle $& $ \Delta\langle 1, 2 \rangle $ \\ \hline
2& $C_1$ branch& $ \langle 2 - 3 \rangle $& $ \Delta_{I_{4}I_{6}}^{1/2}\langle 2 \rangle $ \\ \hline
3& $C_3$ branch& $ \langle 3 - 4 \rangle $& $ \Delta_{I_{2}I_{4}}^{1/2}\langle 3 \rangle $ \\ \hline
4& Tangency between $C_3$ and $C_1$& $ \langle 2 - 3 \rangle $& $ \Delta^{2}\langle 2, 3 \rangle $ \\ \hline
5& Node between $C_1$ and $L_2$& $ \langle 1 - 2 \rangle $& $ \Delta\langle 1, 2 \rangle $ \\ \hline
6& $C_2$ branch& $ \langle 5 - 6 \rangle $& $ \Delta_{\mathbb{R}I_{2}}^{1/2}\langle 5 \rangle $ \\ \hline
7& Tangency between $C_3$ and $L_2$& $ \langle 2 - 3 \rangle $& $ \Delta^{2}\langle 2, 3 \rangle $ \\
\hline \hline
\end{tabular}

{Braids calculation}

\subsubsection{Vertex Number 1}

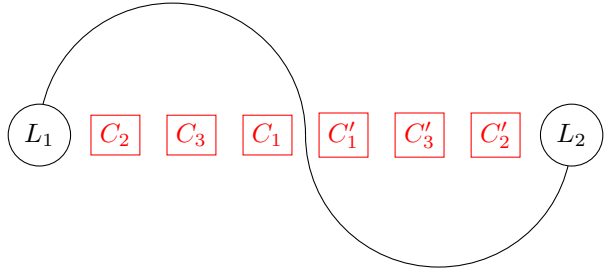
\begin{figure}[H]
\begin{center}
%
\end{center}
\end{figure}

Relation:

\begin{equation*}\begin{split}
\left[\Gamma_{C_1}\Gamma_{C_3}\Gamma_{C_2}\Gamma_{L_1}\Gamma_{C_2}^{-1}\Gamma_{C_3}^{-1}\Gamma_{C_1}^{-1} , \Gamma_{L_2} \right]=e
\end{split}
\end{equation*}
\subsubsection{Vertex Number 2}

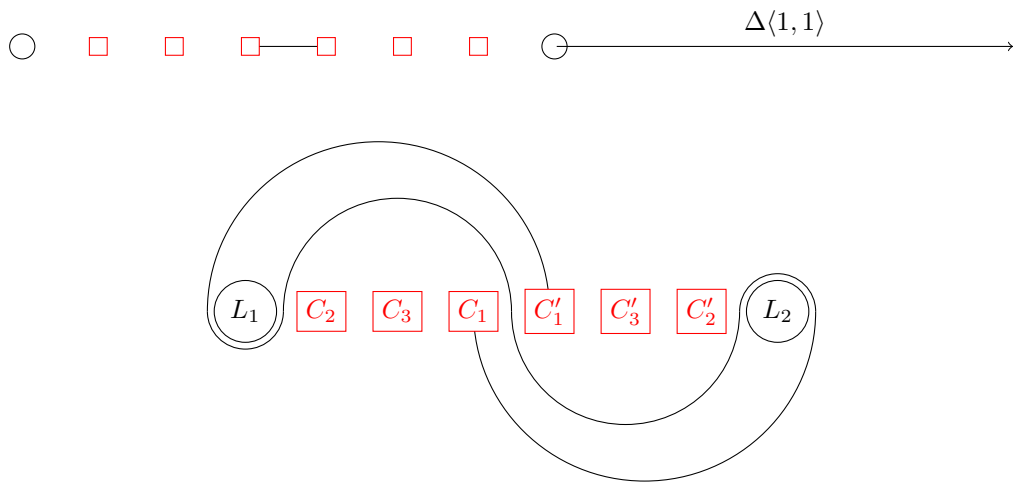
\begin{figure}[H]
\begin{center}
%
\end{center}
\end{figure}

Relation:

\begin{equation*}\begin{split}
\Gamma_{C_1}\Gamma_{C_3}\Gamma_{C_2}\Gamma_{L_1}\Gamma_{C_2}^{-1}\Gamma_{C_3}^{-1}\Gamma_{C_1}^{-1}\Gamma_{L_2}^{-1}\Gamma_{C_1}\Gamma_{L_2}\Gamma_{C_1}\Gamma_{C_3}\Gamma_{C_2}\Gamma_{L_1}^{-1}\Gamma_{C_2}^{-1}\Gamma_{C_3}^{-1}\Gamma_{C_1}^{-1}  = \Gamma_{C_1'}
\end{split}
\end{equation*}
\subsubsection{Vertex Number 3}

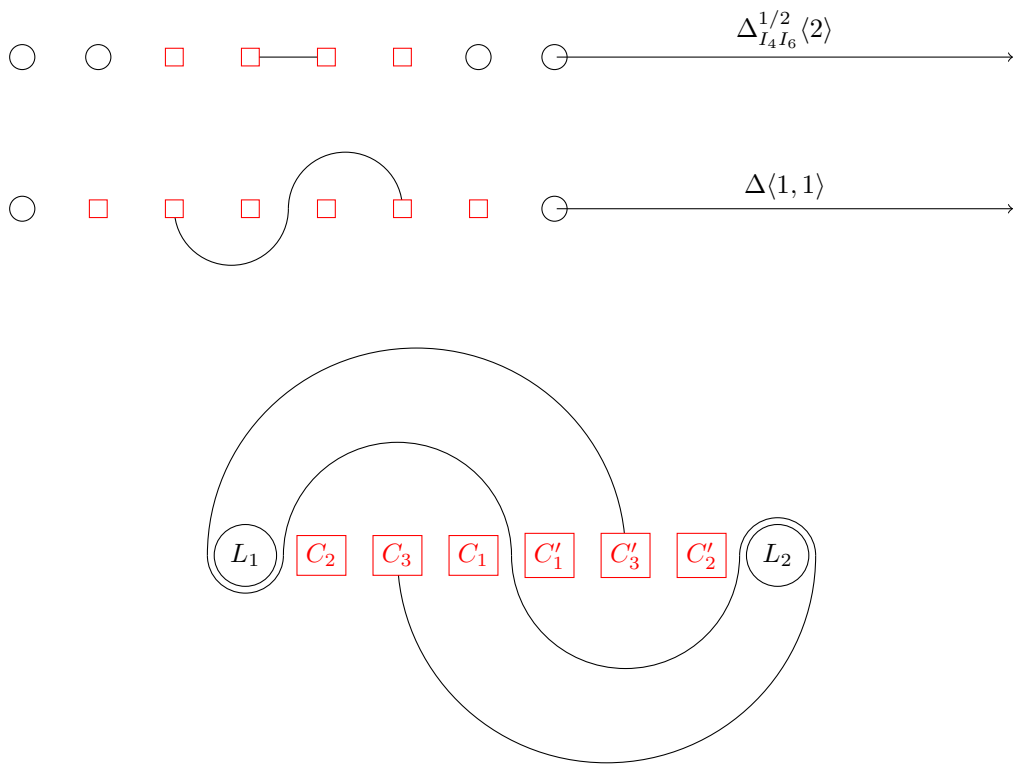
\begin{figure}[H]
\begin{center}
%
\end{center}
\end{figure}

Relation:

\begin{equation*}\begin{split}
\Gamma_{C_1'}\Gamma_{C_1}\Gamma_{C_3}\Gamma_{C_2}\Gamma_{L_1}\Gamma_{C_2}^{-1}\Gamma_{C_3}^{-1}\Gamma_{C_1}^{-1}\Gamma_{L_2}^{-1}\Gamma_{C_3}\Gamma_{L_2}\Gamma_{C_1}\Gamma_{C_3}\Gamma_{C_2}\Gamma_{L_1}^{-1}\Gamma_{C_2}^{-1}\Gamma_{C_3}^{-1}\Gamma_{C_1}^{-1}\Gamma_{C_1'}^{-1}  = \Gamma_{C_3'}
\end{split}
\end{equation*}
\subsubsection{Vertex Number 4}

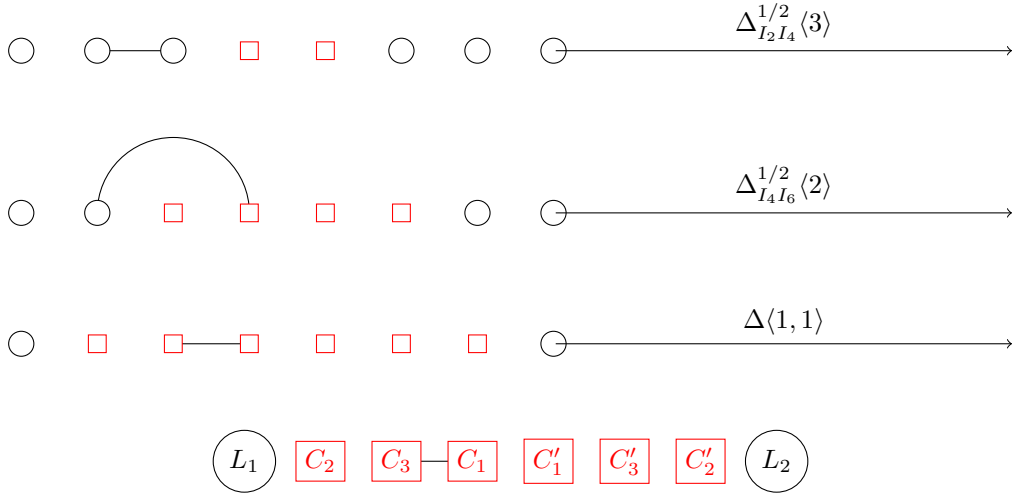
\begin{figure}[H]
\begin{center}
%
\end{center}
\end{figure}

Relation:

\begin{equation*}\begin{split}
\left\{\Gamma_{C_1}  , \Gamma_{C_3} \right\}=e
\end{split}
\end{equation*}
\subsubsection{Vertex Number 5}

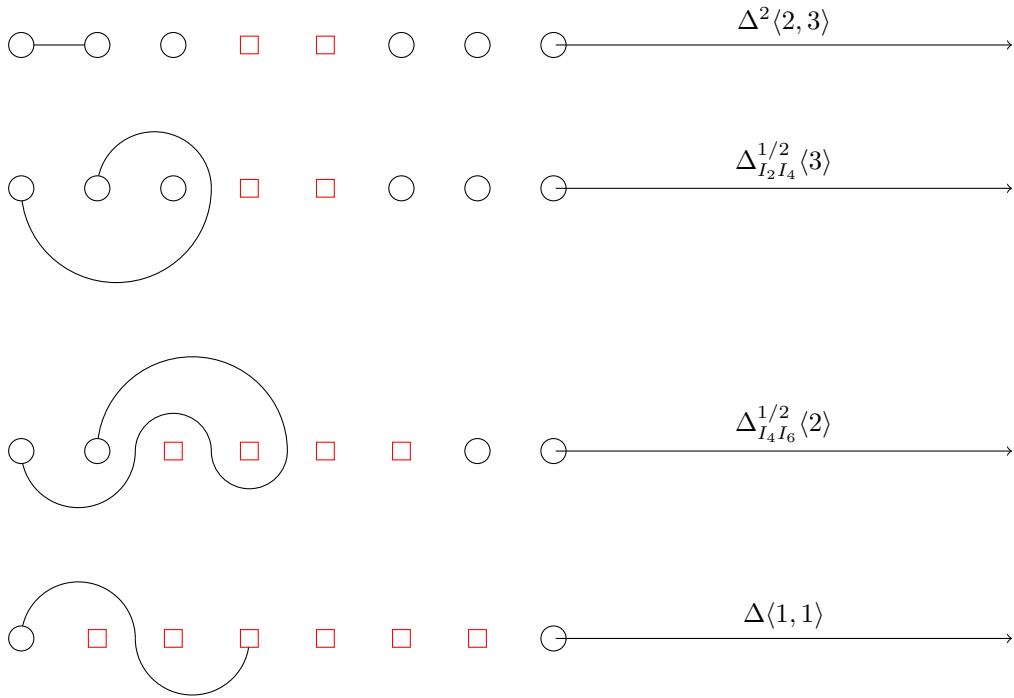
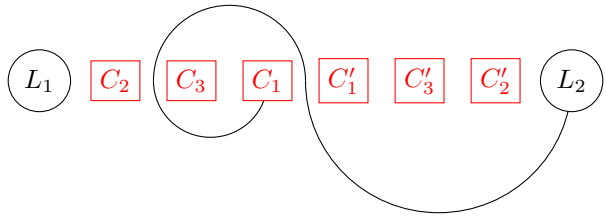
\begin{figure}[H]
\begin{center}
%
\end{center}
\end{figure}

Relation:

\begin{equation*}\begin{split}
\left[\Gamma_{C_3}^{-1}\Gamma_{C_1}^{-1}\Gamma_{L_2}\Gamma_{C_1}\Gamma_{C_3} , \Gamma_{C_1} \right]=e
\end{split}
\end{equation*}
\subsubsection{Vertex Number 6}

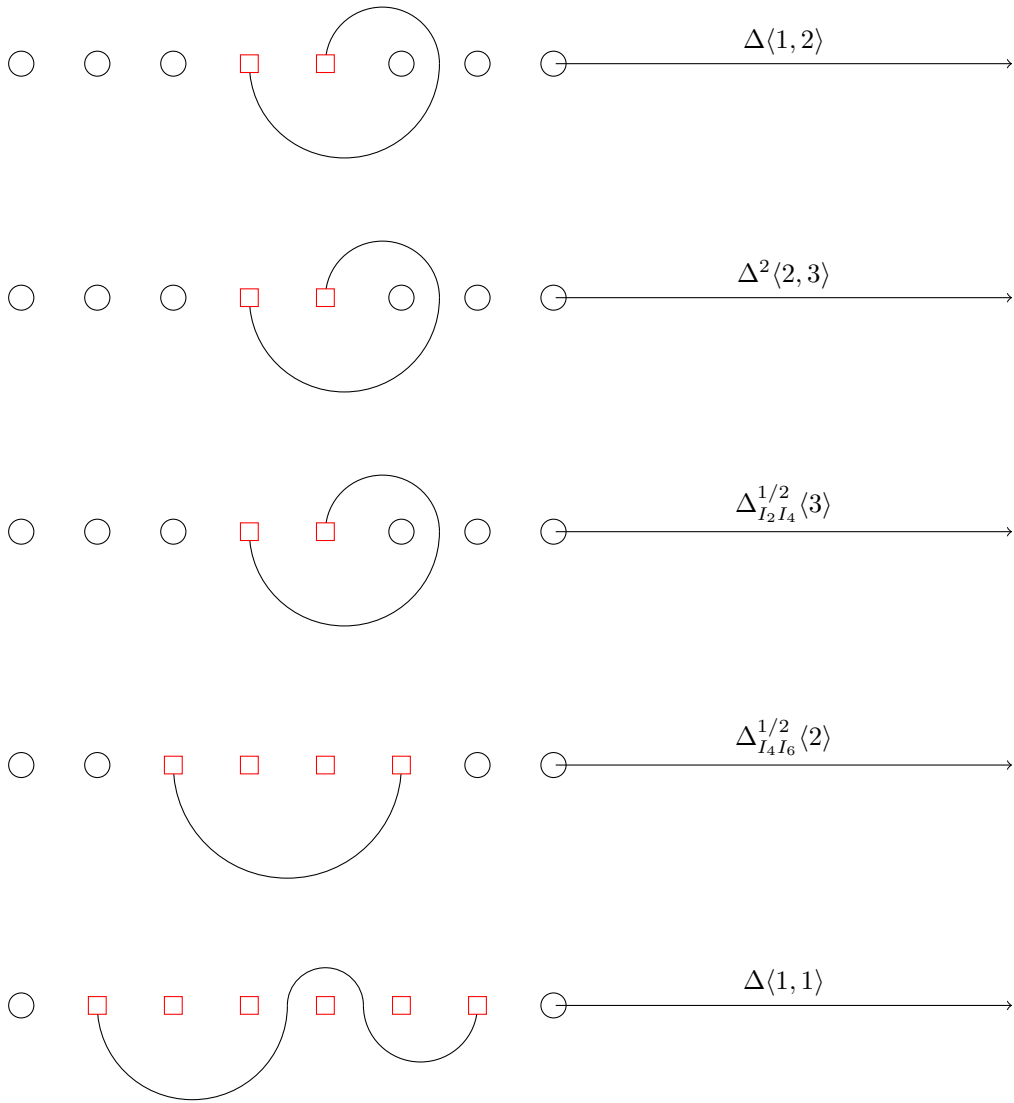
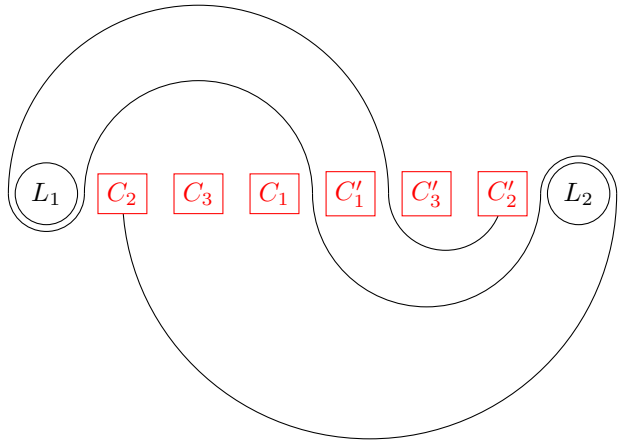
\begin{figure}[H]
\begin{center}
%
\end{center}
\end{figure}

Relation:

\begin{equation*}\begin{split}
\Gamma_{C_1'}\Gamma_{C_1}\Gamma_{C_3}\Gamma_{C_2}\Gamma_{L_1}\Gamma_{C_2}^{-1}\Gamma_{C_3}^{-1}\Gamma_{C_1}^{-1}\Gamma_{L_2}^{-1}\Gamma_{C_2}\Gamma_{L_2}\Gamma_{C_1}\Gamma_{C_3}\Gamma_{C_2}\Gamma_{L_1}^{-1}\Gamma_{C_2}^{-1}\Gamma_{C_3}^{-1}\Gamma_{C_1}^{-1}\Gamma_{C_1'}^{-1}  = \Gamma_{C_2'}
\end{split}
\end{equation*}
\subsubsection{Vertex Number 7}

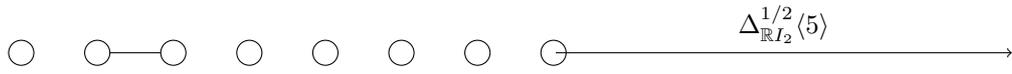
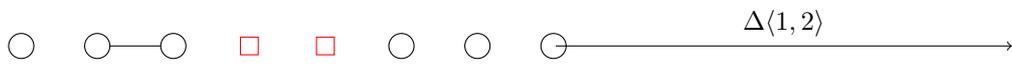
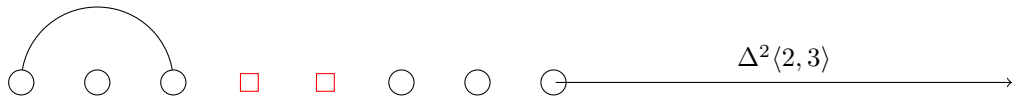
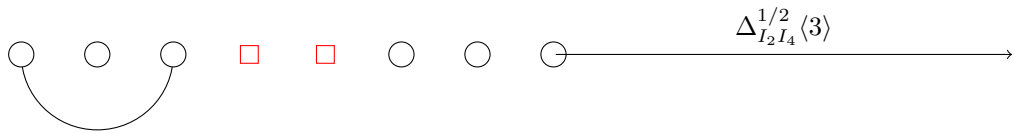
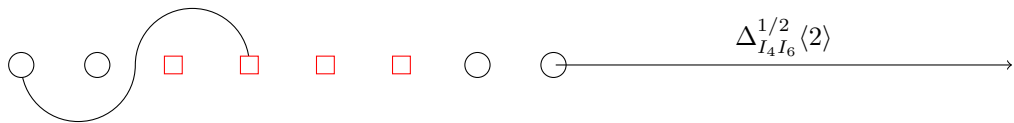
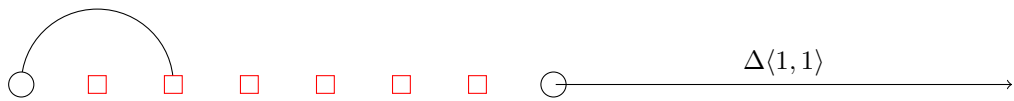
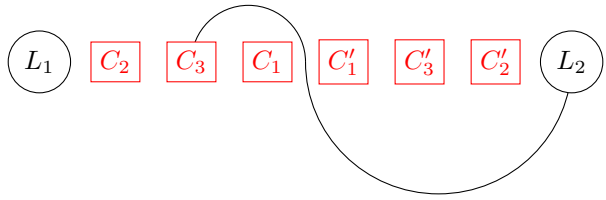
\begin{figure}[H]
\begin{center}
%
\end{center}
\end{figure}

Relation:

\begin{equation*}\begin{split}
\left\{\Gamma_{C_1}^{-1}\Gamma_{L_2}\Gamma_{C_1}  , \Gamma_{C_3} \right\}=e
\end{split}
\end{equation*}
\subsubsection{Vertex Number 8}

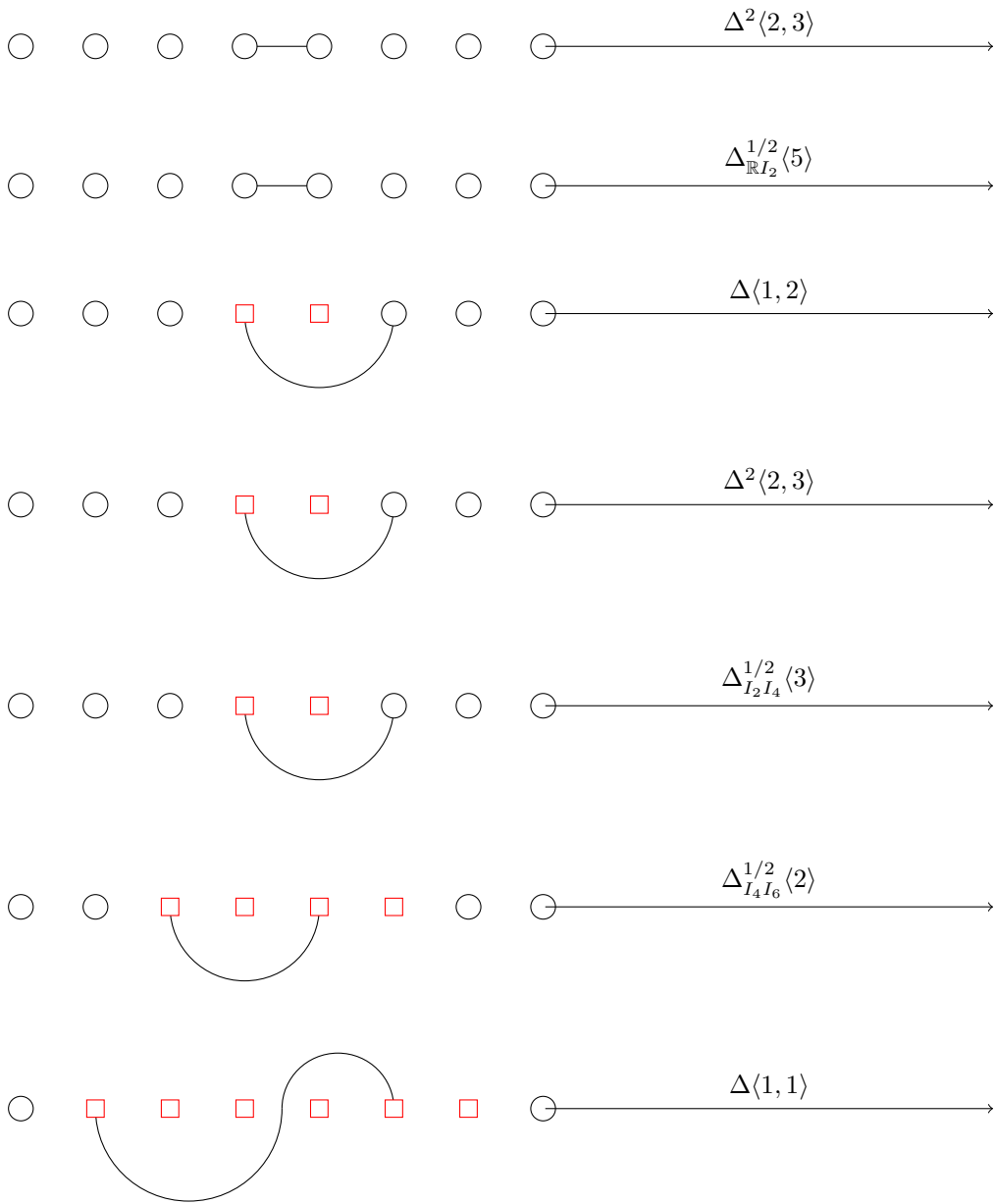
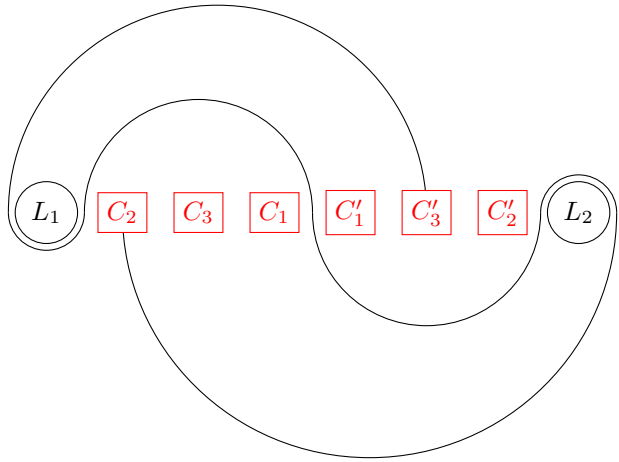
\begin{figure}[H]
\begin{center}
%
\end{center}
\end{figure}

Relation:

\begin{equation*}\begin{split}
\left[\Gamma_{L_2}\Gamma_{C_1}\Gamma_{C_3}\Gamma_{C_2}\Gamma_{L_1}^{-1}\Gamma_{C_2}^{-1}\Gamma_{C_3}^{-1}\Gamma_{C_1}^{-1}\Gamma_{C_1'}^{-1}\Gamma_{C_3'}\Gamma_{C_1'}\Gamma_{C_1}\Gamma_{C_3}\Gamma_{C_2}\Gamma_{L_1}\Gamma_{C_2}^{-1}\Gamma_{C_3}^{-1}\Gamma_{C_1}^{-1}\Gamma_{L_2}^{-1} , \Gamma_{C_2} \right]=e
\end{split}
\end{equation*}
\subsubsection{Vertex Number 9}

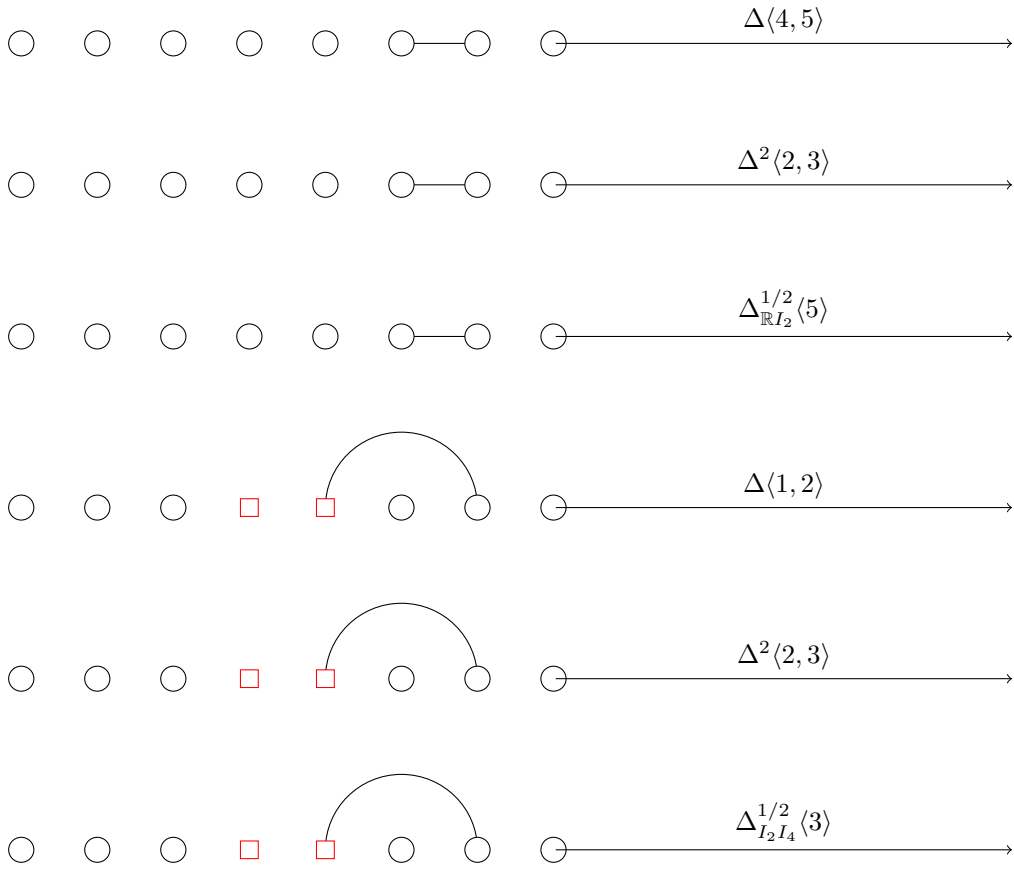
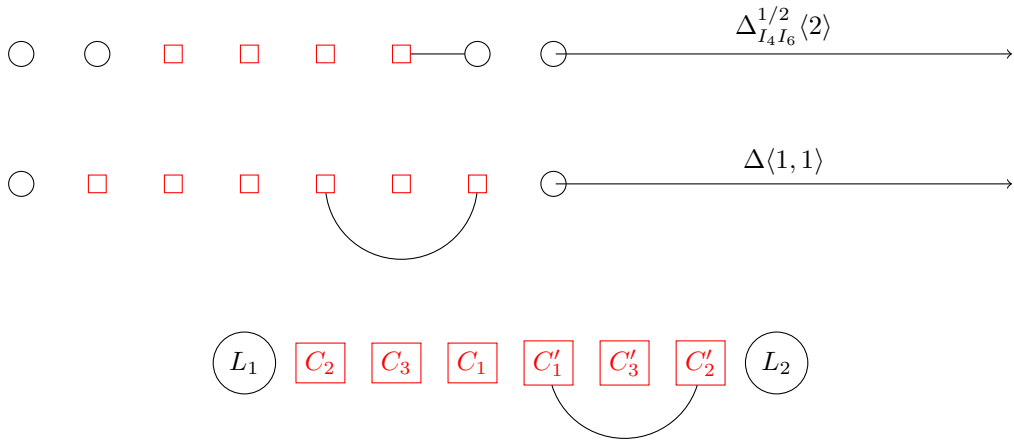
\begin{figure}[H]
\begin{center}
%
\end{center}
\end{figure}

Relation:

\begin{equation*}\begin{split}
\left\{\Gamma_{C_2'}  , \Gamma_{C_1'} \right\}=e
\end{split}
\end{equation*}
\subsubsection{Vertex Number 10}

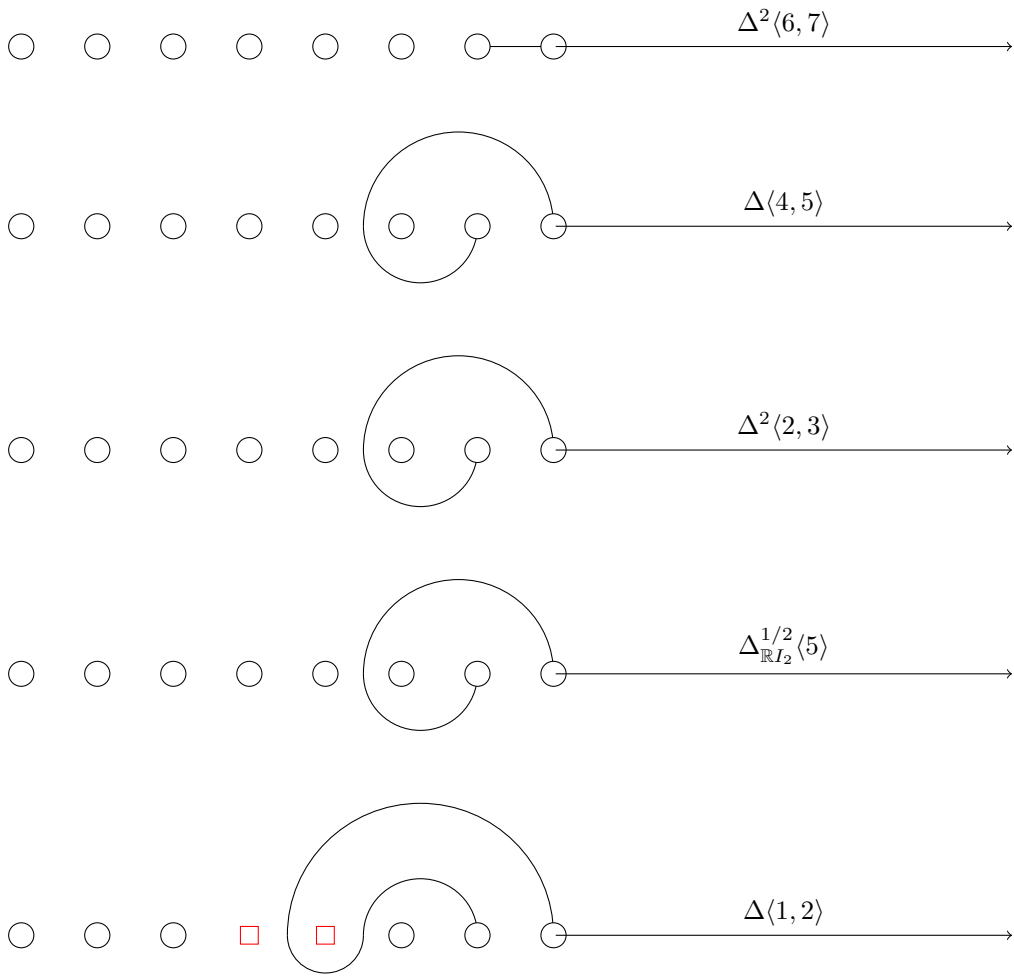
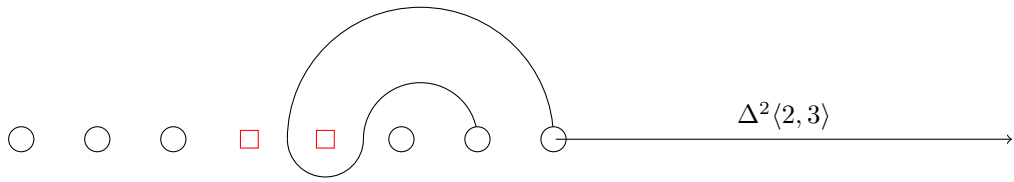
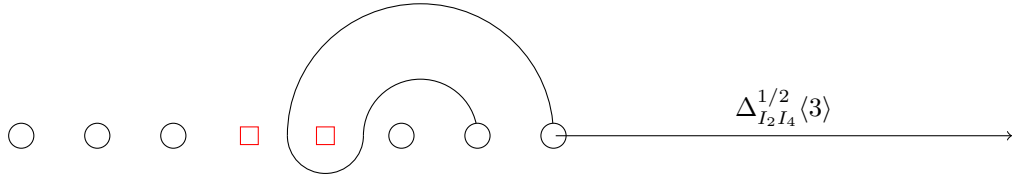
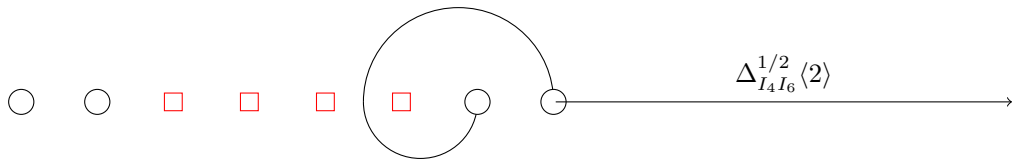
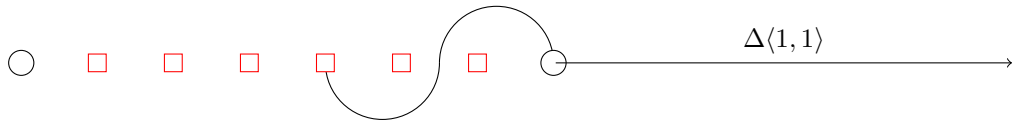
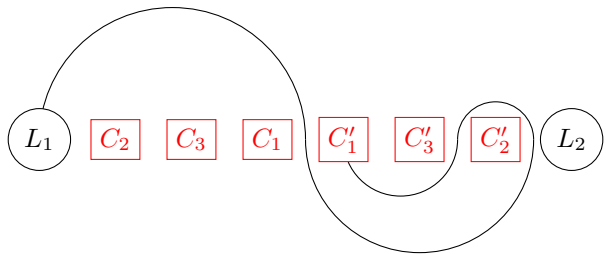
\begin{figure}[H]
\begin{center}
%
\end{center}
\end{figure}

Relation:

\begin{equation*}\begin{split}
\left[\Gamma_{C_2}^{-1}\Gamma_{C_3}^{-1}\Gamma_{C_1}^{-1}\Gamma_{C_2'}\Gamma_{C_1'}\Gamma_{C_2'}^{-1}\Gamma_{C_1}\Gamma_{C_3}\Gamma_{C_2} , \Gamma_{L_1} \right]=e
\end{split}
\end{equation*}
\subsubsection{Vertex Number 11}

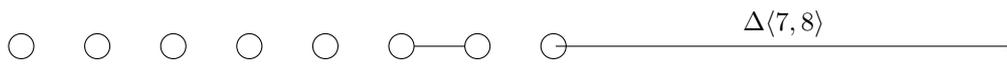
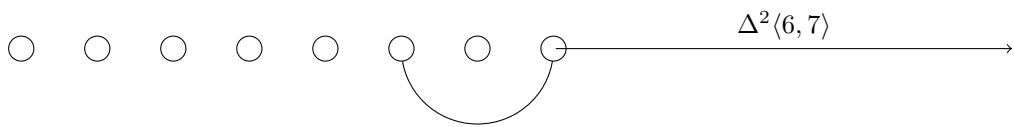
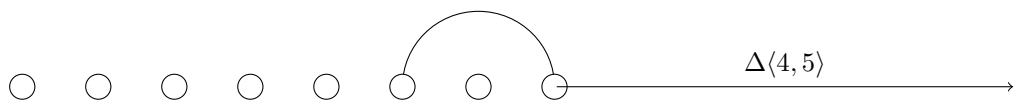
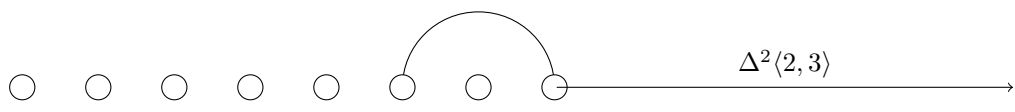
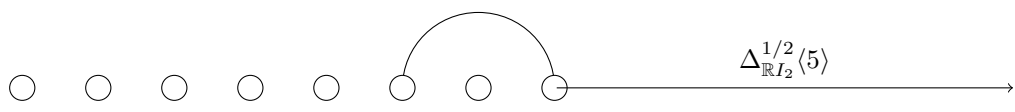
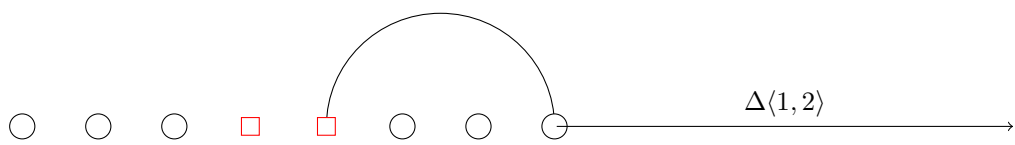
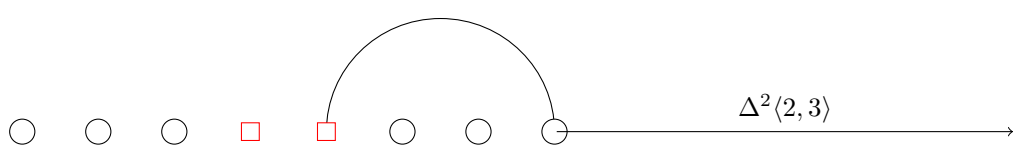
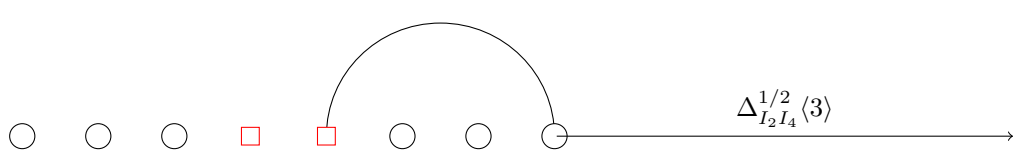
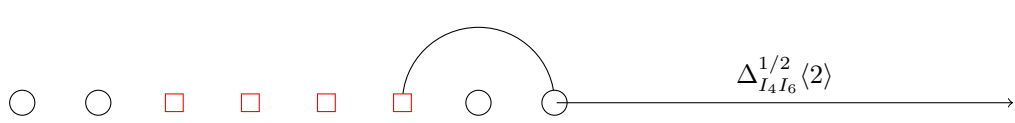
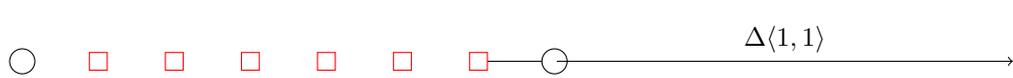
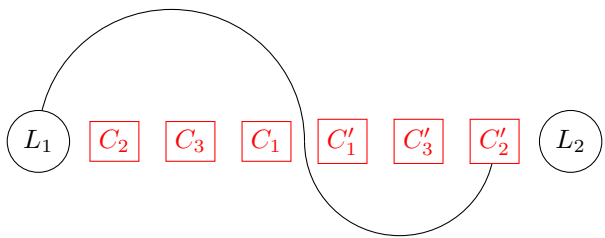
\begin{figure}[H]
\begin{center}
%
\end{center}
\end{figure}

Relation:

\begin{equation*}\begin{split}
\left\{\Gamma_{C_2}^{-1}\Gamma_{C_3}^{-1}\Gamma_{C_1}^{-1}\Gamma_{C_2'}\Gamma_{C_1}\Gamma_{C_3}\Gamma_{C_2}  , \Gamma_{L_1} \right\}=e
\end{split}
\end{equation*}
\subsubsection{Vertex Number 12}

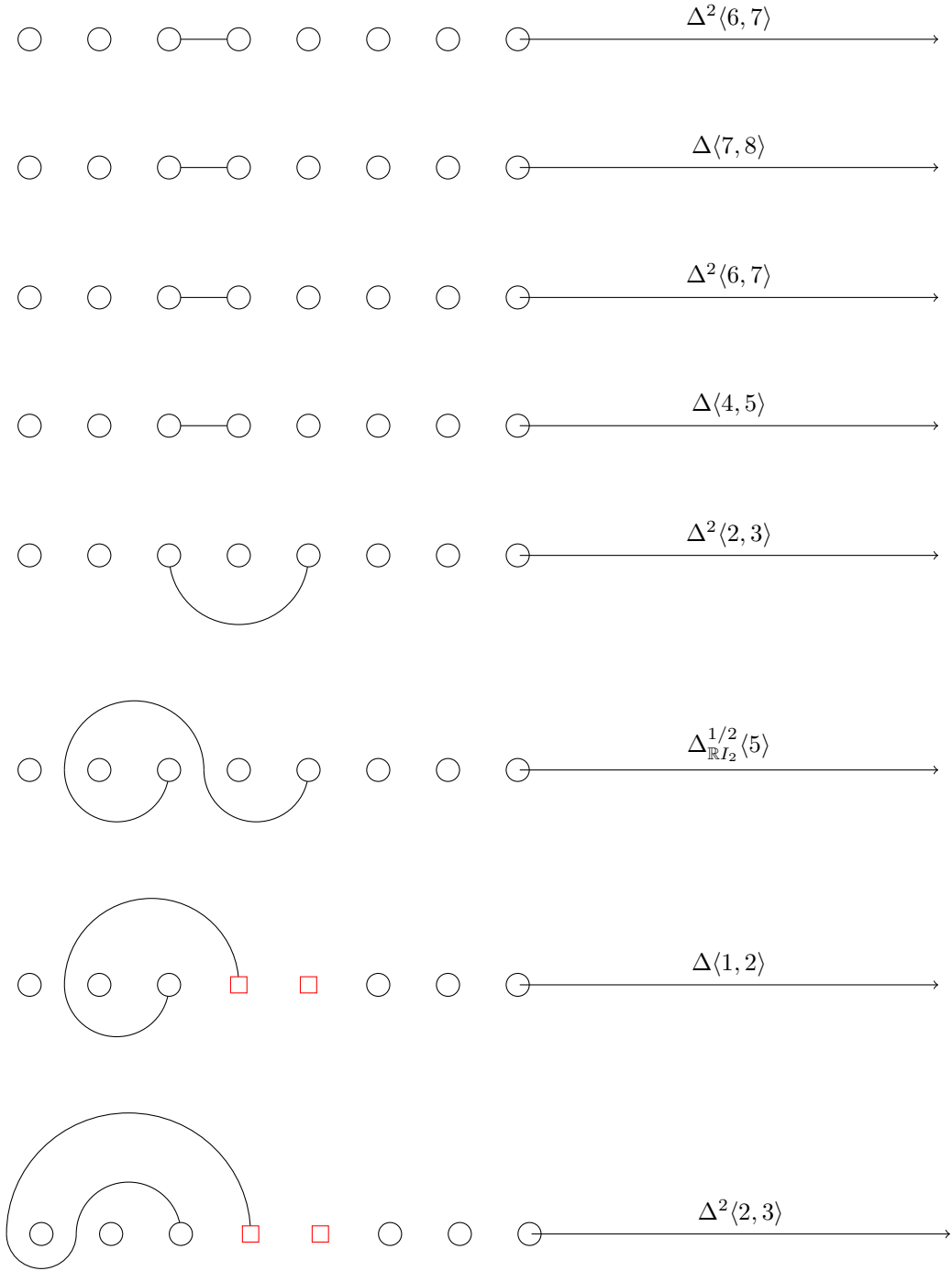
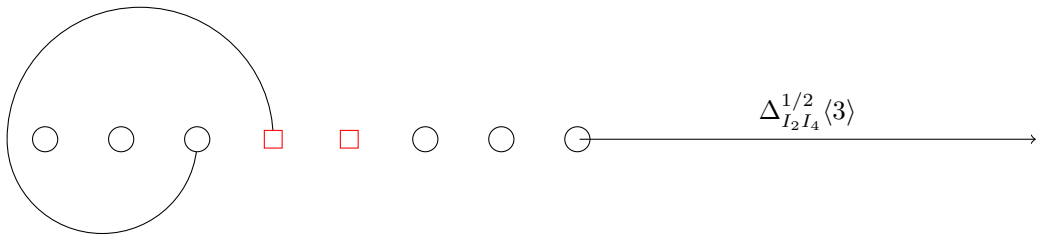
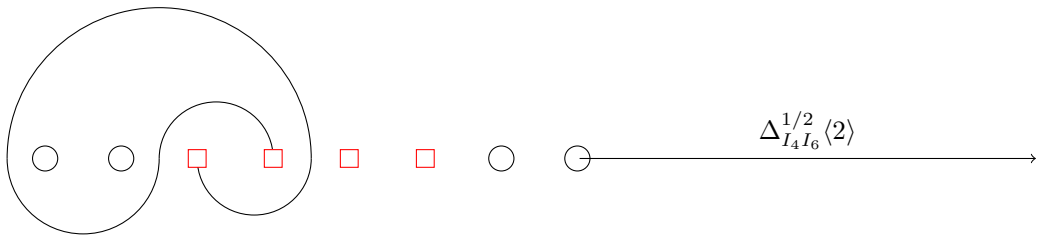
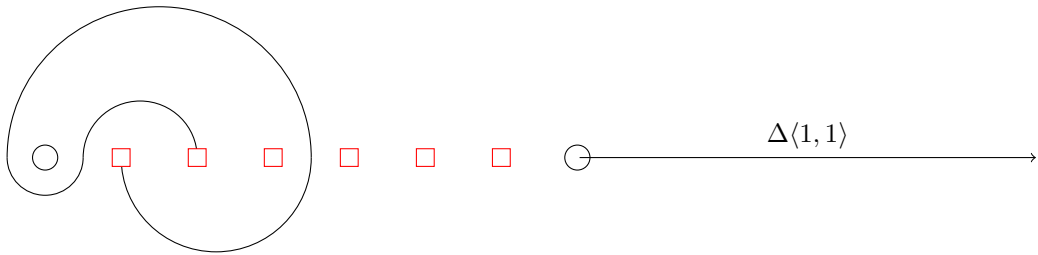
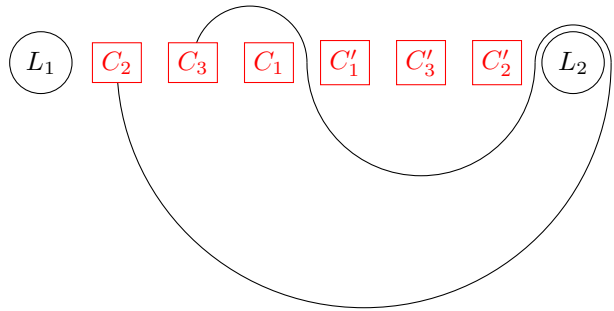
\begin{figure}[H]
\begin{center}
%
\end{center}
\end{figure}

Relation:

\begin{equation*}\begin{split}
\left[\Gamma_{L_2}\Gamma_{C_1}\Gamma_{C_3}\Gamma_{C_1}^{-1}\Gamma_{L_2}^{-1} , \Gamma_{C_2} \right]=e
\end{split}
\end{equation*}
\subsubsection{Vertex Number 13}

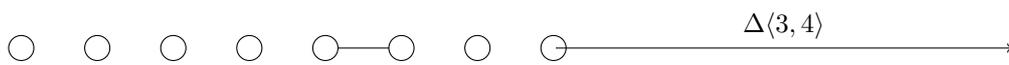
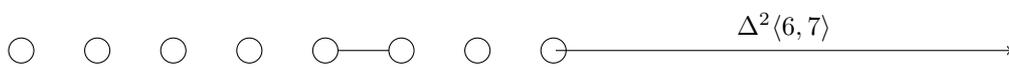
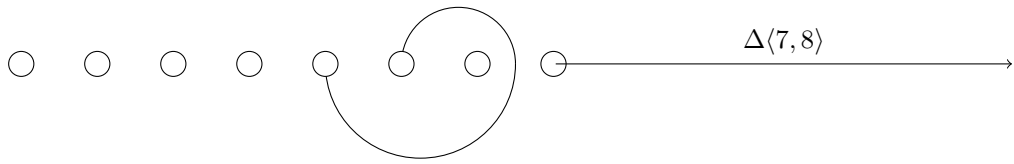
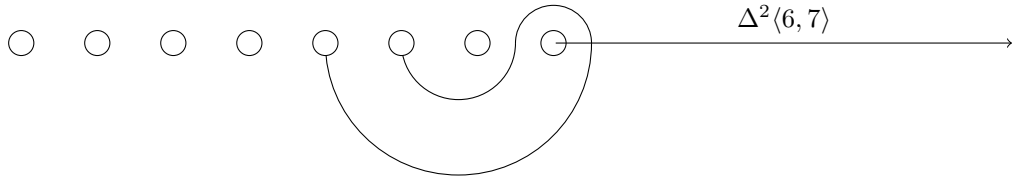
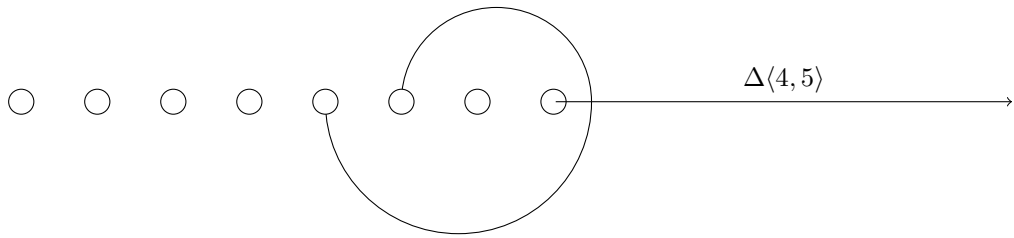
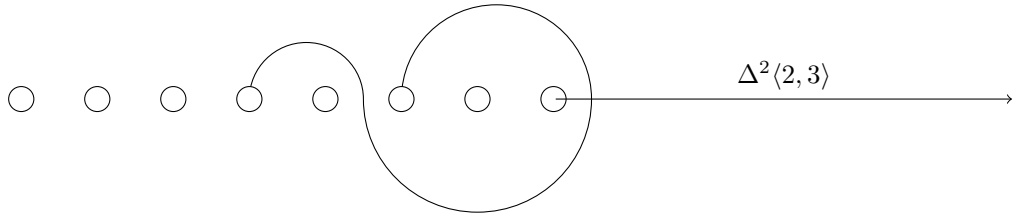
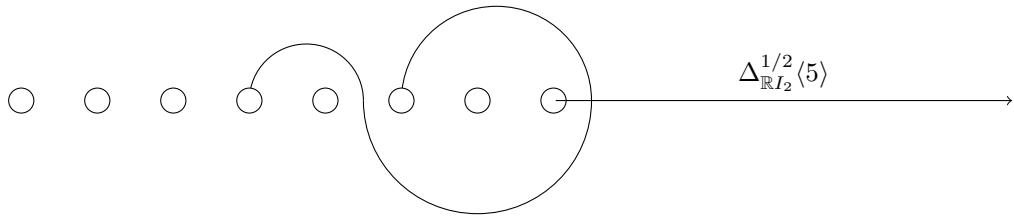
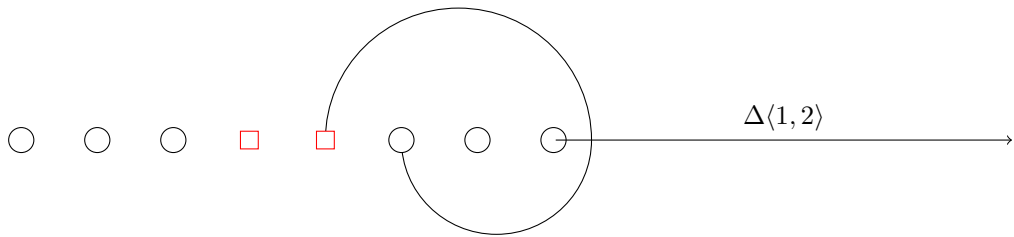
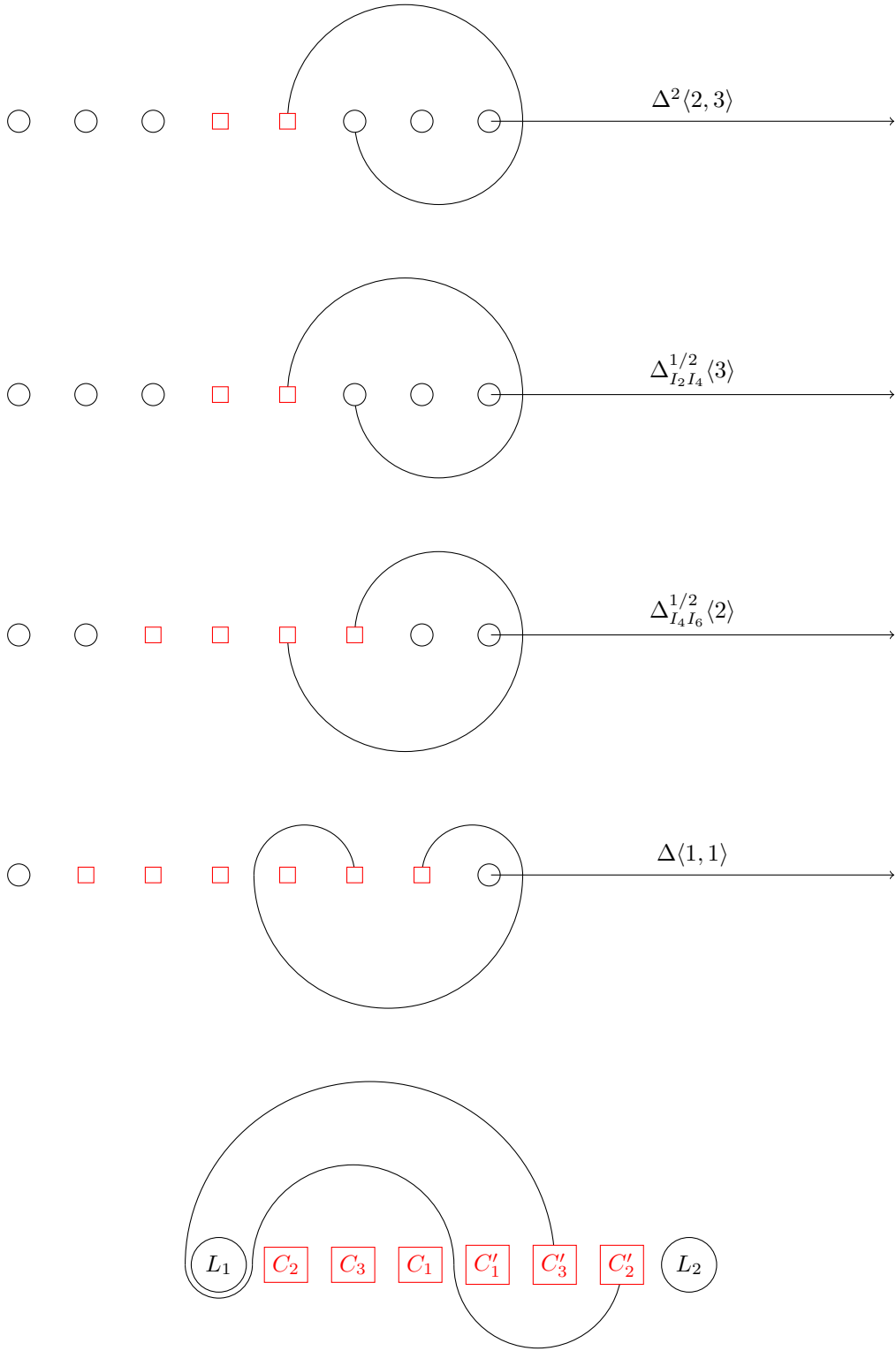
\begin{figure}[H]
\begin{center}
%
\end{center}
\end{figure}

Relation:

\begin{equation*}\begin{split}
\left[\Gamma_{C_1}\Gamma_{C_3}\Gamma_{C_2}\Gamma_{L_1}^{-1}\Gamma_{C_2}^{-1}\Gamma_{C_3}^{-1}\Gamma_{C_1}^{-1}\Gamma_{C_1'}^{-1}\Gamma_{C_3'}\Gamma_{C_1'}\Gamma_{C_1}\Gamma_{C_3}\Gamma_{C_2}\Gamma_{L_1}\Gamma_{C_2}^{-1}\Gamma_{C_3}^{-1}\Gamma_{C_1}^{-1} , \Gamma_{C_2'} \right]=e
\end{split}
\end{equation*}
\subsubsection{Vertex Number 14}

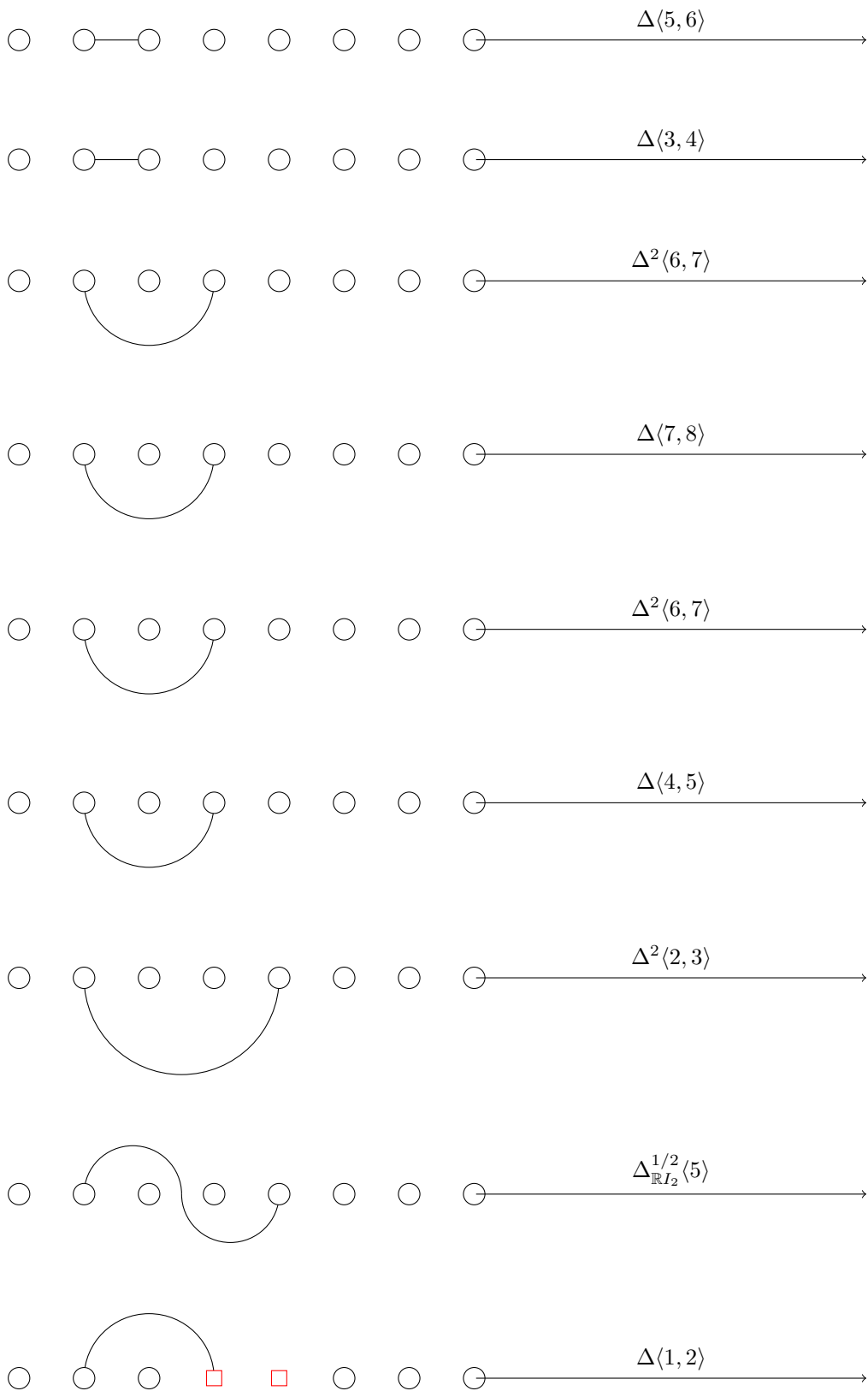
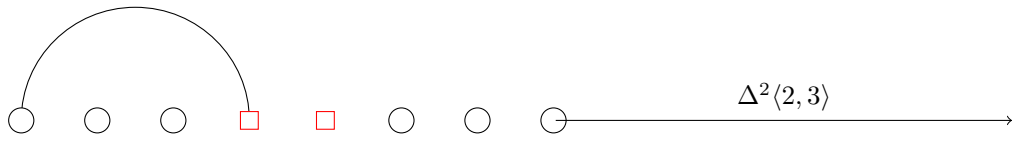
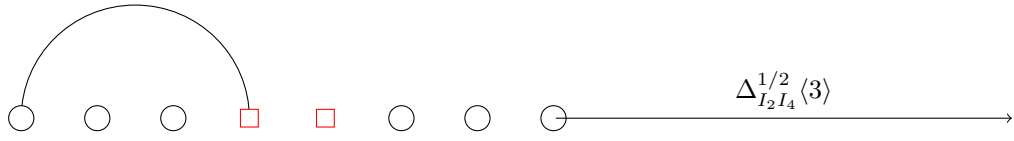
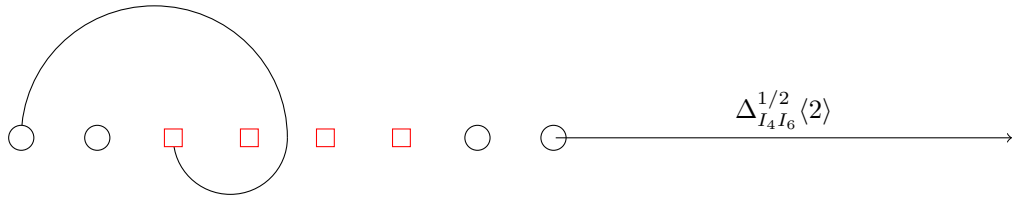
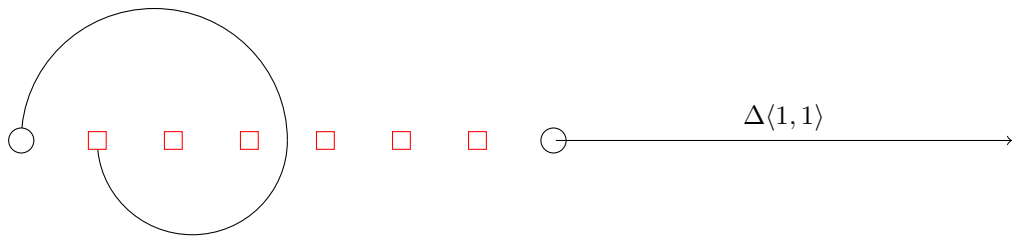
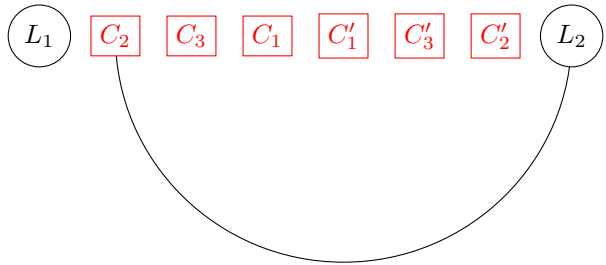
\begin{figure}[H]
\begin{center}
%
\end{center}
\end{figure}

Relation:

\begin{equation*}\begin{split}
\left\{\Gamma_{L_2}  , \Gamma_{C_2} \right\}=e
\end{split}
\end{equation*}
\subsubsection{Vertex Number 15}

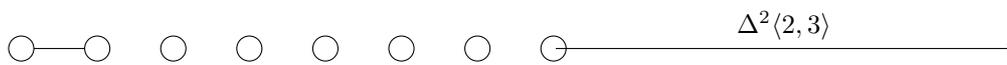
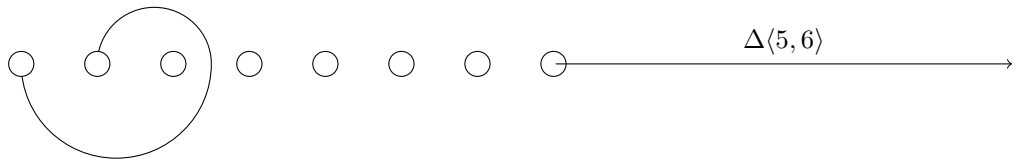
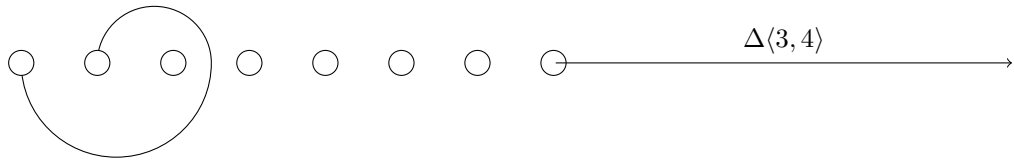
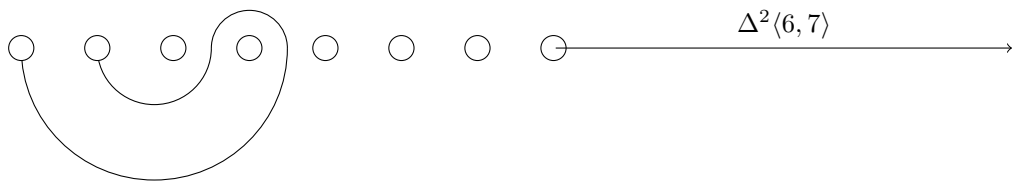
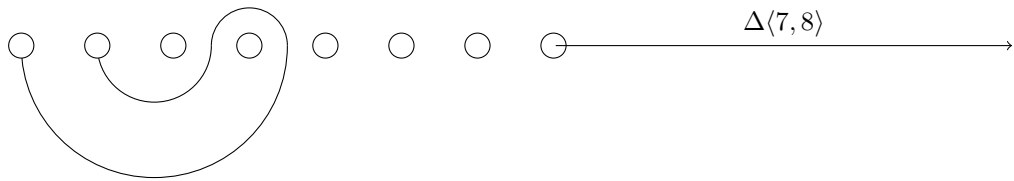
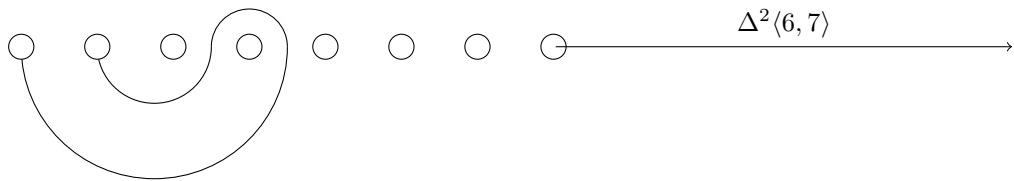
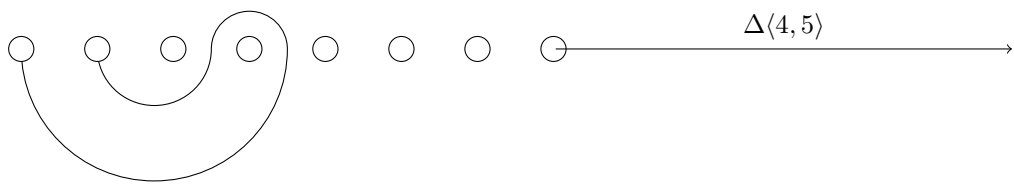
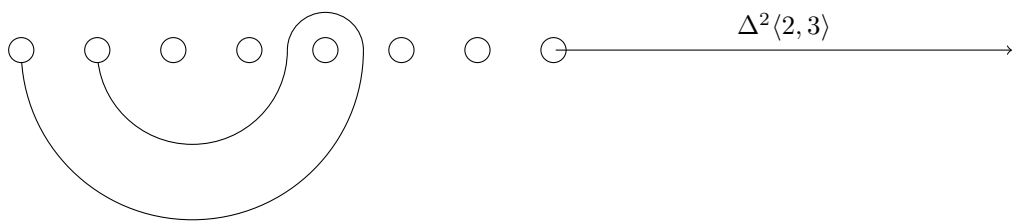
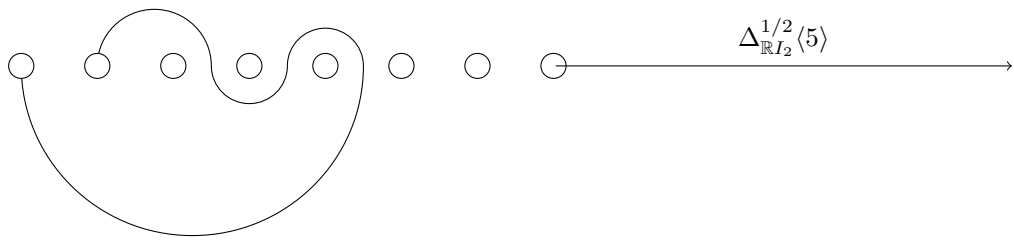
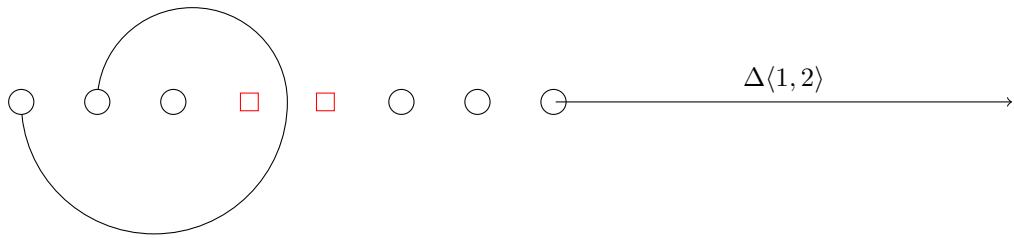
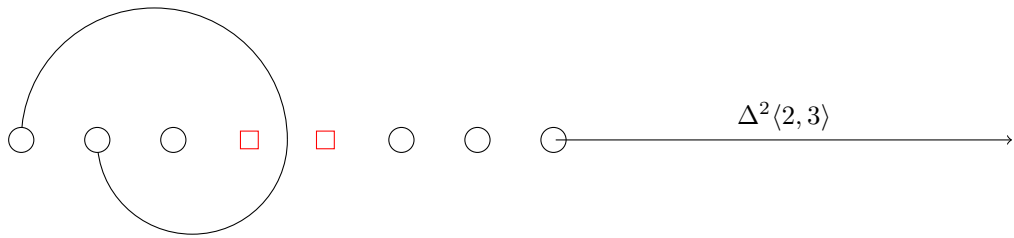
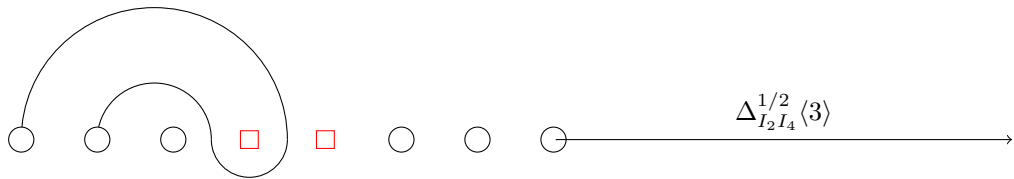
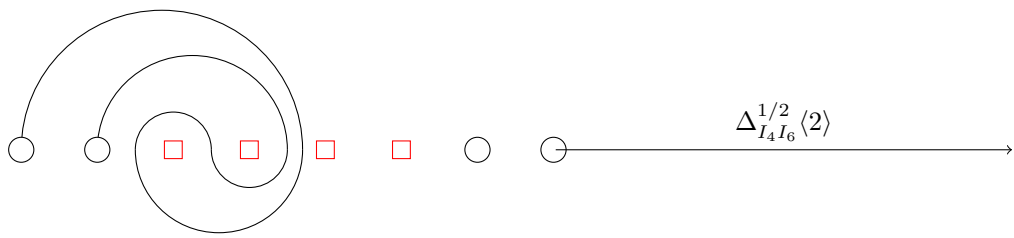
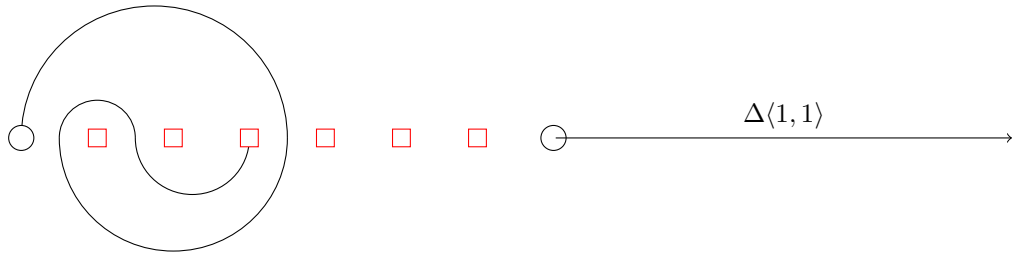
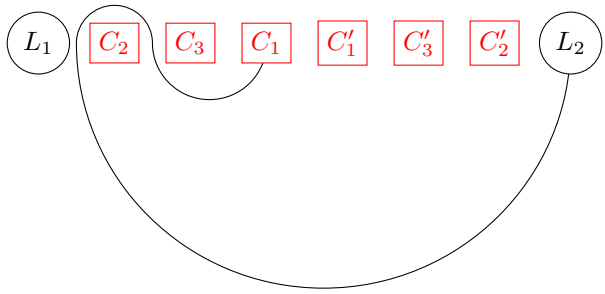
\begin{figure}[H]
\begin{center}
%
\end{center}
\end{figure}

Relation:

\begin{equation*}\begin{split}
\left[\Gamma_{C_2}^{-1}\Gamma_{C_1}\Gamma_{C_2} , \Gamma_{L_2} \right]=e
\end{split}
\end{equation*}
\subsubsection{Vertex Number 16}

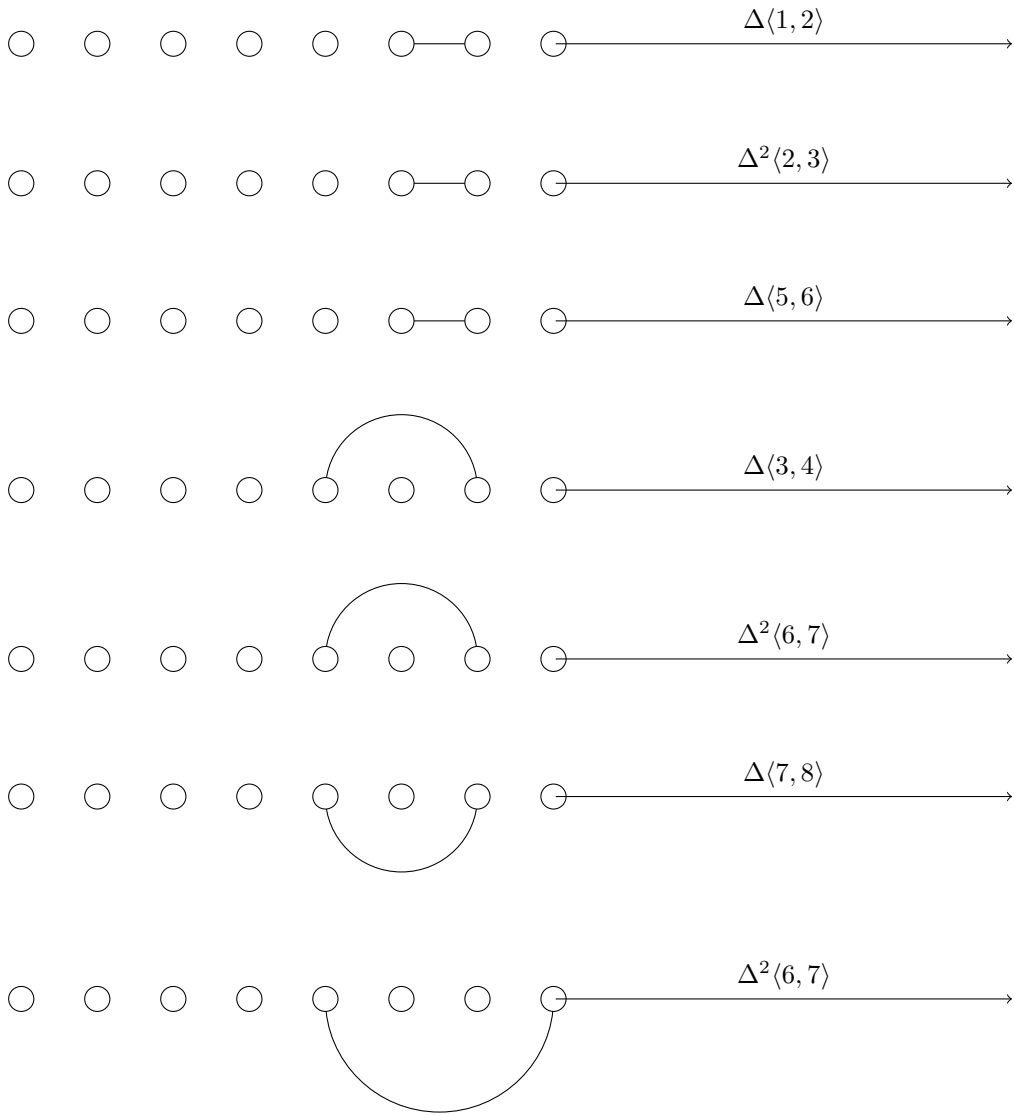
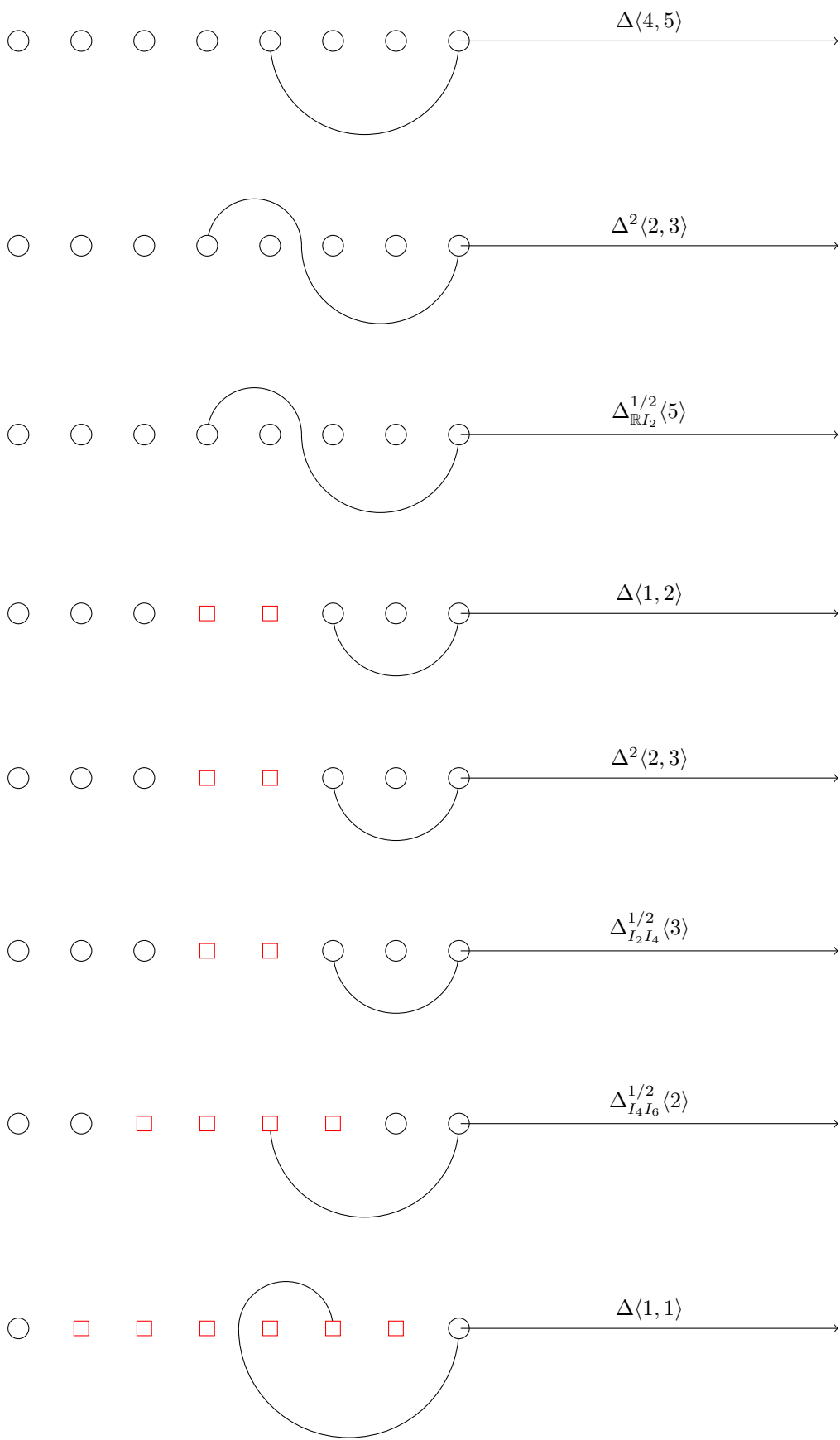
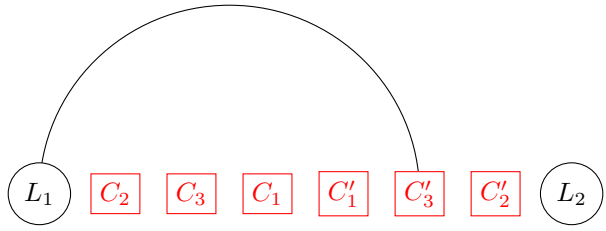
\begin{figure}[H]
\begin{center}
%
\end{center}
\end{figure}

Relation:

\begin{equation*}\begin{split}
\left\{\Gamma_{C_2}^{-1}\Gamma_{C_3}^{-1}\Gamma_{C_1}^{-1}\Gamma_{C_1'}^{-1}\Gamma_{C_3'}\Gamma_{C_1'}\Gamma_{C_1}\Gamma_{C_3}\Gamma_{C_2}  , \Gamma_{L_1} \right\}=e
\end{split}
\end{equation*}
\subsubsection{Vertex Number 17}

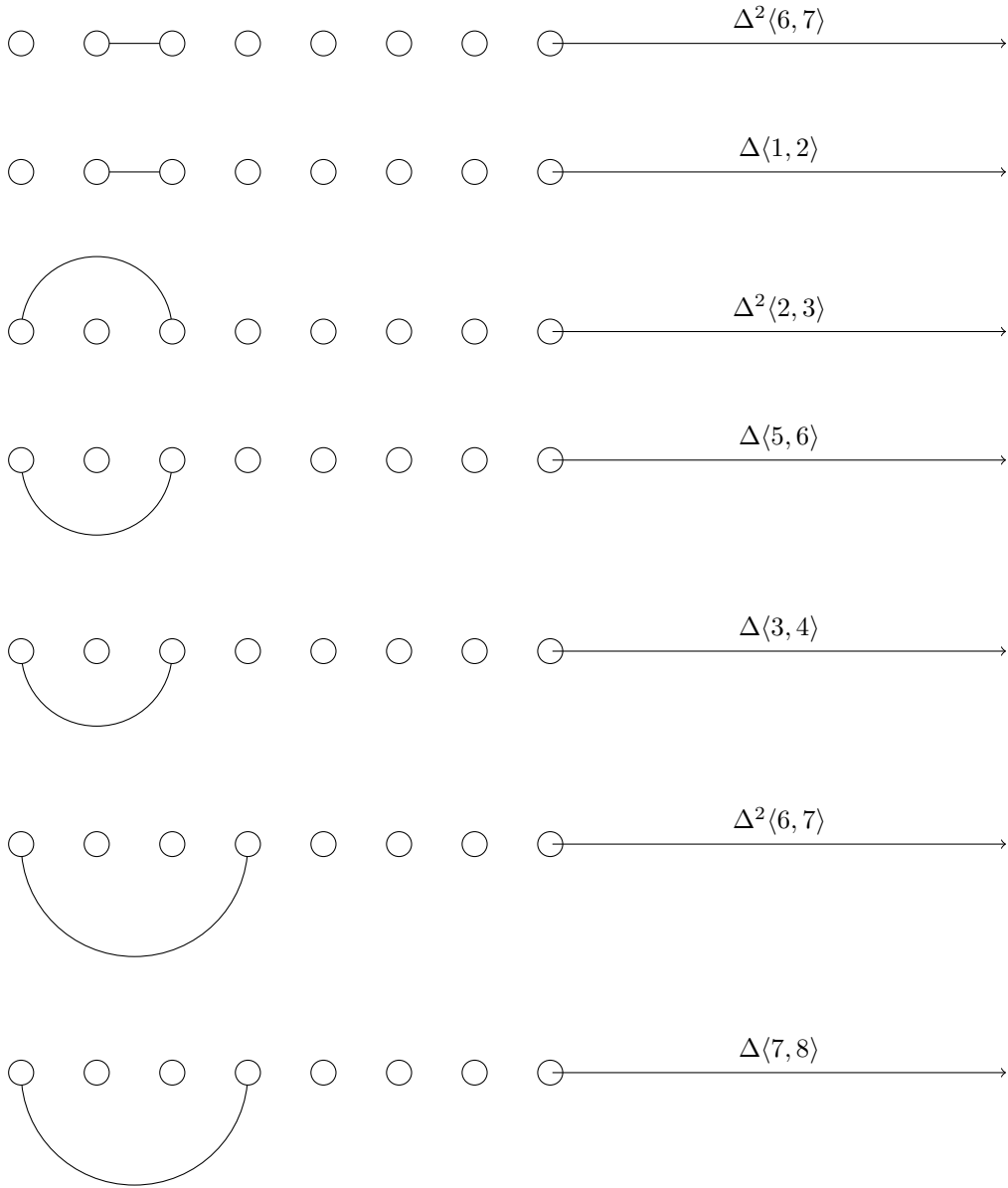
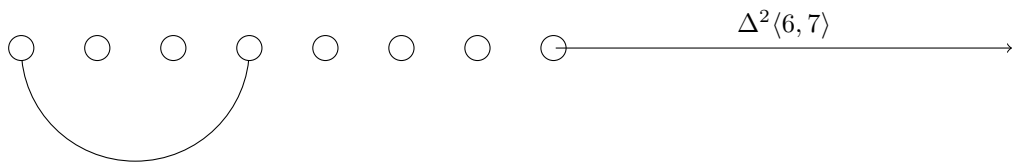
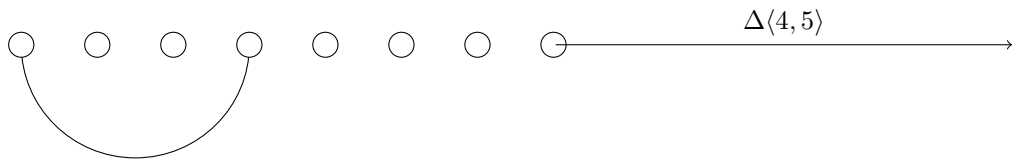
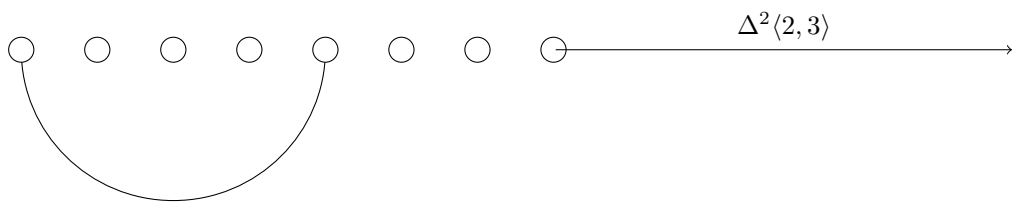
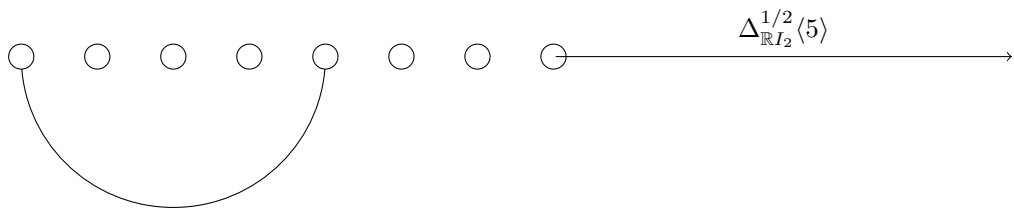
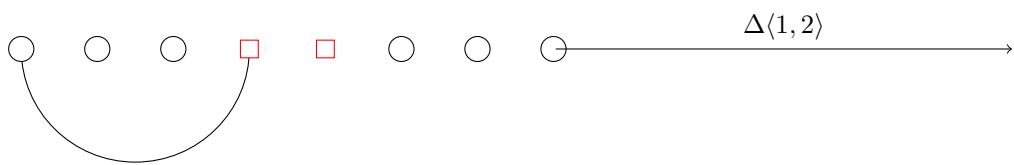
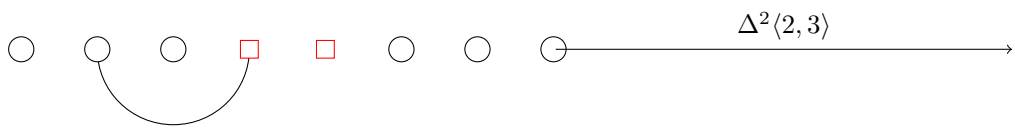
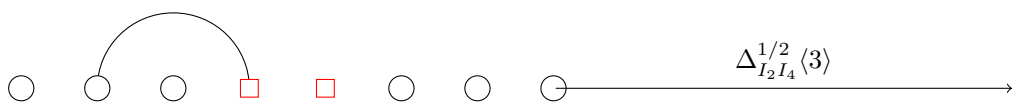
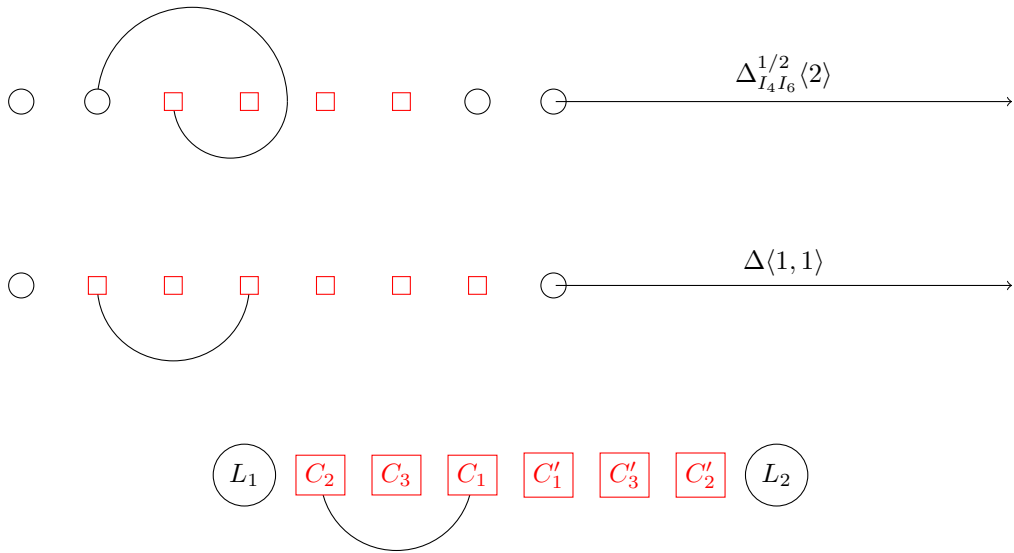
\begin{figure}[H]
\begin{center}
%
\end{center}
\end{figure}

Relation:

\begin{equation*}\begin{split}
\left\{\Gamma_{C_1}  , \Gamma_{C_2} \right\}=e
\end{split}
\end{equation*}
\subsubsection{Vertex Number 18}

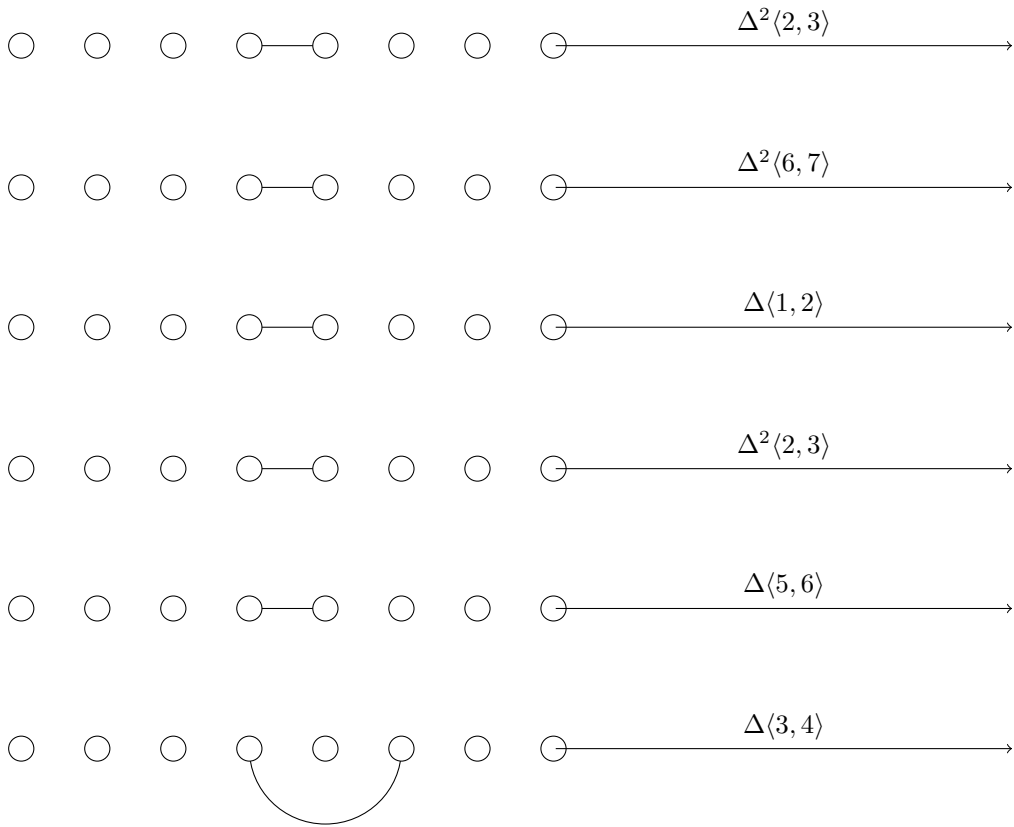
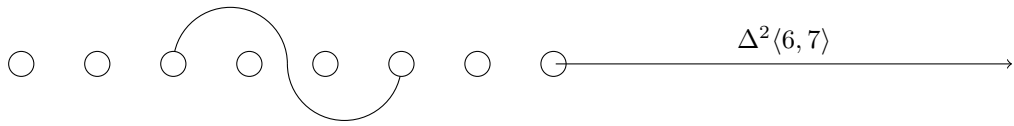
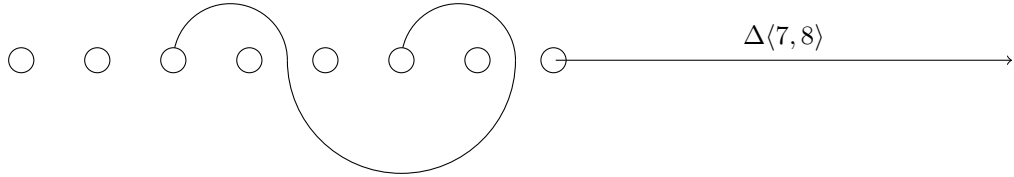
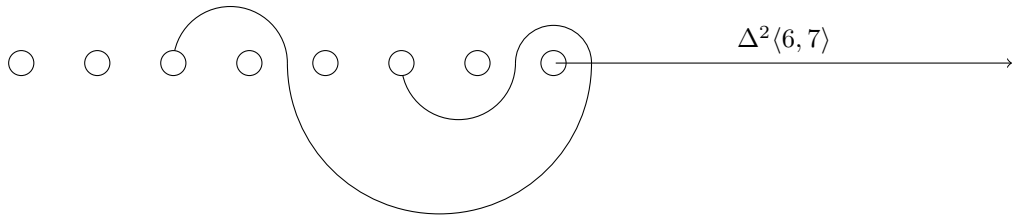
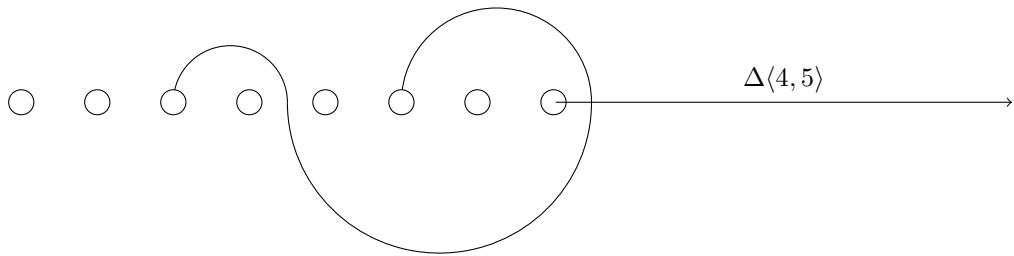
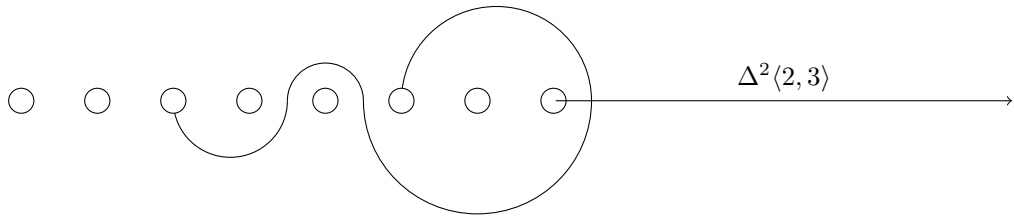
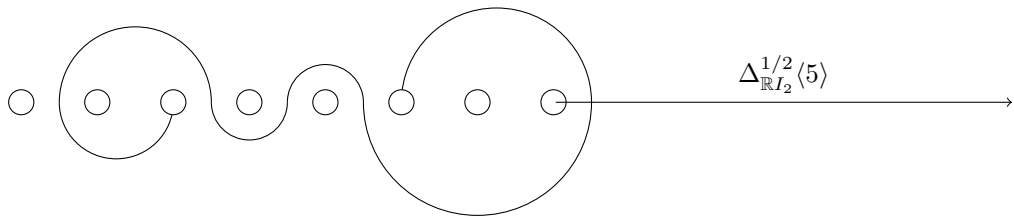
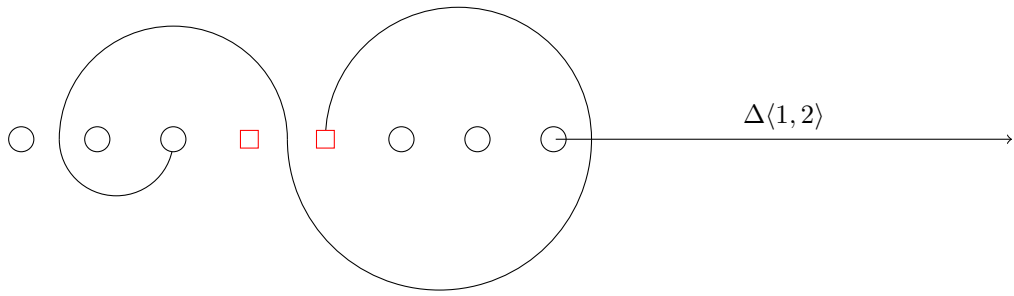
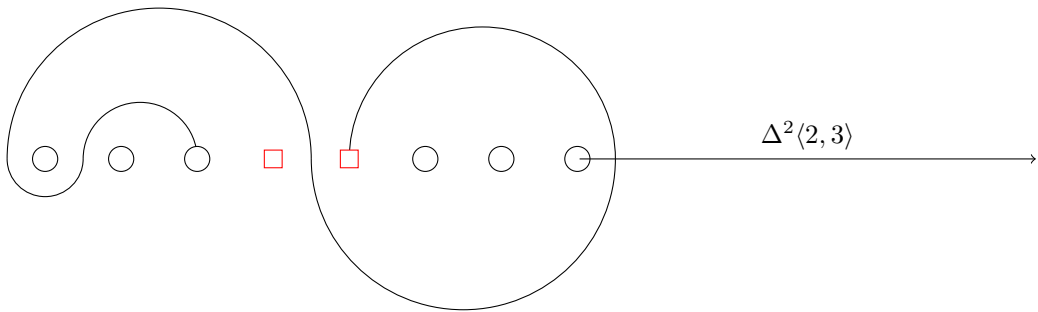
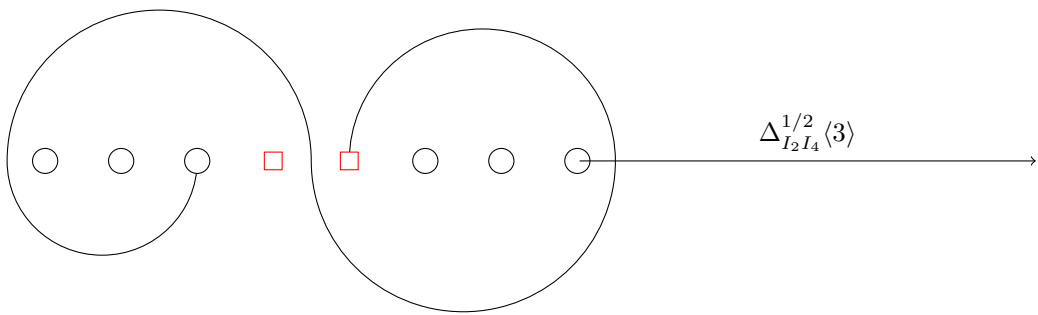
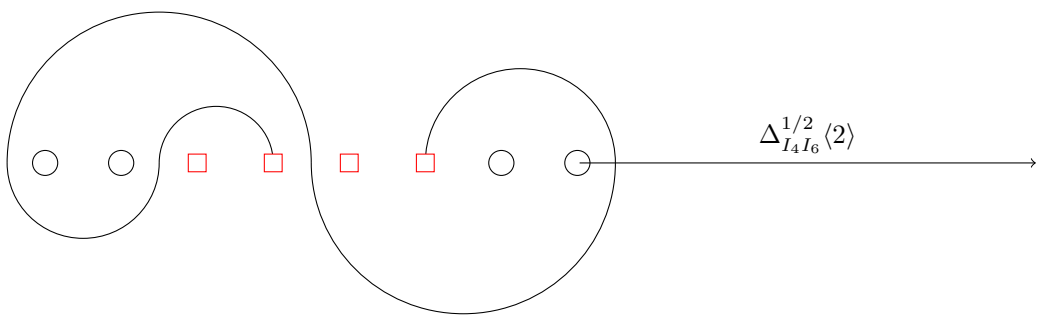
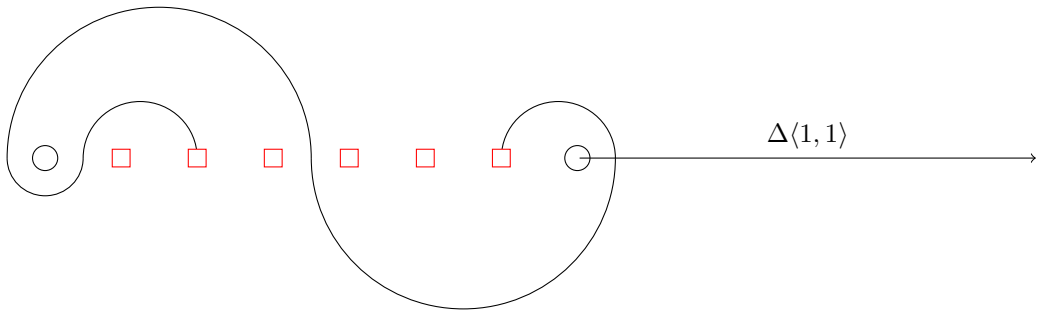
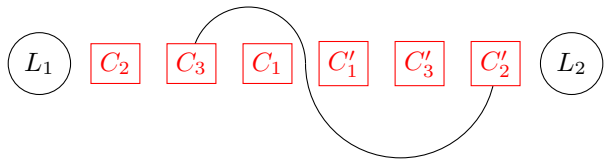
\begin{figure}[H]
\begin{center}
%
\end{center}
\end{figure}

Relation:

\begin{equation*}\begin{split}
\left[\Gamma_{C_1}\Gamma_{C_3}\Gamma_{C_1}^{-1} , \Gamma_{C_2'} \right]=e
\end{split}
\end{equation*}
\subsubsection{Vertex Number 19}

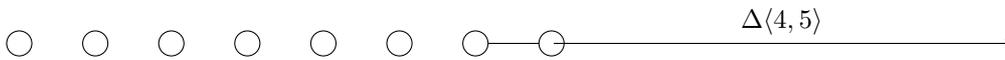
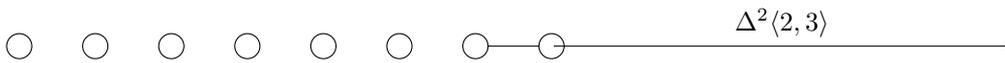
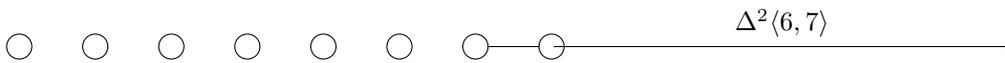
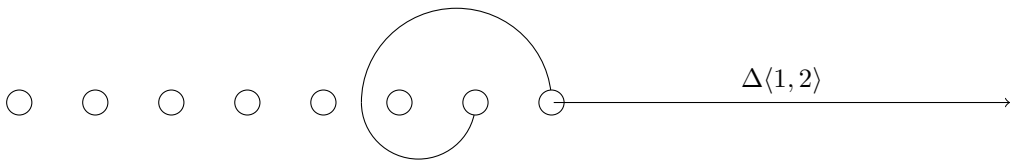
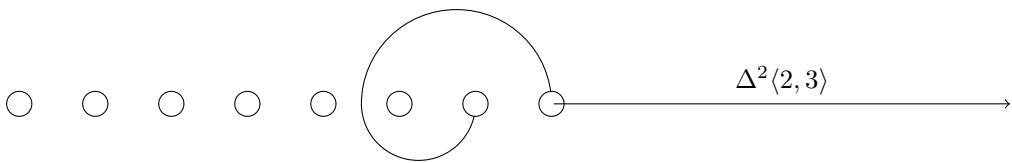
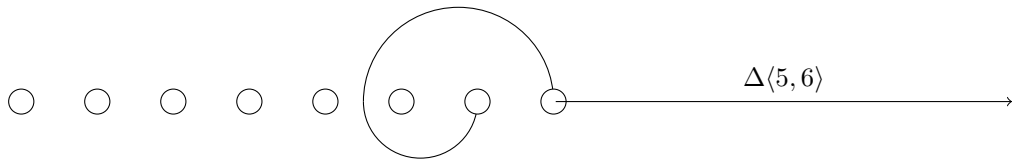
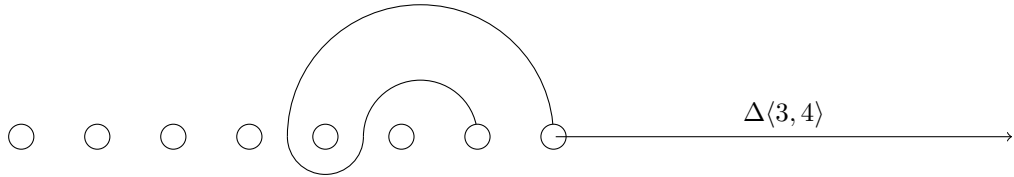
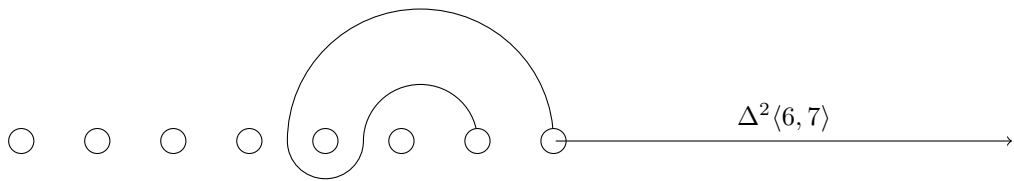
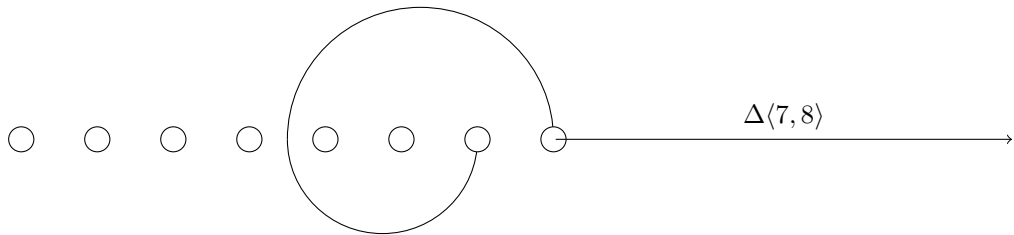
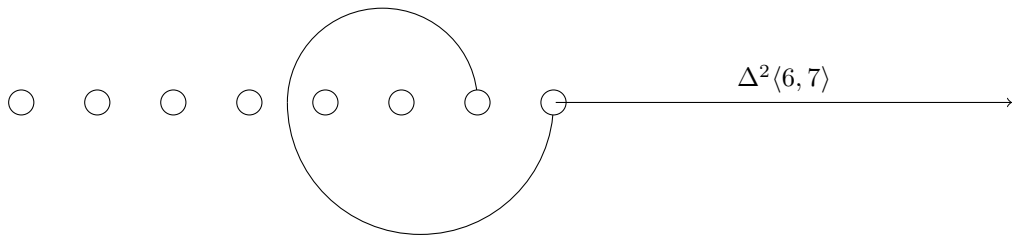
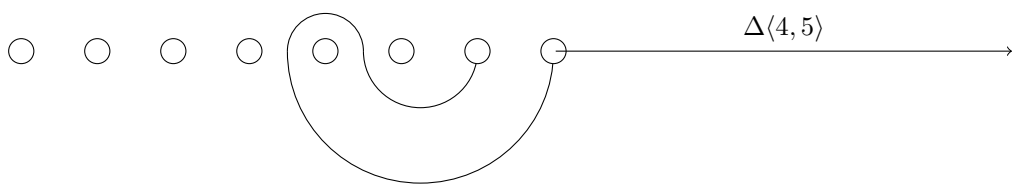
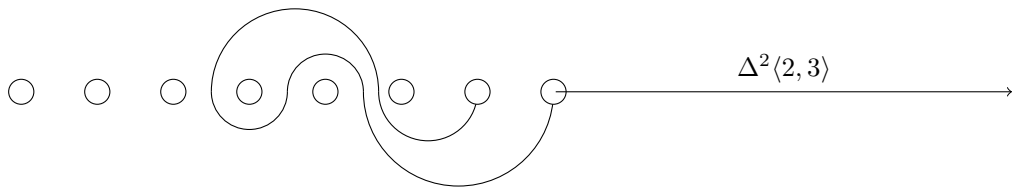
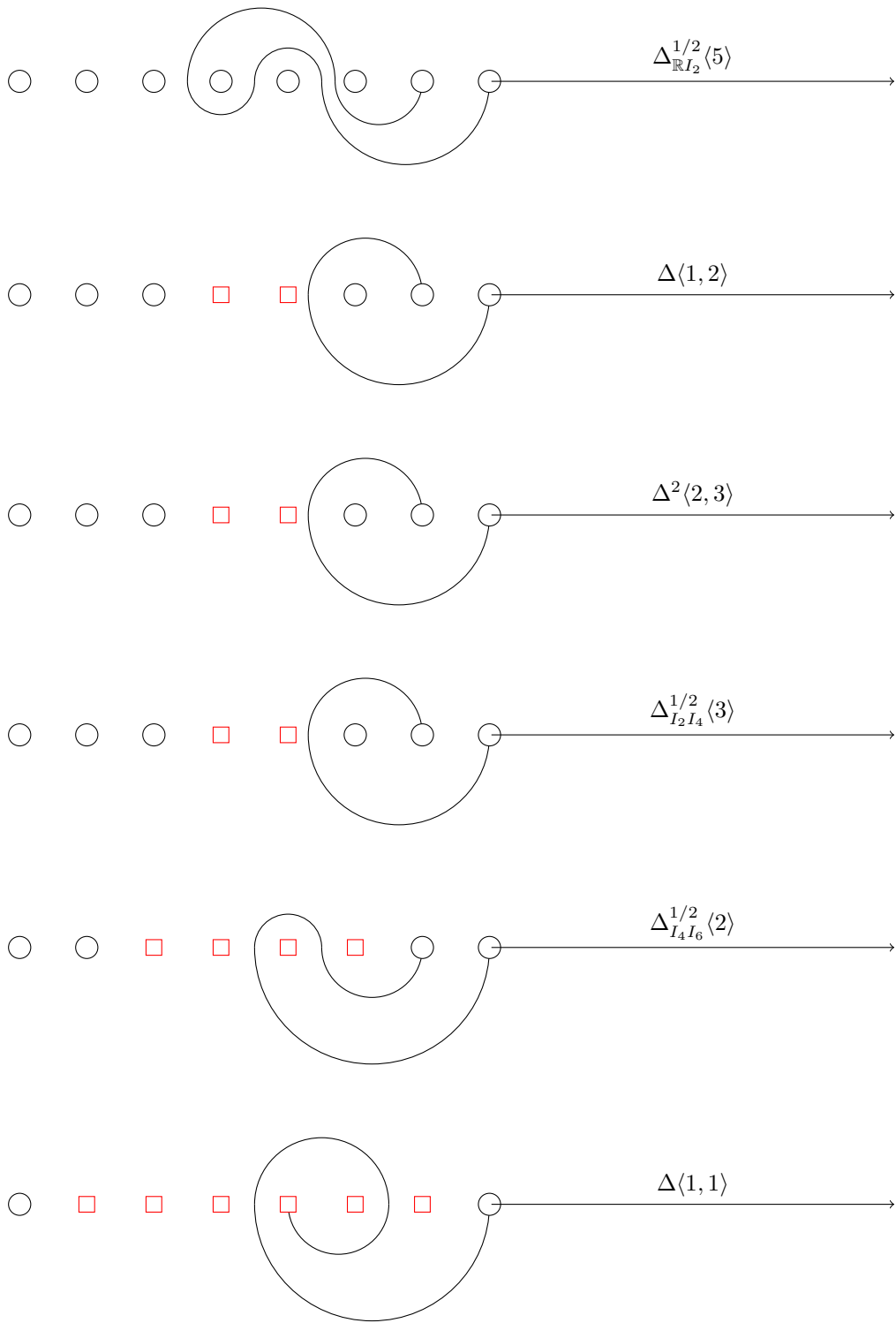
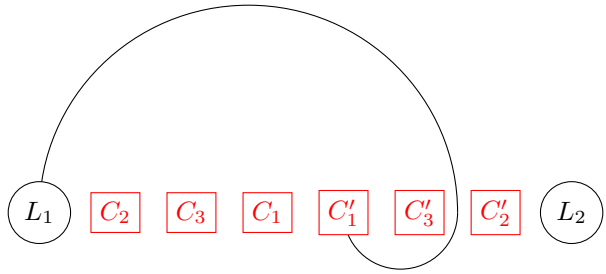
\begin{figure}[H]
\begin{center}
%
\end{center}
\end{figure}

Relation:

\begin{equation*}\begin{split}
\left[\Gamma_{C_3'}\Gamma_{C_1'}\Gamma_{C_1}\Gamma_{C_3}\Gamma_{C_2}\Gamma_{L_1}\Gamma_{C_2}^{-1}\Gamma_{C_3}^{-1}\Gamma_{C_1}^{-1}\Gamma_{C_1'}^{-1}\Gamma_{C_3'}^{-1} , \Gamma_{C_1'} \right]=e
\end{split}
\end{equation*}
\subsubsection{Vertex Number 20}

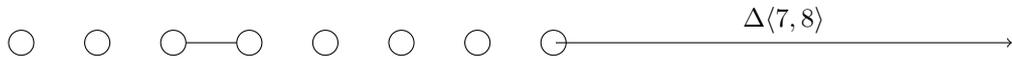
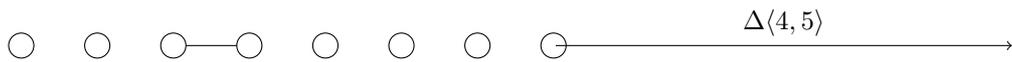
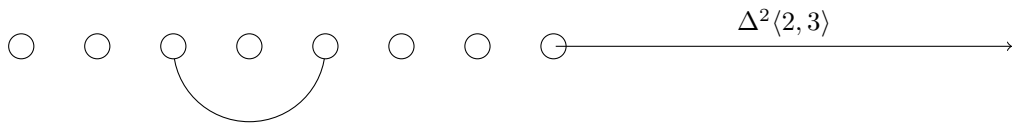
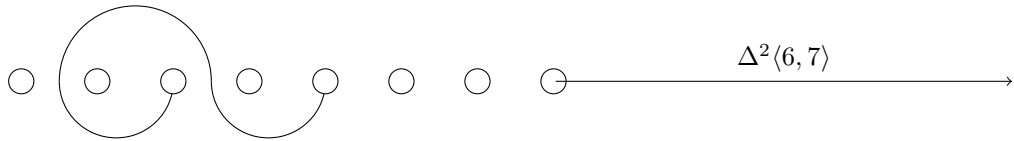
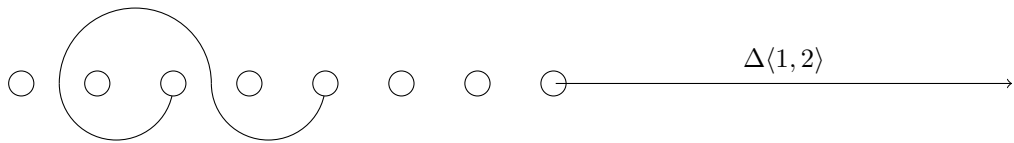
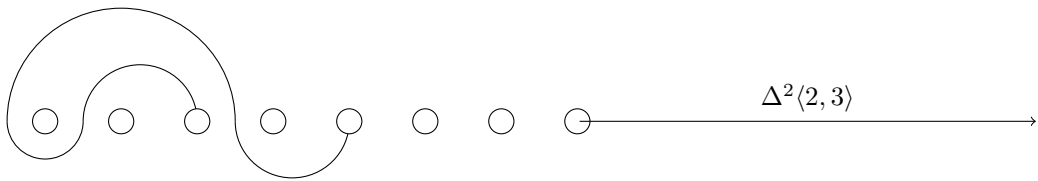
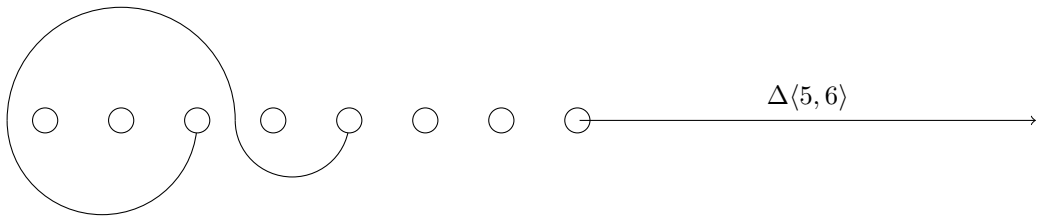
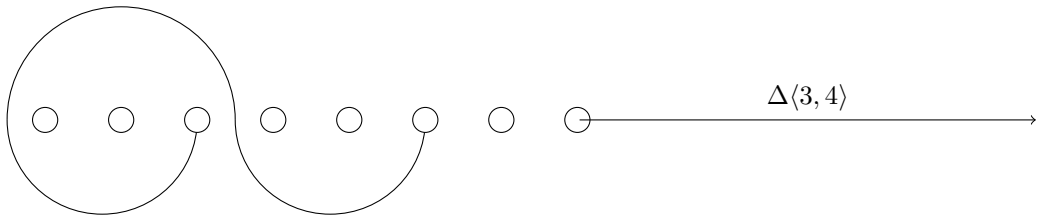
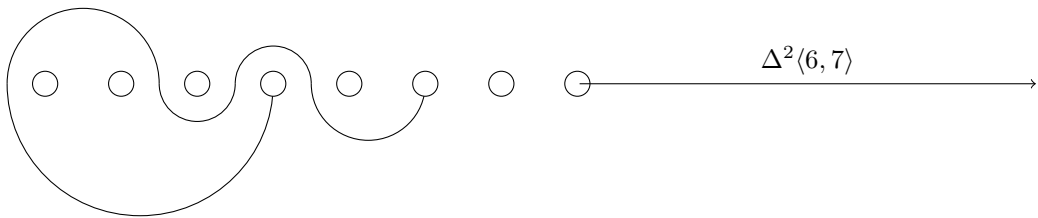
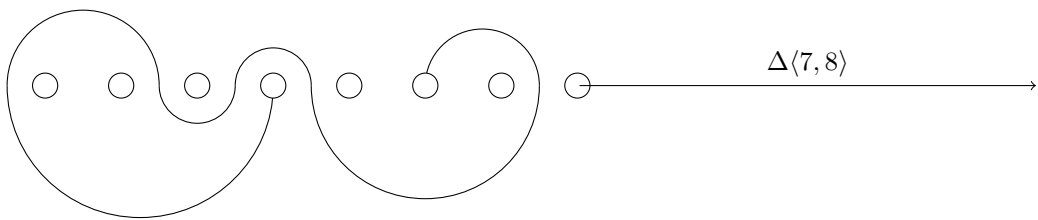
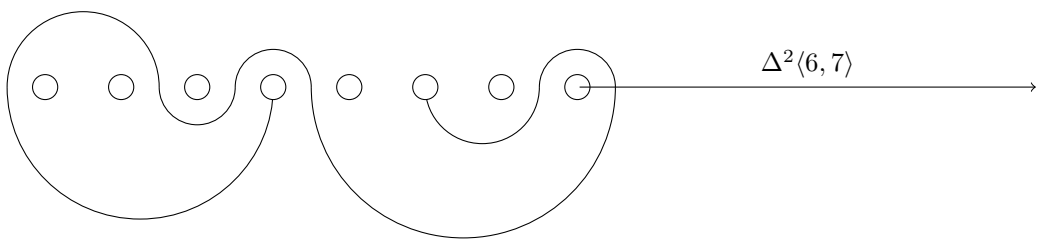
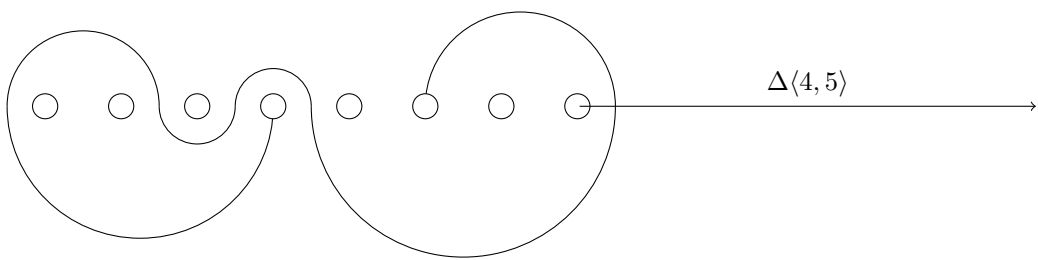
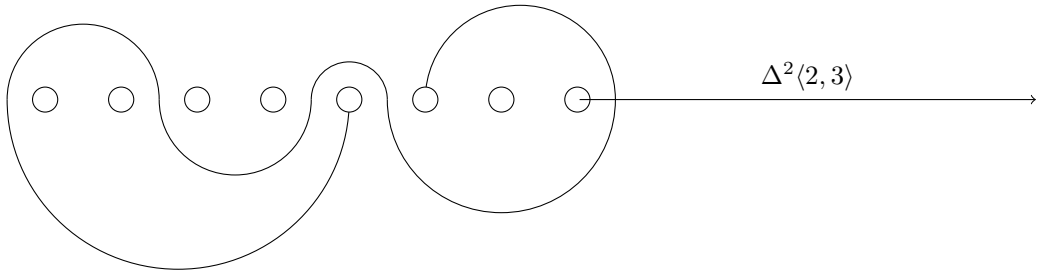
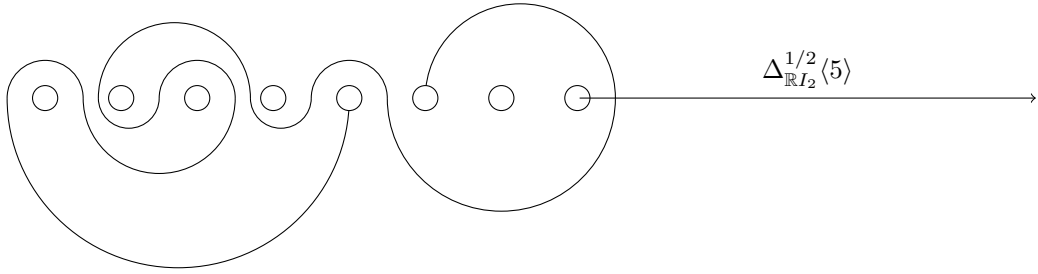
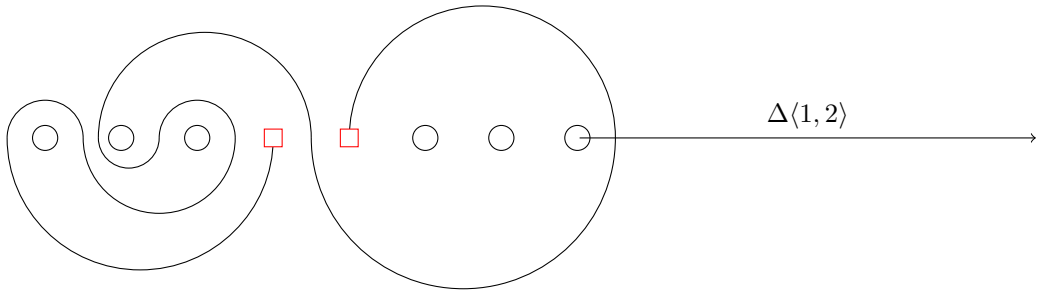
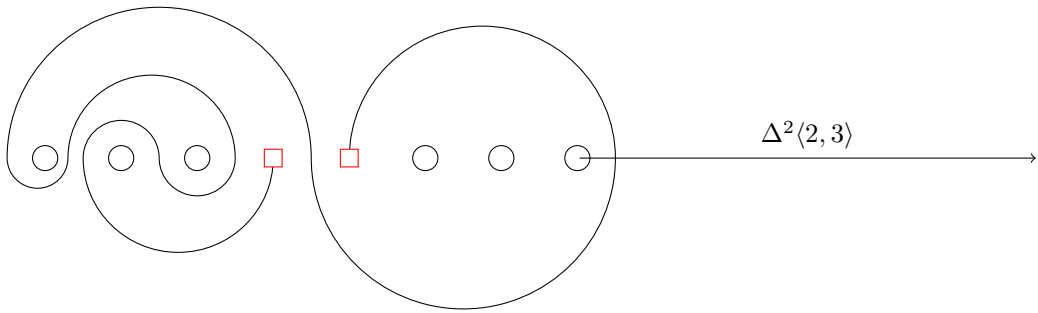
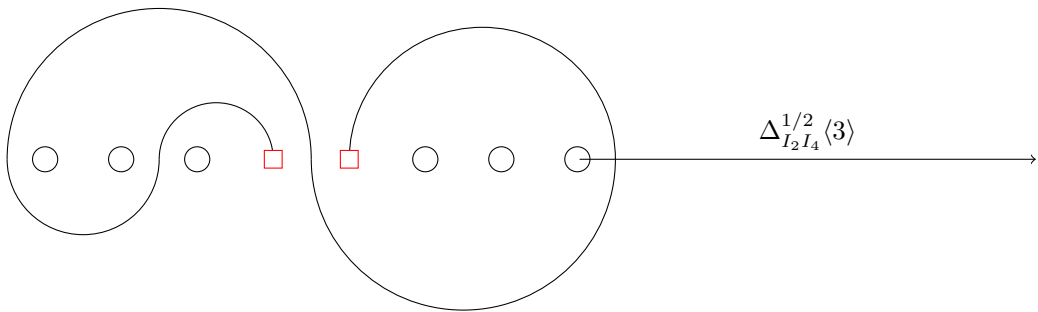
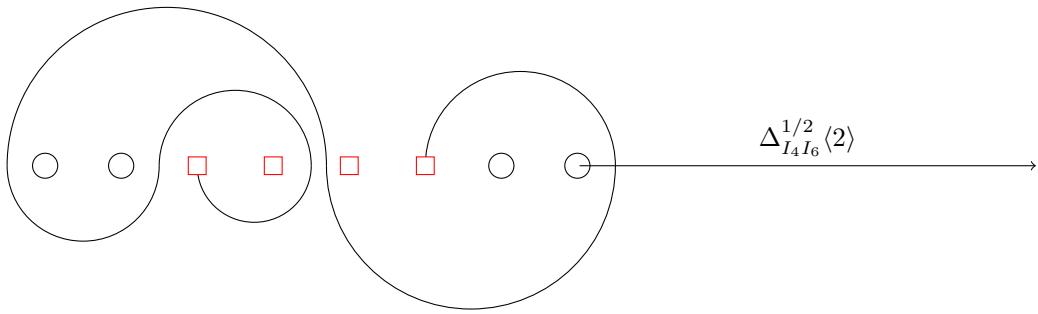
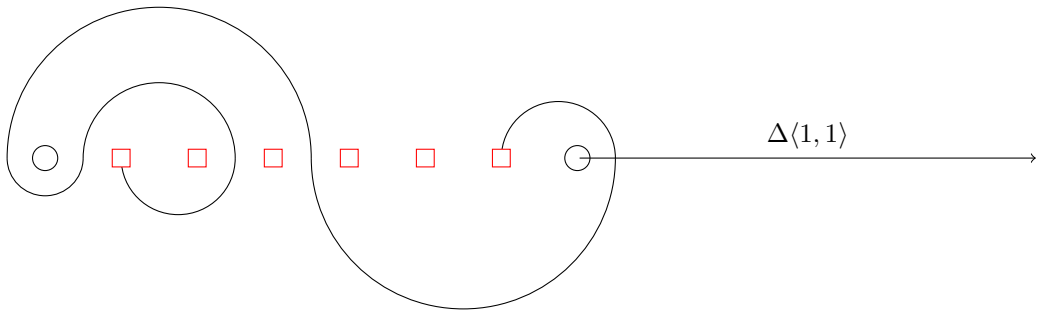
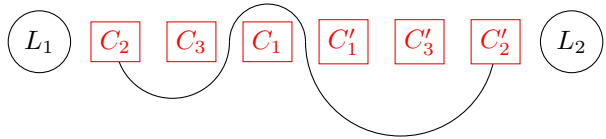
\begin{figure}[H]
\begin{center}
%
\end{center}
\end{figure}

Relation:

\begin{equation*}\begin{split}
\Gamma_{C_1}\Gamma_{C_2}\Gamma_{C_1}^{-1}  = \Gamma_{C_2'}
\end{split}
\end{equation*}
\subsubsection{Vertex Number 21}

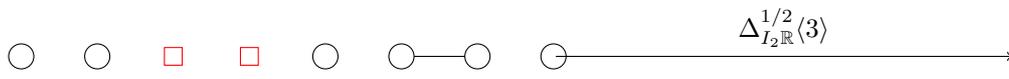
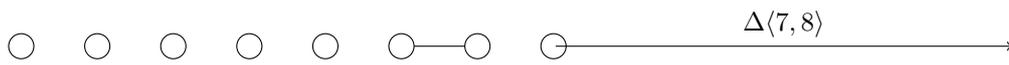
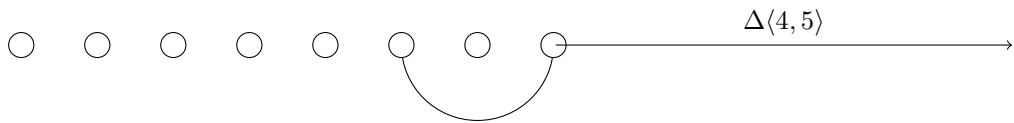
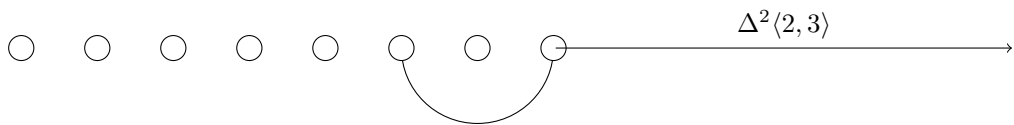
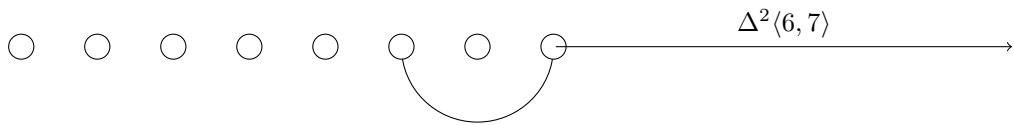
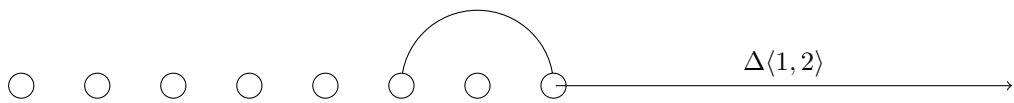
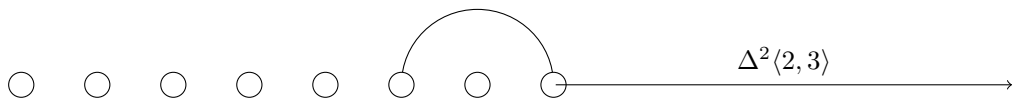
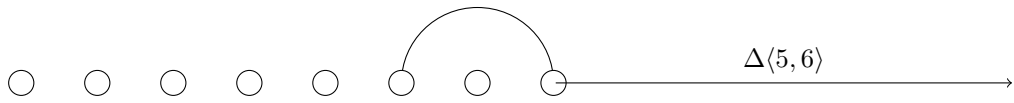
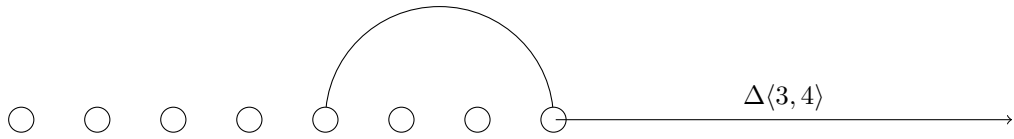
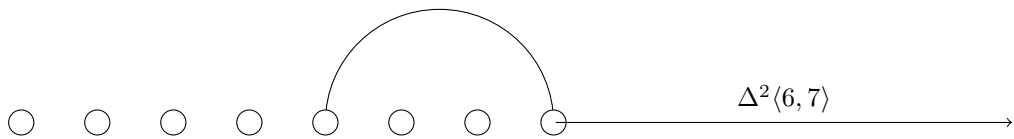
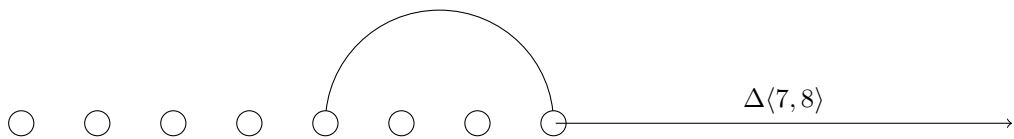
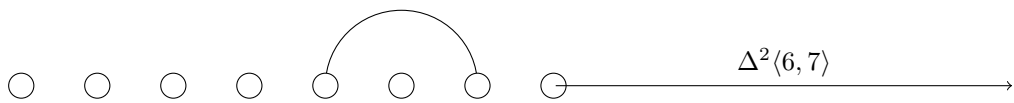
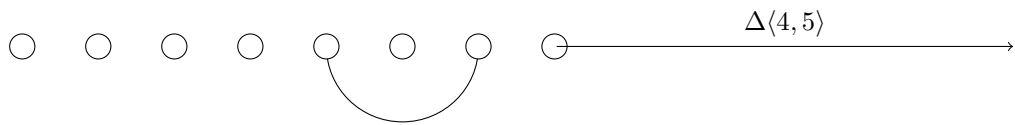
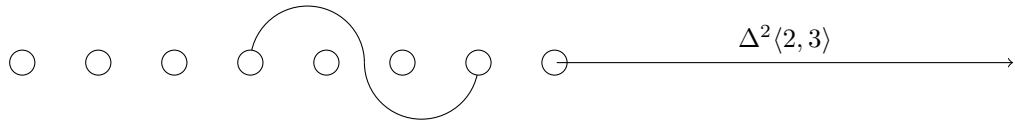
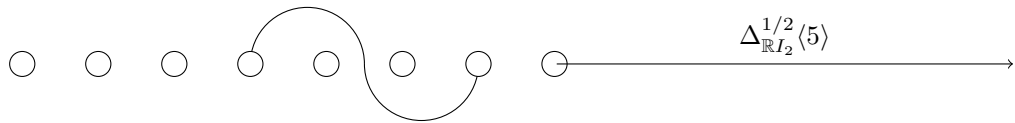
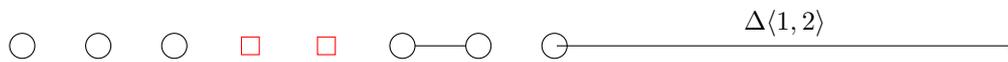
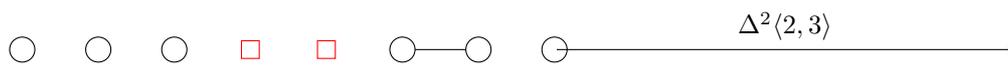
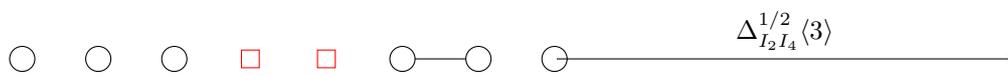
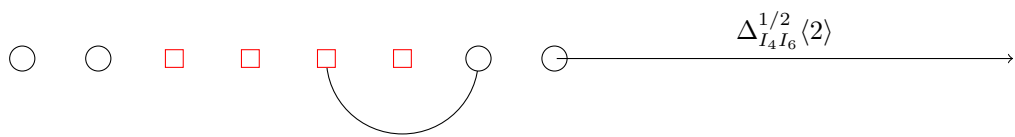
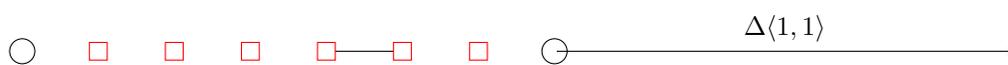
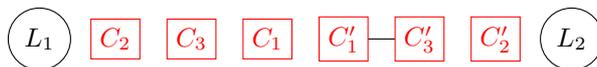
\begin{figure}[H]
\begin{center}
%
\end{center}
\end{figure}

Relation:

\begin{equation*}\begin{split}
\left\{\Gamma_{C_3'}  , \Gamma_{C_1'} \right\}=e
\end{split}
\end{equation*}
\subsubsection{Vertex Number 22}

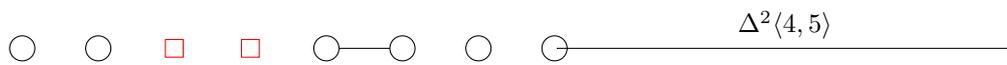
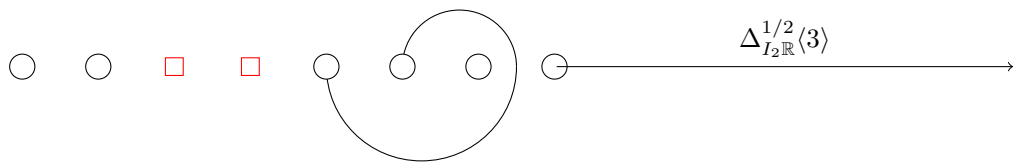
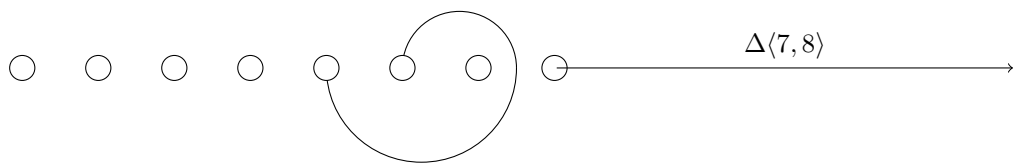
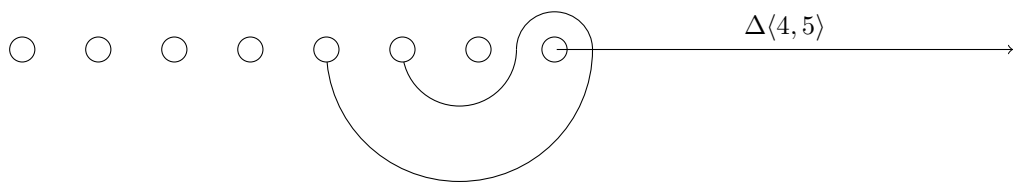
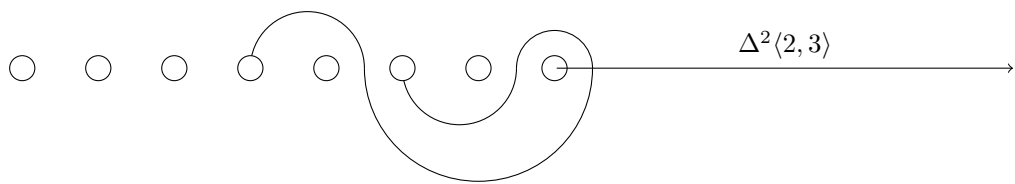
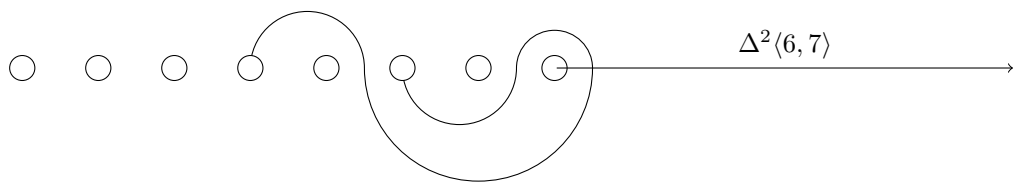
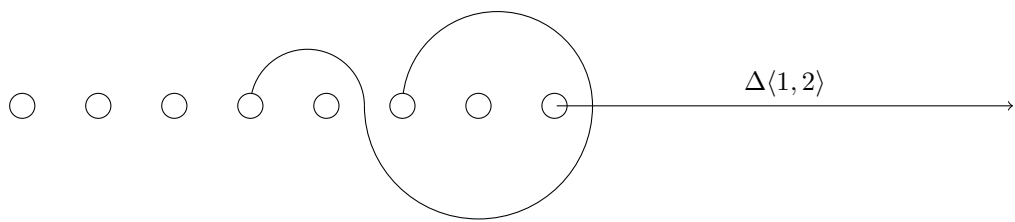
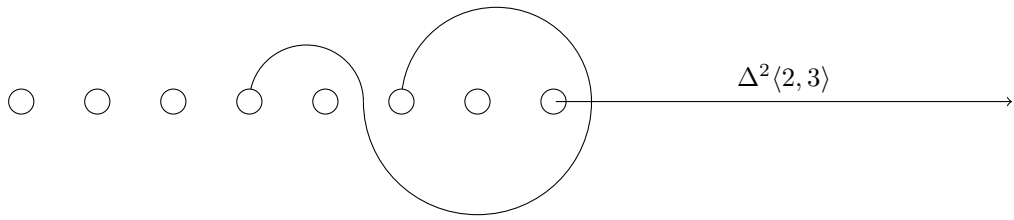
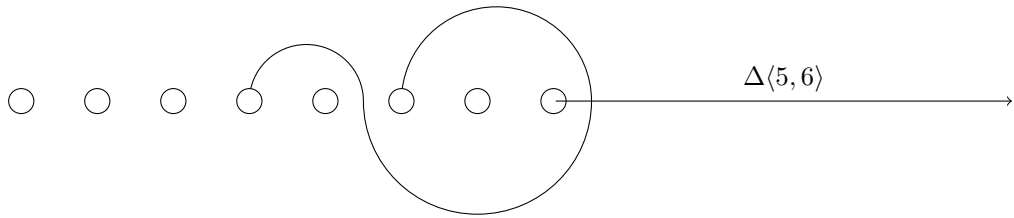
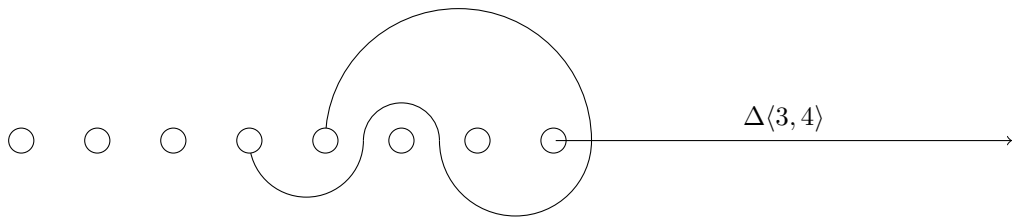
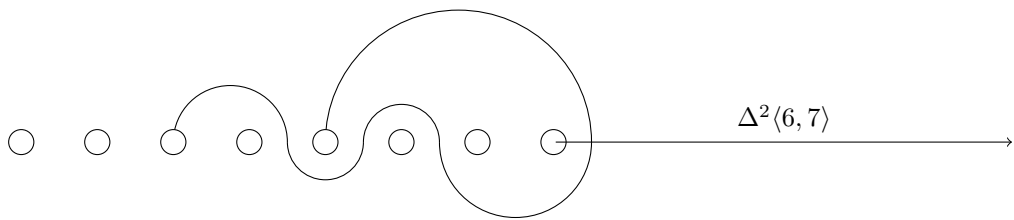
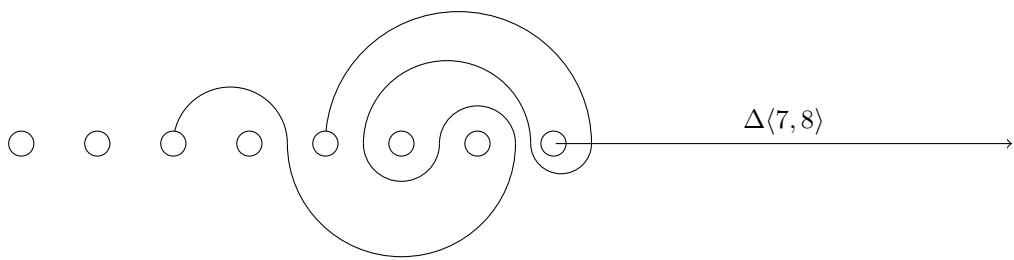
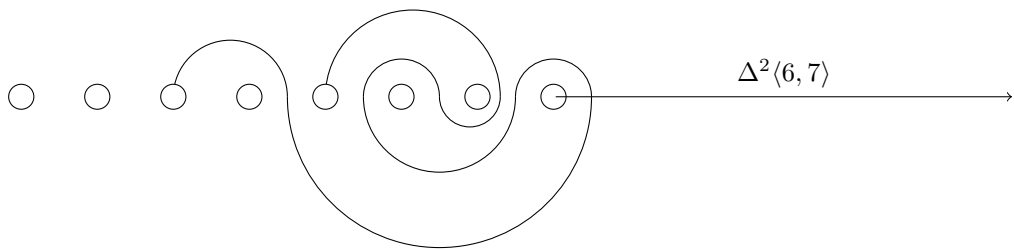
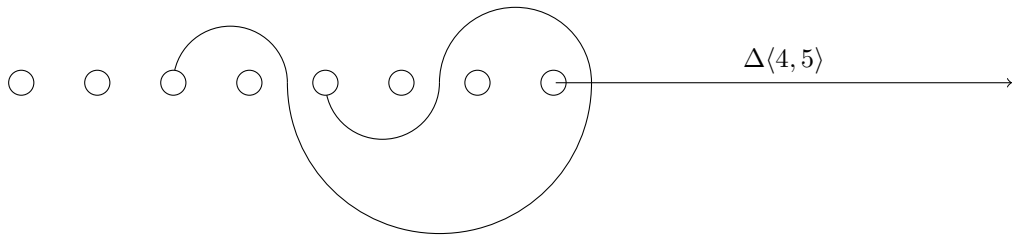
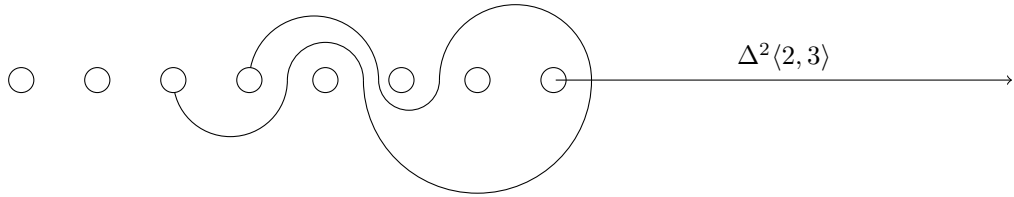
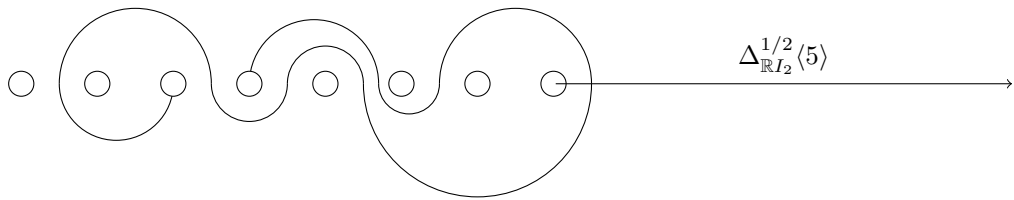
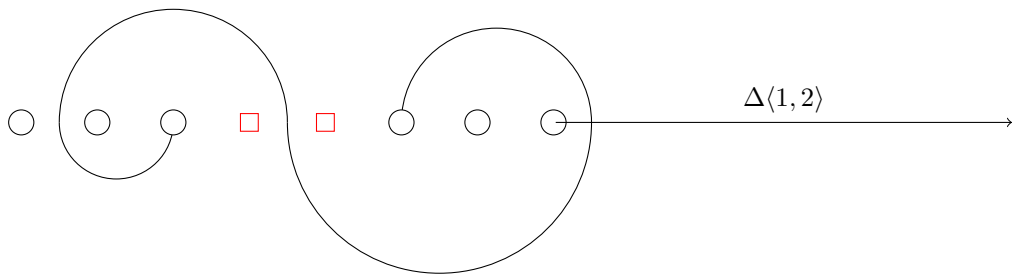
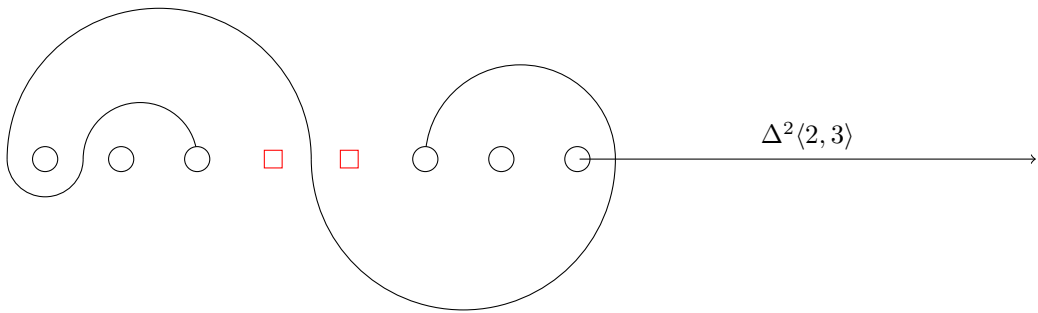
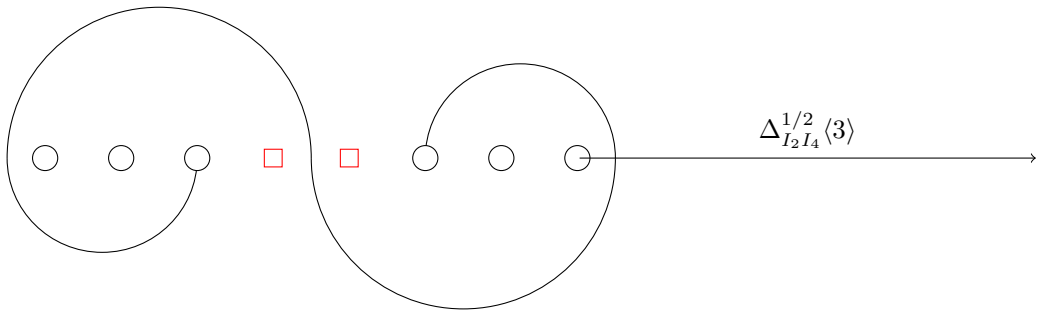
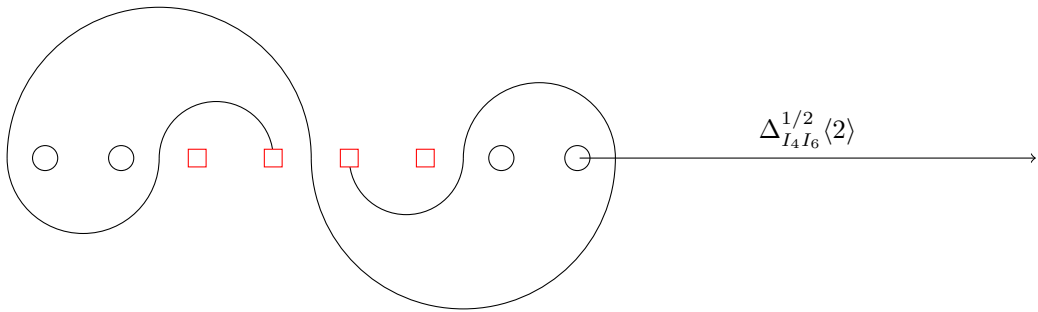
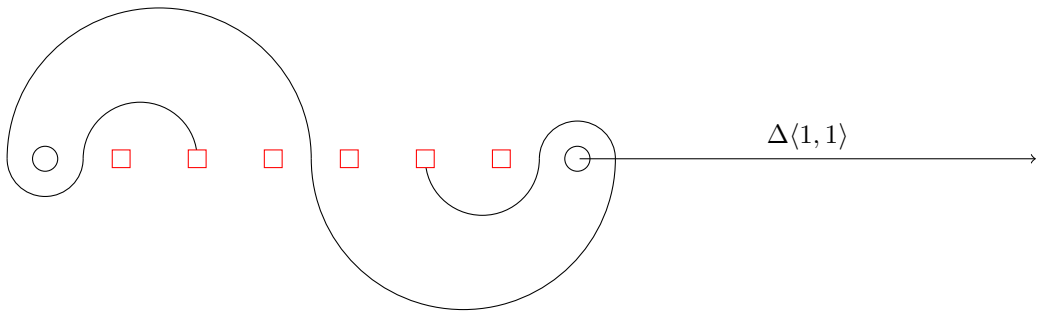
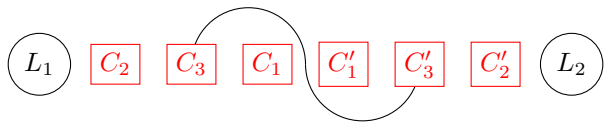
\begin{figure}[H]
\begin{center}
%
\end{center}
\end{figure}

Relation:

\begin{equation*}\begin{split}
\Gamma_{C_1}\Gamma_{C_3}\Gamma_{C_1}^{-1}  = \Gamma_{C_3'}
\end{split}
\end{equation*}
\subsubsection{Vertex Number 23}

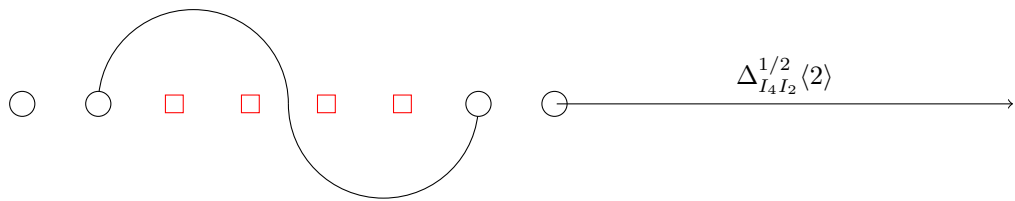
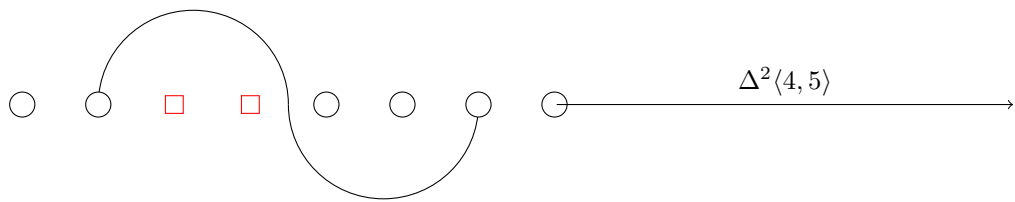
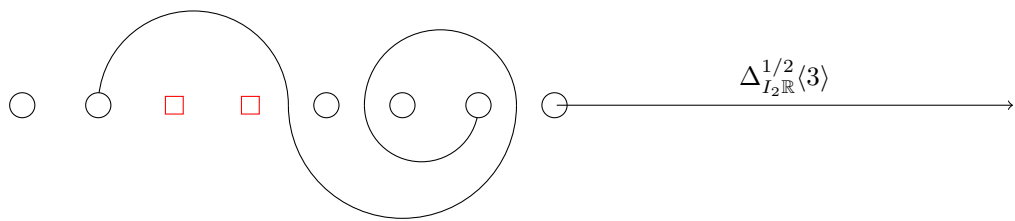
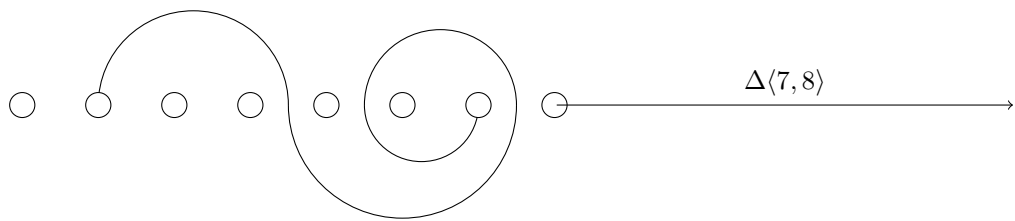
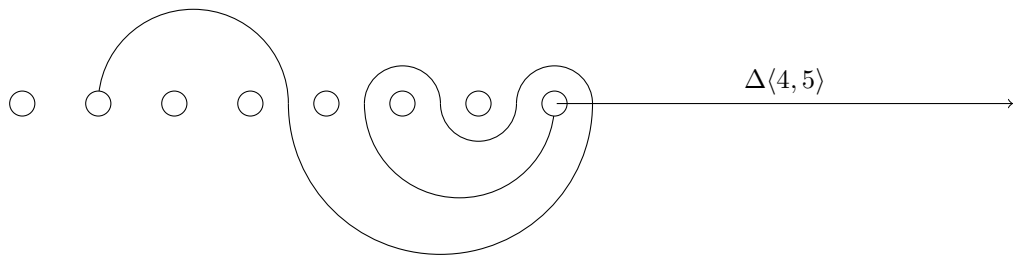
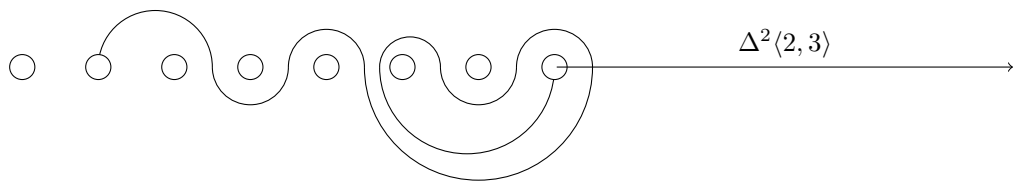
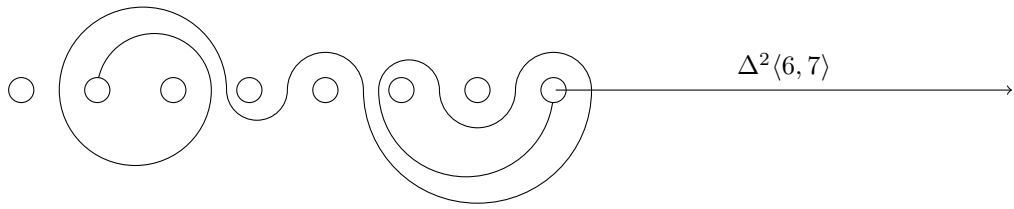
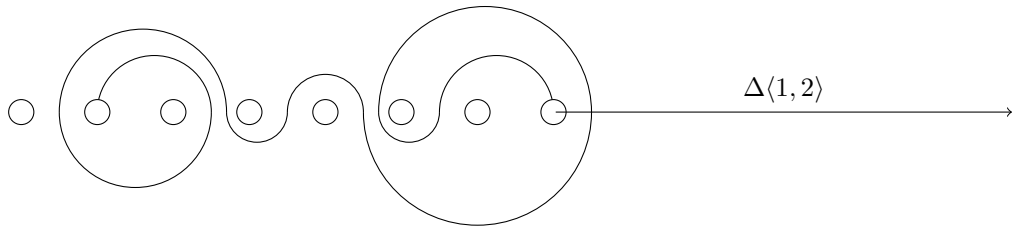
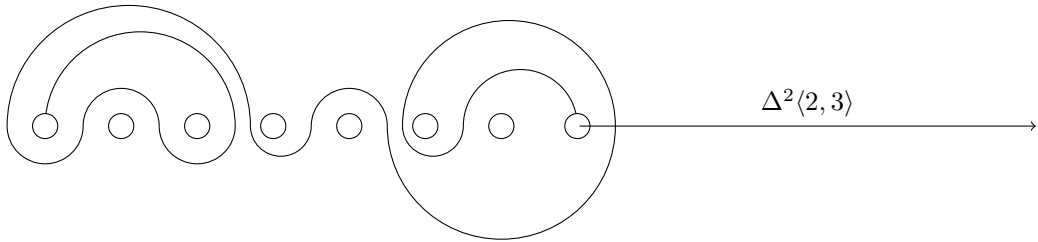
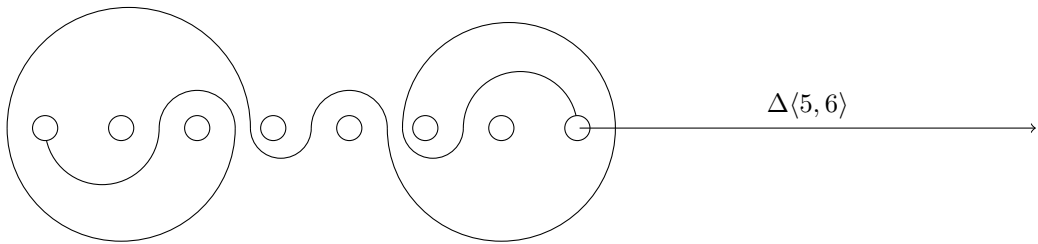
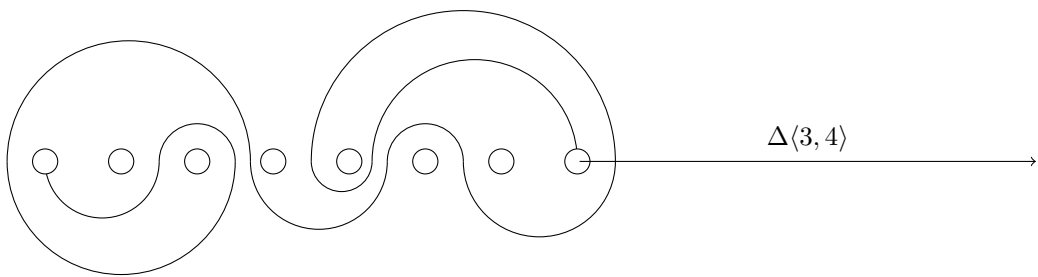
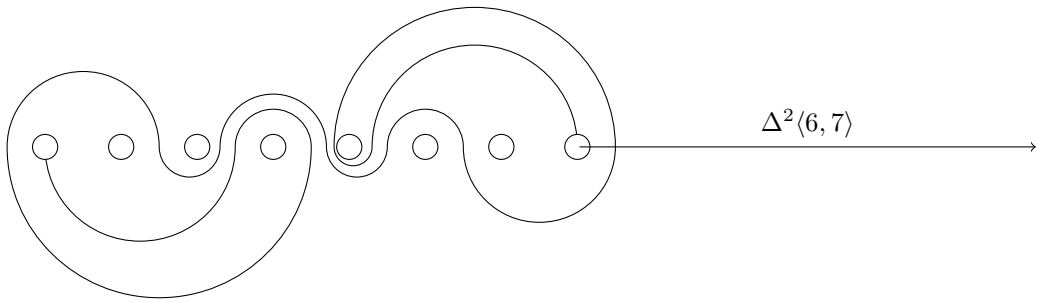
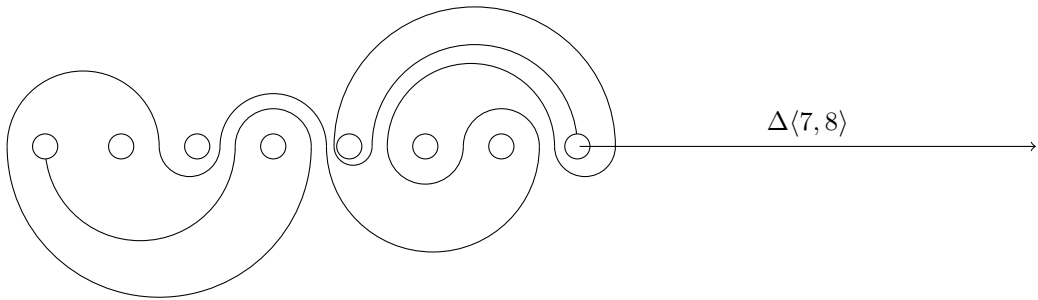
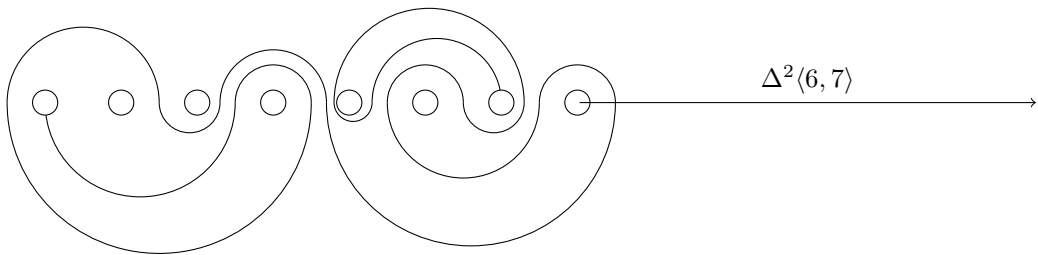
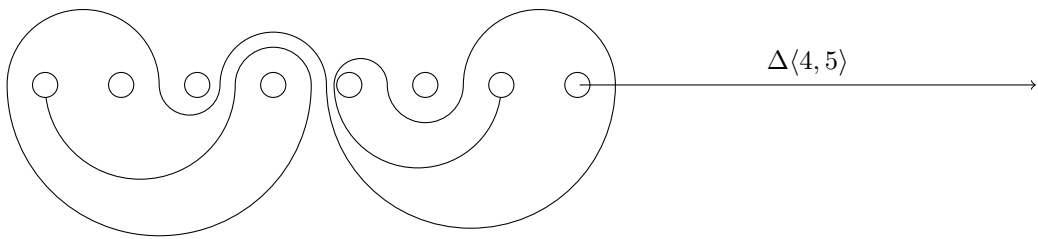
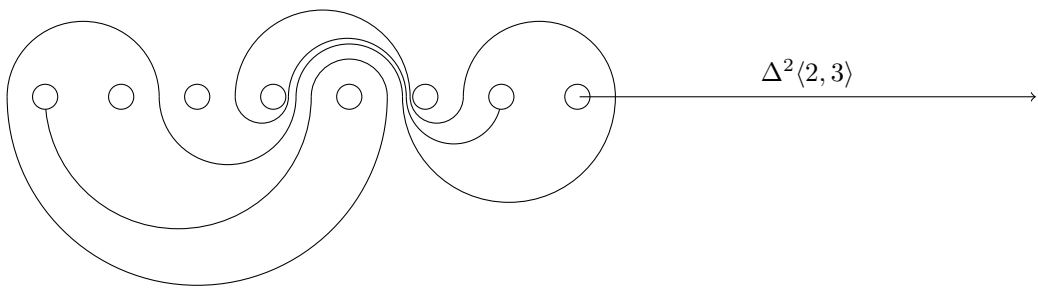
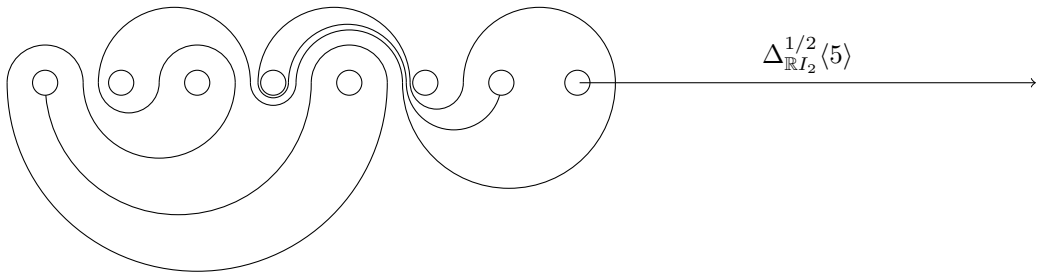
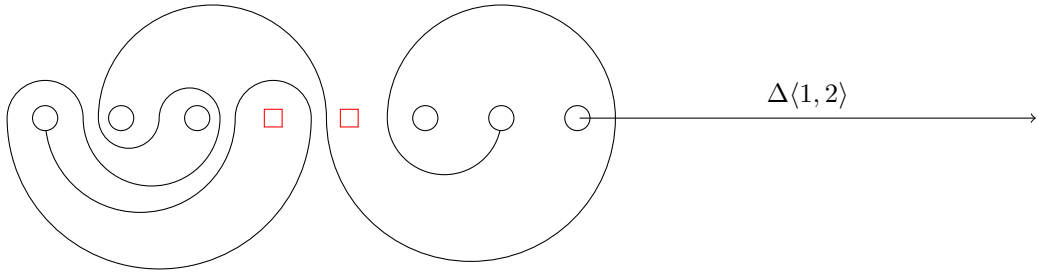
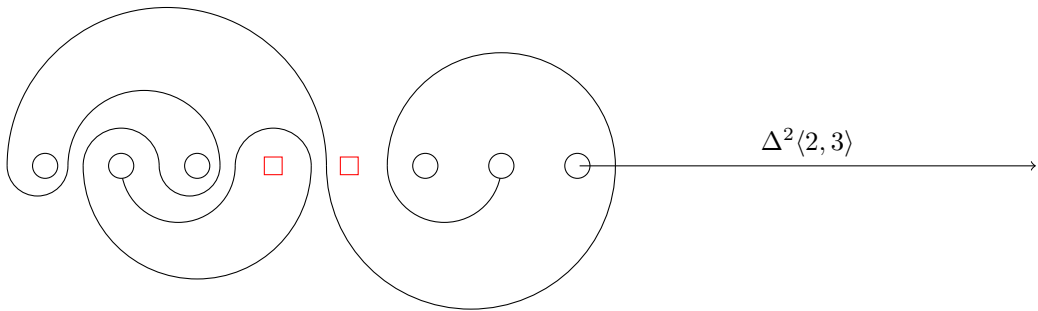
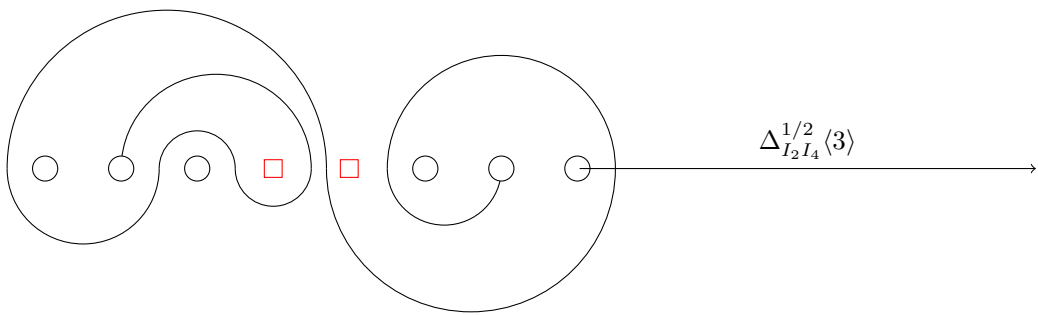
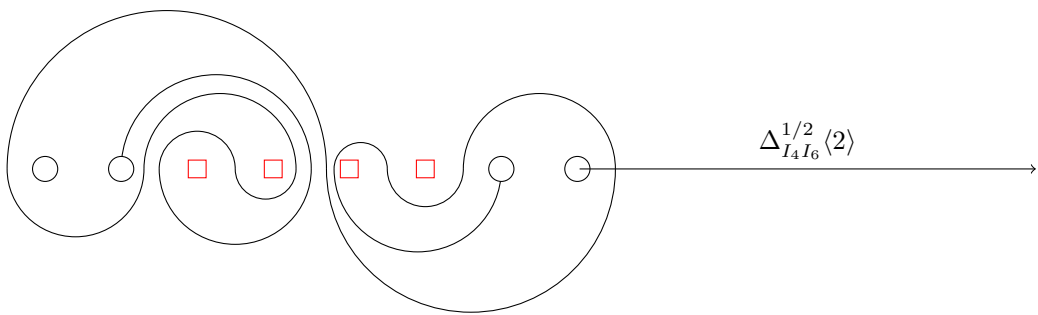
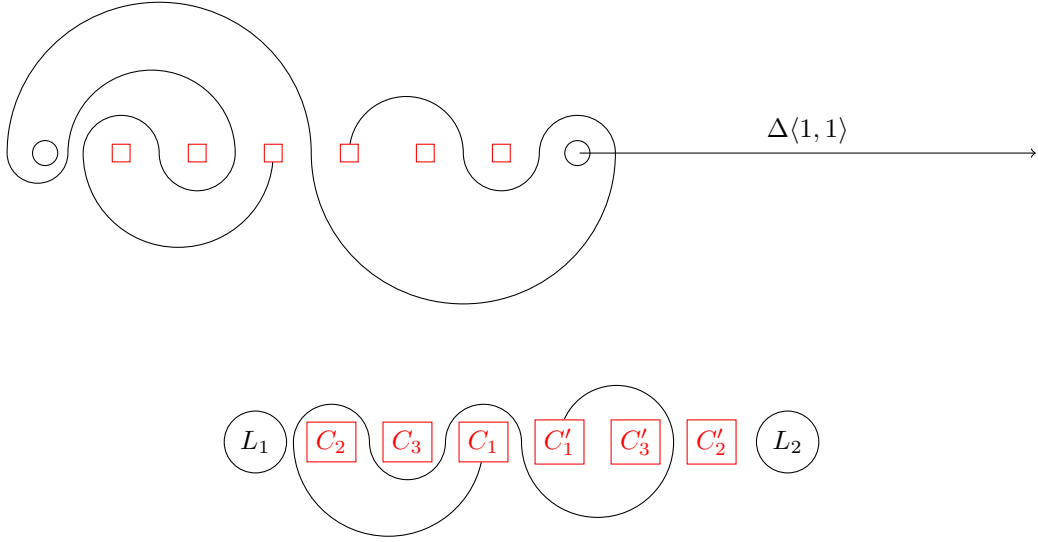
\begin{figure}[H]
\begin{center}
%
\end{center}
\end{figure}

Relation:

\begin{equation*}\begin{split}
\Gamma_{C_3'}^{-1}\Gamma_{C_1}\Gamma_{C_2}\Gamma_{C_1}\Gamma_{C_2}^{-1}\Gamma_{C_1}^{-1}\Gamma_{C_3'}  = \Gamma_{C_1'}
\end{split}
\end{equation*}
\subsubsection{Raw relations}

\begin{equation}\label{arr1auto_calc_vert_1_rel_1}
\begin{split}
\left[\Gamma_{C_1}\Gamma_{C_3}\Gamma_{C_2}\Gamma_{L_1}\Gamma_{C_2}^{-1}\Gamma_{C_3}^{-1}\Gamma_{C_1}^{-1} , \Gamma_{L_2}  \right]=e,
\end{split}
\end{equation}
\begin{equation}\label{arr1auto_calc_vert_2_rel_1}
\begin{split}
\Gamma_{C_1}\Gamma_{C_3}\Gamma_{C_2}\Gamma_{L_1}\Gamma_{C_2}^{-1}\Gamma_{C_3}^{-1}\Gamma_{C_1}^{-1}\Gamma_{L_2}^{-1}\Gamma_{C_1}\Gamma_{L_2}\Gamma_{C_1}\Gamma_{C_3}\Gamma_{C_2}\Gamma_{L_1}^{-1}\Gamma_{C_2}^{-1}\Gamma_{C_3}^{-1}\Gamma_{C_1}^{-1}  = \Gamma_{C_1'} ,
\end{split}
\end{equation}
\begin{equation}\label{arr1auto_calc_vert_3_rel_1}
\begin{split}
\Gamma_{C_1'}\Gamma_{C_1}\Gamma_{C_3}\Gamma_{C_2}\Gamma_{L_1}\Gamma_{C_2}^{-1}\Gamma_{C_3}^{-1}\Gamma_{C_1}^{-1}\Gamma_{L_2}^{-1}\Gamma_{C_3}\Gamma_{L_2}\Gamma_{C_1}\Gamma_{C_3}\Gamma_{C_2}\Gamma_{L_1}^{-1}\Gamma_{C_2}^{-1}\Gamma_{C_3}^{-1}\Gamma_{C_1}^{-1}\Gamma_{C_1'}^{-1}  = \Gamma_{C_3'} ,
\end{split}
\end{equation}
\begin{equation}\label{arr1auto_calc_vert_4_rel_1}
\begin{split}
\left\{\Gamma_{C_1} , \Gamma_{C_3}  \right\}=e,
\end{split}
\end{equation}
\begin{equation}\label{arr1auto_calc_vert_5_rel_1}
\begin{split}
\left[\Gamma_{C_3}^{-1}\Gamma_{C_1}^{-1}\Gamma_{L_2}\Gamma_{C_1}\Gamma_{C_3} , \Gamma_{C_1}  \right]=e,
\end{split}
\end{equation}
\begin{equation}\label{arr1auto_calc_vert_7_rel_1}
\begin{split}
\Gamma_{C_1'}\Gamma_{C_1}\Gamma_{C_3}\Gamma_{C_2}\Gamma_{L_1}\Gamma_{C_2}^{-1}\Gamma_{C_3}^{-1}\Gamma_{C_1}^{-1}\Gamma_{L_2}^{-1}\Gamma_{C_2}\Gamma_{L_2}\Gamma_{C_1}\Gamma_{C_3}\Gamma_{C_2}\Gamma_{L_1}^{-1}\Gamma_{C_2}^{-1}\Gamma_{C_3}^{-1}\Gamma_{C_1}^{-1}\Gamma_{C_1'}^{-1}  = \Gamma_{C_2'} ,
\end{split}
\end{equation}
\begin{equation}\label{arr1auto_calc_vert_8_rel_1}
\begin{split}
\left\{\Gamma_{C_1}^{-1}\Gamma_{L_2}\Gamma_{C_1} , \Gamma_{C_3}  \right\}=e,
\end{split}
\end{equation}
\begin{equation}\label{arr1auto_calc_vert_9_rel_1}
\begin{split}
\left[\Gamma_{L_2}\Gamma_{C_1}\Gamma_{C_3}\Gamma_{C_2}\Gamma_{L_1}^{-1}\Gamma_{C_2}^{-1}\Gamma_{C_3}^{-1}\Gamma_{C_1}^{-1}\Gamma_{C_1'}^{-1}\Gamma_{C_3'}\Gamma_{C_1'}\Gamma_{C_1}\Gamma_{C_3}\Gamma_{C_2}\Gamma_{L_1}\Gamma_{C_2}^{-1}\Gamma_{C_3}^{-1}\Gamma_{C_1}^{-1}\Gamma_{L_2}^{-1} , \Gamma_{C_2}  \right]=e,
\end{split}
\end{equation}
\begin{equation}\label{arr1auto_calc_vert_10_rel_1}
\begin{split}
\left\{\Gamma_{C_2'} , \Gamma_{C_1'}  \right\}=e,
\end{split}
\end{equation}
\begin{equation}\label{arr1auto_calc_vert_11_rel_1}
\begin{split}
\left[\Gamma_{C_2}^{-1}\Gamma_{C_3}^{-1}\Gamma_{C_1}^{-1}\Gamma_{C_2'}\Gamma_{C_1'}\Gamma_{C_2'}^{-1}\Gamma_{C_1}\Gamma_{C_3}\Gamma_{C_2} , \Gamma_{L_1}  \right]=e,
\end{split}
\end{equation}
\begin{equation}\label{arr1auto_calc_vert_12_rel_1}
\begin{split}
\left\{\Gamma_{C_2}^{-1}\Gamma_{C_3}^{-1}\Gamma_{C_1}^{-1}\Gamma_{C_2'}\Gamma_{C_1}\Gamma_{C_3}\Gamma_{C_2} , \Gamma_{L_1}  \right\}=e,
\end{split}
\end{equation}
\begin{equation}\label{arr1auto_calc_vert_13_rel_1}
\begin{split}
\left[\Gamma_{L_2}\Gamma_{C_1}\Gamma_{C_3}\Gamma_{C_1}^{-1}\Gamma_{L_2}^{-1} , \Gamma_{C_2}  \right]=e,
\end{split}
\end{equation}
\begin{equation}\label{arr1auto_calc_vert_14_rel_1}
\begin{split}
\left[\Gamma_{C_1}\Gamma_{C_3}\Gamma_{C_2}\Gamma_{L_1}^{-1}\Gamma_{C_2}^{-1}\Gamma_{C_3}^{-1}\Gamma_{C_1}^{-1}\Gamma_{C_1'}^{-1}\Gamma_{C_3'}\Gamma_{C_1'}\Gamma_{C_1}\Gamma_{C_3}\Gamma_{C_2}\Gamma_{L_1}\Gamma_{C_2}^{-1}\Gamma_{C_3}^{-1}\Gamma_{C_1}^{-1} , \Gamma_{C_2'}  \right]=e,
\end{split}
\end{equation}
\begin{equation}\label{arr1auto_calc_vert_15_rel_1}
\begin{split}
\left\{\Gamma_{L_2} , \Gamma_{C_2}  \right\}=e,
\end{split}
\end{equation}
\begin{equation}\label{arr1auto_calc_vert_16_rel_1}
\begin{split}
\left[\Gamma_{C_2}^{-1}\Gamma_{C_1}\Gamma_{C_2} , \Gamma_{L_2}  \right]=e,
\end{split}
\end{equation}
\begin{equation}\label{arr1auto_calc_vert_17_rel_1}
\begin{split}
\left\{\Gamma_{C_2}^{-1}\Gamma_{C_3}^{-1}\Gamma_{C_1}^{-1}\Gamma_{C_1'}^{-1}\Gamma_{C_3'}\Gamma_{C_1'}\Gamma_{C_1}\Gamma_{C_3}\Gamma_{C_2} , \Gamma_{L_1}  \right\}=e,
\end{split}
\end{equation}
\begin{equation}\label{arr1auto_calc_vert_18_rel_1}
\begin{split}
\left\{\Gamma_{C_1} , \Gamma_{C_2}  \right\}=e,
\end{split}
\end{equation}
\begin{equation}\label{arr1auto_calc_vert_19_rel_1}
\begin{split}
\left[\Gamma_{C_1}\Gamma_{C_3}\Gamma_{C_1}^{-1} , \Gamma_{C_2'}  \right]=e,
\end{split}
\end{equation}
\begin{equation}\label{arr1auto_calc_vert_20_rel_1}
\begin{split}
\left[\Gamma_{C_3'}\Gamma_{C_1'}\Gamma_{C_1}\Gamma_{C_3}\Gamma_{C_2}\Gamma_{L_1}\Gamma_{C_2}^{-1}\Gamma_{C_3}^{-1}\Gamma_{C_1}^{-1}\Gamma_{C_1'}^{-1}\Gamma_{C_3'}^{-1} , \Gamma_{C_1'}  \right]=e,
\end{split}
\end{equation}
\begin{equation}\label{arr1auto_calc_vert_21_rel_1}
\begin{split}
\Gamma_{C_1}\Gamma_{C_2}\Gamma_{C_1}^{-1}  = \Gamma_{C_2'} ,
\end{split}
\end{equation}
\begin{equation}\label{arr1auto_calc_vert_22_rel_1}
\begin{split}
\left\{\Gamma_{C_3'} , \Gamma_{C_1'}  \right\}=e,
\end{split}
\end{equation}
\begin{equation}\label{arr1auto_calc_vert_23_rel_1}
\begin{split}
\Gamma_{C_1}\Gamma_{C_3}\Gamma_{C_1}^{-1}  = \Gamma_{C_3'} ,
\end{split}
\end{equation}
\begin{equation}\label{arr1auto_calc_vert_24_rel_1}
\begin{split}
\Gamma_{C_3'}^{-1}\Gamma_{C_1}\Gamma_{C_2}\Gamma_{C_1}\Gamma_{C_2}^{-1}\Gamma_{C_1}^{-1}\Gamma_{C_3'}  = \Gamma_{C_1'} ,
\end{split}
\end{equation}
\begin{equation}\label{arr1auto_calc_projective_rel}
\begin{split}
\Gamma_{L_2}\Gamma_{C_2'}\Gamma_{C_3'}\Gamma_{C_1'}\Gamma_{C_1}\Gamma_{C_3}\Gamma_{C_2}\Gamma_{L_1} =e.
\end{split}
\end{equation}

\subsection{Computation of $ \pcpt{\mathcal{B}_2} $}

\begin{figure}[H]
\begin{center}
\begin{scaletikzpicturetowidth}{\textwidth}

\end{scaletikzpicturetowidth}
\end{center}
\end{figure}

{Monodromy table}

\begin{tabular}{| c | c | c | c |}
\hline \hline
Vertex number& Vertex description& Skeleton& Diffeomorphism\\ \hline \hline
1& $C_1$ branch& $ \langle 2 - 3 \rangle $& $ \Delta_{I_{4}I_{6}}^{1/2}\langle 2 \rangle $ \\ \hline
2& $C_3$ branch& $ \langle 3 - 4 \rangle $& $ \Delta_{I_{2}I_{4}}^{1/2}\langle 3 \rangle $ \\ \hline
3& Tangency between $C_3$ and $C_1$& $ \langle 2 - 3 \rangle $& $ \Delta^{2}\langle 2, 3 \rangle $ \\ \hline
4& $C_2$ branch& $ \langle 5 - 6 \rangle $& $ \Delta_{\mathbb{R}I_{2}}^{1/2}\langle 5 \rangle $ \\ \hline
5& Node between $C_2$ and $C_3$& $ \langle 4 - 5 \rangle $& $ \Delta\langle 4, 5 \rangle $ \\ \hline
6& Tangency between $C_1$ and $C_2$& $ \langle 6 - 7 \rangle $& $ \Delta^{2}\langle 6, 7 \rangle $ \\ \hline
7& Node between $L_1$ and $C_1$& $ \langle 7 - 8 \rangle $& $ \Delta\langle 7, 8 \rangle $ \\
\hline \hline
\end{tabular}

{Braids calculation}

\subsubsection{Vertex Number 1}

\begin{figure}[H]
\begin{center}
%
\end{center}
\end{figure}

Relation:

\begin{equation*}\begin{split}
\Gamma_{C_1}  = \Gamma_{C_1'}
\end{split}
\end{equation*}
\subsubsection{Vertex Number 2}

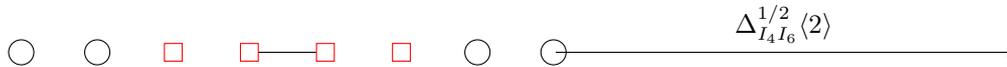
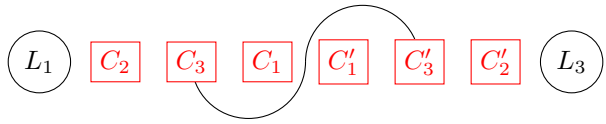
\begin{figure}[H]
\begin{center}
%
\end{center}
\end{figure}

Relation:

\begin{equation*}\begin{split}
\Gamma_{C_1'}\Gamma_{C_3}\Gamma_{C_1'}^{-1}  = \Gamma_{C_3'}
\end{split}
\end{equation*}
\subsubsection{Vertex Number 3}

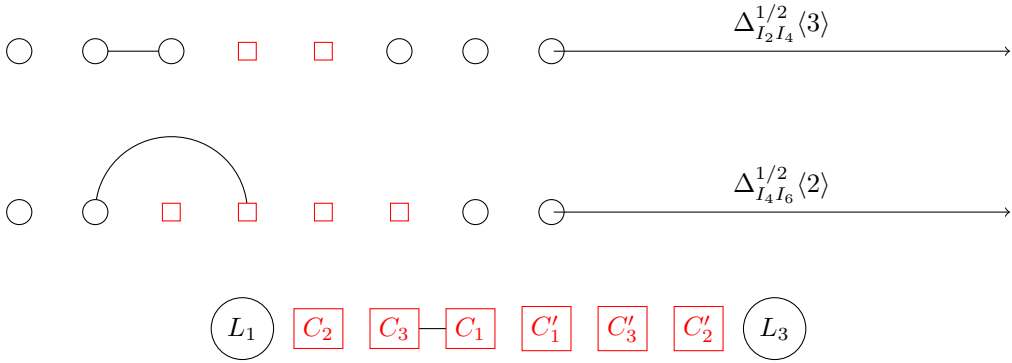
\begin{figure}[H]
\begin{center}
%
\end{center}
\end{figure}

Relation:

\begin{equation*}\begin{split}
\left\{\Gamma_{C_1}  , \Gamma_{C_3} \right\}=e
\end{split}
\end{equation*}
\subsubsection{Vertex Number 4}

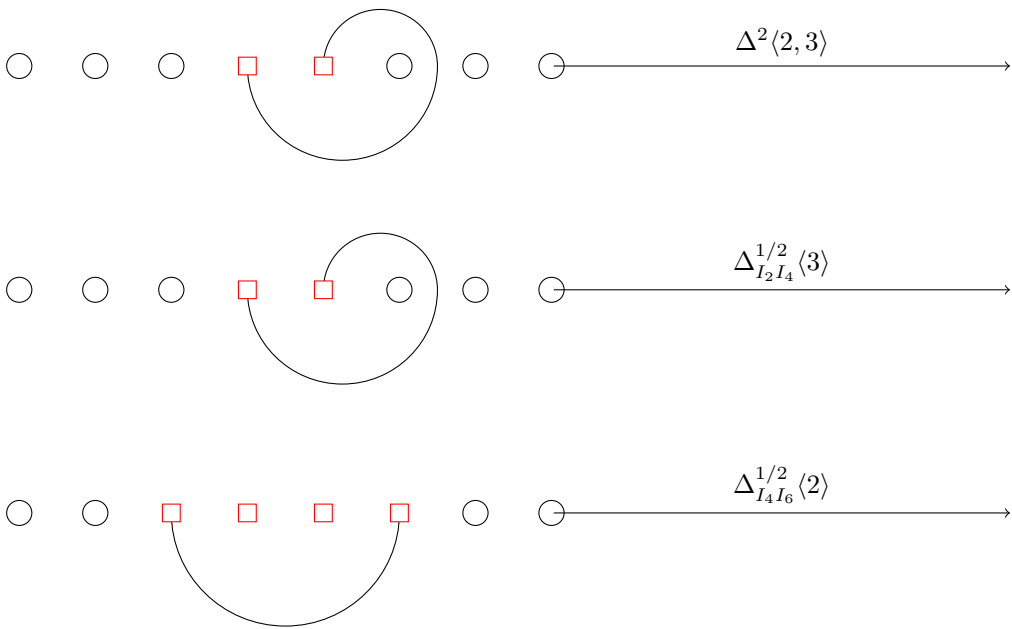
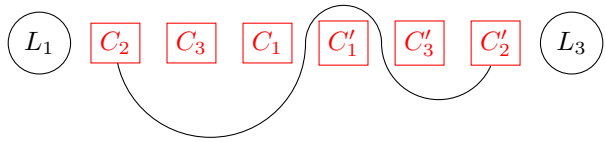
\begin{figure}[H]
\begin{center}
%
\end{center}
\end{figure}

Relation:

\begin{equation*}\begin{split}
\Gamma_{C_1'}\Gamma_{C_2}\Gamma_{C_1'}^{-1}  = \Gamma_{C_2'}
\end{split}
\end{equation*}
\subsubsection{Vertex Number 5}

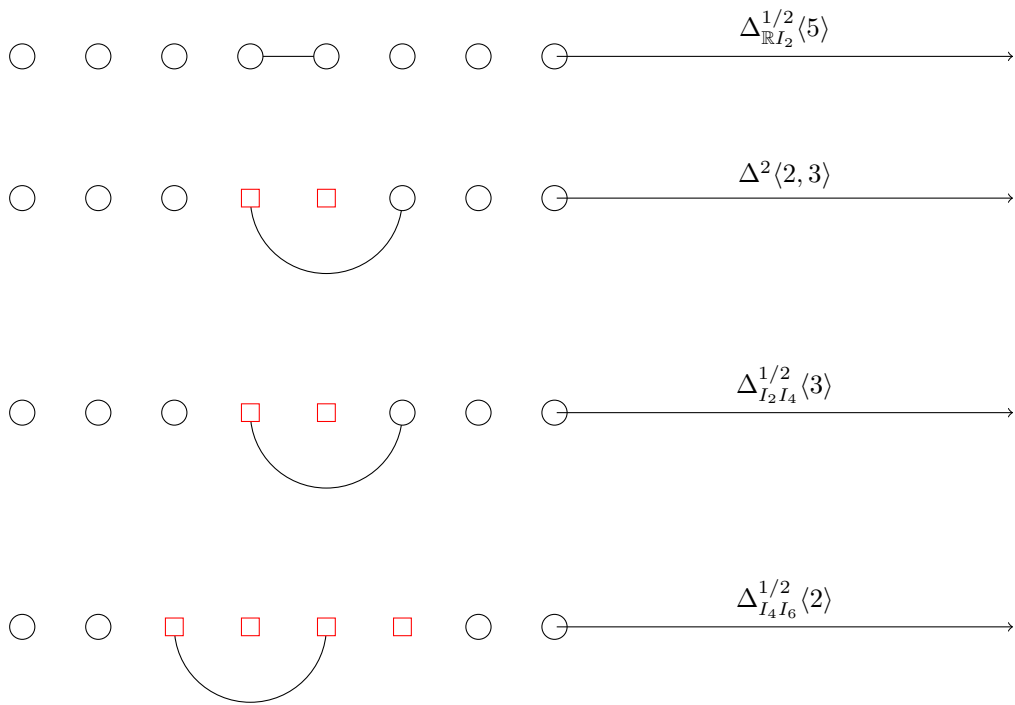
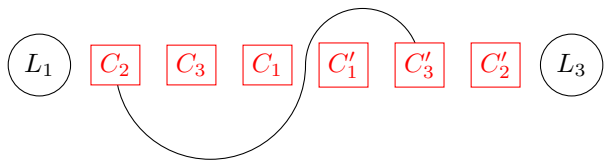
\begin{figure}[H]
\begin{center}
%
\end{center}
\end{figure}

Relation:

\begin{equation*}\begin{split}
\left[\Gamma_{C_1'}^{-1}\Gamma_{C_3'}\Gamma_{C_1'} , \Gamma_{C_2} \right]=e
\end{split}
\end{equation*}
\subsubsection{Vertex Number 6}

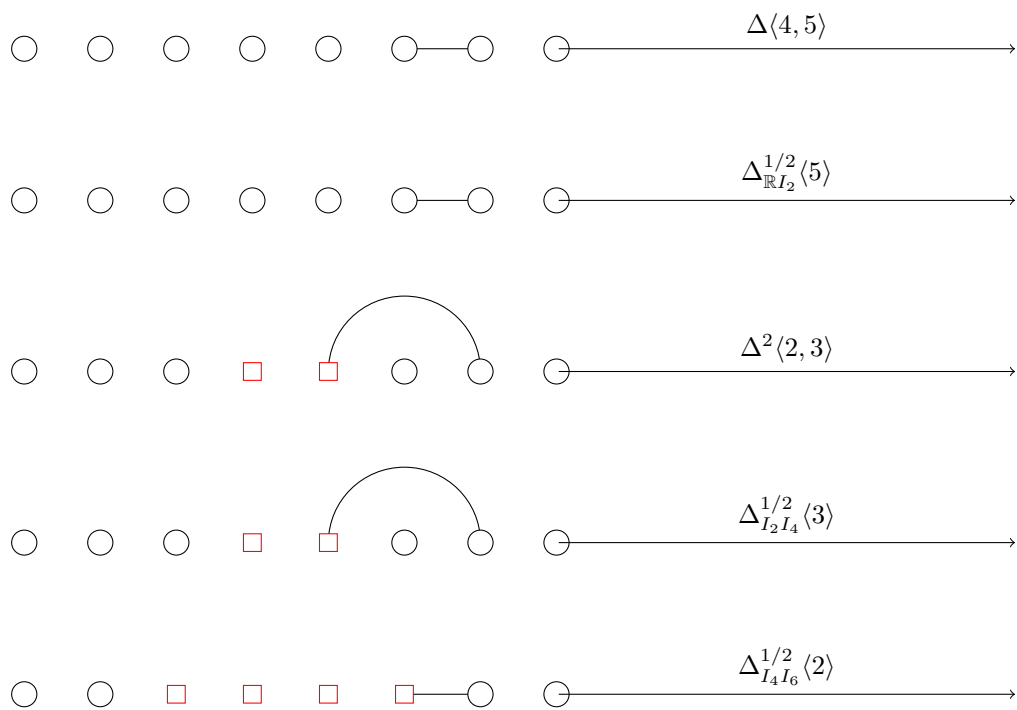
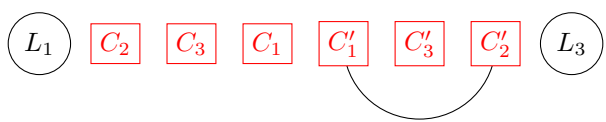
\begin{figure}[H]
\begin{center}
%
\end{center}
\end{figure}

Relation:

\begin{equation*}\begin{split}
\left\{\Gamma_{C_2'}  , \Gamma_{C_1'} \right\}=e
\end{split}
\end{equation*}
\subsubsection{Vertex Number 7}

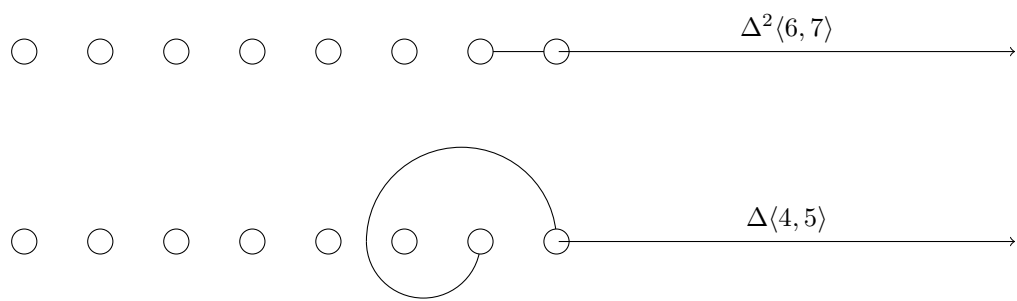
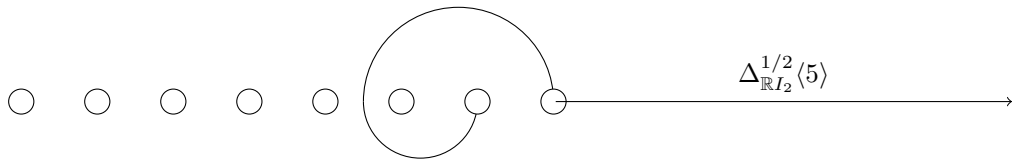
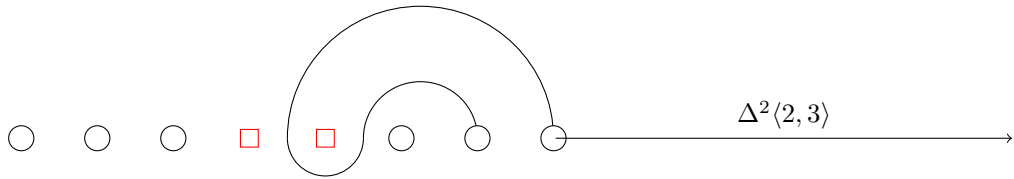
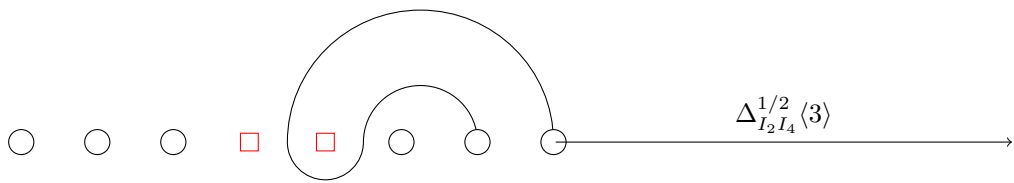
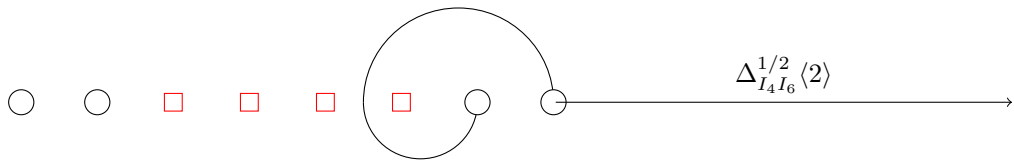
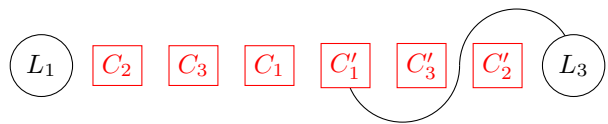
\begin{figure}[H]
\begin{center}
%
\end{center}
\end{figure}

Relation:

\begin{equation*}\begin{split}
\left[\Gamma_{C_2'}\Gamma_{C_1'}\Gamma_{C_2'}^{-1} , \Gamma_{L_3} \right]=e
\end{split}
\end{equation*}
\subsubsection{Vertex Number 8}

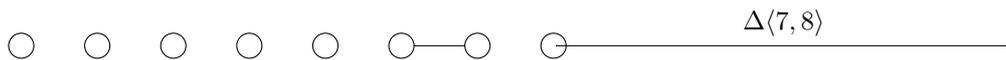
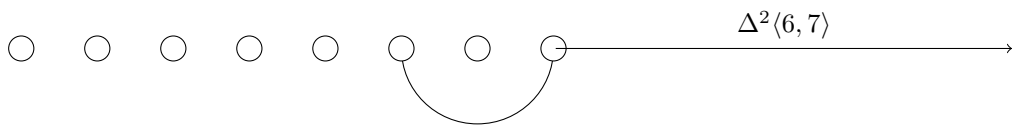
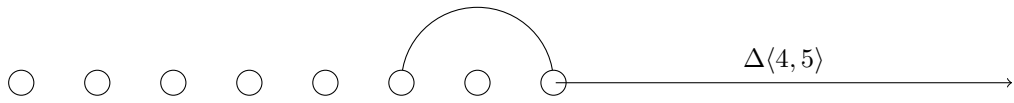
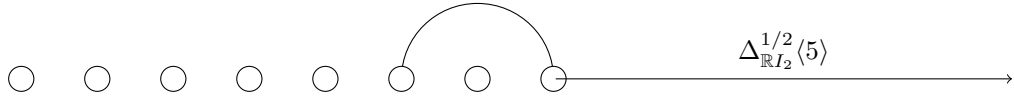
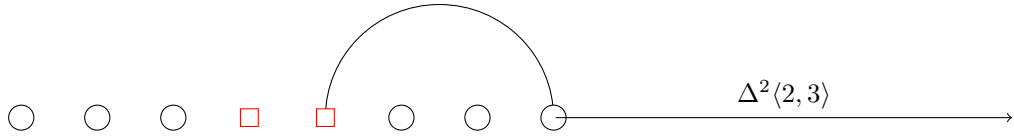
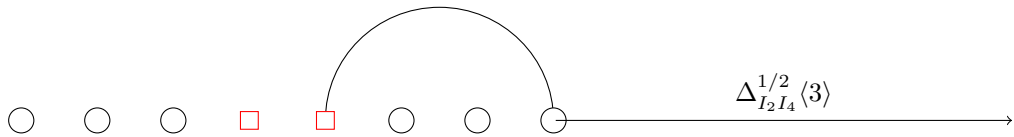
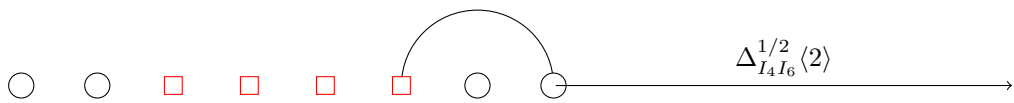
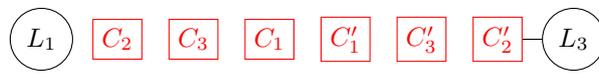
\begin{figure}[H]
\begin{center}
%
\end{center}
\end{figure}

Relation:

\begin{equation*}\begin{split}
\left\{\Gamma_{C_2'}  , \Gamma_{L_3} \right\}=e
\end{split}
\end{equation*}
\subsubsection{Vertex Number 9}

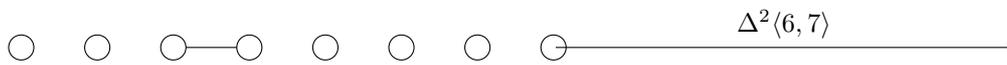
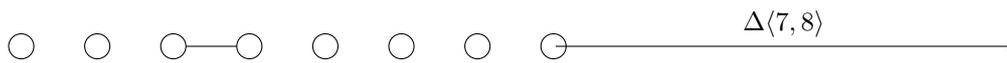
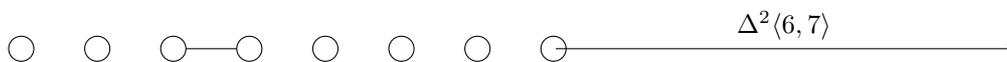
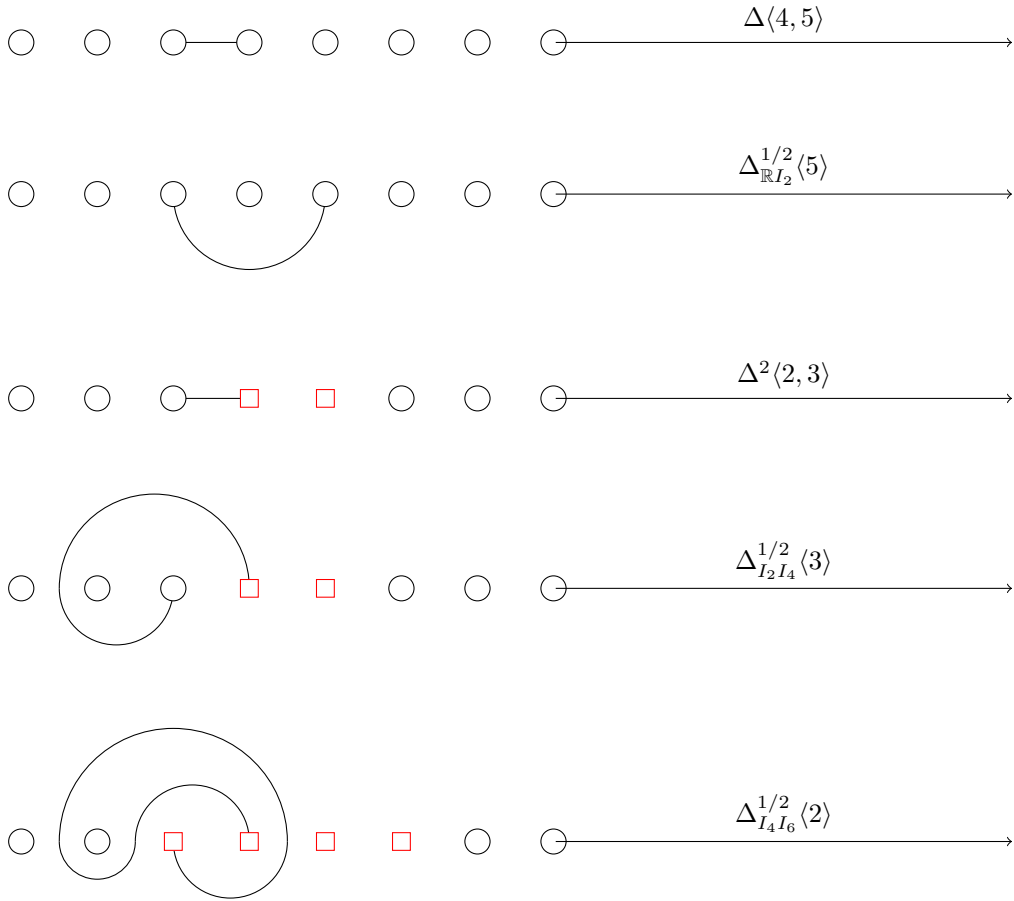
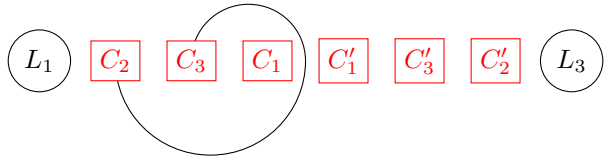
\begin{figure}[H]
\begin{center}
%
\end{center}
\end{figure}

Relation:

\begin{equation*}\begin{split}
\left[\Gamma_{C_1}\Gamma_{C_3}\Gamma_{C_1}^{-1} , \Gamma_{C_2} \right]=e
\end{split}
\end{equation*}
\subsubsection{Vertex Number 10}

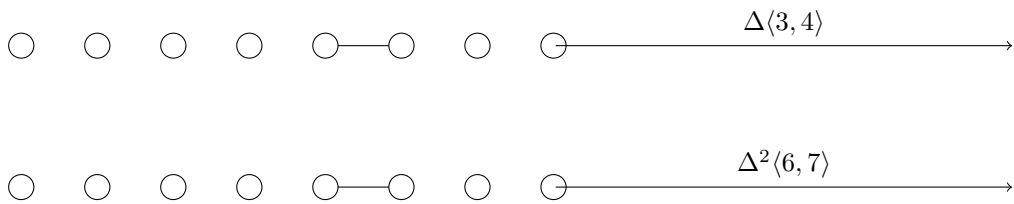
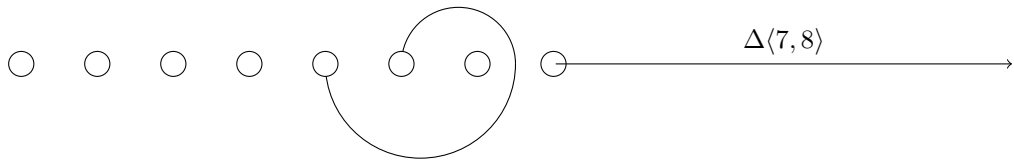
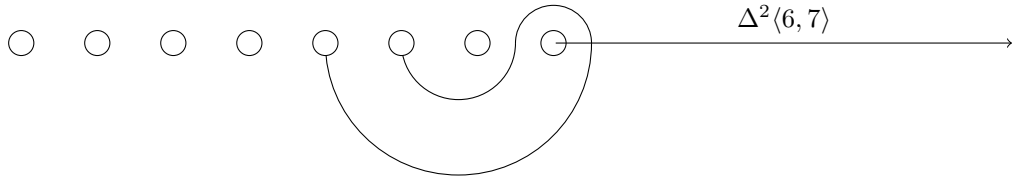
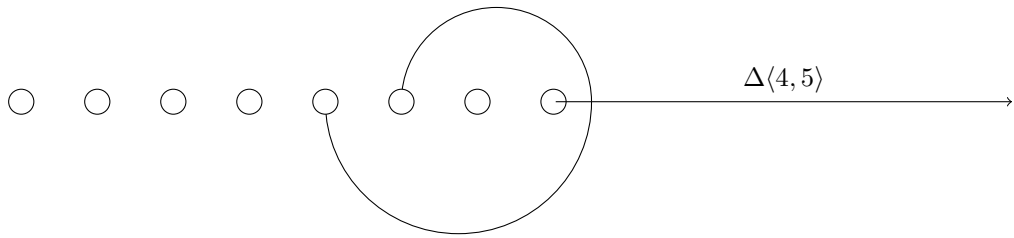
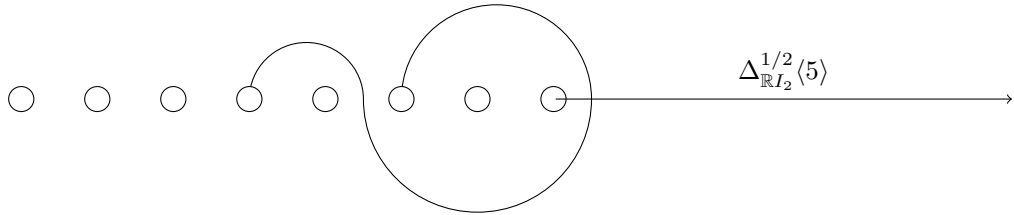
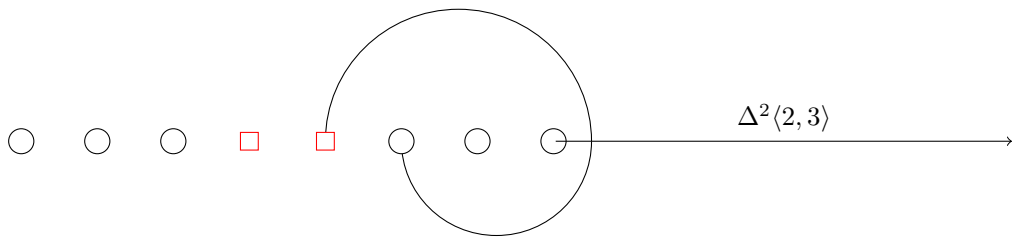
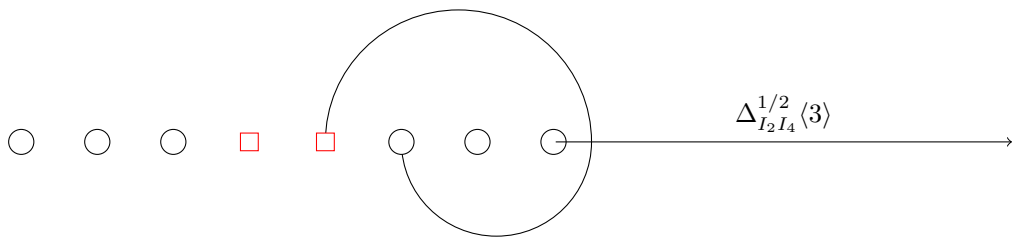
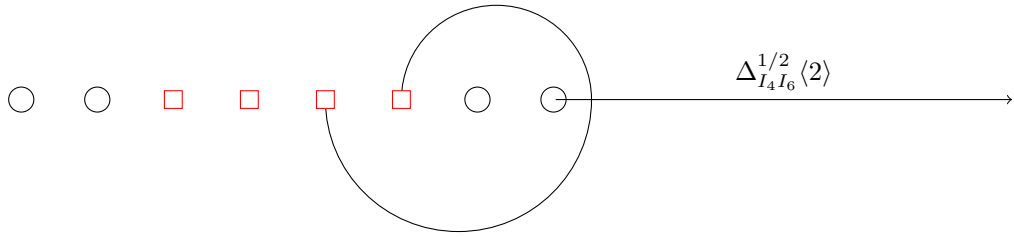
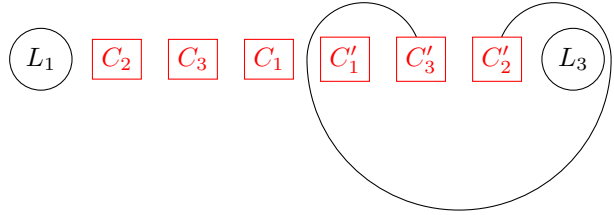
\begin{figure}[H]
\begin{center}
%
\end{center}
\end{figure}

Relation:

\begin{equation*}\begin{split}
\left[\Gamma_{L_3}^{-1}\Gamma_{C_1'}^{-1}\Gamma_{C_3'}\Gamma_{C_1'}\Gamma_{L_3} , \Gamma_{C_2'} \right]=e
\end{split}
\end{equation*}
\subsubsection{Vertex Number 11}

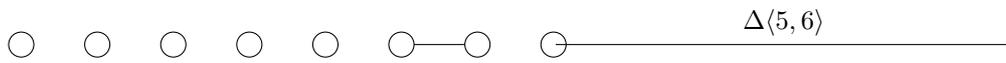
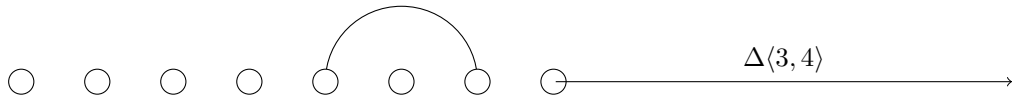
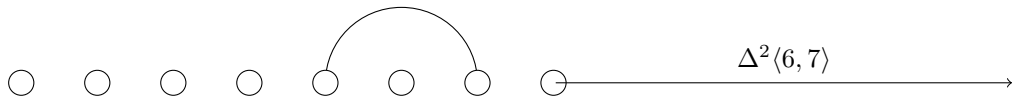
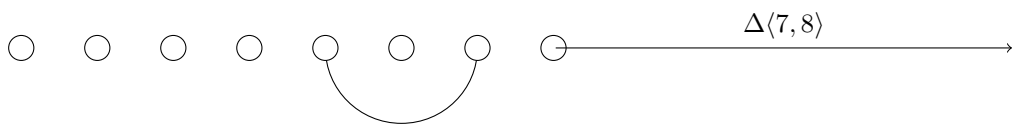
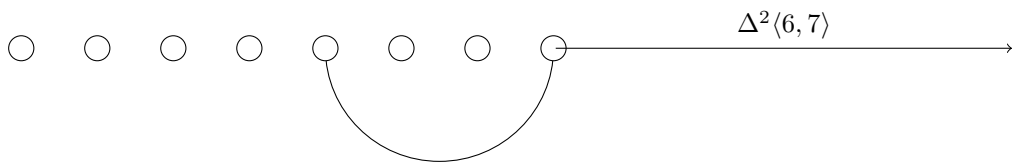
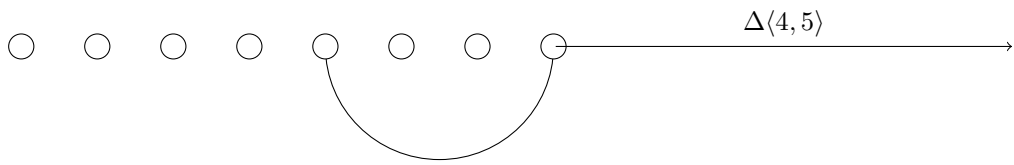
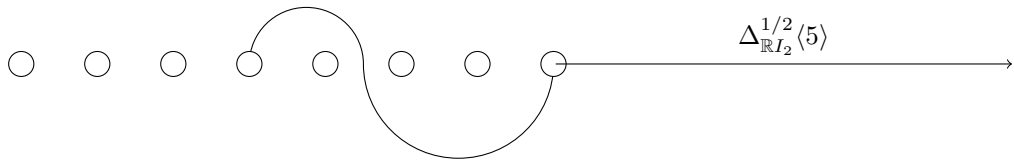
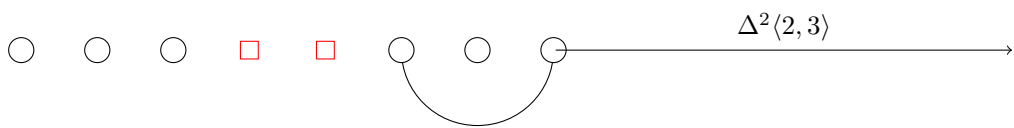
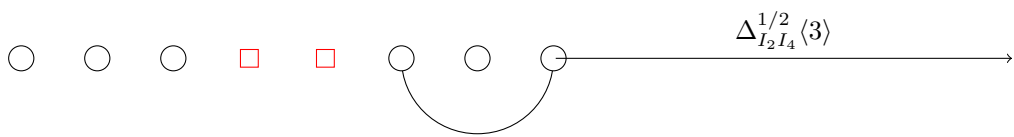
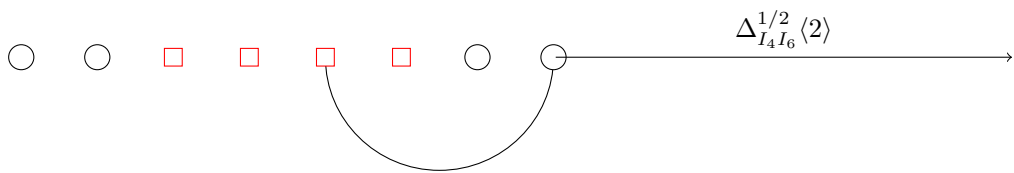
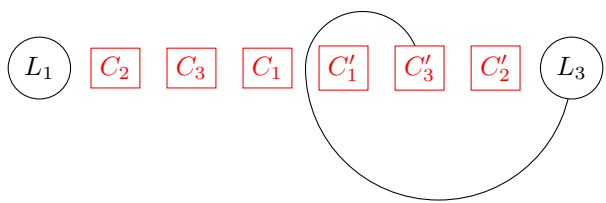
\begin{figure}[H]
\begin{center}
%
\end{center}
\end{figure}

Relation:

\begin{equation*}\begin{split}
\left\{\Gamma_{C_1'}^{-1}\Gamma_{C_3'}\Gamma_{C_1'}  , \Gamma_{L_3} \right\}=e
\end{split}
\end{equation*}
\subsubsection{Vertex Number 12}

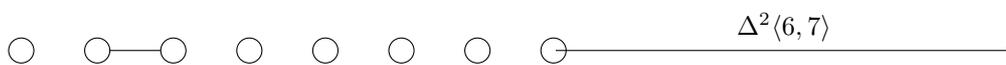
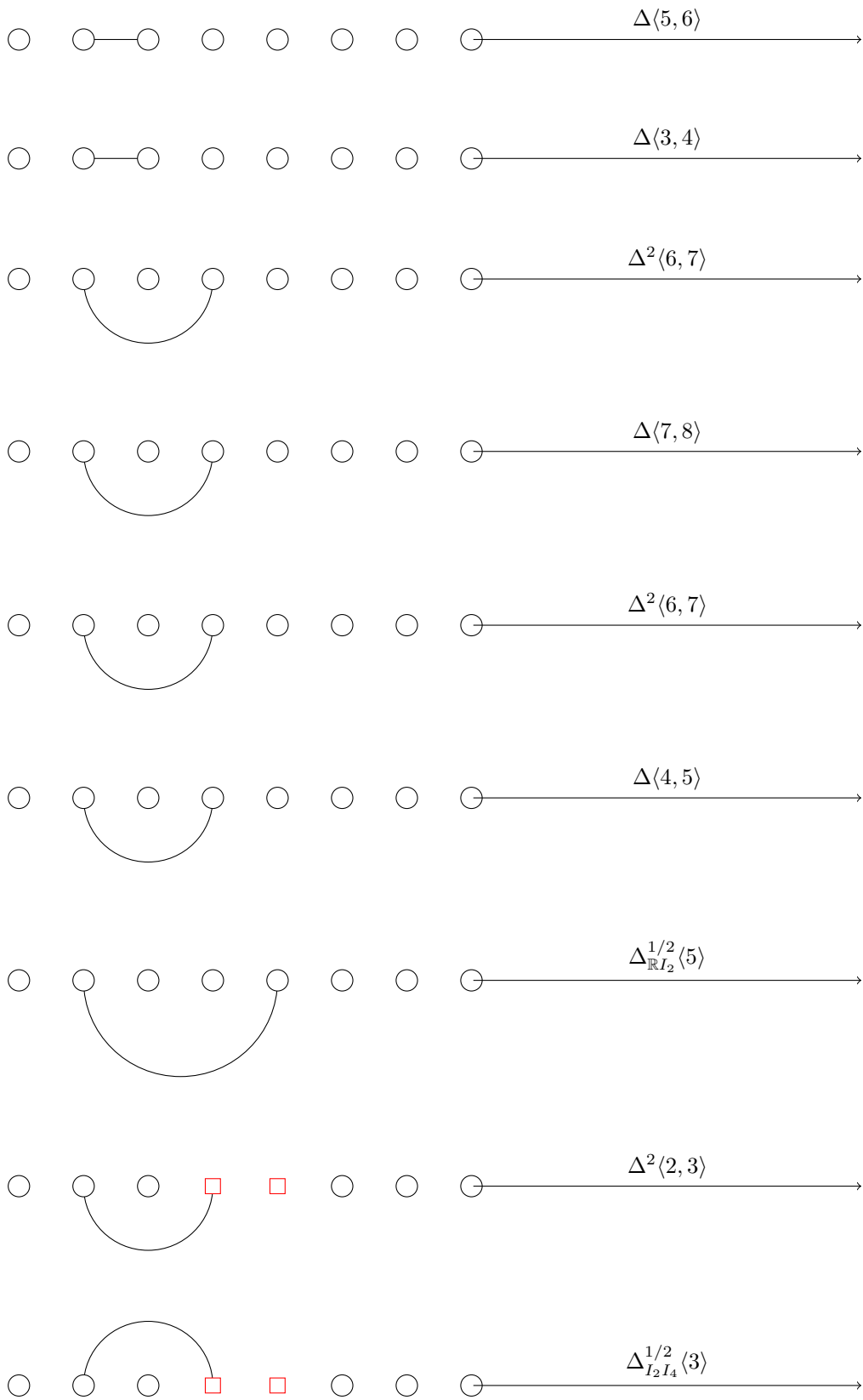
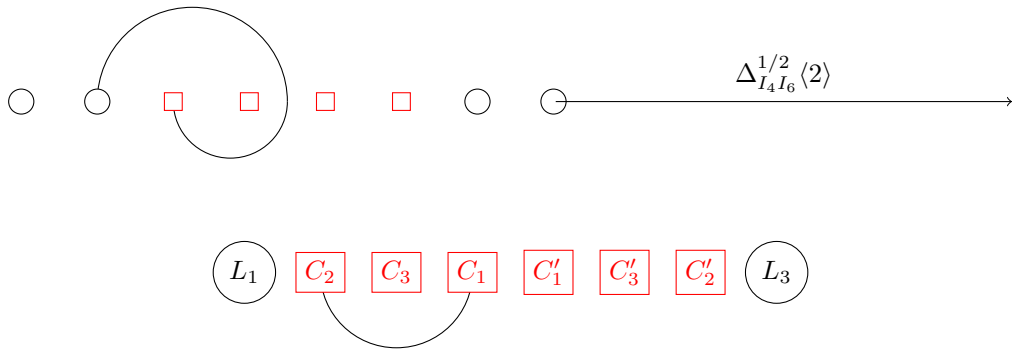
\begin{figure}[H]
\begin{center}
%
\end{center}
\end{figure}

Relation:

\begin{equation*}\begin{split}
\left\{\Gamma_{C_1}  , \Gamma_{C_2} \right\}=e
\end{split}
\end{equation*}
\subsubsection{Vertex Number 13}

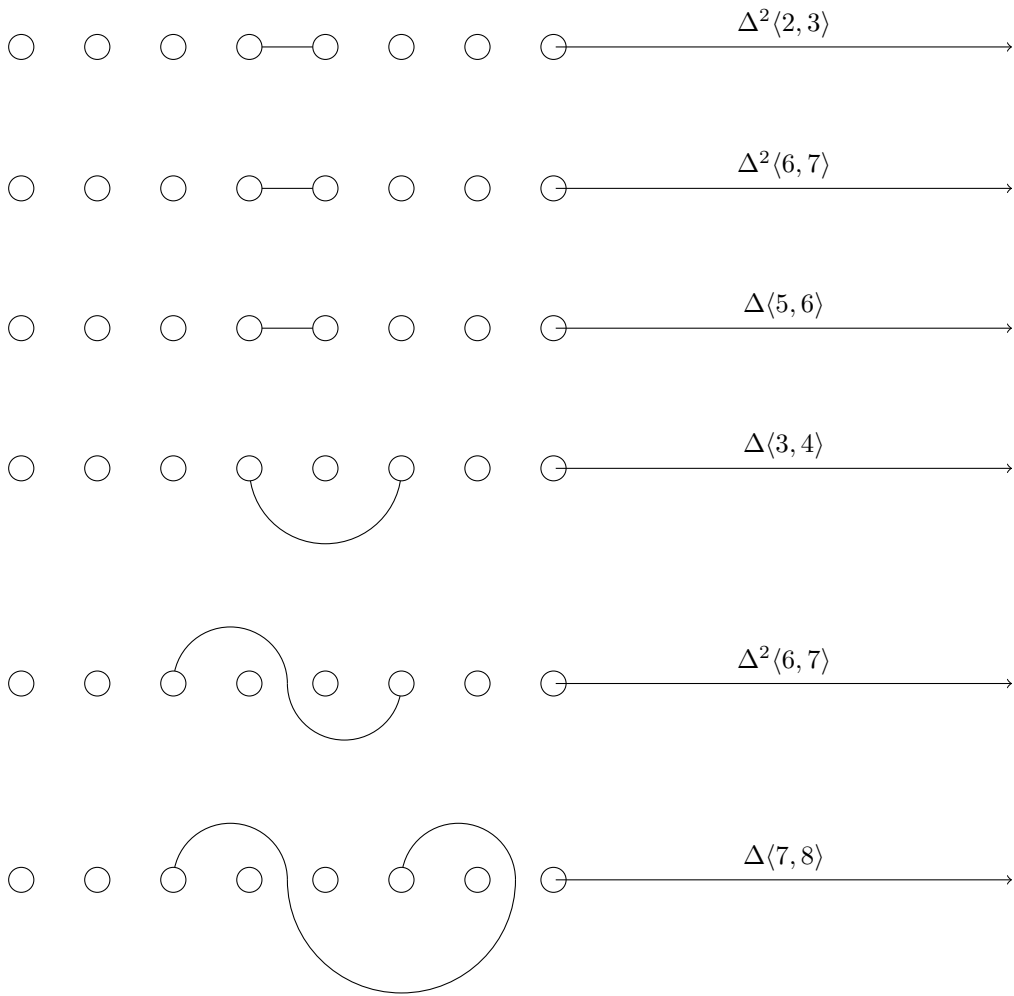
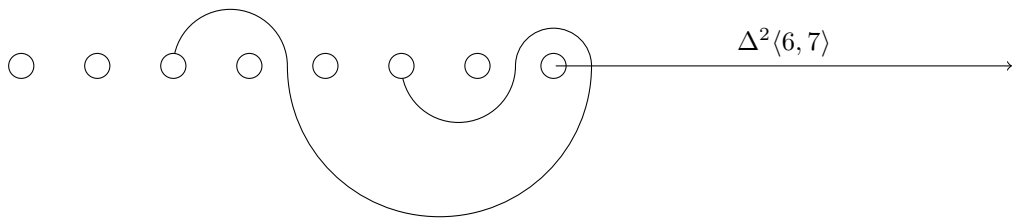
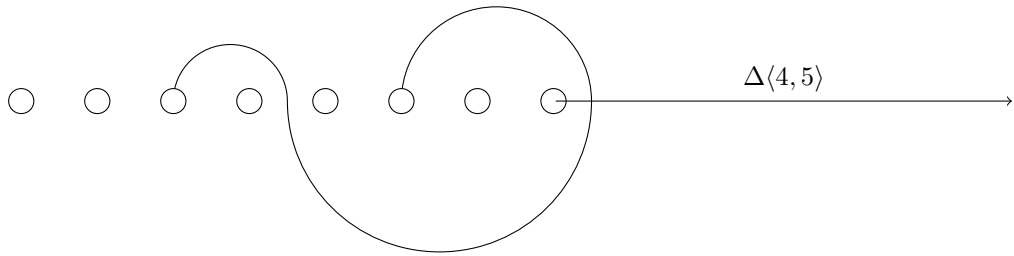
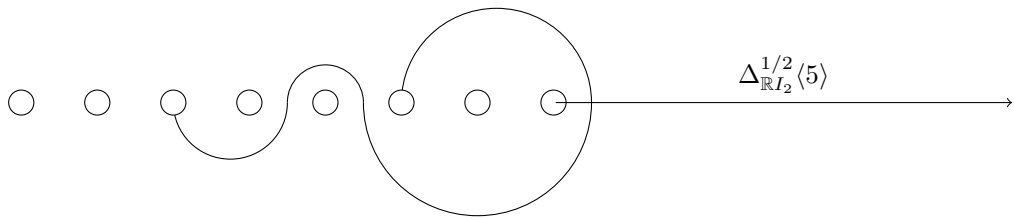
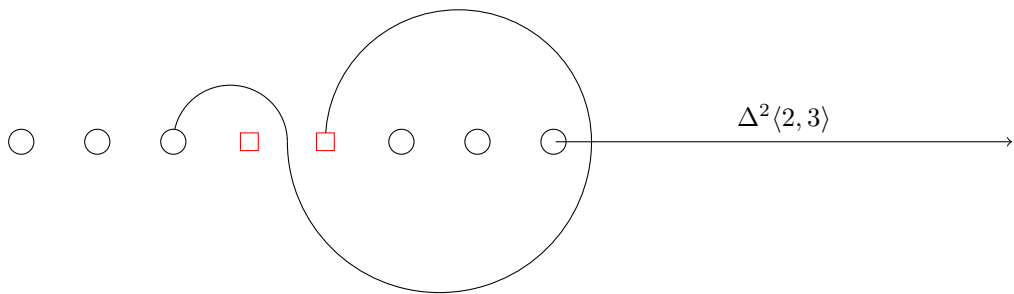
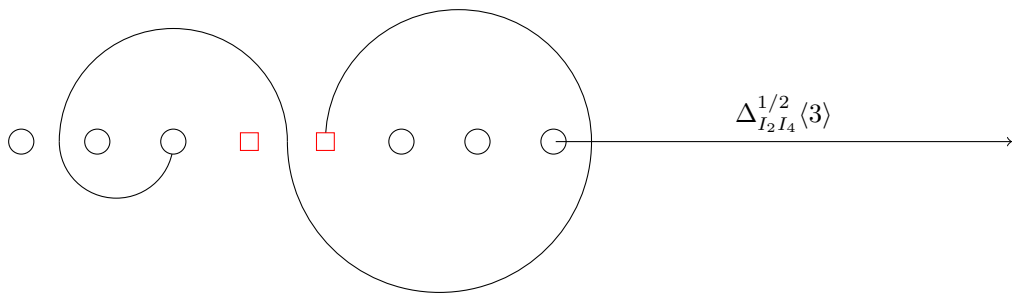
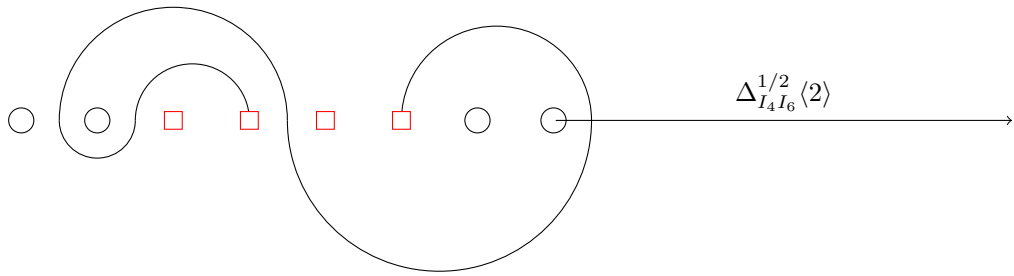
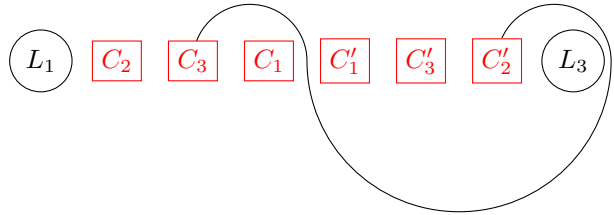
\begin{figure}[H]
\begin{center}
%
\end{center}
\end{figure}

Relation:

\begin{equation*}\begin{split}
\left[\Gamma_{L_3}^{-1}\Gamma_{C_1}\Gamma_{C_3}\Gamma_{C_1}^{-1}\Gamma_{L_3} , \Gamma_{C_2'} \right]=e
\end{split}
\end{equation*}
\subsubsection{Vertex Number 14}

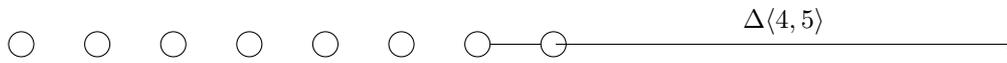
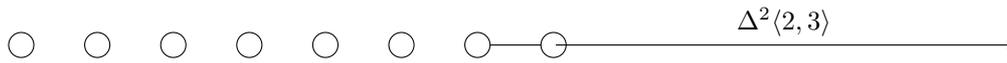
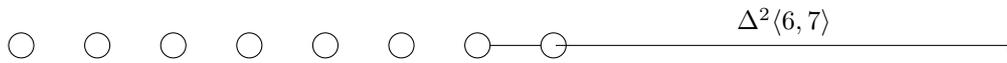
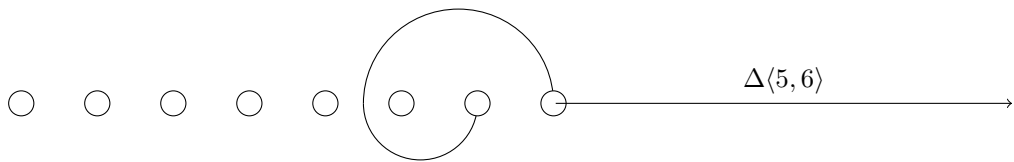
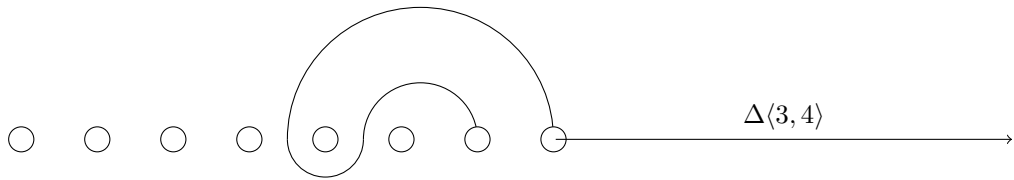
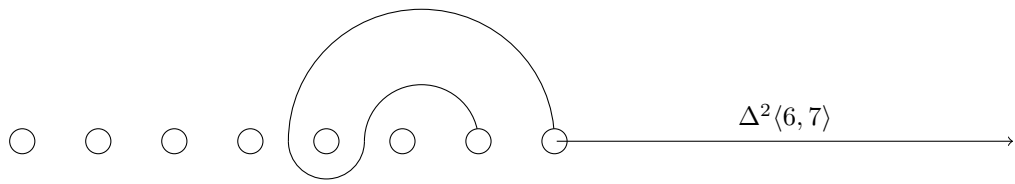
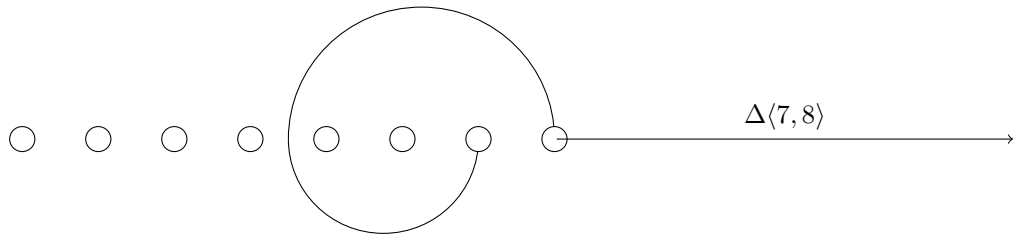
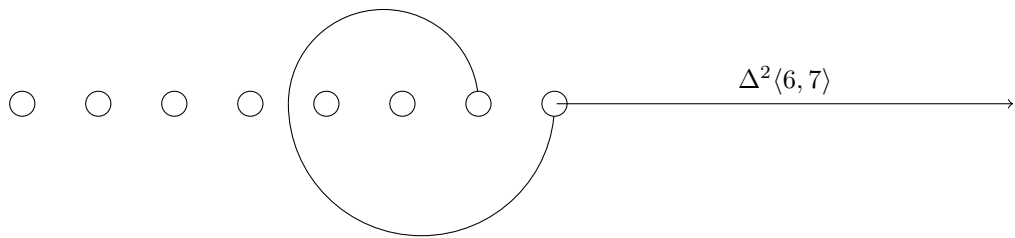
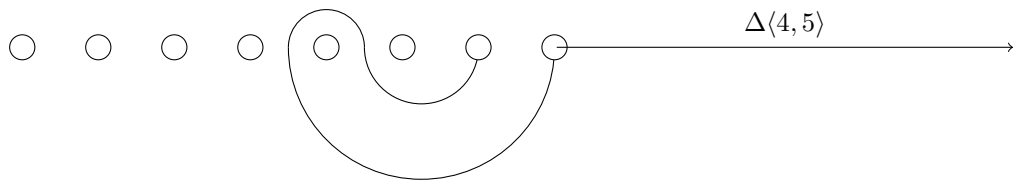
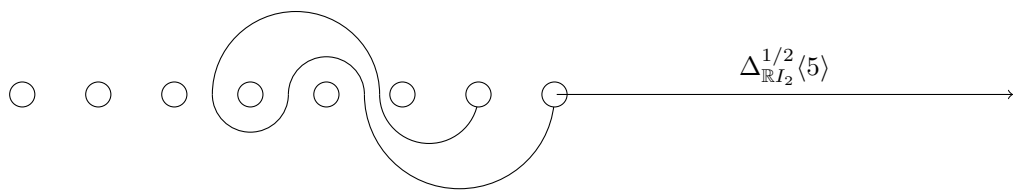
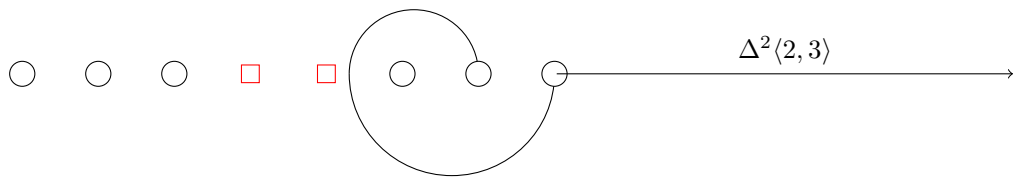
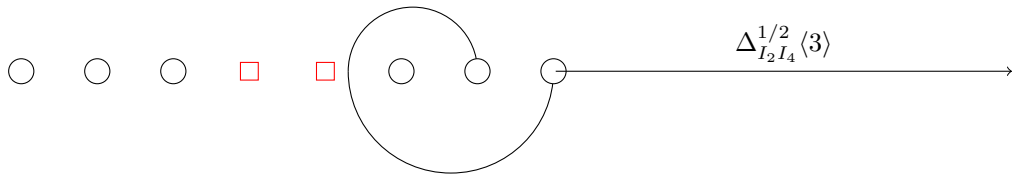
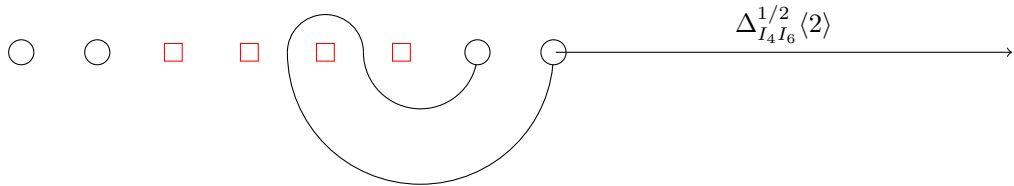
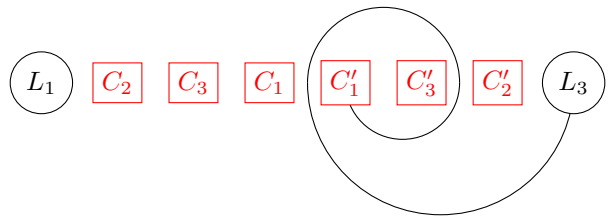
\begin{figure}[H]
\begin{center}
%
\end{center}
\end{figure}

Relation:

\begin{equation*}\begin{split}
\left[\Gamma_{C_3'}\Gamma_{C_1'}\Gamma_{L_3}\Gamma_{C_1'}^{-1}\Gamma_{C_3'}^{-1} , \Gamma_{C_1'} \right]=e
\end{split}
\end{equation*}
\subsubsection{Vertex Number 15}

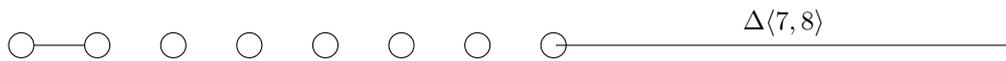
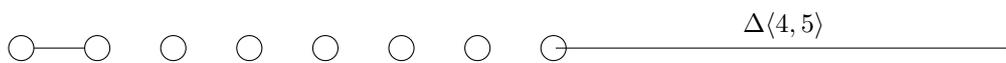
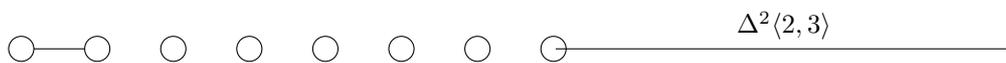
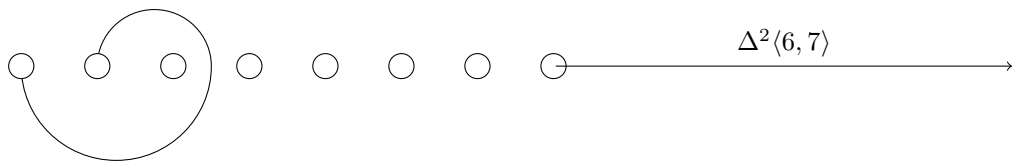
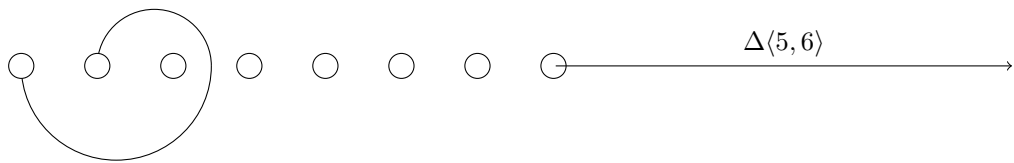
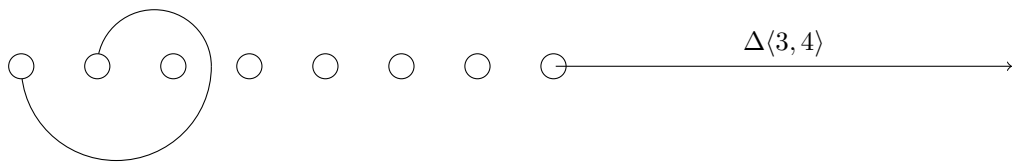
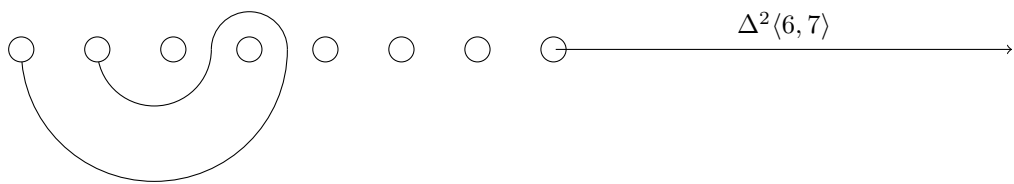
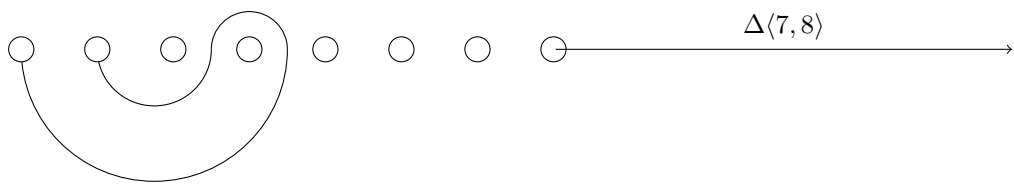
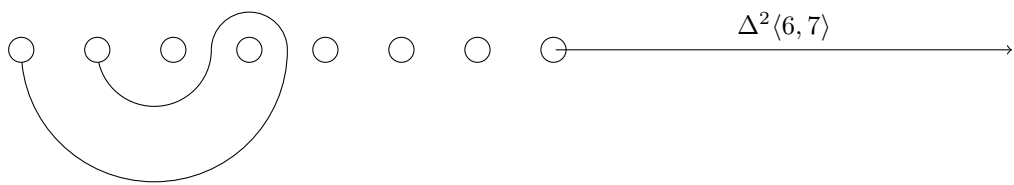
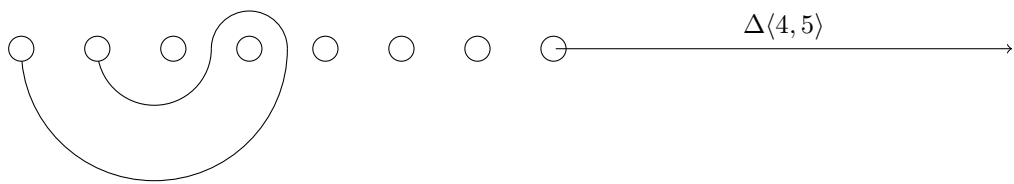
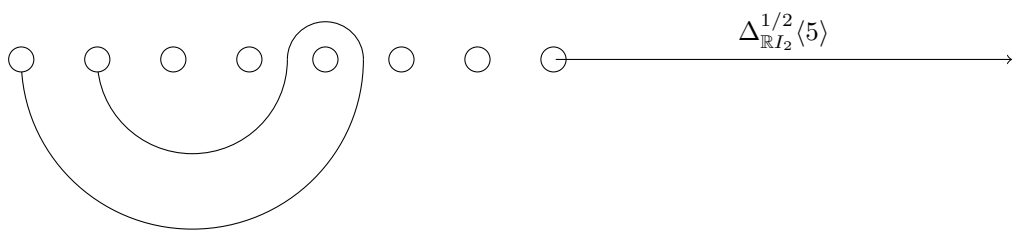
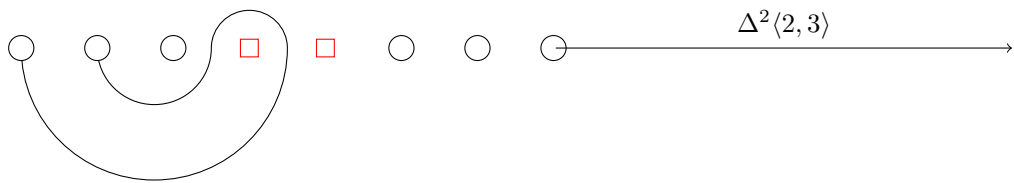
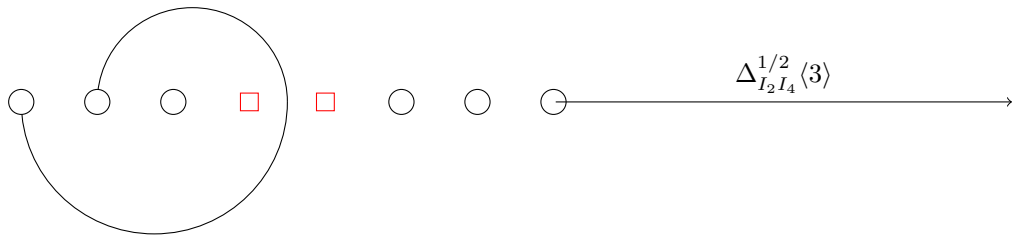
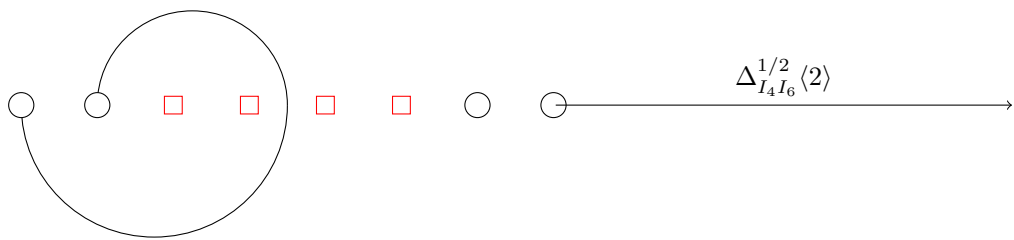
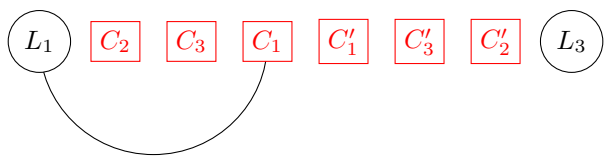
\begin{figure}[H]
\begin{center}
%
\end{center}
\end{figure}

Relation:

\begin{equation*}\begin{split}
\left[\Gamma_{L_1} , \Gamma_{C_1} \right]=e
\end{split}
\end{equation*}
\subsubsection{Vertex Number 16}

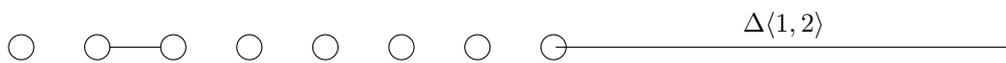
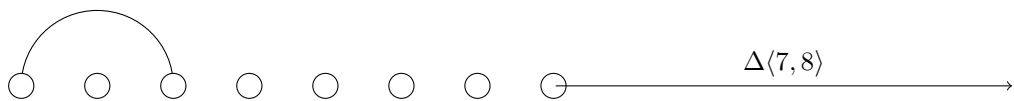
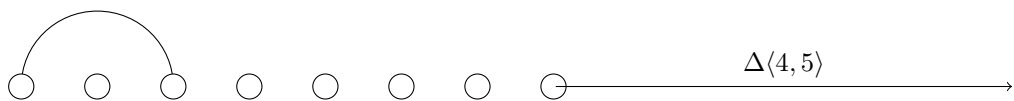
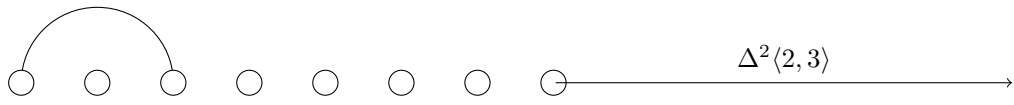
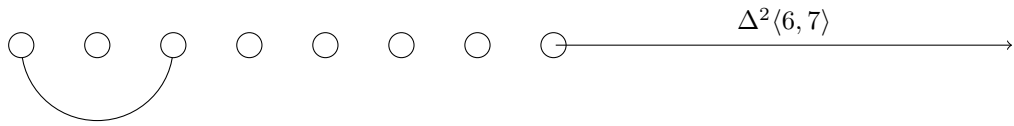
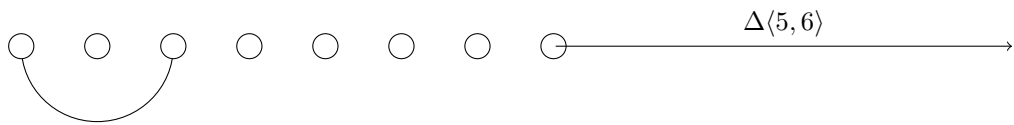
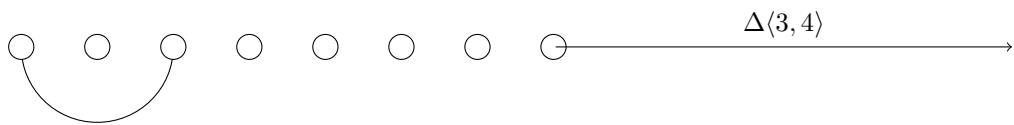
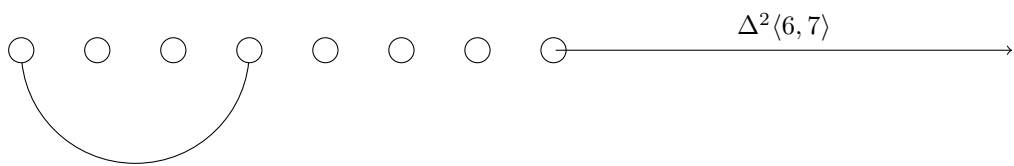
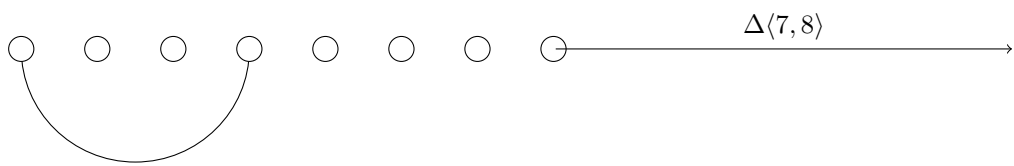
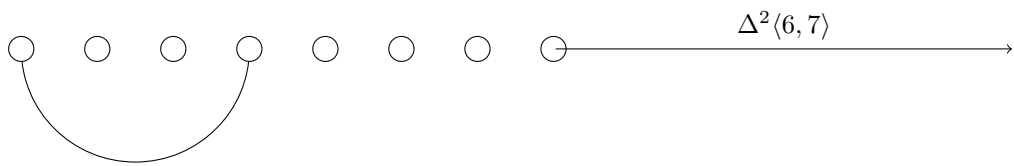
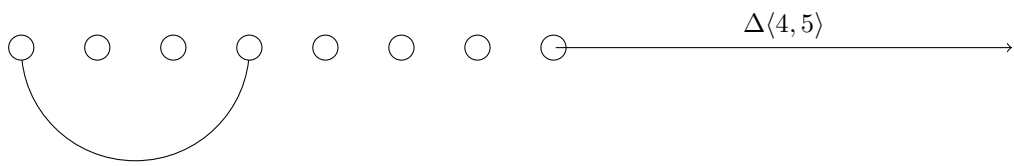
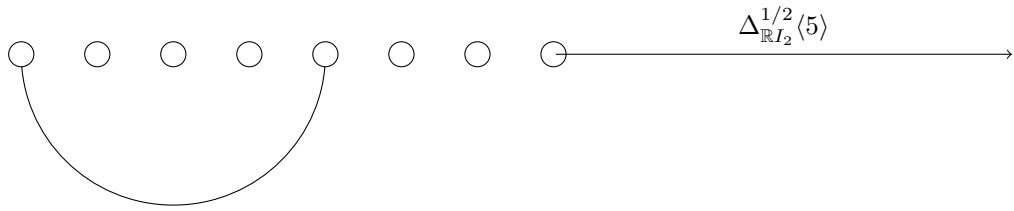
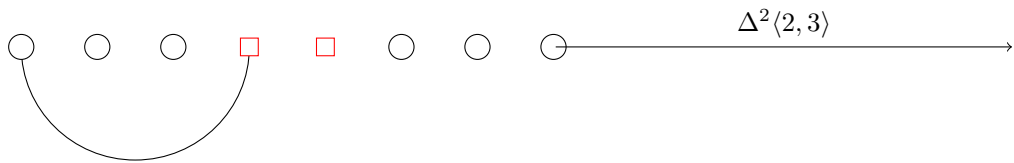
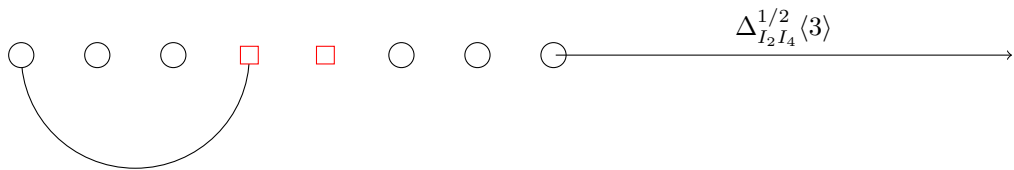
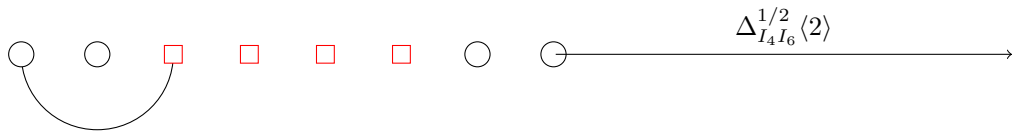
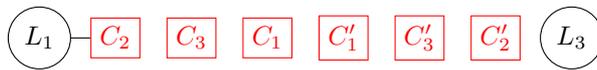
\begin{figure}[H]
\begin{center}
%
\end{center}
\end{figure}

Relation:

\begin{equation*}\begin{split}
\left\{\Gamma_{L_1}  , \Gamma_{C_2} \right\}=e
\end{split}
\end{equation*}
\subsubsection{Vertex Number 17}

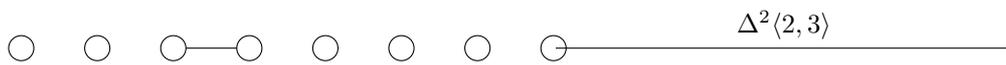
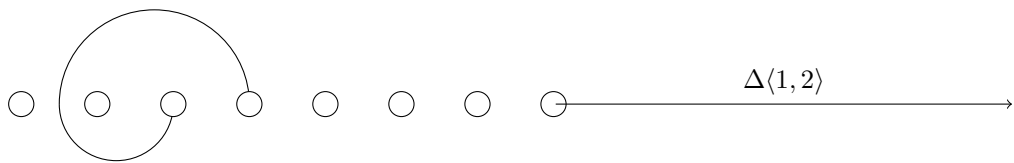
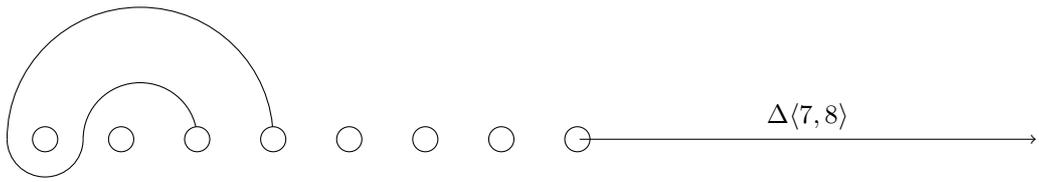
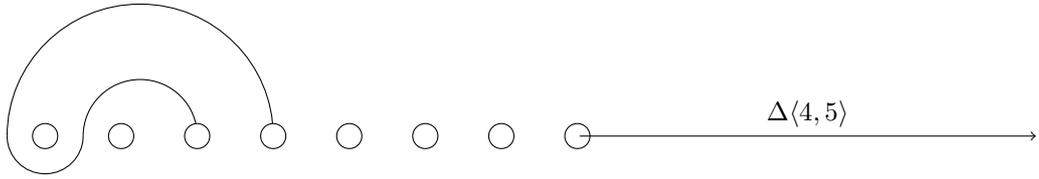
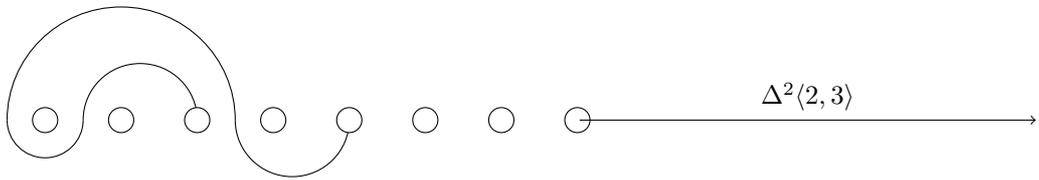
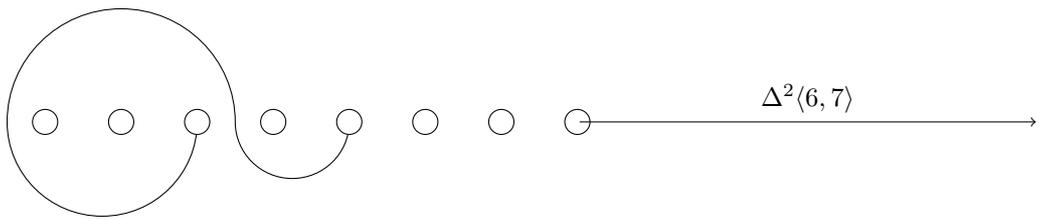
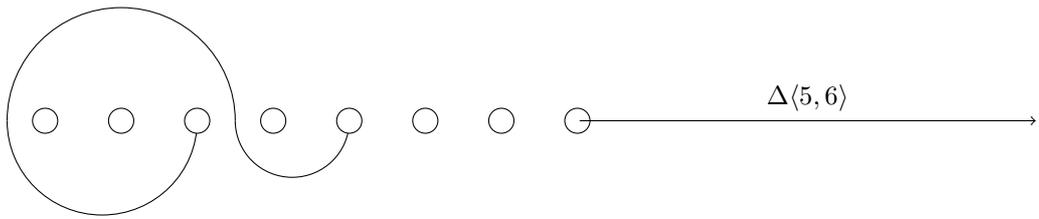
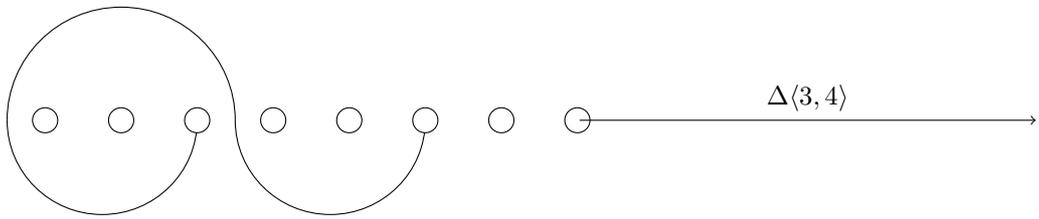
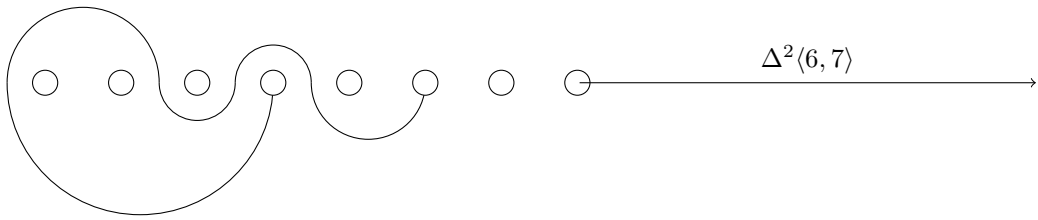
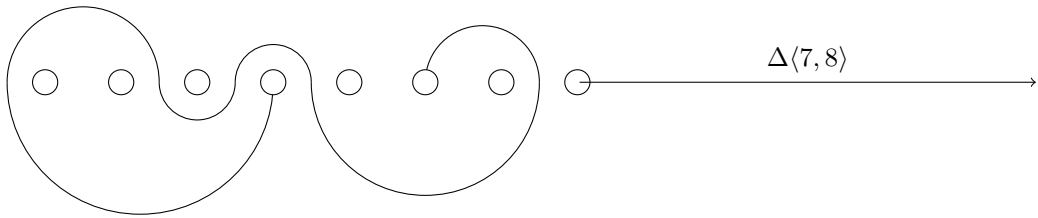
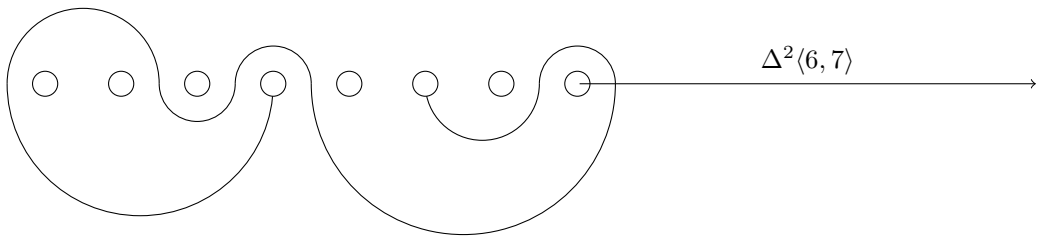
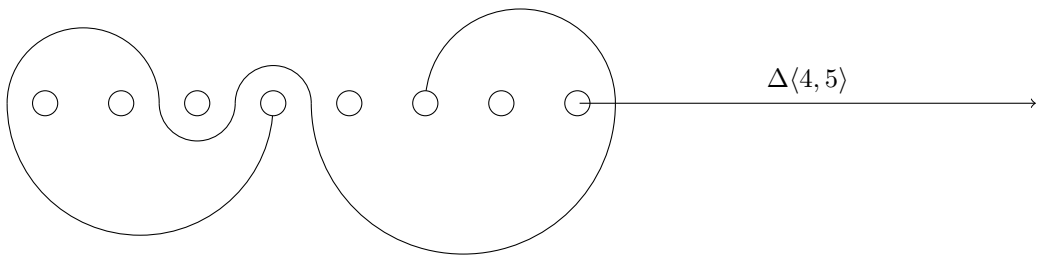
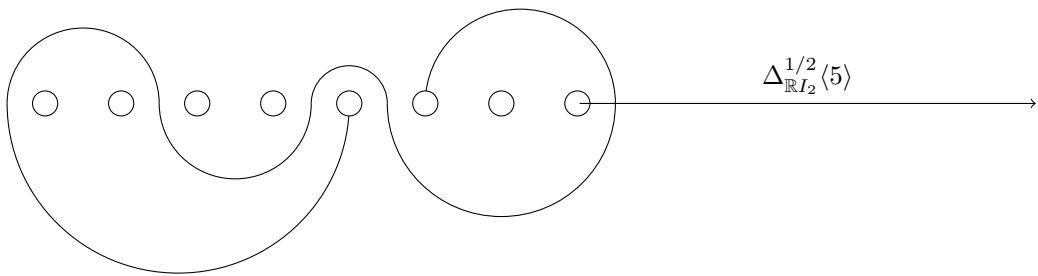
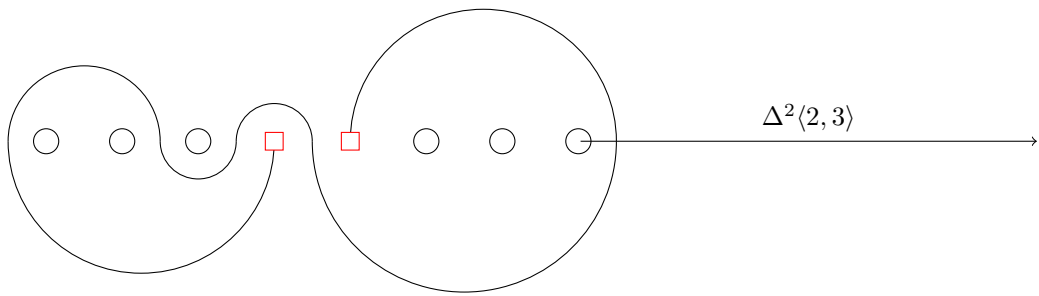
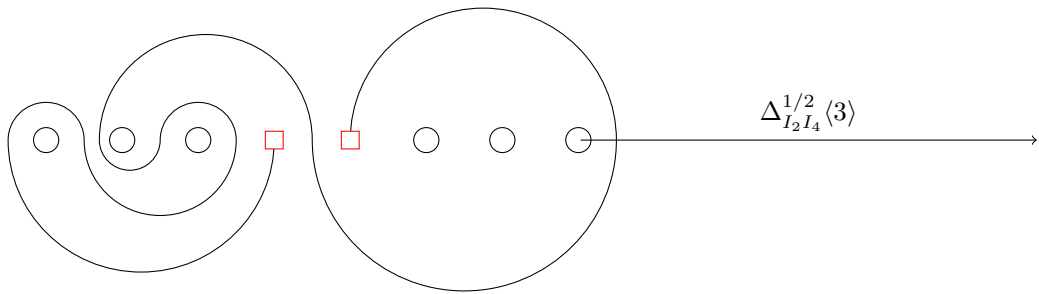
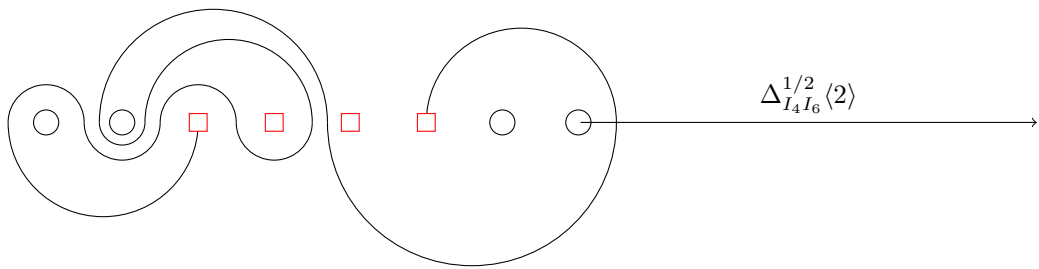
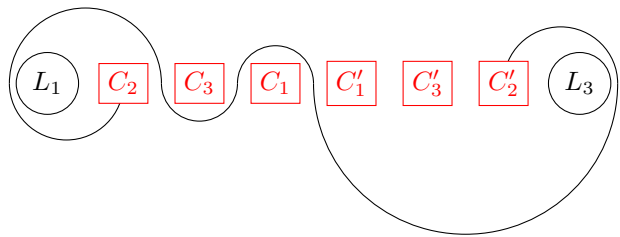
\begin{figure}[H]
\begin{center}
%
\end{center}
\end{figure}

Relation:

\begin{equation*}\begin{split}
\Gamma_{L_3}^{-1}\Gamma_{C_1}\Gamma_{C_2}\Gamma_{L_1}\Gamma_{C_2}\Gamma_{L_1}^{-1}\Gamma_{C_2}^{-1}\Gamma_{C_1}^{-1}\Gamma_{L_3}  = \Gamma_{C_2'}
\end{split}
\end{equation*}
\subsubsection{Vertex Number 18}

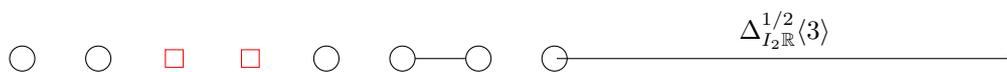
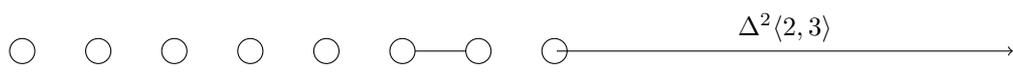
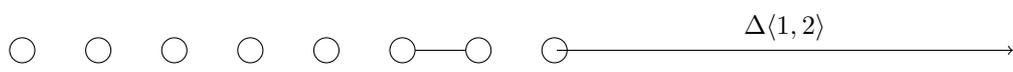
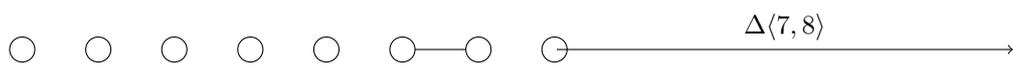
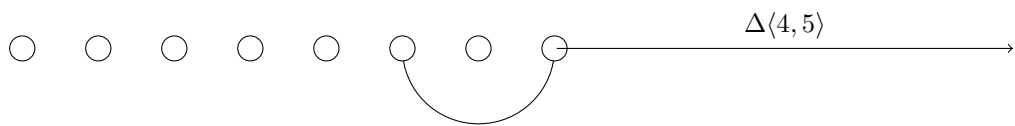
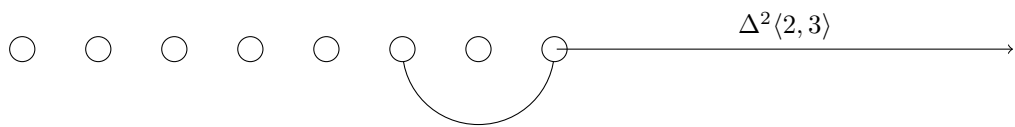
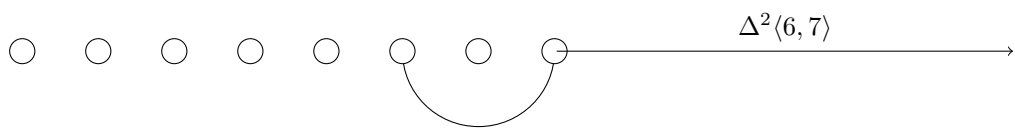
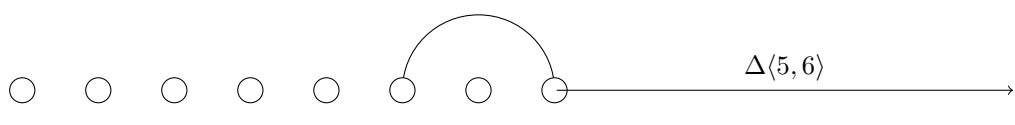
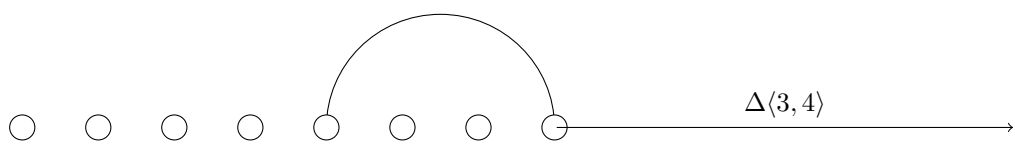
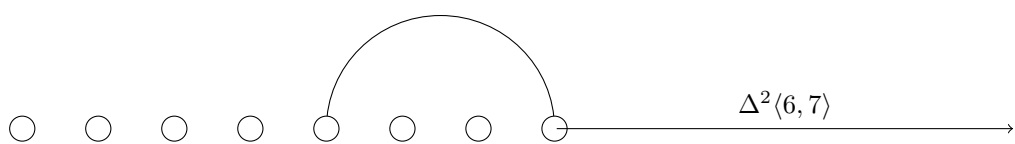
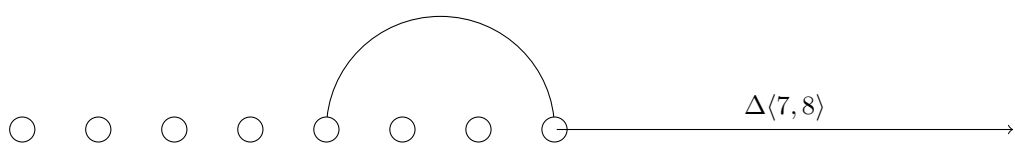
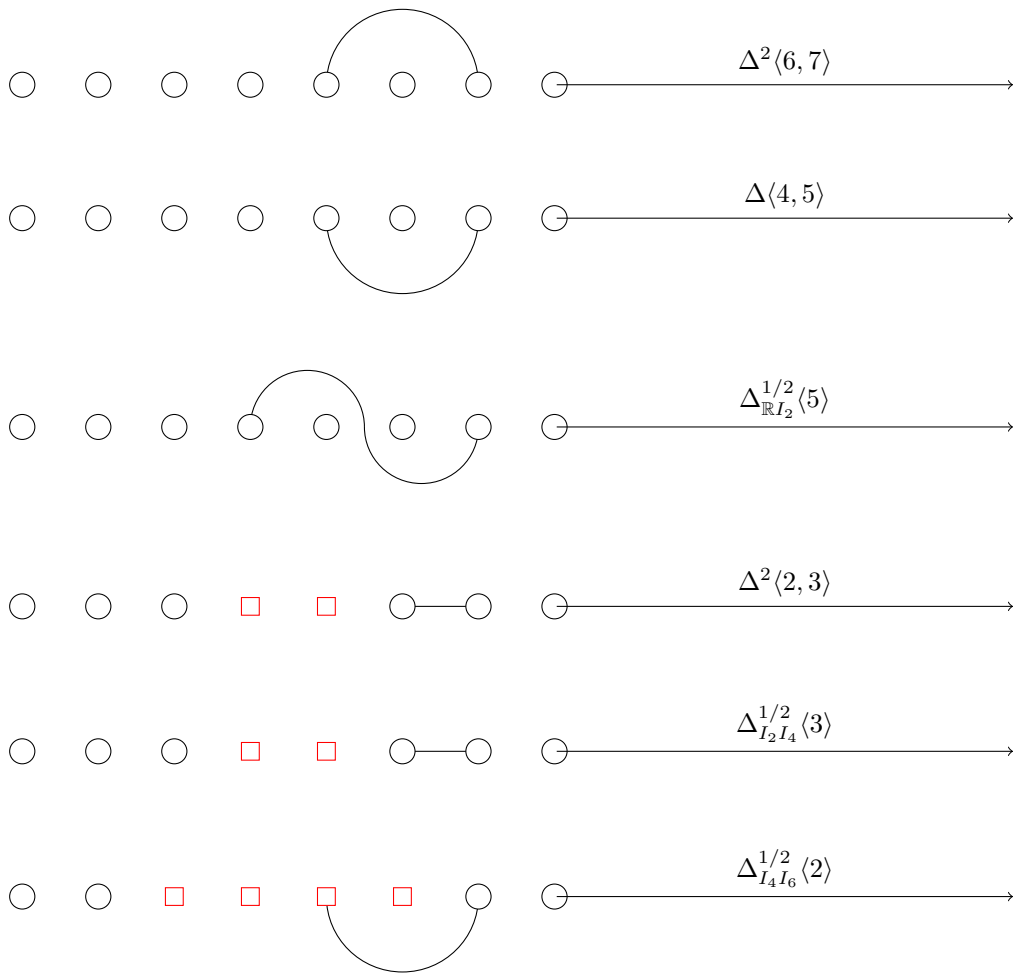
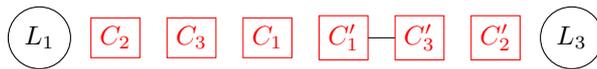
\begin{figure}[H]
\begin{center}
%
\end{center}
\end{figure}

Relation:

\begin{equation*}\begin{split}
\left\{\Gamma_{C_3'}  , \Gamma_{C_1'} \right\}=e
\end{split}
\end{equation*}
\subsubsection{Vertex Number 19}

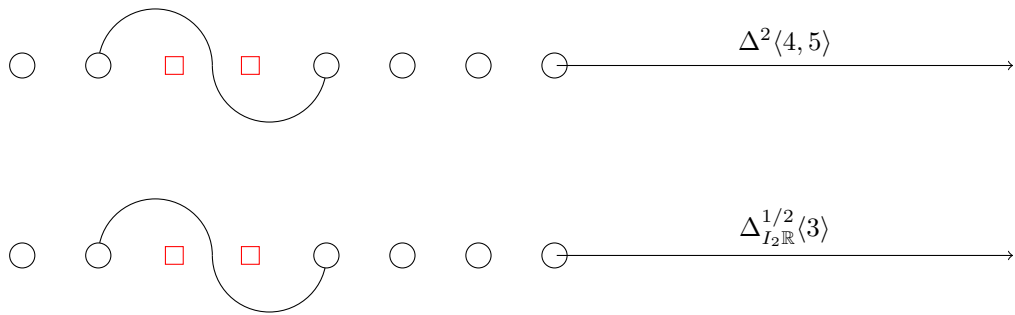
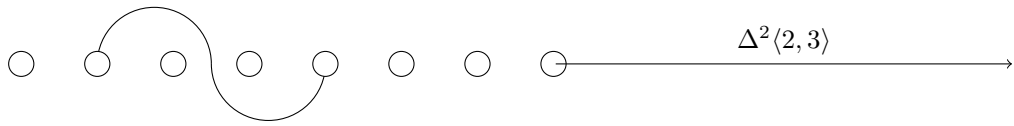
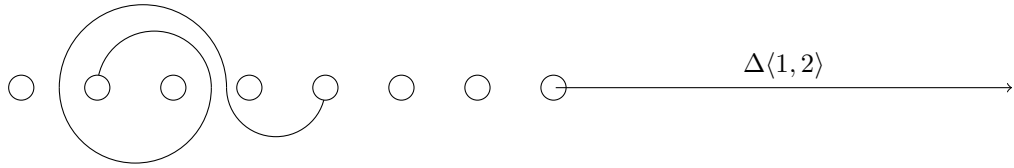
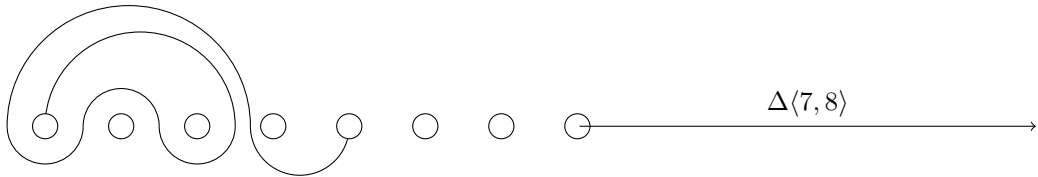
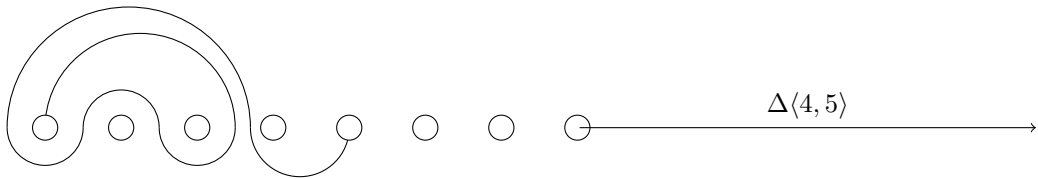
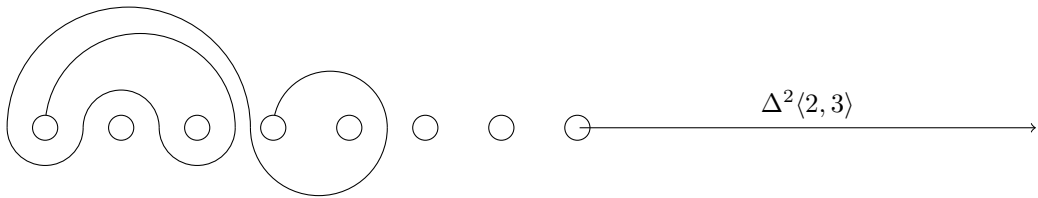
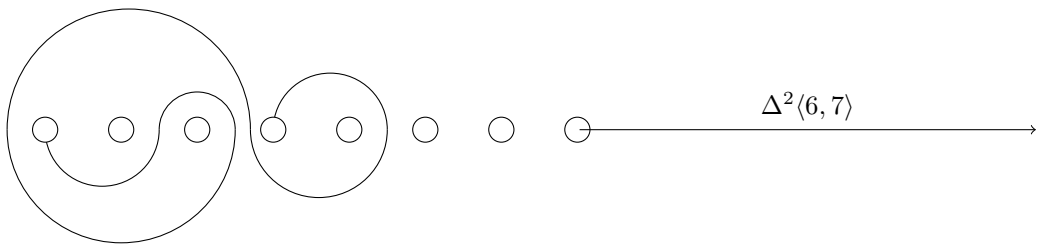
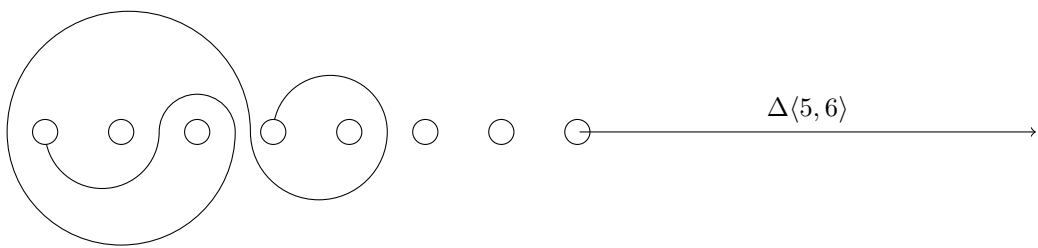
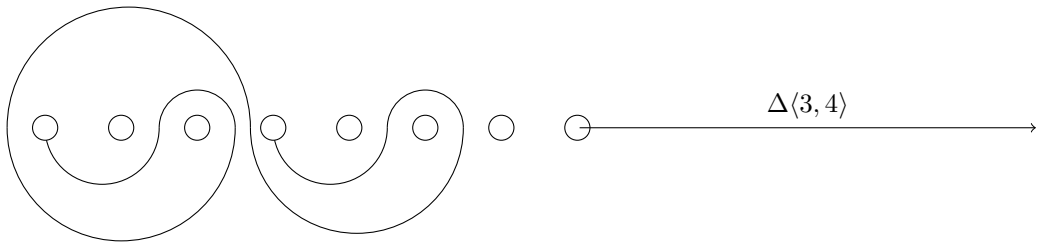
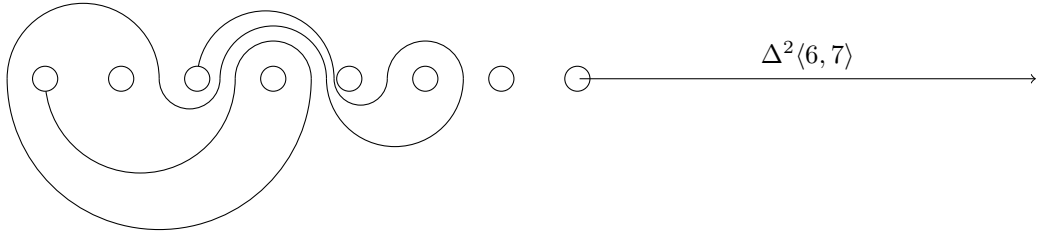
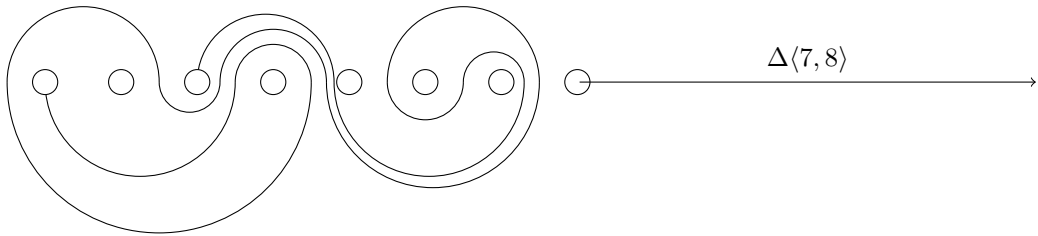
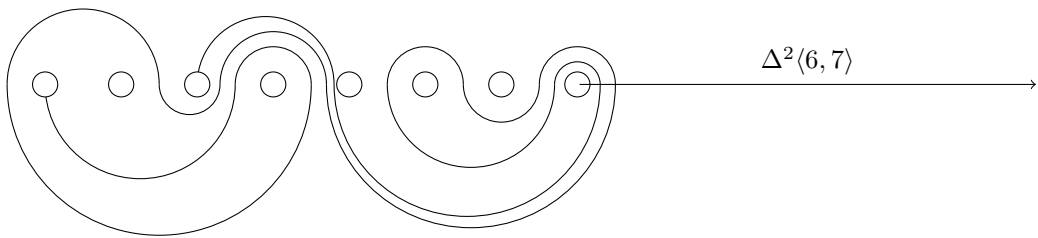
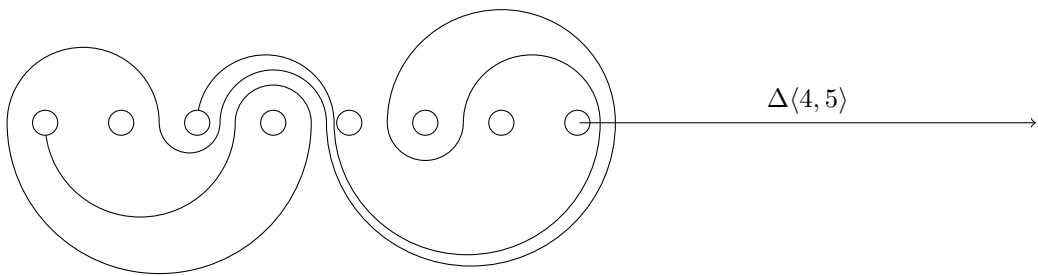
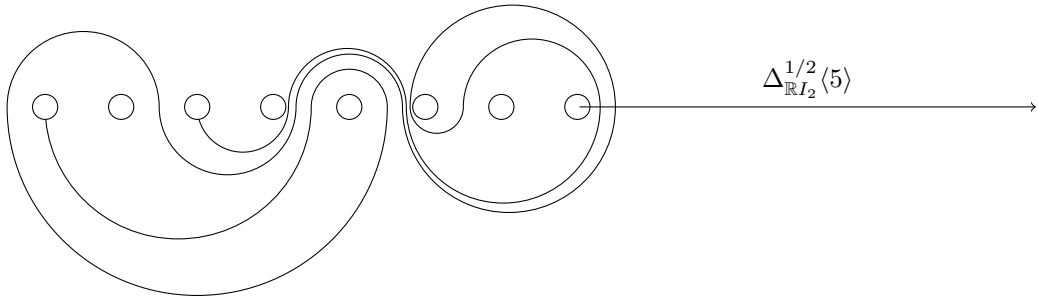
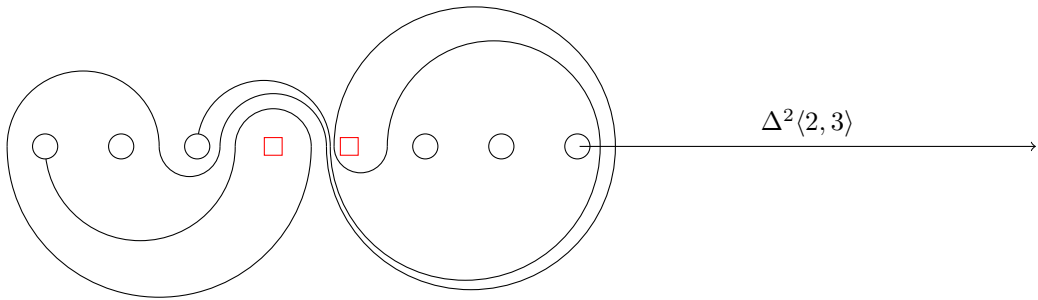
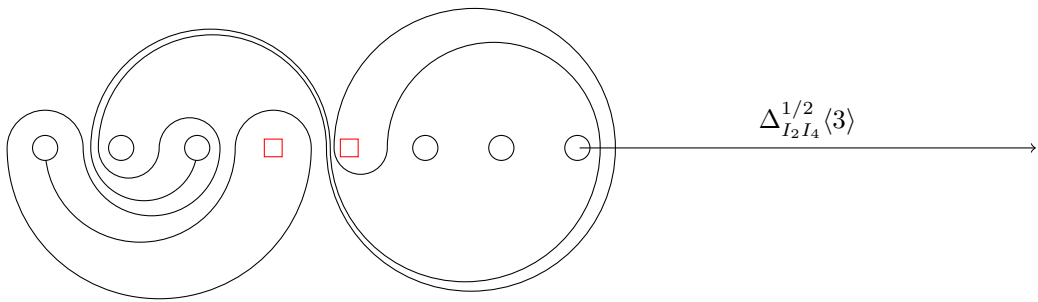
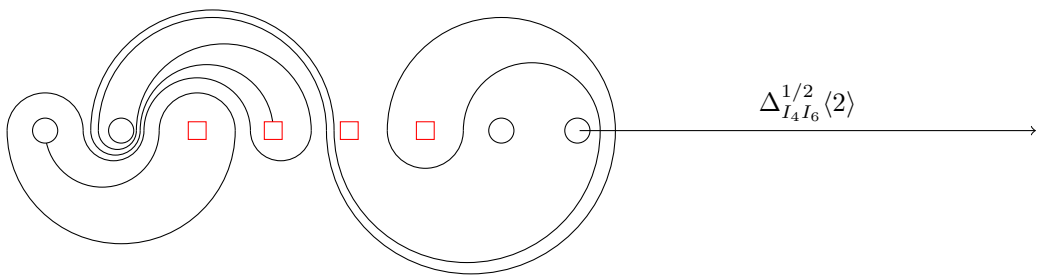
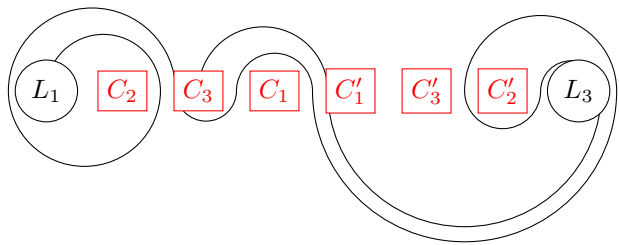
\begin{figure}[H]
\begin{center}
%
\end{center}
\end{figure}

Relation:

\begin{equation*}\begin{split}
\left\{\Gamma_{C_1}^{-1}\Gamma_{L_3}\Gamma_{C_2'}^{-1}\Gamma_{L_3}^{-1}\Gamma_{C_1}\Gamma_{C_2}\Gamma_{L_1}\Gamma_{C_2}\Gamma_{L_1}\Gamma_{C_2}^{-1}\Gamma_{L_1}^{-1}\Gamma_{C_2}^{-1}\Gamma_{C_1}^{-1}\Gamma_{L_3}\Gamma_{C_2'}\Gamma_{L_3}^{-1}\Gamma_{C_1}  , \Gamma_{C_3} \right\}=e
\end{split}
\end{equation*}
\subsubsection{Vertex Number 20}

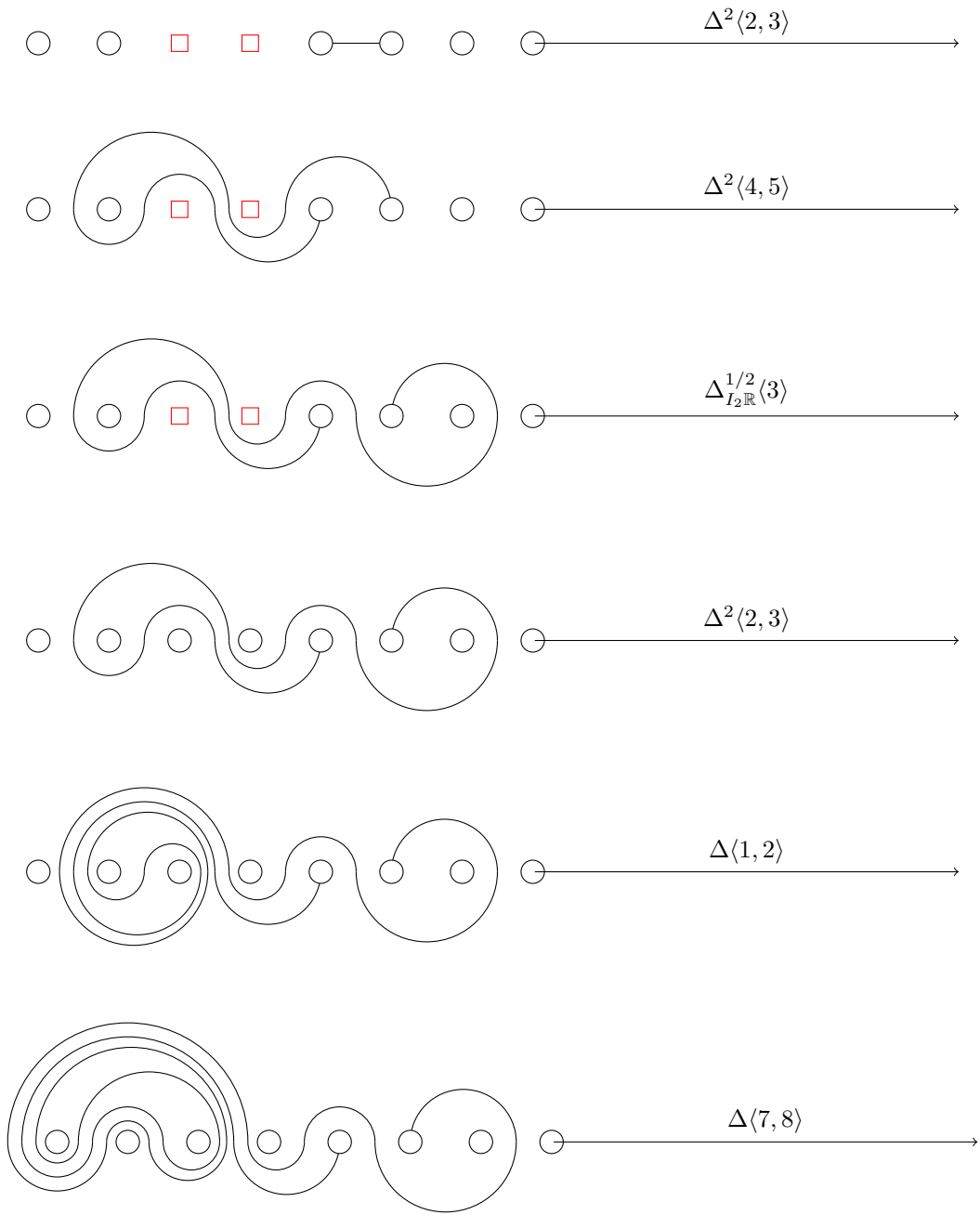
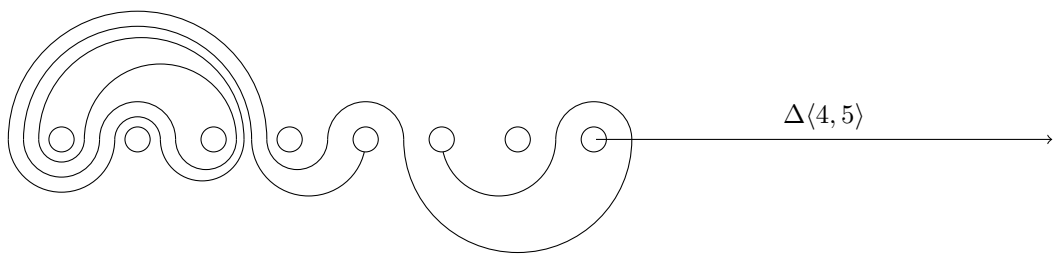
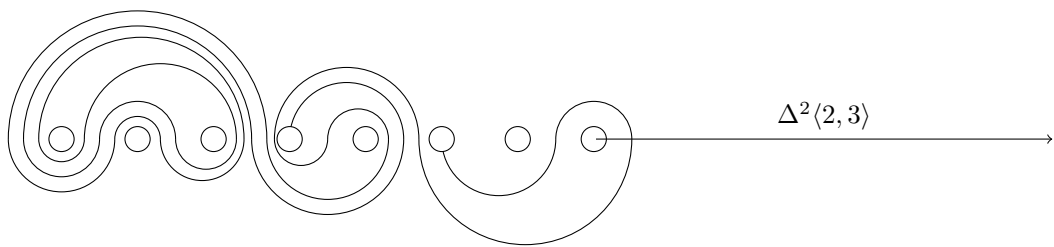
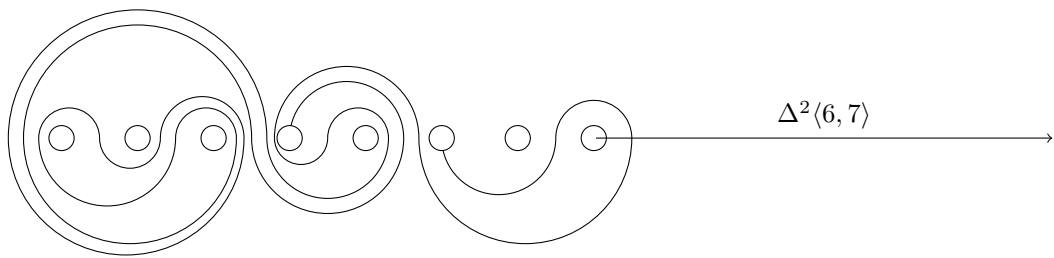
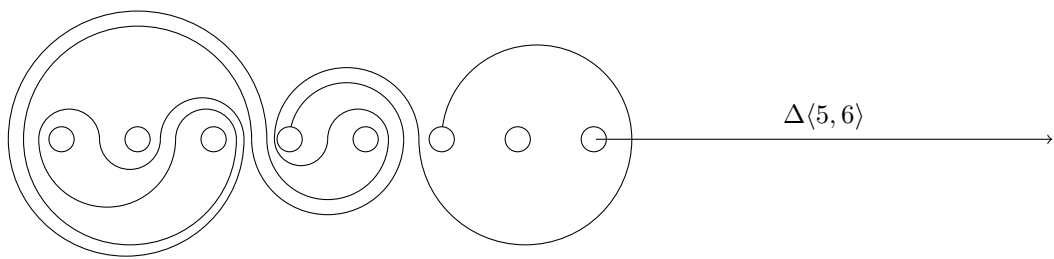
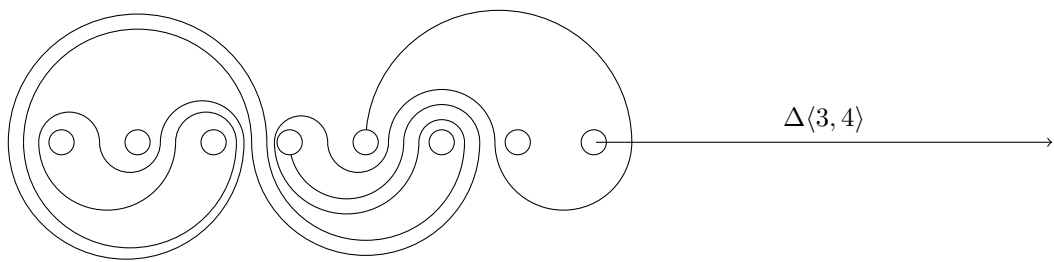
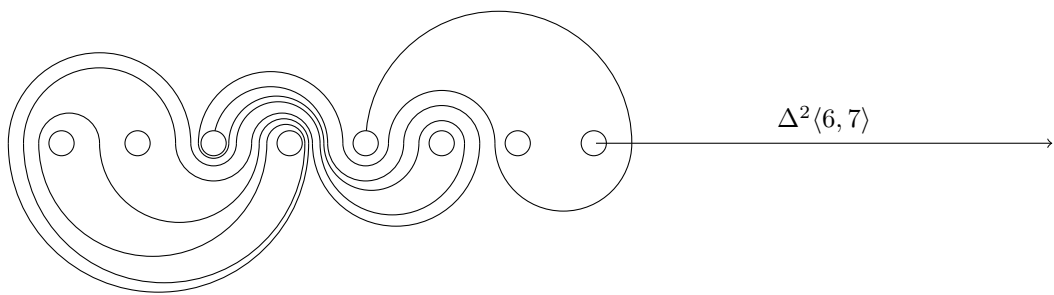
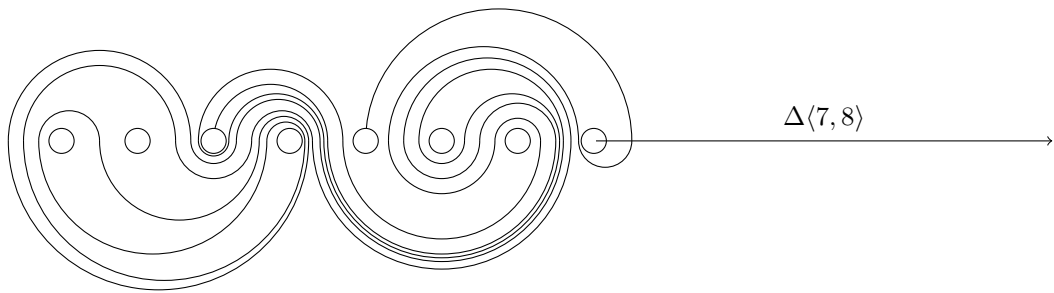
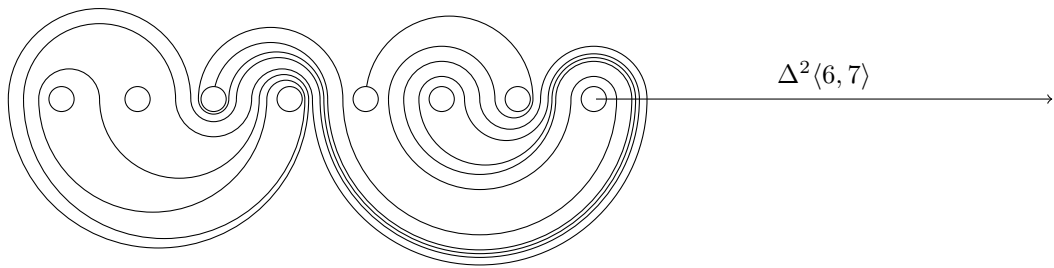
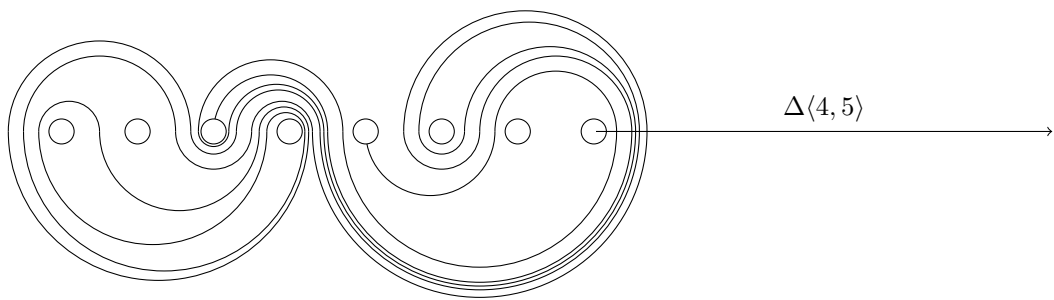
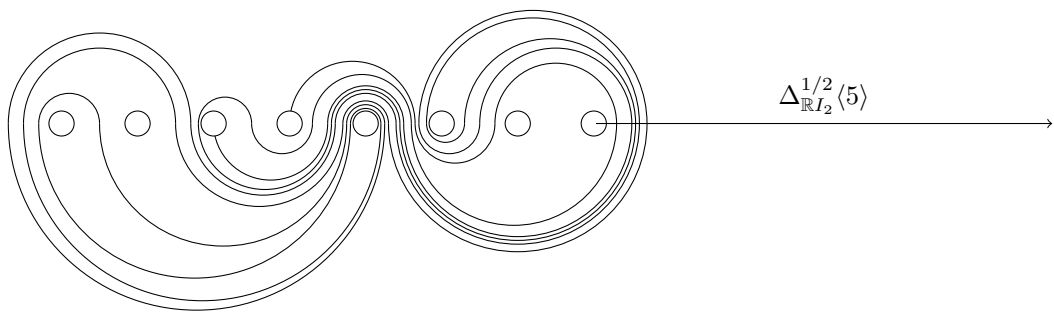
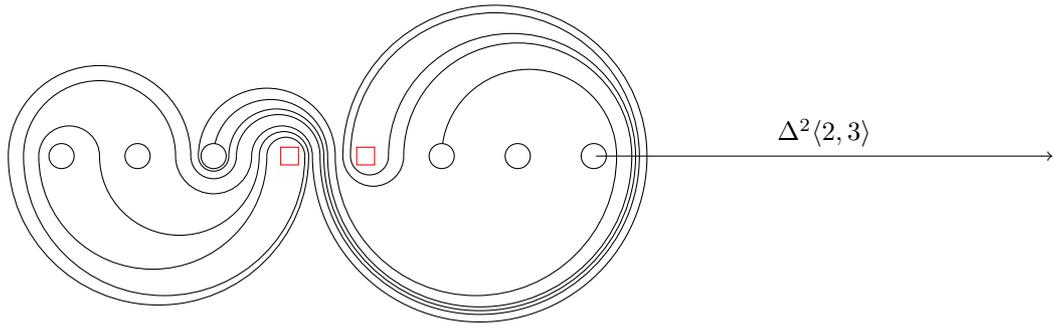
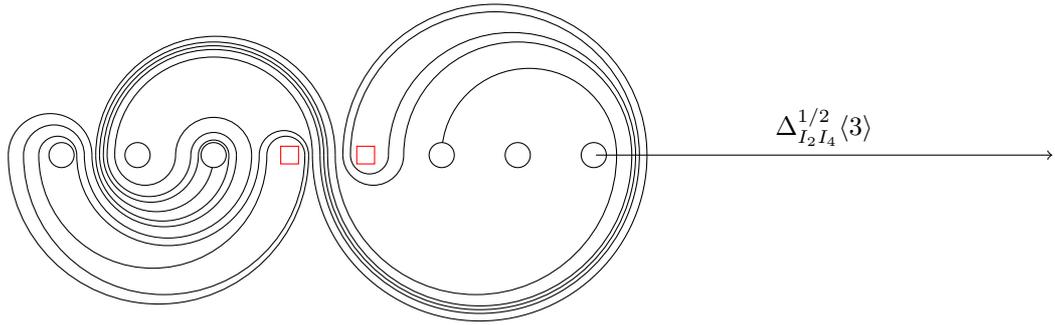
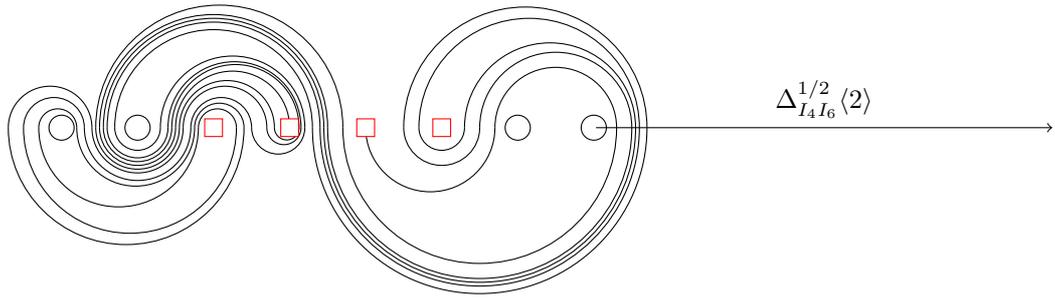
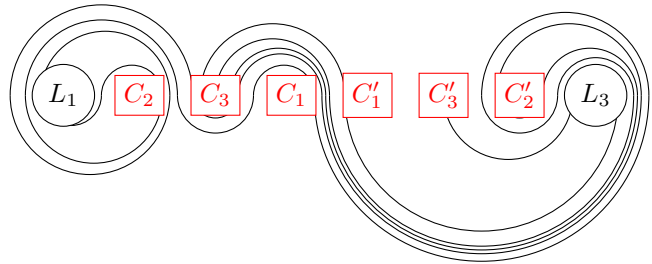
\begin{figure}[H]
\begin{center}
%
\end{center}
\end{figure}

Relation:
\begin{equation*}
\begin{split}
\Gamma_{L_3}^{-1}\Gamma_{C_1}\Gamma_{C_3}\Gamma_{C_1}^{-1}\Gamma_{L_3}\Gamma_{C_2'}^{-1}\Gamma_{L_3}^{-1}  \Gamma_{C_1}\Gamma_{C_2}\Gamma_{L_1}\Gamma_{C_2}\Gamma_{L_1}\Gamma_{C_2}^{-1}\Gamma_{L_1}^{-1}
\Gamma_{C_2}^{-1}\Gamma_{C_1}^{-1}\Gamma_{L_3}\Gamma_{C_2'}\Gamma_{L_3}^{-1}\Gamma_{C_1}\Gamma_{C_3} \times & \\
\times \Gamma_{C_1}^{-1}\Gamma_{L_3}\Gamma_{C_2'}^{-1}\Gamma_{L_3}^{-1}\Gamma_{C_1}\Gamma_{C_2}\Gamma_{L_1}
\Gamma_{C_2}\Gamma_{L_1}^{-1}\Gamma_{C_2}^{-1}\Gamma_{L_1}^{-1}\Gamma_{C_2}^{-1}\Gamma_{C_1}^{-1}
\Gamma_{L_3}\Gamma_{C_2'}\Gamma_{L_3}^{-1}\Gamma_{C_1}\Gamma_{C_3}^{-1}\Gamma_{C_1}^{-1}\Gamma_{L_3}  & = \Gamma_{C_3'} ,
\end{split}
\end{equation*}
\subsubsection{Vertex Number 21}

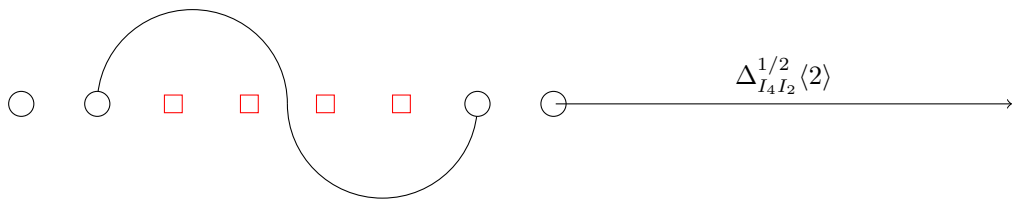
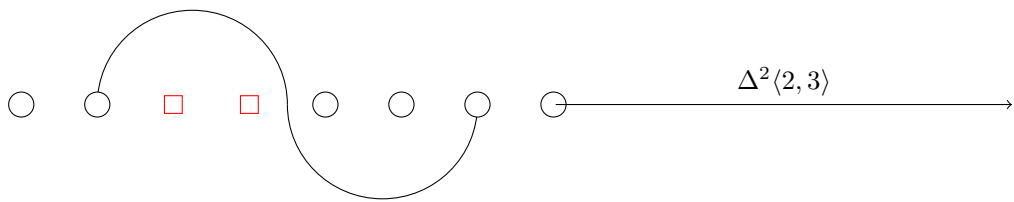
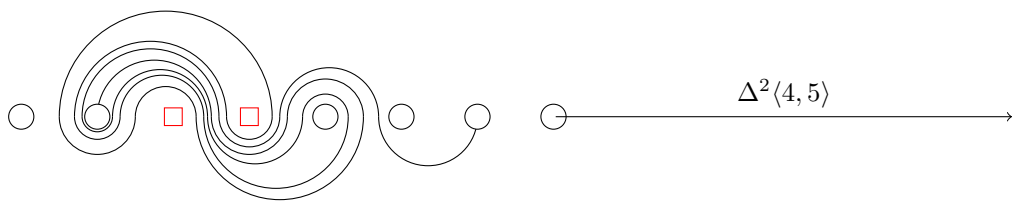
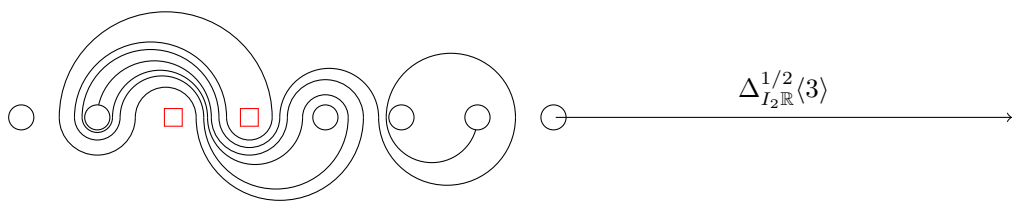
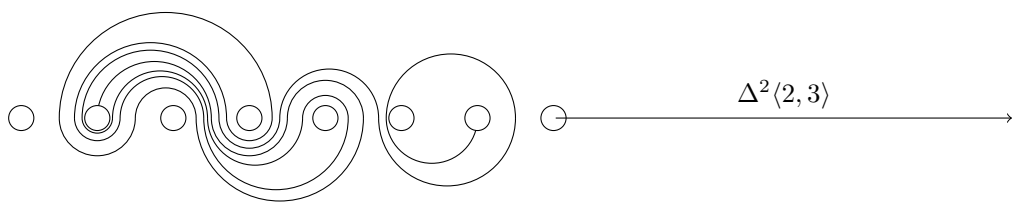
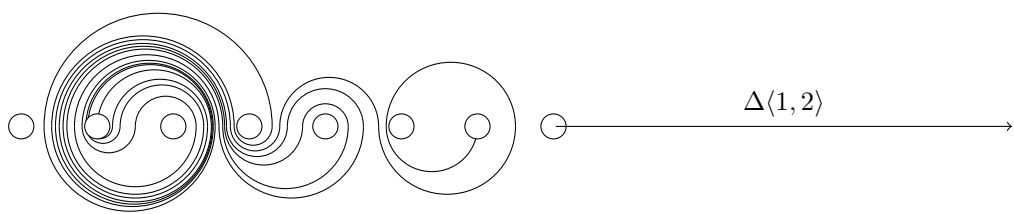
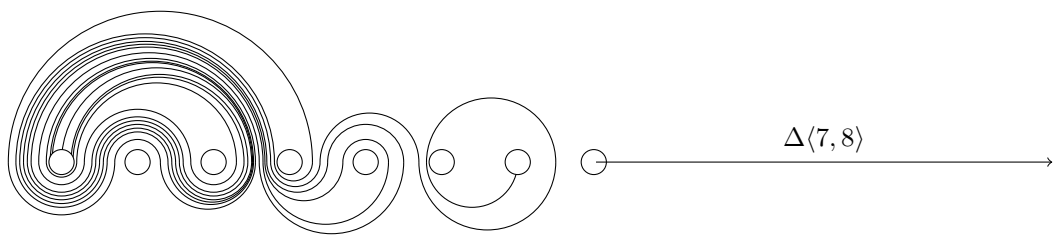
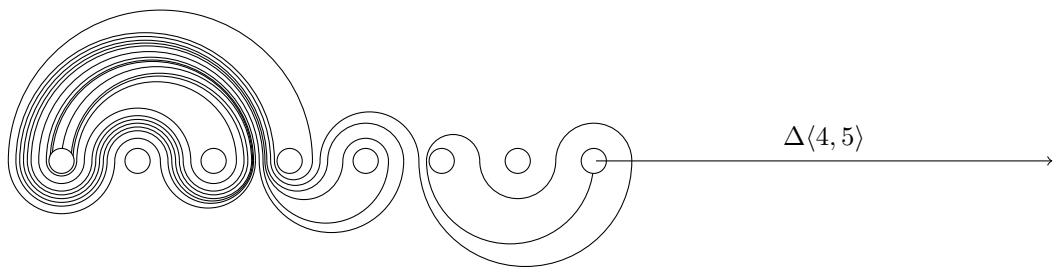
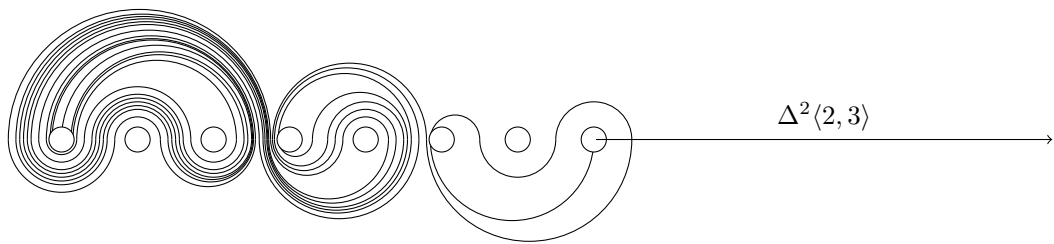
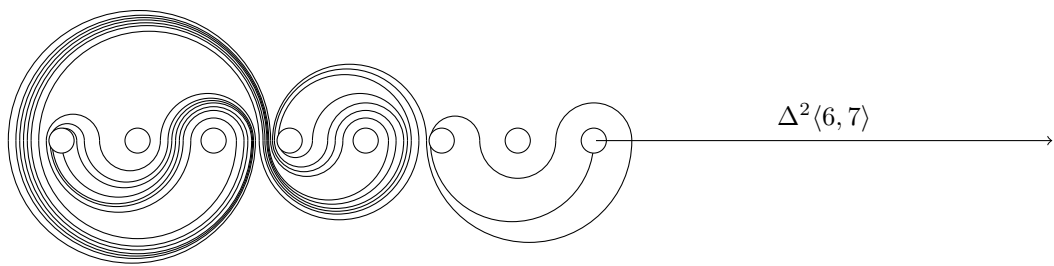
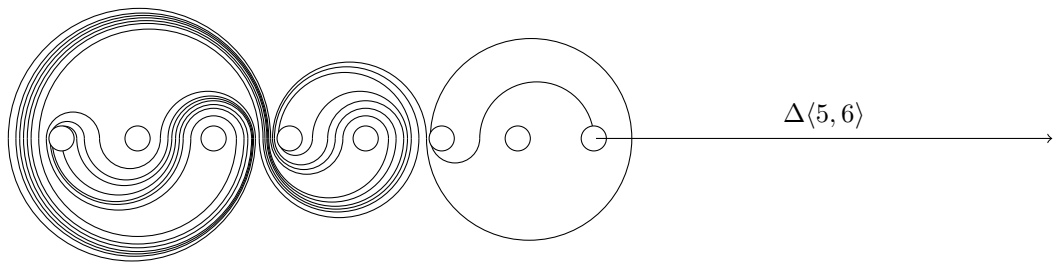
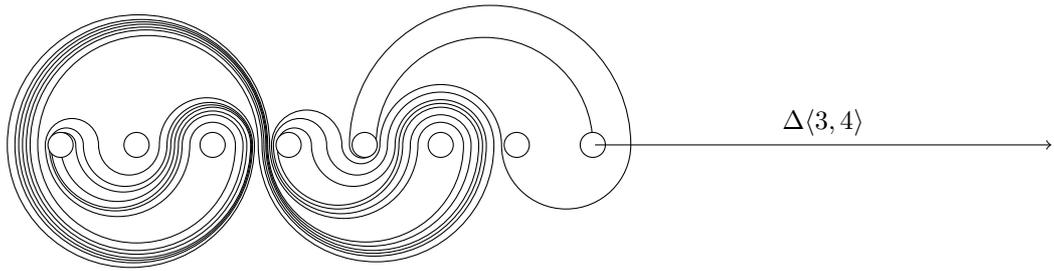
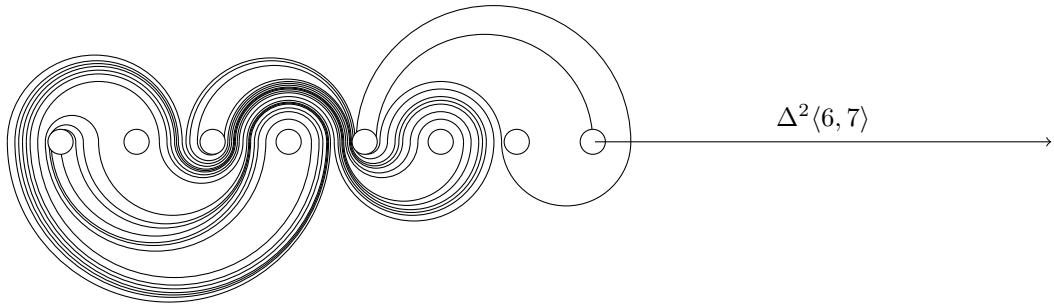
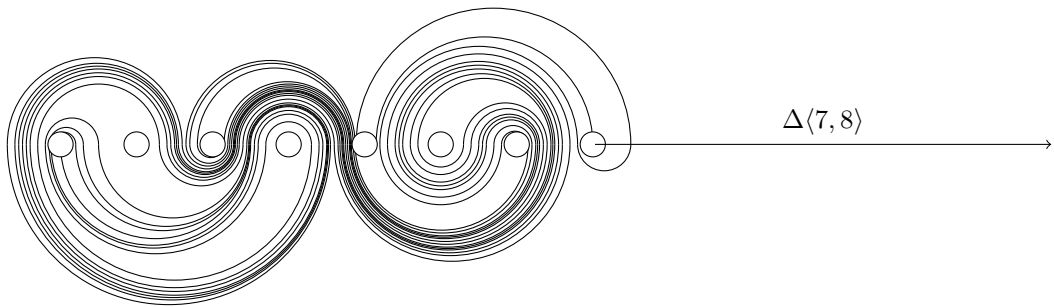
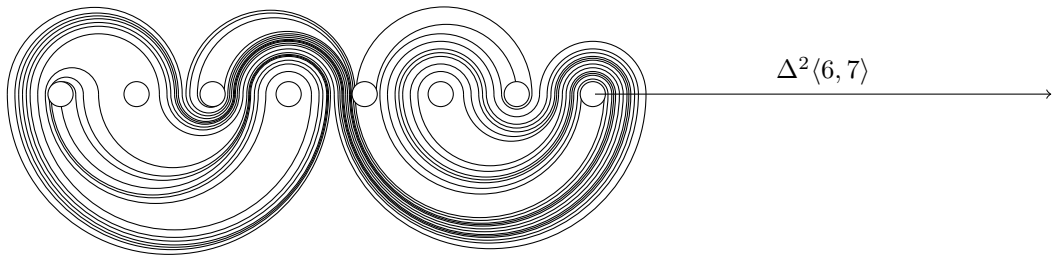
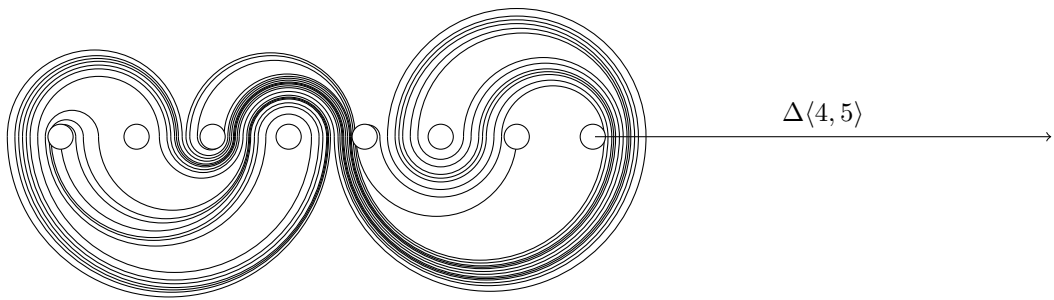
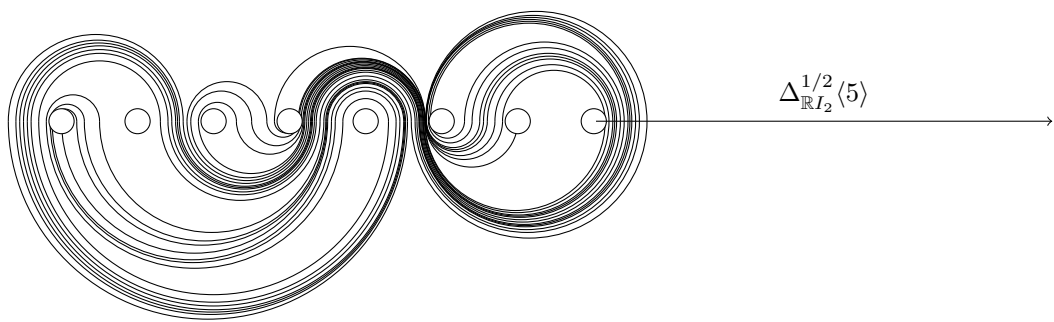
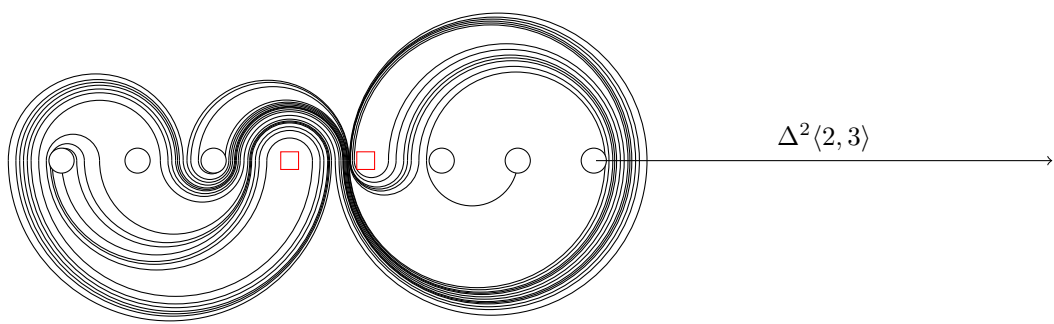
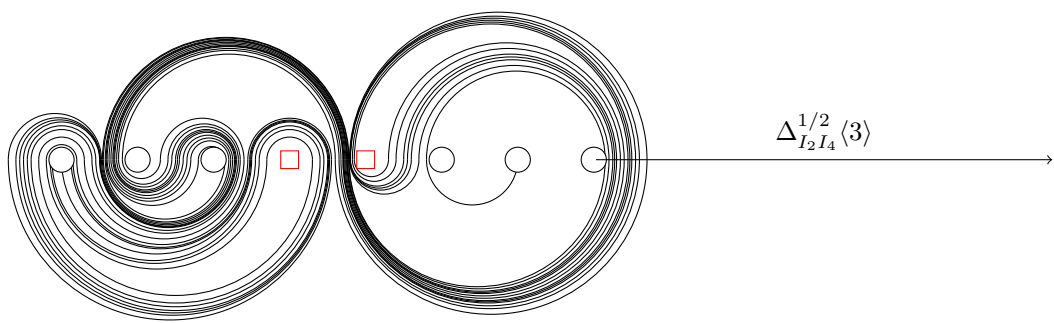
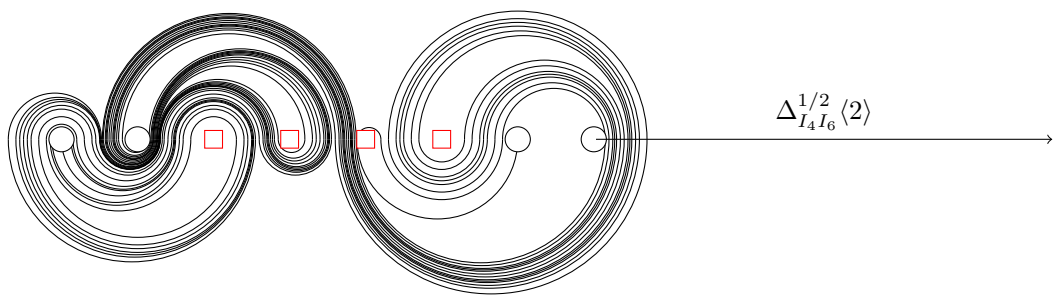
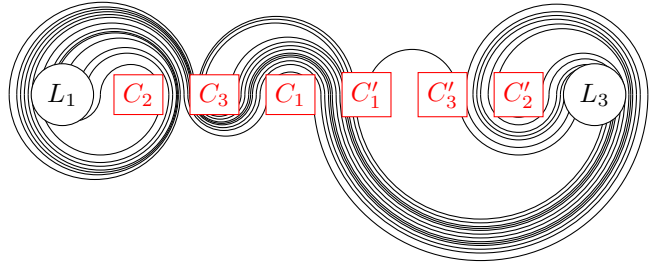
\begin{figure}[H]
\begin{center}

\end{center}
\end{figure}

Relation:
\begin{equation*}\label{arr2auto_calc_vert_22_rel_1}
\begin{split}
\Big[\Gamma_{C_3'}^{-1}\Gamma_{L_3}^{-1}\Gamma_{C_1}\Gamma_{C_3}\Gamma_{C_1}^{-1}\Gamma_{L_3}
\Gamma_{C_2'}^{-1}\Gamma_{L_3}^{-1}\Gamma_{C_1}\Gamma_{C_2}\Gamma_{L_1}\Gamma_{C_2}\Gamma_{L_1}
\Gamma_{C_2}^{-1}\Gamma_{L_1}^{-1}\Gamma_{C_2}^{-1}\Gamma_{C_1}^{-1}\Gamma_{L_3}\Gamma_{C_2'}\Gamma_{L_3}^{-1} \times  & \\
\times \Gamma_{C_1}\Gamma_{C_3}^{-1}\Gamma_{C_1}^{-1}\Gamma_{L_3}\Gamma_{C_2'}^{-1}
\Gamma_{L_3}^{-1}\Gamma_{C_1}\Gamma_{C_2}\Gamma_{L_1}\Gamma_{C_2}\Gamma_{L_1}^{-1}\Gamma_{C_2}^{-1}
\Gamma_{L_1}^{-1}\Gamma_{C_2}^{-1}\Gamma_{C_1}^{-1}\Gamma_{L_3}\Gamma_{C_2'}\Gamma_{L_3}^{-1}\Gamma_{C_1}\Gamma_{C_2}\times  & \\
\times \Gamma_{L_1}\Gamma_{C_2}\Gamma_{L_1}\Gamma_{C_2}^{-1}\Gamma_{L_1}^{-1}
\Gamma_{C_2}^{-1}\Gamma_{C_1}^{-1}\Gamma_{L_3}\Gamma_{C_2'}\Gamma_{L_3}^{-1}\Gamma_{C_1}\Gamma_{C_3}
\Gamma_{C_1}^{-1}\Gamma_{L_3}\Gamma_{C_2'}^{-1}\Gamma_{L_3}^{-1}\Gamma_{C_1}\Gamma_{C_2}\Gamma_{L_1}\Gamma_{C_2} \times & \\
\times \Gamma_{L_1}\Gamma_{C_2}^{-1}\Gamma_{L_1}^{-1}\Gamma_{C_2}^{-1}\Gamma_{C_1}^{-1}
\Gamma_{L_3}\Gamma_{C_2'}\Gamma_{L_3}^{-1}\Gamma_{C_1}\Gamma_{C_3}^{-1}\Gamma_{C_1}^{-1}\Gamma_{L_3}
\Gamma_{C_2'}^{-1}\Gamma_{L_3}^{-1}\Gamma_{C_1}\Gamma_{C_2}\Gamma_{L_1}\Gamma_{C_2}\Gamma_{L_1}^{-1}\Gamma_{C_2}^{-1} \times & \\
\times \Gamma_{L_1}^{-1}\Gamma_{C_2}^{-1}\Gamma_{C_1}^{-1}\Gamma_{L_3}\Gamma_{C_2'}^{-1}
\Gamma_{L_3}^{-1}\Gamma_{C_1}\Gamma_{C_2}\Gamma_{L_1}\Gamma_{C_2}\Gamma_{L_1}\Gamma_{C_2}^{-1}
\Gamma_{L_1}^{-1}\Gamma_{C_2}^{-1}\Gamma_{C_1}^{-1}\Gamma_{L_3}\Gamma_{C_2'}\Gamma_{L_3}^{-1}\Gamma_{C_1}\Gamma_{C_3}\times & \\
 \times  \Gamma_{C_1}^{-1}\Gamma_{L_3}\Gamma_{C_2'}^{-1}\Gamma_{L_3}^{-1}\Gamma_{C_1}\Gamma_{C_2}
\Gamma_{L_1}\Gamma_{C_2}\Gamma_{L_1}^{-1}\Gamma_{C_2}^{-1}\Gamma_{L_1}^{-1}\Gamma_{C_2}^{-1}
\Gamma_{C_1}^{-1}\Gamma_{L_3}\Gamma_{C_2'}\Gamma_{L_3}^{-1}\Gamma_{C_1}\Gamma_{C_3}^{-1}\Gamma_{C_1}^{-1} \Gamma_{L_3}\Gamma_{C_3'} , & \Gamma_{C_1'}  \Big]=e,
\end{split}
\end{equation*}
\subsubsection{Vertex Number 22}

\begin{figure}[H]
\begin{center}

\end{center}
\end{figure}

Relation:
\begin{equation*}\label{arr2auto_calc_vert_23_rel_1}
\begin{split}
\Big[\Gamma_{C_3'}^{-1}\Gamma_{L_3}^{-1}\Gamma_{C_1}\Gamma_{C_3}\Gamma_{C_1}^{-1}\Gamma_{L_3}
\Gamma_{C_2'}^{-1}\Gamma_{L_3}^{-1}\Gamma_{C_1}\Gamma_{C_2}\Gamma_{L_1}\Gamma_{C_2}\Gamma_{L_1}
\Gamma_{C_2}^{-1}\Gamma_{L_1}^{-1}\Gamma_{C_2}^{-1}\Gamma_{C_1}^{-1}\Gamma_{L_3}\Gamma_{C_2'}\Gamma_{L_3}^{-1} \times &\\
\times \Gamma_{C_1}\Gamma_{C_3}^{-1}\Gamma_{C_1}^{-1}\Gamma_{L_3}\Gamma_{C_2'}^{-1}
\Gamma_{L_3}^{-1}\Gamma_{C_1}\Gamma_{C_2}\Gamma_{L_1}\Gamma_{C_2}\Gamma_{L_1}^{-1}\Gamma_{C_2}^{-1}
\Gamma_{L_1}^{-1}\Gamma_{C_2}^{-1}\Gamma_{C_1}^{-1}\Gamma_{L_3}\Gamma_{C_2'}\Gamma_{L_3}^{-1}\Gamma_{C_1}\Gamma_{C_2} \times &\\
\times \Gamma_{L_1}\Gamma_{C_2}\Gamma_{L_1}\Gamma_{C_2}^{-1}\Gamma_{L_1}^{-1}
\Gamma_{C_2}^{-1}\Gamma_{C_1}^{-1}\Gamma_{L_3}\Gamma_{C_2'}\Gamma_{L_3}^{-1}\Gamma_{C_1}\Gamma_{C_3}
\Gamma_{C_1}^{-1}\Gamma_{L_3}\Gamma_{C_2'}^{-1}\Gamma_{L_3}^{-1}\Gamma_{C_1}\Gamma_{C_2}\Gamma_{L_1}
\Gamma_{C_2} \times &\\
\times \Gamma_{L_1}
\Gamma_{C_2}^{-1}\Gamma_{L_1}^{-1}\Gamma_{C_2}^{-1}\Gamma_{C_1}^{-1}\Gamma_{L_3}
\Gamma_{C_2'}\Gamma_{L_3}^{-1}\Gamma_{C_1}\Gamma_{C_3}^{-1}\Gamma_{C_1}^{-1}\Gamma_{L_3}\Gamma_{C_2'}^{-1}
\Gamma_{L_3}^{-1}\Gamma_{C_1}\Gamma_{C_2}\Gamma_{L_1}\Gamma_{C_2}\Gamma_{L_1}^{-1}\Gamma_{C_2}^{-1} \times &\\
\times \Gamma_{L_1}^{-1}\Gamma_{C_2}^{-1}\Gamma_{C_1}^{-1}\Gamma_{L_3}\Gamma_{C_2'}^{-1}\Gamma_{L_3}^{-1}
\Gamma_{C_1}\Gamma_{C_2}\Gamma_{L_1}\Gamma_{C_2}\Gamma_{L_1}\Gamma_{C_2}^{-1}\Gamma_{L_1}^{-1}
\Gamma_{C_2}^{-1}\Gamma_{C_1}^{-1}\Gamma_{L_3}\Gamma_{C_2'}\Gamma_{L_3}^{-1}\Gamma_{C_1}\Gamma_{C_3}\times &\\
 \times \Gamma_{C_1}^{-1}\Gamma_{L_3}\Gamma_{C_2'}^{-1}\Gamma_{L_3}^{-1}\Gamma_{C_1}\Gamma_{C_2}\Gamma_{L_1}
\Gamma_{C_2}\Gamma_{L_1}^{-1}\Gamma_{C_2}^{-1}\Gamma_{L_1}^{-1}\Gamma_{C_2}^{-1}\Gamma_{C_1}^{-1}
\Gamma_{L_3}\Gamma_{C_2'}\Gamma_{L_3}^{-1}\Gamma_{C_1}\Gamma_{C_3}^{-1}\Gamma_{C_1}^{-1}\Gamma_{L_3}
\Gamma_{C_3'}  &, \Gamma_{L_3}  \Big]=e,
\end{split}
\end{equation*}
\subsubsection{Vertex Number 23}

\begin{figure}[H]
\begin{center}

\end{center}
\end{figure}

Relation:
\begin{equation*}\label{arr2auto_calc_vert_24_rel_1}
\begin{split}
\Gamma_{C_3'}^{-1}\Gamma_{L_3}^{-1}\Gamma_{C_1}\Gamma_{C_3}\Gamma_{C_1}^{-1}\Gamma_{L_3}\Gamma_{C_2'}^{-1}
\Gamma_{L_3}^{-1}\Gamma_{C_1}\Gamma_{C_2}\Gamma_{L_1}\Gamma_{C_2}\Gamma_{L_1}\Gamma_{C_2}^{-1}
\Gamma_{L_1}^{-1}\Gamma_{C_2}^{-1}\Gamma_{C_1}^{-1}\Gamma_{L_3}\Gamma_{C_2'}\Gamma_{L_3}^{-1}\Gamma_{C_1} \times &\\
\times \Gamma_{C_3}^{-1}\Gamma_{C_1}^{-1}\Gamma_{L_3}\Gamma_{C_2'}^{-1}\Gamma_{L_3}^{-1}\Gamma_{C_1}\Gamma_{C_2}
\Gamma_{L_1}\Gamma_{C_2}\Gamma_{L_1}^{-1}\Gamma_{C_2}^{-1}\Gamma_{L_1}^{-1}\Gamma_{C_2}^{-1}\Gamma_{C_1}^{-1}
\Gamma_{L_3}\Gamma_{C_2'}\Gamma_{L_3}^{-1}\Gamma_{C_1}\Gamma_{C_2}\Gamma_{L_1} \times & \\
\times \Gamma_{C_2}\Gamma_{L_1}\Gamma_{C_2}^{-1}\Gamma_{L_1}^{-1}\Gamma_{C_2}^{-1}\Gamma_{C_1}^{-1}\Gamma_{L_3}\Gamma_{C_2'}
\Gamma_{L_3}^{-1}\Gamma_{C_1}\Gamma_{C_3}\Gamma_{C_1}^{-1}\Gamma_{L_3}\Gamma_{C_2'}^{-1}\Gamma_{L_3}^{-1}
\Gamma_{C_1}\Gamma_{C_2}\Gamma_{L_1}\Gamma_{C_2}\Gamma_{L_1}^{-1} \times &\\
\times \Gamma_{C_2}^{-1}\Gamma_{L_1}^{-1}\Gamma_{C_2}^{-1}\Gamma_{C_1}^{-1}\Gamma_{L_3}\Gamma_{C_2'}\Gamma_{L_3}^{-1}\Gamma_{C_1}\Gamma_{C_3}^{-1}
\Gamma_{C_1}^{-1}\Gamma_{L_3}\Gamma_{C_2'}^{-1}\Gamma_{L_3}^{-1}\Gamma_{C_1}\Gamma_{C_2}\Gamma_{L_1}
\Gamma_{C_2}\Gamma_{L_1}^{-1}\Gamma_{C_2}^{-1}\Gamma_{C_1} \times &\\
\times \Gamma_{C_2}\Gamma_{L_1}\Gamma_{C_2}^{-1}\Gamma_{L_1}^{-1}\Gamma_{C_2}^{-1}\Gamma_{C_1}^{-1}\Gamma_{L_3}\Gamma_{C_2'}\Gamma_{L_3}^{-1}\Gamma_{C_1}
\Gamma_{C_3}\Gamma_{C_1}^{-1}\Gamma_{L_3}\Gamma_{C_2'}^{-1}\Gamma_{L_3}^{-1}\Gamma_{C_1}\Gamma_{C_2}
\Gamma_{L_1}\Gamma_{C_2}\Gamma_{L_1} \times &\\
\times \Gamma_{C_2}^{-1}\Gamma_{L_1}^{-1}\Gamma_{C_2}^{-1}\Gamma_{C_1}^{-1}\Gamma_{L_3}\Gamma_{C_2'}\Gamma_{L_3}^{-1}\Gamma_{C_1}\Gamma_{C_3}^{-1}\Gamma_{C_1}^{-1}\Gamma_{L_3}
\Gamma_{C_2'}^{-1}\Gamma_{L_3}^{-1}\Gamma_{C_1}\Gamma_{C_2}\Gamma_{L_1}\Gamma_{C_2}\Gamma_{L_1}^{-1}
\Gamma_{C_2}^{-1}\Gamma_{L_1}^{-1}\times &\\
\times \Gamma_{C_2}^{-1}\Gamma_{C_1}^{-1}\Gamma_{L_3}\Gamma_{C_2'}^{-1}\Gamma_{L_3}^{-1}\Gamma_{C_1}\Gamma_{C_2}\Gamma_{L_1}\Gamma_{C_2}\Gamma_{L_1}\Gamma_{C_2}^{-1}
\Gamma_{L_1}^{-1}\Gamma_{C_2}^{-1}\Gamma_{C_1}^{-1}\Gamma_{L_3}\Gamma_{C_2'}\Gamma_{L_3}^{-1}\Gamma_{C_1}
\Gamma_{C_3}\Gamma_{C_1}^{-1} \times &\\
\times \Gamma_{L_3}\Gamma_{C_2'}^{-1}\Gamma_{L_3}^{-1}\Gamma_{C_1}\Gamma_{C_2}\Gamma_{L_1}\Gamma_{C_2}\Gamma_{L_1}^{-1}\Gamma_{C_2}^{-1}\Gamma_{L_1}^{-1}\Gamma_{C_2}^{-1}\Gamma_{C_1}^{-1}
\Gamma_{L_3}\Gamma_{C_2'}\Gamma_{L_3}^{-1}\Gamma_{C_1}\Gamma_{C_3}^{-1}\Gamma_{C_1}^{-1}\Gamma_{L_3}
\Gamma_{C_3'}  &= \Gamma_{C_1'} ,
\end{split}
\end{equation*}
\subsubsection{Raw relations}

\begin{equation}\label{arr2auto_calc_vert_1_rel_1}
\begin{split}
\Gamma_{C_1}  = \Gamma_{C_1'} ,
\end{split}
\end{equation}
\begin{equation}\label{arr2auto_calc_vert_2_rel_1}
\begin{split}
\Gamma_{C_1'}\Gamma_{C_3}\Gamma_{C_1'}^{-1}  = \Gamma_{C_3'} ,
\end{split}
\end{equation}
\begin{equation}\label{arr2auto_calc_vert_3_rel_1}
\begin{split}
\left\{\Gamma_{C_1} , \Gamma_{C_3}  \right\}=e,
\end{split}
\end{equation}
\begin{equation}\label{arr2auto_calc_vert_5_rel_1}
\begin{split}
\Gamma_{C_1'}\Gamma_{C_2}\Gamma_{C_1'}^{-1}  = \Gamma_{C_2'} ,
\end{split}
\end{equation}
\begin{equation}\label{arr2auto_calc_vert_6_rel_1}
\begin{split}
\left[\Gamma_{C_1'}^{-1}\Gamma_{C_3'}\Gamma_{C_1'} , \Gamma_{C_2}  \right]=e,
\end{split}
\end{equation}
\begin{equation}\label{arr2auto_calc_vert_7_rel_1}
\begin{split}
\left\{\Gamma_{C_2'} , \Gamma_{C_1'}  \right\}=e,
\end{split}
\end{equation}
\begin{equation}\label{arr2auto_calc_vert_8_rel_1}
\begin{split}
\left[\Gamma_{C_2'}\Gamma_{C_1'}\Gamma_{C_2'}^{-1} , \Gamma_{L_3}  \right]=e,
\end{split}
\end{equation}
\begin{equation}\label{arr2auto_calc_vert_9_rel_1}
\begin{split}
\left\{\Gamma_{C_2'} , \Gamma_{L_3}  \right\}=e,
\end{split}
\end{equation}
\begin{equation}\label{arr2auto_calc_vert_10_rel_1}
\begin{split}
\left[\Gamma_{C_1}\Gamma_{C_3}\Gamma_{C_1}^{-1} , \Gamma_{C_2}  \right]=e,
\end{split}
\end{equation}
\begin{equation}\label{arr2auto_calc_vert_11_rel_1}
\begin{split}
\left[\Gamma_{L_3}^{-1}\Gamma_{C_1'}^{-1}\Gamma_{C_3'}\Gamma_{C_1'}\Gamma_{L_3} , \Gamma_{C_2'}  \right]=e,
\end{split}
\end{equation}
\begin{equation}\label{arr2auto_calc_vert_12_rel_1}
\begin{split}
\left\{\Gamma_{C_1'}^{-1}\Gamma_{C_3'}\Gamma_{C_1'} , \Gamma_{L_3}  \right\}=e,
\end{split}
\end{equation}
\begin{equation}\label{arr2auto_calc_vert_13_rel_1}
\begin{split}
\left\{\Gamma_{C_1} , \Gamma_{C_2}  \right\}=e,
\end{split}
\end{equation}
\begin{equation}\label{arr2auto_calc_vert_14_rel_1}
\begin{split}
\left[\Gamma_{L_3}^{-1}\Gamma_{C_1}\Gamma_{C_3}\Gamma_{C_1}^{-1}\Gamma_{L_3} , \Gamma_{C_2'}  \right]=e,
\end{split}
\end{equation}
\begin{equation}\label{arr2auto_calc_vert_15_rel_1}
\begin{split}
\left[\Gamma_{C_3'}\Gamma_{C_1'}\Gamma_{L_3}\Gamma_{C_1'}^{-1}\Gamma_{C_3'}^{-1} , \Gamma_{C_1'}  \right]=e,
\end{split}
\end{equation}
\begin{equation}\label{arr2auto_calc_vert_16_rel_1}
\begin{split}
\left[\Gamma_{L_1} , \Gamma_{C_1}  \right]=e,
\end{split}
\end{equation}
\begin{equation}\label{arr2auto_calc_vert_17_rel_1}
\begin{split}
\left\{\Gamma_{L_1} , \Gamma_{C_2}  \right\}=e,
\end{split}
\end{equation}
\begin{equation}\label{arr2auto_calc_vert_18_rel_1}
\begin{split}
\Gamma_{L_3}^{-1}\Gamma_{C_1}\Gamma_{C_2}\Gamma_{L_1}\Gamma_{C_2}\Gamma_{L_1}^{-1}\Gamma_{C_2}^{-1}\Gamma_{C_1}^{-1}\Gamma_{L_3}  = \Gamma_{C_2'} ,
\end{split}
\end{equation}
\begin{equation}\label{arr2auto_calc_vert_19_rel_1}
\begin{split}
\left\{\Gamma_{C_3'} , \Gamma_{C_1'}  \right\}=e,
\end{split}
\end{equation}
\begin{equation}\label{arr2auto_calc_vert_20_rel_1}
\begin{split}
\left\{\Gamma_{C_1}^{-1}\Gamma_{L_3}\Gamma_{C_2'}^{-1}\Gamma_{L_3}^{-1}\Gamma_{C_1}\Gamma_{C_2}\Gamma_{L_1}\Gamma_{C_2}\Gamma_{L_1}\Gamma_{C_2}^{-1}\Gamma_{L_1}^{-1}\Gamma_{C_2}^{-1}\Gamma_{C_1}^{-1}\Gamma_{L_3}\Gamma_{C_2'}\Gamma_{L_3}^{-1}\Gamma_{C_1} , \Gamma_{C_3}  \right\}=e,
\end{split}
\end{equation}
\begin{equation}\label{arr2auto_calc_vert_21_rel_1}
\begin{split}
\Gamma_{L_3}^{-1}\Gamma_{C_1}\Gamma_{C_3}\Gamma_{C_1}^{-1}\Gamma_{L_3}\Gamma_{C_2'}^{-1}\Gamma_{L_3}^{-1}
\Gamma_{C_1}\Gamma_{C_2}\Gamma_{L_1}\Gamma_{C_2}\Gamma_{L_1}\Gamma_{C_2}^{-1}\Gamma_{L_1}^{-1}
\Gamma_{C_2}^{-1}\Gamma_{C_1}^{-1}\Gamma_{L_3}\Gamma_{C_2'}\Gamma_{L_3}^{-1}\Gamma_{C_1}\Gamma_{C_3}\\
\Gamma_{C_1}^{-1}\Gamma_{L_3}\Gamma_{C_2'}^{-1}\Gamma_{L_3}^{-1}\Gamma_{C_1}\Gamma_{C_2}\Gamma_{L_1}
\Gamma_{C_2}\Gamma_{L_1}^{-1}\Gamma_{C_2}^{-1}\Gamma_{L_1}^{-1}\Gamma_{C_2}^{-1}\Gamma_{C_1}^{-1}
\Gamma_{L_3}\Gamma_{C_2'}\Gamma_{L_3}^{-1}\Gamma_{C_1}\Gamma_{C_3}^{-1}\Gamma_{C_1}^{-1}\Gamma_{L_3}  = \Gamma_{C_3'} ,
\end{split}
\end{equation}
\begin{equation}\label{arr2auto_calc_vert_22_rel_1}
\begin{split}
\Big[\Gamma_{C_3'}^{-1}\Gamma_{L_3}^{-1}\Gamma_{C_1}\Gamma_{C_3}\Gamma_{C_1}^{-1}\Gamma_{L_3}
\Gamma_{C_2'}^{-1}\Gamma_{L_3}^{-1}\Gamma_{C_1}\Gamma_{C_2}\Gamma_{L_1}\Gamma_{C_2}\Gamma_{L_1}
\Gamma_{C_2}^{-1}\Gamma_{L_1}^{-1}\Gamma_{C_2}^{-1}\Gamma_{C_1}^{-1}\Gamma_{L_3}\Gamma_{C_2'}\\
\Gamma_{L_3}^{-1}\Gamma_{C_1}\Gamma_{C_3}^{-1}\Gamma_{C_1}^{-1}\Gamma_{L_3}\Gamma_{C_2'}^{-1}
\Gamma_{L_3}^{-1}\Gamma_{C_1}\Gamma_{C_2}\Gamma_{L_1}\Gamma_{C_2}\Gamma_{L_1}^{-1}\Gamma_{C_2}^{-1}
\Gamma_{L_1}^{-1}\Gamma_{C_2}^{-1}\Gamma_{C_1}^{-1}\Gamma_{L_3}\Gamma_{C_2'}\Gamma_{L_3}^{-1}\\
\Gamma_{C_1}\Gamma_{C_2}\Gamma_{L_1}\Gamma_{C_2}\Gamma_{L_1}\Gamma_{C_2}^{-1}\Gamma_{L_1}^{-1}
\Gamma_{C_2}^{-1}\Gamma_{C_1}^{-1}\Gamma_{L_3}\Gamma_{C_2'}\Gamma_{L_3}^{-1}\Gamma_{C_1}\Gamma_{C_3}
\Gamma_{C_1}^{-1}\Gamma_{L_3}\Gamma_{C_2'}^{-1}\Gamma_{L_3}^{-1}\Gamma_{C_1}\Gamma_{C_2}\Gamma_{L_1}\\
\Gamma_{C_2}\Gamma_{L_1}\Gamma_{C_2}^{-1}\Gamma_{L_1}^{-1}\Gamma_{C_2}^{-1}\Gamma_{C_1}^{-1}
\Gamma_{L_3}\Gamma_{C_2'}\Gamma_{L_3}^{-1}\Gamma_{C_1}\Gamma_{C_3}^{-1}\Gamma_{C_1}^{-1}\Gamma_{L_3}
\Gamma_{C_2'}^{-1}\Gamma_{L_3}^{-1}\Gamma_{C_1}\Gamma_{C_2}\Gamma_{L_1}\Gamma_{C_2}\Gamma_{L_1}^{-1}\\
\Gamma_{C_2}^{-1}\Gamma_{L_1}^{-1}\Gamma_{C_2}^{-1}\Gamma_{C_1}^{-1}\Gamma_{L_3}\Gamma_{C_2'}^{-1}
\Gamma_{L_3}^{-1}\Gamma_{C_1}\Gamma_{C_2}\Gamma_{L_1}\Gamma_{C_2}\Gamma_{L_1}\Gamma_{C_2}^{-1}
\Gamma_{L_1}^{-1}\Gamma_{C_2}^{-1}\Gamma_{C_1}^{-1}\Gamma_{L_3}\Gamma_{C_2'}\Gamma_{L_3}^{-1}\Gamma_{C_1}\\
\Gamma_{C_3}\Gamma_{C_1}^{-1}\Gamma_{L_3}\Gamma_{C_2'}^{-1}\Gamma_{L_3}^{-1}\Gamma_{C_1}\Gamma_{C_2}
\Gamma_{L_1}\Gamma_{C_2}\Gamma_{L_1}^{-1}\Gamma_{C_2}^{-1}\Gamma_{L_1}^{-1}\Gamma_{C_2}^{-1}
\Gamma_{C_1}^{-1}\Gamma_{L_3}\Gamma_{C_2'}\Gamma_{L_3}^{-1}\Gamma_{C_1}\Gamma_{C_3}^{-1}\Gamma_{C_1}^{-1}\\
\Gamma_{L_3}\Gamma_{C_3'} , \Gamma_{C_1'}  \Big]=e,
\end{split}
\end{equation}
\begin{equation}\label{arr2auto_calc_vert_23_rel_1}
\begin{split}
\Big[\Gamma_{C_3'}^{-1}\Gamma_{L_3}^{-1}\Gamma_{C_1}\Gamma_{C_3}\Gamma_{C_1}^{-1}\Gamma_{L_3}
\Gamma_{C_2'}^{-1}\Gamma_{L_3}^{-1}\Gamma_{C_1}\Gamma_{C_2}\Gamma_{L_1}\Gamma_{C_2}\Gamma_{L_1}
\Gamma_{C_2}^{-1}\Gamma_{L_1}^{-1}\Gamma_{C_2}^{-1}\Gamma_{C_1}^{-1}\Gamma_{L_3}\Gamma_{C_2'}\\
\Gamma_{L_3}^{-1}\Gamma_{C_1}\Gamma_{C_3}^{-1}\Gamma_{C_1}^{-1}\Gamma_{L_3}\Gamma_{C_2'}^{-1}
\Gamma_{L_3}^{-1}\Gamma_{C_1}\Gamma_{C_2}\Gamma_{L_1}\Gamma_{C_2}\Gamma_{L_1}^{-1}\Gamma_{C_2}^{-1}
\Gamma_{L_1}^{-1}\Gamma_{C_2}^{-1}\Gamma_{C_1}^{-1}\Gamma_{L_3}\Gamma_{C_2'}\Gamma_{L_3}^{-1}\\
\Gamma_{C_1}\Gamma_{C_2}\Gamma_{L_1}\Gamma_{C_2}\Gamma_{L_1}\Gamma_{C_2}^{-1}\Gamma_{L_1}^{-1}
\Gamma_{C_2}^{-1}\Gamma_{C_1}^{-1}\Gamma_{L_3}\Gamma_{C_2'}\Gamma_{L_3}^{-1}\Gamma_{C_1}\Gamma_{C_3}
\Gamma_{C_1}^{-1}\Gamma_{L_3}\Gamma_{C_2'}^{-1}\Gamma_{L_3}^{-1}\Gamma_{C_1}\Gamma_{C_2}\Gamma_{L_1}
\Gamma_{C_2}\\
\Gamma_{L_1}
\Gamma_{C_2}^{-1}\Gamma_{L_1}^{-1}\Gamma_{C_2}^{-1}\Gamma_{C_1}^{-1}\Gamma_{L_3}
\Gamma_{C_2'}\Gamma_{L_3}^{-1}\Gamma_{C_1}\Gamma_{C_3}^{-1}\Gamma_{C_1}^{-1}\Gamma_{L_3}\Gamma_{C_2'}^{-1}
\Gamma_{L_3}^{-1}\Gamma_{C_1}\Gamma_{C_2}\Gamma_{L_1}\Gamma_{C_2}\Gamma_{L_1}^{-1}\Gamma_{C_2}^{-1}\\
\Gamma_{L_1}^{-1}\Gamma_{C_2}^{-1}\Gamma_{C_1}^{-1}\Gamma_{L_3}\Gamma_{C_2'}^{-1}\Gamma_{L_3}^{-1}
\Gamma_{C_1}\Gamma_{C_2}\Gamma_{L_1}\Gamma_{C_2}\Gamma_{L_1}\Gamma_{C_2}^{-1}\Gamma_{L_1}^{-1}
\Gamma_{C_2}^{-1}\Gamma_{C_1}^{-1}\Gamma_{L_3}\Gamma_{C_2'}\Gamma_{L_3}^{-1}\Gamma_{C_1}\Gamma_{C_3}\\
\Gamma_{C_1}^{-1}\Gamma_{L_3}\Gamma_{C_2'}^{-1}\Gamma_{L_3}^{-1}\Gamma_{C_1}\Gamma_{C_2}\Gamma_{L_1}
\Gamma_{C_2}\Gamma_{L_1}^{-1}\Gamma_{C_2}^{-1}\Gamma_{L_1}^{-1}\Gamma_{C_2}^{-1}\Gamma_{C_1}^{-1}
\Gamma_{L_3}\Gamma_{C_2'}\Gamma_{L_3}^{-1}\Gamma_{C_1}\Gamma_{C_3}^{-1}\Gamma_{C_1}^{-1}\Gamma_{L_3}
\Gamma_{C_3'} , \Gamma_{L_3}  \Big]=e,
\end{split}
\end{equation}
\begin{equation}\label{arr2auto_calc_vert_24_rel_1}
\begin{split}
\Gamma_{C_3'}^{-1}\Gamma_{L_3}^{-1}\Gamma_{C_1}\Gamma_{C_3}\Gamma_{C_1}^{-1}\Gamma_{L_3}\Gamma_{C_2'}^{-1}
\Gamma_{L_3}^{-1}\Gamma_{C_1}\Gamma_{C_2}\Gamma_{L_1}\Gamma_{C_2}\Gamma_{L_1}\Gamma_{C_2}^{-1}
\Gamma_{L_1}^{-1}\Gamma_{C_2}^{-1}\Gamma_{C_1}^{-1}\Gamma_{L_3}\Gamma_{C_2'}\Gamma_{L_3}^{-1}\Gamma_{C_1} \times &\\
\times \Gamma_{C_3}^{-1}\Gamma_{C_1}^{-1}\Gamma_{L_3}\Gamma_{C_2'}^{-1}\Gamma_{L_3}^{-1}\Gamma_{C_1}\Gamma_{C_2}
\Gamma_{L_1}\Gamma_{C_2}\Gamma_{L_1}^{-1}\Gamma_{C_2}^{-1}\Gamma_{L_1}^{-1}\Gamma_{C_2}^{-1}\Gamma_{C_1}^{-1}
\Gamma_{L_3}\Gamma_{C_2'}\Gamma_{L_3}^{-1}\Gamma_{C_1}\Gamma_{C_2}\Gamma_{L_1} \times & \\
\times \Gamma_{C_2}\Gamma_{L_1}\Gamma_{C_2}^{-1}\Gamma_{L_1}^{-1}\Gamma_{C_2}^{-1}\Gamma_{C_1}^{-1}\Gamma_{L_3}\Gamma_{C_2'}
\Gamma_{L_3}^{-1}\Gamma_{C_1}\Gamma_{C_3}\Gamma_{C_1}^{-1}\Gamma_{L_3}\Gamma_{C_2'}^{-1}\Gamma_{L_3}^{-1}
\Gamma_{C_1}\Gamma_{C_2}\Gamma_{L_1}\Gamma_{C_2}\Gamma_{L_1}^{-1} \times &\\
\times \Gamma_{C_2}^{-1}\Gamma_{L_1}^{-1}\Gamma_{C_2}^{-1}\Gamma_{C_1}^{-1}\Gamma_{L_3}\Gamma_{C_2'}\Gamma_{L_3}^{-1}\Gamma_{C_1}\Gamma_{C_3}^{-1}
\Gamma_{C_1}^{-1}\Gamma_{L_3}\Gamma_{C_2'}^{-1}\Gamma_{L_3}^{-1}\Gamma_{C_1}\Gamma_{C_2}\Gamma_{L_1}
\Gamma_{C_2}\Gamma_{L_1}^{-1}\Gamma_{C_2}^{-1}\Gamma_{C_1} \times &\\
\times \Gamma_{C_2}\Gamma_{L_1}\Gamma_{C_2}^{-1}\Gamma_{L_1}^{-1}\Gamma_{C_2}^{-1}\Gamma_{C_1}^{-1}\Gamma_{L_3}\Gamma_{C_2'}\Gamma_{L_3}^{-1}\Gamma_{C_1}
\Gamma_{C_3}\Gamma_{C_1}^{-1}\Gamma_{L_3}\Gamma_{C_2'}^{-1}\Gamma_{L_3}^{-1}\Gamma_{C_1}\Gamma_{C_2}
\Gamma_{L_1}\Gamma_{C_2}\Gamma_{L_1} \times &\\
\times \Gamma_{C_2}^{-1}\Gamma_{L_1}^{-1}\Gamma_{C_2}^{-1}\Gamma_{C_1}^{-1}\Gamma_{L_3}\Gamma_{C_2'}\Gamma_{L_3}^{-1}\Gamma_{C_1}\Gamma_{C_3}^{-1}\Gamma_{C_1}^{-1}\Gamma_{L_3}
\Gamma_{C_2'}^{-1}\Gamma_{L_3}^{-1}\Gamma_{C_1}\Gamma_{C_2}\Gamma_{L_1}\Gamma_{C_2}\Gamma_{L_1}^{-1}
\Gamma_{C_2}^{-1}\Gamma_{L_1}^{-1}\times &\\
\times \Gamma_{C_2}^{-1}\Gamma_{C_1}^{-1}\Gamma_{L_3}\Gamma_{C_2'}^{-1}\Gamma_{L_3}^{-1}\Gamma_{C_1}\Gamma_{C_2}\Gamma_{L_1}\Gamma_{C_2}\Gamma_{L_1}\Gamma_{C_2}^{-1}
\Gamma_{L_1}^{-1}\Gamma_{C_2}^{-1}\Gamma_{C_1}^{-1}\Gamma_{L_3}\Gamma_{C_2'}\Gamma_{L_3}^{-1}\Gamma_{C_1}
\Gamma_{C_3}\Gamma_{C_1}^{-1} \times &\\
\times \Gamma_{L_3}\Gamma_{C_2'}^{-1}\Gamma_{L_3}^{-1}\Gamma_{C_1}\Gamma_{C_2}\Gamma_{L_1}\Gamma_{C_2}\Gamma_{L_1}^{-1}\Gamma_{C_2}^{-1}\Gamma_{L_1}^{-1}\Gamma_{C_2}^{-1}\Gamma_{C_1}^{-1}
\Gamma_{L_3}\Gamma_{C_2'}\Gamma_{L_3}^{-1}\Gamma_{C_1}\Gamma_{C_3}^{-1}\Gamma_{C_1}^{-1}\Gamma_{L_3}
\Gamma_{C_3'}  &= \Gamma_{C_1'} ,
\end{split}
\end{equation}
\begin{equation}\label{arr2auto_calc_projective_rel}
\begin{split}
\Gamma_{L_3}\Gamma_{C_2'}\Gamma_{C_3'}\Gamma_{C_1'}\Gamma_{C_1}\Gamma_{C_3}\Gamma_{C_2}\Gamma_{L_1} =e.
\end{split}
\end{equation}

\newpage

\section{Tokunaga's Examples}\label{sec:existing}
In \cite{tokunaga2014}, Tokunaga used the theory of dihedral covers and Mordell-Weil groups to find two examples of Zariski pairs of conic-line arrangements of degree 7. In this section we will show that the fundamental groups of those arrangements differ, so our method can prove the fact that those pairs of arrangements are Zariski pairs with non-isomorphic fundamental groups.

\subsection{First pair}
We will now turn to the first pair $ \{\mathcal{B}_3, \mathcal{B}_4\} $.

\begin{figure}[H]
\centering
\begin{tikzpicture} [scale=0.35, style = very thick]
	\begin{scope}[xshift=-330]
		\draw (0,0) ellipse (7 and 4);
		\draw (-7, -3) node {\footnotesize $C_2$};
		
		\draw (0,0) ellipse (4 and 7);
		\draw (-3, 7) node {\footnotesize $ C_3 $};
		
		\draw (0,0) ellipse (4 and 4);
		\draw (-2, 2) node {\footnotesize $ C_1 $};
		
		\draw (3.47, -7) -- (3.47, 7); 
		\draw (4, -6) node {\footnotesize $ L_1 $};
		
		\draw (-4.7, 0) node {\footnotesize $Q_1$};
		\draw (4.7, 0) node {\footnotesize $Q_2$};
		\draw (0, 4.7) node {\footnotesize $Q_3$};
		\draw (0, -4.7) node {\footnotesize $Q_4$};
		
		\draw (4.4, 3.7) node {\footnotesize $ P_1 $};
		\draw (-4.4, 3.7) node {\footnotesize $ P_2 $};
		\draw (4.4, -3.7) node {\footnotesize $ P_3 $};
		\draw (-4.4, -3.7) node {\footnotesize $ P_4 $};
		\draw (0, -8) node {$\mathcal{B}_3$};
	\end{scope}
	
	\begin{scope}[xshift=330]
		\draw (0,0) ellipse (7 and 4);
		\draw (-7, -3) node {\footnotesize $C_2$};
		
		\draw (0,0) ellipse (4 and 7);
		\draw (-3, 7) node {\footnotesize $ C_3 $};
		
		\draw (0,0) ellipse (4 and 4);
		\draw (-2, 2) node {\footnotesize $ C_1 $};
		
		\draw (-7, -7) -- (7, 7); 
		\draw (-4, -6) node {\footnotesize $ L_2 $};
		
		\draw (-4.7, 0) node {\footnotesize $Q_1$};
		\draw (4.7, 0) node {\footnotesize $Q_2$};
		\draw (0, 4.7) node {\footnotesize $Q_3$};
		\draw (0, -4.7) node {\footnotesize $Q_4$};
		
		\draw (4.4, 3.7) node {\footnotesize $ P_1 $};
		\draw (-4.4, 3.7) node {\footnotesize $ P_2 $};
		\draw (4.4, -3.7) node {\footnotesize $ P_3 $};
		\draw (-4.4, -3.7) node {\footnotesize $ P_4 $};
		
		\draw (0, -8) node {$\mathcal{B}_4$};
	\end{scope}
\end{tikzpicture}

\caption{The arrangements $ \mathcal{B}_3, \mathcal{B}_4 $.}
\label{fig_tokunaga2}
\end{figure}
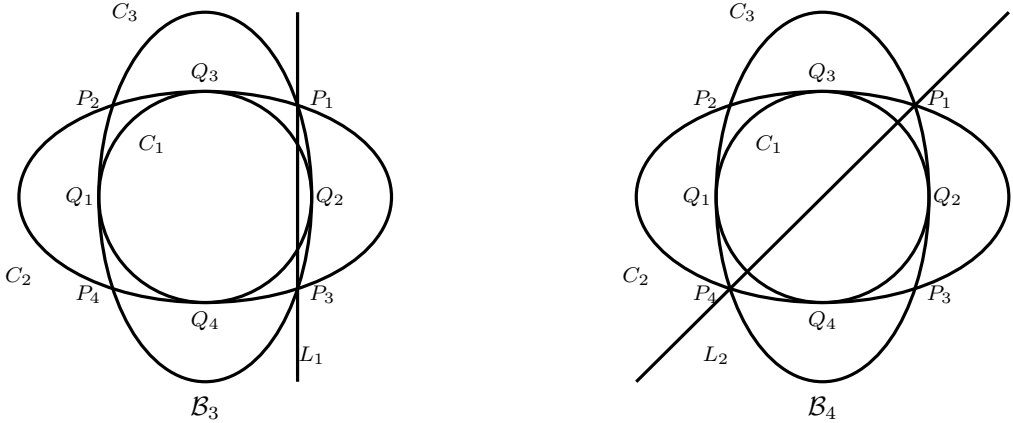

\begin{Construction}\label{construcetion_tokunaga_arrangement2_construction}
Let $ C_i \left(i=1,2,3\right)$ be smooth conics and let $ L_j \left(j=1,2\right) $ be lines as follows (see Figure \ref{fig_tokunaga2}):
\begin{itemize}
	\item $ C_2 $ and $ C_3 $ meet transversely. We put $ C_2\cap C_3=\{ P_1, P_2, P_3, P_4 \} $.
	\item $ C_1 $ is tangent to both $ C_3 $ and $ C_2 $ such that the intersection multiplicities at intersection points are all equal to 2. We set $ C_1\cap C_2=\{ Q_1, Q_2 \} $ and $ C_1\cap C_3 = \{Q_3, Q_4 \} $.
	\item $ L_1 $ passes through $ P_1 $ and $ P_3 $.
	\item $ L_2 $ passes through $ P_1 $ and $ P_4 $.
	\item $ L_1, L_2 $ meet $ C_1 $ transversely.
\end{itemize}
Denote $ \mathcal{B}_3=C_1\cup C_2\cup C_3\cup L_1 $ and $ \mathcal{B}_4=C_1\cup C_2\cup C_3 \cup L_2 $.
\end{Construction}

\subsection{Second pair}
We will turn to the second pair $ \{\mathcal{B}_5, \mathcal{B}_6\} $.
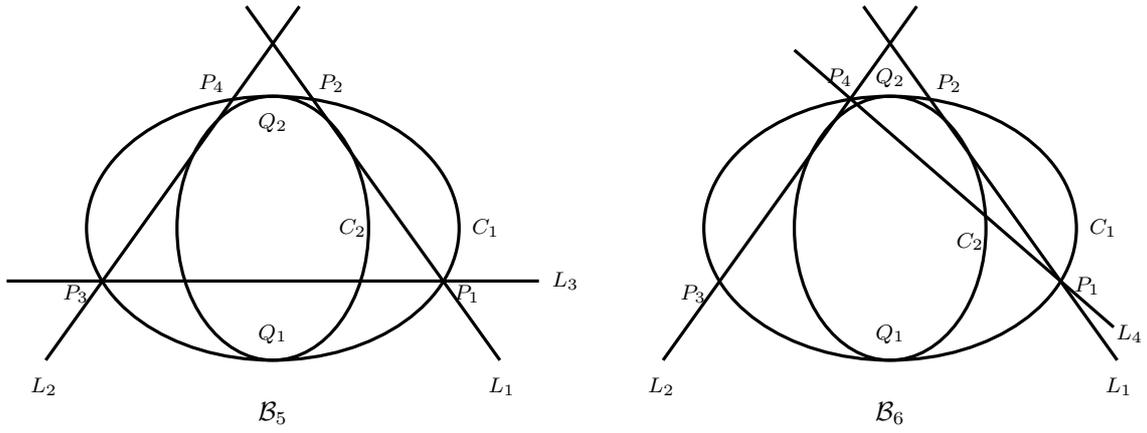
\begin{figure}
\centering
\begin{tikzpicture} [scale=0.35, style = very thick]
	\begin{scope}[xshift=-330]
		\draw (0,0) ellipse (7 and 5);
		\draw (8,0) node {\footnotesize $C_1$};
		\draw (0,0) ellipse (3.6 and 5);
		\draw (3,0) node {\footnotesize $C_2$};
		\draw (-10, -2) -- (10, -2);
		\draw (11, -2) node {\footnotesize $ L_3 $};
		\draw (-1,8.406) -- (8.533, -5);
		\draw (8.6, -6) node {\footnotesize $ L_1 $};
		\draw (1,8.406) -- (-8.533, -5);
		\draw (-8.6, -6) node {\footnotesize $ L_2 $};
		
		\draw (0, -4) node {\footnotesize $ Q_1 $};
		\draw (0, 4) node {\footnotesize $ Q_2 $};
		
		\draw (7.3, -2.5) node {\footnotesize $ P_1 $};
		\draw (-7.4, -2.5) node {\footnotesize $ P_3 $};
		\draw (2.2, 5.5) node {\footnotesize $ P_2 $};
		\draw (-2.3, 5.5) node {\footnotesize $ P_4 $};
		
		\draw (0, -7) node {$\mathcal{B}_5$};
	\end{scope}
	
	\begin{scope}[xshift=330]
		\draw (0,0) ellipse (7 and 5);
		\draw (8,0) node {\footnotesize $C_1$};
		\draw (0,0) ellipse (3.6 and 5);
		\draw (3,-0.5) node {\footnotesize $C_2$};
		\draw (-1,8.406) -- (8.533, -5);
		\draw (8.6, -6) node {\footnotesize $ L_1 $};
		\draw (1,8.406) -- (-8.533, -5);
		\draw (-8.6, -6) node {\footnotesize $ L_2 $};
		\draw (8.4, -3.75) -- (-3.6, 6.75);
		\draw (9, -4) node {\footnotesize $ L_4 $};
		
		\draw (0, -4) node {\footnotesize $ Q_1 $};
		\draw (0, 5.7) node {\footnotesize $ Q_2 $};
		
		\draw (7.4, -2.2) node {\footnotesize $ P_1 $};
		\draw (-7.4, -2.5) node {\footnotesize $ P_3 $};
		\draw (2.2, 5.5) node {\footnotesize $ P_2 $};
		\draw (-1.9, 5.7) node {\footnotesize $ P_4 $};
		
		\draw (0, -7) node {$\mathcal{B}_6$};
	\end{scope}
	
\end{tikzpicture}
\caption{The arrangements $ \mathcal{B}_5, \mathcal{B}_6 $.}
\label{fig_tokunaga1}
\end{figure}
\begin{Construction}\label{construcetion_tokunaga_arrangement1_construction}
Let $ C_i \left(i=1,2\right)$ be smooth conics and let $ L_j \left(j=1,\dots,4\right) $ be lines as follows (see Figure \ref{fig_tokunaga1}):
\begin{itemize}
	\item $ C_1 $ and $ C_2 $ are tangent at $ 2 $ distinct points $ Q_1 $ and $ Q_2 $.
	\item $ C_2 $ is tangent to $ L_1 $ and to $ L_2 $.
	\item $ L_1 $ and $ L_2 $ meet $ C_1 $ transversely. We put $ C_1\cap L_1=\left\{P_1, P_2\right\}, \; C_1\cap L_2=\left\{P_3,P_4\right\} $.
	\item $ L_3 $ is the line connecting $ P_1 $ and $ P_3 $, $ L_4 $ is the line connecting $ P_1 $ and $ P_4 $.
	\item Both $ L_3 $ and $ L_4 $ meet $ C_2 $ transversely.
\end{itemize}
Denote $ \mathcal{B}_5=C_1\cup C_2\cup L_1\cup L_2\cup L_3 $ and $ \mathcal{B}_6=C_1\cup C_2\cup L_1\cup L_2 \cup L_4 $.
\end{Construction}

\newpage

\section{Computations of $ \pcpt{\mathcal{B}_3} $ and $ \pcpt{\mathcal{B}_4} $}

\subsection{$ \pcpt{\mathcal{B}_3} $}
\begin{figure}[H]
	\begin{center}
		\begin{tikzpicture}
		\begin{scope}[ , scale=1]
		\draw [rotate around={90.000000:(0.000000,0.000000)}] (0.000000,0.000000) ellipse (4.000000 and 4.000000);
		\node  (C_1_naming_node) at (0.000000,3.700000) {$ C_1 $};
		\draw [rotate around={-70.000000:(-0.000000,-0.000000)}] (-0.000000,-0.000000) ellipse (5.656854 and 4.000000);
		\node  (C_3_naming_node) at (1.934758,-5.815704) {$ C_3 $};
		\draw [rotate around={-160.000000:(-0.000000,0.000000)}] (-0.000000,0.000000) ellipse (5.656854 and 4.000000);
		\node  (C_2_naming_node) at (-5.815704,-1.934758) {$ C_2 $};
		\draw [] (5.303910, -5.023268) node [below] {$L_1$} -- (1.647270, 5.023268);
		\draw [fill=gray] (5.488930, 0.936849) circle (5.000000pt);
		\node [label=above:1] (sing_pnt_1) at (5.488930,0.936849) {};
		\draw [fill=gray] (4.227487, -1.216397) circle (5.000000pt);
		\node [label=above:2] (sing_pnt_2) at (4.227487,-1.216397) {};
		\draw [fill=gray] (4.186056, -1.951990) circle (5.000000pt);
		\node [label=above:3] (sing_pnt_3) at (4.186056,-1.951990) {};
		\draw [fill=gray] (4.000000, -0.000000) circle (5.000000pt);
		\node [label=above:4] (sing_pnt_4) at (4.000000,-0.000000) {};
		\draw [fill=gray] (3.858885, -1.053094) circle (5.000000pt);
		\node [label=above:5] (sing_pnt_5) at (3.858885,-1.053094) {};
		\draw [fill=gray] (3.758770, 1.368081) circle (5.000000pt);
		\node [label=above:6] (sing_pnt_6) at (3.758770,1.368081) {};
		\draw [fill=gray] (2.279162, 3.287160) circle (5.000000pt);
		\node [label=above:7] (sing_pnt_7) at (2.279162,3.287160) {};
		\draw [fill=gray] (1.951990, 4.186056) circle (5.000000pt);
		\node [label=above:8] (sing_pnt_8) at (1.951990,4.186056) {};
		\draw [fill=gray] (1.368081, -3.758770) circle (5.000000pt);
		\node [label=above:9] (sing_pnt_9) at (1.368081,-3.758770) {};
		\draw [fill=gray] (-1.368081, 3.758770) circle (5.000000pt);
		\node [label=above:10] (sing_pnt_10) at (-1.368081,3.758770) {};
		\draw [fill=gray] (-1.951990, -4.186056) circle (5.000000pt);
		\node [label=above:11] (sing_pnt_11) at (-1.951990,-4.186056) {};
		\draw [fill=gray] (-3.758770, -1.368081) circle (5.000000pt);
		\node [label=above:12] (sing_pnt_12) at (-3.758770,-1.368081) {};
		\draw [fill=gray] (-4.000000, 0.000000) circle (5.000000pt);
		\node [label=above:13] (sing_pnt_13) at (-4.000000,0.000000) {};
		\draw [fill=gray] (-4.186056, 1.951990) circle (5.000000pt);
		\node [label=above:14] (sing_pnt_14) at (-4.186056,1.951990) {};
		\draw [fill=gray] (-4.227487, 1.216397) circle (5.000000pt);
		\node [label=above:15] (sing_pnt_15) at (-4.227487,1.216397) {};
		\draw [fill=gray] (-5.488930, -0.936849) circle (5.000000pt);
		\node [label=above:16] (sing_pnt_16) at (-5.488930,-0.936849) {};
		
		\end{scope}[]
		\end{tikzpicture}
	\end{center}
\end{figure}

{Monodromy table}

\begin{tabular}{| c | c | c | c |}
	\hline \hline
	Vertex number& Vertex description& Skeleton& Diffeomorphism\\ \hline \hline
	1& $C_2$ branch& $ \langle 2 - 3 \rangle $& $ \Delta_{I_{4}I_{6}}^{1/2}\langle 2 \rangle $ \\ \hline
	2& $C_3$ branch& $ \langle 3 - 4 \rangle $& $ \Delta_{I_{2}I_{4}}^{1/2}\langle 3 \rangle $ \\ \hline
	3& Node between $C_3$, $C_2$ and $L_1$& $ \langle 1 - 2 - 3 \rangle $& $ \Delta \langle 1,2,3 \rangle $ \\ \hline
	4& $C_1$ branch& $ \langle 4 - 5 \rangle $& $ \Delta_{\mathbb{R}I_{2}}^{1/2}\langle 4 \rangle $ \\ \hline
	5& Node between $C_1$ and $L_1$& $ \langle 3 - 4 \rangle $& $ \Delta\langle 3, 4 \rangle $ \\ \hline
	6& Tangency between $C_1$ and $C_3$& $ \langle 5 - 6 \rangle $& $ \Delta^{2}\langle 5, 6 \rangle $ \\ \hline
	7& Node between $C_1$ and $L_1$& $ \langle 4 - 5 \rangle $& $ \Delta\langle 4, 5 \rangle $ \\ \hline
	8& Node between $C_2$, $C_3$ and $L_1$& $ \langle 5 - 6 - 7 \rangle $& $ \Delta \langle 5,6,7 \rangle $ \\ \hline
	9& Tangency between $C_1$ and $C_2$& $ \langle 2 - 3 \rangle $& $ \Delta^{2}\langle 2, 3 \rangle $ \\ \hline
	10& Tangency between $C_1$ and $C_2$& $ \langle 4 - 5 \rangle $& $ \Delta^{2}\langle 4, 5 \rangle $ \\ \hline
	11& Node between $C_2$ and $C_3$& $ \langle 1 - 2 \rangle $& $ \Delta\langle 1, 2 \rangle $ \\ \hline
	12& Tangency between $C_1$ and $C_3$& $ \langle 2 - 3 \rangle $& $ \Delta^{2}\langle 2, 3 \rangle $ \\ \hline
	13& $C_1$ branch& $ \langle 3 - 4 \rangle $& $ \Delta_{I_{2}\mathbb{R}}^{1/2}\langle 3 \rangle $ \\ \hline
	14& Node between $C_3$ and $C_2$& $ \langle 3 - 4 \rangle $& $ \Delta\langle 3, 4 \rangle $ \\ \hline
	15& $C_3$ branch& $ \langle 2 - 3 \rangle $& $ \Delta_{I_{4}I_{2}}^{1/2}\langle 2 \rangle $ \\ \hline
	16& $C_2$ branch& $ \langle 1 - 2 \rangle $& $ \Delta_{I_{6}I_{4}}^{1/2}\langle 1 \rangle $
	\\ \hline \hline
\end{tabular}

\subsubsection{Vertex Number 1}

\begin{figure}[H]
	\begin{center}
		\begin{tikzpicture}
		\begin{scope}[]
		\draw [] (3.000000, 0.000000) -- (4.000000, 0.000000) ;
		\node [draw, circle, color=black, fill=white] (vert_0) at (0.000000,0.000000) {$ L_1 $};
		\node [draw, rectangle, color=red, fill=white] (vert_1) at (1.000000,0.000000) {$ C_1 $};
		\node [draw, rectangle, color=red, fill=white] (vert_2) at (2.000000,0.000000) {$ C_3 $};
		\node [draw, rectangle, color=red, fill=white] (vert_3) at (3.000000,0.000000) {$ C_2 $};
		\node [draw, rectangle, color=red, fill=white] (vert_4) at (4.000000,0.000000) {$ C_2' $};
		\node [draw, rectangle, color=red, fill=white] (vert_5) at (5.000000,0.000000) {$ C_3' $};
		\node [draw, rectangle, color=red, fill=white] (vert_6) at (6.000000,0.000000) {$ C_1' $};
		
		\end{scope}[]
		\end{tikzpicture}
	\end{center}
\end{figure}

Relation:

\begin{equation*}\begin{split}
\Gamma_{C_2}  = \Gamma_{C_2'}
\end{split}
\end{equation*}
\subsubsection{Vertex Number 2}

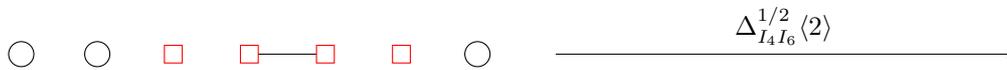
\begin{figure}[H]
	\begin{center}
		\begin{tikzpicture}
		\begin{scope}[scale=1 , xshift=-100]
		\draw [] (3.000000, 0.000000) -- (4.000000, 0.000000) ;
		\node [draw, circle, color=black, fill=white] (vert_0) at (0.000000,0.000000) {};
		\node [draw, circle, color=black, fill=white] (vert_1) at (1.000000,0.000000) {};
		\node [draw, rectangle, color=red, fill=white] (vert_2) at (2.000000,0.000000) {};
		\node [draw, rectangle, color=red, fill=white] (vert_3) at (3.000000,0.000000) {};
		\node [draw, rectangle, color=red, fill=white] (vert_4) at (4.000000,0.000000) {};
		\node [draw, rectangle, color=red, fill=white] (vert_5) at (5.000000,0.000000) {};
		\node [draw, circle, color=black, fill=white] (vert_6) at (6.000000,0.000000) {};
		
		\end{scope}[]
		\begin{scope}[scale=1 , xshift=100]
		\draw [->] (0.000000, 0.000000) -- node [above, midway] {$ \Delta_{I_{4}I_{6}}^{1/2}\langle 2 \rangle $}  (6.000000, 0.000000) ;
		
		\end{scope}[]
		\end{tikzpicture}
	\end{center}
\end{figure}

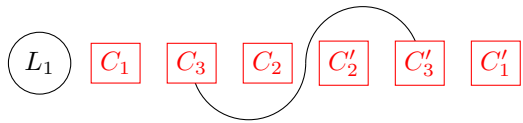
\begin{figure}[H]
	\begin{center}
		\begin{tikzpicture}
		\begin{scope}[]
		\draw [] (2.000000, 0.000000) arc (180.000000:360.000000:0.750000);
		\draw [] (3.500000, 0.000000) arc (180.000000:0.000000:0.750000);
		\node [draw, circle, color=black, fill=white] (vert_0) at (0.000000,0.000000) {$ L_1 $};
		\node [draw, rectangle, color=red, fill=white] (vert_1) at (1.000000,0.000000) {$ C_1 $};
		\node [draw, rectangle, color=red, fill=white] (vert_2) at (2.000000,0.000000) {$ C_3 $};
		\node [draw, rectangle, color=red, fill=white] (vert_3) at (3.000000,0.000000) {$ C_2 $};
		\node [draw, rectangle, color=red, fill=white] (vert_4) at (4.000000,0.000000) {$ C_2' $};
		\node [draw, rectangle, color=red, fill=white] (vert_5) at (5.000000,0.000000) {$ C_3' $};
		\node [draw, rectangle, color=red, fill=white] (vert_6) at (6.000000,0.000000) {$ C_1' $};
		
		\end{scope}[]
		\end{tikzpicture}
	\end{center}
\end{figure}

Relation:

\begin{equation*}\begin{split}
\Gamma_{C_2'}\Gamma_{C_3}\Gamma_{C_2'}^{-1}  = \Gamma_{C_3'}
\end{split}
\end{equation*}
\subsubsection{Vertex Number 3}

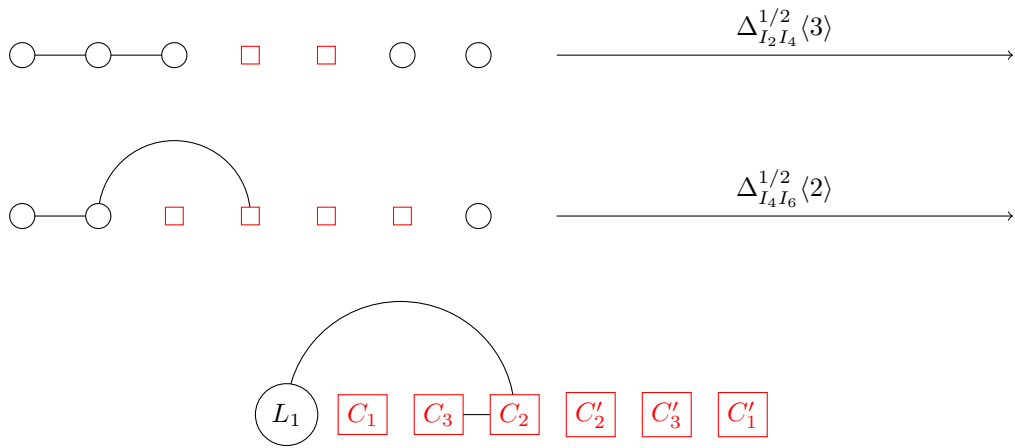
\begin{figure}[H]
	\begin{center}
		\begin{tikzpicture}
		\begin{scope}[scale=1 , xshift=-100]
		\draw [] (0.000000, 0.000000) -- (1.000000, 0.000000) ;
		\draw [] (1.000000, 0.000000) -- (2.000000, 0.000000) ;
		\node [draw, circle, color=black, fill=white] (vert_0) at (0.000000,0.000000) {};
		\node [draw, circle, color=black, fill=white] (vert_1) at (1.000000,0.000000) {};
		\node [draw, circle, color=black, fill=white] (vert_2) at (2.000000,0.000000) {};
		\node [draw, rectangle, color=red, fill=white] (vert_3) at (3.000000,0.000000) {};
		\node [draw, rectangle, color=red, fill=white] (vert_4) at (4.000000,0.000000) {};
		\node [draw, circle, color=black, fill=white] (vert_5) at (5.000000,0.000000) {};
		\node [draw, circle, color=black, fill=white] (vert_6) at (6.000000,0.000000) {};
		
		\end{scope}[]
		\begin{scope}[scale=1 , xshift=100]
		\draw [->] (0.000000, 0.000000) -- node [above, midway] {$ \Delta_{I_{2}I_{4}}^{1/2}\langle 3 \rangle $}  (6.000000, 0.000000) ;
		
		\end{scope}[]
		\end{tikzpicture}
	\end{center}
\end{figure}

\begin{figure}[H]
	\begin{center}
		\begin{tikzpicture}
		\begin{scope}[scale=1 , xshift=-100]
		\draw [] (0.000000, 0.000000) -- (1.000000, 0.000000) ;
		\draw [] (1.000000, 0.000000) arc (180.000000:0.000000:1.000000);
		\node [draw, circle, color=black, fill=white] (vert_0) at (0.000000,0.000000) {};
		\node [draw, circle, color=black, fill=white] (vert_1) at (1.000000,0.000000) {};
		\node [draw, rectangle, color=red, fill=white] (vert_2) at (2.000000,0.000000) {};
		\node [draw, rectangle, color=red, fill=white] (vert_3) at (3.000000,0.000000) {};
		\node [draw, rectangle, color=red, fill=white] (vert_4) at (4.000000,0.000000) {};
		\node [draw, rectangle, color=red, fill=white] (vert_5) at (5.000000,0.000000) {};
		\node [draw, circle, color=black, fill=white] (vert_6) at (6.000000,0.000000) {};
		
		\end{scope}[]
		\begin{scope}[scale=1 , xshift=100]
		\draw [->] (0.000000, 0.000000) -- node [above, midway] {$ \Delta_{I_{4}I_{6}}^{1/2}\langle 2 \rangle $}  (6.000000, 0.000000) ;
		
		\end{scope}[]
		\end{tikzpicture}
	\end{center}
\end{figure}

\begin{figure}[H]
	\begin{center}
		\begin{tikzpicture}
		\begin{scope}[]
		\draw [] (0.000000, 0.000000) arc (180.000000:0.000000:1.500000);
		\draw [] (3.000000, 0.000000) -- (2.000000, 0.000000) ;
		\node [draw, circle, color=black, fill=white] (vert_0) at (0.000000,0.000000) {$ L_1 $};
		\node [draw, rectangle, color=red, fill=white] (vert_1) at (1.000000,0.000000) {$ C_1 $};
		\node [draw, rectangle, color=red, fill=white] (vert_2) at (2.000000,0.000000) {$ C_3 $};
		\node [draw, rectangle, color=red, fill=white] (vert_3) at (3.000000,0.000000) {$ C_2 $};
		\node [draw, rectangle, color=red, fill=white] (vert_4) at (4.000000,0.000000) {$ C_2' $};
		\node [draw, rectangle, color=red, fill=white] (vert_5) at (5.000000,0.000000) {$ C_3' $};
		\node [draw, rectangle, color=red, fill=white] (vert_6) at (6.000000,0.000000) {$ C_1' $};
		
		\end{scope}[]
		\end{tikzpicture}
	\end{center}
\end{figure}

Relation:

\begin{equation*}\begin{split}
\left( \Gamma_{C_3}  \right) \left( \Gamma_{C_2}  \right) \left( \Gamma_{C_3}\Gamma_{C_1}\Gamma_{L_1}\Gamma_{C_1}^{-1}\Gamma_{C_3}^{-1}  \right)& = \\
\left( \Gamma_{C_2}  \right) \left( \Gamma_{C_3}\Gamma_{C_1}\Gamma_{L_1}\Gamma_{C_1}^{-1}\Gamma_{C_3}^{-1}  \right) \left( \Gamma_{C_3}  \right)& = \\
\left( \Gamma_{C_3}\Gamma_{C_1}\Gamma_{L_1}\Gamma_{C_1}^{-1}\Gamma_{C_3}^{-1}  \right) \left( \Gamma_{C_3}  \right) \left( \Gamma_{C_2}  \right)
\end{split}
\end{equation*}
\subsubsection{Vertex Number 4}

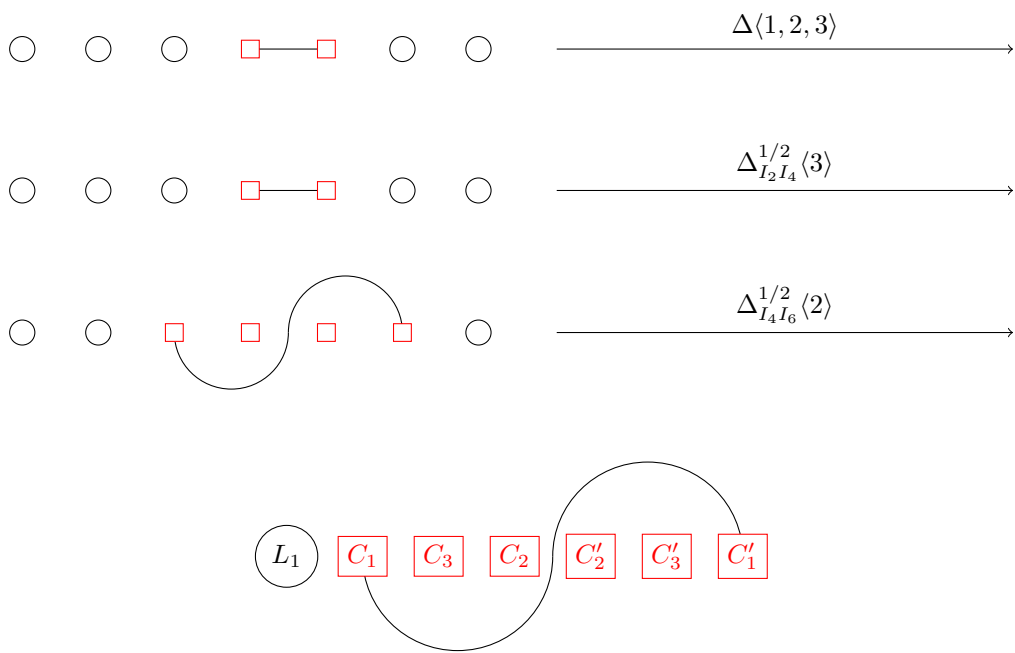
\begin{figure}[H]
	\begin{center}

	\end{center}
\end{figure}

Relation:

\begin{equation*}\begin{split}
\Gamma_{C_3'}\Gamma_{C_2'}\Gamma_{C_1}\Gamma_{C_2'}^{-1}\Gamma_{C_3'}^{-1}  = \Gamma_{C_1'}
\end{split}
\end{equation*}
\subsubsection{Vertex Number 5}

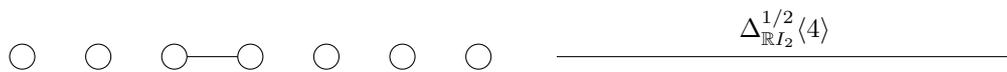
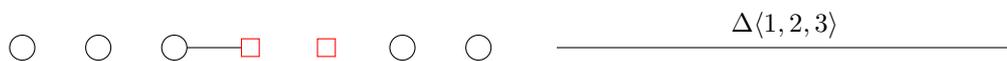
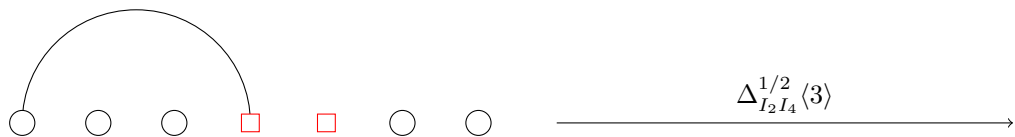
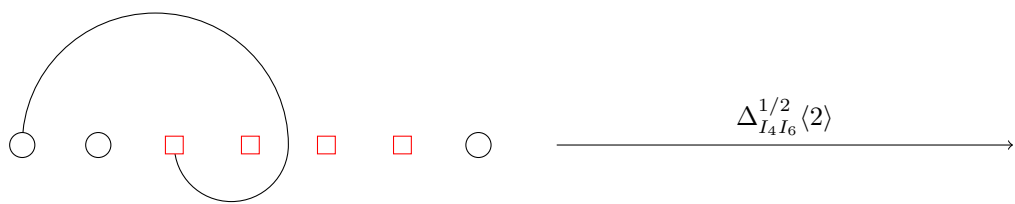
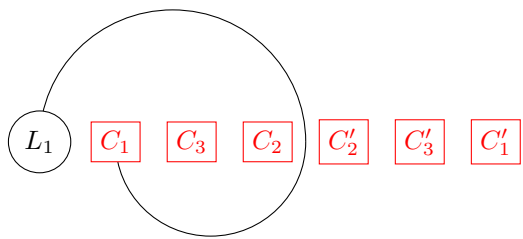
\begin{figure}[H]
	\begin{center}
		%
	\end{center}
\end{figure}

Relation:

\begin{equation*}\begin{split}
\left[\Gamma_{C_2}\Gamma_{C_3}\Gamma_{C_1}\Gamma_{L_1}\Gamma_{C_1}^{-1}\Gamma_{C_3}^{-1}\Gamma_{C_2}^{-1} , \Gamma_{C_1} \right]=e
\end{split}
\end{equation*}
\subsubsection{Vertex Number 6}

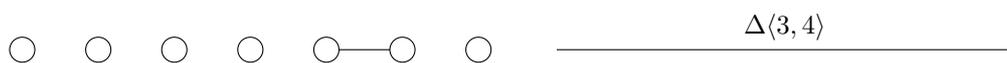
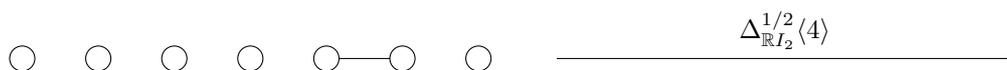
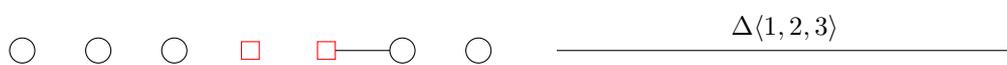
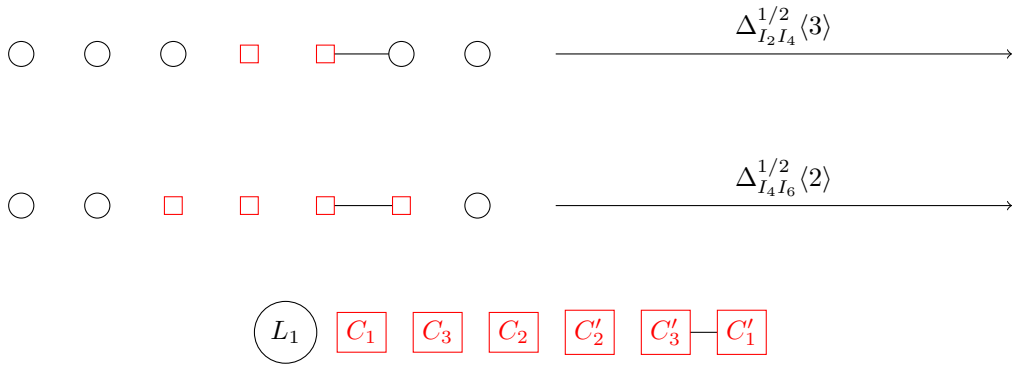
\begin{figure}[H]
	\begin{center}
		%
	\end{center}
\end{figure}

Relation:

\begin{equation*}\begin{split}
\left\{\Gamma_{C_1'}  , \Gamma_{C_3'} \right\}=e
\end{split}
\end{equation*}
\subsubsection{Vertex Number 7}

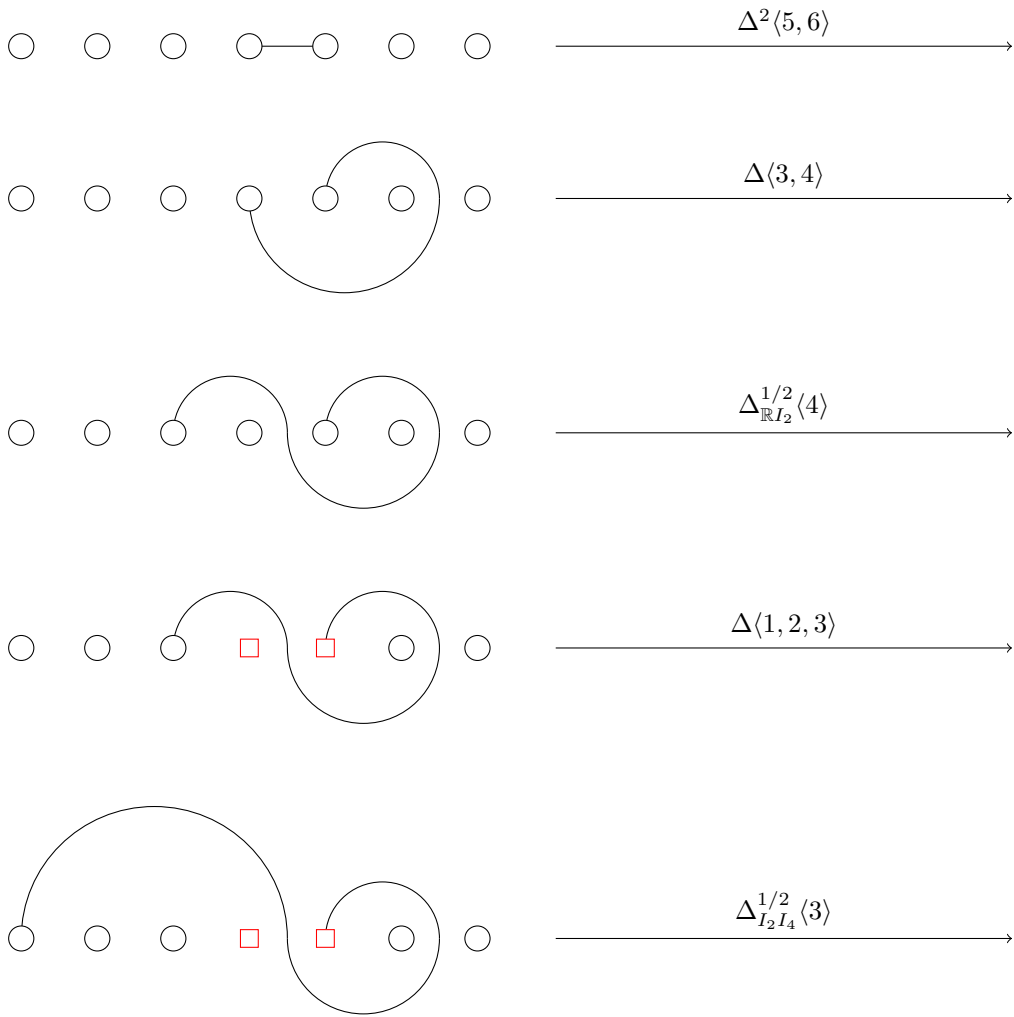
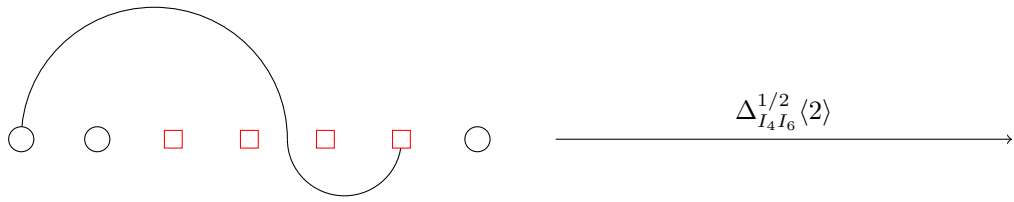
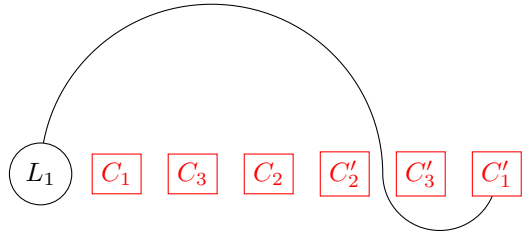
\begin{figure}[H]
	\begin{center}
		%
	\end{center}
\end{figure}

Relation:

\begin{equation*}\begin{split}
\left[\Gamma_{C_2'}\Gamma_{C_2}\Gamma_{C_3}\Gamma_{C_1}\Gamma_{L_1}\Gamma_{C_1}^{-1}\Gamma_{C_3}^{-1}\Gamma_{C_2}^{-1}\Gamma_{C_2'}^{-1} , \Gamma_{C_1'} \right]=e
\end{split}
\end{equation*}
\subsubsection{Vertex Number 8}

\begin{figure}[H]
	\begin{center}
		%
	\end{center}
\end{figure}

Relation:

\begin{equation*}\begin{split}
\left( \Gamma_{C_1'}^{-1}\Gamma_{C_2'}\Gamma_{C_1'}  \right) \left( \Gamma_{C_3'}  \right) \left( \Gamma_{C_2'}\Gamma_{C_2}\Gamma_{C_3}\Gamma_{C_1}\Gamma_{L_1}\Gamma_{C_1}^{-1}\Gamma_{C_3}^{-1}\Gamma_{C_2}^{-1}\Gamma_{C_2'}^{-1}  \right)& = \\
\left( \Gamma_{C_3'}  \right) \left( \Gamma_{C_2'}\Gamma_{C_2}\Gamma_{C_3}\Gamma_{C_1}\Gamma_{L_1}\Gamma_{C_1}^{-1}\Gamma_{C_3}^{-1}\Gamma_{C_2}^{-1}\Gamma_{C_2'}^{-1}  \right) \left( \Gamma_{C_1'}^{-1}\Gamma_{C_2'}\Gamma_{C_1'}  \right)& = \\
\left( \Gamma_{C_2'}\Gamma_{C_2}\Gamma_{C_3}\Gamma_{C_1}\Gamma_{L_1}\Gamma_{C_1}^{-1}\Gamma_{C_3}^{-1}\Gamma_{C_2}^{-1}\Gamma_{C_2'}^{-1}  \right) \left( \Gamma_{C_1'}^{-1}\Gamma_{C_2'}\Gamma_{C_1'}  \right) \left( \Gamma_{C_3'}  \right)
\end{split}
\end{equation*}
\subsubsection{Vertex Number 9}

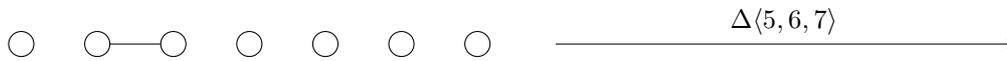
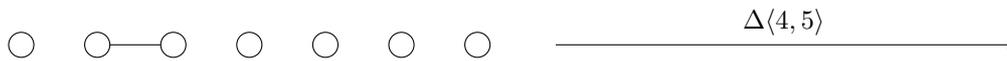
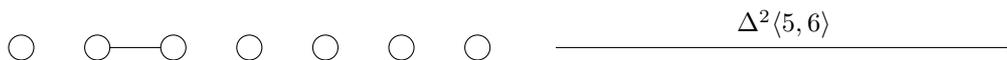
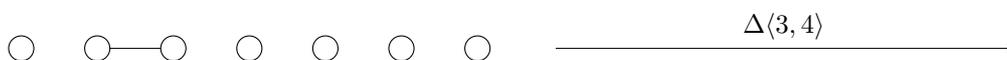
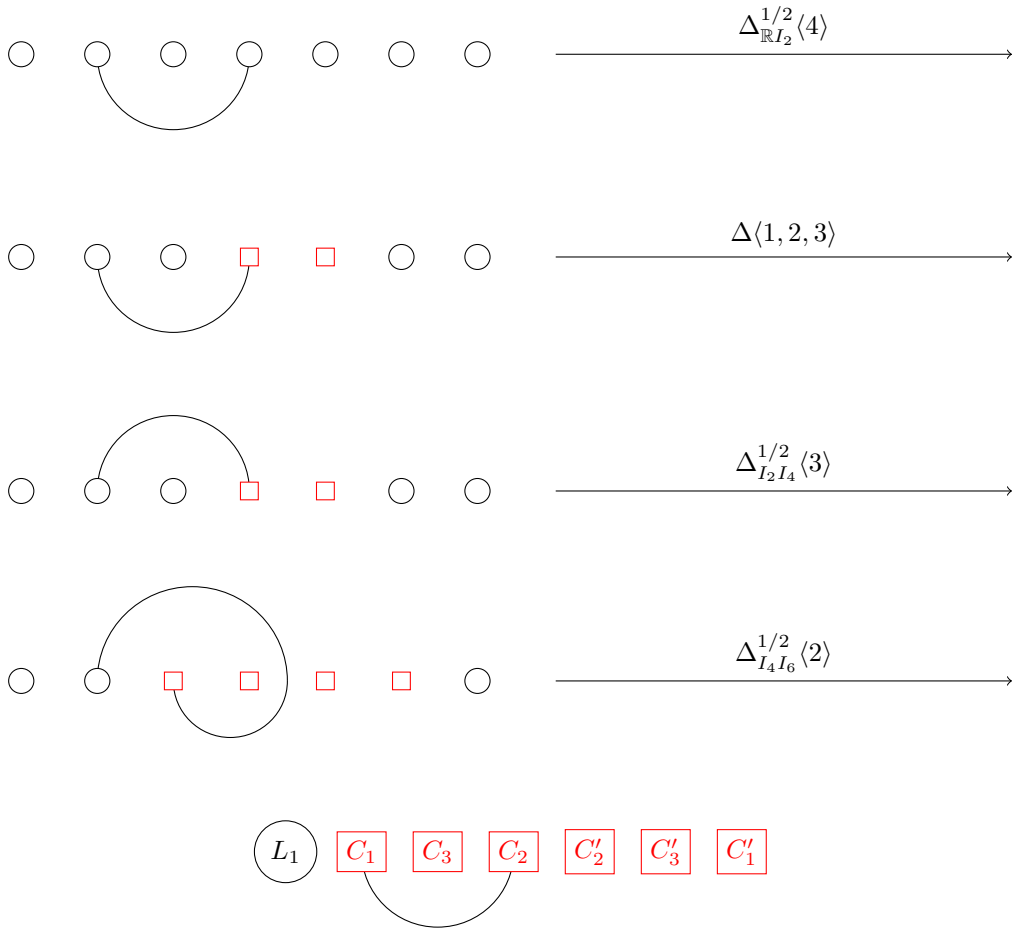
\begin{figure}[H]
	\begin{center}

	\end{center}
\end{figure}

Relation:

\begin{equation*}\begin{split}
\left\{\Gamma_{C_2}  , \Gamma_{C_1} \right\}=e
\end{split}
\end{equation*}
\subsubsection{Vertex Number 10}

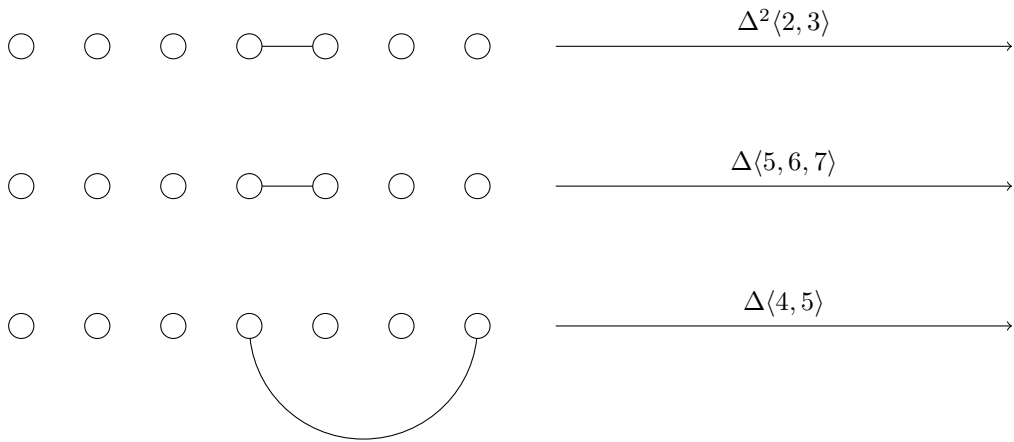
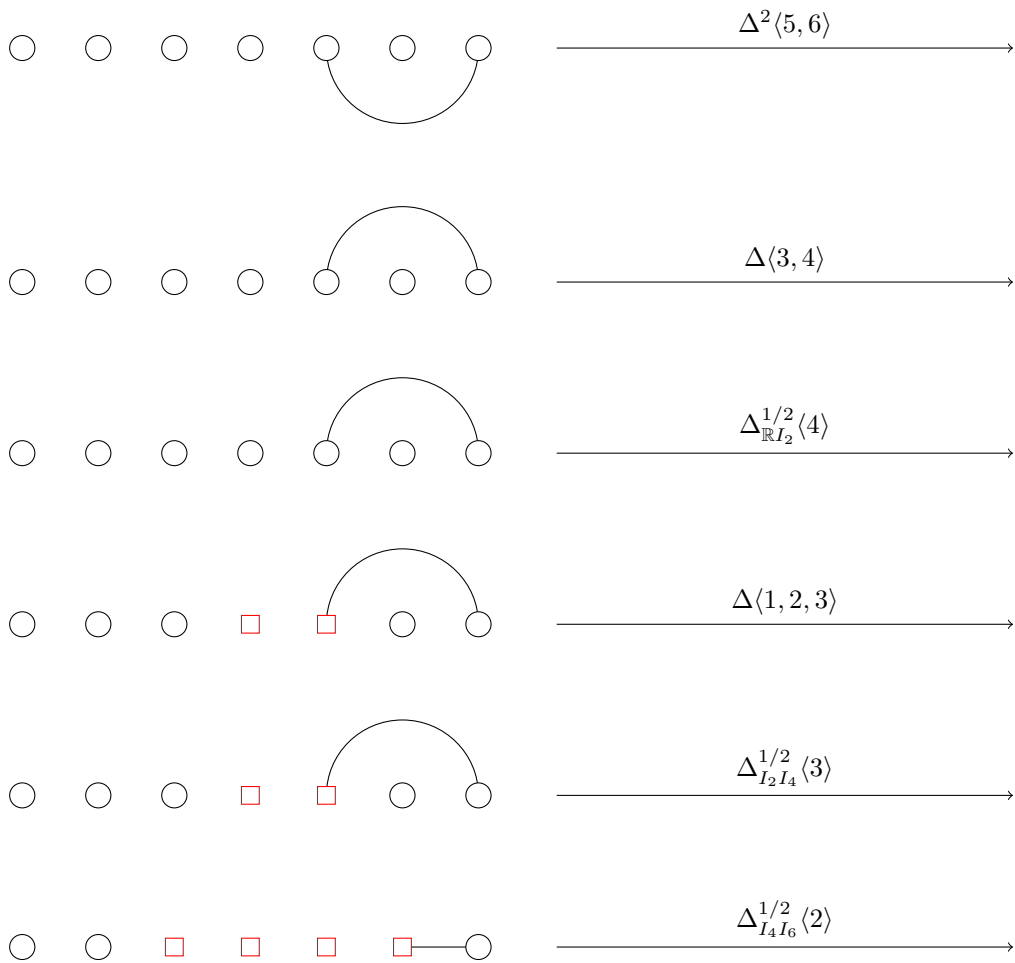
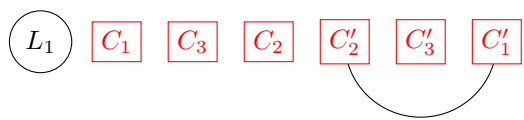
\begin{figure}[H]
	\begin{center}
		%
	\end{center}
\end{figure}

Relation:

\begin{equation*}\begin{split}
\left\{\Gamma_{C_1'}  , \Gamma_{C_2'} \right\}=e
\end{split}
\end{equation*}
\subsubsection{Vertex Number 11}

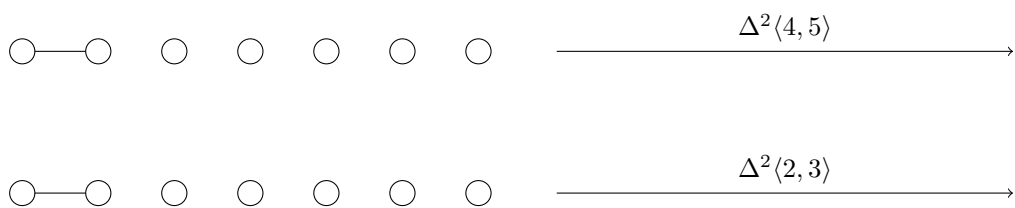
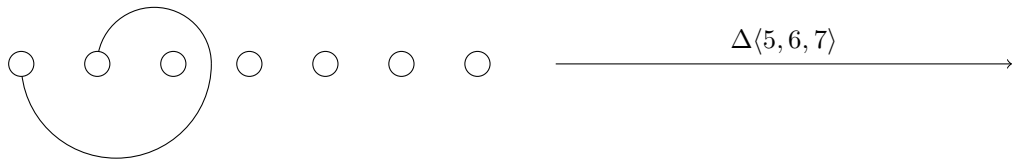
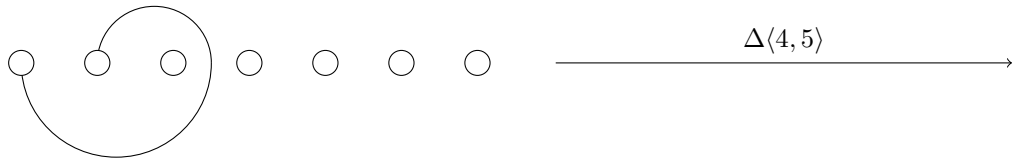
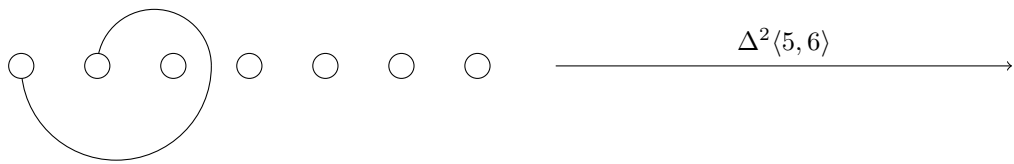
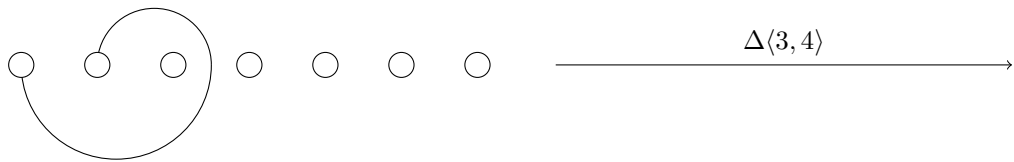
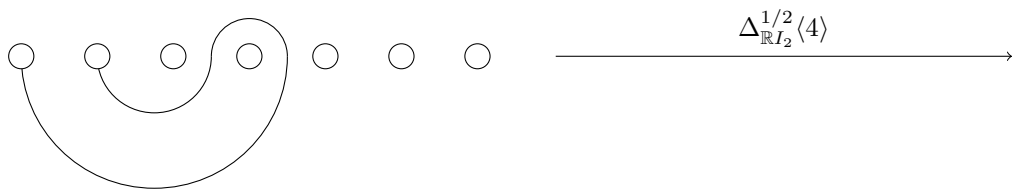
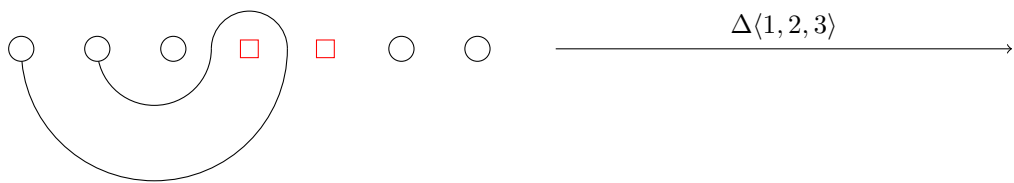
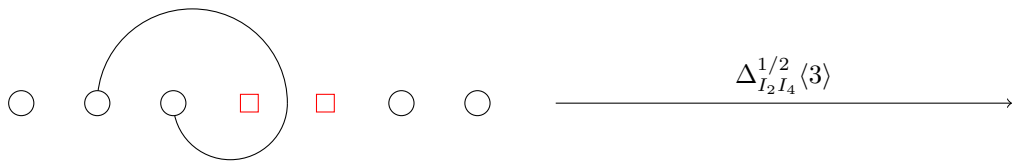
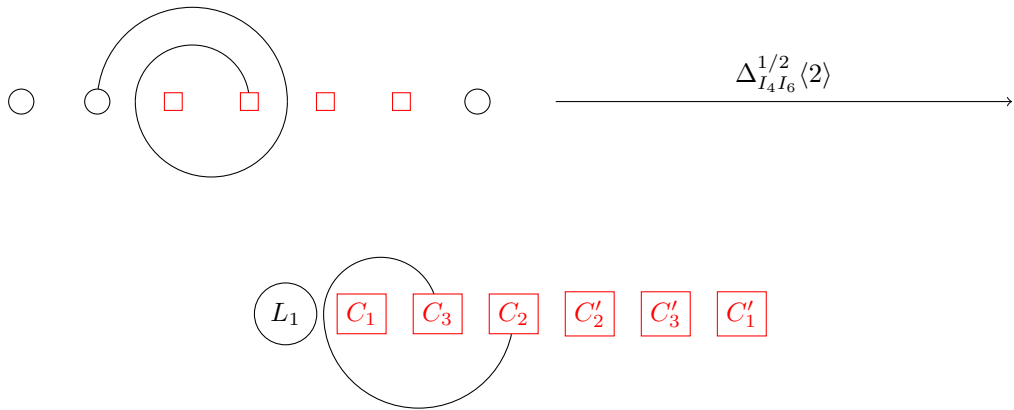
\begin{figure}[H]
	\begin{center}
		%
	\end{center}
\end{figure}

Relation:

\begin{equation*}\begin{split}
\left[\Gamma_{C_1}^{-1}\Gamma_{C_3}\Gamma_{C_1} , \Gamma_{C_2} \right]=e
\end{split}
\end{equation*}
\subsubsection{Vertex Number 12}

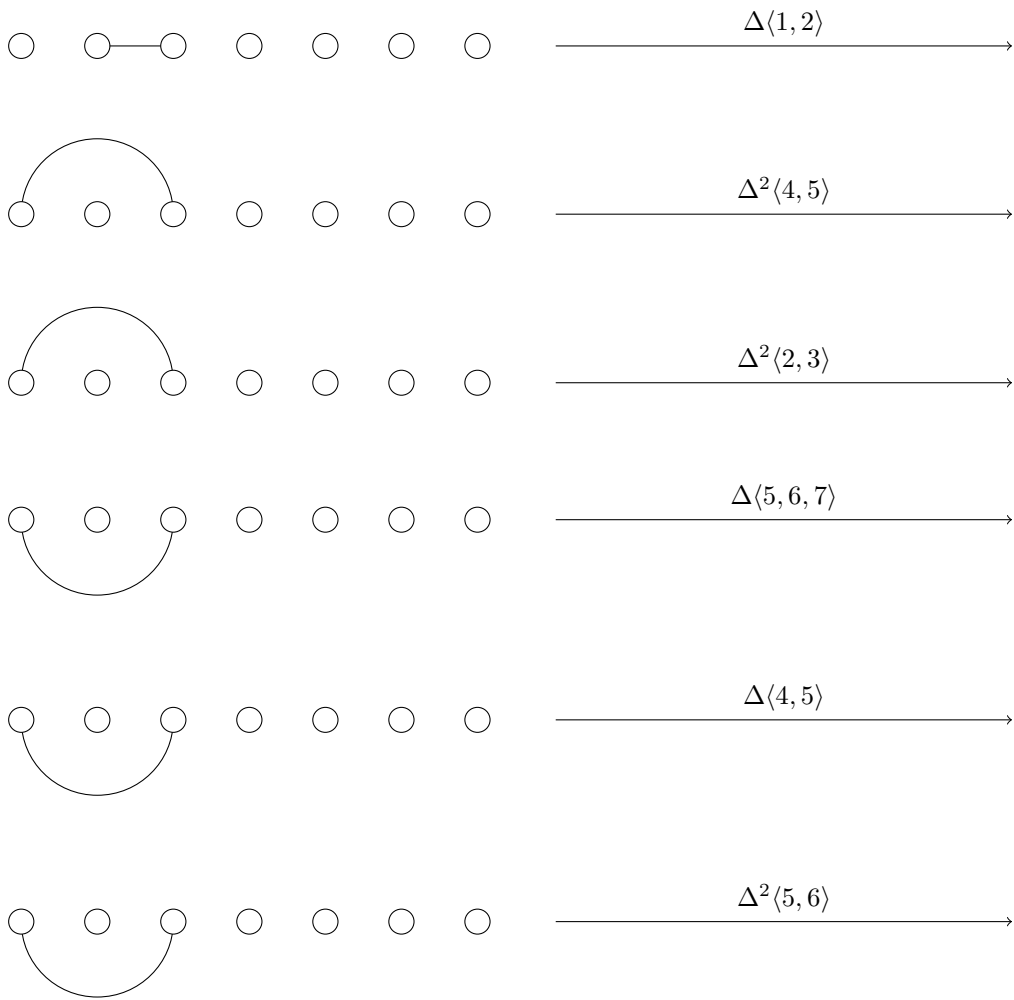
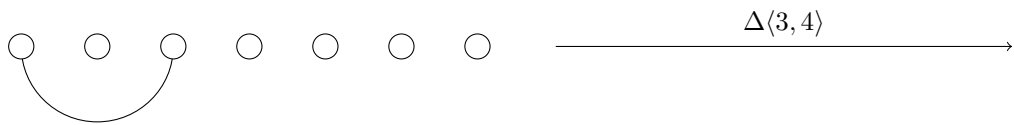
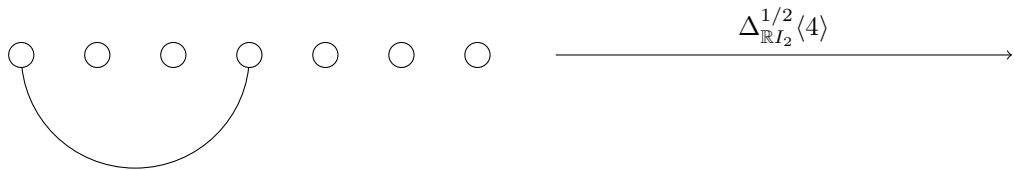
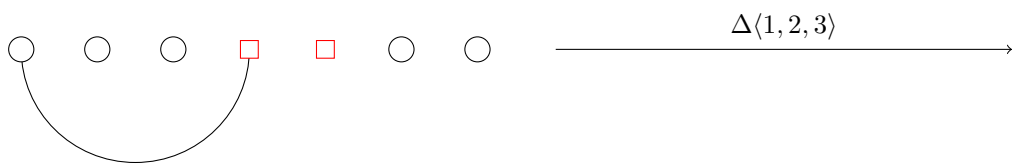
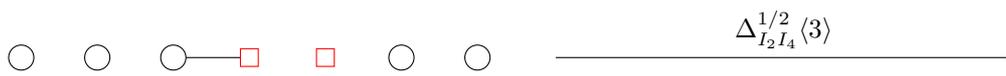
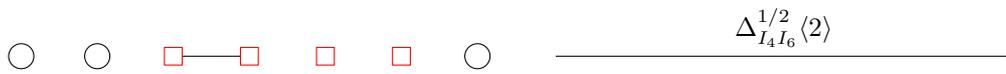
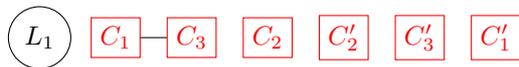
\begin{figure}[H]
	\begin{center}
		%
	\end{center}
\end{figure}

Relation:

\begin{equation*}\begin{split}
\left\{\Gamma_{C_3}  , \Gamma_{C_1} \right\}=e
\end{split}
\end{equation*}
\subsubsection{Vertex Number 13}

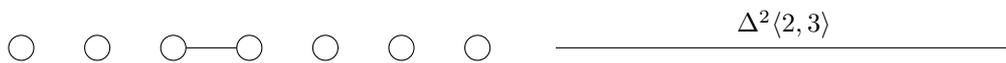
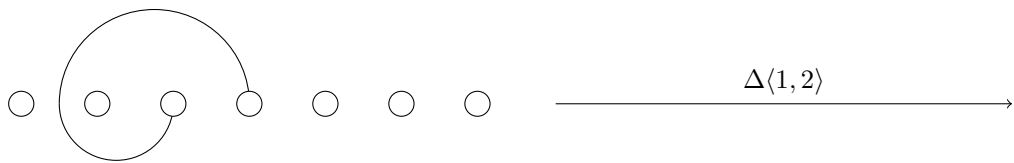
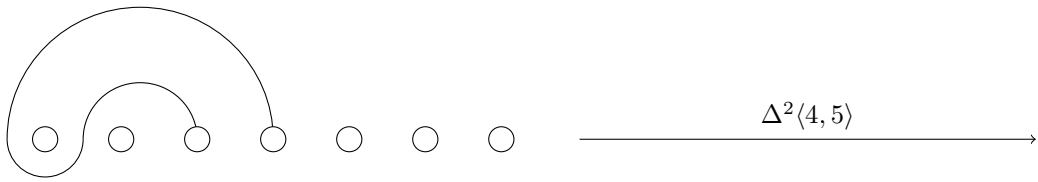
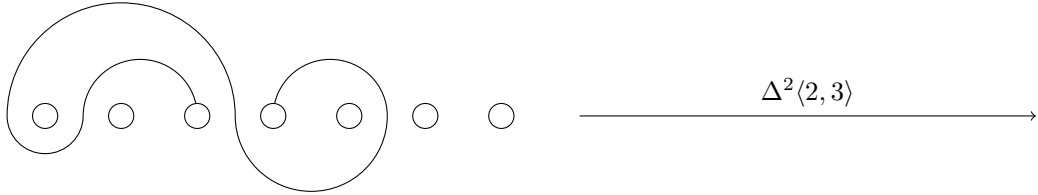
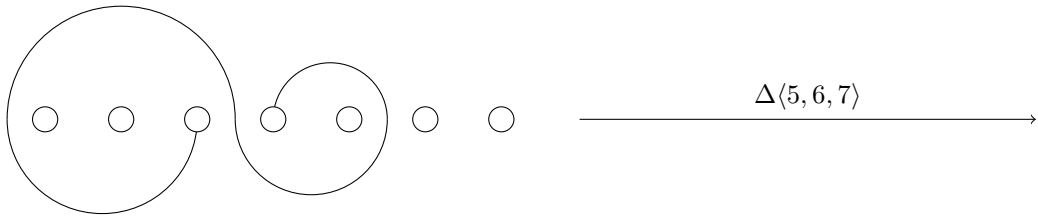
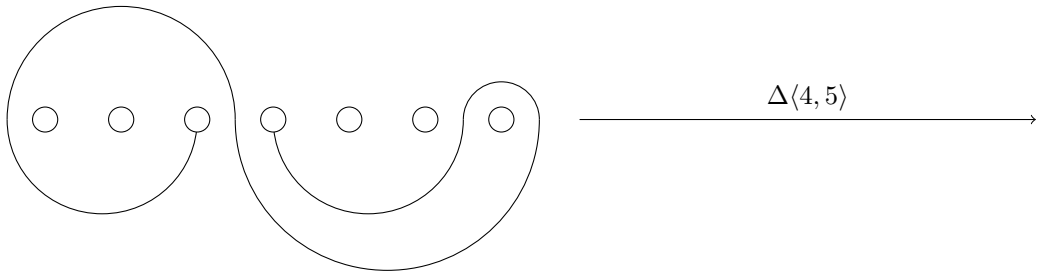
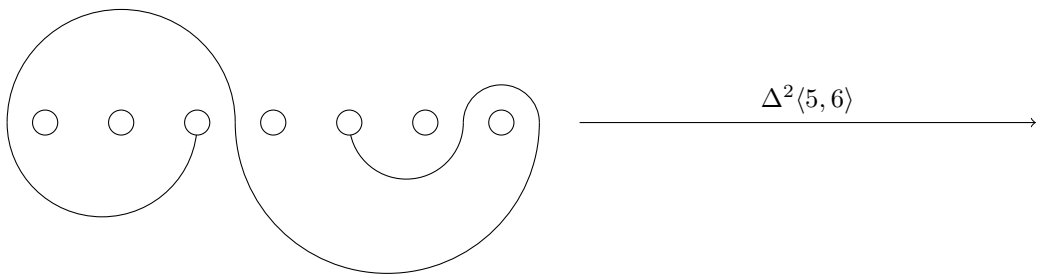
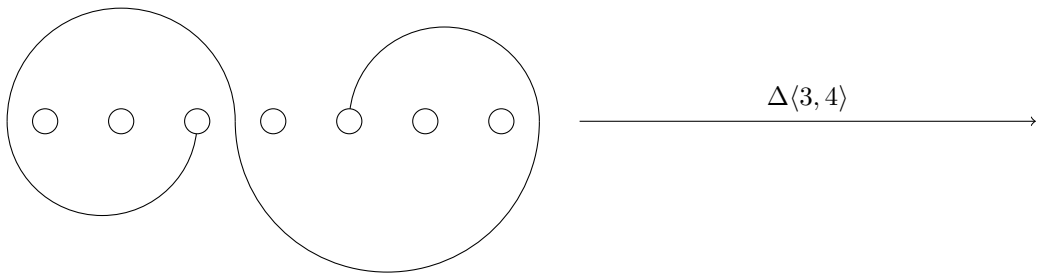
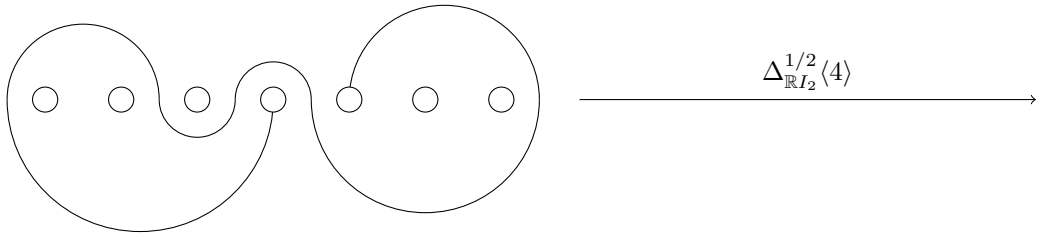
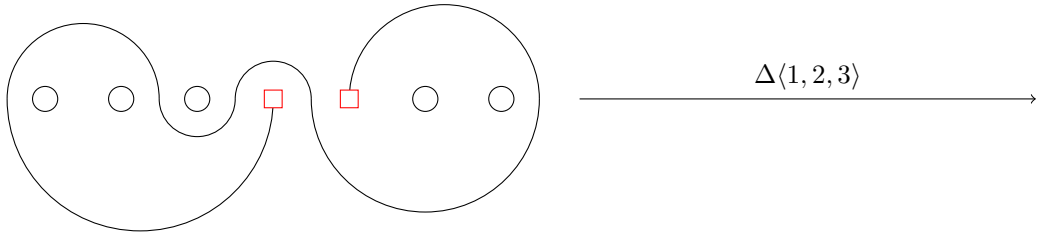
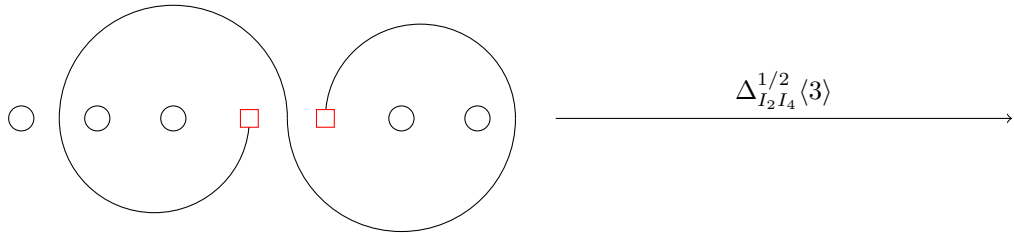
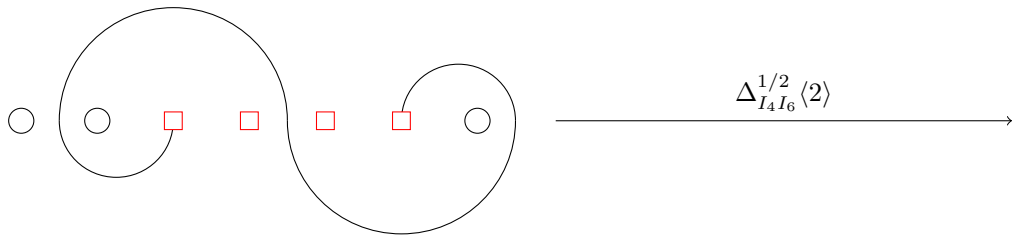
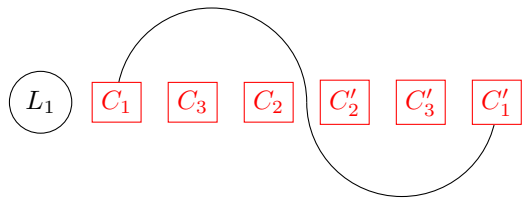
\begin{figure}[H]
	\begin{center}
		%
	\end{center}
\end{figure}

Relation:

\begin{equation*}\begin{split}
\Gamma_{C_2}\Gamma_{C_3}\Gamma_{C_1}\Gamma_{C_3}^{-1}\Gamma_{C_2}^{-1}  = \Gamma_{C_1'}
\end{split}
\end{equation*}
\subsubsection{Vertex Number 14}

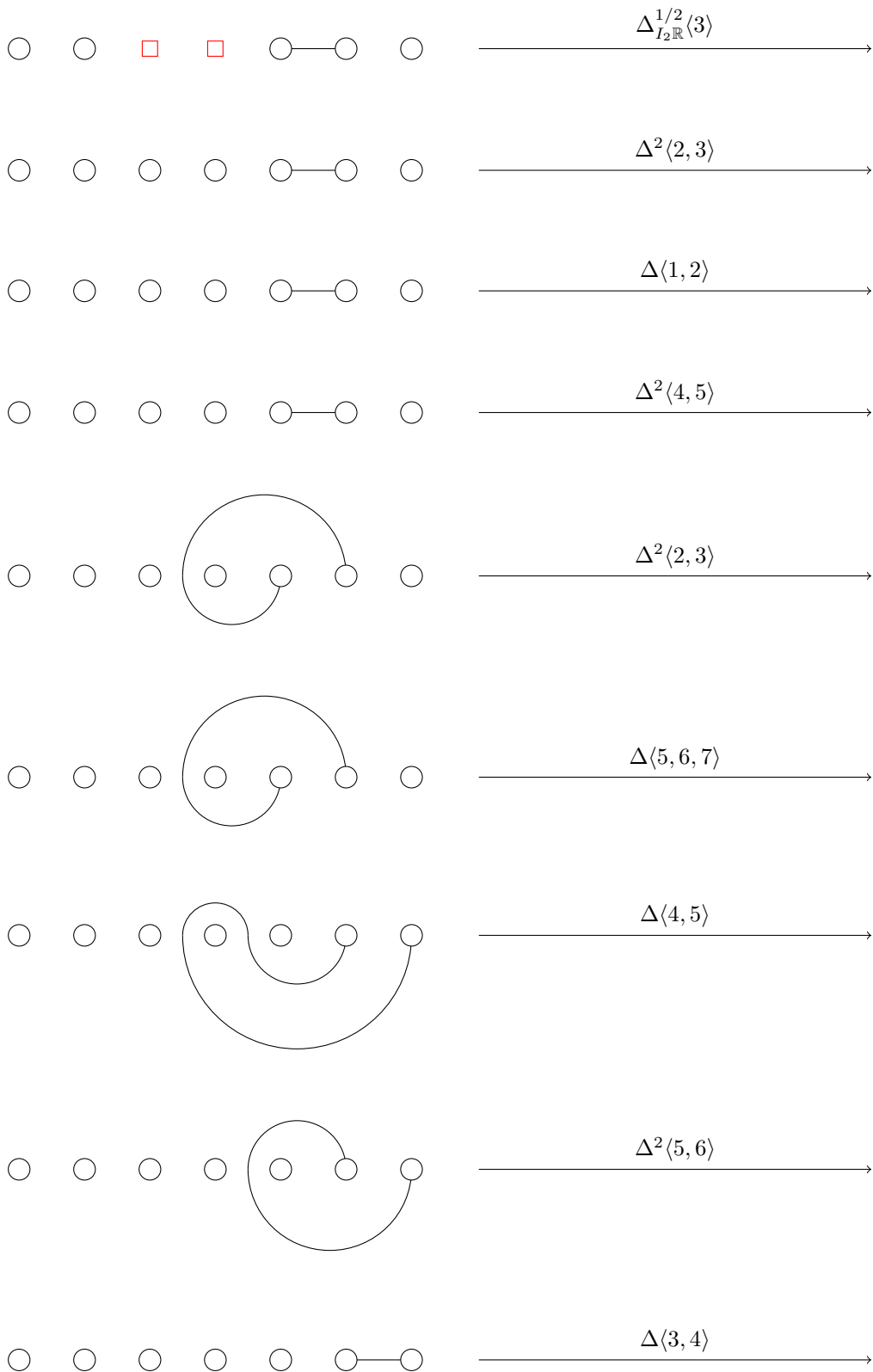
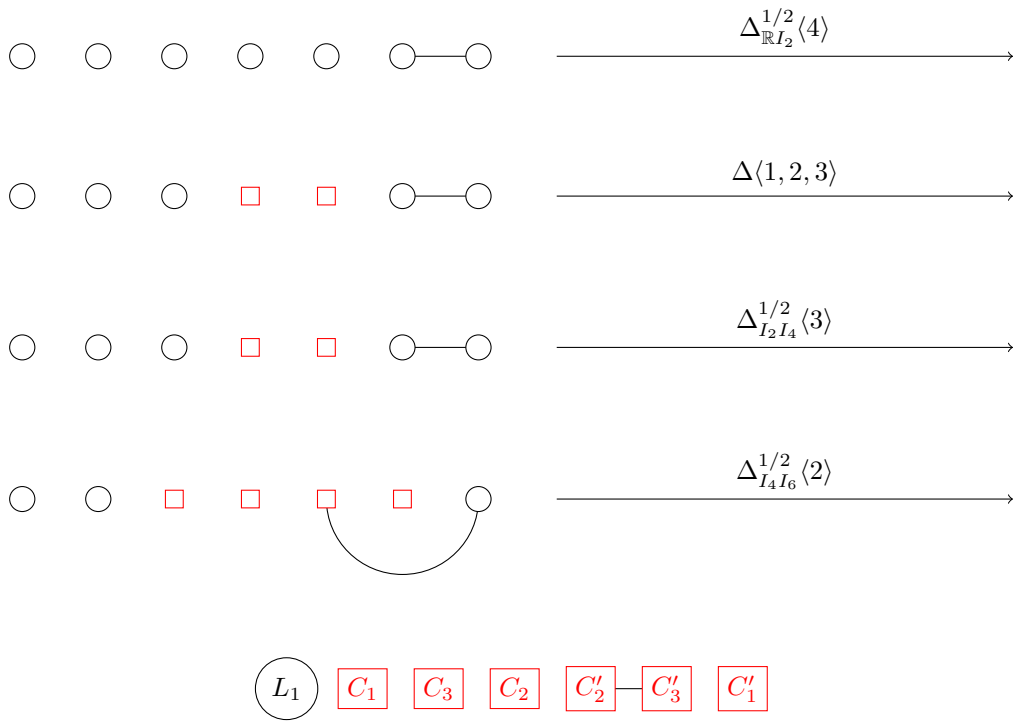
\begin{figure}[H]
	\begin{center}
		%
	\end{center}
\end{figure}

Relation:

\begin{equation*}\begin{split}
\left[\Gamma_{C_2'} , \Gamma_{C_3'} \right]=e
\end{split}
\end{equation*}
\subsubsection{Vertex Number 15}

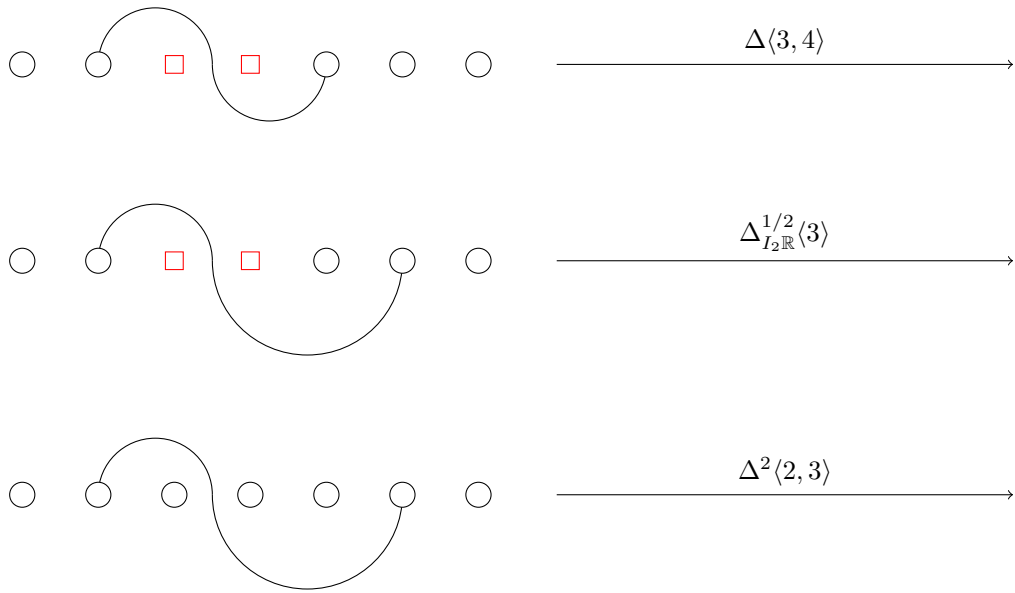
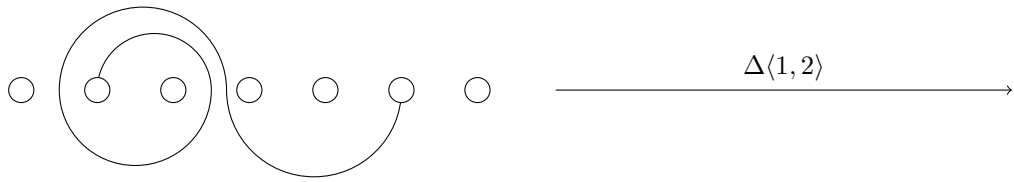
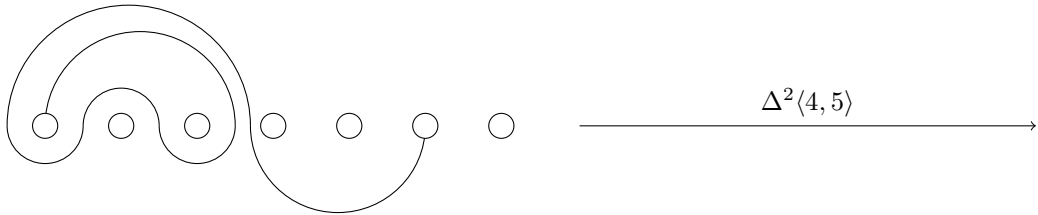
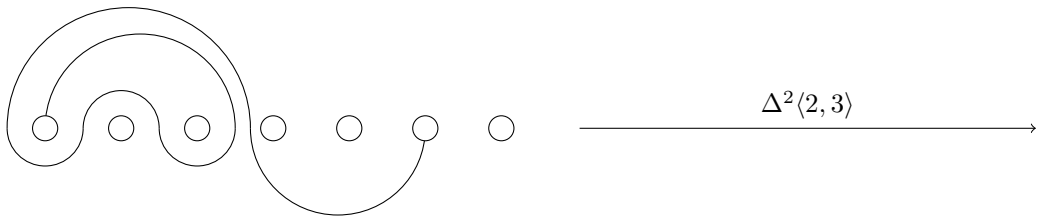
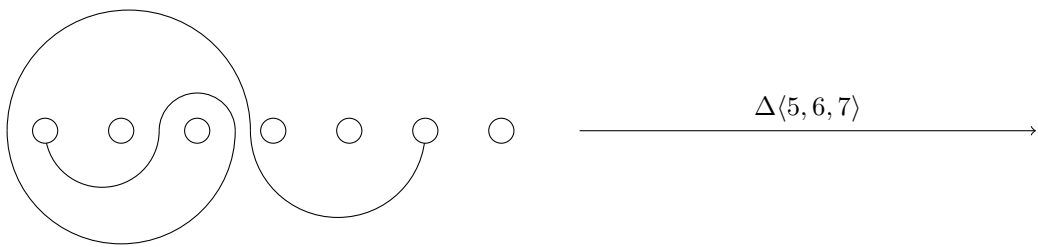
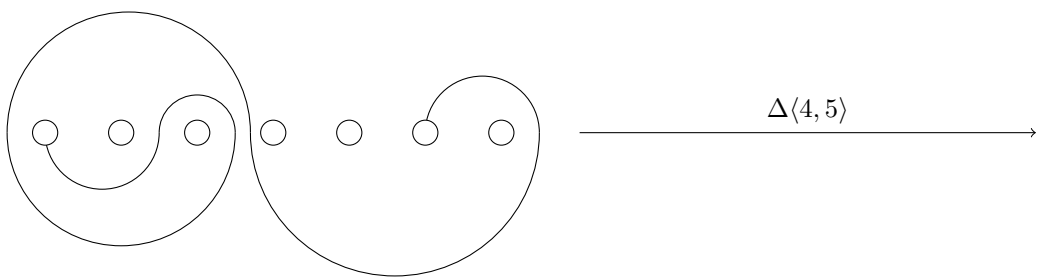
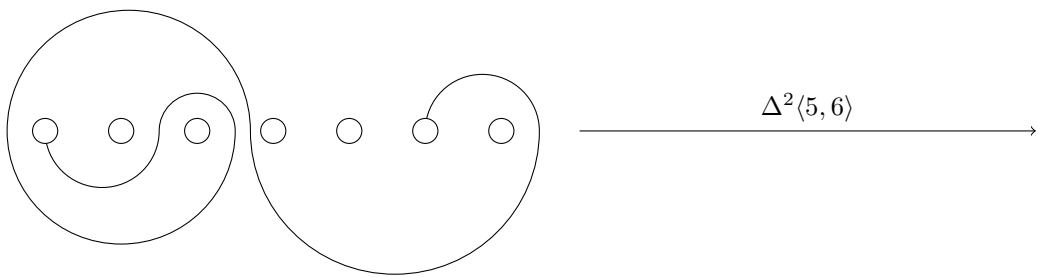
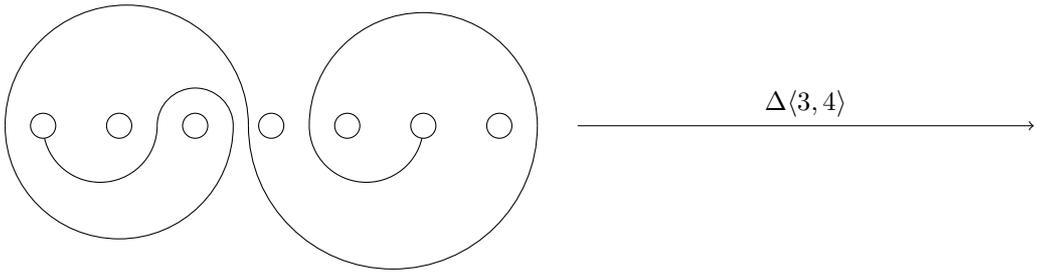
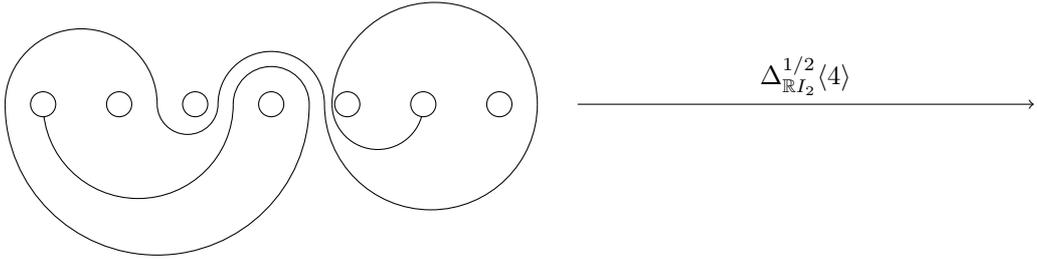
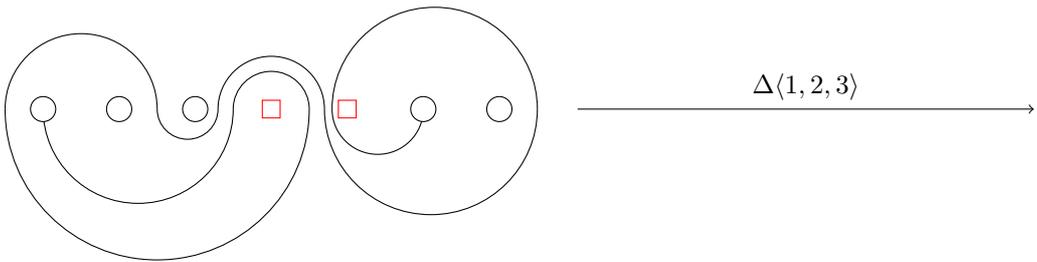
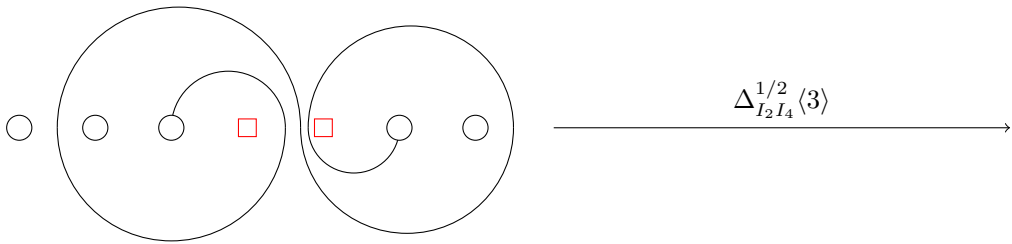
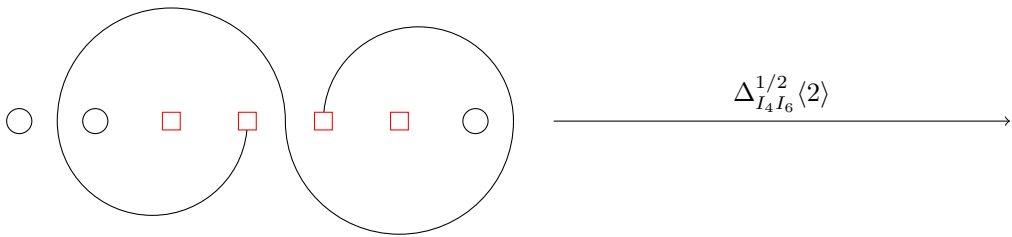
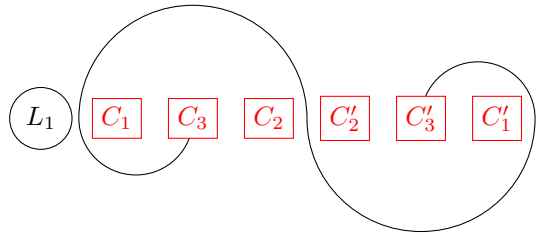
\begin{figure}[H]
	\begin{center}
		%
	\end{center}
\end{figure}

Relation:

\begin{equation*}\begin{split}
\Gamma_{C_1'}^{-1}\Gamma_{C_2}\Gamma_{C_3}\Gamma_{C_1}\Gamma_{C_3}\Gamma_{C_1}^{-1}\Gamma_{C_3}^{-1}\Gamma_{C_2}^{-1}\Gamma_{C_1'}  = \Gamma_{C_3'}
\end{split}
\end{equation*}
\subsubsection{Vertex Number 16}

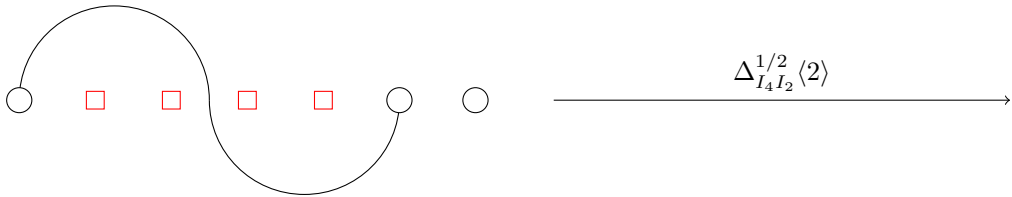
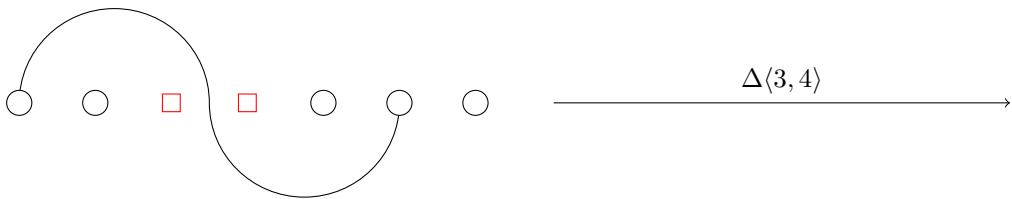
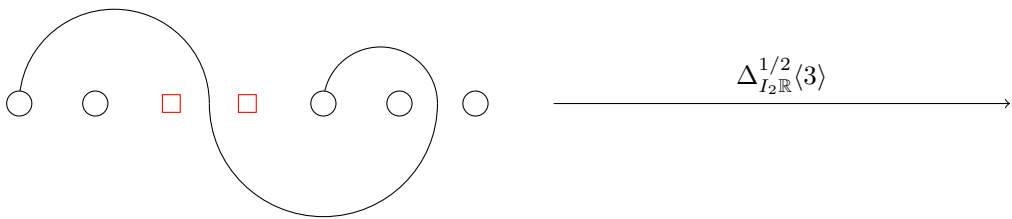
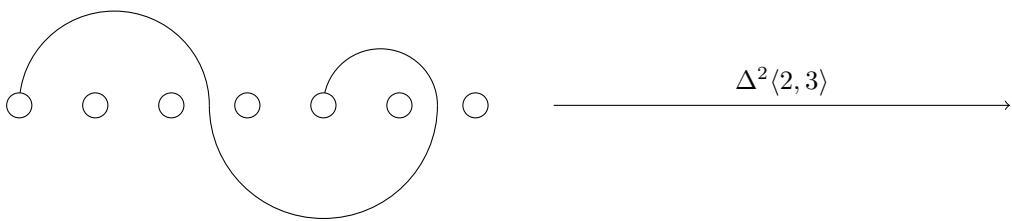
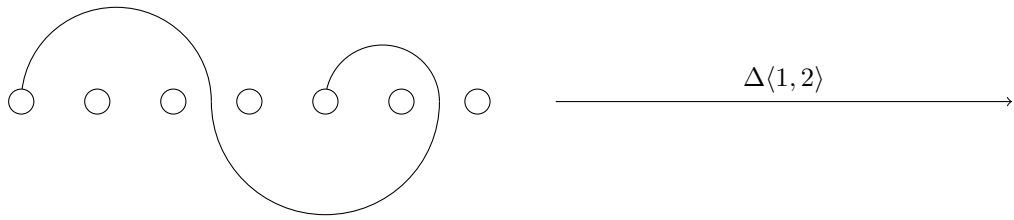
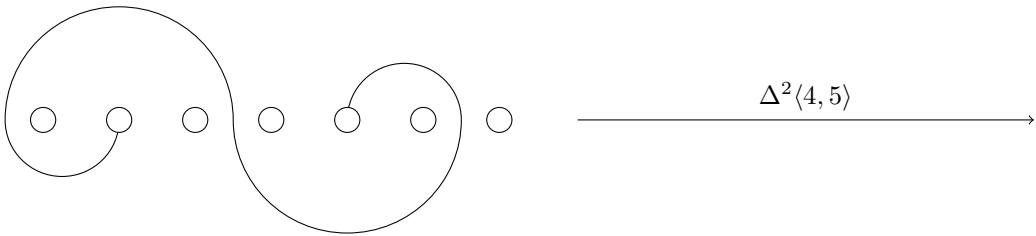
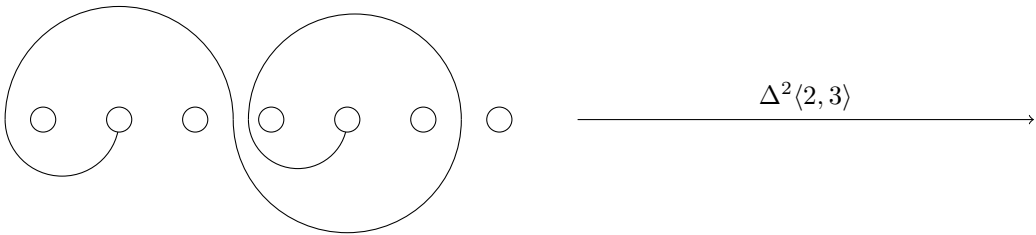
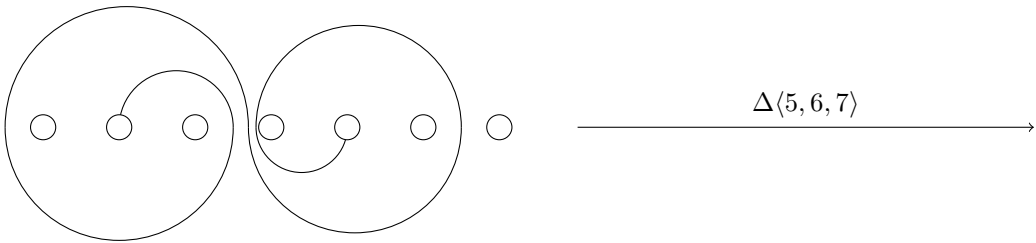
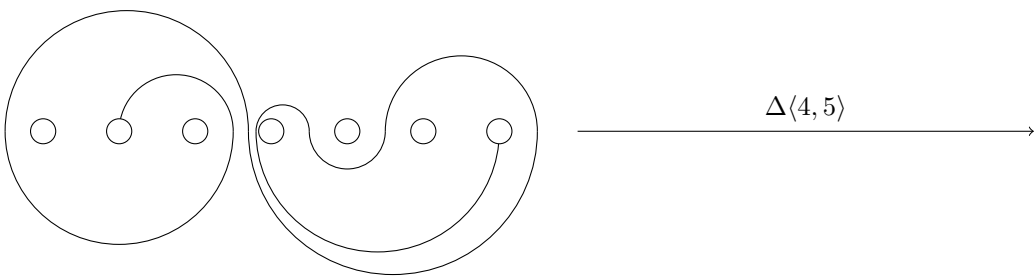
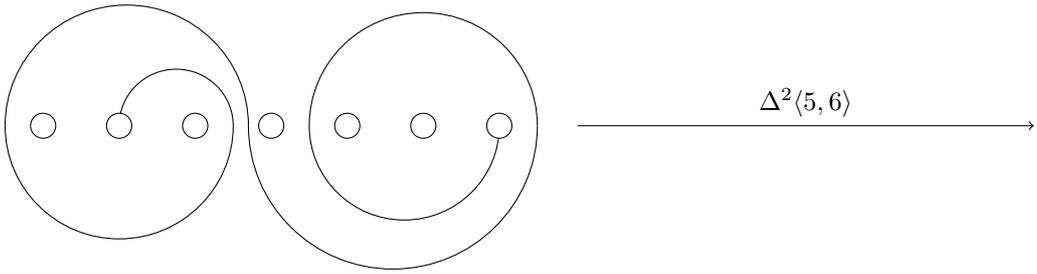
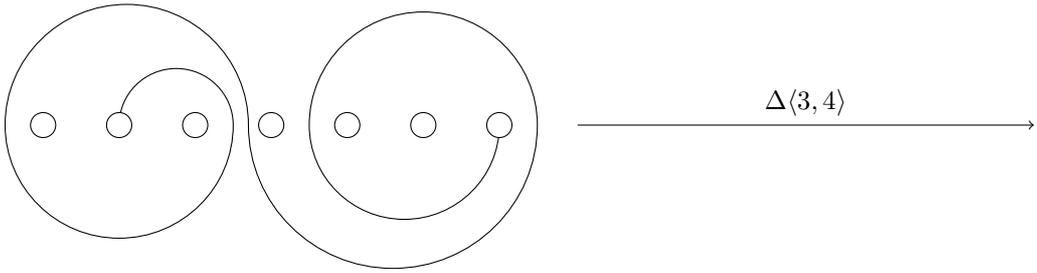
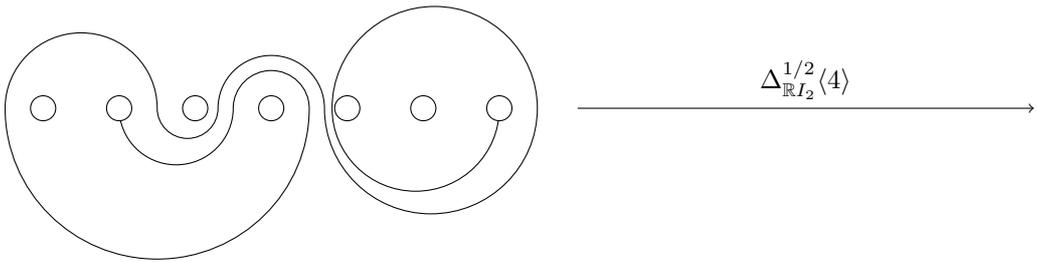
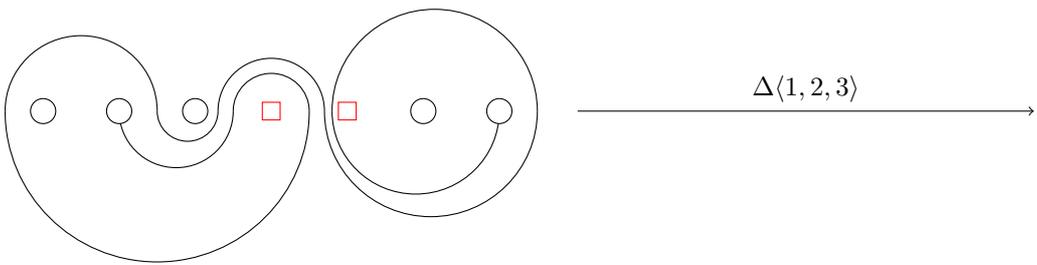
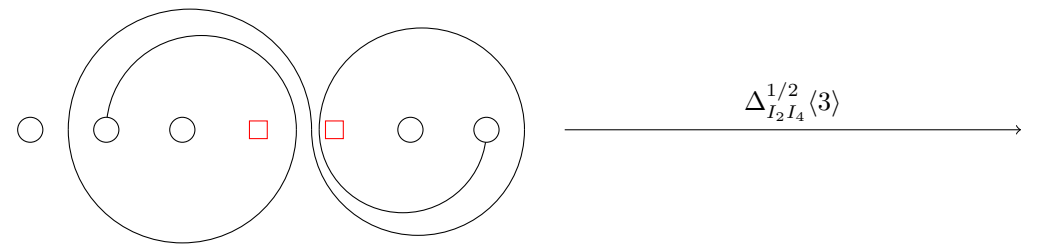
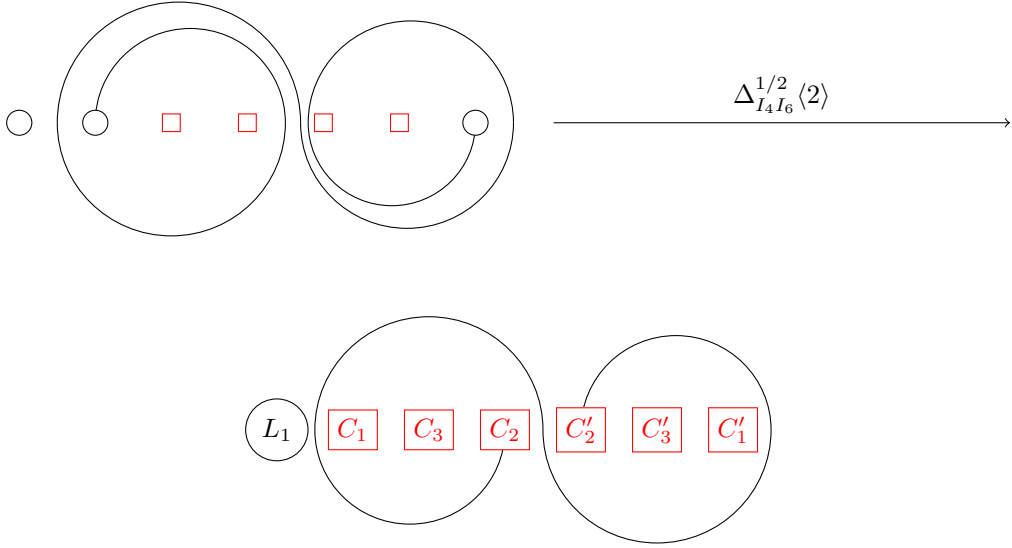
\begin{figure}[H]
	\begin{center}
		%
	\end{center}
\end{figure}

Relation:

\begin{equation*}\begin{split}
\Gamma_{C_3'}^{-1}\Gamma_{C_1'}^{-1}\Gamma_{C_2}\Gamma_{C_3}\Gamma_{C_1}\Gamma_{C_2}\Gamma_{C_1}^{-1}\Gamma_{C_3}^{-1}\Gamma_{C_2}^{-1}\Gamma_{C_1'}\Gamma_{C_3'}  = \Gamma_{C_2'}
\end{split}
\end{equation*}
\subsubsection{Raw relations and simplifications}

\begin{equation}
\begin{split}
\Gamma_{C_2}  = \Gamma_{C_2'} ,
\end{split}
\end{equation}
\begin{equation}
\begin{split}
\Gamma_{C_2'}\Gamma_{C_3}\Gamma_{C_2'}^{-1}  = \Gamma_{C_3'} ,
\end{split}
\end{equation}
\begin{equation}
\begin{split}
\left[\Gamma_{C_3} , \Gamma_{C_2}\Gamma_{C_3}\Gamma_{C_1}\Gamma_{L_1}\Gamma_{C_1}^{-1}\Gamma_{C_3}^{-1}  \right]=e,
\end{split}
\end{equation}
\begin{equation}
\begin{split}
\left[\Gamma_{C_2} , \Gamma_{C_3}\Gamma_{C_1}\Gamma_{L_1}\Gamma_{C_1}^{-1}  \right]=e,
\end{split}
\end{equation}
\begin{equation}
\begin{split}
\Gamma_{C_3'}\Gamma_{C_2'}\Gamma_{C_1}\Gamma_{C_2'}^{-1}\Gamma_{C_3'}^{-1}  = \Gamma_{C_1'} ,
\end{split}
\end{equation}
\begin{equation}
\begin{split}
\left[\Gamma_{C_2}\Gamma_{C_3}\Gamma_{C_1}\Gamma_{L_1}\Gamma_{C_1}^{-1}\Gamma_{C_3}^{-1}\Gamma_{C_2}^{-1} , \Gamma_{C_1}  \right]=e,
\end{split}
\end{equation}
\begin{equation}
\begin{split}
\left\{\Gamma_{C_1'} , \Gamma_{C_3'}  \right\}=e,
\end{split}
\end{equation}
\begin{equation}
\begin{split}
\left[\Gamma_{C_2'}\Gamma_{C_2}\Gamma_{C_3}\Gamma_{C_1}\Gamma_{L_1}\Gamma_{C_1}^{-1}\Gamma_{C_3}^{-1}\Gamma_{C_2}^{-1}\Gamma_{C_2'}^{-1} , \Gamma_{C_1'}  \right]=e,
\end{split}
\end{equation}
\begin{equation}
\begin{split}
\left[\Gamma_{C_1'}^{-1}\Gamma_{C_2'}\Gamma_{C_1'} , \Gamma_{C_3'}\Gamma_{C_2'}\Gamma_{C_2}\Gamma_{C_3}\Gamma_{C_1}\Gamma_{L_1}\Gamma_{C_1}^{-1}\Gamma_{C_3}^{-1}\Gamma_{C_2}^{-1}\Gamma_{C_2'}^{-1}  \right]=e,
\end{split}
\end{equation}
\begin{equation}
\begin{split}
\left[\Gamma_{C_3'} , \Gamma_{C_2'}\Gamma_{C_2}\Gamma_{C_3}\Gamma_{C_1}\Gamma_{L_1}\Gamma_{C_1}^{-1}\Gamma_{C_3}^{-1}\Gamma_{C_2}^{-1}\Gamma_{C_2'}^{-1}\Gamma_{C_1'}^{-1}\Gamma_{C_2'}\Gamma_{C_1'}  \right]=e,
\end{split}
\end{equation}
\begin{equation}
\begin{split}
\left\{\Gamma_{C_2} , \Gamma_{C_1}  \right\}=e,
\end{split}
\end{equation}
\begin{equation}
\begin{split}
\left\{\Gamma_{C_1'} , \Gamma_{C_2'}  \right\}=e,
\end{split}
\end{equation}
\begin{equation}
\begin{split}
\left[\Gamma_{C_1}^{-1}\Gamma_{C_3}\Gamma_{C_1} , \Gamma_{C_2}  \right]=e,
\end{split}
\end{equation}
\begin{equation}
\begin{split}
\left\{\Gamma_{C_3} , \Gamma_{C_1}  \right\}=e,
\end{split}
\end{equation}
\begin{equation}
\begin{split}
\Gamma_{C_2}\Gamma_{C_3}\Gamma_{C_1}\Gamma_{C_3}^{-1}\Gamma_{C_2}^{-1}  = \Gamma_{C_1'} ,
\end{split}
\end{equation}
\begin{equation}
\begin{split}
\left[\Gamma_{C_2'} , \Gamma_{C_3'}  \right]=e,
\end{split}
\end{equation}
\begin{equation}
\begin{split}
\Gamma_{C_1'}^{-1}\Gamma_{C_2}\Gamma_{C_3}\Gamma_{C_1}\Gamma_{C_3}\Gamma_{C_1}^{-1}\Gamma_{C_3}^{-1}\Gamma_{C_2}^{-1}\Gamma_{C_1'}  = \Gamma_{C_3'} ,
\end{split}
\end{equation}
\begin{equation}
\begin{split}
\Gamma_{C_3'}^{-1}\Gamma_{C_1'}^{-1}\Gamma_{C_2}\Gamma_{C_3}\Gamma_{C_1}\Gamma_{C_2}\Gamma_{C_1}^{-1}\Gamma_{C_3}^{-1}\Gamma_{C_2}^{-1}\Gamma_{C_1'}\Gamma_{C_3'}  = \Gamma_{C_2'} ,
\end{split}
\end{equation}
\begin{equation}
\begin{split}
\Gamma_{C_1'}\Gamma_{C_3'}\Gamma_{C_2'}\Gamma_{C_2}\Gamma_{C_3}\Gamma_{C_1}\Gamma_{L_1} =e.
\end{split}
\end{equation}

\bigskip
Here is the simplification of $ G_3 := \pcpt{\mathcal{B}_3}/\langle \Gamma_X^2 \; | \; X\subseteq \mathcal{B}_3 \rangle  $:
\begin{lemma}\label{lemma_tokunaga_arrangement_3_simplification}
Group $ G_3 := \pcpt{\mathcal{B}_3}/\langle \Gamma_X^2 \; | \; X\subseteq \mathcal{B}_3 \rangle  $ is isomorphic to $ \mbb{Z}_2^2\ltimes D_8  $.
\end{lemma}
\begin{proof}
Substituting $ \Gamma_{C_2'}=\Gamma_{C_2} $ we get that $ G_3 $ is generated by $ \Gamma_{L_1}, \Gamma_{C_1}, \Gamma_{C_1'}, \Gamma_{C_2}, \Gamma_{C_3} $, and $ \Gamma_{C_3'} $, subject to the following relations:
\begin{equation}\label{arrangement3_auto_calc_squares}
	\begin{split}
		\Gamma_{L_1}^2 = \Gamma_{C_1}^2 = \Gamma_{C_1'}^2 = \Gamma_{C_2}^2 = \Gamma_{C_3}^2 = \Gamma_{C_3'}^2 = e,
	\end{split}
\end{equation}
\begin{equation}\label{arrangement3_auto_calc_vert_2_rel_1}
	\begin{split}
		\Gamma_{C_2}\Gamma_{C_3}\Gamma_{C_2}  = \Gamma_{C_3'},
	\end{split}
\end{equation}
\begin{equation}\label{arrangement3_auto_calc_vert_3_rel_1}
	\begin{split}
		\left[\Gamma_{C_3} , \Gamma_{C_2}\Gamma_{C_3}\Gamma_{C_1}\Gamma_{L_1}\Gamma_{C_1}\Gamma_{C_3}\right]=e,
	\end{split}
\end{equation}
\begin{equation}\label{arrangement3_auto_calc_vert_3_rel_2}
	\begin{split}
		\left[\Gamma_{C_2} , \Gamma_{C_3}\Gamma_{C_1}\Gamma_{L_1}\Gamma_{C_1}  \right]=e,
	\end{split}
\end{equation}
\begin{equation}\label{arrangement3_auto_calc_vert_4_rel_1}
	\begin{split}
		\Gamma_{C_3'}\Gamma_{C_2}\Gamma_{C_1}\Gamma_{C_2}\Gamma_{C_3'}  = \Gamma_{C_1'},
	\end{split}
\end{equation}
\begin{equation}\label{arrangement3_auto_calc_vert_5_rel_1}
	\begin{split}
		\left[\Gamma_{C_2}\Gamma_{C_3}\Gamma_{C_1}\Gamma_{L_1}\Gamma_{C_1}\Gamma_{C_3}\Gamma_{C_2} , \Gamma_{C_1}  \right]=e,
	\end{split}
\end{equation}
\begin{equation}\label{arrangement3_auto_calc_vert_6_rel_1}
	\begin{split}
		\left\{\Gamma_{C_1'} , \Gamma_{C_3'}  \right\}=e,
	\end{split}
\end{equation}
\begin{equation}\label{arrangement3_auto_calc_vert_7_rel_1}
	\begin{split}
		\left[\Gamma_{C_3}\Gamma_{C_1}\Gamma_{L_1}\Gamma_{C_1}\Gamma_{C_3} , \Gamma_{C_1'}  \right]=e,
	\end{split}
\end{equation}
\begin{equation}\label{arrangement3_auto_calc_vert_8_rel_1}
	\begin{split}
		\left[\Gamma_{C_1'}\Gamma_{C_2}\Gamma_{C_1'} , \Gamma_{C_3'}\Gamma_{C_3}\Gamma_{C_1}\Gamma_{L_1}\Gamma_{C_1}\Gamma_{C_3}  \right]=e,
	\end{split}
\end{equation}
\begin{equation}\label{arrangement3_auto_calc_vert_8_rel_2}
	\begin{split}
		\left[\Gamma_{C_3'} , \Gamma_{C_3}\Gamma_{C_1}\Gamma_{L_1}\Gamma_{C_1}\Gamma_{C_3}
		\Gamma_{C_1'}\Gamma_{C_2}\Gamma_{C_1'}  \right]=e,
	\end{split}
\end{equation}
\begin{equation}\label{arrangement3_auto_calc_vert_9_rel_1}
	\begin{split}
		\left\{\Gamma_{C_2} , \Gamma_{C_1}  \right\}=e,
	\end{split}
\end{equation}
\begin{equation}\label{arrangement3_auto_calc_vert_10_rel_1}
	\begin{split}
		\left\{\Gamma_{C_1'} , \Gamma_{C_2}  \right\}=e,
	\end{split}
\end{equation}
\begin{equation}\label{arrangement3_auto_calc_vert_11_rel_1}
	\begin{split}
		\left[\Gamma_{C_1}\Gamma_{C_3}\Gamma_{C_1} , \Gamma_{C_2}  \right]=e,
	\end{split}
\end{equation}
\begin{equation}\label{arrangement3_auto_calc_vert_12_rel_1}
	\begin{split}
		\left\{\Gamma_{C_3} , \Gamma_{C_1}  \right\}=e,
	\end{split}
\end{equation}
\begin{equation}\label{arrangement3_auto_calc_vert_13_rel_1}
	\begin{split}
		\Gamma_{C_2}\Gamma_{C_3}\Gamma_{C_1}\Gamma_{C_3}\Gamma_{C_2}  = \Gamma_{C_1'},
	\end{split}
\end{equation}
\begin{equation}\label{arrangement3_auto_calc_vert_14_rel_1}
	\begin{split}
		\left[\Gamma_{C_2} , \Gamma_{C_3'}  \right]=e,
	\end{split}
\end{equation}
\begin{equation}\label{arrangement3_auto_calc_vert_15_rel_1}
	\begin{split}	\Gamma_{C_1'}\Gamma_{C_2}\Gamma_{C_3}\Gamma_{C_1}\Gamma_{C_3}\Gamma_{C_1}\Gamma_{C_3}\Gamma_{C_2}\Gamma_{C_1'}  = \Gamma_{C_3'},
	\end{split}
\end{equation}
\begin{equation}\label{arrangement3_auto_calc_vert_16_rel_1}
	\begin{split} \Gamma_{C_3'}\Gamma_{C_1'}\Gamma_{C_2}\Gamma_{C_3}\Gamma_{C_1}\Gamma_{C_2}\Gamma_{C_1}\Gamma_{C_3}
		\Gamma_{C_2}\Gamma_{C_1'}\Gamma_{C_3'}  = \Gamma_{C_2},
	\end{split}
\end{equation}
\begin{equation}\label{arrangement3_auto_calc_projective_rel}
	\begin{split}
		\Gamma_{C_1'}\Gamma_{C_3'}\Gamma_{C_3}\Gamma_{C_1}\Gamma_{L_1} =e.
	\end{split}
\end{equation}

Because $ \Gamma_{C_2} $ and $ \Gamma_{C_3'} $ commute by \eqref{arrangement3_auto_calc_vert_14_rel_1}, we get that   $\Gamma_{C_3}=\Gamma_{C_3'}$ from \eqref{arrangement3_auto_calc_vert_2_rel_1}.

From this we simplify \eqref{arrangement3_auto_calc_vert_3_rel_1} to $ \left[ \Gamma_{C_3}, \Gamma_{C_1}\Gamma_{L_1}\Gamma_{C_1} \right]=e  $ and  \eqref{arrangement3_auto_calc_vert_3_rel_2} to $ \left[ \Gamma_{C_2}, \Gamma_{C_1}\Gamma_{L_1}\Gamma_{C_1} \right]=e  $.
This allows us to rewrite \eqref{arrangement3_auto_calc_vert_5_rel_1} as $ [\Gamma_{C_1}, \Gamma_{L_1}]=e $ and thus, we get $\left[\Gamma_{L_1} , \Gamma_{C_3}  \right]=e$ and $\left[\Gamma_{L_1} , \Gamma_{C_2}  \right]=e$.
These commutations simplify easily \eqref{arrangement3_auto_calc_vert_7_rel_1} to $\left[\Gamma_{L_1} , \Gamma_{C_1}  \right]=e$  and it is redundant.

Now we simplify \eqref{arrangement3_auto_calc_vert_8_rel_2} using the above resulting commutations:
\begin{align*} \left[\Gamma_{C_3},\Gamma_{C_3}\Gamma_{C_1}\Gamma_{L_1}\Gamma_{C_1}\Gamma_{C_3}\Gamma_{C_1'}\Gamma_{C_2}
	\Gamma_{C_1'}  \right]=e \Rightarrow \left[ \Gamma_{C_3}, \Gamma_{C_1'}\Gamma_{C_2}\Gamma_{C_1'} \right] = e \Rightarrow  \left[\Gamma_{C_3},\Gamma_{C_1}\Gamma_{C_2}\Gamma_{C_1}  \right]=e.
\end{align*}
In a similar way we get from \eqref{arrangement3_auto_calc_vert_8_rel_1} the same relation and it is redundant.

Because $ \Gamma_{C_2} $ and $ \Gamma_{C_3} $ commute, substituting \eqref{arrangement3_auto_calc_vert_13_rel_1} into \eqref{arrangement3_auto_calc_vert_10_rel_1} gives us $ \{ \Gamma_{C_1}, \Gamma_{C_2} \} = e $, and it is redundant because we have it already in  \eqref{arrangement3_auto_calc_vert_9_rel_1}.
Similarly we can simplify $\left\{\Gamma_{C_1'} , \Gamma_{C_3'}  \right\}=e$ (from \eqref{arrangement3_auto_calc_vert_6_rel_1}) to $\left\{\Gamma_{C_1} , \Gamma_{C_3}  \right\}=e$, and it is redundant because we have it already in \eqref{arrangement3_auto_calc_vert_12_rel_1}.

In \eqref{arrangement3_auto_calc_vert_15_rel_1} we substitute $\Gamma_{C_2}\Gamma_{C_3}\Gamma_{C_1}\Gamma_{C_3}\Gamma_{C_2}  = \Gamma_{C_1'}$ and it is redundant: $$\Gamma_{C_1'}\Gamma_{C_2}\Gamma_{C_3}\Gamma_{C_1}\Gamma_{C_3}\Gamma_{C_1}\Gamma_{C_3}\Gamma_{C_2}\Gamma_{C_1'}  = \Gamma_{C_3'} \Rightarrow$$ $$\Gamma_{C_2}\Gamma_{C_3}\Gamma_{C_1}\Gamma_{C_3}\Gamma_{C_2}\Gamma_{C_2}\Gamma_{C_3}\Gamma_{C_1}\Gamma_{C_3}\Gamma_{C_1}\Gamma_{C_3}\Gamma_{C_2} \Gamma_{C_2}\Gamma_{C_3}\Gamma_{C_1}\Gamma_{C_3}\Gamma_{C_2}  = \Gamma_{C_3} \Rightarrow$$
$$\Gamma_{C_3}  = \Gamma_{C_3}.$$
In a similar way,  relation  \eqref{arrangement3_auto_calc_vert_16_rel_1} is also redundant.

The last relation to simplify is \eqref{arrangement3_auto_calc_projective_rel}. We substitute
$\Gamma_{C_1'}=\Gamma_{C_2}\Gamma_{C_3}\Gamma_{C_1}\Gamma_{C_3}\Gamma_{C_2}$ in it and get
$$\Gamma_{L_1} =(\Gamma_{C_2}\Gamma_{C_3}\Gamma_{C_1})^{2}.$$

Following is the list of all simplified relations:
\begin{equation}
	\Gamma_{L_1}^2 = \Gamma_{C_1}^2 = \Gamma_{C_2}^2 = \Gamma_{C_3}^2 = e,
\end{equation}
\begin{equation}\label{arrangement3_simp_comm}
	\begin{split}
		\left[\Gamma_{C_2} , \Gamma_{C_3}  \right]=\left[\Gamma_{C_1} , \Gamma_{L_1}  \right]= \left[\Gamma_{C_2} , \Gamma_{L_1}  \right]= \left[\Gamma_{C_3} , \Gamma_{L_1}  \right]=e,
	\end{split}
\end{equation}
\begin{equation}\label{arrangement3_simp_fork}
	\begin{split}
		\left[\Gamma_{C_3} , \Gamma_{C_1}\Gamma_{C_2}\Gamma_{C_1} \right] =e,
	\end{split}
\end{equation}
\begin{equation}\label{arrangement3_simp_quad}
	\begin{split}
		\left\{\Gamma_{C_1} , \Gamma_{C_2}  \right\} = \left\{\Gamma_{C_1} , \Gamma_{C_3}  \right\}=e,
	\end{split}
\end{equation}
\begin{equation}\label{arrangement3_simp_l1}
	\begin{split}
		\Gamma_{L_1} = (\Gamma_{C_2}\Gamma_{C_3}\Gamma_{C_1})^{2}.
	\end{split}
\end{equation}

By substituting \eqref{arrangement3_simp_l1} into the relation $\left[ \Gamma_{C_3}, \Gamma_{L_1} \right] = e  $, we show that it is redundant:
$$\left[\Gamma_{C_3} , \Gamma_{L_1}  \right]=e \Rightarrow$$
$$\left[\Gamma_{C_3} , (\Gamma_{C_2}\Gamma_{C_3}\Gamma_{C_1})^2 \right]=e \Rightarrow$$
$$\Gamma_{C_3} \Gamma_{C_2}\Gamma_{C_3}\Gamma_{C_1} \Gamma_{C_2}\Gamma_{C_3}\Gamma_{C_1} = \Gamma_{C_2}\Gamma_{C_3}\Gamma_{C_1} \Gamma_{C_2}\Gamma_{C_3}\Gamma_{C_1} \Gamma_{C_3} \Rightarrow$$
$$\Gamma_{C_3}\Gamma_{C_1} \Gamma_{C_2}\Gamma_{C_3}\Gamma_{C_1} = \Gamma_{C_1} \Gamma_{C_2}\Gamma_{C_3}\Gamma_{C_1} \Gamma_{C_3} \Rightarrow$$
$$\Gamma_{C_3}\Gamma_{C_1} \Gamma_{C_2} = \Gamma_{C_1} \Gamma_{C_2}(\Gamma_{C_1}\Gamma_{C_1})\Gamma_{C_3}\Gamma_{C_1} \Gamma_{C_3} \Gamma_{C_1}\Gamma_{C_3}\Rightarrow$$
$$\Gamma_{C_3}\Gamma_{C_1} \Gamma_{C_2} = \Gamma_{C_1} \Gamma_{C_2}\Gamma_{C_1}\Gamma_{C_3}\Gamma_{C_1}\Gamma_{C_3}\Gamma_{C_1} \Gamma_{C_1}\Gamma_{C_3}\Rightarrow$$
$$\Gamma_{C_3}\Gamma_{C_1} \Gamma_{C_2} = \Gamma_{C_1} \Gamma_{C_2}\Gamma_{C_1}\Gamma_{C_3}\Gamma_{C_1}\Rightarrow$$
$$\Gamma_{C_3}\Gamma_{C_1} \Gamma_{C_2} = \Gamma_{C_3}\Gamma_{C_1} \Gamma_{C_2}\Gamma_{C_1}\Gamma_{C_1}\Rightarrow$$
$$\Gamma_{C_1} =  \Gamma_{C_1}.$$
In a very similar way, $ \left[ \Gamma_{C_1}, \Gamma_{L_1} \right]=e$ and $ \left[ \Gamma_{C_2}, \Gamma_{L_1} \right]=e  $ are also redundant.

Therefore, $\Gamma_{L_1}$ is eliminated, and group $ G_3 $ is now generated by $\Gamma_{C_1}, \Gamma_{C_2}$, and $\Gamma_{C_3}$ and admits the following relations:
\begin{equation}\label{eq:arrangement3_final1}
	\Gamma_{C_1}^2 = \Gamma_{C_2}^2 = \Gamma_{C_3}^2 = e,
\end{equation}
\begin{equation}\label{eq:arrangement3_final2}
	\begin{split}
		\left[\Gamma_{C_2} , \Gamma_{C_3}  \right]=  \left[ \Gamma_{C_1}\Gamma_{C_2}\Gamma_{C_1}, \Gamma_{C_3} \right]=e,
	\end{split}
\end{equation}
\begin{equation}\label{eq:arrangement3_final3}
	\begin{split}
		\left\{\Gamma_{C_1} , \Gamma_{C_2}  \right\}= \left\{\Gamma_{C_1} , \Gamma_{C_3}  \right\}=e.
	\end{split}
\end{equation}

The group generated by $ \Gamma_{C_1},\Gamma_{C_2} $, and $ \Gamma_{C_3} $, subject to the relations \eqref{eq:arrangement3_final1}-\eqref{eq:arrangement3_final3},  without the relation $ \left[\Gamma_{C_1}\Gamma_{C_2}\Gamma_{C_1}, \Gamma_{C_3}\right] = e $, is a Coxeter group of type $ \widetilde{C}_2 $ and its elements can be represented by signed permutations
$$  \begin{pmatrix}
	1 & 2 \\
	\sigma(1)^{\epsilon_1, n_1} & \sigma(2)^{\epsilon_2, n_2}
\end{pmatrix}. $$
The generators correspond to
$$ \Gamma_{C_1} \mapsto \begin{pmatrix} 1 & 2 \\ 2 & 1 \end{pmatrix} \; ; \;  \Gamma_{C_2} \mapsto \begin{pmatrix} 1 & 2 \\ \overline{1} & 2 \end{pmatrix} \; ; \; \Gamma_{C_1} \mapsto \begin{pmatrix} 1 & 2 \\ 1 & \overline{2}^1 \end{pmatrix}. $$
Now computing the relation  $ \left[\Gamma_{C_1}\Gamma_{C_2}\Gamma_{C_1}, \Gamma_{C_3}\right] = e $ we get
$$ \Gamma_{C_1}\Gamma_{C_2}\Gamma_{C_1}\Gamma_{C_3} \mapsto \begin{pmatrix} 1 & 2 \\ 1 & 2^{-1} \end{pmatrix} $$
and
$$ \Gamma_{C_3}\Gamma_{C_1}\Gamma_{C_2}\Gamma_{C_1} \mapsto \begin{pmatrix} 1 & 2 \\ 1 & 2^{1} \end{pmatrix} $$
meaning that the elements of $ G_3 $ are signed permutations that are identified if the underlying permutations are the same, the signs are the same, and the exponents have the same parity.
Because the subgroup of a Coxeter group of type $ \widetilde{C}_2 $ of elements without exponent is exactly $ D_8 $ (see \cite{Bjorner2005}), we get that $ G_3\cong \mbb{Z}_2^2\ltimes D_8  $ as needed.
\end{proof}

\newpage

\subsection{$ \pcpt{\mathcal{B}_4} $}

\begin{figure}[H]
	\begin{center}
		\begin{tikzpicture}
		\begin{scope}[ , scale=1]
		\draw [rotate around={90.000000:(0.000000,0.000000)}] (0.000000,0.000000) ellipse (4.000000 and 4.000000);
		\node  (C_1_naming_node) at (0.000000,3.800000) {$ C_1 $};
		\draw [] (6.586716, -3.071436) -- (-6.586716, 3.071436) node[very near end, below] {$L_2$};
		\draw [rotate around={-70.000000:(-0.000000,-0.000000)}] (-0.000000,-0.000000) ellipse (5.656854 and 4.000000);
		\node  (C_3_naming_node) at (1.934758,-5.815704) {$ C_3 $};
		\draw [rotate around={-160.000000:(-0.000000,0.000000)}] (-0.000000,0.000000) ellipse (5.656854 and 4.000000);
		\node  (C_2_naming_node) at (-5.815704,-1.934758) {$ C_2 $};
		\draw [fill=gray] (5.488930, 0.936849) circle (5.000000pt);
		\node [label=above:1] (sing_pnt_1) at (5.488930,0.936849) {};
		\draw [fill=gray] (4.227487, -1.216397) circle (5.000000pt);
		\node [label=above:2] (sing_pnt_2) at (4.227487,-1.216397) {};
		\draw [fill=gray] (4.186056, -1.951990) circle (5.000000pt);
		\node [label=above:3] (sing_pnt_3) at (4.186056,-1.951990) {};
		\draw [fill=gray] (4.000000, -0.000000) circle (5.000000pt);
		\node [label=above:4] (sing_pnt_4) at (4.000000,-0.000000) {};
		\draw [fill=gray] (3.758770, 1.368081) circle (5.000000pt);
		\node [label=above:5] (sing_pnt_5) at (3.758770,1.368081) {};
		\draw [fill=gray] (3.625231, -1.690473) circle (5.000000pt);
		\node [label=above:6] (sing_pnt_6) at (3.625231,-1.690473) {};
		\draw [fill=gray] (1.951990, 4.186056) circle (5.000000pt);
		\node [label=above:7] (sing_pnt_7) at (1.951990,4.186056) {};
		\draw [fill=gray] (1.368081, -3.758770) circle (5.000000pt);
		\node [label=above:8] (sing_pnt_8) at (1.368081,-3.758770) {};
		\draw [fill=gray] (-1.368081, 3.758770) circle (5.000000pt);
		\node [label=above:9] (sing_pnt_9) at (-1.368081,3.758770) {};
		\draw [fill=gray] (-1.951990, -4.186056) circle (5.000000pt);
		\node [label=above:10] (sing_pnt_10) at (-1.951990,-4.186056) {};
		\draw [fill=gray] (-3.625231, 1.690473) circle (5.000000pt);
		\node [label=above:11] (sing_pnt_11) at (-3.625231,1.690473) {};
		\draw [fill=gray] (-3.758770, -1.368081) circle (5.000000pt);
		\node [label=above:12] (sing_pnt_12) at (-3.758770,-1.368081) {};
		\draw [fill=gray] (-4.000000, 0.000000) circle (5.000000pt);
		\node [label=above:13] (sing_pnt_13) at (-4.000000,0.000000) {};
		\draw [fill=gray] (-4.186056, 1.951990) circle (5.000000pt);
		\node [label=above:14] (sing_pnt_14) at (-4.186056,1.951990) {};
		\draw [fill=gray] (-4.227487, 1.216397) circle (5.000000pt);
		\node [label=above:15] (sing_pnt_15) at (-4.227487,1.216397) {};
		\draw [fill=gray] (-5.488930, -0.936849) circle (5.000000pt);
		\node [label=above:16] (sing_pnt_16) at (-5.488930,-0.936849) {};
		
		\end{scope}[]
		\end{tikzpicture}
	\end{center}
\end{figure}

{Monodromy table}

\begin{tabular}{| c | c | c | c |}
	\hline \hline
	Vertex number& Vertex description& Skeleton& Diffeomorphism\\ \hline \hline
	1& $C_2$ branch& $ \langle 2 - 3 \rangle $& $ \Delta_{I_{4}I_{6}}^{1/2}\langle 2 \rangle $ \\ \hline
	2& $C_3$ branch& $ \langle 3 - 4 \rangle $& $ \Delta_{I_{2}I_{4}}^{1/2}\langle 3 \rangle $ \\ \hline
	3& Node between $C_3$, $C_2$ and $L_2$& $ \langle 1 - 2 - 3 \rangle $& $ \Delta \langle 1,2,3 \rangle $ \\ \hline
	4& $C_1$ branch& $ \langle 4 - 5 \rangle $& $ \Delta_{\mathbb{R}I_{2}}^{1/2}\langle 4 \rangle $ \\ \hline
	5& Tangency between $C_1$ and $C_3$& $ \langle 5 - 6 \rangle $& $ \Delta^{2}\langle 5, 6 \rangle $ \\ \hline
	6& Node between $C_1$ and $L_2$& $ \langle 3 - 4 \rangle $& $ \Delta\langle 3, 4 \rangle $ \\ \hline
	7& Node between $C_2$ and $C_3$& $ \langle 6 - 7 \rangle $& $ \Delta\langle 6, 7 \rangle $ \\ \hline
	8& Tangency between $C_1$ and $C_2$& $ \langle 2 - 3 \rangle $& $ \Delta^{2}\langle 2, 3 \rangle $ \\ \hline
	9& Tangency between $C_1$ and $C_2$& $ \langle 5 - 6 \rangle $& $ \Delta^{2}\langle 5, 6 \rangle $ \\ \hline
	10& Node between $C_2$ and $C_3$& $ \langle 1 - 2 \rangle $& $ \Delta\langle 1, 2 \rangle $ \\ \hline
	11& Node between $C_1$ and $L_2$& $ \langle 4 - 5 \rangle $& $ \Delta\langle 4, 5 \rangle $ \\ \hline
	12& Tangency between $C_1$ and $C_3$& $ \langle 2 - 3 \rangle $& $ \Delta^{2}\langle 2, 3 \rangle $ \\ \hline
	13& $C_1$ branch& $ \langle 3 - 4 \rangle $& $ \Delta_{I_{2}\mathbb{R}}^{1/2}\langle 3 \rangle $ \\ \hline
	14& Node between $C_3$, $C_2$ and $L_2$& $ \langle 3 - 4 - 5 \rangle $& $ \Delta \langle 3,4,5 \rangle $ \\ \hline
	15& $C_3$ branch& $ \langle 2 - 3 \rangle $& $ \Delta_{I_{4}I_{2}}^{1/2}\langle 2 \rangle $ \\ \hline
	16& $C_2$ branch& $ \langle 1 - 2 \rangle $& $ \Delta_{I_{6}I_{4}}^{1/2}\langle 1 \rangle $
	\\ \hline \hline
\end{tabular}

\subsubsection{Vertex Number 1}

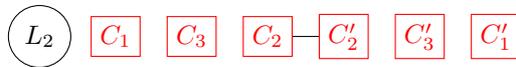
\begin{figure}[H]
	\begin{center}
		\begin{tikzpicture}
		\begin{scope}[]
		\draw [] (3.000000, 0.000000) -- (4.000000, 0.000000) ;
		\node [draw, circle, color=black, fill=white] (vert_0) at (0.000000,0.000000) {$ L_2 $};
		\node [draw, rectangle, color=red, fill=white] (vert_1) at (1.000000,0.000000) {$ C_1 $};
		\node [draw, rectangle, color=red, fill=white] (vert_2) at (2.000000,0.000000) {$ C_3 $};
		\node [draw, rectangle, color=red, fill=white] (vert_3) at (3.000000,0.000000) {$ C_2 $};
		\node [draw, rectangle, color=red, fill=white] (vert_4) at (4.000000,0.000000) {$ C_2' $};
		\node [draw, rectangle, color=red, fill=white] (vert_5) at (5.000000,0.000000) {$ C_3' $};
		\node [draw, rectangle, color=red, fill=white] (vert_6) at (6.000000,0.000000) {$ C_1' $};
		
		\end{scope}[]
		\end{tikzpicture}
	\end{center}
\end{figure}

Relation:

\begin{equation*}\begin{split}
\Gamma_{C_2}  = \Gamma_{C_2'}
\end{split}
\end{equation*}
\subsubsection{Vertex Number 2}

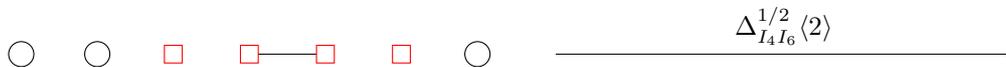
\begin{figure}[H]
	\begin{center}
		\begin{tikzpicture}
		\begin{scope}[scale=1 , xshift=-100]
		\draw [] (3.000000, 0.000000) -- (4.000000, 0.000000) ;
		\node [draw, circle, color=black, fill=white] (vert_0) at (0.000000,0.000000) {};
		\node [draw, circle, color=black, fill=white] (vert_1) at (1.000000,0.000000) {};
		\node [draw, rectangle, color=red, fill=white] (vert_2) at (2.000000,0.000000) {};
		\node [draw, rectangle, color=red, fill=white] (vert_3) at (3.000000,0.000000) {};
		\node [draw, rectangle, color=red, fill=white] (vert_4) at (4.000000,0.000000) {};
		\node [draw, rectangle, color=red, fill=white] (vert_5) at (5.000000,0.000000) {};
		\node [draw, circle, color=black, fill=white] (vert_6) at (6.000000,0.000000) {};
		
		\end{scope}[]
		\begin{scope}[scale=1 , xshift=100]
		\draw [->] (0.000000, 0.000000) -- node [above, midway] {$ \Delta_{I_{4}I_{6}}^{1/2}\langle 2 \rangle $}  (6.000000, 0.000000) ;
		
		\end{scope}[]
		\end{tikzpicture}
	\end{center}
\end{figure}

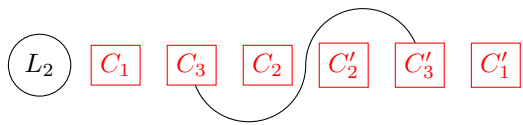
\begin{figure}[H]
	\begin{center}
		\begin{tikzpicture}
		\begin{scope}[]
		\draw [] (2.000000, 0.000000) arc (180.000000:360.000000:0.750000);
		\draw [] (3.500000, 0.000000) arc (180.000000:0.000000:0.750000);
		\node [draw, circle, color=black, fill=white] (vert_0) at (0.000000,0.000000) {$ L_2 $};
		\node [draw, rectangle, color=red, fill=white] (vert_1) at (1.000000,0.000000) {$ C_1 $};
		\node [draw, rectangle, color=red, fill=white] (vert_2) at (2.000000,0.000000) {$ C_3 $};
		\node [draw, rectangle, color=red, fill=white] (vert_3) at (3.000000,0.000000) {$ C_2 $};
		\node [draw, rectangle, color=red, fill=white] (vert_4) at (4.000000,0.000000) {$ C_2' $};
		\node [draw, rectangle, color=red, fill=white] (vert_5) at (5.000000,0.000000) {$ C_3' $};
		\node [draw, rectangle, color=red, fill=white] (vert_6) at (6.000000,0.000000) {$ C_1' $};
		
		\end{scope}[]
		\end{tikzpicture}
	\end{center}
\end{figure}

Relation:

\begin{equation*}\begin{split}
\Gamma_{C_2'}\Gamma_{C_3}\Gamma_{C_2'}^{-1}  = \Gamma_{C_3'}
\end{split}
\end{equation*}
\subsubsection{Vertex Number 3}

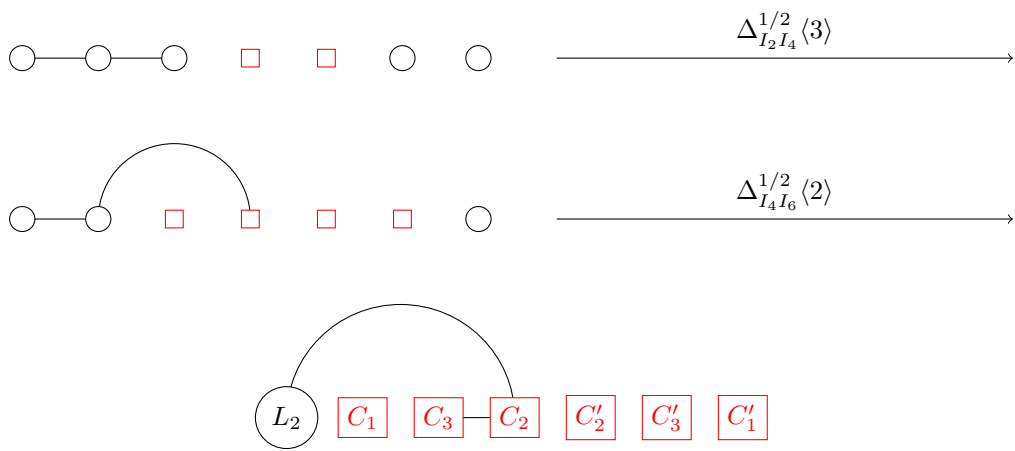
\begin{figure}[H]
	\begin{center}
		\begin{tikzpicture}
		\begin{scope}[scale=1 , xshift=-100]
		\draw [] (0.000000, 0.000000) -- (1.000000, 0.000000) ;
		\draw [] (1.000000, 0.000000) -- (2.000000, 0.000000) ;
		\node [draw, circle, color=black, fill=white] (vert_0) at (0.000000,0.000000) {};
		\node [draw, circle, color=black, fill=white] (vert_1) at (1.000000,0.000000) {};
		\node [draw, circle, color=black, fill=white] (vert_2) at (2.000000,0.000000) {};
		\node [draw, rectangle, color=red, fill=white] (vert_3) at (3.000000,0.000000) {};
		\node [draw, rectangle, color=red, fill=white] (vert_4) at (4.000000,0.000000) {};
		\node [draw, circle, color=black, fill=white] (vert_5) at (5.000000,0.000000) {};
		\node [draw, circle, color=black, fill=white] (vert_6) at (6.000000,0.000000) {};
		
		\end{scope}[]
		\begin{scope}[scale=1 , xshift=100]
		\draw [->] (0.000000, 0.000000) -- node [above, midway] {$ \Delta_{I_{2}I_{4}}^{1/2}\langle 3 \rangle $}  (6.000000, 0.000000) ;
		
		\end{scope}[]
		\end{tikzpicture}
	\end{center}
\end{figure}

\begin{figure}[H]
	\begin{center}
		\begin{tikzpicture}
		\begin{scope}[scale=1 , xshift=-100]
		\draw [] (0.000000, 0.000000) -- (1.000000, 0.000000) ;
		\draw [] (1.000000, 0.000000) arc (180.000000:0.000000:1.000000);
		\node [draw, circle, color=black, fill=white] (vert_0) at (0.000000,0.000000) {};
		\node [draw, circle, color=black, fill=white] (vert_1) at (1.000000,0.000000) {};
		\node [draw, rectangle, color=red, fill=white] (vert_2) at (2.000000,0.000000) {};
		\node [draw, rectangle, color=red, fill=white] (vert_3) at (3.000000,0.000000) {};
		\node [draw, rectangle, color=red, fill=white] (vert_4) at (4.000000,0.000000) {};
		\node [draw, rectangle, color=red, fill=white] (vert_5) at (5.000000,0.000000) {};
		\node [draw, circle, color=black, fill=white] (vert_6) at (6.000000,0.000000) {};
		
		\end{scope}[]
		\begin{scope}[scale=1 , xshift=100]
		\draw [->] (0.000000, 0.000000) -- node [above, midway] {$ \Delta_{I_{4}I_{6}}^{1/2}\langle 2 \rangle $}  (6.000000, 0.000000) ;
		
		\end{scope}[]
		\end{tikzpicture}
	\end{center}
\end{figure}

\begin{figure}[H]
	\begin{center}
		\begin{tikzpicture}
		\begin{scope}[]
		\draw [] (0.000000, 0.000000) arc (180.000000:0.000000:1.500000);
		\draw [] (3.000000, 0.000000) -- (2.000000, 0.000000) ;
		\node [draw, circle, color=black, fill=white] (vert_0) at (0.000000,0.000000) {$ L_2 $};
		\node [draw, rectangle, color=red, fill=white] (vert_1) at (1.000000,0.000000) {$ C_1 $};
		\node [draw, rectangle, color=red, fill=white] (vert_2) at (2.000000,0.000000) {$ C_3 $};
		\node [draw, rectangle, color=red, fill=white] (vert_3) at (3.000000,0.000000) {$ C_2 $};
		\node [draw, rectangle, color=red, fill=white] (vert_4) at (4.000000,0.000000) {$ C_2' $};
		\node [draw, rectangle, color=red, fill=white] (vert_5) at (5.000000,0.000000) {$ C_3' $};
		\node [draw, rectangle, color=red, fill=white] (vert_6) at (6.000000,0.000000) {$ C_1' $};
		
		\end{scope}[]
		\end{tikzpicture}
	\end{center}
\end{figure}

Relation:

\begin{equation*}\begin{split}
\left( \Gamma_{C_3}  \right) \left( \Gamma_{C_2}  \right) \left( \Gamma_{C_3}\Gamma_{C_1}\Gamma_{L_2}\Gamma_{C_1}^{-1}\Gamma_{C_3}^{-1}  \right)& = \\
\left( \Gamma_{C_2}  \right) \left( \Gamma_{C_3}\Gamma_{C_1}\Gamma_{L_2}\Gamma_{C_1}^{-1}\Gamma_{C_3}^{-1}  \right) \left( \Gamma_{C_3}  \right)& = \\
\left( \Gamma_{C_3}\Gamma_{C_1}\Gamma_{L_2}\Gamma_{C_1}^{-1}\Gamma_{C_3}^{-1}  \right) \left( \Gamma_{C_3}  \right) \left( \Gamma_{C_2}  \right)
\end{split}
\end{equation*}
\subsubsection{Vertex Number 4}

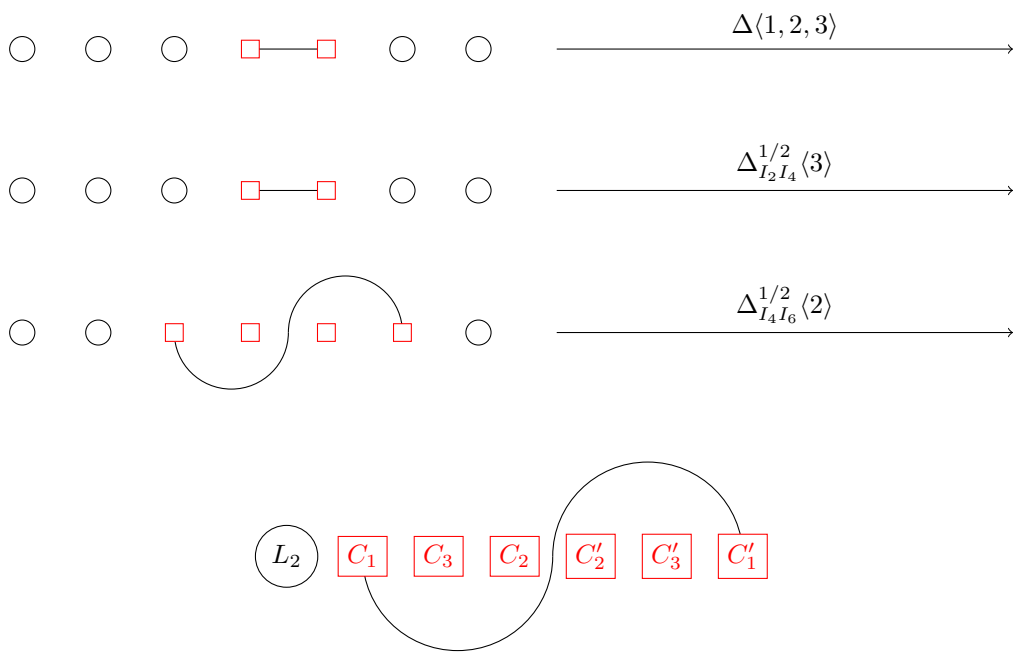
\begin{figure}[H]
	\begin{center}

	\end{center}
\end{figure}

Relation:

\begin{equation*}\begin{split}
\Gamma_{C_3'}\Gamma_{C_2'}\Gamma_{C_1}\Gamma_{C_2'}^{-1}\Gamma_{C_3'}^{-1}  = \Gamma_{C_1'}
\end{split}
\end{equation*}
\subsubsection{Vertex Number 5}

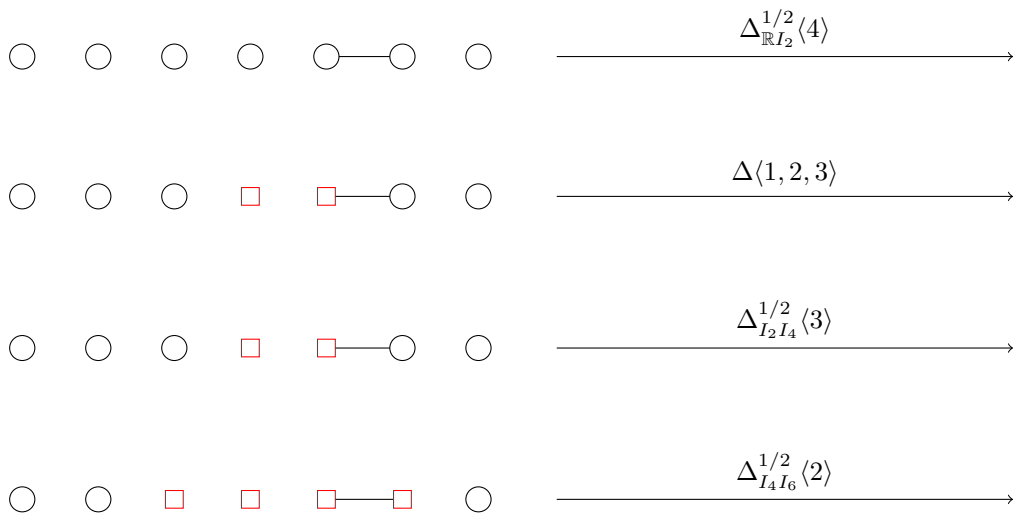
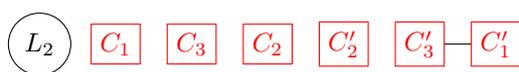
\begin{figure}[H]
	\begin{center}
		%
	\end{center}
\end{figure}

Relation:

\begin{equation*}\begin{split}
\left\{\Gamma_{C_1'}  , \Gamma_{C_3'} \right\}=e
\end{split}
\end{equation*}
\subsubsection{Vertex Number 6}

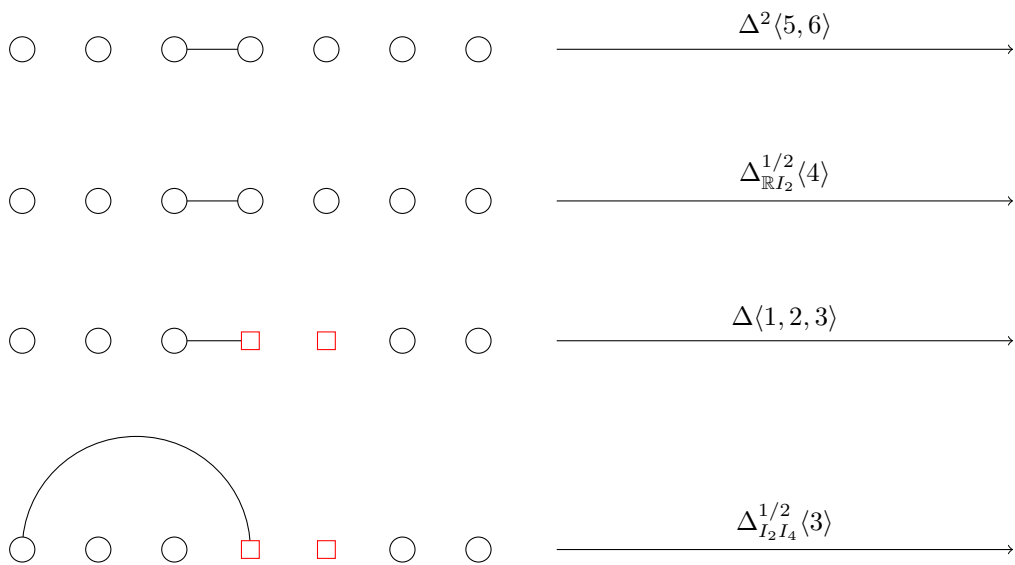
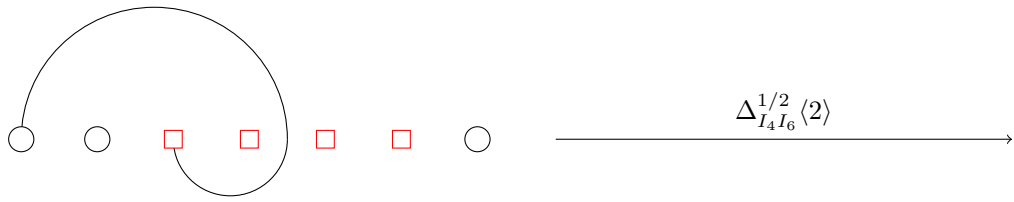
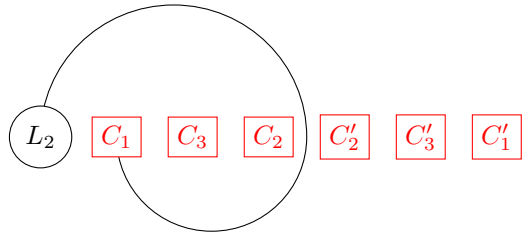
\begin{figure}[H]
	\begin{center}
		%
	\end{center}
\end{figure}

Relation:

\begin{equation*}\begin{split}
\left[\Gamma_{C_2}\Gamma_{C_3}\Gamma_{C_1}\Gamma_{L_2}\Gamma_{C_1}^{-1}\Gamma_{C_3}^{-1}\Gamma_{C_2}^{-1} , \Gamma_{C_1} \right]=e
\end{split}
\end{equation*}
\subsubsection{Vertex Number 7}

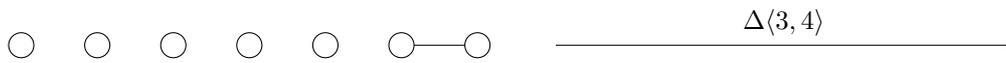
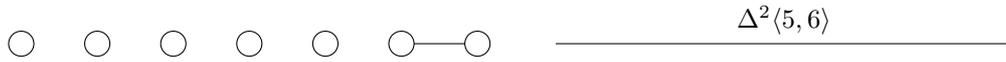
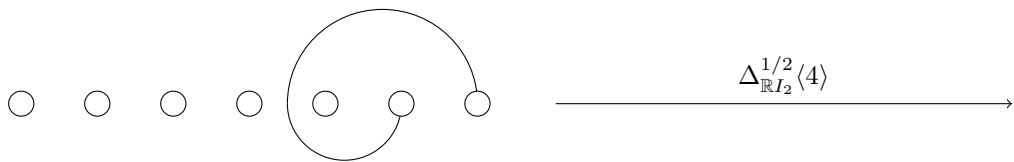
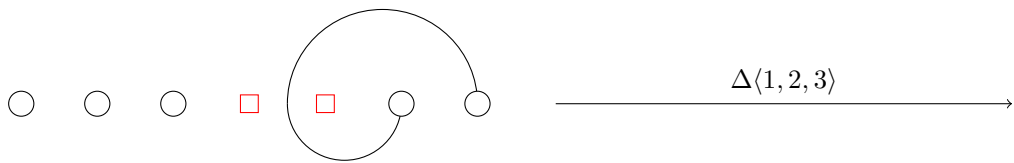
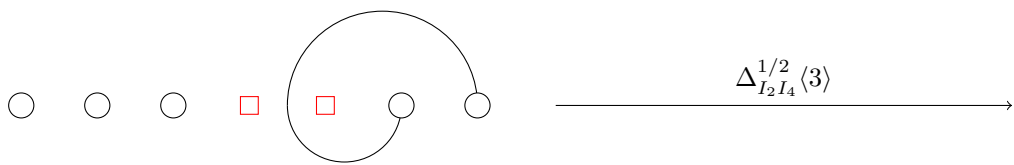
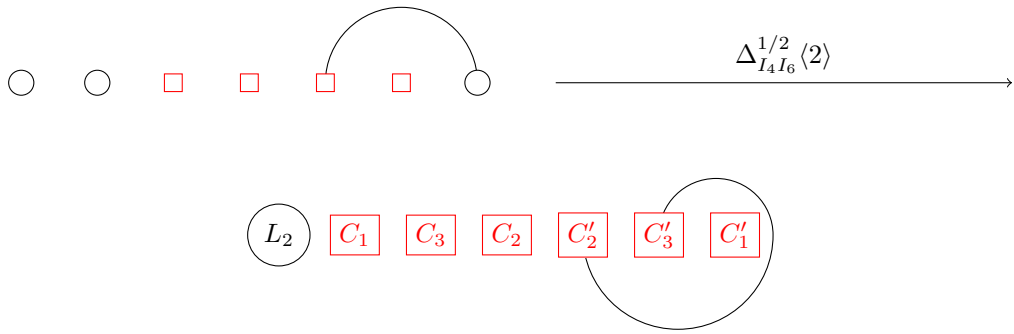
\begin{figure}[H]
	\begin{center}
		%
	\end{center}
\end{figure}

Relation:

\begin{equation*}\begin{split}
\left[\Gamma_{C_1'}\Gamma_{C_3'}\Gamma_{C_1'}^{-1} , \Gamma_{C_2'} \right]=e
\end{split}
\end{equation*}
\subsubsection{Vertex Number 8}

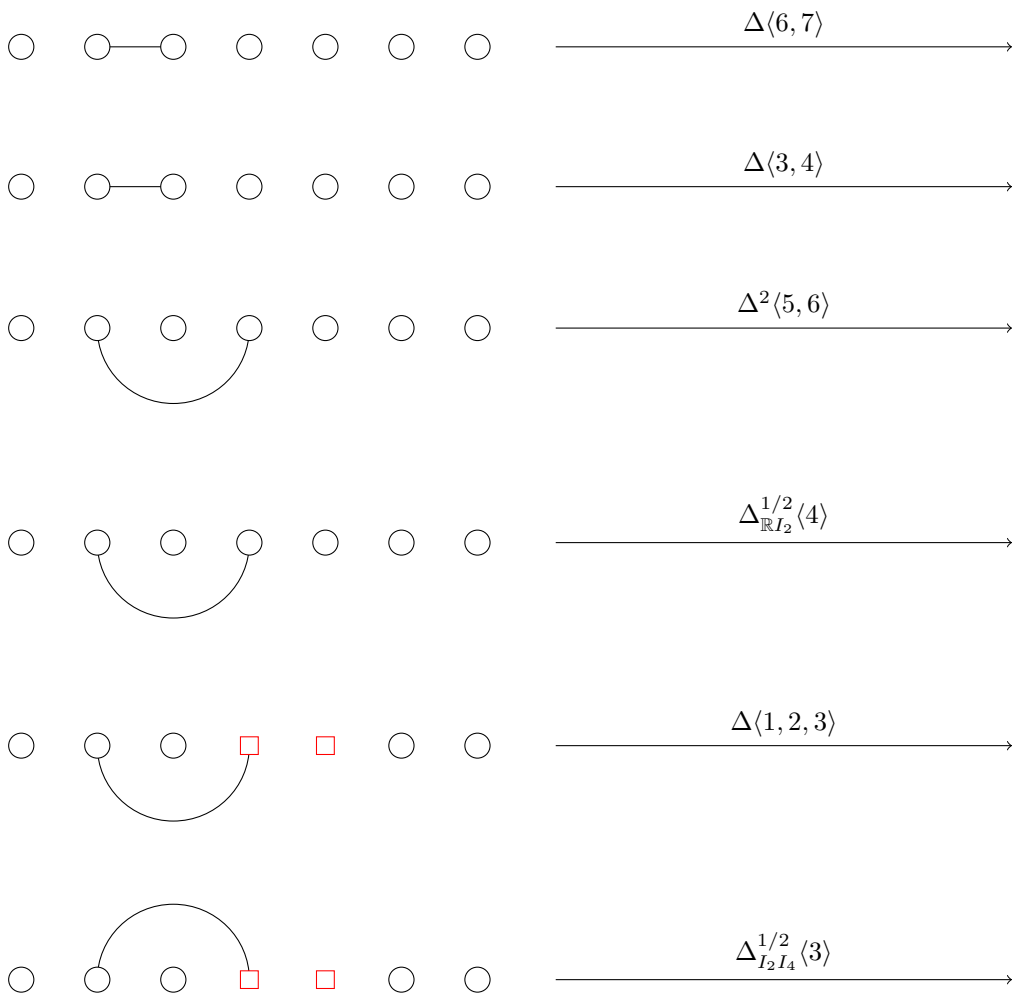
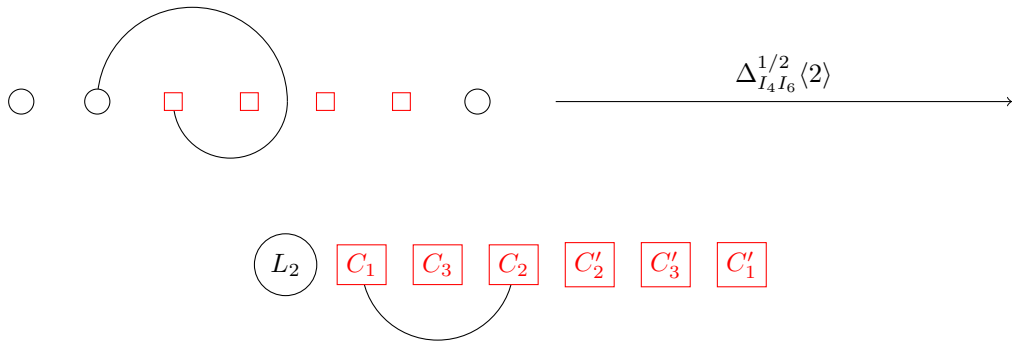
\begin{figure}[H]
	\begin{center}
		%
	\end{center}
\end{figure}

Relation:

\begin{equation*}\begin{split}
\left\{\Gamma_{C_2}  , \Gamma_{C_1} \right\}=e
\end{split}
\end{equation*}
\subsubsection{Vertex Number 9}

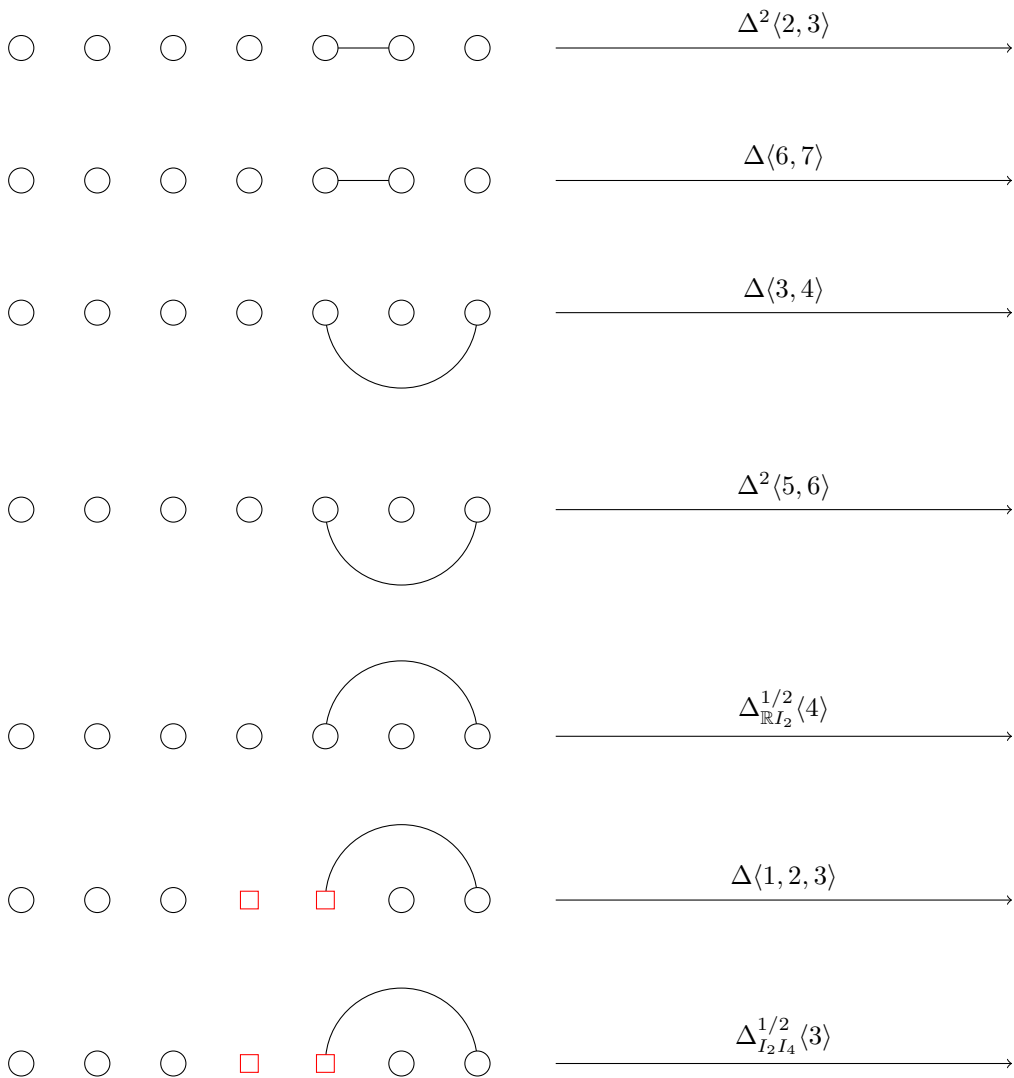
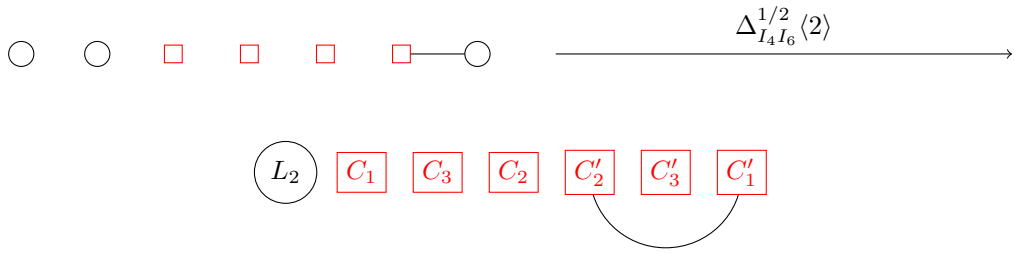
\begin{figure}[H]
	\begin{center}
		%
	\end{center}
\end{figure}

Relation:

\begin{equation*}\begin{split}
\left\{\Gamma_{C_1'}  , \Gamma_{C_2'} \right\}=e
\end{split}
\end{equation*}
\subsubsection{Vertex Number 10}

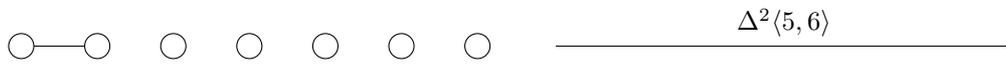
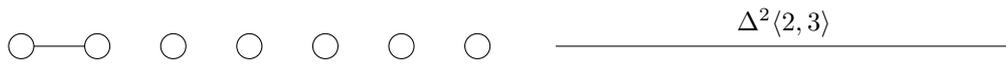
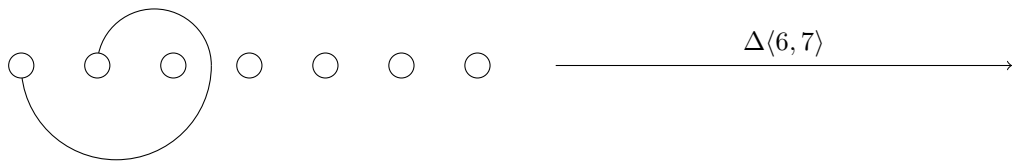
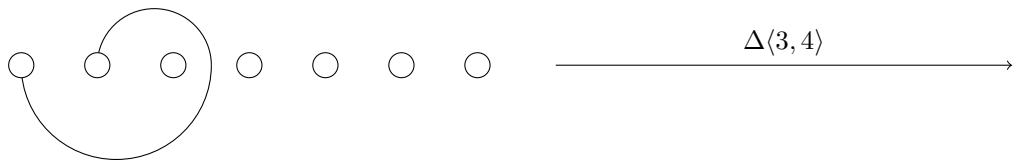
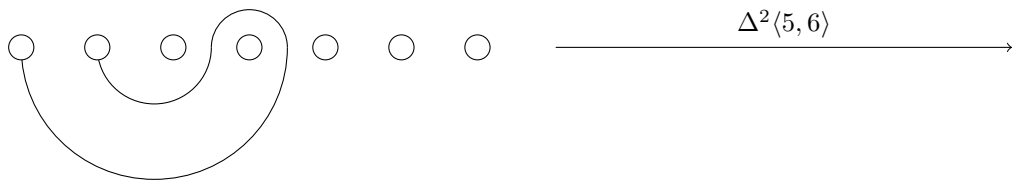
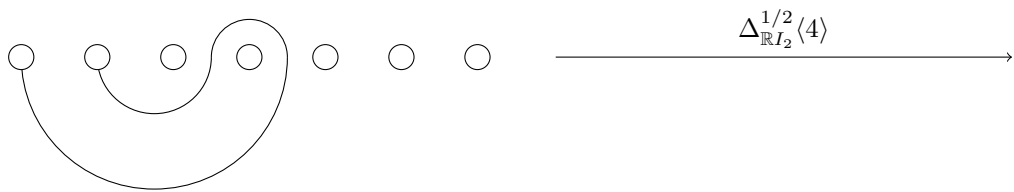
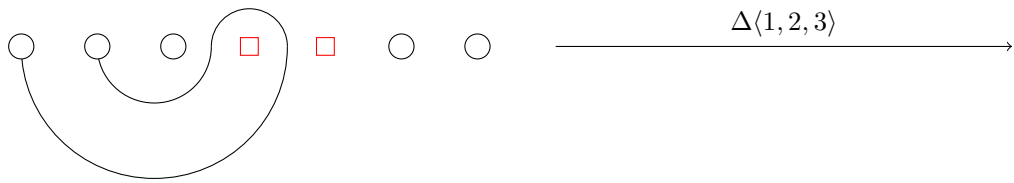
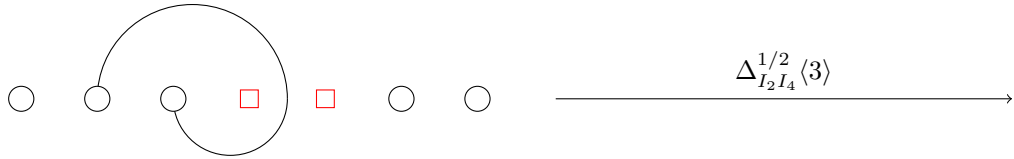
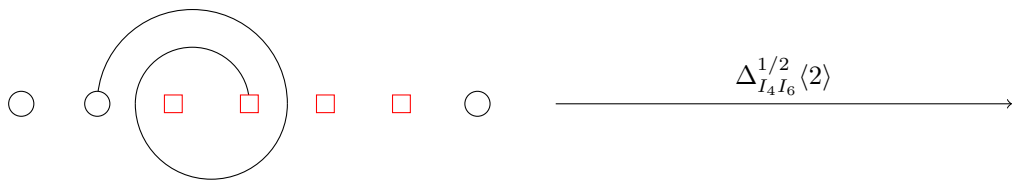
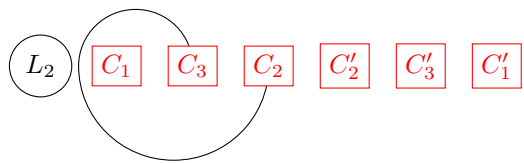
\begin{figure}[H]
	\begin{center}
		%
	\end{center}
\end{figure}

Relation:

\begin{equation*}\begin{split}
\left[\Gamma_{C_1}^{-1}\Gamma_{C_3}\Gamma_{C_1} , \Gamma_{C_2} \right]=e
\end{split}
\end{equation*}
\subsubsection{Vertex Number 11}

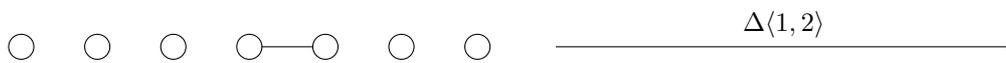
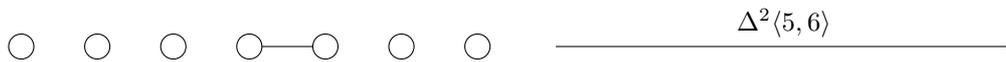
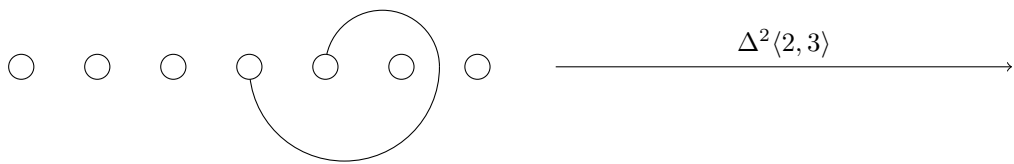
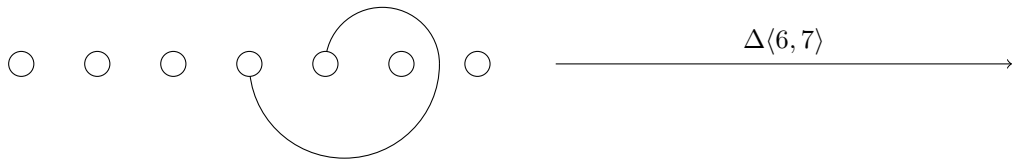
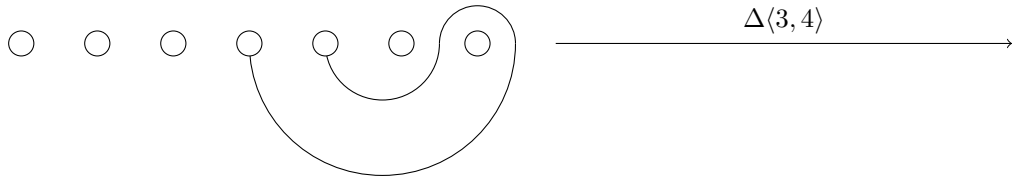
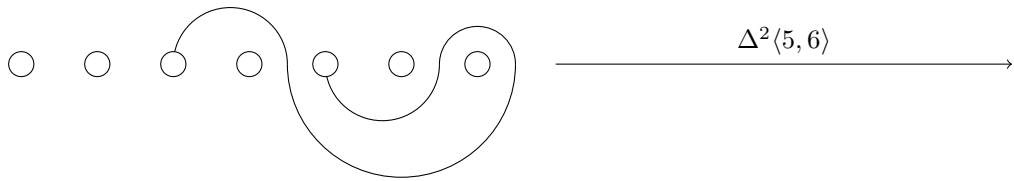
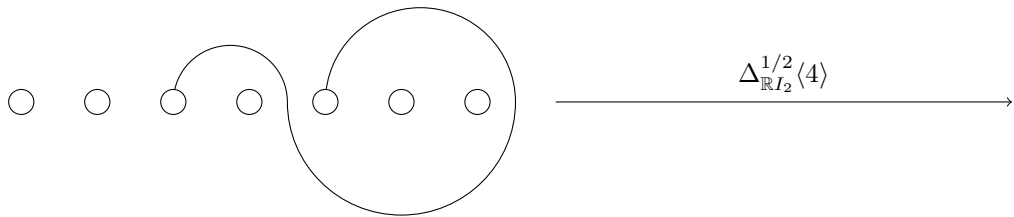
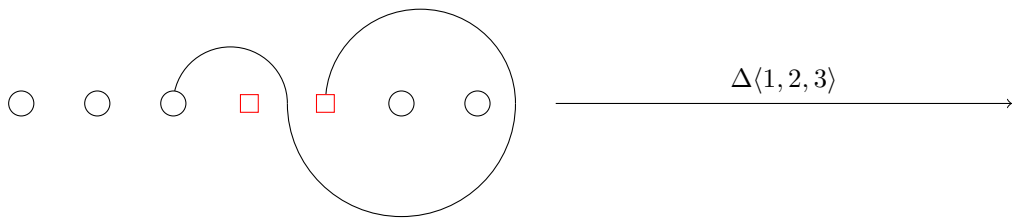
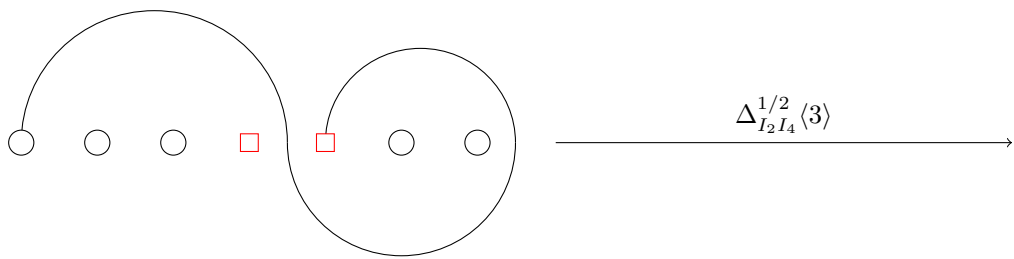
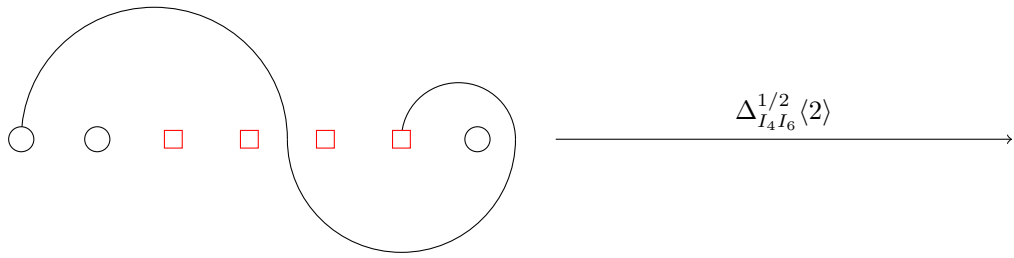
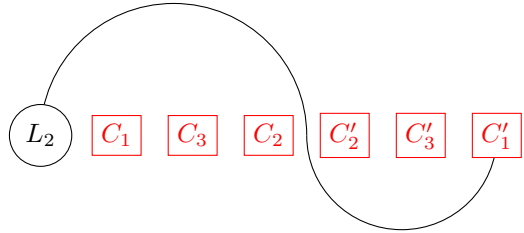
\begin{figure}[H]
	\begin{center}
		%
	\end{center}
\end{figure}

Relation:

\begin{equation*}\begin{split}
\left[\Gamma_{C_2}\Gamma_{C_3}\Gamma_{C_1}\Gamma_{L_2}\Gamma_{C_1}^{-1}\Gamma_{C_3}^{-1}\Gamma_{C_2}^{-1} , \Gamma_{C_1'} \right]=e
\end{split}
\end{equation*}
\subsubsection{Vertex Number 12}

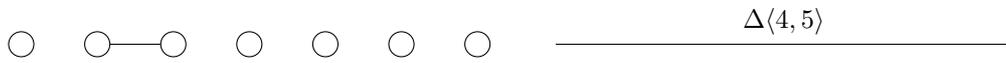
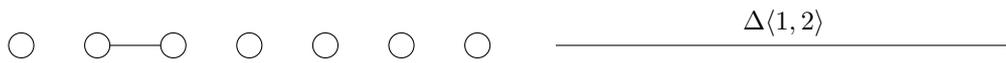
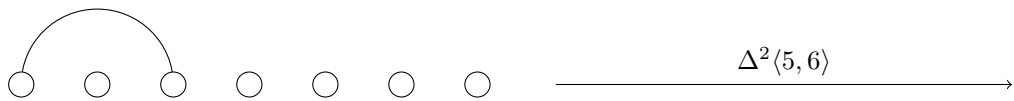
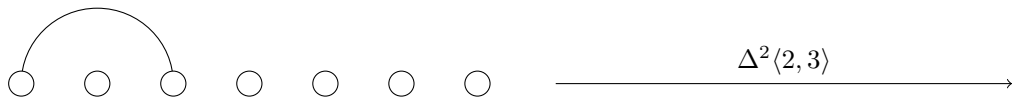
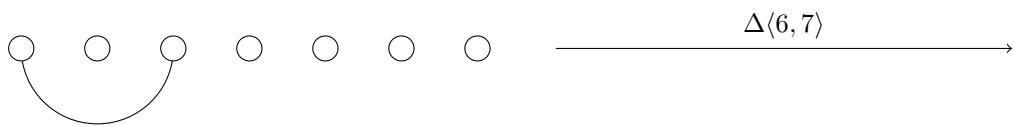
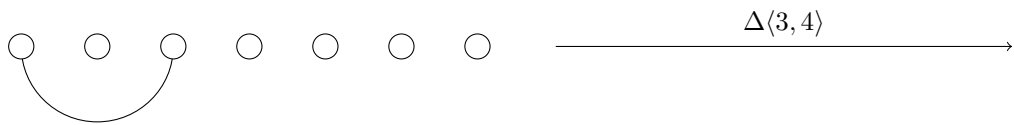
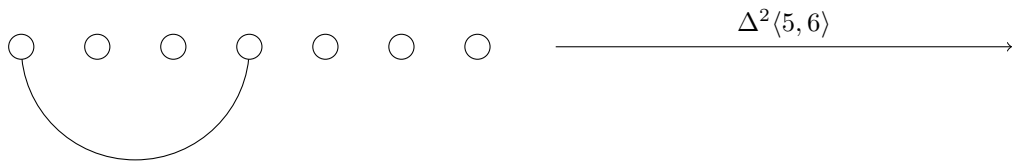
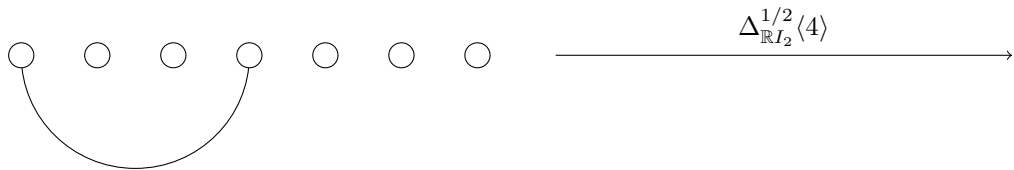
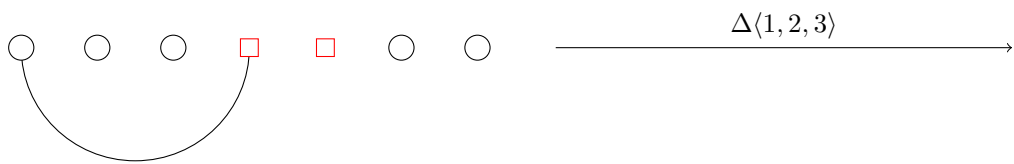
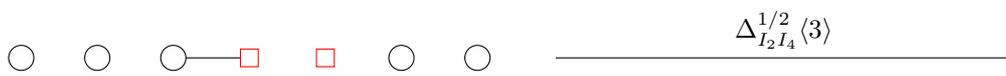
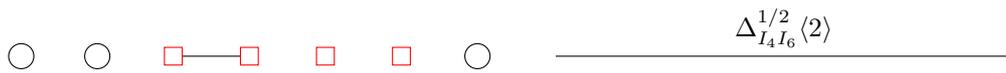
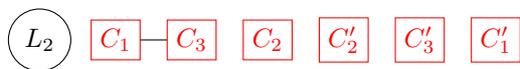
\begin{figure}[H]
	\begin{center}
		%
	\end{center}
\end{figure}

Relation:

\begin{equation*}\begin{split}
\left\{\Gamma_{C_3}  , \Gamma_{C_1} \right\}=e
\end{split}
\end{equation*}
\subsubsection{Vertex Number 13}

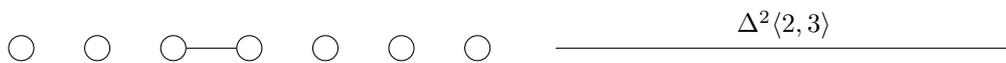
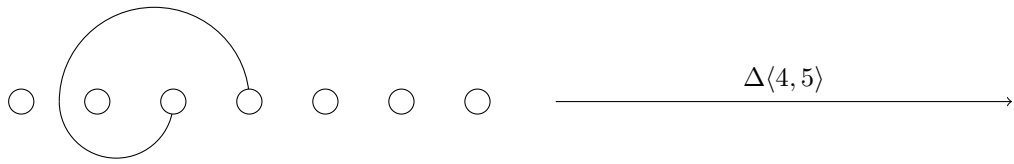
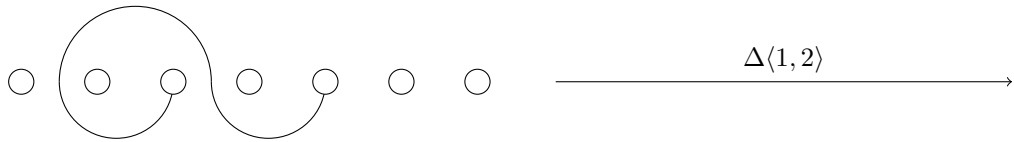
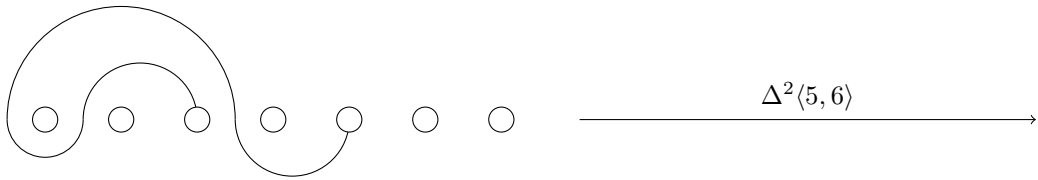
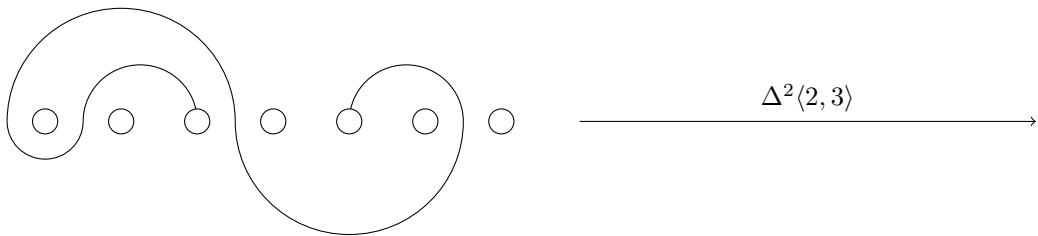
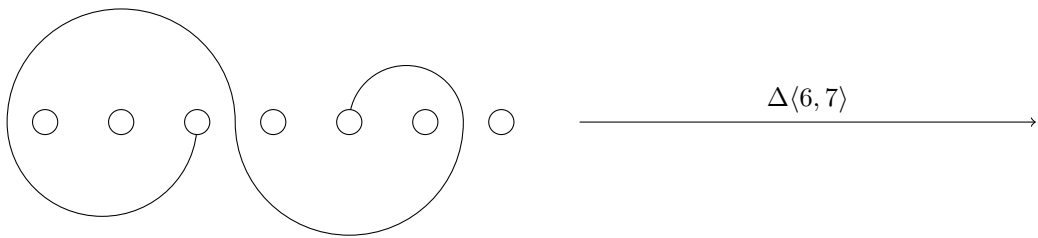
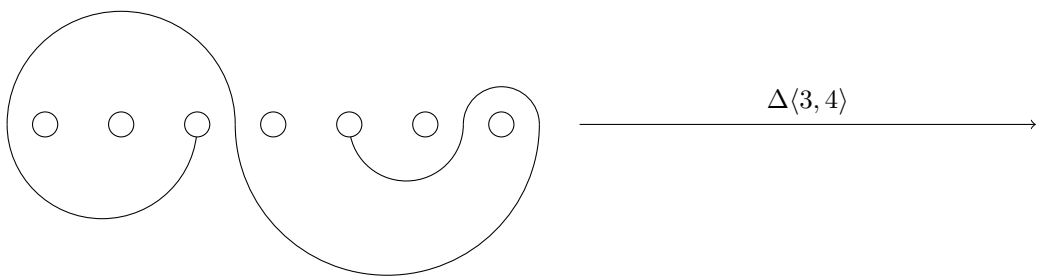
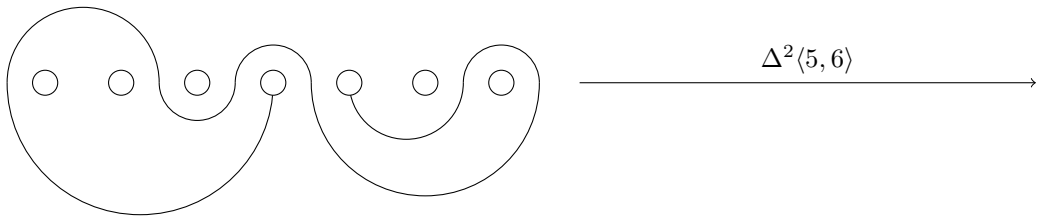
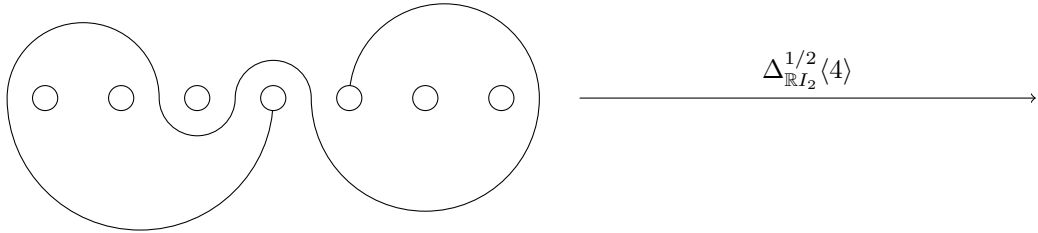
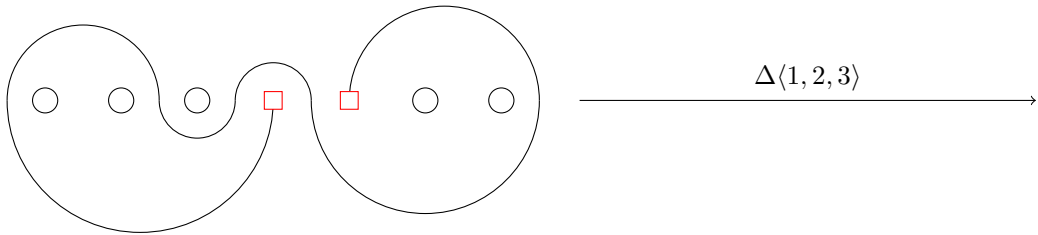
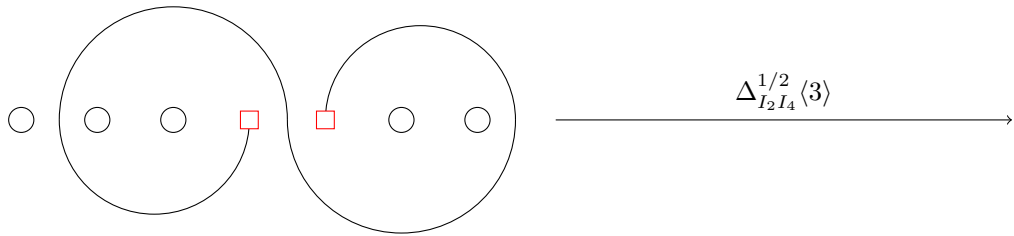
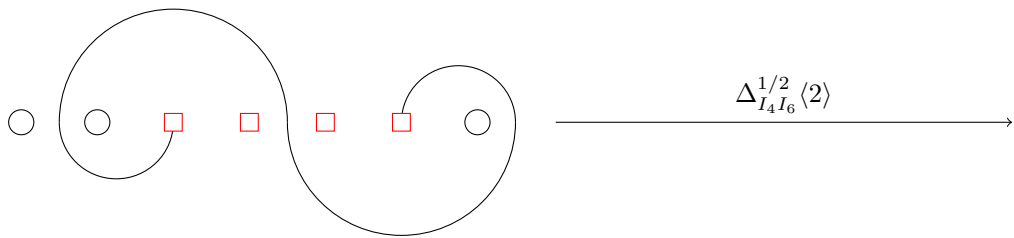
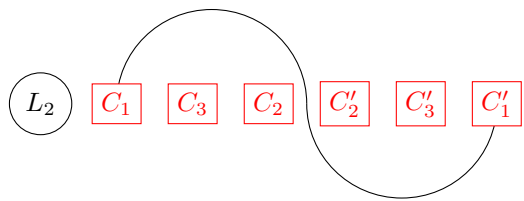
\begin{figure}[H]
	\begin{center}
		%
	\end{center}
\end{figure}

Relation:

\begin{equation*}\begin{split}
\Gamma_{C_2}\Gamma_{C_3}\Gamma_{C_1}\Gamma_{C_3}^{-1}\Gamma_{C_2}^{-1}  = \Gamma_{C_1'}
\end{split}
\end{equation*}
\subsubsection{Vertex Number 14}

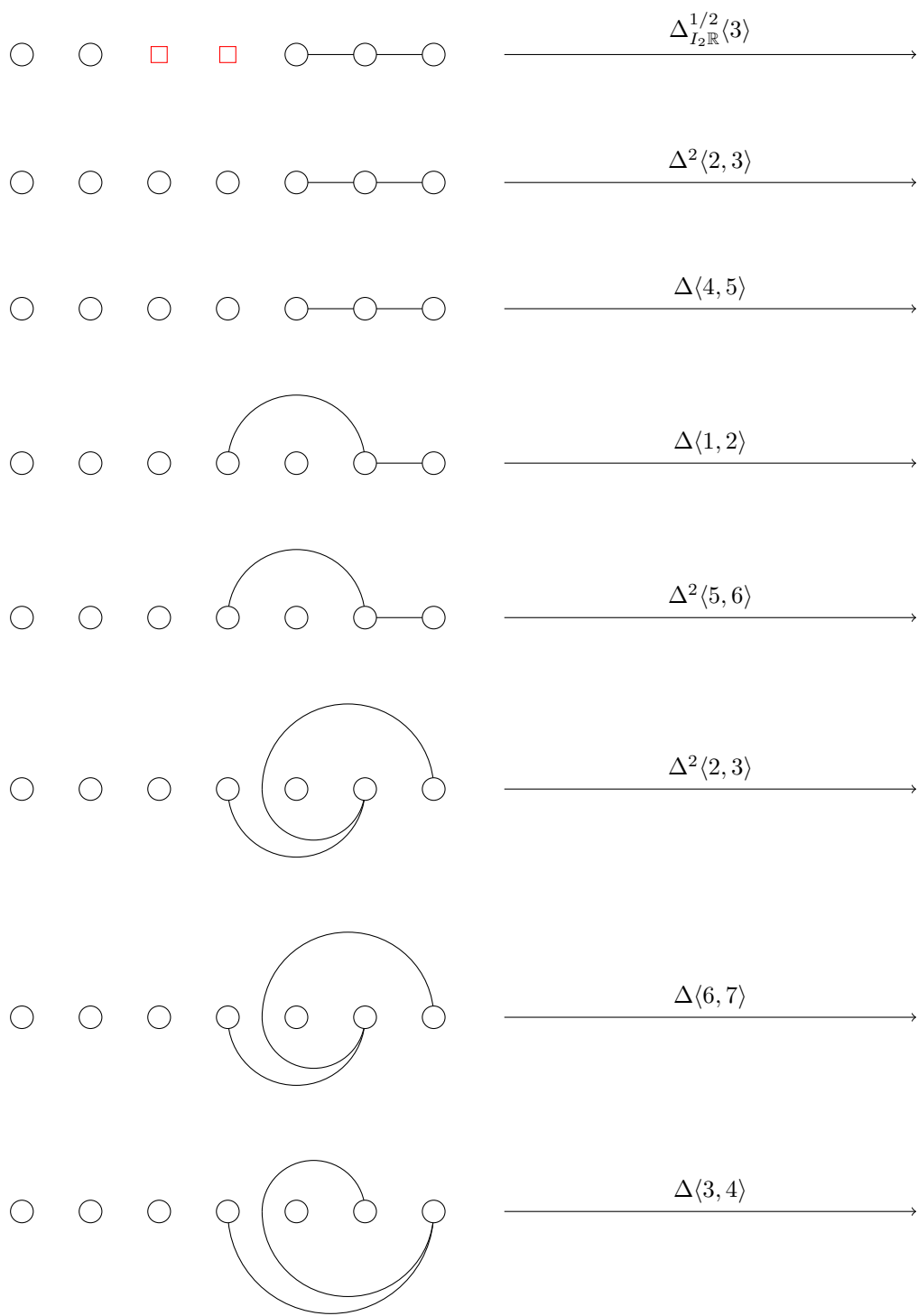
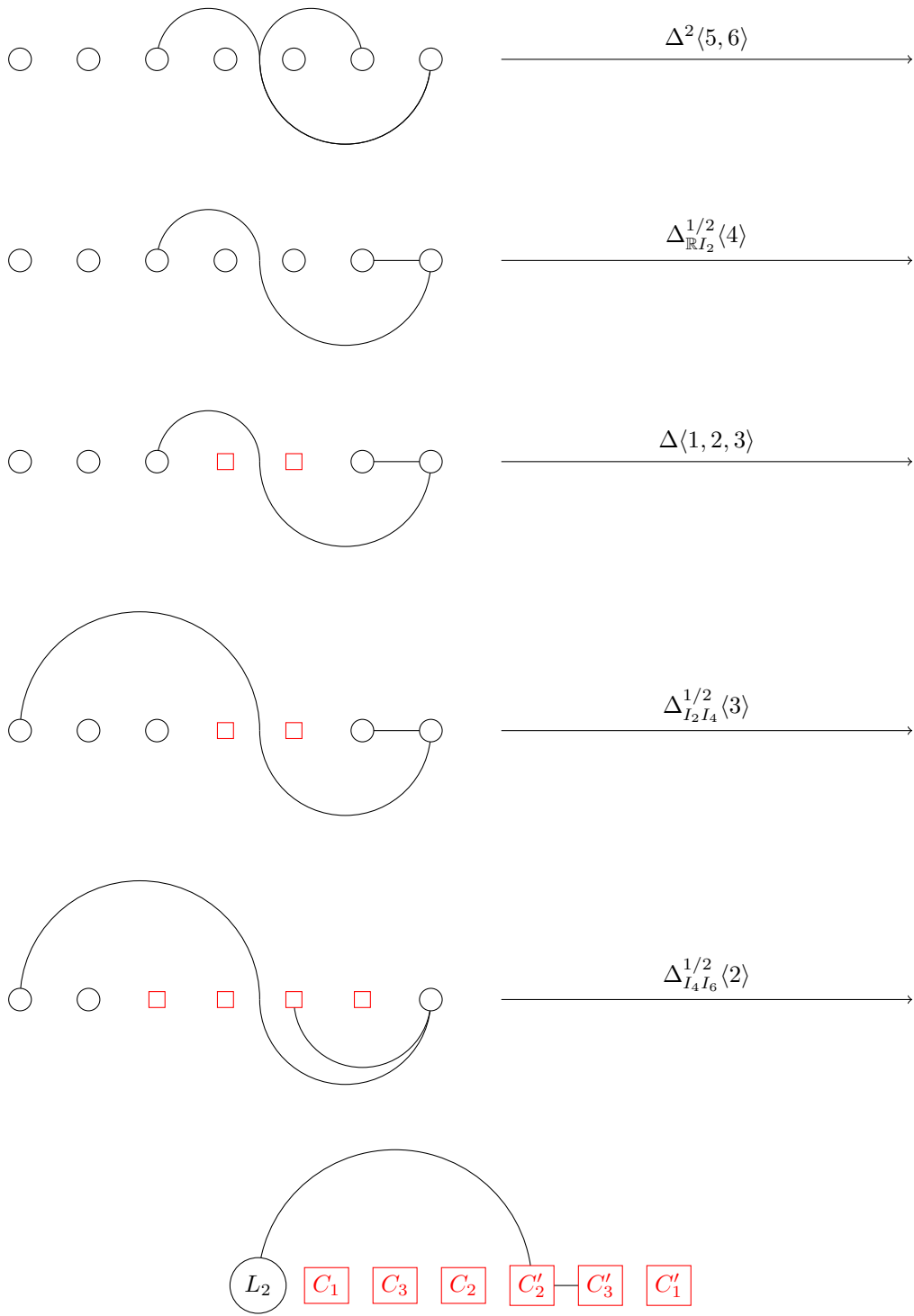
\begin{figure}[H]
	\begin{center}
		%
	\end{center}
\end{figure}

Relation:

\begin{equation*}\begin{split}
\left( \Gamma_{C_3'}  \right) \left( \Gamma_{C_2'}  \right) \left( \Gamma_{C_2}\Gamma_{C_3}\Gamma_{C_1}\Gamma_{L_2}\Gamma_{C_1}^{-1}\Gamma_{C_3}^{-1}\Gamma_{C_2}^{-1}  \right)& = \\
\left( \Gamma_{C_2'}  \right) \left( \Gamma_{C_2}\Gamma_{C_3}\Gamma_{C_1}\Gamma_{L_2}\Gamma_{C_1}^{-1}\Gamma_{C_3}^{-1}\Gamma_{C_2}^{-1}  \right) \left( \Gamma_{C_3'}  \right)& = \\
\left( \Gamma_{C_2}\Gamma_{C_3}\Gamma_{C_1}\Gamma_{L_2}\Gamma_{C_1}^{-1}\Gamma_{C_3}^{-1}\Gamma_{C_2}^{-1}  \right) \left( \Gamma_{C_3'}  \right) \left( \Gamma_{C_2'}  \right)
\end{split}
\end{equation*}
\subsubsection{Vertex Number 15}

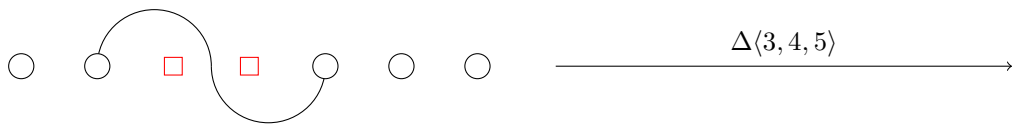
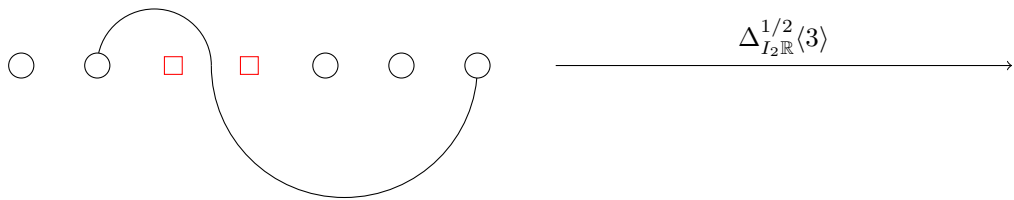
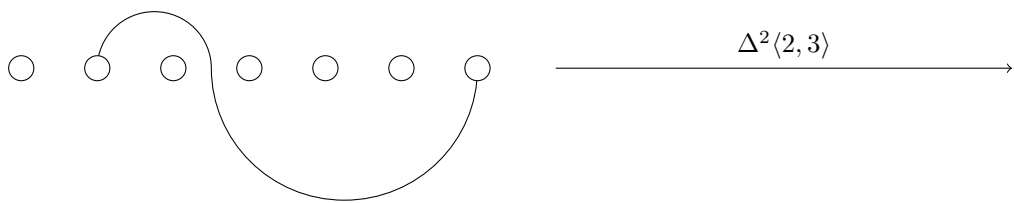
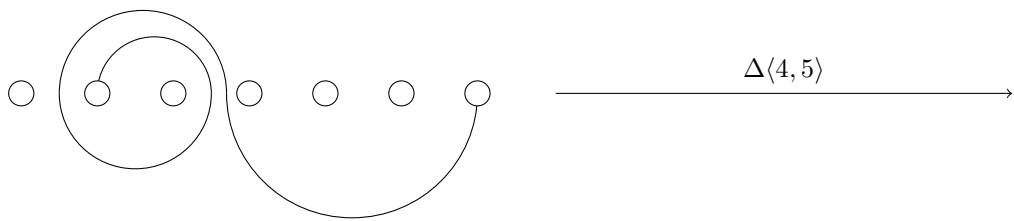
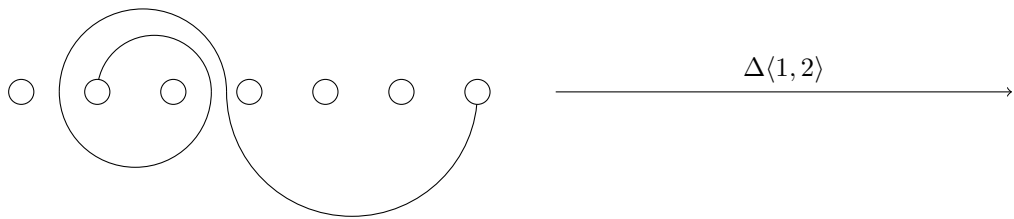
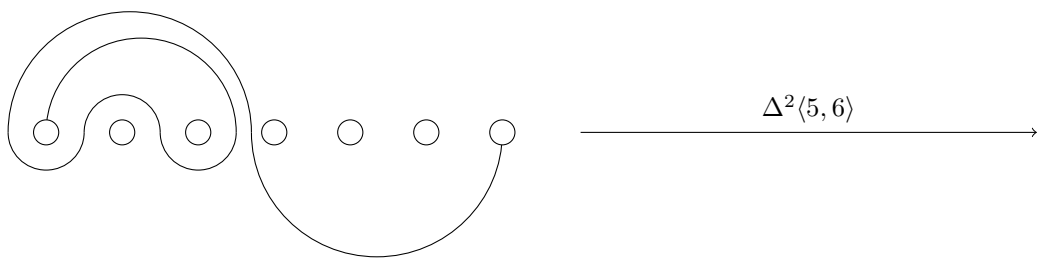
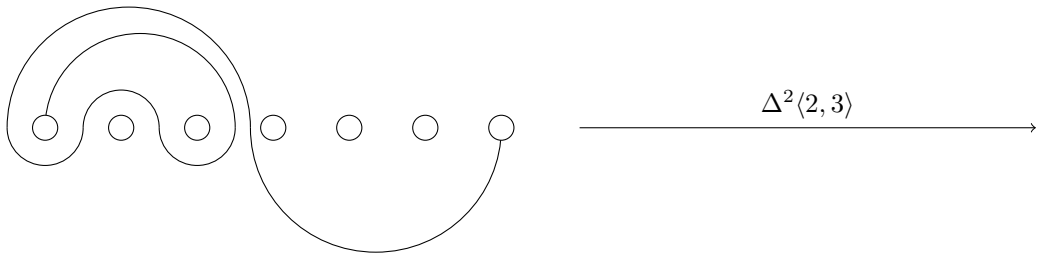
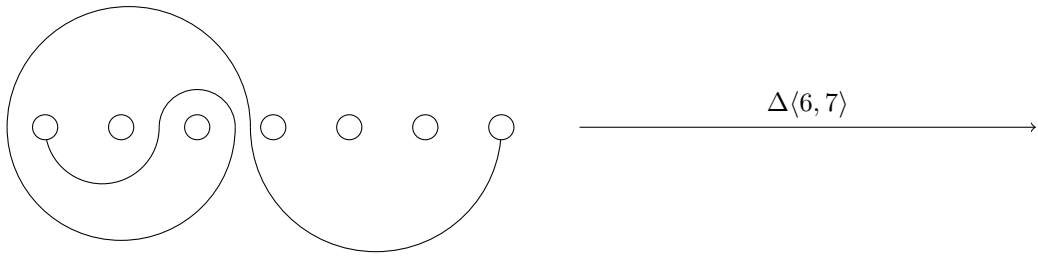
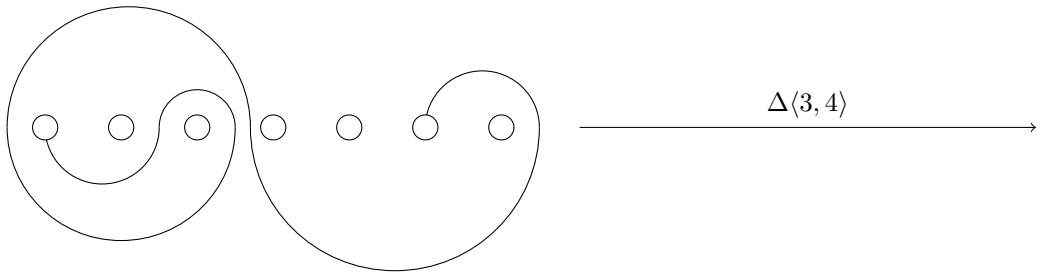
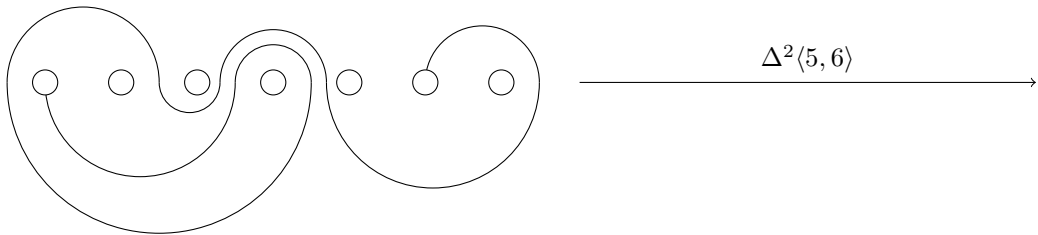
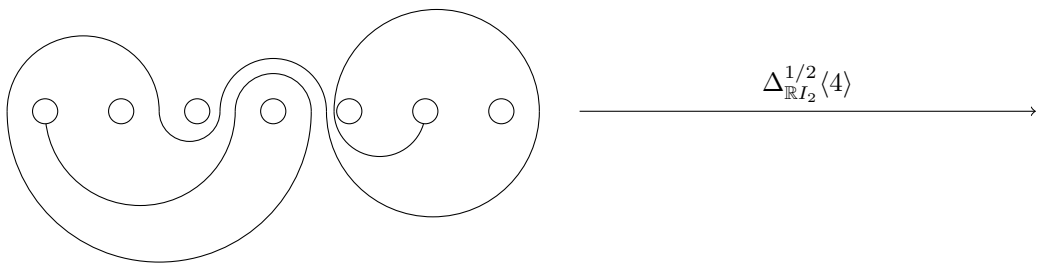
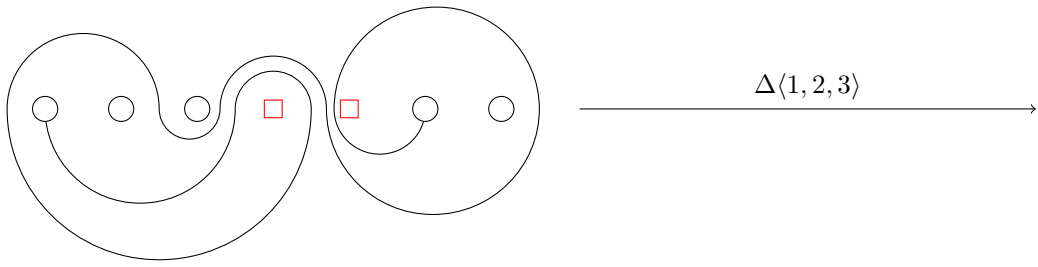
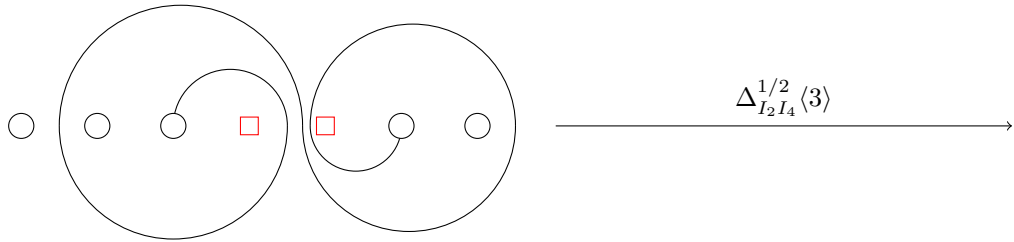
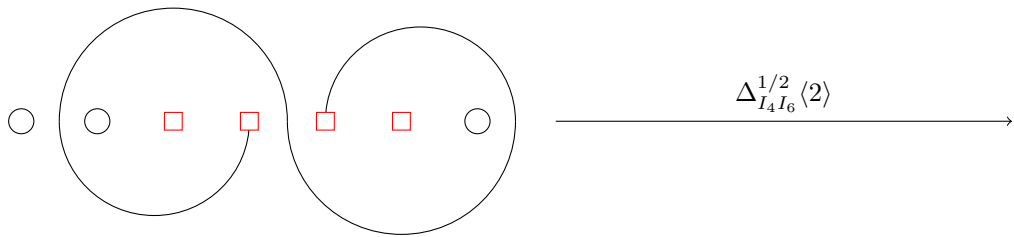
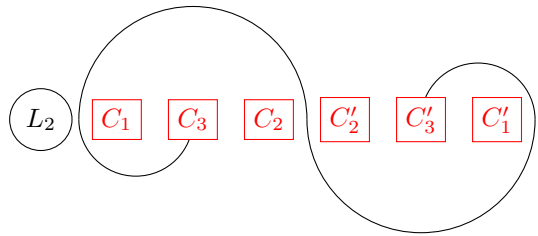
\begin{figure}[H]
	\begin{center}

	\end{center}
\end{figure}

Relation:

\begin{equation*}\begin{split}
\Gamma_{C_1'}^{-1}\Gamma_{C_2}\Gamma_{C_3}\Gamma_{C_1}\Gamma_{C_3}\Gamma_{C_1}^{-1}\Gamma_{C_3}^{-1}\Gamma_{C_2}^{-1}\Gamma_{C_1'}  = \Gamma_{C_3'}
\end{split}
\end{equation*}
\subsubsection{Vertex Number 16}

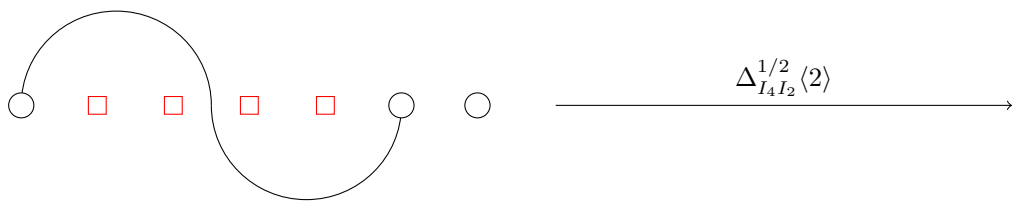
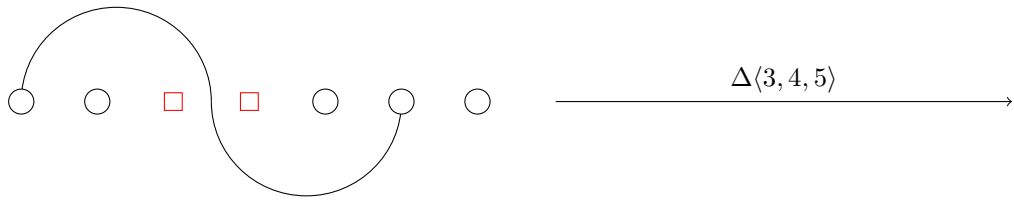
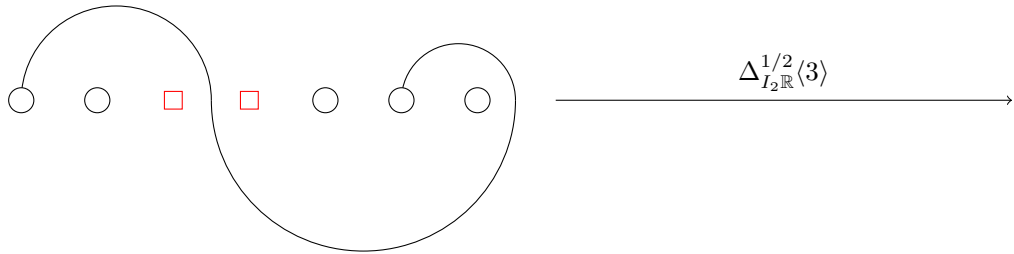
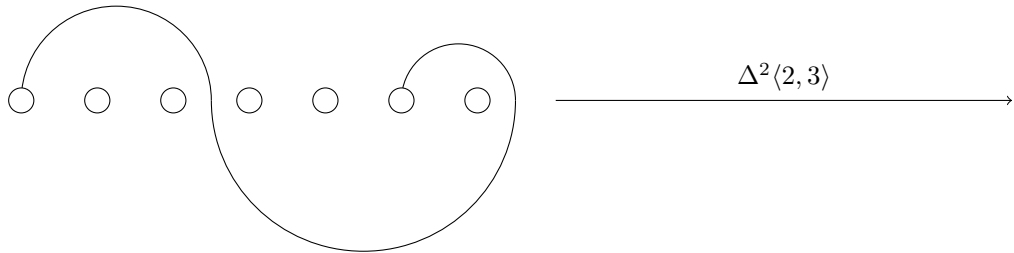
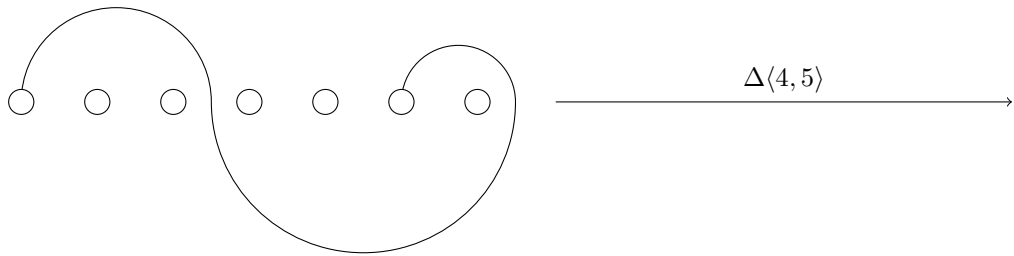
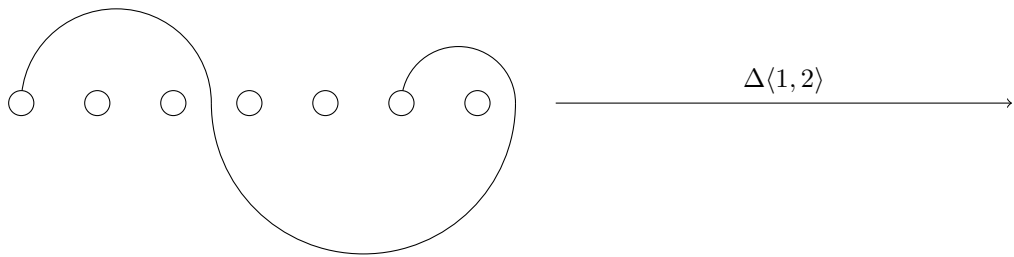
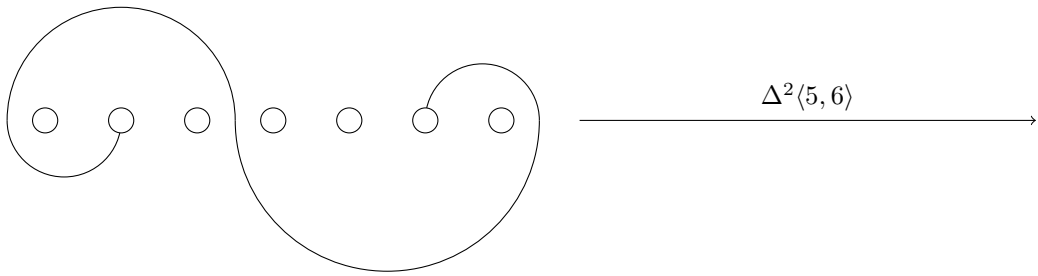
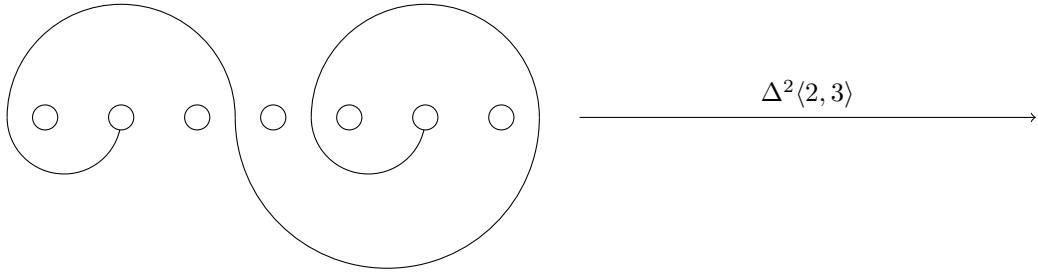
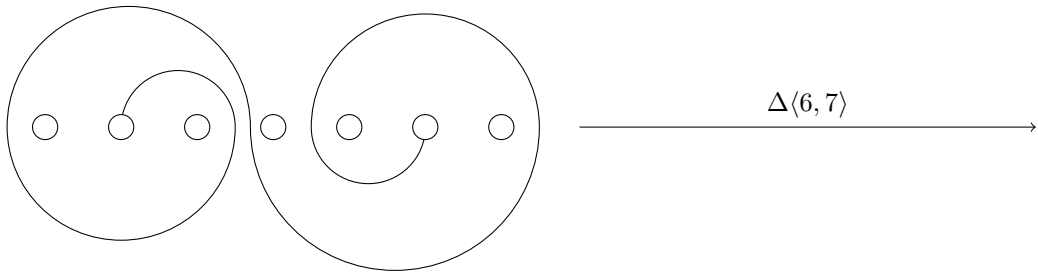
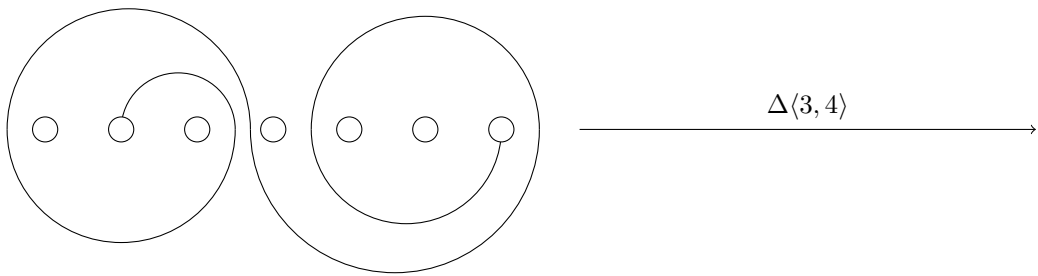
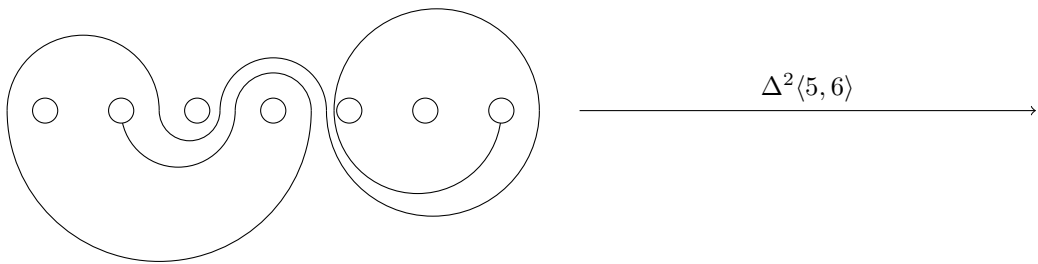
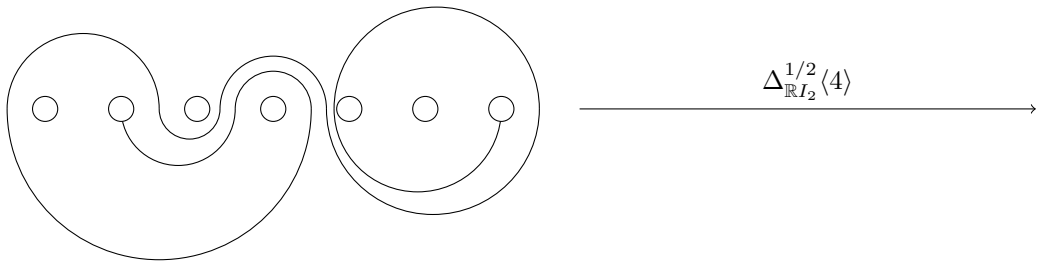
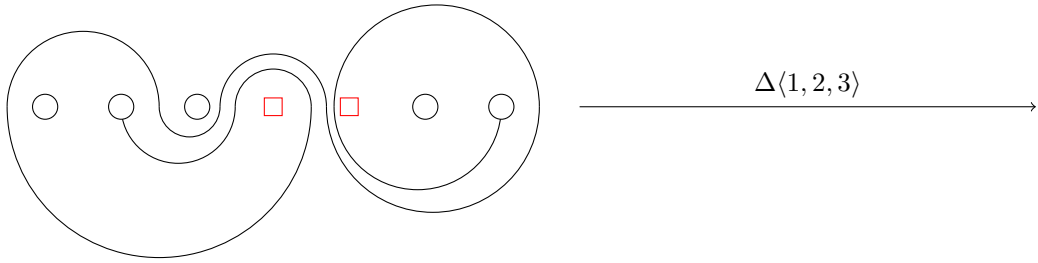
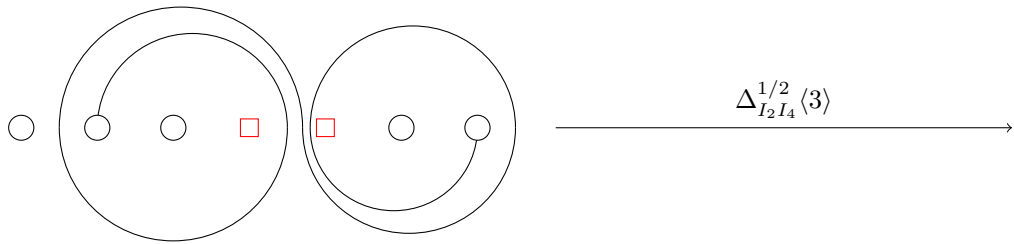
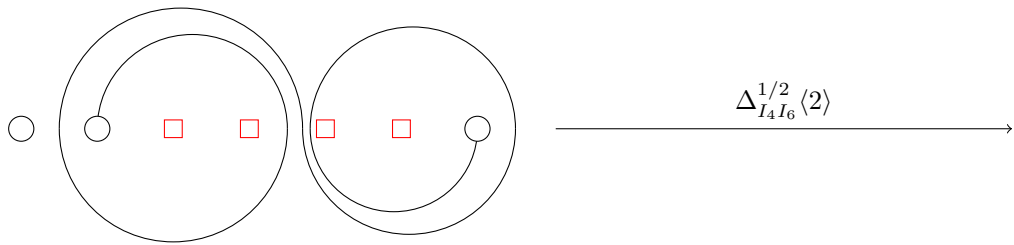
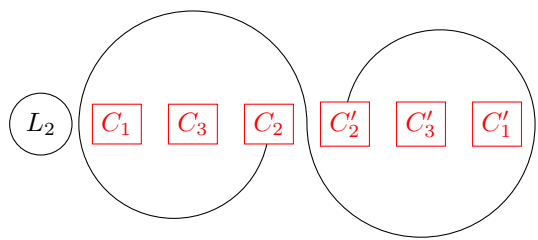
\begin{figure}[H]
	\begin{center}
		%
	\end{center}
\end{figure}

Relation:

\begin{equation*}\begin{split}
\Gamma_{C_3'}^{-1}\Gamma_{C_1'}^{-1}\Gamma_{C_2}\Gamma_{C_3}\Gamma_{C_1}\Gamma_{C_2}\Gamma_{C_1}^{-1}\Gamma_{C_3}^{-1}\Gamma_{C_2}^{-1}\Gamma_{C_1'}\Gamma_{C_3'}  = \Gamma_{C_2'}
\end{split}
\end{equation*}
\subsubsection{Raw relations and simplifications}

\begin{equation}
\begin{split}
\Gamma_{C_2}  = \Gamma_{C_2'} ,
\end{split}
\end{equation}
\begin{equation}
\begin{split}
\Gamma_{C_2'}\Gamma_{C_3}\Gamma_{C_2'}^{-1}  = \Gamma_{C_3'} ,
\end{split}
\end{equation}
\begin{equation}
\begin{split}
\left[\Gamma_{C_3} , \Gamma_{C_2}\Gamma_{C_3}\Gamma_{C_1}\Gamma_{L_2}\Gamma_{C_1}^{-1}\Gamma_{C_3}^{-1}  \right]=e,
\end{split}
\end{equation}
\begin{equation}
\begin{split}
\left[\Gamma_{C_2} , \Gamma_{C_3}\Gamma_{C_1}\Gamma_{L_2}\Gamma_{C_1}^{-1}  \right]=e,
\end{split}
\end{equation}
\begin{equation}
\begin{split}
\Gamma_{C_3'}\Gamma_{C_2'}\Gamma_{C_1}\Gamma_{C_2'}^{-1}\Gamma_{C_3'}^{-1}  = \Gamma_{C_1'} ,
\end{split}
\end{equation}
\begin{equation}
\begin{split}
\left\{\Gamma_{C_1'} , \Gamma_{C_3'}  \right\}=e,
\end{split}
\end{equation}
\begin{equation}
\begin{split}
\left[\Gamma_{C_2}\Gamma_{C_3}\Gamma_{C_1}\Gamma_{L_2}\Gamma_{C_1}^{-1}\Gamma_{C_3}^{-1}\Gamma_{C_2}^{-1} , \Gamma_{C_1}  \right]=e,
\end{split}
\end{equation}
\begin{equation}
\begin{split}
\left[\Gamma_{C_1'}\Gamma_{C_3'}\Gamma_{C_1'}^{-1} , \Gamma_{C_2'}  \right]=e,
\end{split}
\end{equation}
\begin{equation}
\begin{split}
\left\{\Gamma_{C_2} , \Gamma_{C_1}  \right\}=e,
\end{split}
\end{equation}
\begin{equation}
\begin{split}
\left\{\Gamma_{C_1'} , \Gamma_{C_2'}  \right\}=e,
\end{split}
\end{equation}
\begin{equation}
\begin{split}
\left[\Gamma_{C_1}^{-1}\Gamma_{C_3}\Gamma_{C_1} , \Gamma_{C_2}  \right]=e,
\end{split}
\end{equation}
\begin{equation}
\begin{split}
\left[\Gamma_{C_2}\Gamma_{C_3}\Gamma_{C_1}\Gamma_{L_2}\Gamma_{C_1}^{-1}\Gamma_{C_3}^{-1}\Gamma_{C_2}^{-1} , \Gamma_{C_1'}  \right]=e,
\end{split}
\end{equation}
\begin{equation}
\begin{split}
\left\{\Gamma_{C_3} , \Gamma_{C_1}  \right\}=e,
\end{split}
\end{equation}
\begin{equation}
\begin{split}
\Gamma_{C_2}\Gamma_{C_3}\Gamma_{C_1}\Gamma_{C_3}^{-1}\Gamma_{C_2}^{-1}  = \Gamma_{C_1'} ,
\end{split}
\end{equation}
\begin{equation}
\begin{split}
\left[\Gamma_{C_3'} , \Gamma_{C_2'}\Gamma_{C_2}\Gamma_{C_3}\Gamma_{C_1}\Gamma_{L_2}\Gamma_{C_1}^{-1}\Gamma_{C_3}^{-1}\Gamma_{C_2}^{-1}  \right]=e,
\end{split}
\end{equation}
\begin{equation}
\begin{split}
\left[\Gamma_{C_2'} , \Gamma_{C_2}\Gamma_{C_3}\Gamma_{C_1}\Gamma_{L_2}\Gamma_{C_1}^{-1}\Gamma_{C_3}^{-1}\Gamma_{C_2}^{-1}\Gamma_{C_3'}  \right]=e,
\end{split}
\end{equation}
\begin{equation}
\begin{split}
\Gamma_{C_1'}^{-1}\Gamma_{C_2}\Gamma_{C_3}\Gamma_{C_1}\Gamma_{C_3}\Gamma_{C_1}^{-1}\Gamma_{C_3}^{-1}\Gamma_{C_2}^{-1}\Gamma_{C_1'}  = \Gamma_{C_3'} ,
\end{split}
\end{equation}
\begin{equation}
\begin{split}
\Gamma_{C_3'}^{-1}\Gamma_{C_1'}^{-1}\Gamma_{C_2}\Gamma_{C_3}\Gamma_{C_1}\Gamma_{C_2}\Gamma_{C_1}^{-1}\Gamma_{C_3}^{-1}\Gamma_{C_2}^{-1}\Gamma_{C_1'}\Gamma_{C_3'}  = \Gamma_{C_2'} ,
\end{split}
\end{equation}
\begin{equation}
\begin{split}
\Gamma_{C_1'}\Gamma_{C_3'}\Gamma_{C_2'}\Gamma_{C_2}\Gamma_{C_3}\Gamma_{C_1}\Gamma_{L_2} =e.
\end{split}
\end{equation}

\bigskip
Here are the computations of $ \pcpt{\mathcal{B}_4} $:
\begin{lemma}\label{lemma:b4_computation}
\begin{enumerate}
	\item
	$ \pcpt{\mathcal{B}_4} $ is generated by $ \Gamma_{C_1},\Gamma_{C_2},\Gamma_{C_3} $ subject to:
	\begin{equation}\label{case4-simpl3-3}
		\begin{split}
			\left[\Gamma_{C_2} , \Gamma_{C_1}^{-1}\Gamma_{C_3}\Gamma_{C_1} \right]=e,
		\end{split}
	\end{equation}
	\begin{equation}\label{case4-simpl3-4}
		\left\{\Gamma_{C_1} , \Gamma_{C_2}  \right\}= \left\{\Gamma_{C_1} , \Gamma_{C_3}  \right\} =e.
	\end{equation}
	
	\item
	$ \pcpt{\mathcal{B}_4}  $ has Coxeter group of type $ \widetilde{C}_2 $ as a quotient.
	
	\item
	The elements $ \Gamma_{C_2} $ and $ \Gamma_{C_3} $ in $ \pcpt{\mathcal{B}_4} $ do not commute.
\end{enumerate}
\end{lemma}
\begin{proof}
Group $G:=\pcpt{\mathcal{B}_4} $ is generated by $\Gamma_{L_1},\Gamma_{C_1},\Gamma_{C_1'}, \Gamma_{C_2}, \Gamma_{C_2'}, \Gamma_{C_3}, \Gamma_{C_3'}$, subject to the following relations:
\begin{equation}\label{arrangement4_auto_calc_vert_1_rel_1}
	\begin{split}
		\Gamma_{C_2}  = \Gamma_{C_2'} ,
	\end{split}
\end{equation}
\begin{equation}\label{arrangement4_auto_calc_vert_2_rel_1}
	\begin{split}
		\Gamma_{C_2'}\Gamma_{C_3}\Gamma_{C_2'}^{-1}  = \Gamma_{C_3'},
	\end{split}
\end{equation}
\begin{equation}\label{arrangement4_auto_calc_vert_3_rel_1}
	\begin{split}
		\left[\Gamma_{C_3} , \Gamma_{C_2}\Gamma_{C_3}\Gamma_{C_1}\Gamma_{L_2}\Gamma_{C_1}^{-1}\Gamma_{C_3}^{-1}  \right]=e,
	\end{split}
\end{equation}
\begin{equation}\label{arrangement4_auto_calc_vert_3_rel_2}
	\begin{split}
		\left[\Gamma_{C_2} , \Gamma_{C_3}\Gamma_{C_1}\Gamma_{L_2}\Gamma_{C_1}^{-1}  \right]=e,
	\end{split}
\end{equation}
\begin{equation}\label{arrangement4_auto_calc_vert_4_rel_1}
	\begin{split}
		\Gamma_{C_3'}\Gamma_{C_2'}\Gamma_{C_1}\Gamma_{C_2'}^{-1}\Gamma_{C_3'}^{-1}  = \Gamma_{C_1'} ,
	\end{split}
\end{equation}
\begin{equation}\label{arrangement4_auto_calc_vert_5_rel_1}
	\begin{split}
		\left\{\Gamma_{C_1'} , \Gamma_{C_3'}  \right\}=e,
	\end{split}
\end{equation}
\begin{equation}\label{arrangement4_auto_calc_vert_6_rel_1}
	\begin{split} \left[\Gamma_{C_2}\Gamma_{C_3}\Gamma_{C_1}\Gamma_{L_2}\Gamma_{C_1}^{-1}\Gamma_{C_3}^{-1}\Gamma_{C_2}^{-1} , \Gamma_{C_1}  \right]=e,
	\end{split}
\end{equation}
\begin{equation}\label{arrangement4_auto_calc_vert_7_rel_1}
	\begin{split}
		\left[\Gamma_{C_1'}\Gamma_{C_3'}\Gamma_{C_1'}^{-1} , \Gamma_{C_2'}  \right]=e,
	\end{split}
\end{equation}
\begin{equation}\label{arrangement4_auto_calc_vert_8_rel_1}
	\begin{split}
		\left\{\Gamma_{C_2} , \Gamma_{C_1}  \right\}=e,
	\end{split}
\end{equation}
\begin{equation}\label{arrangement4_auto_calc_vert_9_rel_1}
	\begin{split}
		\left\{\Gamma_{C_1'} , \Gamma_{C_2'}  \right\}=e,
	\end{split}
\end{equation}
\begin{equation}\label{arrangement4_auto_calc_vert_10_rel_1}
	\begin{split}
		\left[\Gamma_{C_1}^{-1}\Gamma_{C_3}\Gamma_{C_1} , \Gamma_{C_2}  \right]=e,
	\end{split}
\end{equation}
\begin{equation}\label{arrangement4_auto_calc_vert_11_rel_1}
	\begin{split} \left[\Gamma_{C_2}\Gamma_{C_3}\Gamma_{C_1}\Gamma_{L_2}\Gamma_{C_1}^{-1}\Gamma_{C_3}^{-1}\Gamma_{C_2}^{-1} , \Gamma_{C_1'}  \right]=e,
	\end{split}
\end{equation}
\begin{equation}\label{arrangement4_auto_calc_vert_12_rel_1}
	\begin{split}
		\left\{\Gamma_{C_3} , \Gamma_{C_1}  \right\}=e,
	\end{split}
\end{equation}
\begin{equation}\label{arrangement4_auto_calc_vert_13_rel_1}
	\begin{split}
		\Gamma_{C_2}\Gamma_{C_3}\Gamma_{C_1}\Gamma_{C_3}^{-1}\Gamma_{C_2}^{-1}  = \Gamma_{C_1'} ,
	\end{split}
\end{equation}
\begin{equation}\label{arrangement4_auto_calc_vert_14_rel_1}
	\begin{split}
		\left[\Gamma_{C_3'} , \Gamma_{C_2'}\Gamma_{C_2}\Gamma_{C_3}\Gamma_{C_1}\Gamma_{L_2}\Gamma_{C_1}^{-1}\Gamma_{C_3}^{-1}\Gamma_{C_2}^{-1}  \right]=e,
	\end{split}
\end{equation}
\begin{equation}\label{arrangement4_auto_calc_vert_14_rel_2}
	\begin{split}
		\left[\Gamma_{C_2'} , \Gamma_{C_2}\Gamma_{C_3}\Gamma_{C_1}\Gamma_{L_2}\Gamma_{C_1}^{-1}\Gamma_{C_3}^{-1}\Gamma_{C_2}^{-1}
		\Gamma_{C_3'}  \right]=e,
	\end{split}
\end{equation}
\begin{equation}\label{arrangement4_auto_calc_vert_15_rel_1}
	\begin{split} \Gamma_{C_1'}^{-1}\Gamma_{C_2}\Gamma_{C_3}\Gamma_{C_1}\Gamma_{C_3}\Gamma_{C_1}^{-1}\Gamma_{C_3}^{-1}\Gamma_{C_2}^{-1}\Gamma_{C_1'}  = \Gamma_{C_3'} ,
	\end{split}
\end{equation}
\begin{equation}\label{arrangement4_auto_calc_vert_16_rel_1}
	\begin{split} \Gamma_{C_3'}^{-1}\Gamma_{C_1'}^{-1}\Gamma_{C_2}\Gamma_{C_3}\Gamma_{C_1}\Gamma_{C_2}\Gamma_{C_1}^{-1}\Gamma_{C_3}^{-1}\Gamma_{C_2}^{-1}\Gamma_{C_1'}\Gamma_{C_3'}  = \Gamma_{C_2'} ,
	\end{split}
\end{equation}
\begin{equation}\label{arrangement4_auto_calc_projective_rel}
	\begin{split}
		\Gamma_{C_1'}\Gamma_{C_3'}\Gamma_{C_2'}\Gamma_{C_2}\Gamma_{C_3}\Gamma_{C_1}\Gamma_{L_2} =e.
	\end{split}
\end{equation}

By \eqref{arrangement4_auto_calc_vert_1_rel_1}, relations \eqref{arrangement4_auto_calc_vert_2_rel_1} and \eqref{arrangement4_auto_calc_vert_4_rel_1} are rewritten as $\Gamma_{C_3'}=\Gamma_{C_2}\Gamma_{C_3}\Gamma_{C_2}^{-1}$ and $\Gamma_{C_1'}=\Gamma_{C_2}\Gamma_{C_3}\Gamma_{C_1}\Gamma_{C_3}^{-1}\Gamma_{C_2}^{-1}$, respectively.

Using the above expressions and \eqref{arrangement4_auto_calc_vert_12_rel_1},
relation \eqref{arrangement4_auto_calc_vert_7_rel_1} is simplified to
$$\left[\Gamma_{C_1'}\Gamma_{C_3'}\Gamma_{C_1'}^{-1} , \Gamma_{C_2'}  \right]=e \Rightarrow$$ $$\left[\Gamma_{C_2}\Gamma_{C_3}\Gamma_{C_1}\Gamma_{C_3}^{-1}\Gamma_{C_2}^{-1}\Gamma_{C_2}\Gamma_{C_3}\Gamma_{C_2}^{-1}\Gamma_{C_2}\Gamma_{C_3}
\Gamma_{C_1}^{-1}\Gamma_{C_3}^{-1}\Gamma_{C_2}^{-1} , \Gamma_{C_2}  \right]=e \Rightarrow$$
$$\left[\Gamma_{C_3}\Gamma_{C_1}\Gamma_{C_3}\Gamma_{C_1}^{-1}\Gamma_{C_3}^{-1} , \Gamma_{C_2}  \right]=e \Rightarrow$$
$$\left[\Gamma_{C_1}^{-1}\Gamma_{C_3}\Gamma_{C_1} , \Gamma_{C_2}  \right]=e.$$

Now we simplify \eqref{arrangement4_auto_calc_vert_9_rel_1}, using $\Gamma_{C_1'}=\Gamma_{C_2}\Gamma_{C_3}\Gamma_{C_1}\Gamma_{C_3}^{-1}\Gamma_{C_2}^{-1}$ again, as follows:
$$\Gamma_{C_1'}  \Gamma_{C_2'} \Gamma_{C_1'}  \Gamma_{C_2'} =   \Gamma_{C_2'} \Gamma_{C_1'}  \Gamma_{C_2'} \Gamma_{C_1'} \Rightarrow$$
$$\Gamma_{C_2}\Gamma_{C_3}\Gamma_{C_1}\Gamma_{C_3}^{-1}\Gamma_{C_2}^{-1}  \Gamma_{C_2} \Gamma_{C_2}\Gamma_{C_3}\Gamma_{C_1}\Gamma_{C_3}^{-1}\Gamma_{C_2}^{-1}  \Gamma_{C_2} =
\Gamma_{C_2} \Gamma_{C_2}\Gamma_{C_3}\Gamma_{C_1}\Gamma_{C_3}^{-1}\Gamma_{C_2}^{-1}
\Gamma_{C_2} \Gamma_{C_2}\Gamma_{C_3}\Gamma_{C_1}\Gamma_{C_3}^{-1}\Gamma_{C_2}^{-1} \Rightarrow$$
$$\Gamma_{C_3}\Gamma_{C_1}\Gamma_{C_3}^{-1} \Gamma_{C_2}\Gamma_{C_3}\Gamma_{C_1}\Gamma_{C_3}^{-1} =
\Gamma_{C_2}\Gamma_{C_3}\Gamma_{C_1}\Gamma_{C_3}^{-1} \Gamma_{C_2}\Gamma_{C_3}\Gamma_{C_1}\Gamma_{C_3}^{-1}\Gamma_{C_2}^{-1} \Rightarrow$$
$$\left\{\Gamma_{C_2} , \Gamma_{C_3} \Gamma_{C_1} \Gamma_{C_3}^{-1}  \right\}=e.$$

Relation  \eqref{arrangement4_auto_calc_vert_5_rel_1} is simplified by substituting $\Gamma_{C_3'}=\Gamma_{C_2}\Gamma_{C_3}\Gamma_{C_2}^{-1}$ and $\Gamma_{C_1'}=\Gamma_{C_2}\Gamma_{C_3}\Gamma_{C_1}\Gamma_{C_3}^{-1}\Gamma_{C_2}^{-1}$ to
$\left\{\Gamma_{C_1} , \Gamma_{C_3}  \right\}=e$. But this substitution is redundant, as we have already this relation in \eqref{arrangement4_auto_calc_vert_12_rel_1}.

Note that relations \eqref{arrangement4_auto_calc_vert_15_rel_1} and \eqref{arrangement4_auto_calc_vert_16_rel_1} are easily redundant, using $\Gamma_{C_3'}=\Gamma_{C_2}\Gamma_{C_3}\Gamma_{C_2}^{-1}$ and $\Gamma_{C_1'}=\Gamma_{C_2}\Gamma_{C_3}\Gamma_{C_1}\Gamma_{C_3}^{-1}\Gamma_{C_2}^{-1}$.

Relation \eqref{arrangement4_auto_calc_vert_14_rel_1} is simplified as follows, using $\Gamma_{C_3'}=\Gamma_{C_2}\Gamma_{C_3}\Gamma_{C_2}^{-1}$ and \eqref{arrangement4_auto_calc_vert_3_rel_2}:
$$\left[\Gamma_{C_3'},\Gamma_{C_2'}\Gamma_{C_2}\Gamma_{C_3}\Gamma_{C_1}\Gamma_{L_2}\Gamma_{C_1}^{-1}
\Gamma_{C_3}^{-1}\Gamma_{C_2}^{-1}  \right]=e \Rightarrow$$
$$\left[\Gamma_{C_2}\Gamma_{C_3}\Gamma_{C_2}^{-1},\Gamma_{C_2}\Gamma_{C_2}\Gamma_{C_3}\Gamma_{C_1}
\Gamma_{L_2}\Gamma_{C_1}^{-1}\Gamma_{C_3}^{-1}\Gamma_{C_2}^{-1}  \right]=e \Rightarrow$$ $$\left[\Gamma_{C_3},\Gamma_{C_2}(\Gamma_{C_3}\Gamma_{C_1}\Gamma_{L_2}\Gamma_{C_1}^{-1})\Gamma_{C_3}^{-1}  \right]=e \Rightarrow$$ $$\left[\Gamma_{C_3},(\Gamma_{C_3}\Gamma_{C_1}\Gamma_{L_2}\Gamma_{C_1}^{-1})\Gamma_{C_2}\Gamma_{C_3}^{-1}  \right]=e \Rightarrow$$
$$\left[\Gamma_{C_3},\Gamma_{C_1}\Gamma_{L_2}\Gamma_{C_1}^{-1}\Gamma_{C_2} \right]=e \Rightarrow$$
$$\Gamma_{C_3}\Gamma_{C_1}\Gamma_{L_2}\Gamma_{C_1}^{-1}\Gamma_{C_2}= \Gamma_{C_1}\Gamma_{L_2}\Gamma_{C_1}^{-1}\Gamma_{C_2}\Gamma_{C_3} \Rightarrow$$
$$\Gamma_{C_2}\Gamma_{C_3}\Gamma_{C_1}\Gamma_{L_2}\Gamma_{C_1}^{-1}= \Gamma_{C_1}\Gamma_{L_2}\Gamma_{C_1}^{-1}\Gamma_{C_2}\Gamma_{C_3}\Rightarrow$$
\begin{equation}\label{1temp}
	\left[\Gamma_{C_2}\Gamma_{C_3},\Gamma_{C_1}\Gamma_{L_2}\Gamma_{C_1}^{-1} \right]=e.
\end{equation}

Relation \eqref{arrangement4_auto_calc_vert_14_rel_2} is redundant after an easy simplification, because it appears already in \eqref{arrangement4_auto_calc_vert_3_rel_2}.

Now we simplify \eqref{arrangement4_auto_calc_vert_6_rel_1}, using (\ref{1temp}): $$\Gamma_{C_2}\Gamma_{C_3}(\Gamma_{C_1}\Gamma_{L_2}\Gamma_{C_1}^{-1})\Gamma_{C_3}^{-1}\Gamma_{C_2}^{-1} \Gamma_{C_1} = \Gamma_{C_1}\Gamma_{C_2}\Gamma_{C_3}(\Gamma_{C_1}\Gamma_{L_2}\Gamma_{C_1}^{-1})\Gamma_{C_3}^{-1}
\Gamma_{C_2}^{-1} \Rightarrow$$ $$(\Gamma_{C_1}\Gamma_{L_2}\Gamma_{C_1}^{-1})\Gamma_{C_2}\Gamma_{C_3}\Gamma_{C_3}^{-1}\Gamma_{C_2}^{-1} \Gamma_{C_1} = \Gamma_{C_1}(\Gamma_{C_1}\Gamma_{L_2}\Gamma_{C_1}^{-1})\Gamma_{C_2}\Gamma_{C_3}\Gamma_{C_3}^{-1}
\Gamma_{C_2}^{-1} \Rightarrow$$
$$\Gamma_{L_2}= \Gamma_{C_1}\Gamma_{L_2}\Gamma_{C_1}^{-1} \Rightarrow$$
$$\left[\Gamma_{C_1}\,\Gamma_{L_2}\right]=e.$$
This resulting commutation simplifies \eqref{arrangement4_auto_calc_vert_3_rel_2}
to $\left[\Gamma_{C_2} , \Gamma_{C_3}\Gamma_{L_2} \right]=e$, and  (\ref{1temp}) to $\left[\Gamma_{C_2}\Gamma_{C_3},\Gamma_{L_2} \right]=e$.
Also, relation \eqref{arrangement4_auto_calc_vert_3_rel_1} is easily redundant.

The last relation to simplify is \eqref{arrangement4_auto_calc_vert_11_rel_1}: $$\left[\Gamma_{C_2}\Gamma_{C_3}\Gamma_{C_1}\Gamma_{L_2}\Gamma_{C_1}^{-1}\Gamma_{C_3}^{-1}\Gamma_{C_2}^{-1}
, \Gamma_{C_1'}  \right]=e \Rightarrow$$ $$\left[\Gamma_{C_2}\Gamma_{C_3}\Gamma_{C_1}\Gamma_{L_2}\Gamma_{C_1}^{-1}\Gamma_{C_3}^{-1}\Gamma_{C_2}^{-1}
, \Gamma_{C_2}\Gamma_{C_3}\Gamma_{C_1}\Gamma_{C_3}^{-1}\Gamma_{C_2}^{-1}  \right]=e \Rightarrow$$
$$\left[\Gamma_{L_2}, \Gamma_{C_1}  \right]=e,$$
and it is redundant.

The projective relation can be rewritten as $\Gamma_{L_2} = (\Gamma_{C_2}\Gamma_{C_3}\Gamma_{C_1})^{-2}$. Now we gather all simplified relations in $G$, as follows: 	
\begin{equation}\label{arrangement4_simp1}
	\begin{split}
		\left[\Gamma_{C_2} \Gamma_{C_3}, \Gamma_{L_2}  \right]=e,
	\end{split}
\end{equation}
\begin{equation}\label{arrangement4_simp2}
	\begin{split}
		\left[\Gamma_{C_2}, \Gamma_{C_3} \Gamma_{L_2}  \right]=e,
	\end{split}
\end{equation}
\begin{equation}\label{arrangement4_simp3}
	\begin{split}
		\left[\Gamma_{C_1} , \Gamma_{L_2}  \right]=e,
	\end{split}
\end{equation}
\begin{equation}\label{arrangement4_simp4}
	\begin{split}
		\left[\Gamma_{C_2} , \Gamma_{C_1}^{-1}\Gamma_{C_3}\Gamma_{C_1} \right]=e,
	\end{split}
\end{equation}
\begin{equation}\label{arrangement4_simp5}
	\begin{split}
		\left\{\Gamma_{C_1} , \Gamma_{C_2}  \right\}=e,
	\end{split}
\end{equation}
\begin{equation}\label{arrangement4_simp6}
	\begin{split}
		\left\{\Gamma_{C_1} , \Gamma_{C_3}  \right\}=e,
	\end{split}
\end{equation}
\begin{equation}\label{arrangement4_simp7}
	\begin{split}
		\left\{\Gamma_{C_2} , \Gamma_{C_3}\Gamma_{C_1}\Gamma_{C_3}^{-1}  \right\}=e,
	\end{split}
\end{equation}
\begin{equation}\label{arrangement4_simp8}
	\begin{split}
		\Gamma_{L_2} = (\Gamma_{C_2}\Gamma_{C_3}\Gamma_{C_1})^{-2}.
	\end{split}
\end{equation}

By  \eqref{arrangement4_simp8}, \eqref{arrangement4_simp4}, and \eqref{arrangement4_simp5}, relation \eqref{arrangement4_simp3} is simplified to:
$$\Gamma_{C_1}  \Gamma_{L_2}  = \Gamma_{L_2}  \Gamma_{C_1} \Rightarrow$$
$$\Gamma_{C_1}  \Gamma_{C_2}\Gamma_{C_3}\Gamma_{C_1}\Gamma_{C_2}\Gamma_{C_3}\Gamma_{C_1}  = \Gamma_{C_2}\Gamma_{C_3}\Gamma_{C_1}\Gamma_{C_2}\Gamma_{C_3}\Gamma_{C_1}  \Gamma_{C_1} \Rightarrow$$
$$\Gamma_{C_1}  \Gamma_{C_2}\Gamma_{C_3}\Gamma_{C_1}\Gamma_{C_2} (\Gamma_{C_1}^{-1}\Gamma_{C_1}) \Gamma_{C_3}  = \Gamma_{C_2}\Gamma_{C_3}\Gamma_{C_1}\Gamma_{C_2} (\Gamma_{C_1}^{-1}\Gamma_{C_1}) \Gamma_{C_3}\Gamma_{C_1}  \Rightarrow$$
$$\Gamma_{C_1}  \Gamma_{C_2}\Gamma_{C_1}\Gamma_{C_2} \Gamma_{C_1}^{-1}\Gamma_{C_3}\Gamma_{C_1} \Gamma_{C_3}  = \Gamma_{C_2}\Gamma_{C_1}\Gamma_{C_2} \Gamma_{C_1}^{-1}\Gamma_{C_3}\Gamma_{C_1} \Gamma_{C_3}\Gamma_{C_1}  \Rightarrow$$
$$\Gamma_{C_2}  \Gamma_{C_1}\Gamma_{C_2}\Gamma_{C_1} \Gamma_{C_1}^{-1}\Gamma_{C_3}\Gamma_{C_1} \Gamma_{C_3}  = \Gamma_{C_2}\Gamma_{C_1}\Gamma_{C_2} \Gamma_{C_1}^{-1}\Gamma_{C_3}\Gamma_{C_1} \Gamma_{C_3}\Gamma_{C_1}  \Rightarrow$$
$$\Gamma_{C_3}\Gamma_{C_1} \Gamma_{C_3}  = \Gamma_{C_1}^{-1}\Gamma_{C_3}\Gamma_{C_1} \Gamma_{C_3}\Gamma_{C_1}  \Rightarrow$$
$$\Gamma_{C_1}\Gamma_{C_3}\Gamma_{C_1} \Gamma_{C_3}  = \Gamma_{C_3}\Gamma_{C_1} \Gamma_{C_3}\Gamma_{C_1},$$
and as we see, it already exists in \eqref{arrangement4_simp6}, and it is therefore redundant.
Similarly, relations \eqref{arrangement4_simp1} and \eqref{arrangement4_simp2} are easily redundant.

Now we show how \eqref{arrangement4_simp7} is redundant, using \eqref{arrangement4_simp4}, \eqref{arrangement4_simp5}, and \eqref{arrangement4_simp6}:
$$\Gamma_{C_2} \Gamma_{C_3}\Gamma_{C_1}\Gamma_{C_3}^{-1} \Gamma_{C_2} \Gamma_{C_3}\Gamma_{C_1}\Gamma_{C_3}^{-1} = \Gamma_{C_3}\Gamma_{C_1}\Gamma_{C_3}^{-1} \Gamma_{C_2} \Gamma_{C_3}\Gamma_{C_1}\Gamma_{C_3}^{-1} \Gamma_{C_2} \Rightarrow$$
$$\Gamma_{C_2} \Gamma_{C_3}\Gamma_{C_1}\Gamma_{C_3}^{-1} (\Gamma_{C_1}^{-1}\Gamma_{C_1}) \Gamma_{C_2} (\Gamma_{C_1}^{-1}\Gamma_{C_1})\Gamma_{C_3}\Gamma_{C_1}\Gamma_{C_3}^{-1} = \Gamma_{C_3}\Gamma_{C_1}\Gamma_{C_3}^{-1} \Gamma_{C_2} (\Gamma_{C_1}^{-1}\Gamma_{C_1})\Gamma_{C_3}\Gamma_{C_1}\Gamma_{C_3}^{-1} \Gamma_{C_2} \Rightarrow$$
$$\Gamma_{C_2} (\Gamma_{C_1}^{-1}\Gamma_{C_3}^{-1} \Gamma_{C_1}\Gamma_{C_3})\Gamma_{C_1} \Gamma_{C_2} \Gamma_{C_1}^{-1}(\Gamma_{C_3}^{-1}\Gamma_{C_1}\Gamma_{C_3}\Gamma_{C_1}) = (\Gamma_{C_1}^{-1}\Gamma_{C_3}^{-1}\Gamma_{C_1}\Gamma_{C_3}\Gamma_{C_1}) \Gamma_{C_2} \Gamma_{C_1}^{-1}\Gamma_{C_1}\Gamma_{C_3}\Gamma_{C_1}\Gamma_{C_3}^{-1} \Gamma_{C_2} \Rightarrow$$
$$\Gamma_{C_2} \Gamma_{C_1}^{-1}\Gamma_{C_3}^{-1} \Gamma_{C_1}\Gamma_{C_1} \Gamma_{C_2}\Gamma_{C_3}\Gamma_{C_1} = \Gamma_{C_1}^{-1}\Gamma_{C_3}^{-1}\Gamma_{C_1}\Gamma_{C_3}\Gamma_{C_1} \Gamma_{C_2} \Gamma_{C_1}^{-1}\Gamma_{C_1}\Gamma_{C_3}\Gamma_{C_1}\Gamma_{C_3}^{-1} \Gamma_{C_2} \Rightarrow$$
$$\Gamma_{C_2} \Gamma_{C_1} \Gamma_{C_2} \Gamma_{C_3}\Gamma_{C_1} = \Gamma_{C_1} \Gamma_{C_2} \Gamma_{C_1}^{-1}\Gamma_{C_3}\Gamma_{C_1}\Gamma_{C_3}\Gamma_{C_1}\Gamma_{C_3}^{-1} \Gamma_{C_2} \Rightarrow$$
$$\Gamma_{C_2}^{-1}\Gamma_{C_1}^{-1}\Gamma_{C_2} \Gamma_{C_1} \Gamma_{C_2} \Gamma_{C_3}\Gamma_{C_1} = \Gamma_{C_1}^{-1}\Gamma_{C_3}\Gamma_{C_1}\Gamma_{C_3}\Gamma_{C_1}\Gamma_{C_3}^{-1} \Gamma_{C_2}\Rightarrow $$
$$\Gamma_{C_1} \Gamma_{C_2} \Gamma_{C_1}^{-1} \Gamma_{C_3}\Gamma_{C_1} = \Gamma_{C_1}^{-1}\Gamma_{C_1}\Gamma_{C_3}\Gamma_{C_1}\Gamma_{C_3}\Gamma_{C_3}^{-1} \Gamma_{C_2} \Rightarrow$$
$$\Gamma_{C_1} \Gamma_{C_2} \Gamma_{C_1}^{-1} \Gamma_{C_3}\Gamma_{C_1} = \
\Gamma_{C_3}\Gamma_{C_1} \Gamma_{C_2} \Rightarrow$$
$$\Gamma_{C_3}\Gamma_{C_1} \Gamma_{C_2} \Gamma_{C_1}^{-1} \Gamma_{C_1} = \
\Gamma_{C_3}\Gamma_{C_1} \Gamma_{C_2} \Rightarrow$$
$$\Gamma_{C_2} = \Gamma_{C_2}.$$

Therefore, group $ \pi_1(\mathbb{CP}^2-\mathcal{B}_4)$ is generated by $\Gamma_{C_1}, \Gamma_{C_2}$, and $\Gamma_{C_3}$ and admits the following relations:
\begin{equation}\label{arrangement4_simp9}
	\begin{split}
		\left[\Gamma_{C_2} , \Gamma_{C_1}^{-1}\Gamma_{C_3}\Gamma_{C_1} \right]=e,
	\end{split}
\end{equation}
\begin{equation}\label{arrangement4_simp10}
	\begin{split}
		\left\{\Gamma_{C_1} , \Gamma_{C_2}  \right\}=e,
	\end{split}
\end{equation}
\begin{equation}\label{arrangement4_simp11}
	\begin{split}
		\left\{\Gamma_{C_1} , \Gamma_{C_3}  \right\}=e.
	\end{split}
\end{equation}

The quotient $G_4 := \pcpt{\mathcal{B}_4}/\langle \Gamma_X^2 \; | \; X \subseteq \mathcal{B}_4 \rangle $ is then generated by $\Gamma_{C_1}, \Gamma_{C_2}, \Gamma_{C_3}$ and has the following relations:
\begin{equation}\label{arrangement4_simpl12}
	\Gamma_{C_1}^2 = \Gamma_{C_2}^2 = \Gamma_{C_3}^2 = e,
\end{equation}
\begin{equation}\label{arrangement4_simp13}
	\begin{split}
		\left[\Gamma_{C_2} , \Gamma_{C_1}\Gamma_{C_3}\Gamma_{C_1} \right]=e,
	\end{split}
\end{equation}
\begin{equation}\label{arrangement4_simp14}
	\begin{split}
		\left\{\Gamma_{C_1} , \Gamma_{C_2}  \right\}=e,
	\end{split}
\end{equation}
\begin{equation}\label{arrangement4_simp15}
	\begin{split}
		\left\{\Gamma_{C_1} , \Gamma_{C_3}  \right\}=e.
	\end{split}
\end{equation}

Obviously, $ G_4 $ is also generated by $ \Gamma_{C_1}, \Gamma_{C_2} $, and $ \widetilde{\Gamma}:=\Gamma_{C_1}\Gamma_{C_3}\Gamma_{C_1} $.
Examining the relations between $ \Gamma_{C_1}, \Gamma_{C_2} $, and $ \widetilde{\Gamma} $, we see that it is a set of Coxeter generators, and that $ G_4 $ is a Coxeter group of type $\widetilde{C}_2$.
Computing the commutator $ [\Gamma_{C_2}, \Gamma_{C_3}] $ in $ \widetilde{C}_2 $, we get
$$ [\Gamma_{C_2}, \Gamma_{C_3}] = \begin{pmatrix} 1 & 2 \\ 1^{-2} & 2 \end{pmatrix}, $$ so $ \Gamma_{C_2} $
and $ \Gamma_{C_3} $ do not commute in $ G_4 $, thus, they do not commute in $ \pcpt{\mathcal{B}_4} $, either.

\end{proof}

\subsection{$\mathcal{B}_3$ and $\mathcal{B}_4$ are Zariski pair with non-isomorphic groups}

As a corollary, we get the following theorem:
\begin{thm}\label{thm:tokunaga_pair2}
$ \mathcal{B}_3 $ and $ \mathcal{B}_4 $ form a Zariski pair with non-isomorphic fundamental groups.
\end{thm}

\begin{proof}
Any homeomorphism from $ (\cpt, \mathcal{B}_3) $ to $ (\cpt, \mathcal{B}_4) $ must send $ C_2 $ to either $ C_2 $ or $ C_3 $ (because their combinatorics in $ \mathcal{B}_4 $ are identical).
In $ G_3 $ we have $ \left[ \Gamma_{C_2}, \Gamma_{C_3} \right] = e  $, while in $ G_4$ the generators $ \Gamma_{C_2} $ and $ \Gamma_{C_3} $ do not commute, so such an isomorphism cannot exist.

\end{proof}

\newpage

\section{Computations of $ \pcpt{\mathcal{B}_5} $ and $ \pcpt{\mathcal{B}_6}$}
\label{section:braids_computation_case1}
\subsection{$ \pcpt{\mathcal{B}_5} $}
\begin{figure}[H]
	\begin{center}

	\end{center}
\end{figure}

{Monodromy table}

\begin{tabular}{| c | c | c | c |}
	\hline \hline
	Vertex number& Vertex description& Skeleton& Diffeomorphism\\ \hline \hline
	1& $C_1$ branch& $ \langle 3 - 4 \rangle $& $ \Delta_{I_{2}I_{4}}^{1/2}\langle 3 \rangle $ \\ \hline
	2& Node between $C_1$, $L_3$ and $L_2$& $ \langle 1 - 2 - 3 \rangle $& $ \Delta \langle 1,2,3 \rangle $ \\ \hline
	3& $C_2$ branch& $ \langle 3 - 4 \rangle $& $ \Delta_{\mathbb{R}I_{2}}^{1/2}\langle 3 \rangle $ \\ \hline
	4& Node between $C_2$ and $L_3$& $ \langle 2 - 3 \rangle $& $ \Delta\langle 2, 3 \rangle $ \\ \hline
	5& Tangency between $L_2$ and $C_2$& $ \langle 4 - 5 \rangle $& $ \Delta^{2}\langle 4, 5 \rangle $ \\ \hline
	6& Node between $C_1$ and $L_2$& $ \langle 5 - 6 \rangle $& $ \Delta\langle 5, 6 \rangle $ \\ \hline
	7& Tangency between $C_1$ and $C_2$& $ \langle 1 - 2 \rangle $& $ \Delta^{2}\langle 1, 2 \rangle $ \\ \hline
	8& Tangency between $C_1$ and $C_2$& $ \langle 4 - 5 \rangle $& $ \Delta^{2}\langle 4, 5 \rangle $ \\ \hline
	9& Node between $L_1$ and $L_2$& $ \langle 6 - 7 \rangle $& $ \Delta\langle 6, 7 \rangle $ \\ \hline
	10& Node between $L_1$ and $C_1$& $ \langle 5 - 6 \rangle $& $ \Delta\langle 5, 6 \rangle $ \\ \hline
	11& Tangency between $L_1$ and $C_2$& $ \langle 4 - 5 \rangle $& $ \Delta^{2}\langle 4, 5 \rangle $ \\ \hline
	12& Node between $L_3$ and $C_2$& $ \langle 2 - 3 \rangle $& $ \Delta\langle 2, 3 \rangle $ \\ \hline
	13& $C_2$ branch& $ \langle 3 - 4 \rangle $& $ \Delta_{I_{2}\mathbb{R}}^{1/2}\langle 3 \rangle $ \\ \hline
	14& Node between $L_1$, $L_3$ and $C_1$& $ \langle 1 - 2 - 3 \rangle $& $ \Delta\langle 1, 2 , 3\rangle $ \\ \hline
	15& $C_1$ branch& $ \langle 3 - 4 \rangle $& $ \Delta_{I_{4}I_{2}}^{1/2}\langle 3 \rangle $
	\\ \hline \hline
\end{tabular}

\subsubsection{Vertex Number 1}


\begin{figure}[H]
	\begin{center}
		%
	\end{center}
\end{figure}

Relation:

\begin{equation*}\begin{split}
\Gamma_{C_1}  = \Gamma_{C_1'}
\end{split}
\end{equation*}
\subsubsection{Vertex Number 2}

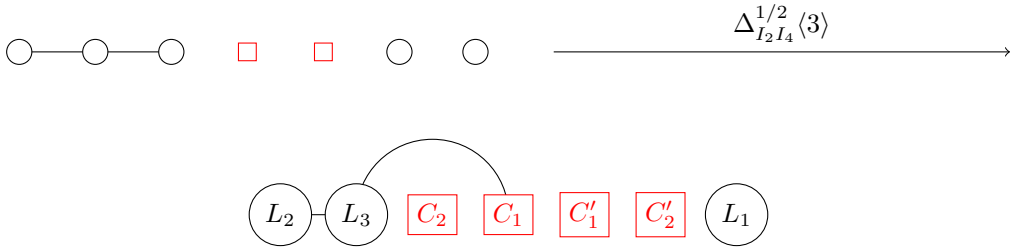
\begin{figure}[H]
	\begin{center}
		%
	\end{center}
\end{figure}

Relation:

\begin{equation*}\begin{split}
\left( \Gamma_{C_2}^{-1}\Gamma_{C_1}\Gamma_{C_2}  \right) \left( \Gamma_{L_3}  \right) \left( \Gamma_{L_2}  \right)& = \\
\left( \Gamma_{L_3}  \right) \left( \Gamma_{L_2}  \right) \left( \Gamma_{C_2}^{-1}\Gamma_{C_1}\Gamma_{C_2}  \right)& = \\
\left( \Gamma_{L_2}  \right) \left( \Gamma_{C_2}^{-1}\Gamma_{C_1}\Gamma_{C_2}  \right) \left( \Gamma_{L_3}  \right)
\end{split}
\end{equation*}
\subsubsection{Vertex Number 3}

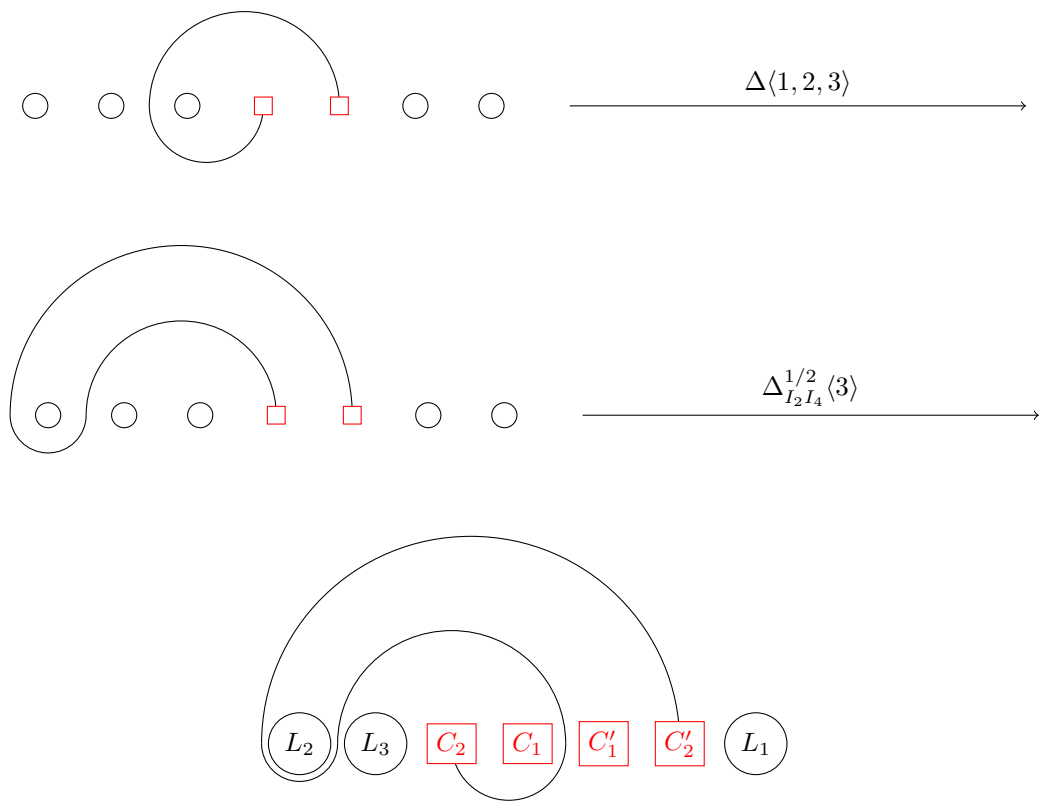
\begin{figure}[H]
	\begin{center}
		%
	\end{center}
\end{figure}

Relation:

\begin{equation*}\begin{split}
\Gamma_{C_1'}\Gamma_{C_1}\Gamma_{C_2}\Gamma_{L_3}\Gamma_{L_2}\Gamma_{L_3}^{-1}\Gamma_{C_2}^{-1}\Gamma_{C_1}^{-1}\Gamma_{C_2}\Gamma_{C_1}\Gamma_{C_2}\Gamma_{L_3}\Gamma_{L_2}^{-1}\Gamma_{L_3}^{-1}\Gamma_{C_2}^{-1}\Gamma_{C_1}^{-1}\Gamma_{C_1'}^{-1}  = \Gamma_{C_2'}
\end{split}
\end{equation*}
\subsubsection{Vertex Number 4}

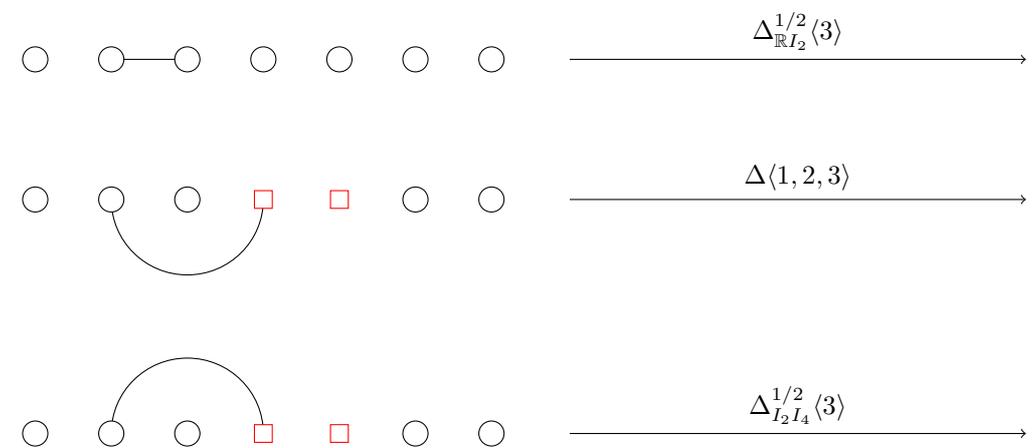
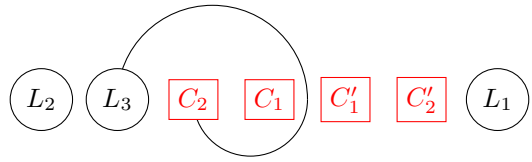
\begin{figure}[H]
	\begin{center}
		%
	\end{center}
\end{figure}

Relation:

\begin{equation*}\begin{split}
\left[\Gamma_{C_1}\Gamma_{C_2}\Gamma_{L_3}\Gamma_{C_2}^{-1}\Gamma_{C_1}^{-1} , \Gamma_{C_2} \right]=e
\end{split}
\end{equation*}
\subsubsection{Vertex Number 5}

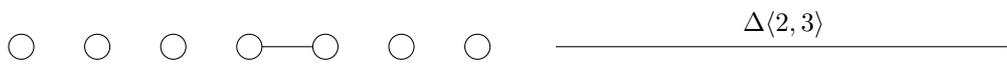
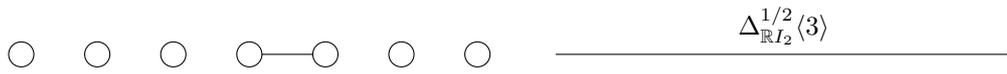
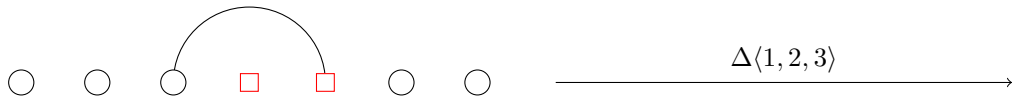
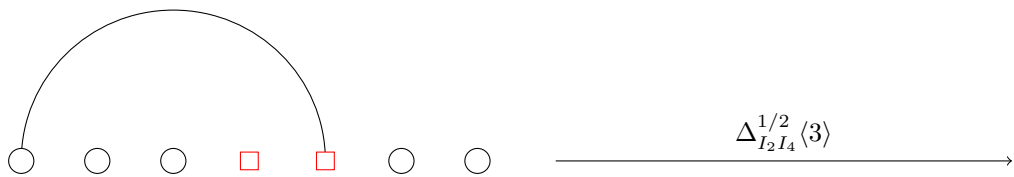
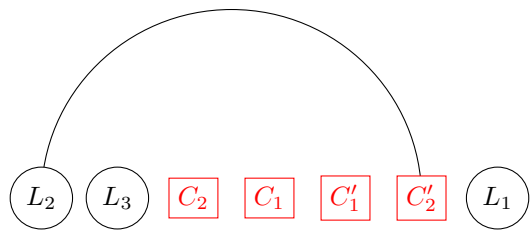
\begin{figure}[H]
	\begin{center}
		%
	\end{center}
\end{figure}

Relation:

\begin{equation*}\begin{split}
\left\{\Gamma_{L_3}^{-1}\Gamma_{C_2}^{-1}\Gamma_{C_1}^{-1}\Gamma_{C_1'}^{-1}\Gamma_{C_2'}\Gamma_{C_1'}\Gamma_{C_1}\Gamma_{C_2}\Gamma_{L_3}  , \Gamma_{L_2} \right\}=e
\end{split}
\end{equation*}
\subsubsection{Vertex Number 6}

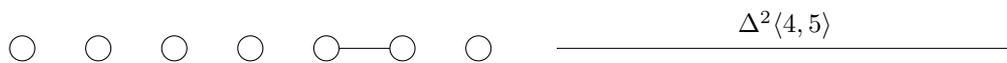
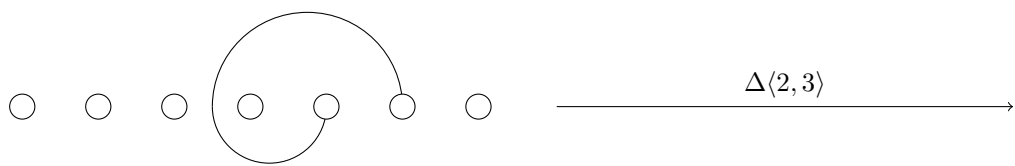
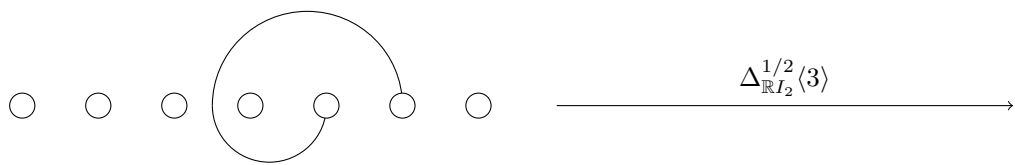
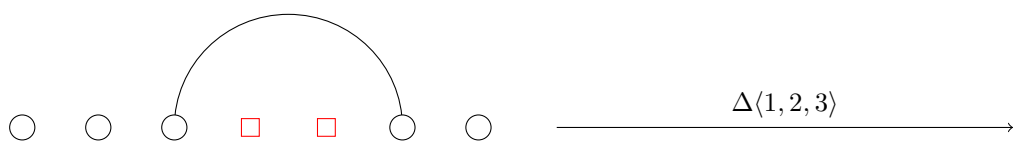
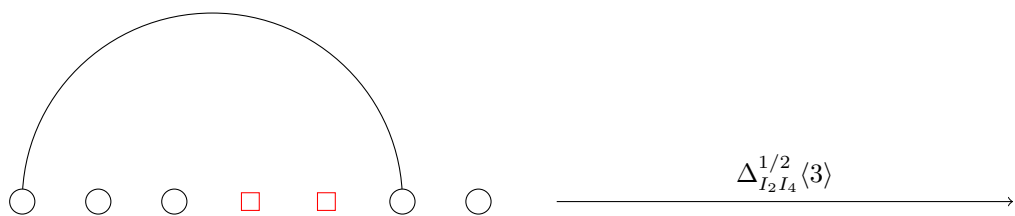
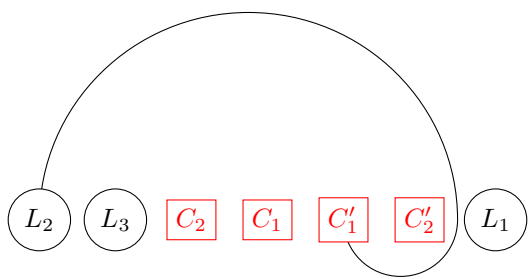
\begin{figure}[H]
	\begin{center}
		%
	\end{center}
\end{figure}

Relation:

\begin{equation*}\begin{split}
\left[\Gamma_{C_2'}\Gamma_{C_1'}\Gamma_{C_1}\Gamma_{C_2}\Gamma_{L_3}\Gamma_{L_2}\Gamma_{L_3}^{-1}\Gamma_{C_2}^{-1}\Gamma_{C_1}^{-1}\Gamma_{C_1'}^{-1}\Gamma_{C_2'}^{-1} , \Gamma_{C_1'} \right]=e
\end{split}
\end{equation*}
\subsubsection{Vertex Number 7}

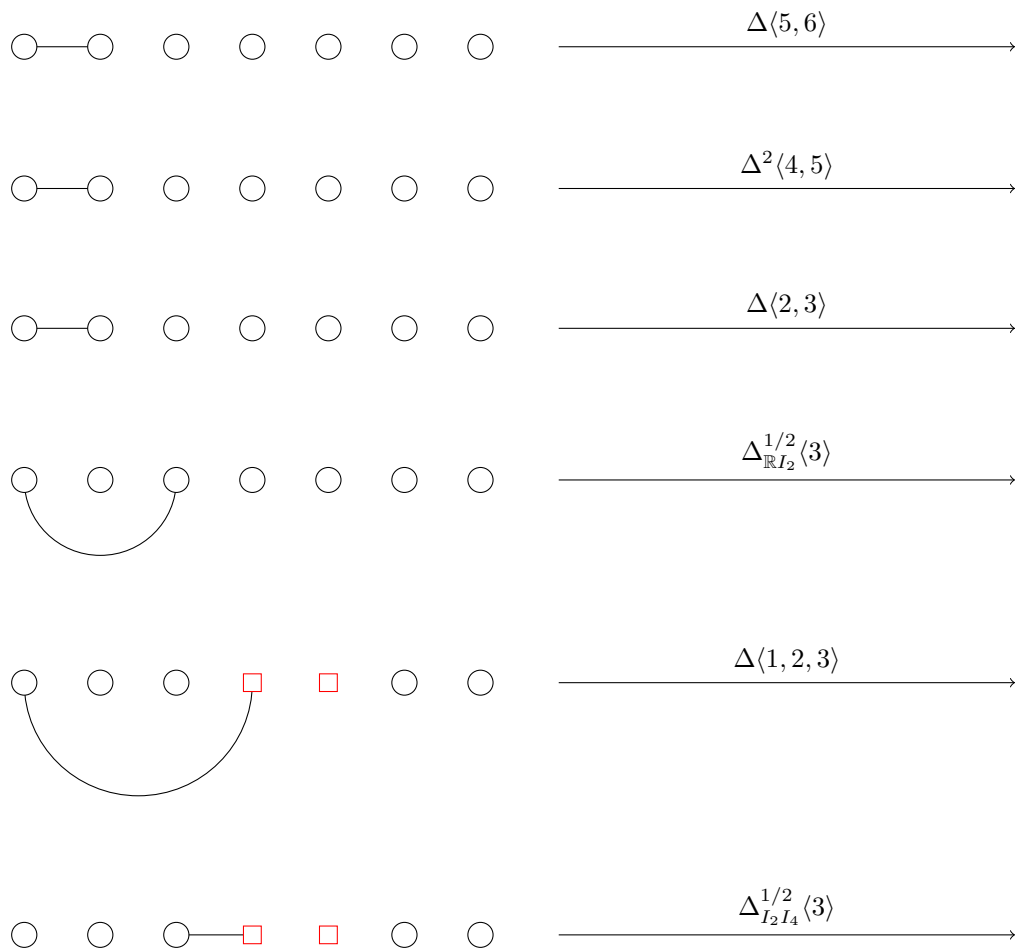
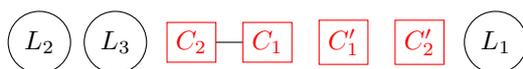
\begin{figure}[H]
	\begin{center}
		%
	\end{center}
\end{figure}

Relation:

\begin{equation*}\begin{split}
\left\{\Gamma_{C_1}  , \Gamma_{C_2} \right\}=e
\end{split}
\end{equation*}
\subsubsection{Vertex Number 8}

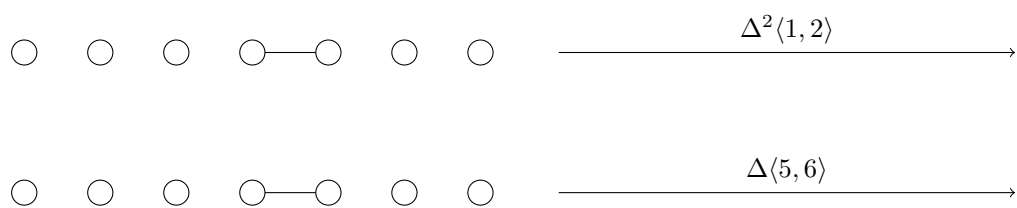
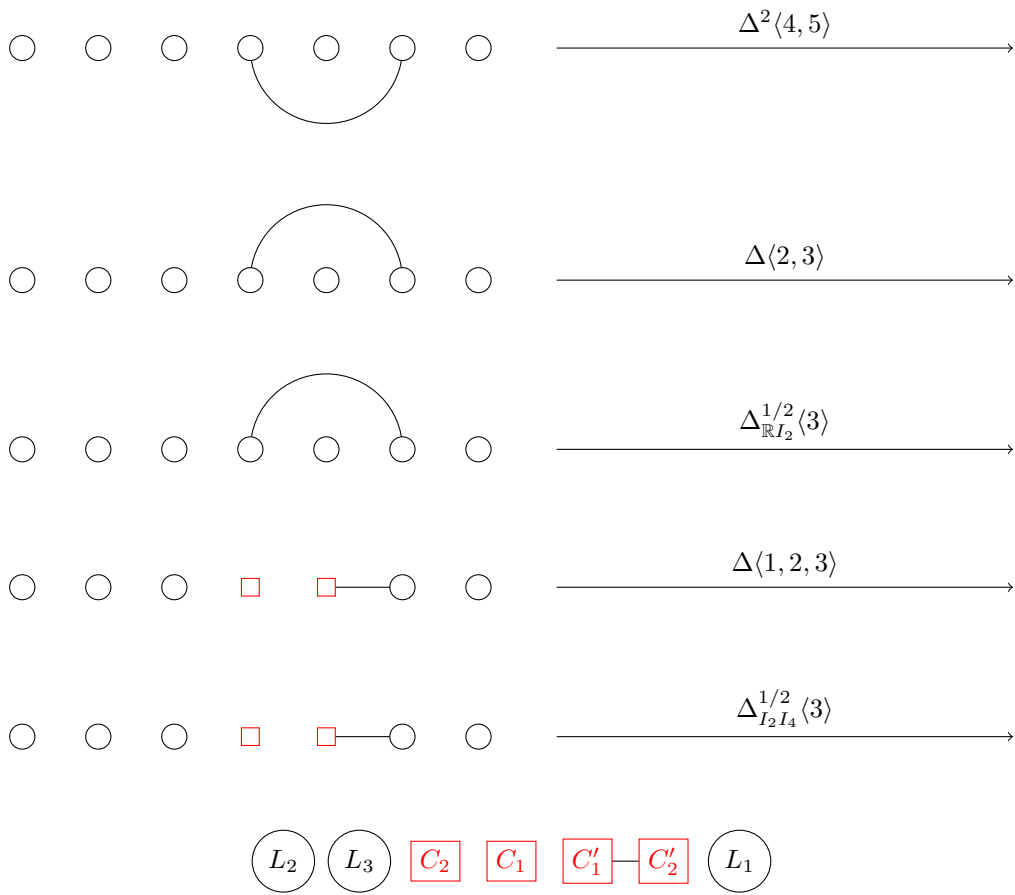
\begin{figure}[H]
	\begin{center}
		%
	\end{center}
\end{figure}

Relation:

\begin{equation*}\begin{split}
\left\{\Gamma_{C_2'}  , \Gamma_{C_1'} \right\}=e
\end{split}
\end{equation*}
\subsubsection{Vertex Number 9}

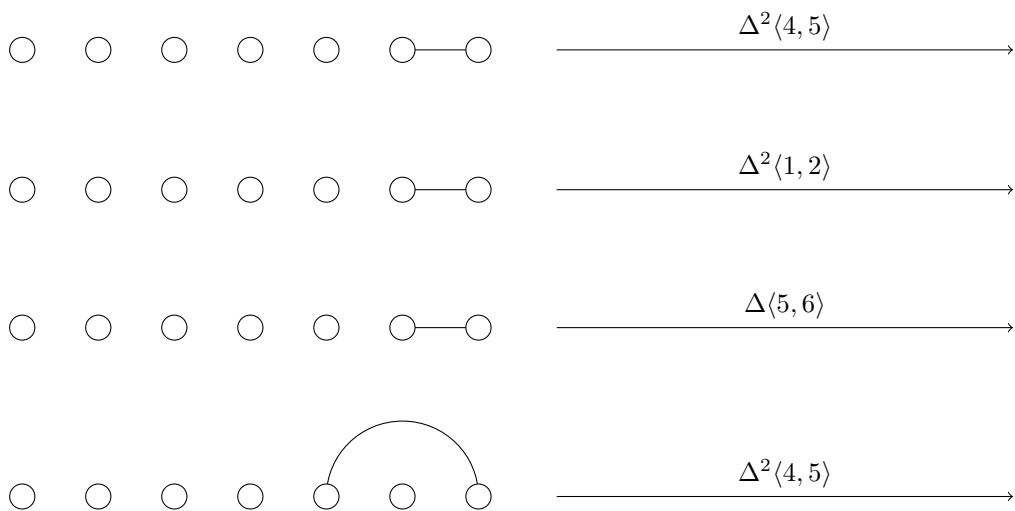
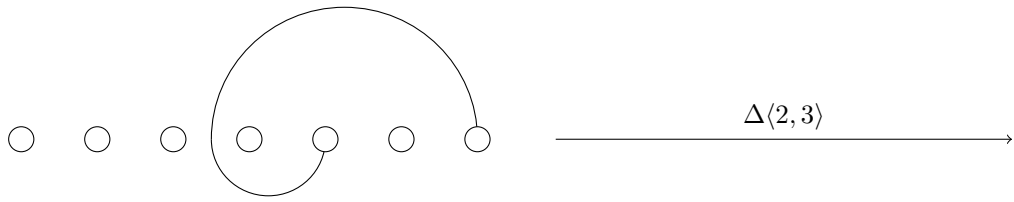
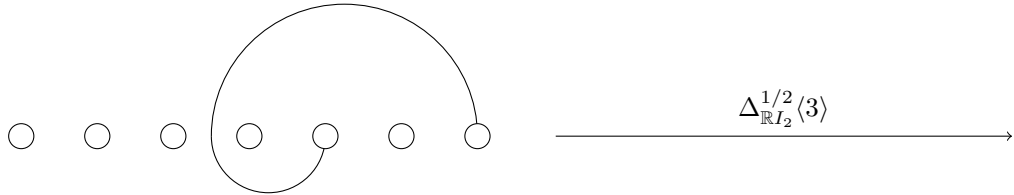
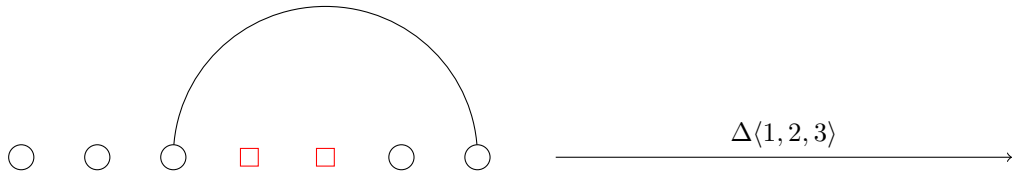
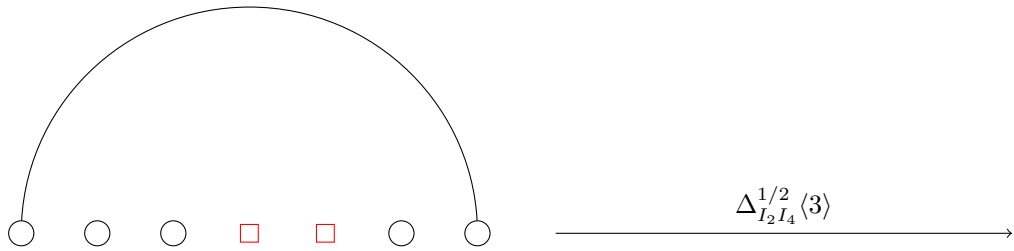
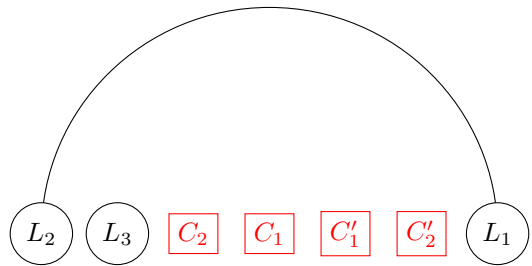
\begin{figure}[H]
	\begin{center}
		%
	\end{center}
\end{figure}

Relation:

\begin{equation*}\begin{split}
\left[\Gamma_{C_2'}\Gamma_{C_1'}\Gamma_{C_1}\Gamma_{C_2}\Gamma_{L_3}\Gamma_{L_2}\Gamma_{L_3}^{-1}\Gamma_{C_2}^{-1}\Gamma_{C_1}^{-1}\Gamma_{C_1'}^{-1}\Gamma_{C_2'}^{-1} , \Gamma_{L_1} \right]=e
\end{split}
\end{equation*}
\subsubsection{Vertex Number 10}

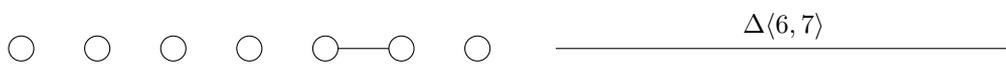
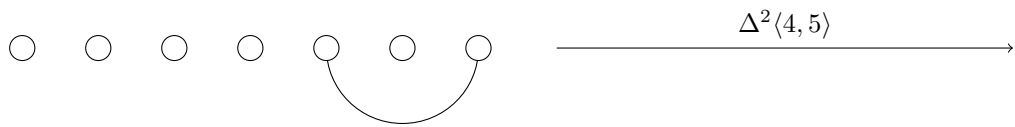
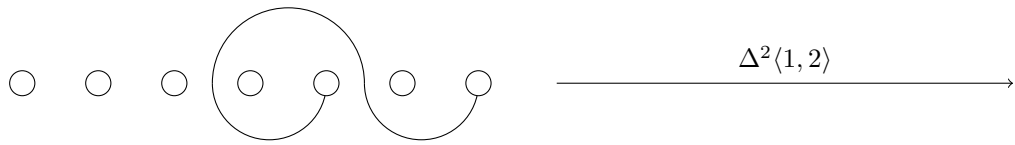
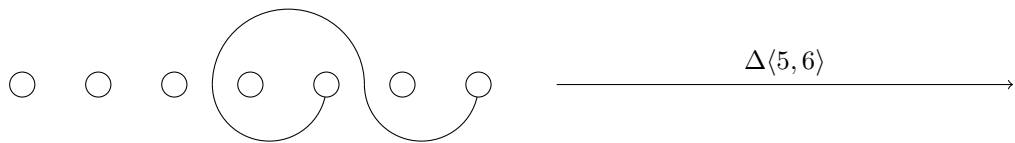
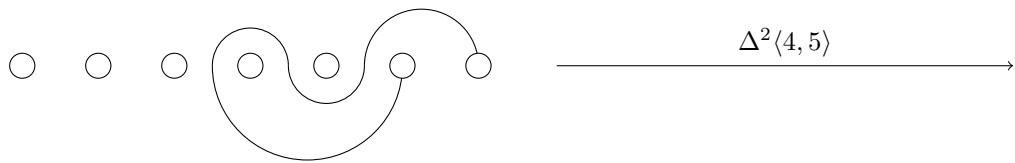
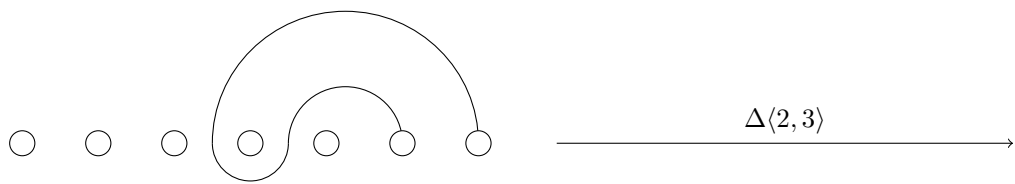
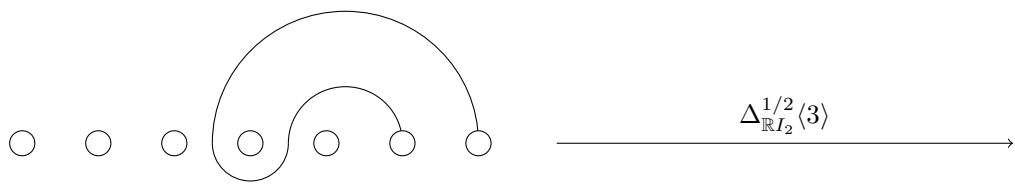
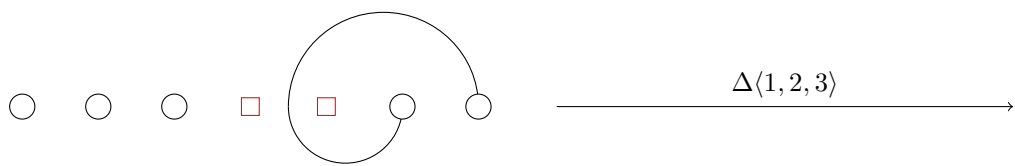
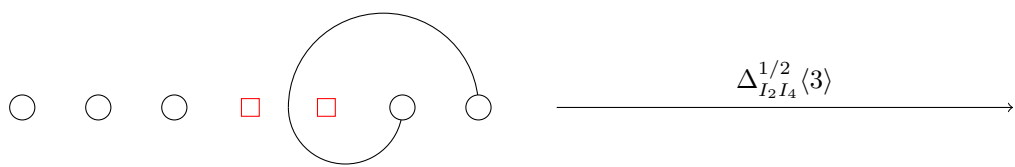
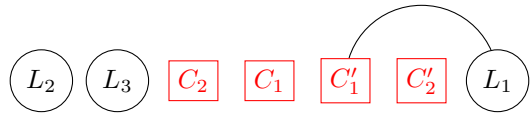
\begin{figure}[H]
	\begin{center}
		%
	\end{center}
\end{figure}

Relation:

\begin{equation*}\begin{split}
\left[\Gamma_{C_2'}\Gamma_{C_1'}\Gamma_{C_2'}^{-1} , \Gamma_{L_1} \right]=e
\end{split}
\end{equation*}
\subsubsection{Vertex Number 11}

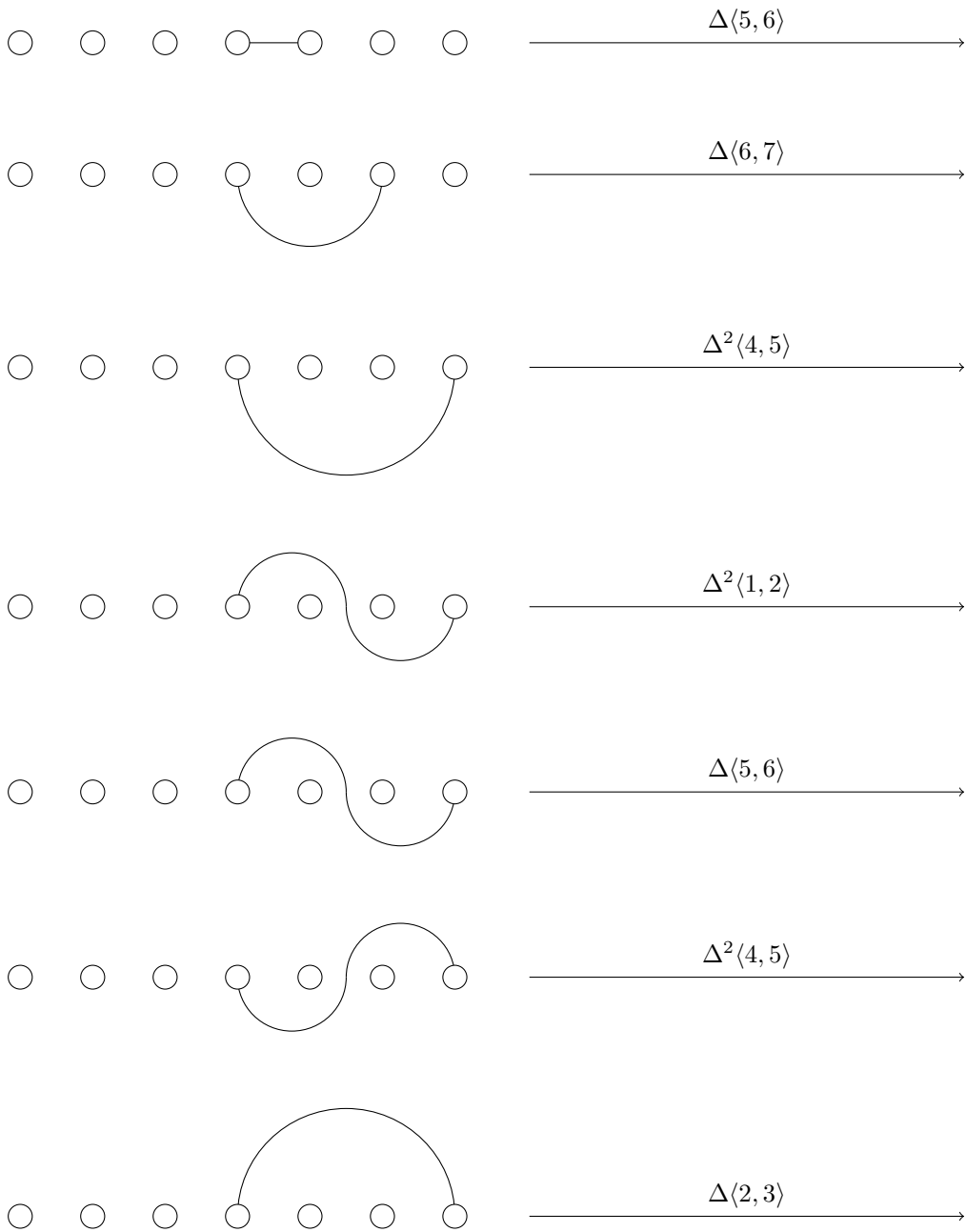
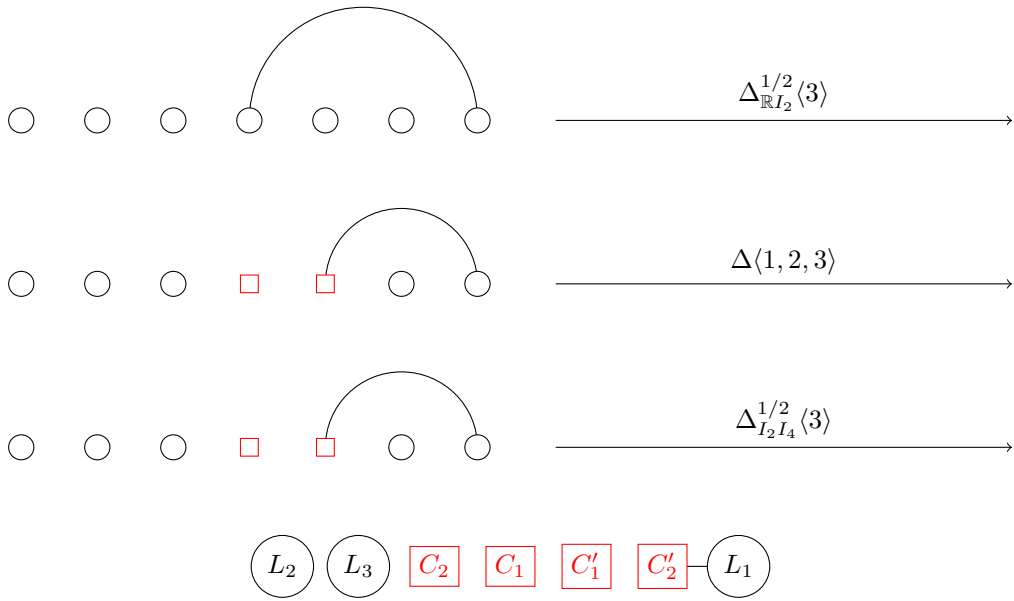
\begin{figure}[H]
	\begin{center}
		%
	\end{center}
\end{figure}

Relation:

\begin{equation*}\begin{split}
\left\{\Gamma_{C_2'}  , \Gamma_{L_1} \right\}=e
\end{split}
\end{equation*}
\subsubsection{Vertex Number 12}

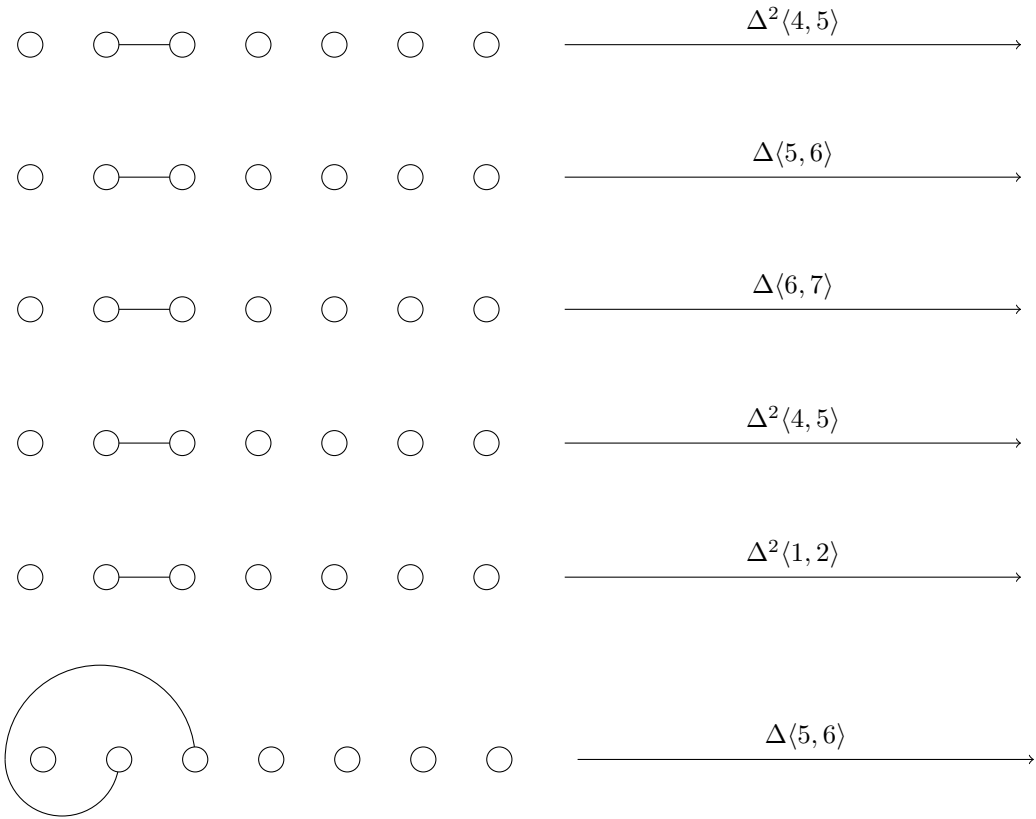
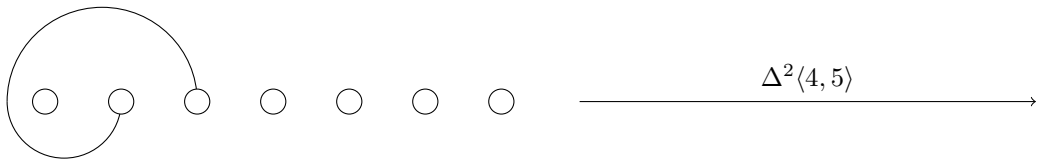
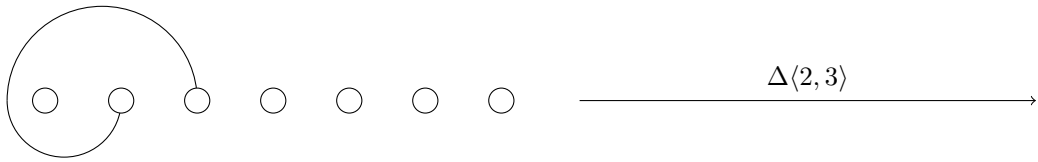
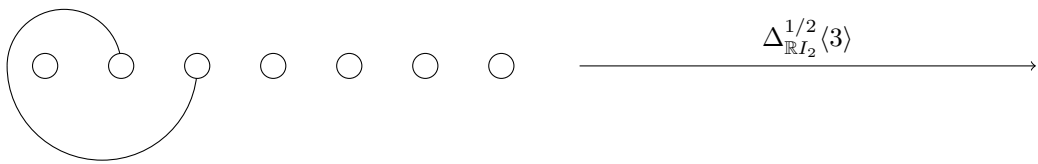
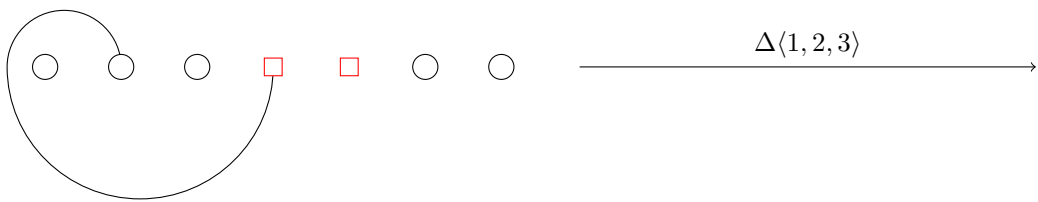
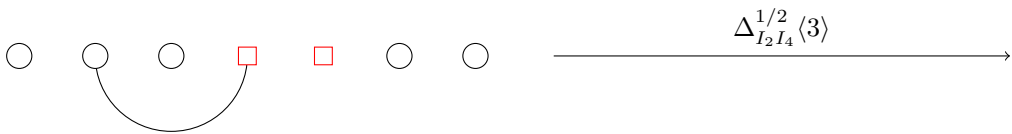
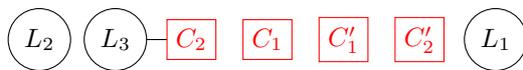
\begin{figure}[H]
	\begin{center}
		%
	\end{center}
\end{figure}

Relation:

\begin{equation*}\begin{split}
\left[\Gamma_{C_2} , \Gamma_{L_3} \right]=e
\end{split}
\end{equation*}
\subsubsection{Vertex Number 13}

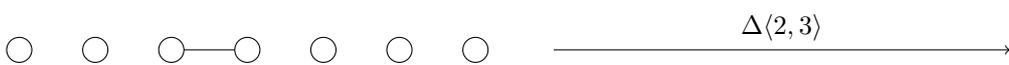
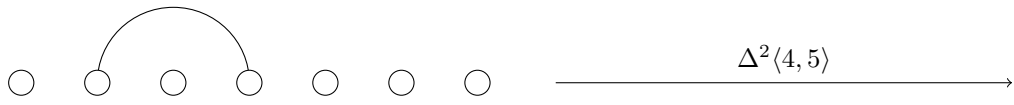
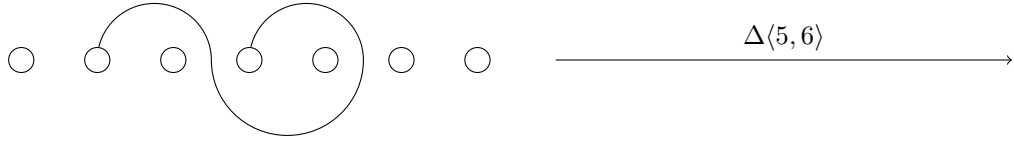
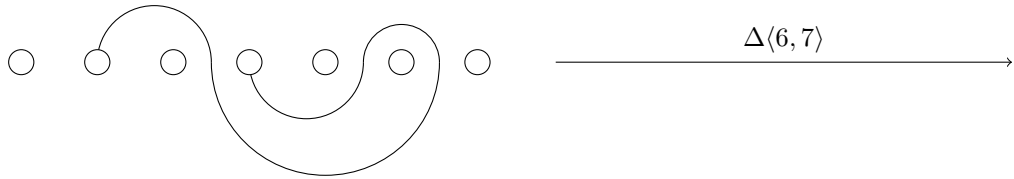
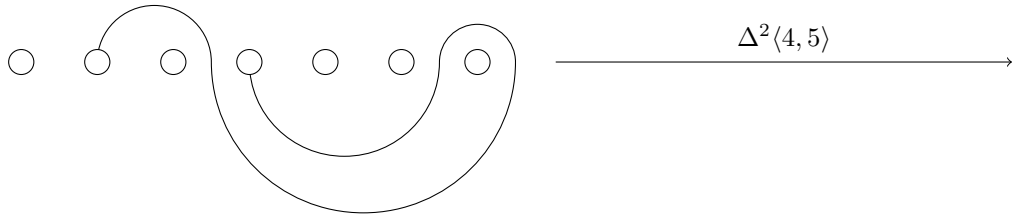
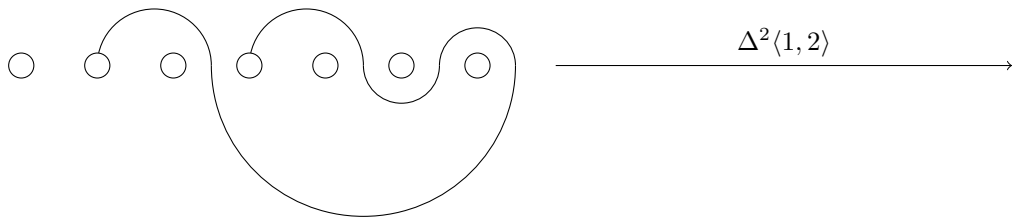
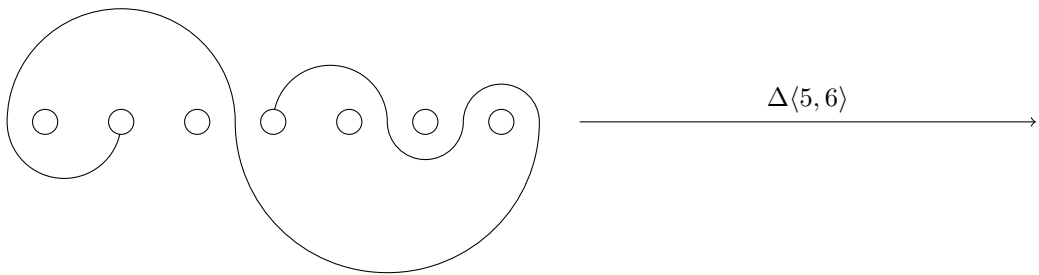
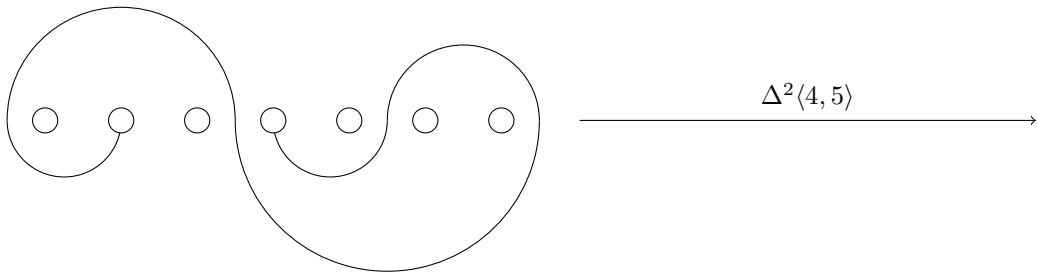
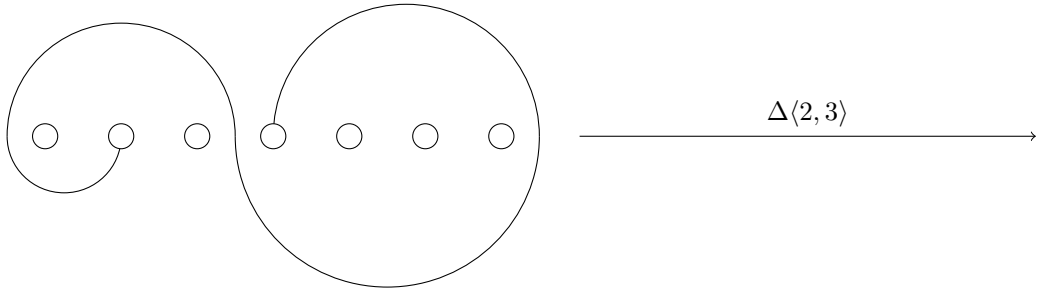
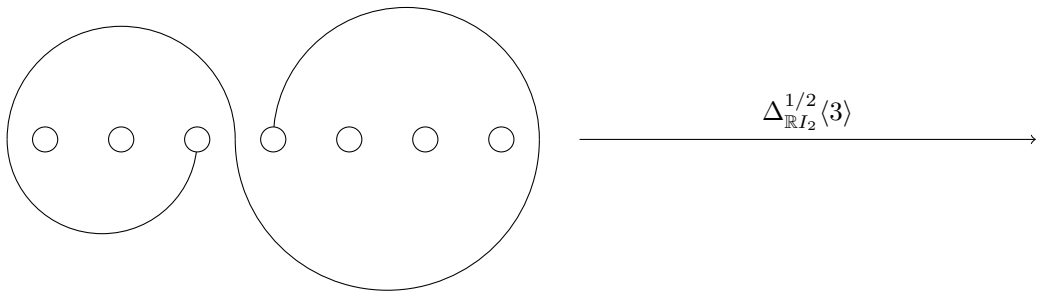
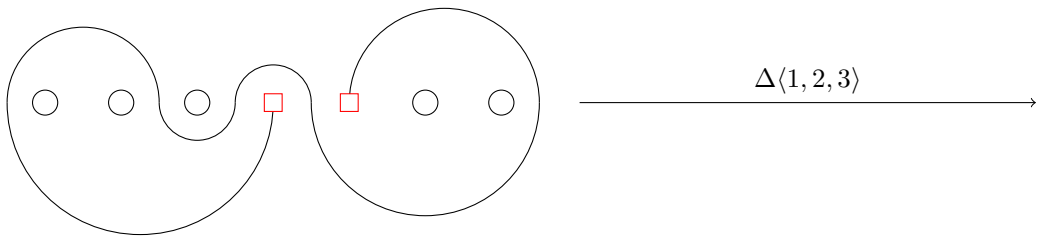
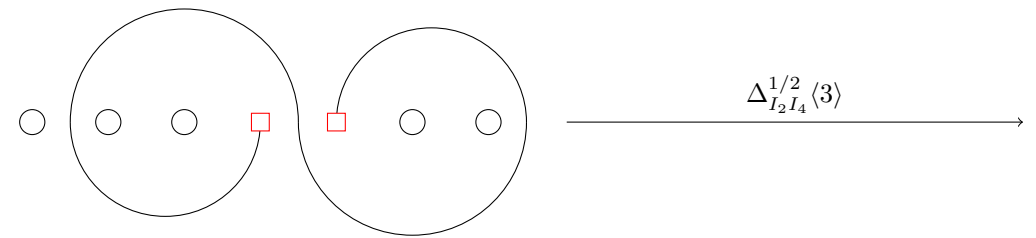
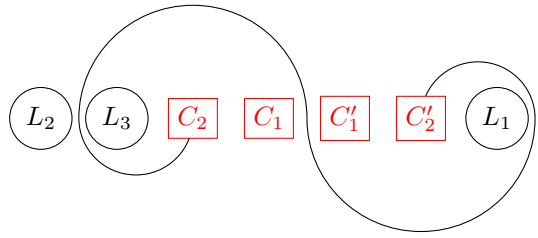
\begin{figure}[H]
	\begin{center}
		%
	\end{center}
\end{figure}

Relation:

\begin{equation*}\begin{split}
\Gamma_{L_1}^{-1}\Gamma_{C_1}\Gamma_{C_2}\Gamma_{L_3}\Gamma_{C_2}\Gamma_{L_3}^{-1}\Gamma_{C_2}^{-1}\Gamma_{C_1}^{-1}\Gamma_{L_1}  = \Gamma_{C_2'}
\end{split}
\end{equation*}
\subsubsection{Vertex Number 14}

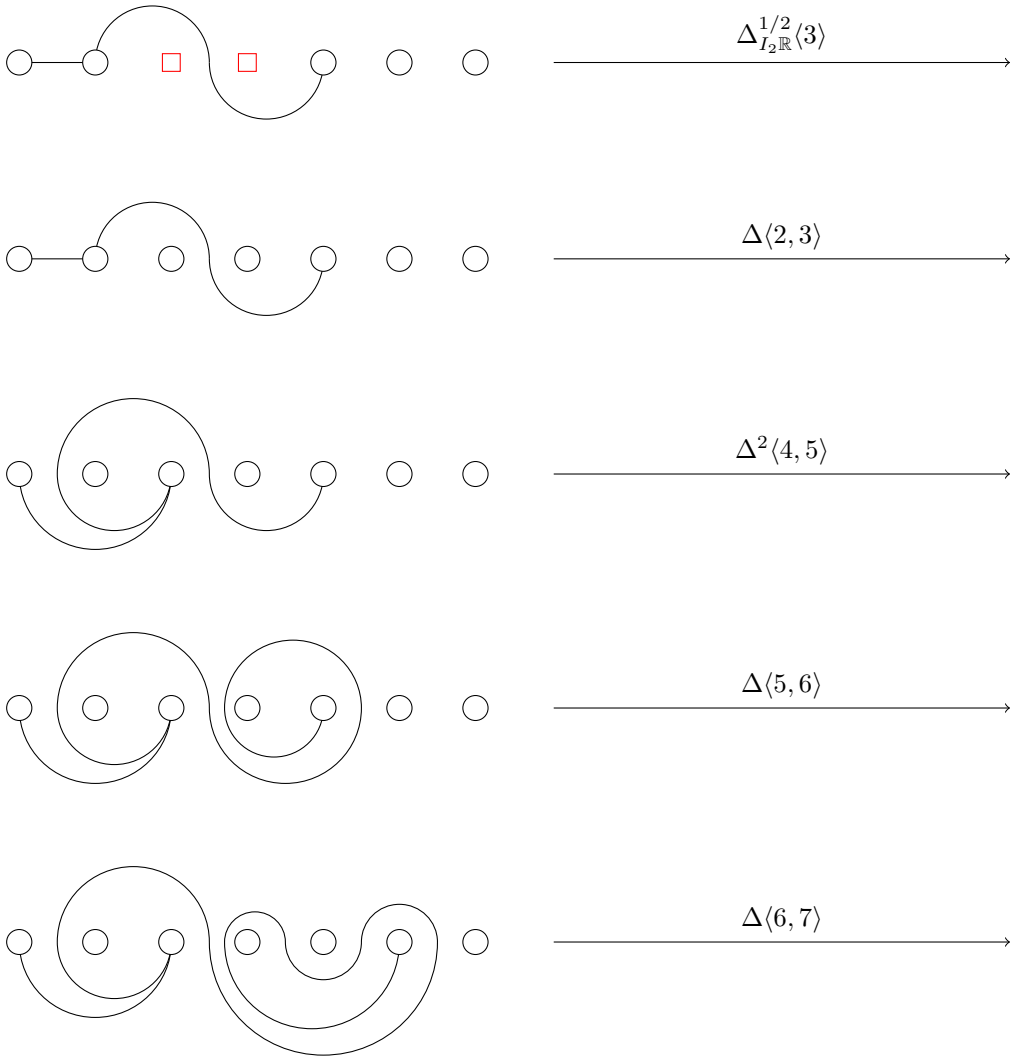
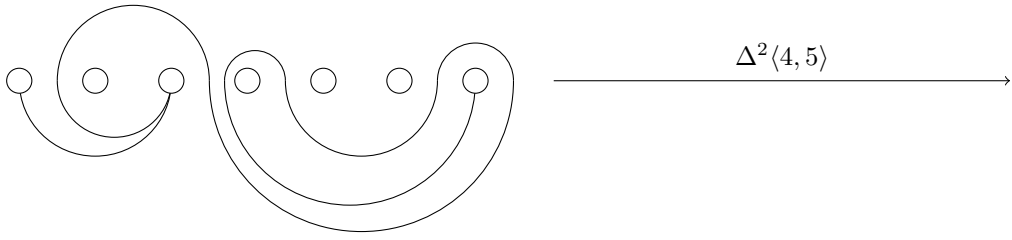
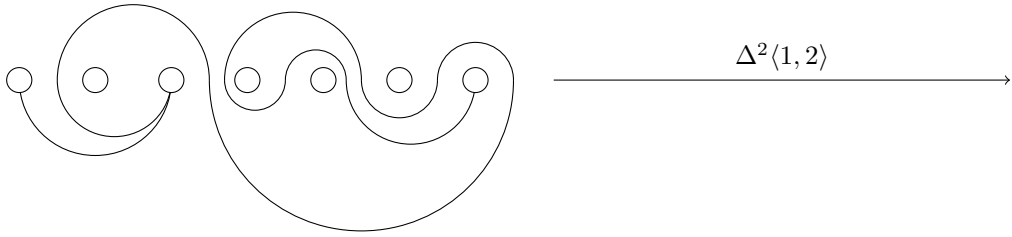
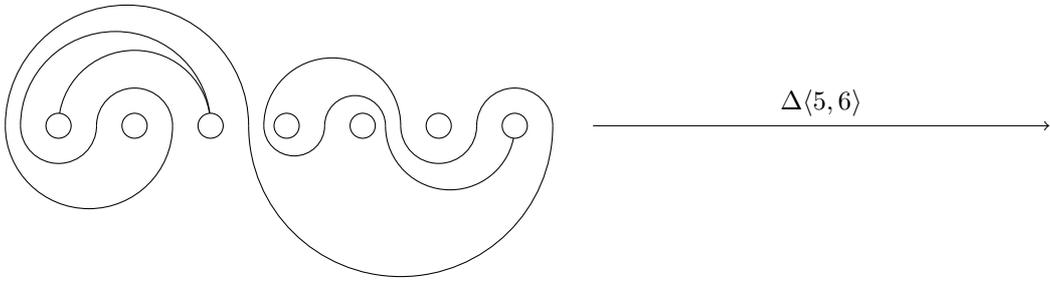
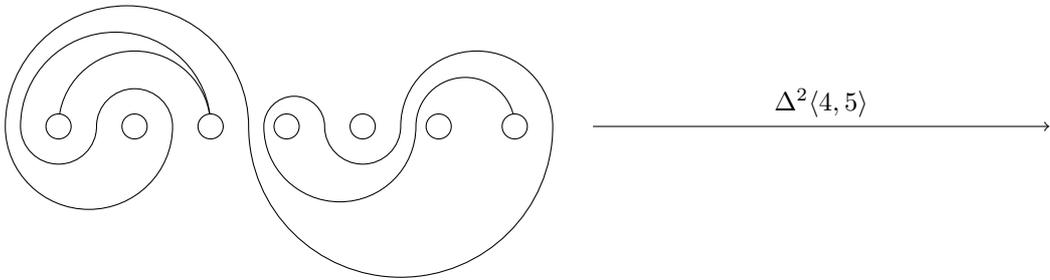
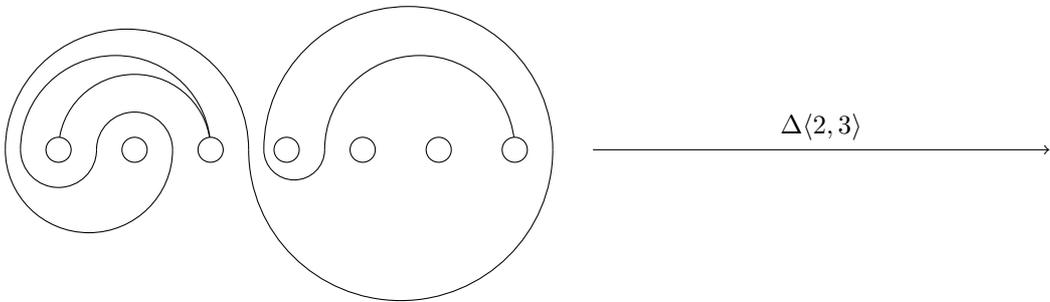
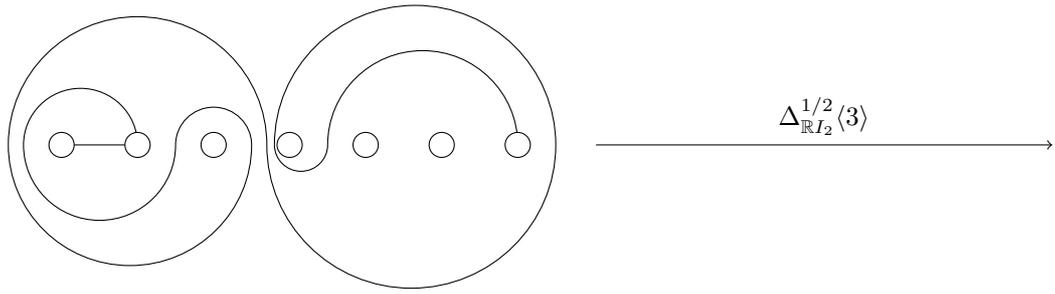
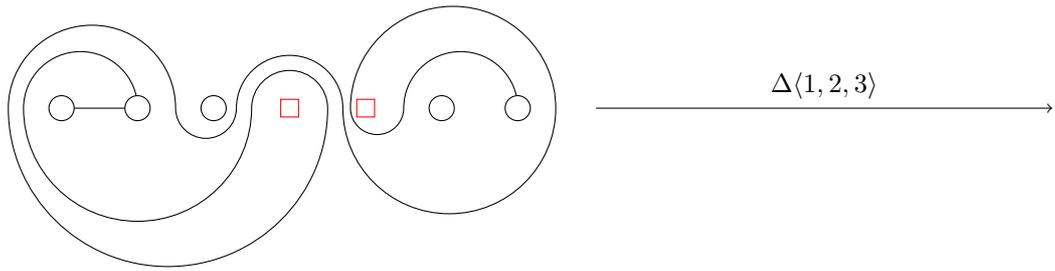
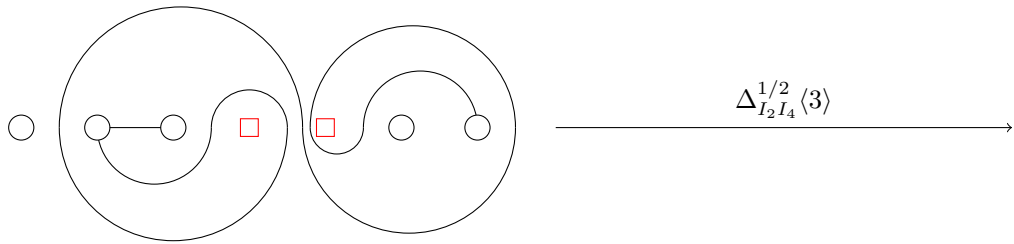
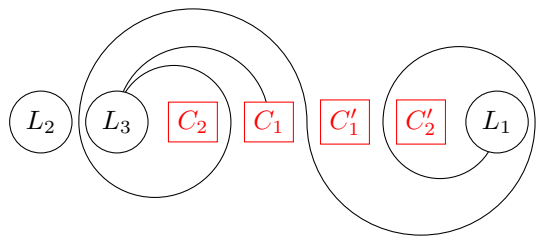
\begin{figure}[H]
	\begin{center}
		%
	\end{center}
\end{figure}

Relation:

\begin{equation*}\begin{split}
\left( \Gamma_{C_2}^{-1}\Gamma_{L_3}^{-1}\Gamma_{C_2}^{-1}\Gamma_{C_1}^{-1}\Gamma_{L_1}\Gamma_{C_2'}\Gamma_{L_1}\Gamma_{C_2'}^{-1}\Gamma_{L_1}^{-1}\Gamma_{C_1}\Gamma_{C_2}\Gamma_{L_3}\Gamma_{C_2}  \right) \left( \Gamma_{L_3}  \right) \left( \Gamma_{C_2}^{-1}\Gamma_{C_1}\Gamma_{C_2}  \right)& = \\
\left( \Gamma_{L_3}  \right) \left( \Gamma_{C_2}^{-1}\Gamma_{C_1}\Gamma_{C_2}  \right) \left( \Gamma_{C_2}^{-1}\Gamma_{L_3}^{-1}\Gamma_{C_2}^{-1}\Gamma_{C_1}^{-1}\Gamma_{L_1}\Gamma_{C_2'}\Gamma_{L_1}\Gamma_{C_2'}^{-1}\Gamma_{L_1}^{-1}\Gamma_{C_1}\Gamma_{C_2}\Gamma_{L_3}\Gamma_{C_2}  \right)& = \\
\left( \Gamma_{C_2}^{-1}\Gamma_{C_1}\Gamma_{C_2}  \right) \left( \Gamma_{C_2}^{-1}\Gamma_{L_3}^{-1}\Gamma_{C_2}^{-1}\Gamma_{C_1}^{-1}\Gamma_{L_1}\Gamma_{C_2'}\Gamma_{L_1}\Gamma_{C_2'}^{-1}\Gamma_{L_1}^{-1}\Gamma_{C_1}\Gamma_{C_2}\Gamma_{L_3}\Gamma_{C_2}  \right) \left( \Gamma_{L_3}  \right)
\end{split}
\end{equation*}
\subsubsection{Vertex Number 15}

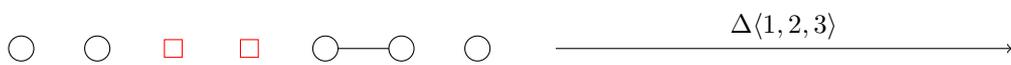
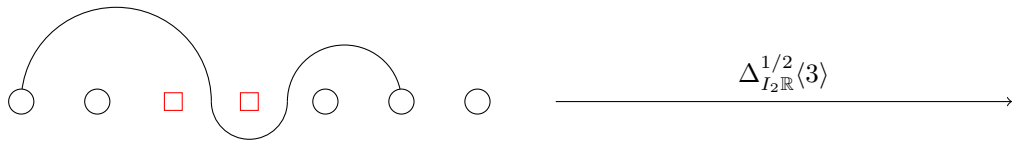
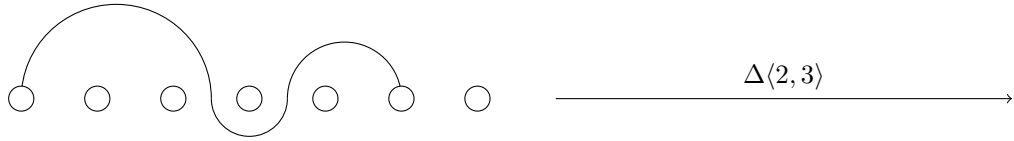
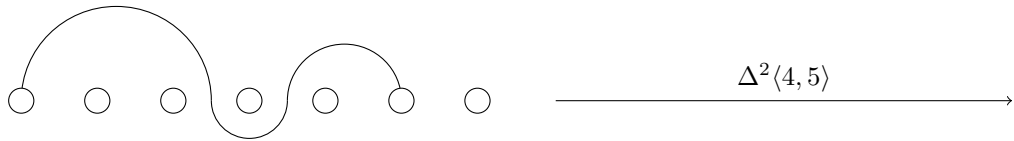
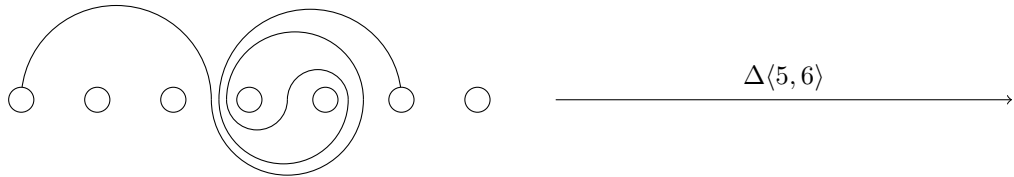
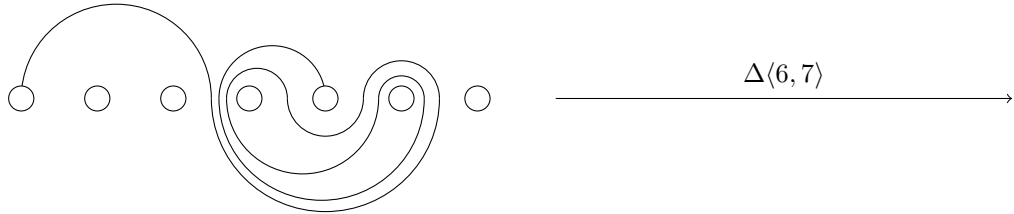
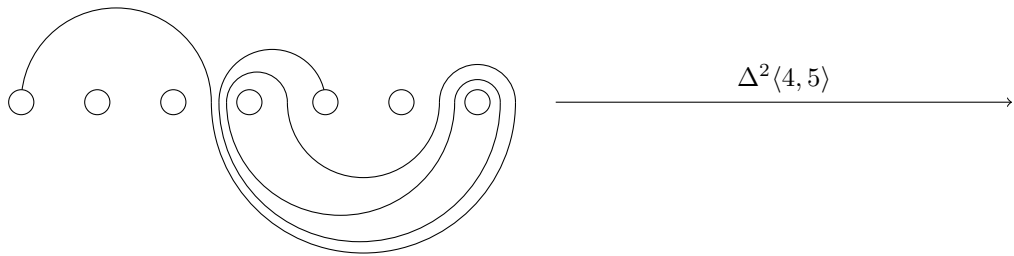
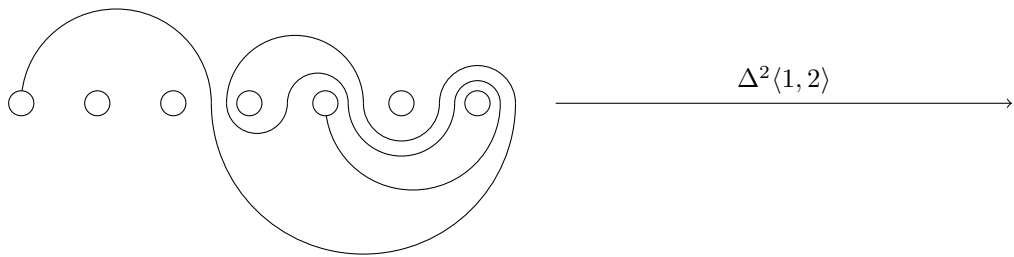
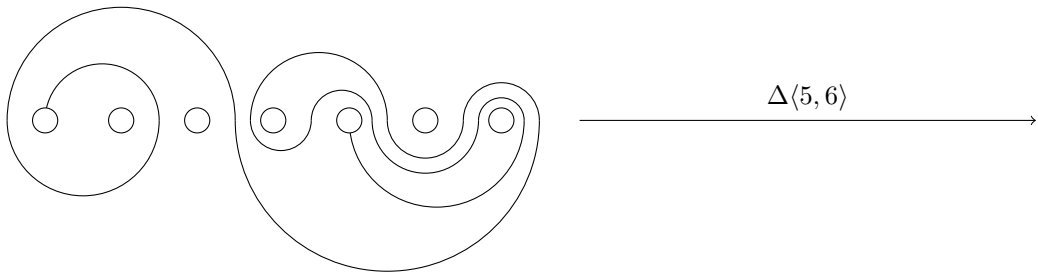
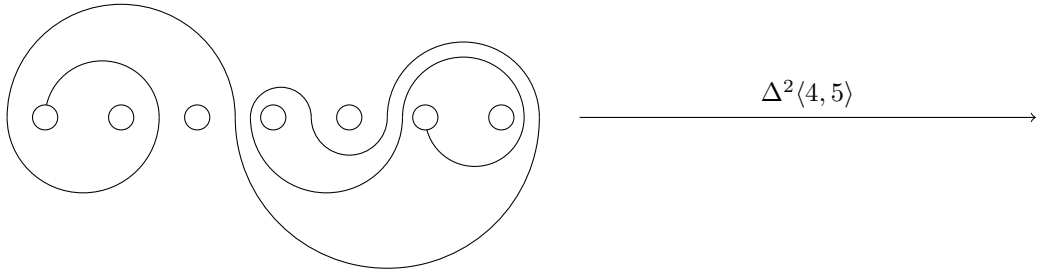
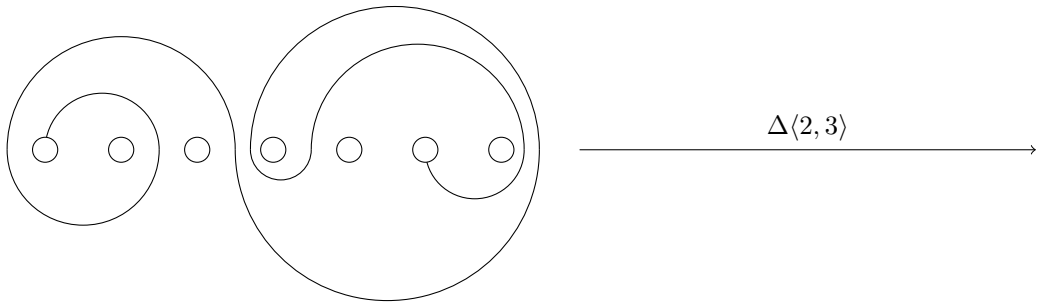
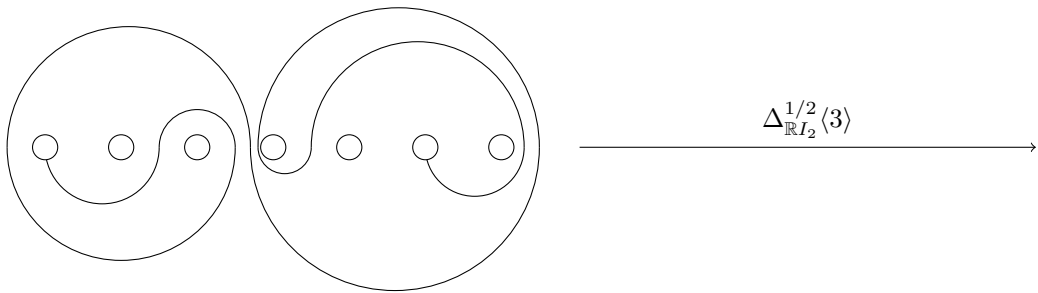
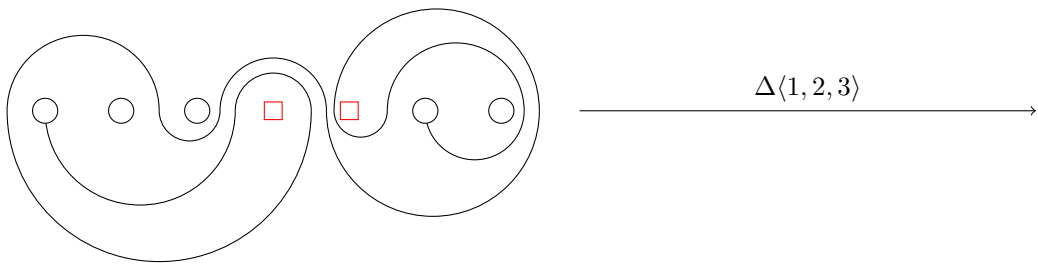
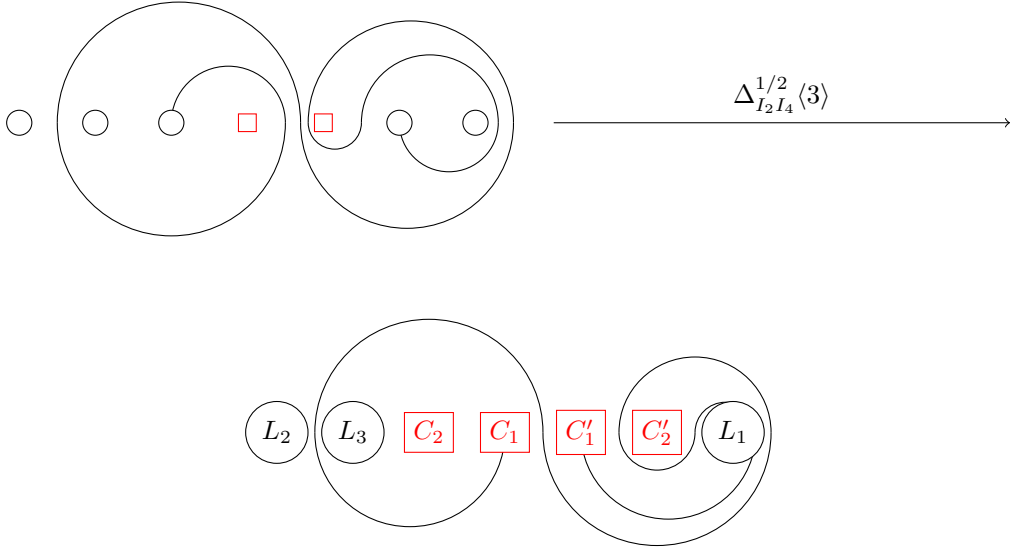
\begin{figure}[H]
	\begin{center}

	\end{center}
\end{figure}

Relation:

\begin{equation*}\begin{split}
\Gamma_{L_1}\Gamma_{C_2'}^{-1}\Gamma_{L_1}^{-1}\Gamma_{C_1}\Gamma_{C_2}\Gamma_{L_3}\Gamma_{C_1}\Gamma_{L_3}^{-1}\Gamma_{C_2}^{-1}\Gamma_{C_1}^{-1}\Gamma_{L_1}\Gamma_{C_2'}\Gamma_{L_1}^{-1}  = \Gamma_{C_1'}
\end{split}
\end{equation*}
\subsubsection{Raw relations and simplifications}

\begin{equation}\label{arrangement1_auto_calc_vert_1_rel_1}
\begin{split}
\Gamma_{C_1}  = \Gamma_{C_1'} ,
\end{split}
\end{equation}
\begin{equation}\label{arrangement1_auto_calc_vert_2_rel_1}
\begin{split}
\left[\Gamma_{C_2}^{-1}\Gamma_{C_1}\Gamma_{C_2} , \Gamma_{L_3}\Gamma_{L_2}  \right]=e,
\end{split}
\end{equation}
\begin{equation}\label{arrangement1_auto_calc_vert_2_rel_2}
\begin{split}
\left[\Gamma_{L_3} , \Gamma_{L_2}\Gamma_{C_2}^{-1}\Gamma_{C_1}\Gamma_{C_2}  \right]=e,
\end{split}
\end{equation}
\begin{equation}\label{arrangement1_auto_calc_vert_4_rel_1}
\begin{split}
\Gamma_{C_1'}\Gamma_{C_1}\Gamma_{C_2}\Gamma_{L_3}\Gamma_{L_2}\Gamma_{L_3}^{-1}\Gamma_{C_2}^{-1}\Gamma_{C_1}^{-1}\Gamma_{C_2}\Gamma_{C_1}\Gamma_{C_2}\Gamma_{L_3}\Gamma_{L_2}^{-1}\Gamma_{L_3}^{-1}\Gamma_{C_2}^{-1}\Gamma_{C_1}^{-1}\Gamma_{C_1'}^{-1}  = \Gamma_{C_2'} ,
\end{split}
\end{equation}
\begin{equation}\label{arrangement1_auto_calc_vert_5_rel_1}
\begin{split}
\left[\Gamma_{C_1}\Gamma_{C_2}\Gamma_{L_3}\Gamma_{C_2}^{-1}\Gamma_{C_1}^{-1} , \Gamma_{C_2}  \right]=e,
\end{split}
\end{equation}
\begin{equation}\label{arrangement1_auto_calc_vert_6_rel_1}
\begin{split}
\left\{\Gamma_{L_3}^{-1}\Gamma_{C_2}^{-1}\Gamma_{C_1}^{-1}\Gamma_{C_1'}^{-1}\Gamma_{C_2'}\Gamma_{C_1'}\Gamma_{C_1}\Gamma_{C_2}\Gamma_{L_3} , \Gamma_{L_2}  \right\}=e,
\end{split}
\end{equation}
\begin{equation}\label{arrangement1_auto_calc_vert_7_rel_1}
\begin{split}
\left[\Gamma_{C_2'}\Gamma_{C_1'}\Gamma_{C_1}\Gamma_{C_2}\Gamma_{L_3}\Gamma_{L_2}\Gamma_{L_3}^{-1}\Gamma_{C_2}^{-1}\Gamma_{C_1}^{-1}\Gamma_{C_1'}^{-1}\Gamma_{C_2'}^{-1} , \Gamma_{C_1'}  \right]=e,
\end{split}
\end{equation}
\begin{equation}\label{arrangement1_auto_calc_vert_8_rel_1}
\begin{split}
\left\{\Gamma_{C_1} , \Gamma_{C_2}  \right\}=e,
\end{split}
\end{equation}
\begin{equation}\label{arrangement1_auto_calc_vert_9_rel_1}
\begin{split}
\left\{\Gamma_{C_2'} , \Gamma_{C_1'}  \right\}=e,
\end{split}
\end{equation}
\begin{equation}\label{arrangement1_auto_calc_vert_10_rel_1}
\begin{split}
\left[\Gamma_{C_2'}\Gamma_{C_1'}\Gamma_{C_1}\Gamma_{C_2}\Gamma_{L_3}\Gamma_{L_2}\Gamma_{L_3}^{-1}\Gamma_{C_2}^{-1}\Gamma_{C_1}^{-1}\Gamma_{C_1'}^{-1}\Gamma_{C_2'}^{-1} , \Gamma_{L_1}  \right]=e,
\end{split}
\end{equation}
\begin{equation}\label{arrangement1_auto_calc_vert_11_rel_1}
\begin{split}
\left[\Gamma_{C_2'}\Gamma_{C_1'}\Gamma_{C_2'}^{-1} , \Gamma_{L_1}  \right]=e,
\end{split}
\end{equation}
\begin{equation}\label{arrangement1_auto_calc_vert_12_rel_1}
\begin{split}
\left\{\Gamma_{C_2'} , \Gamma_{L_1}  \right\}=e,
\end{split}
\end{equation}
\begin{equation}\label{arrangement1_auto_calc_vert_13_rel_1}
\begin{split}
\left[\Gamma_{C_2} , \Gamma_{L_3}  \right]=e,
\end{split}
\end{equation}
\begin{equation}\label{arrangement1_auto_calc_vert_14_rel_1}
\begin{split}
\Gamma_{L_1}^{-1}\Gamma_{C_1}\Gamma_{C_2}\Gamma_{L_3}\Gamma_{C_2}\Gamma_{L_3}^{-1}\Gamma_{C_2}^{-1}\Gamma_{C_1}^{-1}\Gamma_{L_1}  = \Gamma_{C_2'} ,
\end{split}
\end{equation}
\begin{equation}\label{arrangement1_auto_calc_vert_15_rel_1}
\begin{split}
\left[\Gamma_{C_2}^{-1}\Gamma_{L_3}^{-1}\Gamma_{C_2}^{-1}\Gamma_{C_1}^{-1}\Gamma_{L_1}\Gamma_{C_2'}\Gamma_{L_1}\Gamma_{C_2'}^{-1}\Gamma_{L_1}^{-1}\Gamma_{C_1}\Gamma_{C_2}\Gamma_{L_3}\Gamma_{C_2} , \Gamma_{L_3}\Gamma_{C_2}^{-1}\Gamma_{C_1}\Gamma_{C_2}  \right]=e,
\end{split}
\end{equation}
\begin{equation}\label{arrangement1_auto_calc_vert_15_rel_2}
\begin{split}
\left[\Gamma_{L_3} , \Gamma_{C_2}^{-1}\Gamma_{C_1}\Gamma_{L_3}^{-1}\Gamma_{C_2}^{-1}\Gamma_{C_1}^{-1}\Gamma_{L_1}\Gamma_{C_2'}\Gamma_{L_1}\Gamma_{C_2'}^{-1}\Gamma_{L_1}^{-1}\Gamma_{C_1}\Gamma_{C_2}\Gamma_{L_3}\Gamma_{C_2}  \right]=e,
\end{split}
\end{equation}
\begin{equation}\label{arrangement1_auto_calc_vert_16_rel_1}
\begin{split}
\Gamma_{L_1}\Gamma_{C_2'}^{-1}\Gamma_{L_1}^{-1}\Gamma_{C_1}\Gamma_{C_2}\Gamma_{L_3}\Gamma_{C_1}\Gamma_{L_3}^{-1}\Gamma_{C_2}^{-1}\Gamma_{C_1}^{-1}\Gamma_{L_1}\Gamma_{C_2'}\Gamma_{L_1}^{-1}  = \Gamma_{C_1'} ,
\end{split}
\end{equation}
\begin{equation}\label{arrangement1_auto_calc_projective_rel}
\begin{split}
\Gamma_{L_1}\Gamma_{C_2'}\Gamma_{C_1'}\Gamma_{C_1}\Gamma_{C_2}\Gamma_{L_3}\Gamma_{L_2} =e.
\end{split}
\end{equation}

\bigskip
Here are the computations of group $ G_5 := \pcpt{\mathcal{B}_5} / \langle \Gamma_X^2 \; | \; X\subseteq \mathcal{B}_5 \rangle $:
\begin{lemma}\label{tokunaga-arrangement1-b1-group-calc}
The quotient $ G_5 := \pcpt{\mathcal{B}_5} / \langle \Gamma_X^2 \; | \; X\subseteq \mathcal{B}_5 \rangle $ is generated by $ \Gamma_{C_1}, \Gamma_{C_2}, \Gamma_{L_1}, \Gamma_{L_2} $, and $ \Gamma_{L_3} $, subject to the following relations:
\begin{equation}\label{arrangement1_final_squares}
	\Gamma_{C_1}^2 = \Gamma_{C_2}^2 = \Gamma_{L_1}^2 = \Gamma_{L_2}^2 = \Gamma_{L_3}^2 = e,
\end{equation}
\begin{equation}\label{arrangement1_final1}
	[\Gamma_{L_3},\Gamma_{L_1}] = [\Gamma_{L_3}, \Gamma_{L_2}] = [\Gamma_{L_3}, \Gamma_{C_1}] = [\Gamma_{L_3}, \Gamma_{C_2}] = [\Gamma_{L_1}, \Gamma_{L_2}] = [\Gamma_{C_1}, \Gamma_{L_1}] = [\Gamma_{C_1}, \Gamma_{L_2}] = e,
\end{equation}
\begin{equation}\label{arrangement1_final2}
	[\Gamma_{L_1}, \Gamma_{C_2}\Gamma_{C_1}\Gamma_{C_2}] = [\Gamma_{L_2}, \Gamma_{C_2}\Gamma_{C_1}\Gamma_{C_2}] = [\Gamma_{C_2}, \Gamma_{L_1}\Gamma_{L_2}] = e,
\end{equation}
\begin{equation}\label{arrangement1_final3}
	\{\Gamma_{C_1}, \Gamma_{C_2}\} = \{ \Gamma_{C_2}, \Gamma_{L_1} \} = \{ \Gamma_{C_2}, \Gamma_{L_2}\} = e,
\end{equation}
\begin{equation}\label{arrangement1_final4}
	\Gamma_{C_1}\Gamma_{C_2}\Gamma_{C_1}\Gamma_{L_1}\Gamma_{C_2}\Gamma_{L_3}\Gamma_{L_2} = e.
\end{equation}
In particular, $ \Gamma_{L_3} $ is in the center of $ G_5 $.
\end{lemma}

\begin{proof}
After the substitutions $ \Gamma_{C_1'} = \Gamma_{C_1} $ and $ \Gamma_{C_2'} = \Gamma_{L_1}\Gamma_{C_1}\Gamma_{C_2}\Gamma_{C_1}\Gamma_{L_1} $ we derive that $ G_5 $ is generated by $ \Gamma_{C_1},\Gamma_{C_2},\Gamma_{L_1},\Gamma_{L_2} $, and $ \Gamma_{L_3} $, subject to:
\begin{equation}\label{arrangement1_simpl_1}
	\begin{split}
		\left[\Gamma_{C_2}\Gamma_{C_1}\Gamma_{C_2} , \Gamma_{L_3}\Gamma_{L_2}  \right]=e,
	\end{split}
\end{equation}
\begin{equation}\label{arrangement1_simpl_2}
	\begin{split}
		\left[\Gamma_{L_3} , \Gamma_{L_2}\Gamma_{C_2}\Gamma_{C_1}\Gamma_{C_2}  \right]=e,
	\end{split}
\end{equation}
\begin{equation}\label{arrangement1_simpl_3}
	\begin{split} \Gamma_{C_2}\Gamma_{L_3}\Gamma_{L_2}\Gamma_{L_3}\Gamma_{C_2}\Gamma_{C_1}\Gamma_{C_2}\Gamma_{C_1}\Gamma_{C_2}\Gamma_{L_3}\Gamma_{L_2}\Gamma_{L_3}\Gamma_{C_2}  = \Gamma_{L_1}\Gamma_{C_1}\Gamma_{C_2}\Gamma_{C_1}\Gamma_{L_1} ,
	\end{split}
\end{equation}
\begin{equation}\label{arrangement1_simpl_4}
	\begin{split}
		\left[\Gamma_{C_1}\Gamma_{C_2}\Gamma_{L_3}\Gamma_{C_2}\Gamma_{C_1} , \Gamma_{C_2}  \right]=e,
	\end{split}
\end{equation}
\begin{equation}\label{arrangement1_simpl_5}
	\begin{split} \left\{\Gamma_{L_3}\Gamma_{C_2}\Gamma_{L_1}\Gamma_{C_1}\Gamma_{C_2}\Gamma_{C_1}\Gamma_{L_1}\Gamma_{C_2}\Gamma_{L_3} , \Gamma_{L_2}  \right\}=e,
	\end{split}
\end{equation}
\begin{equation}\label{arrangement1_simpl_6}
	\begin{split} \left[\Gamma_{L_1}\Gamma_{C_1}\Gamma_{C_2}\Gamma_{C_1}\Gamma_{L_1}\Gamma_{C_2}\Gamma_{L_3}\Gamma_{L_2}\Gamma_{L_3}\Gamma_{C_2}\Gamma_{L_1}\Gamma_{C_1}\Gamma_{C_2}\Gamma_{C_1}\Gamma_{L_1} , \Gamma_{C_1}  \right]=e,
	\end{split}
\end{equation}
\begin{equation}\label{arrangement1_simpl_7}
	\begin{split}
		\left\{\Gamma_{C_1} , \Gamma_{C_2}  \right\}=e,
	\end{split}
\end{equation}
\begin{equation}\label{arrangement1_simpl_8}
	\begin{split}
		\left\{\Gamma_{L_1}\Gamma_{C_1}\Gamma_{C_2}\Gamma_{C_1}\Gamma_{L_1} , \Gamma_{C_1}  \right\}=e,
	\end{split}
\end{equation}
\begin{equation}\label{arrangement1_simpl_9}
	\begin{split} \left[\Gamma_{L_1}\Gamma_{C_1}\Gamma_{C_2}\Gamma_{C_1}\Gamma_{L_1}\Gamma_{C_2}\Gamma_{L_3}\Gamma_{L_2}\Gamma_{L_3}\Gamma_{C_2}\Gamma_{L_1}\Gamma_{C_1}\Gamma_{C_2}\Gamma_{C_1}\Gamma_{L_1} , \Gamma_{L_1}  \right]=e,
	\end{split}
\end{equation}
\begin{equation}\label{arrangement1_simpl_10}
	\begin{split} \left[\Gamma_{L_1}\Gamma_{C_1}\Gamma_{C_2}\Gamma_{C_1}\Gamma_{L_1}\Gamma_{C_1}\Gamma_{L_1}\Gamma_{C_1}\Gamma_{C_2}\Gamma_{C_1}\Gamma_{L_1} , \Gamma_{L_1}  \right]=e,
	\end{split}
\end{equation}
\begin{equation}\label{arrangement1_simpl_11}
	\begin{split}
		\left\{\Gamma_{L_1}\Gamma_{C_1}\Gamma_{C_2}\Gamma_{C_1}\Gamma_{L_1} , \Gamma_{L_1}  \right\}=e,
	\end{split}
\end{equation}
\begin{equation}\label{arrangement1_simpl_12}
	\begin{split}
		\left[\Gamma_{C_2} , \Gamma_{L_3}  \right]=e,
	\end{split}
\end{equation}
\begin{equation}\label{arrangement1_simpl_13}
	\begin{split}
		\left[\Gamma_{C_2}\Gamma_{L_3}\Gamma_{C_1}\Gamma_{L_1}\Gamma_{C_1}\Gamma_{L_3}\Gamma_{C_2} , \Gamma_{L_3}\Gamma_{C_2}\Gamma_{C_1}\Gamma_{C_2}  \right]=e,
	\end{split}
\end{equation}
\begin{equation}\label{arrangement1_simpl_14}
	\begin{split}
		\left[\Gamma_{L_3} , \Gamma_{C_2}\Gamma_{C_1}\Gamma_{L_3}\Gamma_{C_1}\Gamma_{L_1}\Gamma_{C_1}\Gamma_{L_3}\Gamma_{C_2}  \right]=e,
	\end{split}
\end{equation}
\begin{equation}\label{arrangement1_simpl_15}
	\begin{split}
		\Gamma_{C_1}\Gamma_{L_3}\Gamma_{C_1}\Gamma_{L_3}\Gamma_{C_1}  = \Gamma_{C_1} ,
	\end{split}
\end{equation}
\begin{equation}\label{arrangement1_simpl_16}
	\begin{split}
		\Gamma_{C_1}\Gamma_{C_2}\Gamma_{C_1}\Gamma_{L_1}\Gamma_{C_2}\Gamma_{L_3}\Gamma_{L_2} =e.
	\end{split}
\end{equation}

From \eqref{arrangement1_simpl_15} we get that $ \Gamma_{L_3} $ commutes with $ \Gamma_{C_1} $.
Substituting this, together with \eqref{arrangement1_simpl_12}, into \eqref{arrangement1_simpl_2} and \eqref{arrangement1_simpl_14} gives us that $ \Gamma_{L_3} $ commutes with $ \Gamma_{L_2} $ and $ \Gamma_{L_1} $, so $ \Gamma_{L_3} $ is in the center of $ G_5 $.
Omitting $ \Gamma_{L_3} $ from \eqref{arrangement1_simpl_13} yields $ [\Gamma_{C_1}, \Gamma_{L_1}] = e $.
By applying this we get $ \{ \Gamma_{L_1}, \Gamma_{C_2} \} = e $ from \eqref{arrangement1_simpl_11} and $ [\Gamma_{L_1}, \Gamma_{C_2}\Gamma_{C_1}\Gamma_{C_2}]=e $ from \eqref{arrangement1_simpl_10}.
By isolating $ \Gamma_{L_2} $ in \eqref{arrangement1_simpl_9} we get
\begin{align*}
	e & = [ \Gamma_{L_2}, \Gamma_{L_3}\Gamma_{C_2}\Gamma_{L_1}\Gamma_{C_1}\Gamma_{C_2}\Gamma_{C_1}\Gamma_{L_1}\Gamma_{C_1}\Gamma_{C_2}\Gamma_{C_1}\Gamma_{L_1}\Gamma_{C_2}\Gamma_{L_3} ] = \\
	& = [ \Gamma_{L_2}, \Gamma_{C_2}\Gamma_{C_1}\Gamma_{L_1}\Gamma_{C_2}\Gamma_{L_1}\Gamma_{C_2}\Gamma_{L_1}\Gamma_{C_1}\Gamma_{C_2} ] = \\
	& = [ \Gamma_{L_2}, \Gamma_{C_2}\Gamma_{C_1}\Gamma_{C_2}\Gamma_{L_1}\Gamma_{C_2}\Gamma_{C_1}\Gamma_{C_2} ] = \\
	& = [ \Gamma_{L_2}, \Gamma_{L_1} ],
\end{align*}
where in the first equality we used the fact that $ \Gamma_{L_3} $ is in the center, as well as the fact that $ \Gamma_{C_1} $ and $ \Gamma_{L_1} $ commute; in the second equality we used $ \{ \Gamma_{C_2}, \Gamma_{L_1} \} = e $; and in the third equality we used $ [\Gamma_{L_1}, \Gamma_{C_2}\Gamma_{C_1}\Gamma_{C_2}] = e $.

By isolating $ \Gamma_{L_2} $ in \eqref{arrangement1_simpl_6} we get
\begin{align*}
	e &= [\Gamma_{L_2}, \Gamma_{L_3}\Gamma_{C_2}\Gamma_{L_1}\Gamma_{C_1}\Gamma_{C_2}\Gamma_{C_1}\Gamma_{L_1}\Gamma_{C_1}\Gamma_{L_1}\Gamma_{C_1}\Gamma_{C_2}\Gamma_{C_1}\Gamma_{L_1}\Gamma_{C_2}\Gamma_{L_3}] = \\
	& =  [\Gamma_{L_2}, \Gamma_{C_2}\Gamma_{L_1}\Gamma_{C_1}\Gamma_{C_2}\Gamma_{C_1}\Gamma_{C_2}\Gamma_{C_1}\Gamma_{L_1}\Gamma_{C_2}] = \\
	& =  [\Gamma_{L_2}, \Gamma_{C_2}\Gamma_{L_1}\Gamma_{C_2}\Gamma_{C_1}\Gamma_{C_2}\Gamma_{L_1}\Gamma_{C_2}] = \\
	& =  [\Gamma_{L_2}, \Gamma_{C_1}],
\end{align*}
where again we used known commutation relations in the first equality, $ {\{\Gamma_{C_1}, \Gamma_{C_2}\} = e} $ in the second equality, and $ [\Gamma_{L_1}, \Gamma_{C_2}\Gamma_{C_1}\Gamma_{C_2}] = e $ in the third equality.

Next we deal with \eqref{arrangement1_simpl_3}. By eliminating $ \Gamma_{L_3} $ and rearranging, we get
\begin{align*}
	& \Gamma_{C_2}\Gamma_{C_1}\Gamma_{C_2}\Gamma_{C_1}\Gamma_{C_2} = \\
	=& \Gamma_{L_2}\Gamma_{C_2}\Gamma_{L_1}\Gamma_{C_1}\Gamma_{C_2}\Gamma_{C_1}\Gamma_{L_1}\Gamma_{C_2}\Gamma_{L_2} = \\
	=& \Gamma_{L_2}\Gamma_{C_2}\Gamma_{C_1}\Gamma_{L_1}\Gamma_{C_2}\Gamma_{L_1}\Gamma_{C_1}\Gamma_{C_2}\Gamma_{L_2} = \\
	=& \Gamma_{L_2}\Gamma_{C_2}\Gamma_{C_1}\Gamma_{C_2}\Gamma_{L_1}\Gamma_{C_2}\Gamma_{L_1}\Gamma_{C_2}\Gamma_{C_1}\Gamma_{C_2}\Gamma_{L_2} = \\
	=& \Gamma_{L_2}\Gamma_{L_1}\Gamma_{C_2}\Gamma_{C_1}\Gamma_{C_2}\Gamma_{C_1}\Gamma_{C_2}\Gamma_{L_1}\Gamma_{L_2} = \\
	=& \Gamma_{L_2}\Gamma_{L_1}\Gamma_{C_1}\Gamma_{C_2}\Gamma_{C_1}\Gamma_{L_1}\Gamma_{L_2}.
\end{align*}
By $ \{ \Gamma_{C_1}, \Gamma_{C_2} \} =e $, the left side of the above equation is equal to $ \Gamma_{C_1}\Gamma_{C_2}\Gamma_{C_1} $, and by using the commutation relations $ [\Gamma_{C_1},\Gamma_{L_1}]=[\Gamma_{C_1}, \Gamma_{L_2}]=e $, the above equation is equivalent to $ [\Gamma_{C_2}, \Gamma_{L_1}\Gamma_{L_2}]=e $.

From \eqref{arrangement1_simpl_1} we get $ [\Gamma_{L_2}, \Gamma_{C_2}\Gamma_{C_1}\Gamma_{C_2}]=e $. Using this together with $ \{ \Gamma_{C_2}, \Gamma_{L_1} \}=e $ and $ [\Gamma_{L_1}, \Gamma_{L_2}]=e $, we can rearrange \eqref{arrangement1_simpl_5} and derive
\begin{align*}
	e & = \{ \Gamma_{C_1}\Gamma_{C_2}\Gamma_{L_2}\Gamma_{C_2}\Gamma_{C_1}, \Gamma_{L_1}\Gamma_{C_2}\Gamma_{L_1} \} = \\
	& = \{ \Gamma_{C_2}\Gamma_{C_1}\Gamma_{C_2}\Gamma_{L_2}\Gamma_{C_2}\Gamma_{C_1}\Gamma_{C_2}, \Gamma_{C_2}\Gamma_{L_1}\Gamma_{C_2}\Gamma_{L_1}\Gamma_{C_2} \} = \\
	& = \{ \Gamma_{L_2}, \Gamma_{L_1}\Gamma_{C_2}\Gamma_{L_1} \} =
	  \{ \Gamma_{L_2}, \Gamma_{C_2} \}.
\end{align*}

Combining those derivations we get the desired result.
\end{proof}

\newpage

\subsection{$ \pcpt{\mathcal{B}_6} $}

\label{section:braids_computation_case2}

\begin{figure}[H]
	\begin{center}

	\end{center}
\end{figure}

{Monodromy table}

\begin{tabular}{| c | c | c | c |}
	\hline \hline
	Vertex number& Vertex description& Skeleton& Diffeomorphism\\ \hline \hline
	1& $C_1$ branch& $ \langle 3 - 4 \rangle $& $ \Delta_{I_{2}I_{4}}^{1/2}\langle 3 \rangle $ \\ \hline
	2& Node between $C_1$, $L_4$ and $L_2$& $ \langle 1 - 2 - 3 \rangle $& $ \Delta \langle 1,2,3 \rangle $ \\ \hline
	3& $C_2$ branch& $ \langle 3 - 4 \rangle $& $ \Delta_{\mathbb{R}I_{2}}^{1/2}\langle 3 \rangle $ \\ \hline
	4& Node between $C_2$ and $L_4$& $ \langle 2 - 3 \rangle $& $ \Delta\langle 2, 3 \rangle $ \\ \hline
	5& Tangency between $L_2$ and $C_2$& $ \langle 4 - 5 \rangle $& $ \Delta^{2}\langle 4, 5 \rangle $ \\ \hline
	6& Node between $C_1$ and $L_2$& $ \langle 5 - 6 \rangle $& $ \Delta\langle 5, 6 \rangle $ \\ \hline
	7& Tangency between $C_1$ and $C_2$& $ \langle 1 - 2 \rangle $& $ \Delta^{2}\langle 1, 2 \rangle $ \\ \hline
	8& Tangency between $C_1$ and $C_2$& $ \langle 4 - 5 \rangle $& $ \Delta^{2}\langle 4, 5 \rangle $ \\ \hline
	9& Node between $L_1$ and $L_2$& $ \langle 6 - 7 \rangle $& $ \Delta\langle 6, 7 \rangle $ \\ \hline
	10& Node between $C_2$ and $L_4$& $ \langle 3 - 4 \rangle $& $ \Delta\langle 3, 4 \rangle $ \\ \hline
	11& Node between $L_1$, $C_1$ and $L_4$& $ \langle 4 - 5 - 6 \rangle $& $ \Delta \langle 4,5,6 \rangle $ \\ \hline
	12& Tangency between $L_1$ and $C_2$& $ \langle 3 - 4 \rangle $& $ \Delta^{2}\langle 3, 4 \rangle $ \\ \hline
	13& $C_2$ branch& $ \langle 2 - 3 \rangle $& $ \Delta_{I_{2}\mathbb{R}}^{1/2}\langle 2 \rangle $ \\ \hline
	14& Node between $L_1$ and $C_1$& $ \langle 1 - 2 \rangle $& $ \Delta\langle 1, 2 \rangle $ \\ \hline
	15& $C_1$ branch& $ \langle 2 - 3 \rangle $& $ \Delta_{I_{4}I_{2}}^{1/2}\langle 2 \rangle $
	\\ \hline \hline
\end{tabular}

\subsubsection{Vertex Number 1}


\begin{figure}[H]
	\begin{center}
		%
	\end{center}
\end{figure}

Relation:

\begin{equation*}\begin{split}
\Gamma_{C_1}  = \Gamma_{C_1'}
\end{split}
\end{equation*}
\subsubsection{Vertex Number 2}

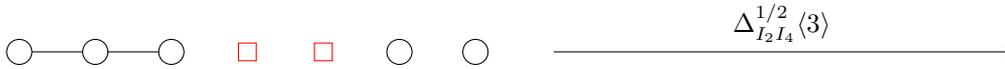
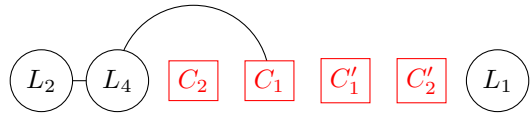
\begin{figure}[H]
	\begin{center}
		%
	\end{center}
\end{figure}

Relation:

\begin{equation*}\begin{split}
\left( \Gamma_{C_2}^{-1}\Gamma_{C_1}\Gamma_{C_2}  \right) \left( \Gamma_{L_4}  \right) \left( \Gamma_{L_2}  \right)& = \\
\left( \Gamma_{L_4}  \right) \left( \Gamma_{L_2}  \right) \left( \Gamma_{C_2}^{-1}\Gamma_{C_1}\Gamma_{C_2}  \right)& = \\
\left( \Gamma_{L_2}  \right) \left( \Gamma_{C_2}^{-1}\Gamma_{C_1}\Gamma_{C_2}  \right) \left( \Gamma_{L_4}  \right)
\end{split}
\end{equation*}
\subsubsection{Vertex Number 3}

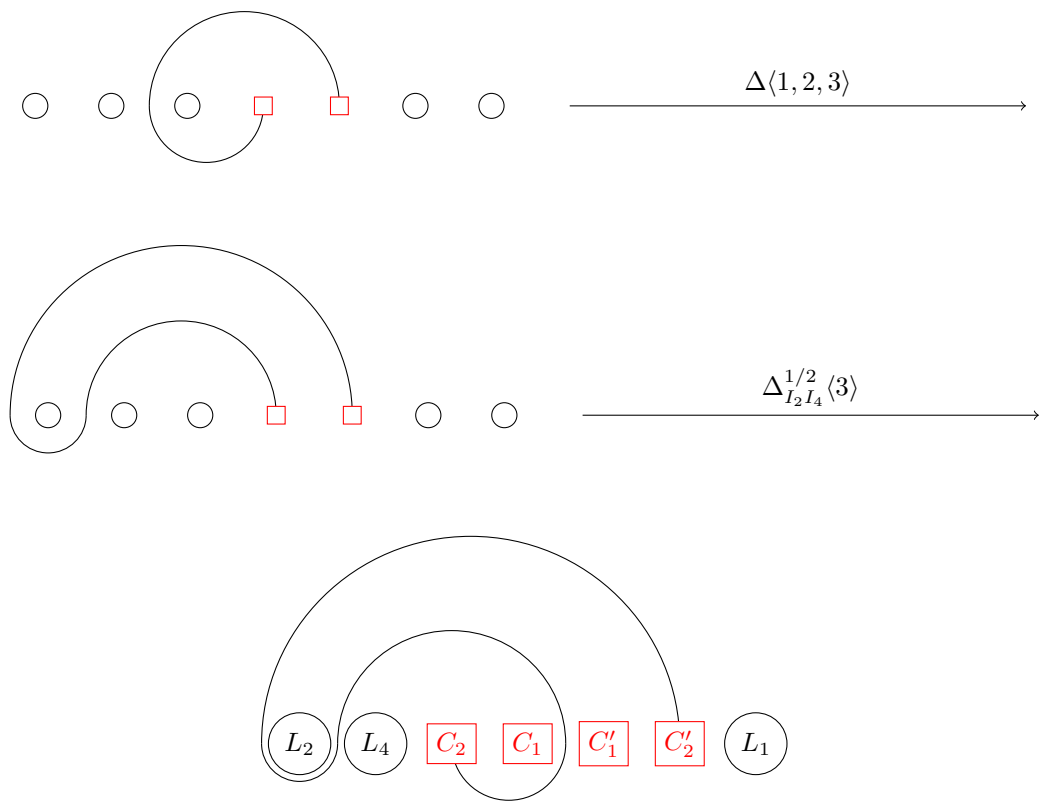
\begin{figure}[H]
	\begin{center}
		%
	\end{center}
\end{figure}

Relation:

\begin{equation*}\begin{split}
\Gamma_{C_1'}\Gamma_{C_1}\Gamma_{C_2}\Gamma_{L_4}\Gamma_{L_2}\Gamma_{L_4}^{-1}\Gamma_{C_2}^{-1}\Gamma_{C_1}^{-1}\Gamma_{C_2}\Gamma_{C_1}\Gamma_{C_2}\Gamma_{L_4}\Gamma_{L_2}^{-1}\Gamma_{L_4}^{-1}\Gamma_{C_2}^{-1}\Gamma_{C_1}^{-1}\Gamma_{C_1'}^{-1}  = \Gamma_{C_2'}
\end{split}
\end{equation*}
\subsubsection{Vertex Number 4}

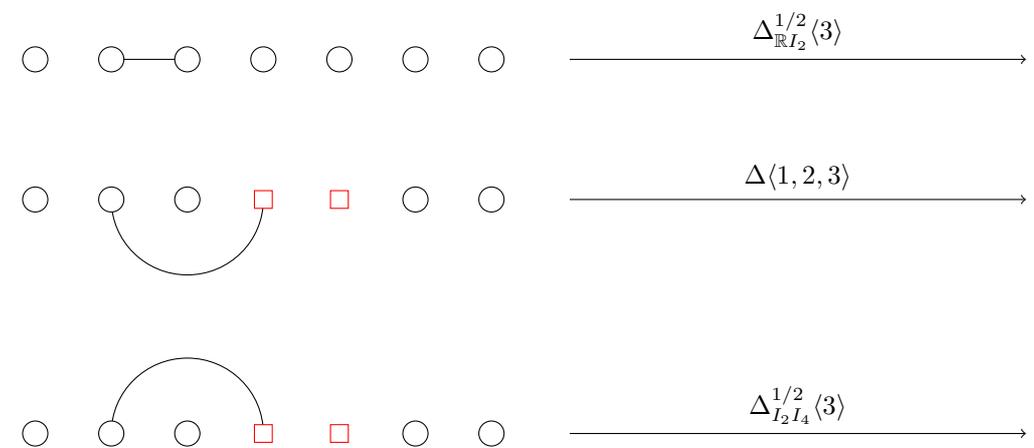
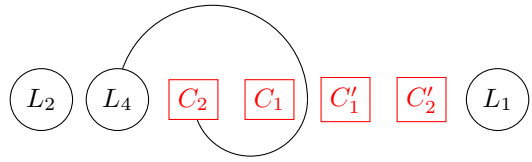
\begin{figure}[H]
	\begin{center}
		%
	\end{center}
\end{figure}

Relation:

\begin{equation*}\begin{split}
\left[\Gamma_{C_1}\Gamma_{C_2}\Gamma_{L_4}\Gamma_{C_2}^{-1}\Gamma_{C_1}^{-1} , \Gamma_{C_2} \right]=e
\end{split}
\end{equation*}
\subsubsection{Vertex Number 5}

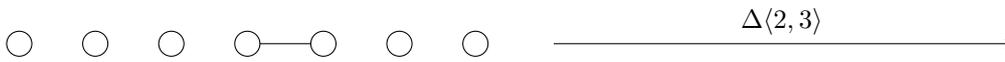
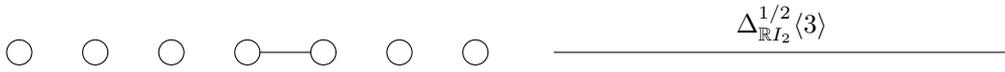
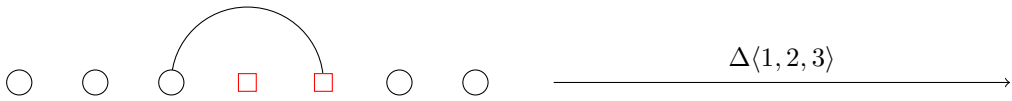
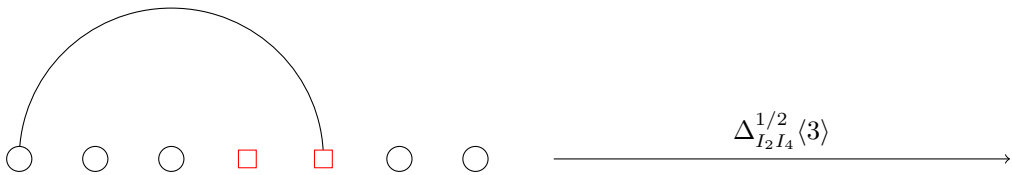
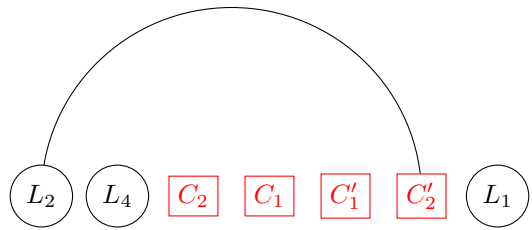
\begin{figure}[H]
	\begin{center}
		%
	\end{center}
\end{figure}

Relation:

\begin{equation*}\begin{split}
\left\{\Gamma_{L_4}^{-1}\Gamma_{C_2}^{-1}\Gamma_{C_1}^{-1}\Gamma_{C_1'}^{-1}\Gamma_{C_2'}\Gamma_{C_1'}\Gamma_{C_1}\Gamma_{C_2}\Gamma_{L_4}  , \Gamma_{L_2} \right\}=e
\end{split}
\end{equation*}
\subsubsection{Vertex Number 6}

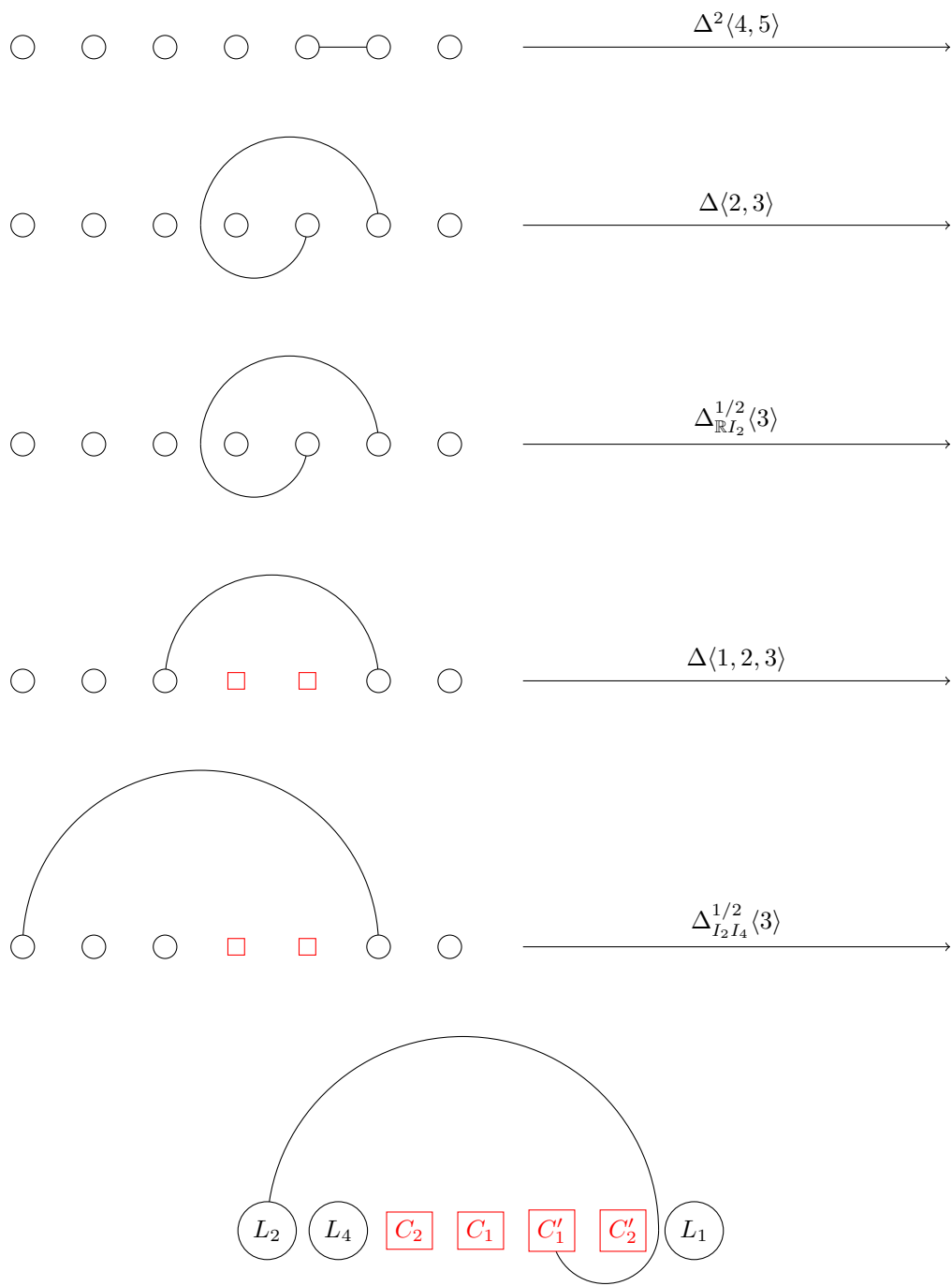
\begin{figure}[H]
	\begin{center}
		%
	\end{center}
\end{figure}

Relation:

\begin{equation*}\begin{split}
\left[\Gamma_{C_2'}\Gamma_{C_1'}\Gamma_{C_1}\Gamma_{C_2}\Gamma_{L_4}\Gamma_{L_2}\Gamma_{L_4}^{-1}\Gamma_{C_2}^{-1}\Gamma_{C_1}^{-1}\Gamma_{C_1'}^{-1}\Gamma_{C_2'}^{-1} , \Gamma_{C_1'} \right]=e
\end{split}
\end{equation*}
\subsubsection{Vertex Number 7}

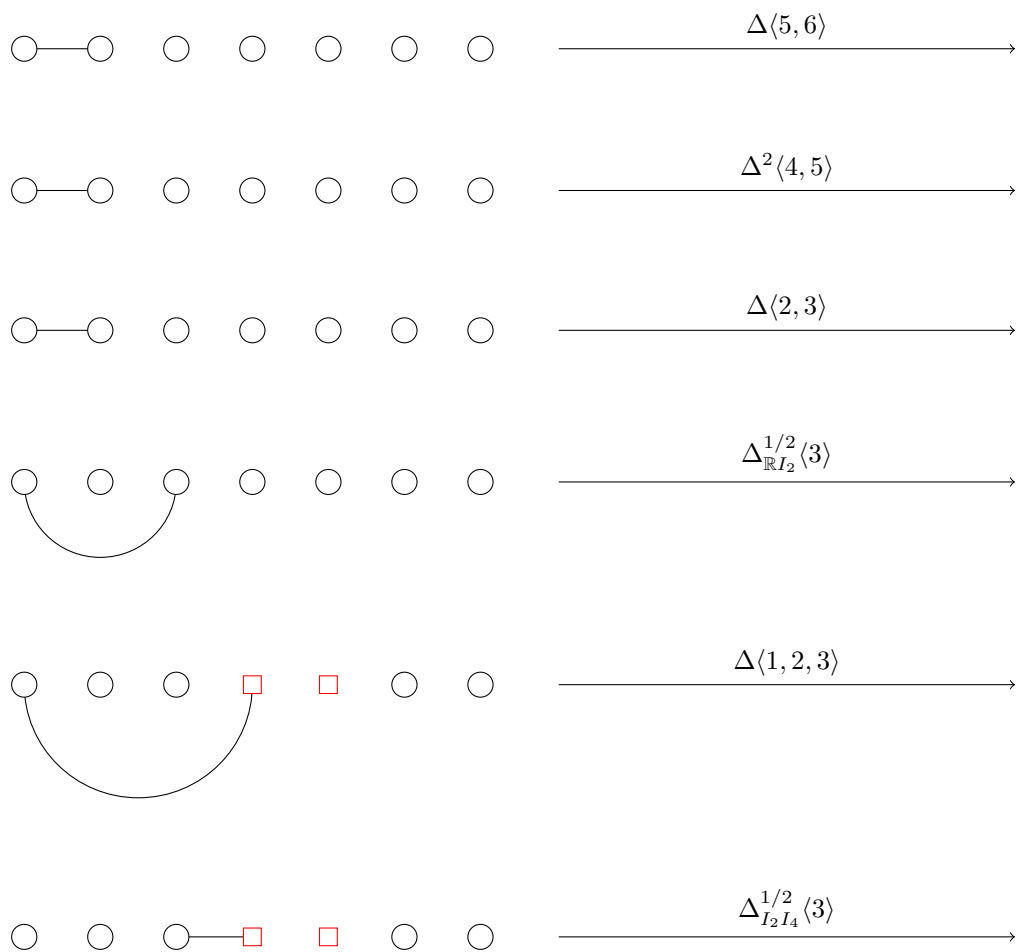
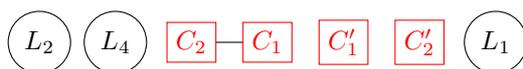
\begin{figure}[H]
	\begin{center}
		%
	\end{center}
\end{figure}

Relation:

\begin{equation*}\begin{split}
\left\{\Gamma_{C_1}  , \Gamma_{C_2} \right\}=e
\end{split}
\end{equation*}
\subsubsection{Vertex Number 8}

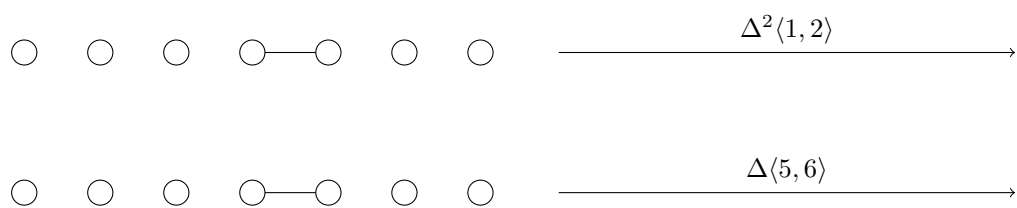
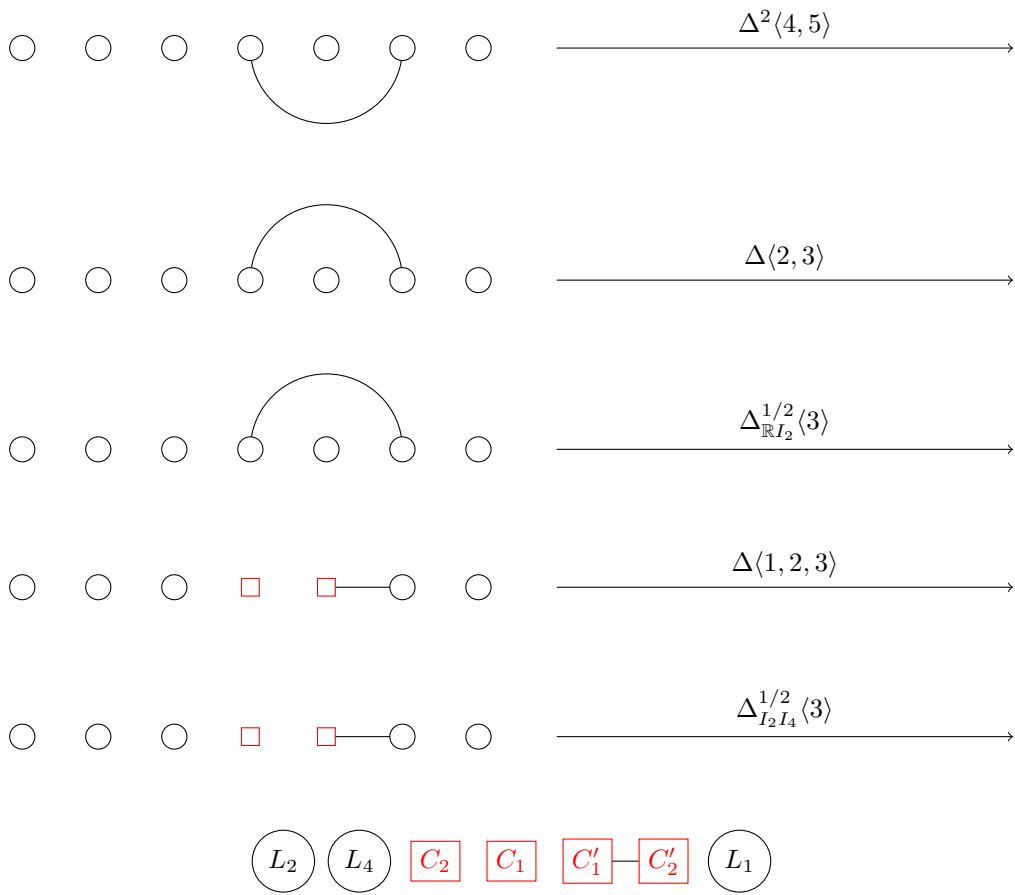
\begin{figure}[H]
	\begin{center}
		%
	\end{center}
\end{figure}

Relation:

\begin{equation*}\begin{split}
\left\{\Gamma_{C_2'}  , \Gamma_{C_1'} \right\}=e
\end{split}
\end{equation*}
\subsubsection{Vertex Number 9}

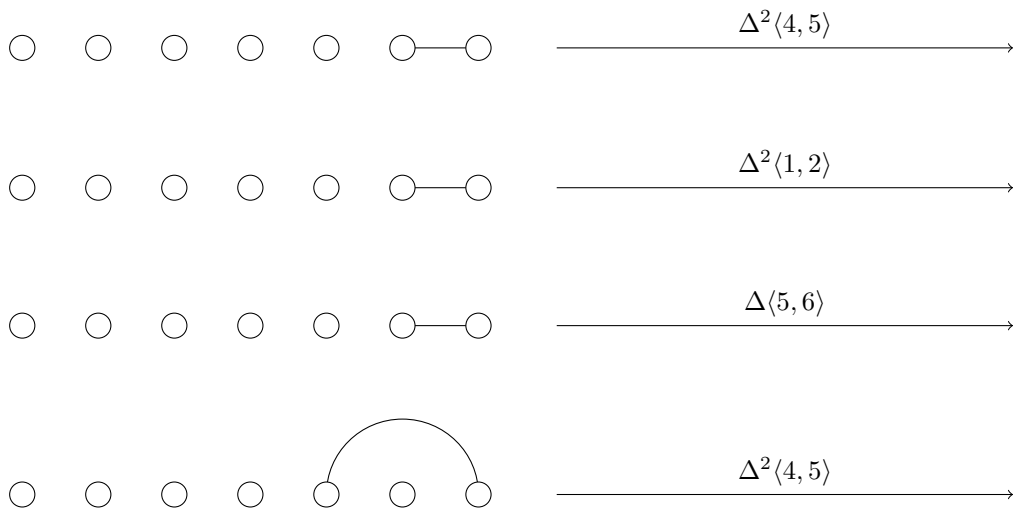
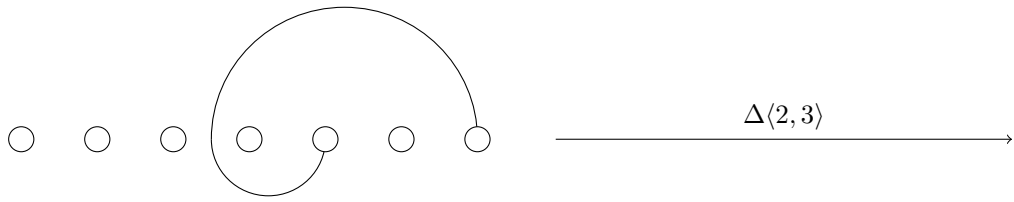
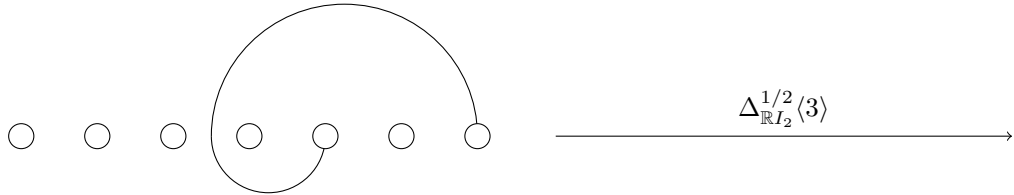
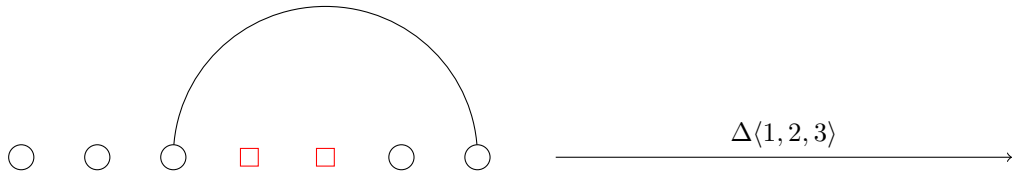
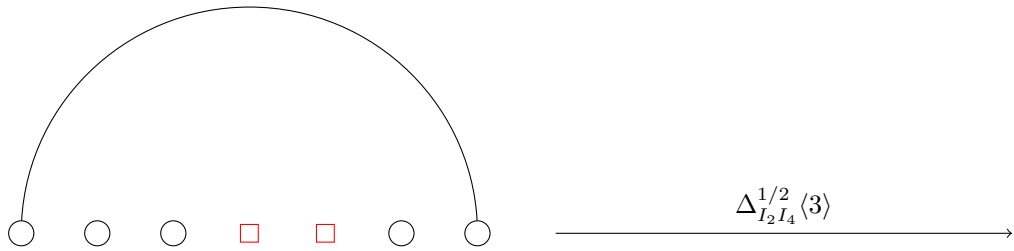
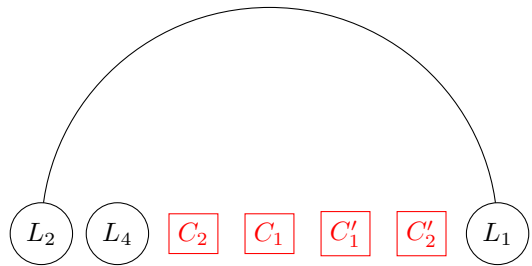
\begin{figure}[H]
	\begin{center}
		%
	\end{center}
\end{figure}

Relation:

\begin{equation*}\begin{split}
\left[\Gamma_{C_2'}\Gamma_{C_1'}\Gamma_{C_1}\Gamma_{C_2}\Gamma_{L_4}\Gamma_{L_2}\Gamma_{L_4}^{-1}\Gamma_{C_2}^{-1}\Gamma_{C_1}^{-1}\Gamma_{C_1'}^{-1}\Gamma_{C_2'}^{-1} , \Gamma_{L_1} \right]=e
\end{split}
\end{equation*}
\subsubsection{Vertex Number 10}

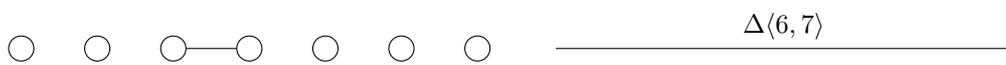
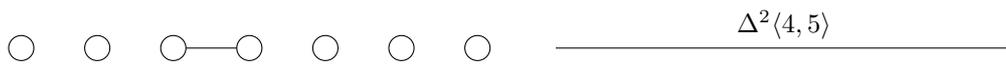
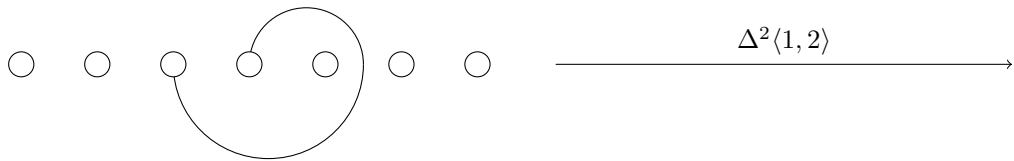
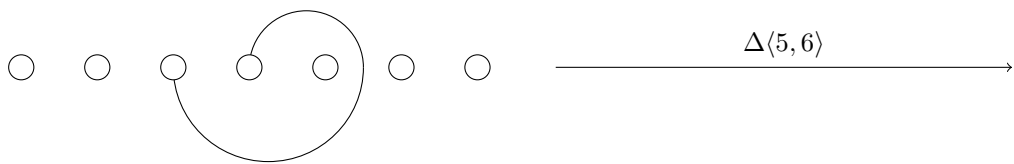
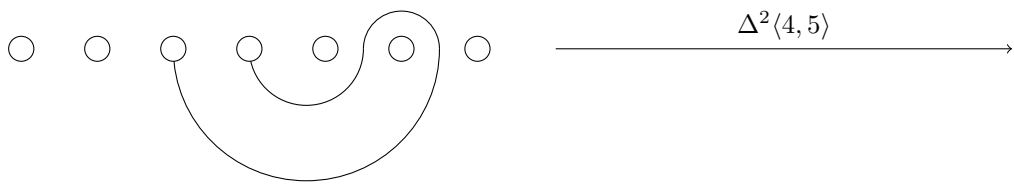
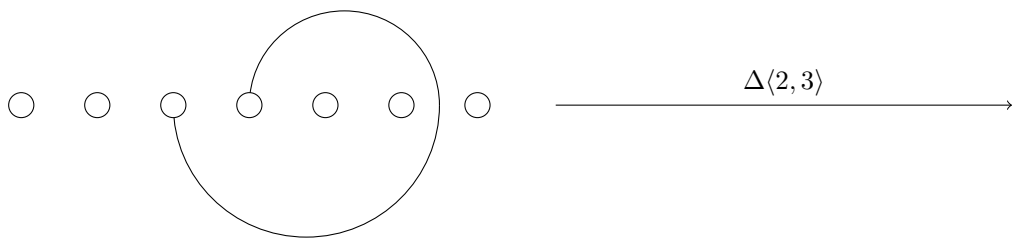
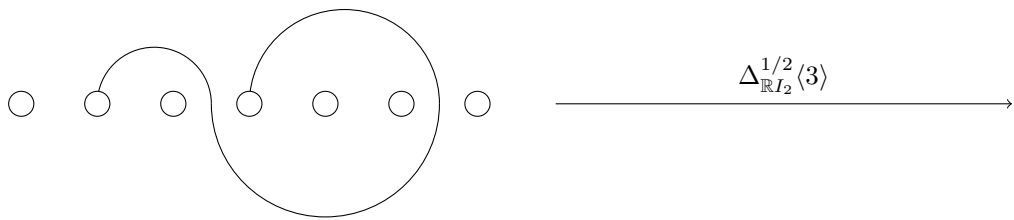
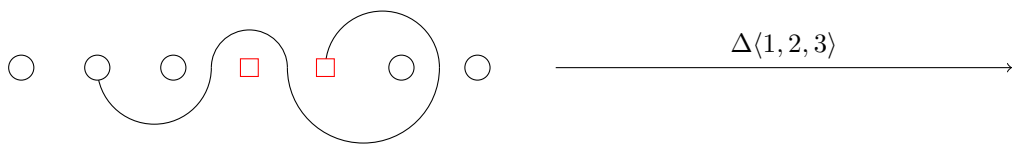
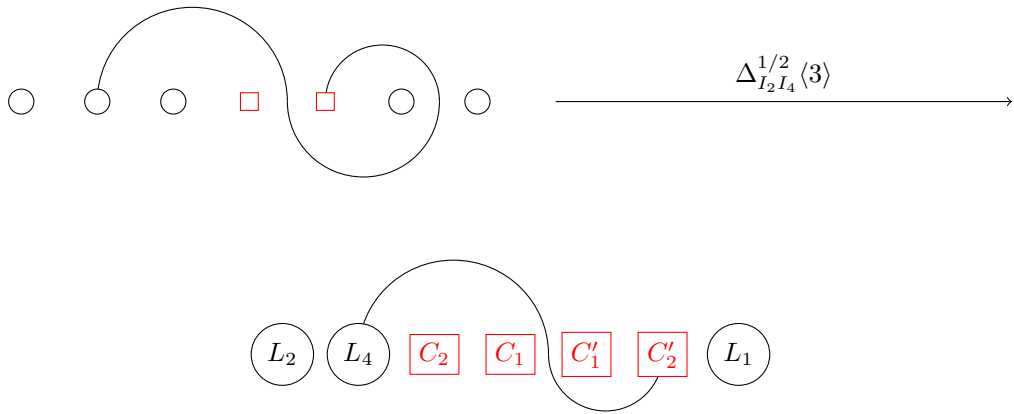
\begin{figure}[H]
	\begin{center}
		%
	\end{center}
\end{figure}

Relation:

\begin{equation*}\begin{split}
\left[\Gamma_{C_1}\Gamma_{C_2}\Gamma_{L_4}\Gamma_{C_2}^{-1}\Gamma_{C_1}^{-1} , \Gamma_{C_2'} \right]=e
\end{split}
\end{equation*}
\subsubsection{Vertex Number 11}

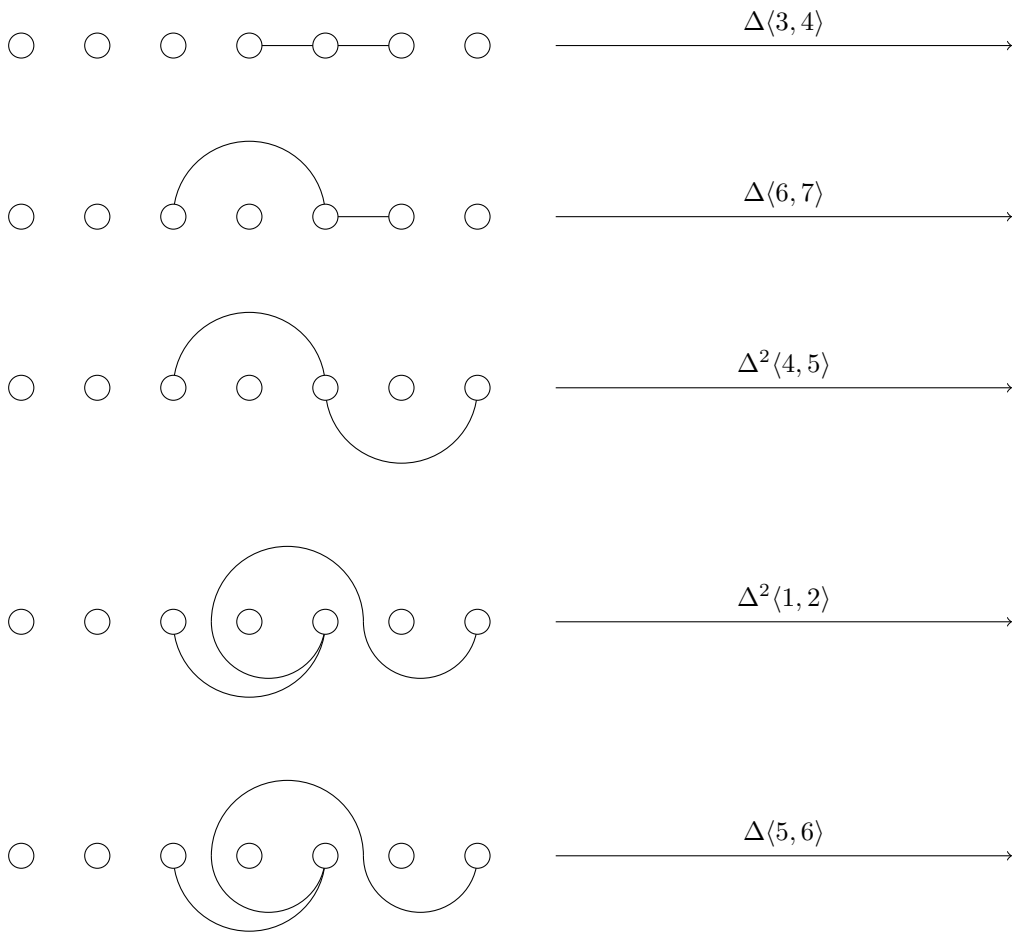
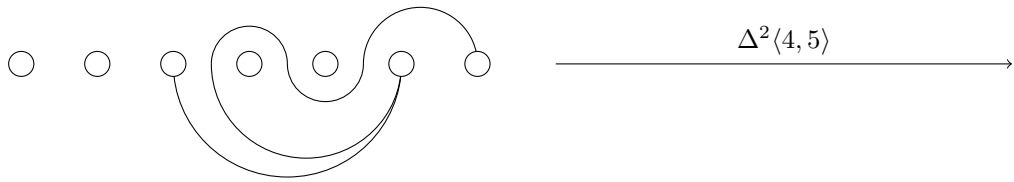
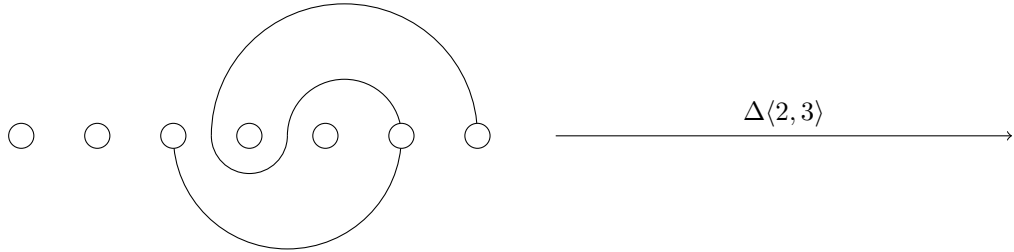
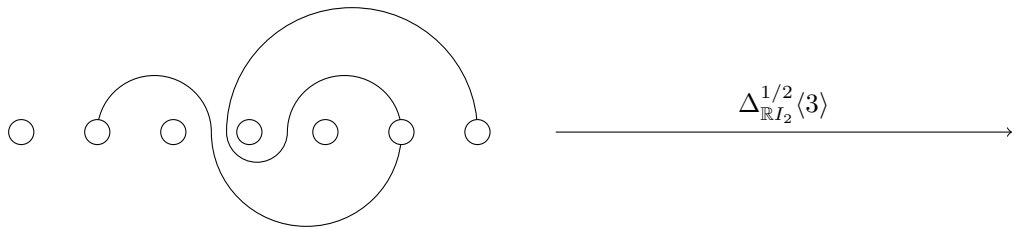
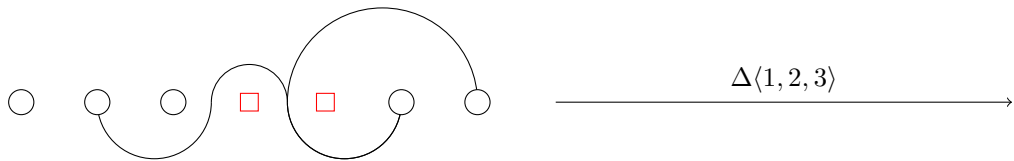
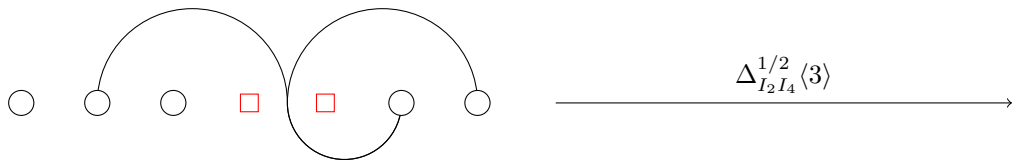
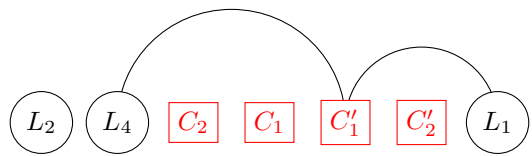
\begin{figure}[H]
	\begin{center}
		%
	\end{center}
\end{figure}

Relation:

\begin{equation*}\begin{split}
\left( \Gamma_{C_2'}^{-1}\Gamma_{L_1}\Gamma_{C_2'}  \right) \left( \Gamma_{C_1'}  \right) \left( \Gamma_{C_1}\Gamma_{C_2}\Gamma_{L_4}\Gamma_{C_2}^{-1}\Gamma_{C_1}^{-1}  \right)& = \\
\left( \Gamma_{C_1'}  \right) \left( \Gamma_{C_1}\Gamma_{C_2}\Gamma_{L_4}\Gamma_{C_2}^{-1}\Gamma_{C_1}^{-1}  \right) \left( \Gamma_{C_2'}^{-1}\Gamma_{L_1}\Gamma_{C_2'}  \right)& = \\
\left( \Gamma_{C_1}\Gamma_{C_2}\Gamma_{L_4}\Gamma_{C_2}^{-1}\Gamma_{C_1}^{-1}  \right) \left( \Gamma_{C_2'}^{-1}\Gamma_{L_1}\Gamma_{C_2'}  \right) \left( \Gamma_{C_1'}  \right)
\end{split}
\end{equation*}
\subsubsection{Vertex Number 12}

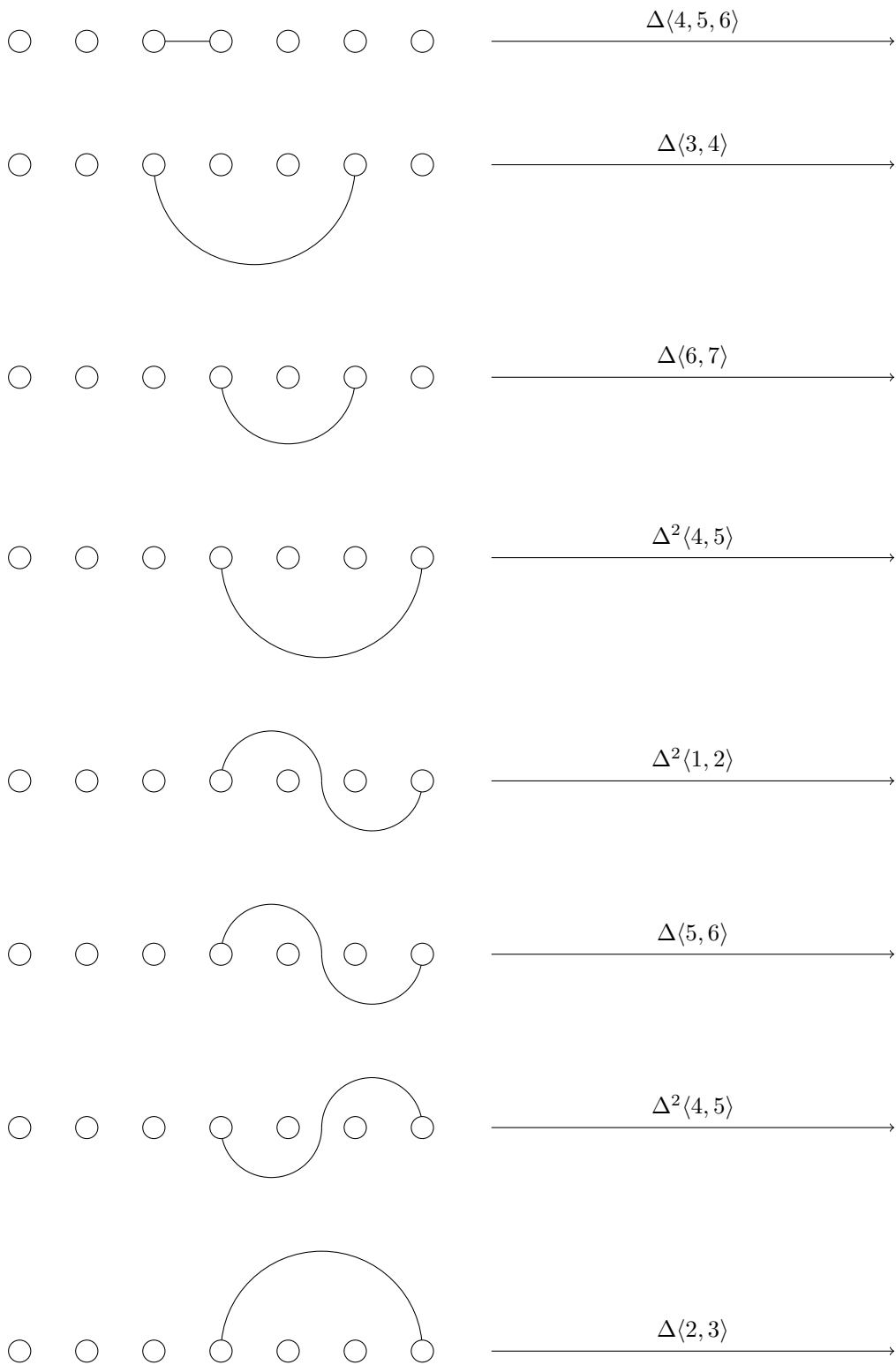
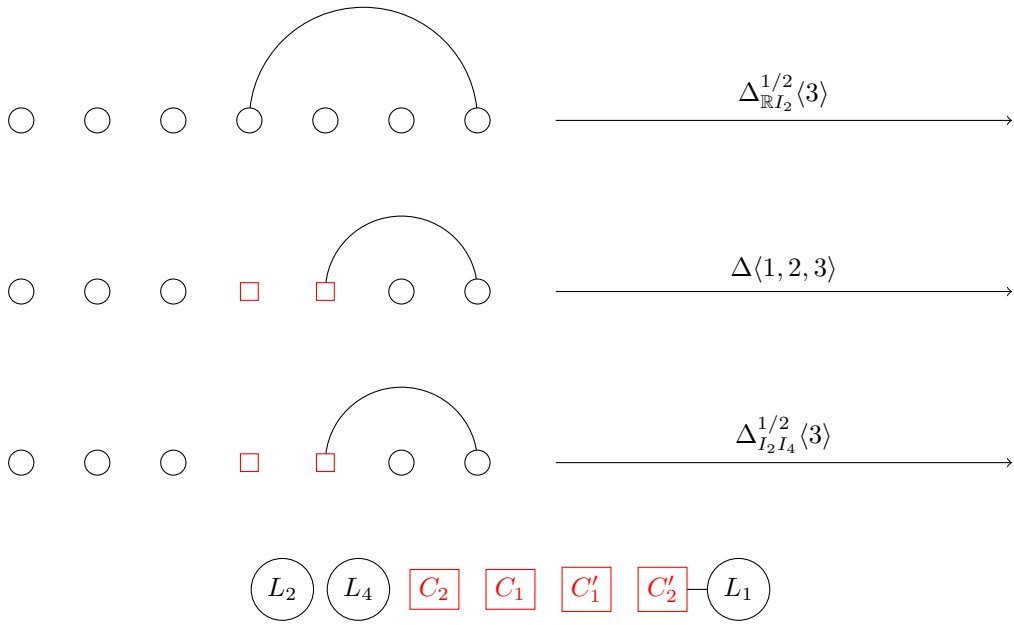
\begin{figure}[H]
	\begin{center}

	\end{center}
\end{figure}

Relation:

\begin{equation*}\begin{split}
\left\{\Gamma_{C_2'}  , \Gamma_{L_1} \right\}=e
\end{split}
\end{equation*}
\subsubsection{Vertex Number 13}

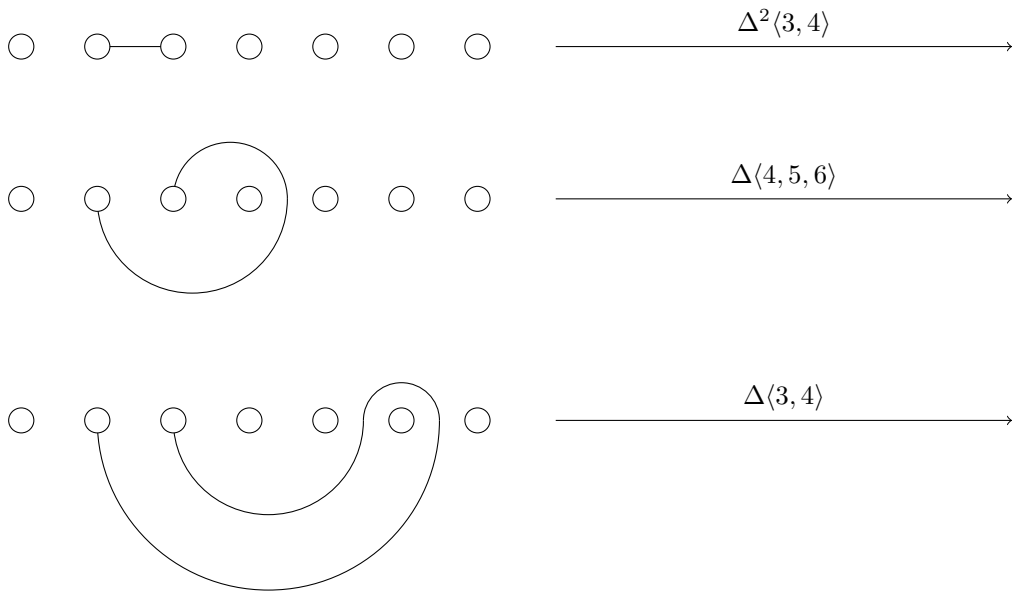
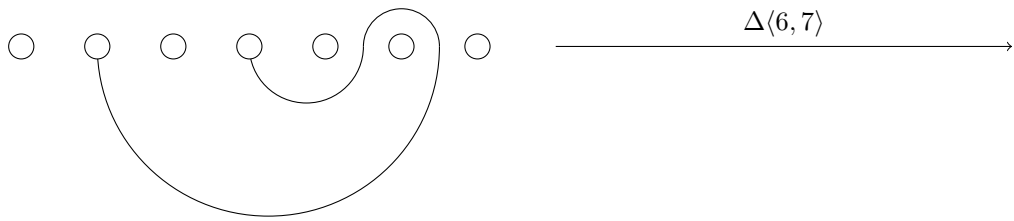
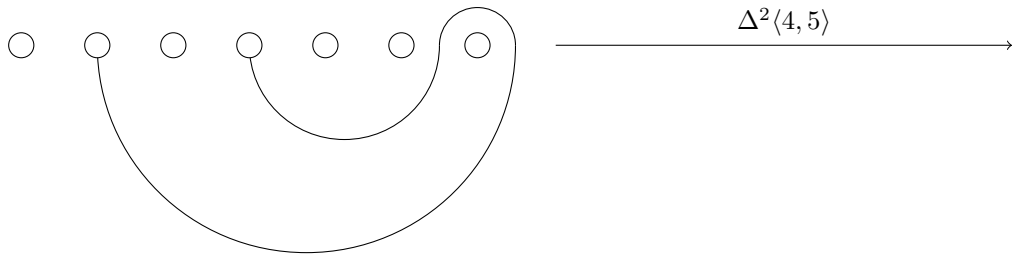
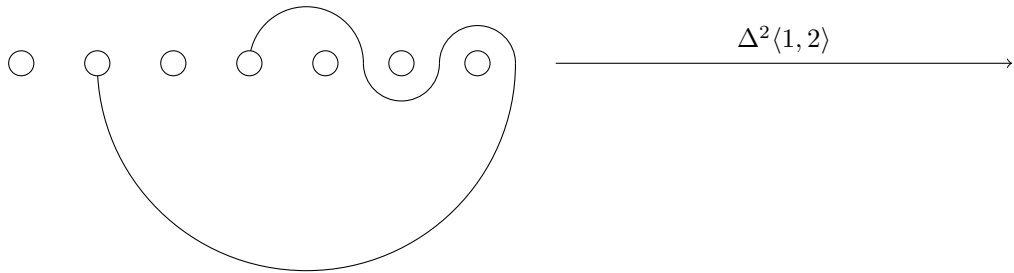
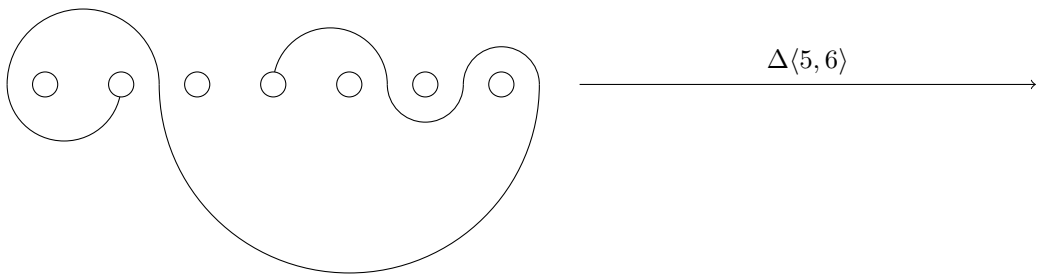
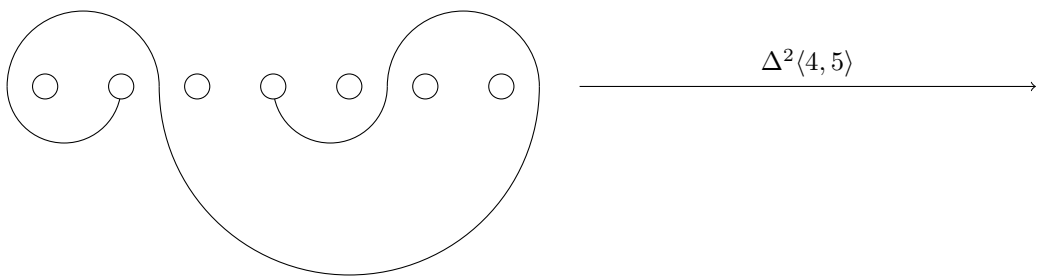
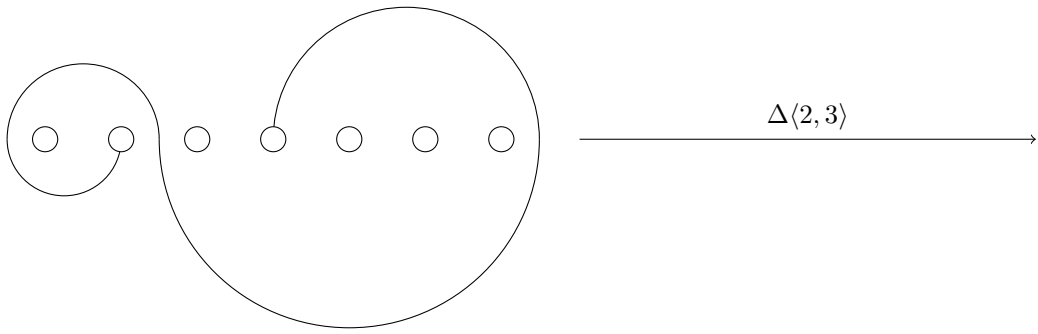
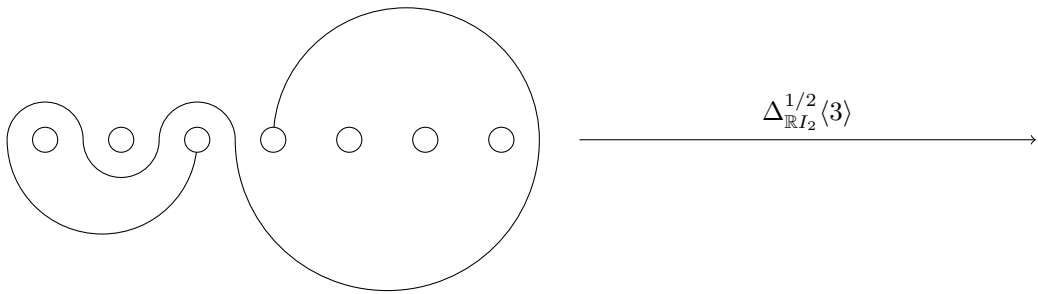
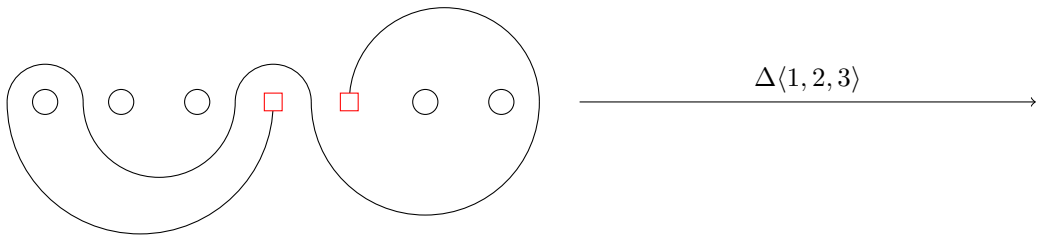
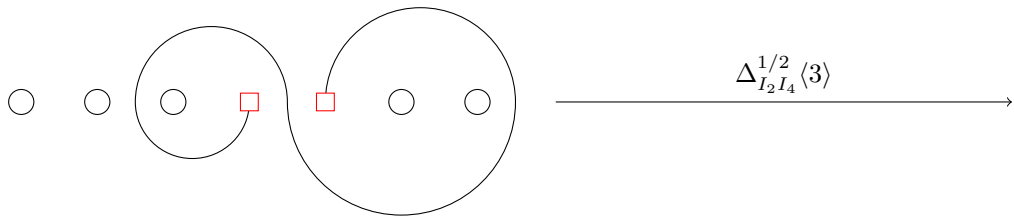
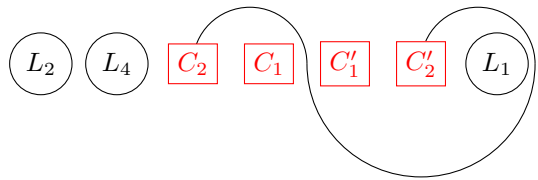
\begin{figure}[H]
	\begin{center}
		%
	\end{center}
\end{figure}

Relation:

\begin{equation*}\begin{split}
\Gamma_{L_1}^{-1}\Gamma_{C_1}\Gamma_{C_2}\Gamma_{C_1}^{-1}\Gamma_{L_1}  = \Gamma_{C_2'}
\end{split}
\end{equation*}
\subsubsection{Vertex Number 14}

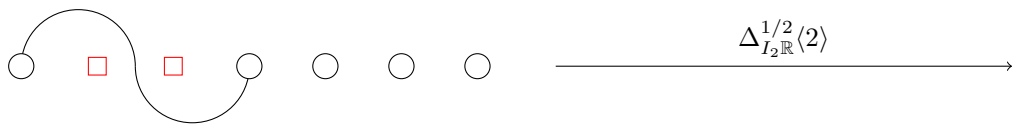
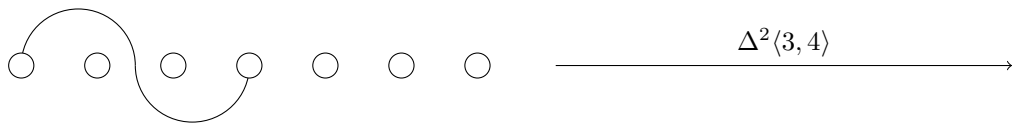
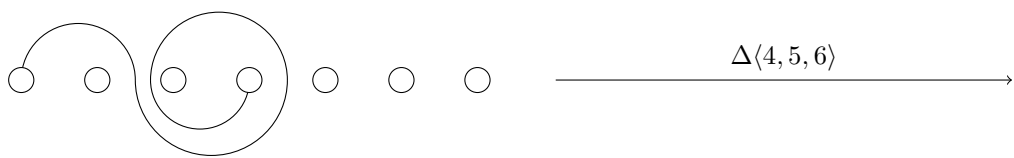
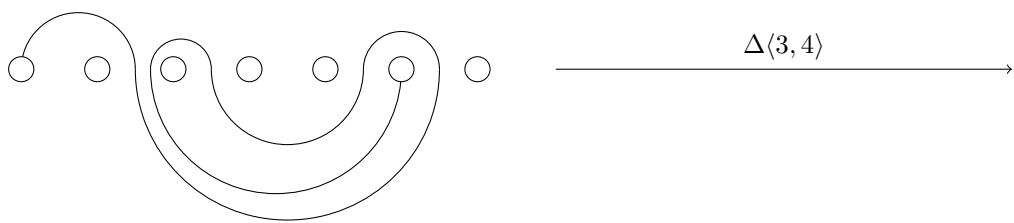
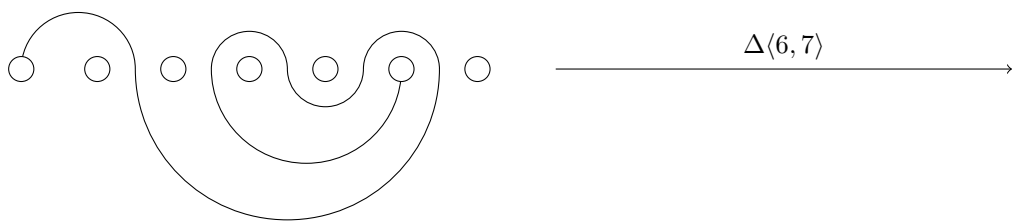
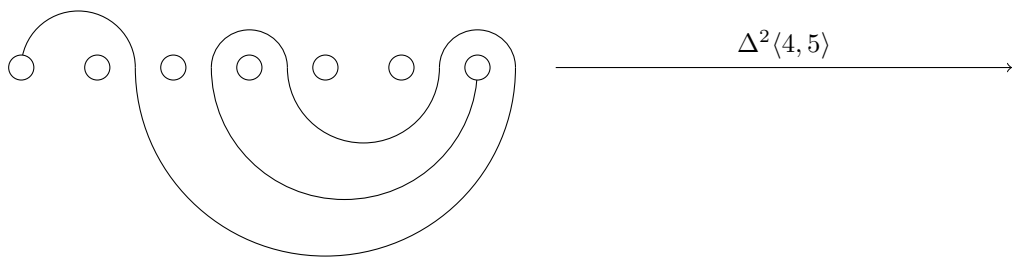
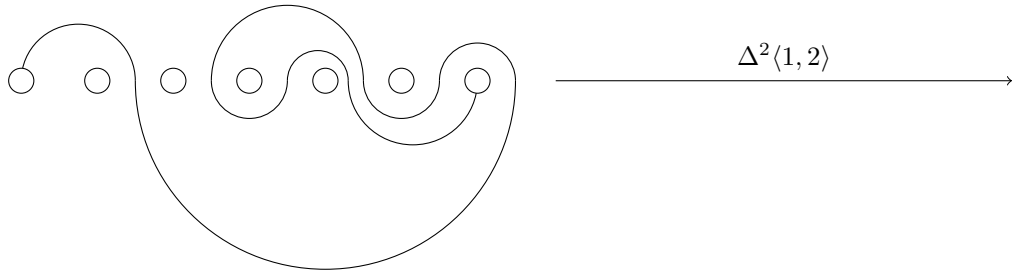
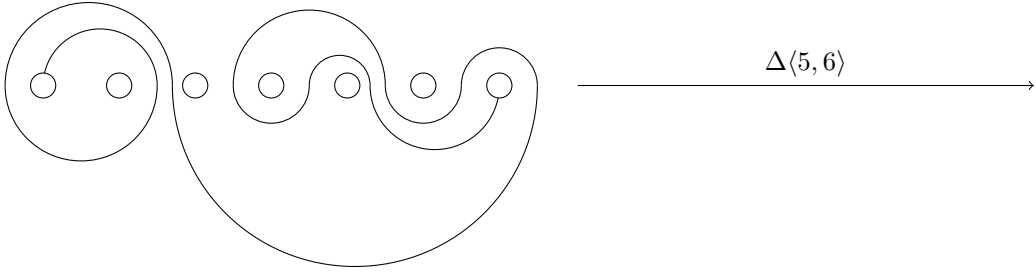
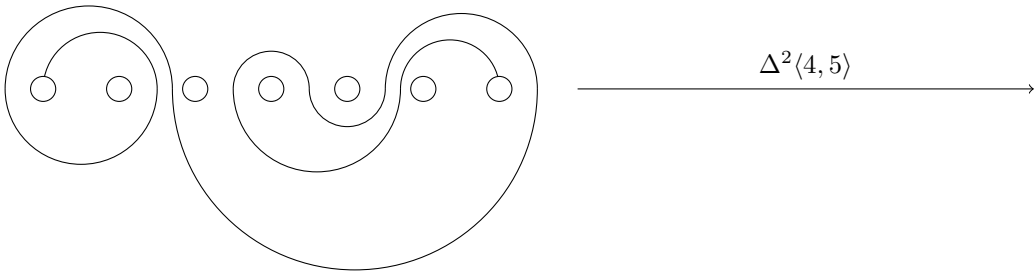
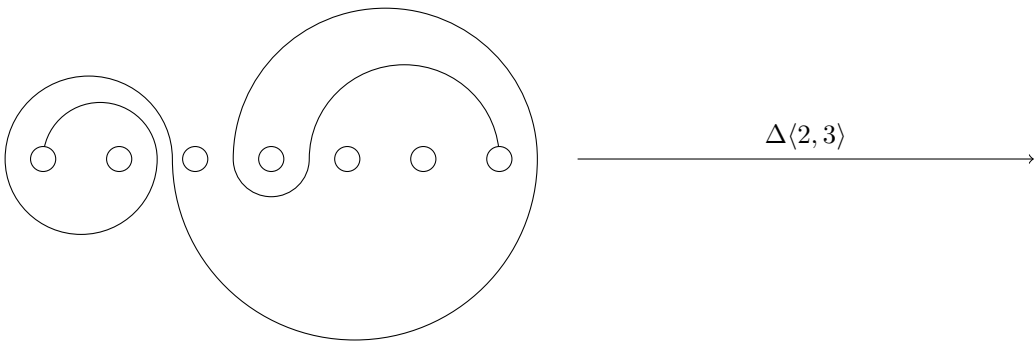
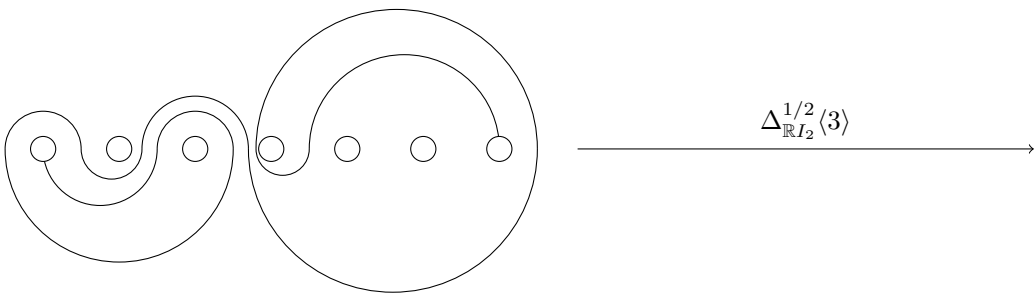
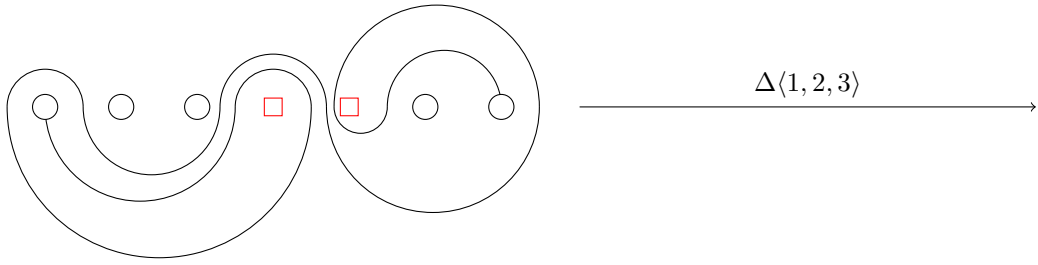
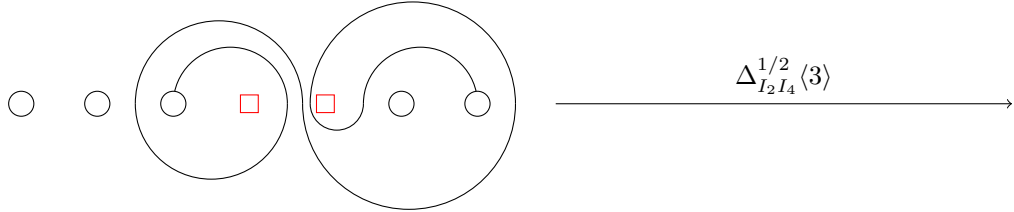
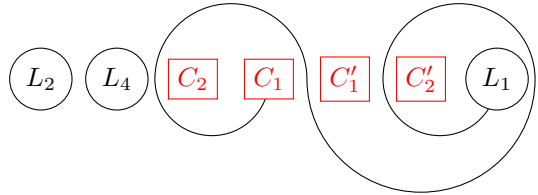
\begin{figure}[H]
	\begin{center}
		%
	\end{center}
\end{figure}

Relation:

\begin{equation*}\begin{split}
\left[\Gamma_{C_2'}^{-1}\Gamma_{L_1}^{-1}\Gamma_{C_1}\Gamma_{C_2}\Gamma_{C_1}\Gamma_{C_2}^{-1}\Gamma_{C_1}^{-1}\Gamma_{L_1}\Gamma_{C_2'} , \Gamma_{L_1} \right]=e
\end{split}
\end{equation*}
\subsubsection{Vertex Number 15}

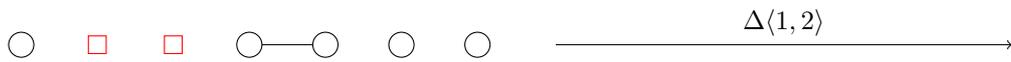
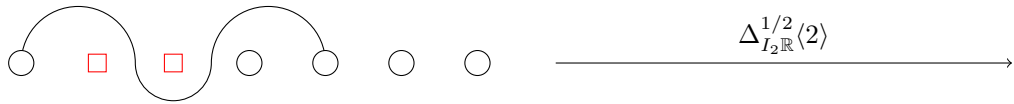
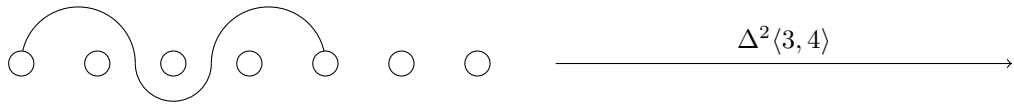
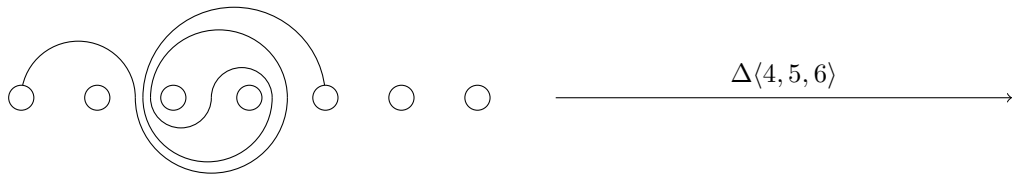
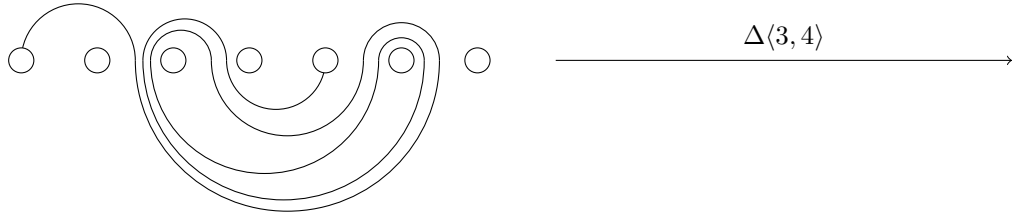
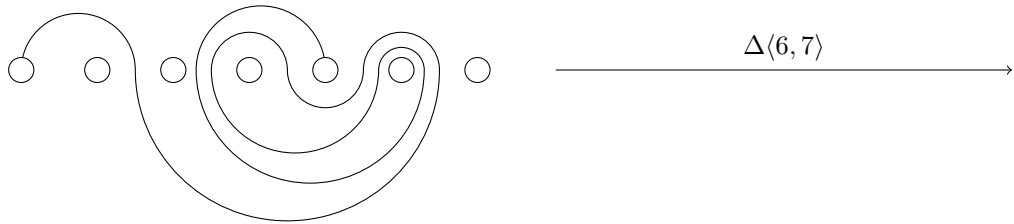
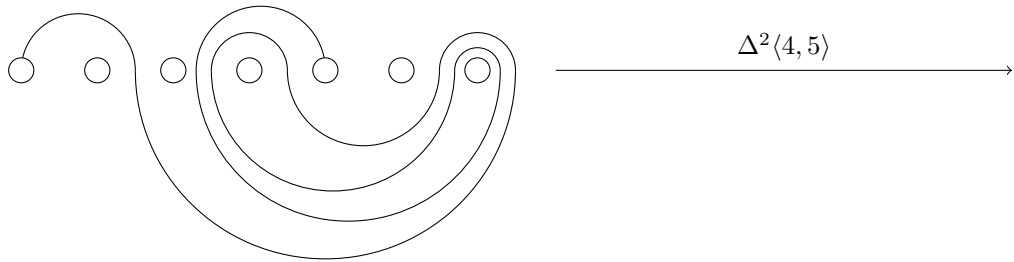
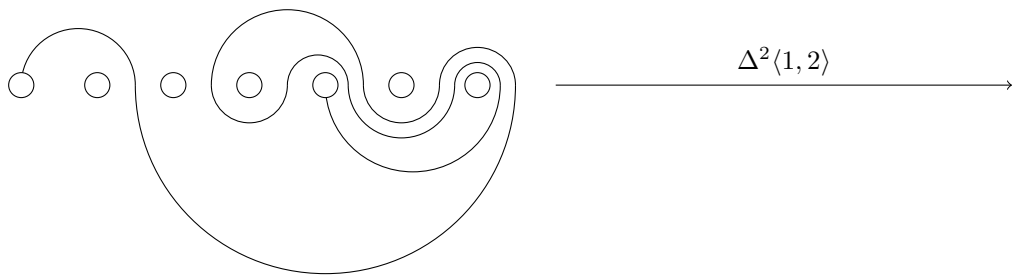
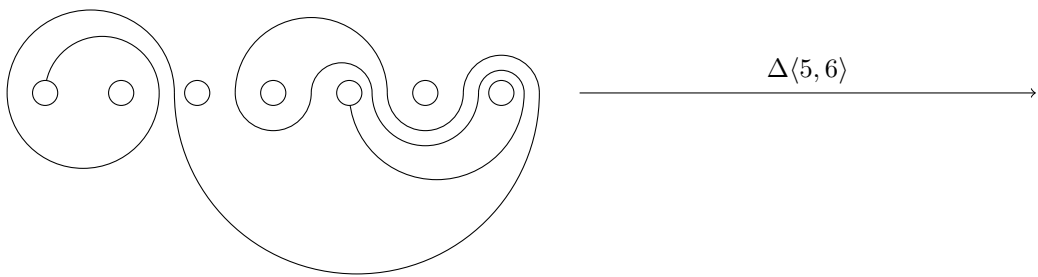
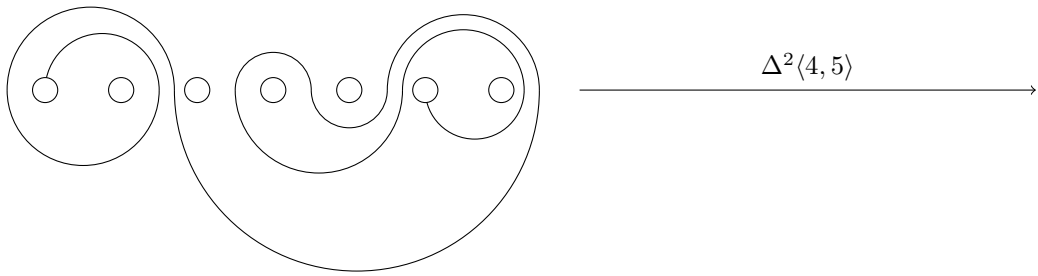
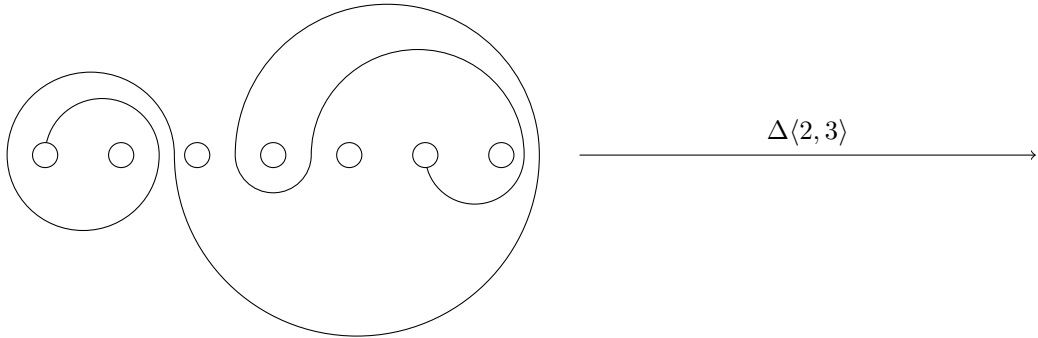
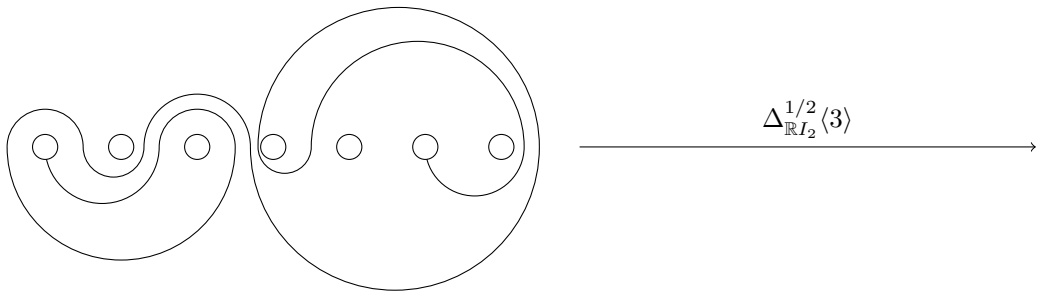
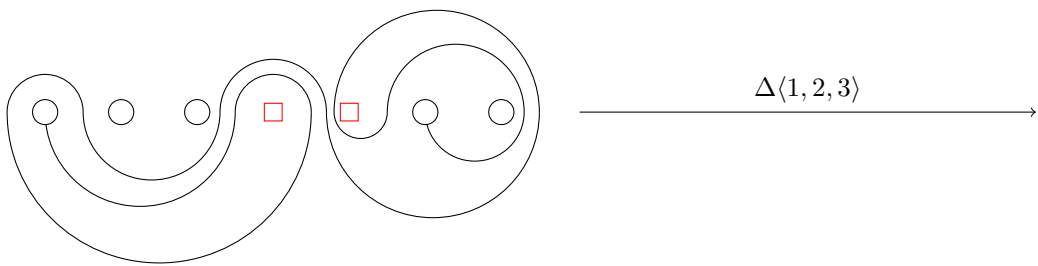
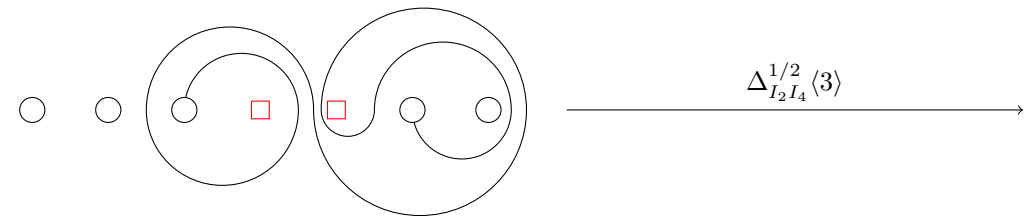
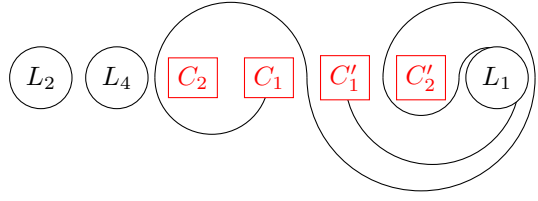
\begin{figure}[H]
	\begin{center}
		%
	\end{center}
\end{figure}

Relation:

\begin{equation*}\begin{split}
\Gamma_{L_1}\Gamma_{C_2'}^{-1}\Gamma_{L_1}^{-1}\Gamma_{C_1}\Gamma_{C_2}\Gamma_{C_1}\Gamma_{C_2}^{-1}\Gamma_{C_1}^{-1}\Gamma_{L_1}\Gamma_{C_2'}\Gamma_{L_1}^{-1}  = \Gamma_{C_1'}
\end{split}
\end{equation*}
\subsubsection{Raw relations and simplifications}

\begin{equation}\label{arrangement2_auto_calc_vert_1_rel_1}
\begin{split}
\Gamma_{C_1}  = \Gamma_{C_1'} ,
\end{split}
\end{equation}
\begin{equation}\label{arrangement2_auto_calc_vert_2_rel_1}
\begin{split}
\left[\Gamma_{C_2}^{-1}\Gamma_{C_1}\Gamma_{C_2} , \Gamma_{L_4}\Gamma_{L_2}  \right]=e,
\end{split}
\end{equation}
\begin{equation}\label{arrangement2_auto_calc_vert_2_rel_2}
\begin{split}
\left[\Gamma_{L_4} , \Gamma_{L_2}\Gamma_{C_2}^{-1}\Gamma_{C_1}\Gamma_{C_2}  \right]=e,
\end{split}
\end{equation}
\begin{equation}\label{arrangement2_auto_calc_vert_4_rel_1}
\begin{split}
\Gamma_{C_1'}\Gamma_{C_1}\Gamma_{C_2}\Gamma_{L_4}\Gamma_{L_2}\Gamma_{L_4}^{-1}\Gamma_{C_2}^{-1}\Gamma_{C_1}^{-1}\Gamma_{C_2}\Gamma_{C_1}\Gamma_{C_2}\Gamma_{L_4}\Gamma_{L_2}^{-1}\Gamma_{L_4}^{-1}\Gamma_{C_2}^{-1}\Gamma_{C_1}^{-1}\Gamma_{C_1'}^{-1}  = \Gamma_{C_2'} ,
\end{split}
\end{equation}
\begin{equation}\label{arrangement2_auto_calc_vert_5_rel_1}
\begin{split}
\left[\Gamma_{C_1}\Gamma_{C_2}\Gamma_{L_4}\Gamma_{C_2}^{-1}\Gamma_{C_1}^{-1} , \Gamma_{C_2}  \right]=e,
\end{split}
\end{equation}
\begin{equation}\label{arrangement2_auto_calc_vert_6_rel_1}
\begin{split}
\left\{\Gamma_{L_4}^{-1}\Gamma_{C_2}^{-1}\Gamma_{C_1}^{-1}\Gamma_{C_1'}^{-1}\Gamma_{C_2'}\Gamma_{C_1'}\Gamma_{C_1}\Gamma_{C_2}\Gamma_{L_4} , \Gamma_{L_2}  \right\}=e,
\end{split}
\end{equation}
\begin{equation}\label{arrangement2_auto_calc_vert_7_rel_1}
\begin{split}
\left[\Gamma_{C_2'}\Gamma_{C_1'}\Gamma_{C_1}\Gamma_{C_2}\Gamma_{L_4}\Gamma_{L_2}\Gamma_{L_4}^{-1}\Gamma_{C_2}^{-1}\Gamma_{C_1}^{-1}\Gamma_{C_1'}^{-1}\Gamma_{C_2'}^{-1} , \Gamma_{C_1'}  \right]=e,
\end{split}
\end{equation}
\begin{equation}\label{arrangement2_auto_calc_vert_8_rel_1}
\begin{split}
\left\{\Gamma_{C_1} , \Gamma_{C_2}  \right\}=e,
\end{split}
\end{equation}
\begin{equation}\label{arrangement2_auto_calc_vert_9_rel_1}
\begin{split}
\left\{\Gamma_{C_2'} , \Gamma_{C_1'}  \right\}=e,
\end{split}
\end{equation}
\begin{equation}\label{arrangement2_auto_calc_vert_10_rel_1}
\begin{split}
\left[\Gamma_{C_2'}\Gamma_{C_1'}\Gamma_{C_1}\Gamma_{C_2}\Gamma_{L_4}\Gamma_{L_2}\Gamma_{L_4}^{-1}\Gamma_{C_2}^{-1}\Gamma_{C_1}^{-1}\Gamma_{C_1'}^{-1}\Gamma_{C_2'}^{-1} , \Gamma_{L_1}  \right]=e,
\end{split}
\end{equation}
\begin{equation}\label{arrangement2_auto_calc_vert_11_rel_1}
\begin{split}
\left[\Gamma_{C_1}\Gamma_{C_2}\Gamma_{L_4}\Gamma_{C_2}^{-1}\Gamma_{C_1}^{-1} , \Gamma_{C_2'}  \right]=e,
\end{split}
\end{equation}
\begin{equation}\label{arrangement2_auto_calc_vert_12_rel_1}
\begin{split}
\left[\Gamma_{C_2'}^{-1}\Gamma_{L_1}\Gamma_{C_2'} , \Gamma_{C_1'}\Gamma_{C_1}\Gamma_{C_2}\Gamma_{L_4}\Gamma_{C_2}^{-1}\Gamma_{C_1}^{-1}  \right]=e,
\end{split}
\end{equation}
\begin{equation}\label{arrangement2_auto_calc_vert_12_rel_2}
\begin{split}
\left[\Gamma_{C_1'} , \Gamma_{C_1}\Gamma_{C_2}\Gamma_{L_4}\Gamma_{C_2}^{-1}\Gamma_{C_1}^{-1}\Gamma_{C_2'}^{-1}\Gamma_{L_1}\Gamma_{C_2'}  \right]=e,
\end{split}
\end{equation}
\begin{equation}\label{arrangement2_auto_calc_vert_13_rel_1}
\begin{split}
\left\{\Gamma_{C_2'} , \Gamma_{L_1}  \right\}=e,
\end{split}
\end{equation}
\begin{equation}\label{arrangement2_auto_calc_vert_14_rel_1}
\begin{split}
\Gamma_{L_1}^{-1}\Gamma_{C_1}\Gamma_{C_2}\Gamma_{C_1}^{-1}\Gamma_{L_1}  = \Gamma_{C_2'} ,
\end{split}
\end{equation}
\begin{equation}\label{arrangement2_auto_calc_vert_15_rel_1}
\begin{split}
\left[\Gamma_{C_2'}^{-1}\Gamma_{L_1}^{-1}\Gamma_{C_1}\Gamma_{C_2}\Gamma_{C_1}\Gamma_{C_2}^{-1}\Gamma_{C_1}^{-1}\Gamma_{L_1}\Gamma_{C_2'} , \Gamma_{L_1}  \right]=e,
\end{split}
\end{equation}
\begin{equation}\label{arrangement2_auto_calc_vert_16_rel_1}
\begin{split}
\Gamma_{L_1}\Gamma_{C_2'}^{-1}\Gamma_{L_1}^{-1}\Gamma_{C_1}\Gamma_{C_2}\Gamma_{C_1}\Gamma_{C_2}^{-1}\Gamma_{C_1}^{-1}\Gamma_{L_1}\Gamma_{C_2'}\Gamma_{L_1}^{-1}  = \Gamma_{C_1'} ,
\end{split}
\end{equation}
\begin{equation}\label{arrangement2_auto_calc_projective_rel}
\begin{split}
\Gamma_{L_1}\Gamma_{C_2'}\Gamma_{C_1'}\Gamma_{C_1}\Gamma_{C_2}\Gamma_{L_4}\Gamma_{L_2} =e.
\end{split}
\end{equation}

Here are the computations for group $ G_6 := \pcpt{\mathcal{B}_6} / \langle \Gamma_{X}^2 \; | \; X\subseteq \mathcal{B}_6 \rangle $:
\begin{lemma}\label{tokunaga-arrangement1-b2-group-calc}
\begin{enumerate}
	\item Group $ G_6 := \pcpt{\mathcal{B}_6} / \langle \Gamma_{X}^2 \; | \; X\subseteq \mathcal{B}_6 \rangle $ is generated by $ \Gamma_{C_1},\Gamma_{C_2}, \Gamma_{L_1} $, and $ \Gamma_{L_2} $, subject to the following relations:
	\begin{equation}\label{arrangement2_second_final_1}
		\Gamma_{C_1}^2 = \Gamma_{C_2}^2 = \Gamma_{L_1}^2 = \Gamma_{L_2}^2 = e,
	\end{equation}
	\begin{equation}\label{arrangement2_second_final_2}
		[\Gamma_{L_1}, \Gamma_{L_2} ] = [\Gamma_{C_1}, \Gamma_{L_1}] = [\Gamma_{C_1}, \Gamma_{L_2}] = [\Gamma_{L_1}\Gamma_{L_2}, \Gamma_{C_2}] = e,
	\end{equation}
	\begin{equation}\label{arrangement2_second_final_3}
		\{ \Gamma_{C_1}, \Gamma_{C_2} \} = \{ \Gamma_{C_2}, \Gamma_{L_1} \} = \{ \Gamma_{C_2}, \Gamma_{L_2} \} = e,
	\end{equation}
	\begin{equation}\label{arrangement2_second_final_4}
		[\Gamma_{C_2}, \Gamma_{L_1}\Gamma_{C_1}\Gamma_{C_2}\Gamma_{C_1}\Gamma_{L_1}] = e.
	\end{equation}
	
	\item Group $ G_6 $ is isomorphic to the direct product $ G_6'\times \mbb{Z}_2 $, where $ G_6' $ is generated by $ \Gamma_{C_1},\Gamma_{C_2} $, and $ \Gamma_{L_1} $, subject to the following relations:
	\begin{equation}\label{arrangement2_without_L2_1}
		\Gamma_{C_1}^2 = \Gamma_{C_2}^2 = \Gamma_{L_1}^2 = e,
	\end{equation}
	\begin{equation}\label{arrangement2_without_L2_2}
		[\Gamma_{C_1}, \Gamma_{L_1}] = e,
	\end{equation}
	\begin{equation}\label{arrangement2_without_L2_3}
		\{ \Gamma_{C_1}, \Gamma_{C_2} \} = \{ \Gamma_{C_2}, \Gamma_{L_1} \} = e,
	\end{equation}
	\begin{equation}\label{arrangement2_without_L2_4}
		[\Gamma_{C_2}, \Gamma_{L_1}\Gamma_{C_1}\Gamma_{C_2}\Gamma_{C_1}\Gamma_{L_1}] = e.
	\end{equation}		
	
	\item The element $ \Gamma_{L_4} \in \pcpt{\mathcal{B}_6} $ is  not contained in the center of $ \pcpt{\mathcal{B}_6} $.
\end{enumerate}
\end{lemma}
\begin{proof}
Substituting $ \Gamma_{C_1'}=\Gamma_{C_1} $ and $ \Gamma_{C_2'} = \Gamma_{L_1}\Gamma_{C_1}\Gamma_{C_2}\Gamma_{C_1}\Gamma_{L_1} $ we get that $ G_6 $ is generated by $ \Gamma_{C_1}, \Gamma_{C_2}, \Gamma_{L_1},\Gamma_{L_2} $ and $ \Gamma_{L_4} $ subject to:	
\begin{equation}\label{arrangement2_simpl_1}
	\begin{split}
		(\Gamma_{C_2}\Gamma_{C_1}\Gamma_{C_2})  \Gamma_{L_4}\Gamma_{L_2} = \Gamma_{L_4}\Gamma_{L_2}(\Gamma_{C_2}\Gamma_{C_1}\Gamma_{C_2}) =\Gamma_{L_2}(\Gamma_{C_2}\Gamma_{C_1}\Gamma_{C_2})\Gamma_{L_4},
	\end{split}
\end{equation}
\begin{equation}\label{arrangement2_simpl_3}
	\begin{split} \Gamma_{C_2}\Gamma_{L_4}\Gamma_{L_2}\Gamma_{L_4}\Gamma_{C_2}\Gamma_{C_1}\Gamma_{C_2}\Gamma_{C_1}\Gamma_{C_2}\Gamma_{L_4}\Gamma_{L_2}\Gamma_{L_4}\Gamma_{C_2}  = \Gamma_{L_1}\Gamma_{C_1}\Gamma_{C_2}\Gamma_{C_1}\Gamma_{L_1} ,
	\end{split}
\end{equation}
\begin{equation}\label{arrangement2_simpl_4}
	\begin{split}
		\left[\Gamma_{C_1}\Gamma_{C_2}\Gamma_{L_4}\Gamma_{C_2}\Gamma_{C_1} , \Gamma_{C_2}  \right]=e,
	\end{split}
\end{equation}
\begin{equation}\label{arrangement2_simpl_5}
	\begin{split} \left\{\Gamma_{L_4}\Gamma_{C_2}\Gamma_{L_1}\Gamma_{C_1}\Gamma_{C_2}\Gamma_{C_1}\Gamma_{L_1}\Gamma_{C_2}\Gamma_{L_4} , \Gamma_{L_2}  \right\}=e,
	\end{split}
\end{equation}
\begin{equation}\label{arrangement2_simpl_6}
	\begin{split} \left[\Gamma_{L_1}\Gamma_{C_1}\Gamma_{C_2}\Gamma_{C_1}\Gamma_{L_1}\Gamma_{C_2}\Gamma_{L_4}\Gamma_{L_2}\Gamma_{L_4}\Gamma_{C_2}\Gamma_{L_1}\Gamma_{C_1}\Gamma_{C_2}\Gamma_{C_1}\Gamma_{L_1} , \Gamma_{C_1}  \right]=e,
	\end{split}
\end{equation}
\begin{equation}\label{arrangement2_simpl_7}
	\begin{split}
		\left\{\Gamma_{C_1} , \Gamma_{C_2}  \right\}=e,
	\end{split}
\end{equation}
\begin{equation}\label{arrangement2_simpl_8}
	\begin{split}
		\left\{\Gamma_{L_1}\Gamma_{C_1}\Gamma_{C_2}\Gamma_{C_1}\Gamma_{L_1} , \Gamma_{C_1}  \right\}=e,
	\end{split}
\end{equation}
\begin{equation}\label{arrangement2_simpl_9}
	\begin{split} \left[\Gamma_{L_1}\Gamma_{C_1}\Gamma_{C_2}\Gamma_{C_1}\Gamma_{L_1}\Gamma_{C_2}\Gamma_{L_4}\Gamma_{L_2}\Gamma_{L_4}\Gamma_{C_2}\Gamma_{L_1}\Gamma_{C_1}\Gamma_{C_2}\Gamma_{C_1}\Gamma_{L_1} , \Gamma_{L_1}  \right]=e,
	\end{split}
\end{equation}
\begin{equation}\label{arrangement2_simpl_10}
	\begin{split}
		\left[\Gamma_{C_1}\Gamma_{C_2}\Gamma_{L_4}\Gamma_{C_2}\Gamma_{C_1} , \Gamma_{L_1}\Gamma_{C_1}\Gamma_{C_2}\Gamma_{C_1}\Gamma_{L_1}  \right]=e,
	\end{split}
\end{equation}
\begin{equation}\label{arrangement2_simpl_11}
	\begin{split} \left[\Gamma_{L_1}\Gamma_{C_1}\Gamma_{C_2}\Gamma_{C_1}\Gamma_{L_1}\Gamma_{C_1}\Gamma_{C_2}\Gamma_{C_1}\Gamma_{L_1} , \Gamma_{C_2}\Gamma_{L_4}\Gamma_{C_2}\Gamma_{C_1}  \right]=e,
	\end{split}
\end{equation}
\begin{equation}\label{arrangement2_simpl_12}
	\begin{split}
		\left[\Gamma_{C_1} , \Gamma_{C_1}\Gamma_{C_2}\Gamma_{L_4}\Gamma_{C_2}\Gamma_{C_1}\Gamma_{L_1}\Gamma_{C_1}\Gamma_{C_2}\Gamma_{C_1}\Gamma_{L_1}\Gamma_{C_1}\Gamma_{C_2}\Gamma_{C_1}\Gamma_{L_1}  \right]=e,
	\end{split}
\end{equation}
\begin{equation}\label{arrangement2_simpl_13}
	\begin{split}
		\left\{\Gamma_{L_1}\Gamma_{C_1}\Gamma_{C_2}\Gamma_{C_1}\Gamma_{L_1} , \Gamma_{L_1}  \right\}=e,
	\end{split}
\end{equation}
\begin{equation}\label{arrangement2_simpl_14}
	\begin{split}
		\left[\Gamma_{L_1}\Gamma_{C_1}\Gamma_{L_1} , \Gamma_{L_1}  \right]=e,
	\end{split}
\end{equation}
\begin{equation}\label{arrangement2_simpl_15}
	\begin{split}
		\Gamma_{C_1}\Gamma_{C_2}\Gamma_{C_1}\Gamma_{L_1}\Gamma_{C_2}\Gamma_{L_4}\Gamma_{L_2} =e.
	\end{split}
\end{equation}

From \eqref{arrangement2_simpl_14} we get that $ \Gamma_{C_1} $ and $ \Gamma_{L_1} $ commute, substituting this into \eqref{arrangement2_simpl_13} we get that $ \{ \Gamma_{C_2}, \Gamma_{L_1}\}=e $.
Using those relations in \eqref{arrangement2_simpl_10} we get that
\begin{align*}
	e & = [ \Gamma_{L_4}, \Gamma_{C_2}\Gamma_{C_1}\Gamma_{L_1}\Gamma_{C_1}\Gamma_{C_2}\Gamma_{C_1}\Gamma_{L_1}\Gamma_{C_1}
	\Gamma_{C_2} ] = \\
	& = [ \Gamma_{L_4}, \Gamma_{C_2}\Gamma_{L_1}\Gamma_{C_2}\Gamma_{L_1}
	\Gamma_{C_2} ] = [ \Gamma_{L_4}, \Gamma_{L_1}\Gamma_{C_2}\Gamma_{L_1} ].
\end{align*}

Now, from \eqref{arrangement2_simpl_4}, using \eqref{arrangement2_simpl_7}, we get
\[ e = [\Gamma_{L_4}, \Gamma_{C_2}\Gamma_{C_1}\Gamma_{C_2}\Gamma_{C_1}\Gamma_{C_2}] = [\Gamma_{L_4}, \Gamma_{C_1}\Gamma_{C_2}\Gamma_{C_1}].  \]

Using the above, we can simplify \eqref{arrangement2_simpl_3} as follows:
\begin{align*} \Gamma_{C_2}\Gamma_{L_4}\Gamma_{L_2}\Gamma_{L_4}\Gamma_{C_2}\Gamma_{C_1}\Gamma_{C_2}\Gamma_{C_1}\Gamma_{C_2}\Gamma_{L_4}\Gamma_{L_2}\Gamma_{L_4}\Gamma_{C_2} & = \Gamma_{L_1}\Gamma_{C_1}\Gamma_{C_2}\Gamma_{C_1}\Gamma_{L_1} \\ \Gamma_{C_2}\Gamma_{L_4}\Gamma_{L_2}\Gamma_{C_1}\Gamma_{C_2}\Gamma_{C_1}\Gamma_{L_2}\Gamma_{L_4}\Gamma_{C_2} & = \Gamma_{L_1}\Gamma_{C_1}\Gamma_{C_2}\Gamma_{C_1}\Gamma_{L_1} \\ \Gamma_{C_1}\Gamma_{C_2}\Gamma_{L_4}\Gamma_{L_2}\Gamma_{C_1}\Gamma_{C_2}\Gamma_{C_1}\Gamma_{L_2}\Gamma_{L_4}\Gamma_{C_2}\Gamma_{C_1} & = \Gamma_{C_2}\Gamma_{L_1}\Gamma_{C_2}\Gamma_{L_1}\Gamma_{C_2} \\ \Gamma_{L_4}\Gamma_{L_2}\Gamma_{C_2}\Gamma_{C_1}\Gamma_{C_2}\Gamma_{C_1}\Gamma_{C_2}\Gamma_{C_1}\Gamma_{C_2}\Gamma_{C_1}\Gamma_{C_2}\Gamma_{L_2}\Gamma_{L_4} & = \Gamma_{L_1}\Gamma_{C_2}\Gamma_{L_1}\\
	\Gamma_{L_2}\Gamma_{C_2}\Gamma_{L_2} = \Gamma_{L_4}\Gamma_{L_1}\Gamma_{C_2}\Gamma_{L_1}\Gamma_{L_4} & = \Gamma_{L_1}\Gamma_{C_2}\Gamma_{L_1} \\
	[\Gamma_{C_2}, \Gamma_{L_1}\Gamma_{L_2}] &= e,
\end{align*}
where we first  used $ \{\Gamma_{C_1}, \Gamma_{C_2} \}=e $ and $ [\Gamma_{L_4}, \Gamma_{C_1}\Gamma_{C_2}\Gamma_{C_1}]=e $; then $ [\Gamma_{C_1}, \Gamma_{L_1}]=e $ and $ \{\Gamma_{C_2}, \Gamma_{L_1}\}=e $; next we applied \eqref{arrangement2_simpl_1}; and finally used $ \{\Gamma_{C_1},\Gamma_{C_2}\}=e $ followed by $ [\Gamma_{L_4}, \Gamma_{L_1}\Gamma_{C_2}\Gamma_{L_1}]=e $.

We now turn to \eqref{arrangement2_simpl_9}. By using derived commutation relations and rearranging we get that it is equivalent to
\[ \left[\Gamma_{C_1}\Gamma_{C_2}\Gamma_{L_4}\Gamma_{L_2}\Gamma_{L_4}\Gamma_{C_2}\Gamma_{C_1} , \Gamma_{L_1}\Gamma_{C_2}\Gamma_{L_1}\Gamma_{C_2}\Gamma_{L_1}  \right]=e. \]
By applying $ \{ \Gamma_{C_2}, \Gamma_{L_1}\} = e $ and rearranging again,
\[ \left[\Gamma_{C_2}\Gamma_{C_1}\Gamma_{C_2}\Gamma_{L_4}\Gamma_{L_2}\Gamma_{L_4}\Gamma_{C_2}\Gamma_{C_1}\Gamma_{C_2} , \Gamma_{L_1}  \right]=e \]
which, by \eqref{arrangement2_simpl_1}, is equivalent to $ [\Gamma_{L_1}, \Gamma_{L_2}]=e $.

Using $ [\Gamma_{C_1}, \Gamma_{L_1}]=e $ and $ \{\Gamma_{C_2}, \Gamma_{L_1}\}=e $ we get that \eqref{arrangement2_simpl_5} is equivalent to
\[ \{ \Gamma_{L_1}\Gamma_{C_2}\Gamma_{L_1}, \Gamma_{C_2}\Gamma_{C_1}\Gamma_{C_2}\Gamma_{L_4}\Gamma_{L_2}\Gamma_{L_4}\Gamma_{C_2}\Gamma_{C_1}\Gamma_{C_2} \} = e \]
which is, by \eqref{arrangement2_simpl_1} and $ [\Gamma_{L_1}, \Gamma_{L_2}]=e $, equivalent to $ \{ \Gamma_{C_2}, \Gamma_{L_2}\}=e $.

Next we deal with \eqref{arrangement2_simpl_6}. Applying derived commutations, $ \{ \Gamma_{C_1}, \Gamma_{C_2} \}=e $ and rearranging, we get
\[ [\Gamma_{L_1}\Gamma_{C_2}\Gamma_{L_4}\Gamma_{L_2}\Gamma_{L_4}\Gamma_{C_2}\Gamma_{L_1}, \Gamma_{C_2}\Gamma_{C_1}\Gamma_{C_2}] = e. \]
Using $ \{\Gamma_{C_2}, \Gamma_{L_1} \} =e $, $ [\Gamma_{C_1}, \Gamma_{L_1}]=e $ and $ {[\Gamma_{L_4}, \Gamma_{L_1}\Gamma_{C_2}\Gamma_{L_1}] = e} $ yields
\[ [\Gamma_{C_2}\Gamma_{L_4}\Gamma_{L_1}\Gamma_{C_2}\Gamma_{L_1}\Gamma_{L_2}\Gamma_{L_1}\Gamma_{C_2}\Gamma_{L_1}\Gamma_{L_4}\Gamma_{C_2}, \Gamma_{C_1}] = e. \]
Applying $ [\Gamma_{C_2}, \Gamma_{L_1}\Gamma_{L_2}]=e $ and rearranging once more, gives us
\[ [\Gamma_{C_2}\Gamma_{L_2}\Gamma_{C_2}, \Gamma_{L_2}\Gamma_{L_4}\Gamma_{C_2}\Gamma_{C_1}\Gamma_{C_2}\Gamma_{L_4}\Gamma_{L_2} ] \]
which then, by using \eqref{arrangement2_simpl_1}, becomes $ [\Gamma_{L_2}, \Gamma_{C_1}] = e $.

From \eqref{arrangement2_simpl_11} we get, by conjugating with $ \Gamma_{C_2}\Gamma_{C_1} $ the following:
\begin{align*}
	e & = [\Gamma_{C_2}\Gamma_{C_1}\Gamma_{L_1}\Gamma_{C_1}\Gamma_{C_2}\Gamma_{C_1}\Gamma_{L_1}\Gamma_{C_1}\Gamma_{C_2}\Gamma_{C_1}\Gamma_{L_1}\Gamma_{C_1}\Gamma_{C_2}, \Gamma_{C_2}\Gamma_{C_1}\Gamma_{C_2}\Gamma_{L_4}] = \\
	& = [\Gamma_{C_2}\Gamma_{L_1}\Gamma_{C_2}\Gamma_{L_1}\Gamma_{C_2}\Gamma_{L_1}\Gamma_{C_2}, \Gamma_{C_2}\Gamma_{C_1}\Gamma_{C_2}\Gamma_{L_4}] = \\
	& = [\Gamma_{L_1}, \Gamma_{C_2}\Gamma_{C_1}\Gamma_{C_2}\Gamma_{L_4}],
\end{align*}
and from \eqref{arrangement2_simpl_12} we get:
\begin{align*}
	e & =
	\left[\Gamma_{C_1} , \Gamma_{C_1}\Gamma_{C_2}\Gamma_{L_4}\Gamma_{C_2}\Gamma_{C_1}\Gamma_{L_1}\Gamma_{C_1}\Gamma_{C_2}\Gamma_{C_1}\Gamma_{L_1}\Gamma_{C_1}\Gamma_{C_2}\Gamma_{C_1}\Gamma_{L_1}  \right] = \\
	& =
	\left[\Gamma_{C_1} , \Gamma_{C_2}\Gamma_{L_4}\Gamma_{C_2}\Gamma_{L_1}\Gamma_{C_2}\Gamma_{L_1}\Gamma_{C_2}\Gamma_{L_1}  \right] = \\
	& =
	\left[\Gamma_{C_1} , \Gamma_{C_2}\Gamma_{L_4}\Gamma_{L_1}\Gamma_{C_2}  \right].
\end{align*}
Together those relations can be written as
\[ \Gamma_{L_1}(\Gamma_{C_2}\Gamma_{C_1}\Gamma_{C_2})\Gamma_{L_4} = (\Gamma_{C_2}\Gamma_{C_1}\Gamma_{C_2})\Gamma_{L_4}\Gamma_{L_1} = \Gamma_{L_4}\Gamma_{L_1}(\Gamma_{C_2}\Gamma_{C_1}\Gamma_{C_2}). \]

Summarizing those simplifications we get that $ G_6 $ is generated by $ \Gamma_{C_1},\Gamma_{C_2}, \Gamma_{L_1}, \Gamma_{L_2} $, and $ \Gamma_{L_4} $, subject to the following relations:
\begin{equation}\label{arrangement2_final_1}
	\Gamma_{C_1}^2 = \Gamma_{C_2}^2 = \Gamma_{L_1}^2 = \Gamma_{L_2}^2 = \Gamma_{L_4}^2 =e,
\end{equation}
\begin{equation}\label{arrangement2_final_2}
	[\Gamma_{L_1}, \Gamma_{L_2} ] = [\Gamma_{C_1}, \Gamma_{L_1}] = [\Gamma_{C_1}, \Gamma_{L_2}] = e,
\end{equation}
\begin{equation}\label{arrangement2_final_3}
	\{ \Gamma_{C_1}, \Gamma_{C_2} \} = \{ \Gamma_{C_2}, \Gamma_{L_1} \} = \{ \Gamma_{C_2}, \Gamma_{L_2} \} = e,
\end{equation}
\begin{equation}\label{arrangement2_final_4}
	\Gamma_{L_1}(\Gamma_{C_2}\Gamma_{C_1}\Gamma_{C_2})\Gamma_{L_4} = (\Gamma_{C_2}\Gamma_{C_1}\Gamma_{C_2})\Gamma_{L_4}\Gamma_{L_1} = \Gamma_{L_4}\Gamma_{L_1}(\Gamma_{C_2}\Gamma_{C_1}\Gamma_{C_2}),
\end{equation}
\begin{equation}\label{arrangement2_final_5}
	\Gamma_{L_2}(\Gamma_{C_2}\Gamma_{C_1}\Gamma_{C_2})\Gamma_{L_4} = (\Gamma_{C_2}\Gamma_{C_1}\Gamma_{C_2})\Gamma_{L_4}\Gamma_{L_2} = \Gamma_{L_4}\Gamma_{L_2}(\Gamma_{C_2}\Gamma_{C_1}\Gamma_{C_2}),
\end{equation}
\begin{equation}\label{arrangement2_final_6}
	[\Gamma_{C_2}, \Gamma_{L_1}\Gamma_{L_2}] = [\Gamma_{L_4}, \Gamma_{C_1}\Gamma_{C_2}\Gamma_{C_1}] = [\Gamma_{L_4}, \Gamma_{L_1}\Gamma_{C_2}\Gamma_{L_1}] = e,
\end{equation}
\begin{equation}\label{arrangement2_final_7}
	\Gamma_{C_1}\Gamma_{C_2}\Gamma_{C_1}\Gamma_{L_1}\Gamma_{C_2}\Gamma_{L_4}\Gamma_{L_2} =e.
\end{equation}

Note that because $ [ \Gamma_{C_2}, \Gamma_{L_1}\Gamma_{L_2} ] = e $, we have
\begin{equation}\label{arrangement2_second_L1_L2}
	\Gamma_{L_2}\Gamma_{C_2}\Gamma_{L_2} = \Gamma_{L_1}\Gamma_{C_2}\Gamma_{L_1},
\end{equation}
and
\begin{equation}
	\Gamma_{L_2}\Gamma_{C_2}\Gamma_{L_1} = \Gamma_{L_1}\Gamma_{C_2}\Gamma_{L_2}.
\end{equation}

We now use \eqref{arrangement2_final_7} to eliminate $ \Gamma_{L_4} $. To that end, write
$$ \Gamma_{L_2}\Gamma_{C_1}\Gamma_{C_2}\Gamma_{L_1}\Gamma_{C_1}\Gamma_{C_2} = \Gamma_{L_4} = \Gamma_{L_4}^{-1} = \Gamma_{C_2}\Gamma_{C_1}\Gamma_{L_1}\Gamma_{C_2}\Gamma_{C_1}\Gamma_{L_2}, $$
and substitute it in relations \eqref{arrangement2_final_1}-\eqref{arrangement2_final_6}.

First, using \eqref{arrangement2_second_L1_L2}, $ \Gamma_{L_4}^2 = e $ is equivalent to
\begin{equation}\label{arrangement2_second1}
	\begin{split}
		e = \Gamma_{L_2}\Gamma_{C_1}\Gamma_{C_2}\Gamma_{L_1}\Gamma_{C_1}\Gamma_{C_2}\Gamma_{L_2}\Gamma_{C_1}\Gamma_{C_2}\Gamma_{L_1}\Gamma_{C_1}\Gamma_{C_2} &= \Gamma_{C_1}\Gamma_{L_2}\Gamma_{C_2}\Gamma_{C_1}\Gamma_{L_1}\Gamma_{C_2}\Gamma_{L_2}\Gamma_{C_1}\Gamma_{C_2}\Gamma_{L_1}\Gamma_{C_1}\Gamma_{C_2} = \\
		= \Gamma_{C_1}\Gamma_{L_2}\Gamma_{C_2}\Gamma_{C_1}\Gamma_{L_2}\Gamma_{C_2}\Gamma_{L_1}\Gamma_{C_1}\Gamma_{C_2}\Gamma_{L_1}\Gamma_{C_1}\Gamma_{C_2} &= \Gamma_{C_1}\Gamma_{L_1}\Gamma_{C_2}\Gamma_{L_1}\Gamma_{C_1}\Gamma_{C_2}\Gamma_{C_1}\Gamma_{L_1}\Gamma_{C_2}\Gamma_{L_1}\Gamma_{C_1}\Gamma_{C_2} = \\
		&= [ \Gamma_{C_2}, \Gamma_{L_1}\Gamma_{C_1}\Gamma_{C_2}\Gamma_{C_1}\Gamma_{L_1}].
	\end{split}
\end{equation}

Next we deal with \eqref{arrangement2_final_4}. After the substitution, the first part reads
\begin{equation*}
	\Gamma_{L_1}\Gamma_{C_2}\Gamma_{L_1}\Gamma_{C_2}\Gamma_{C_1}\Gamma_{L_2} = \Gamma_{C_2}\Gamma_{L_1}\Gamma_{C_2}\Gamma_{C_1}\Gamma_{L_2}\Gamma_{L_1}
\end{equation*}
and because $ \Gamma_{L_1} $ commutes with both $ \Gamma_{C_1} $ and $ \Gamma_{L_2} $, it is equivalent to $ \{ \Gamma_{C_2}, \Gamma_{L_1} \} = e $, which already appears in  \eqref{arrangement2_final_3}.
The second part of \eqref{arrangement2_final_4} reads
\begin{equation*}
	\Gamma_{C_2}\Gamma_{L_1}\Gamma_{C_2}\Gamma_{C_1}\Gamma_{L_2}\Gamma_{L_1} = \Gamma_{L_2}\Gamma_{C_1}\Gamma_{C_2}\Gamma_{L_1}\Gamma_{C_1}\Gamma_{C_2}\Gamma_{L_1}\Gamma_{C_2}\Gamma_{C_1}\Gamma_{C_2}.
\end{equation*}
Multiplying by $ \Gamma_{C_1}\Gamma_{C_2} $ from the right and using the fact that both $ \Gamma_{C_1} $ and $\Gamma_{C_2}$ commute with $ \Gamma_{L_2}\Gamma_{L_1} $, $ [\Gamma_{L_1},\Gamma_{L_2}]=e $ and $ \{ \Gamma_{C_1},\Gamma_{C_2} \}=e $, we get
\begin{equation*}
	\Gamma_{C_2}\Gamma_{L_2} = \Gamma_{L_2}\Gamma_{C_1}\Gamma_{C_2}\Gamma_{L_1}\Gamma_{C_1}\Gamma_{C_2}\Gamma_{L_1}\Gamma_{C_1}\Gamma_{C_2}\Gamma_{C_1}.
\end{equation*}
This equation is equivalent to
\begin{equation}\label{arrangement2_second2}
	\Gamma_{C_2}\Gamma_{L_2}\Gamma_{C_2}\Gamma_{L_2}\Gamma_{C_2} = \Gamma_{C_2}\Gamma_{C_1}\Gamma_{C_2}\Gamma_{L_1}\Gamma_{C_1}\Gamma_{C_2}\Gamma_{L_1}\Gamma_{C_1}\Gamma_{C_2}\Gamma_{C_1}\Gamma_{C_2}.
\end{equation}
Using \eqref{arrangement2_second_L1_L2} and  \eqref{arrangement2_final_3} we see that
\begin{equation*}
	\Gamma_{L_1}\Gamma_{C_2}\Gamma_{L_1} = \Gamma_{L_2}\Gamma_{C_2}\Gamma_{L_2} = \Gamma_{C_2}\Gamma_{L_2}\Gamma_{C_2}\Gamma_{L_2}\Gamma_{C_2}
\end{equation*}
so, the left side of \eqref{arrangement2_second2} is equal to $\Gamma_{L_1}\Gamma_{C_2}\Gamma_{L_1}$.
Using \eqref{arrangement2_final_3} and $ [\Gamma_{C_1}, \Gamma_{L_1} ] =e $, we get that the right side of \eqref{arrangement2_second2} is equal to
\begin{align*} \Gamma_{C_2}\Gamma_{C_1}\Gamma_{C_2}\Gamma_{L_1}\Gamma_{C_1}\Gamma_{C_2}\Gamma_{L_1}\Gamma_{C_1}\Gamma_{C_2}\Gamma_{C_1}\Gamma_{C_2} &	= \Gamma_{C_2}\Gamma_{C_1}\Gamma_{C_2}\Gamma_{C_1}\Gamma_{L_1}\Gamma_{C_2}\Gamma_{L_1}\Gamma_{C_1}\Gamma_{C_2}\Gamma_{C_1}\Gamma_{C_2} = \\
	= \Gamma_{C_1}\Gamma_{C_2}\Gamma_{C_1}\Gamma_{C_2}\Gamma_{L_1}\Gamma_{C_2}\Gamma_{L_1}\Gamma_{C_2}\Gamma_{C_1}\Gamma_{C_2}\Gamma_{C_1}	&= \Gamma_{C_1}\Gamma_{C_2}\Gamma_{C_1}\Gamma_{L_1}\Gamma_{C_2}\Gamma_{L_1}\Gamma_{C_1}\Gamma_{C_2}\Gamma_{C_1},
\end{align*}
meaning that \eqref{arrangement2_second2} reads
\begin{equation*}
	\Gamma_{L_1}\Gamma_{C_2}\Gamma_{L_1} = \Gamma_{C_1}\Gamma_{C_2}\Gamma_{C_1}\Gamma_{L_1}\Gamma_{C_2}\Gamma_{L_1}\Gamma_{C_1}\Gamma_{C_2}\Gamma_{C_1},
\end{equation*}
which is equivalent to \eqref{arrangement2_second1}.

We now turn our attention to \eqref{arrangement2_final_5}.
The first part reads
\begin{equation*}
	\Gamma_{L_2}\Gamma_{C_2}\Gamma_{L_1}\Gamma_{C_2}\Gamma_{C_1}\Gamma_{L_2} = \Gamma_{C_2}\Gamma_{L_1}\Gamma_{C_2}\Gamma_{C_1}.
\end{equation*}
Because $ \Gamma_{C_1} $ commutes with $ \Gamma_{L_2} $ it is equivalent to $ [\Gamma_{L_2}, \Gamma_{C_2}\Gamma_{L_1}\Gamma_{C_2}]=e $. We assert that this relation follows from existing ones. Indeed, because $ [\Gamma_{L_1}\Gamma_{L_2}, \Gamma_{C_2}]=e $ and $ [\Gamma_{L_1}\Gamma_{L_2}, \Gamma_{L_1}]=e $, we have $ {[\Gamma_{L_1}\Gamma_{L_2}, \Gamma_{C_2}\Gamma_{L_1}\Gamma_{C_2}]=e} $, and the relation $ \{ \Gamma_{L_1}, \Gamma_{C_2}\}=e $ is equivalent to $ [\Gamma_{L_1}, \Gamma_{C_2}\Gamma_{L_1}\Gamma_{C_2}] = e $, which together give us the desired $ [\Gamma_{L_2}, \Gamma_{C_2}\Gamma_{L_1}\Gamma_{C_2}] = e $.

The second part of \eqref{arrangement2_final_5} becomes after the substitution
\begin{equation*}
	\Gamma_{L_1}\Gamma_{C_2}\Gamma_{C_1} = \Gamma_{C_1}\Gamma_{L_1}\Gamma_{C_2}\Gamma_{C_1}\Gamma_{C_2}\Gamma_{C_1}\Gamma_{C_2}.
\end{equation*}
Obviously $ \Gamma_{L_1} $ cancels out and we get $ \{\Gamma_{C_1}, \Gamma_{C_2}\}=e $.

We are left with the two commutation relations in \eqref{arrangement2_final_6} that contain $ \Gamma_{L_4} $.
First, using $ \{ \Gamma_{C_1}, \Gamma_{C_2} \} = \{ \Gamma_{L_1}, \Gamma_{C_2} \} = e $ and $ [\Gamma_{C_1}, \Gamma_{L_2} ] =e $ (in this order), we get
\begin{equation*}
	\begin{split}
		[\Gamma_{L_4}, \Gamma_{C_1}\Gamma_{C_2}\Gamma_{C_1}] & = \Gamma_{L_2}\Gamma_{C_1}\Gamma_{C_2}\Gamma_{L_1}\Gamma_{C_1}\Gamma_{C_2}\Gamma_{C_1}\Gamma_{C_2}\Gamma_{C_1}\Gamma_{C_2}\Gamma_{C_1}\Gamma_{L_1}\Gamma_{C_2}\Gamma_{C_1}\Gamma_{L_2}\Gamma_{C_1}\Gamma_{C_2}\Gamma_{C_1} = \\
		& = \Gamma_{L_2}\Gamma_{C_1}\Gamma_{C_2}\Gamma_{L_1}\Gamma_{C_2}\Gamma_{L_1}\Gamma_{C_2}\Gamma_{C_1}\Gamma_{L_2}\Gamma_{C_1}\Gamma_{C_2}\Gamma_{C_1} = \\
		& = \Gamma_{L_2}\Gamma_{C_1}\Gamma_{L_1}\Gamma_{C_2}\Gamma_{L_1}\Gamma_{C_1}\Gamma_{L_2}\Gamma_{C_1}\Gamma_{C_2}\Gamma_{C_1} = \\
		& = \Gamma_{C_1}\Gamma_{L_2}\Gamma_{L_1}\Gamma_{C_2}\Gamma_{L_1}\Gamma_{L_2}\Gamma_{C_2}\Gamma_{C_1},
	\end{split}
\end{equation*}
which is trivial by $ [\Gamma_{L_1}\Gamma_{L_2}, \Gamma_{C_2}] = e $.
Finally,
\begin{equation*}
	\begin{split}
		[\Gamma_{L_4}, \Gamma_{L_1}\Gamma_{C_2}\Gamma_{L_1}] &= \Gamma_{L_2}\Gamma_{C_1}\Gamma_{C_2}\Gamma_{L_1}\Gamma_{C_1}\Gamma_{C_2}\Gamma_{L_1}\Gamma_{C_2}\Gamma_{L_1}\Gamma_{C_2}\Gamma_{C_1}\Gamma_{L_1}\Gamma_{C_2}\Gamma_{C_1}\Gamma_{L_2}\Gamma_{L_1}\Gamma_{C_2}\Gamma_{L_1} = \\
		&=  \Gamma_{L_2}\Gamma_{C_1}\Gamma_{C_2}\Gamma_{C_1}\Gamma_{C_2}\Gamma_{C_1}\Gamma_{C_2}\Gamma_{C_1}\Gamma_{L_2}\Gamma_{L_1}\Gamma_{C_2}\Gamma_{L_1} = \\
		&=  \Gamma_{L_2}\Gamma_{C_2}\Gamma_{L_2}\Gamma_{L_1}\Gamma_{C_2}\Gamma_{L_1},
	\end{split}
\end{equation*}
which is redundant because $ \Gamma_{C_2} $ commutes with $ \Gamma_{L_1}\Gamma_{L_2} $. This finishes the proof of the first assertion in the lemma.

\def\additionalgen{\widetilde{\Gamma}_L}

To prove item $ (2) $ of the lemma, define $ \additionalgen = \Gamma_{L_1}\Gamma_{L_2} \in G_6 $.
Because $ \Gamma_{L_2} = \Gamma_{L_1}\additionalgen $, group $ G_6 $ is generated by $ \Gamma_{C_1},\Gamma_{C_2},\Gamma_{L_1},\additionalgen $, we compute the relations on those generators.
By $ [\Gamma_{L_1}, \Gamma_{L_2}] = e $ we immediately get that $ \additionalgen $ and $ \Gamma_{L_1} $ commute.
Thus, $ \Gamma_{L_2}^2 = e $ becomes $ \additionalgen^2 = e $.
Because $ \Gamma_{C_1} $ commutes with both $ \Gamma_{L_1} $ and $ \Gamma_{L_2} $, it commutes with $ \additionalgen $ as well.
We get that $ \additionalgen $ is in the center of $ G_6 $ and so the relation $ \{ \Gamma_{C_2}, \Gamma_{L_2} \} = e $ follows from $ \{ \Gamma_{C_2}, \Gamma_{L_1} \} =e $.
The above discussion implies that the generators $ \Gamma_{C_1}, \Gamma_{C_2} $, and $ \Gamma_{L_1} $ satisfy exactly the relations among \eqref{arrangement2_second_final_1}-\eqref{arrangement2_second_final_4} that do not involve $ \Gamma_{L_2} $, which is precisely the assertion of item $ (2) $ of the lemma.

Now, consider group $G_6' = G_6/\langle\Gamma_{\additionalgen} \rangle $.
The group generated by $ \Gamma_{C_1},\Gamma_{C_2} $, and $ \Gamma_{L_1} $, subject to the relations $ [\Gamma_{C_1},\Gamma_{L_1}] = \{ \Gamma_{C_1}, \Gamma_{C_2} \} = \{ \Gamma_{C_2}, \Gamma_{L_1} \} = e $, is the affine Coxeter group of type $ \widetilde{C}_2 $.
Thus, $ G_6' $ is its quotient by the relation $ [\Gamma_{C_2}, \Gamma_{L_1}\Gamma_{C_1}\Gamma_{C_2}\Gamma_{C_1}\Gamma_{L_1}]=e $.
The generators are mapped to this representation of $ \widetilde{C}_2 $ as follows:
\begin{equation*}
	\Gamma_{C_1} \mapsto \begin{pmatrix}
		1 & 2 \\ \overline{1} & 2
	\end{pmatrix} \; ; \; \Gamma_{C_2} \mapsto \begin{pmatrix}
		1 & 2 \\ 2 & 1
	\end{pmatrix} \; ; \; \Gamma_{L_1} \mapsto \begin{pmatrix}
		1 & 2 \\ 1 & \overline{2}^1
	\end{pmatrix}
\end{equation*}
where an over line denotes $ \epsilon=-1 $ and we omit positive signs and zero exponents.

We now consider the additional relation, by computing it explicitly in $ \widetilde{C}_2 $, as follows:
\begin{equation*}
	\Gamma_{C_2}\Gamma_{L_1}\Gamma_{C_1}\Gamma_{C_2}\Gamma_{C_1}\Gamma_{L_1} \mapsto \begin{pmatrix}
		1 & 2 \\ 1^{-1} & 2^1
	\end{pmatrix}
\end{equation*}
\begin{equation*}
	\Gamma_{L_1}\Gamma_{C_1}\Gamma_{C_2}\Gamma_{C_1}\Gamma_{L_1}\Gamma_{C_2} \mapsto \begin{pmatrix}
		1 & 2 \\ 1^1 & 2^{-1}
	\end{pmatrix}.
\end{equation*}
This means that the elements of $ G_6' $ are signed permutations, and we identify two signed permutations if the underlying regular permutation (in $ S_2 $) is the same, they have the same signs, the same sum of exponents, and the same parity of exponents.

We are now ready to prove that $ [\Gamma_{L_1}, \Gamma_{C_2}\Gamma_{C_1}\Gamma_{C_2}] \ne e $.
For that, compare the images of $ \Gamma_{L_1}\Gamma_{C_2}\Gamma_{C_1}
\Gamma_{C_2} $ and $ \Gamma_{C_2}\Gamma_{C_1}\Gamma_{C_2}\Gamma_{L_1} $ in $ G_6' $:
\begin{equation*}
	\Gamma_{L_1}\Gamma_{C_2}\Gamma_{C_1}\Gamma_{C_2} \mapsto \begin{pmatrix}
		1 & 2 \\ 1 & \overline{2}^1
	\end{pmatrix}\cdot \begin{pmatrix}
		1 & 2 \\ 1 & \overline{2}
	\end{pmatrix} = \begin{pmatrix}
		1 & 2 \\ 1 & 2^{-1}
	\end{pmatrix}
\end{equation*}
\begin{equation*}
	\Gamma_{C_2}\Gamma_{C_1}\Gamma_{C_2}\Gamma_{L_1} = \begin{pmatrix}
		1 & 2 \\ 1 & \overline{2}
	\end{pmatrix} \cdot \begin{pmatrix}
		1 & 2 \\ 1 & \overline{2}^1
	\end{pmatrix} = \begin{pmatrix}
		1 & 2 \\ 1 & 2^1
	\end{pmatrix}.
\end{equation*}
Because the sums of exponents are different, the images of $ \Gamma_{L_1}\Gamma_{C_2}\Gamma_{C_1}\Gamma_{C_2} $ and $ \Gamma_{C_2}\Gamma_{C_1}\Gamma_{C_2}\Gamma_{L_1} $ in
$$ G_6' \cong \left. \widetilde{C}_2 \middle/ \left\langle \begin{pmatrix}
	1 & 2 \\ 1^2 & 2^{-2}
\end{pmatrix} \right\rangle \right. $$
are different, meaning that $ [\Gamma_{L_1}, \Gamma_{C_2}\Gamma_{C_1}\Gamma_{C_2}] \ne e $ in $ G_6 $.

Finally, we deduce that $ L_4 $ does not commute with $ L_1 $ in $ G_6 $ and thus, it is not contained in the center of $ G_6 $, meaning that it is not in the center of $ \pi_1(\mathbb{CP}^2 - \mathcal{B}_6) $, as well.
Indeed, because $ [\Gamma_{L_1}, \Gamma_{C_2}\Gamma_{C_1}\Gamma_{C_2}] \ne e $ and $ [\Gamma_{L_1}, \Gamma_{C_2}\Gamma_{L_1}\Gamma_{C_2}] = e $ we get that $ [\Gamma_{L_1}, \Gamma_{C_2}\Gamma_{C_1}\Gamma_{L_1}\Gamma_{C_2}] \ne e $, that is,
$$ [\Gamma_{L_1}, \Gamma_{L_4} ] = [\Gamma_{L_1}, \Gamma_{C_2}\Gamma_{C_1}\Gamma_{L_1}\Gamma_{C_2}\Gamma_{C_1}\Gamma_{L_2}] \ne e. $$
\end{proof}

\subsection{$\mathcal{B}_5$ and $\mathcal{B}_6$ are Zariski pair with non-isomorphic groups}

\begin{thm}\label{thm:tokunaga_pair1}
$ \left(\mathcal{B}_5, \mathcal{B}_6\right) $ is a Zariski pair with non-isomorphic fundamental groups.
\end{thm}
\begin{proof}
It is a direct and easy verification that $ \mathcal{B}_5 $ and $ \mathcal{B}_6 $ has the same combinatorics. Assume, towards a contradiction, that $ \varphi: \left(\mathbb{CP}^2, \mathcal{B}_5\right) \rightarrow \left(\mathbb{CP}^2, \mathcal{B}_6\right) $ is an homeomorphism. Then $ \varphi\left(\left\{P_1,P_3\right\}\right)=\left\{P_1, P_4\right\} $ so $ \varphi\left(L_3\cup C_1\right)=L_4\cup C_1 $ meaning that $ \varphi\left(C_1\right)=C_1 $ and $ \varphi\left( L_3 \right) = L_4  $.

Now, consider the induced map on fundamental groups $ \varphi_* : \pcpt{\mathcal{B}_5} \to \pcpt{\mathcal{B}_6} $.
It must send a loop around $ L_3 $ to a loop around $ L_4 $, meaning that $ \varphi_*(\Gamma_{L_3}) $ is conjugate either to $ L_4 $ or to $ L_4^{-1} $.
Moreover, $ \varphi_* $ sends squares of the generators to squares of the generators, that is, it induces isomorphism $
\varphi_*:G_5\to G_6 $ and this isomorphism sends $ \Gamma_{L_3} $ to a conjugate of $ \Gamma_{L_4} $.
Now, because $ L_3 $ is in the center of $ G_5 $ by Lemma \ref{tokunaga-arrangement1-b1-group-calc}, $ \varphi_*(L_3) $ must be in the center of $ G_6 $ and thus, also $ \Gamma_{L_4} $, in contradiction to Lemma \ref{tokunaga-arrangement1-b2-group-calc}.
\end{proof}


\begin{thebibliography}{10}
	\providecommand{\url}[1]{{#1}}
	\providecommand{\urlprefix}{URL }
	\expandafter\ifx\csname urlstyle\endcsname\relax
	  \providecommand{\doi}[1]{\discretionary{}{}{}#1}\else
	  \providecommand{\doi}{\discretionary{}{}{}\begingroup \urlstyle{rm}\Url}\fi
	
	\bibitem{new}
	Amram, M., Bannai, S., Shirane, T., Sinichkin, U., Tokunaga, H.o.: The
	  realization space of a certain conic line arrangement of degree 7 and a
	  $\pi_1$-equivalent {Z}ariski pair.
	\newblock To appear in Isr. J. Math.  (2025).
	\newblock \doi{arxiv.org/abs/2307.01736}
	
	\bibitem{AGT_order6_conic_line_arr}
	Amram, M., Garber, D., Teicher, M.: Fundamental groups of tangent conic-line
	  arrangements with singularities up to order 6.
	\newblock Math. Z. \textbf{256}, 837--870 (2007).
	\newblock \doi{10.1007/s00209-007-0109-4}
	
	\bibitem{arxiv_version}
	Amram, M., Shwartz, R., Sinichkin, U., Tan, S.L., Tokunaga, H.o.: Zariski pairs
	  of conic-line arrangements of degrees 7 and 8 via fundamental groups.
	\newblock Preprint at https://arxiv.org/abs/2106.03507 (2023)
	
	\bibitem{Bartolo94}
	Artal-Bartolo, E.: On {Z}ariski pairs (sur les couples de {Z}ariski).
	\newblock J. Algebraic Geom. \textbf{3}, 223--247 (1994).
	\newblock \doi{Zbl 0524.14026}
	
	\bibitem{bartolo_tokunaga2020torsion}
	Artal-Bartolo, E., Bannai, S., Shirane, T., Tokunaga, H.o.: {Torsion divisors
	  of plane curves and {Z}ariski pairs}.
	\newblock Algebra i Analiz \textbf{34}(5), 1--22 (2022).
	\newblock \doi{10.1090/spmj/1776}
	
	\bibitem{bartolo_tokunaga2020torsion2}
	Artal-Bartolo, E., Bannai, S., Shirane, T., Tokunaga, H.o.: {Torsion divisors
	  of plane curves with maximal flexes and {Z}ariski pairs}.
	\newblock Math. Nachr. \textbf{296}(6), 2214--2235 (2023).
	\newblock \doi{10.1002/mana.202000319}
	
	\bibitem{Artal2}
	Artal-Bartolo, E., Carmona~Ruber, J., Cogolludo-Agustin, J.: Braid monodromy
	  and topology of plane curves.
	\newblock Duke Math. J. \textbf{118} (2003).
	\newblock \doi{DOI: 10.1215/S0012-7094-03-11823-2}
	
	\bibitem{Artal3}
	Artal-Bartolo, E., Carmona~Ruber, J., Cogolludo~Agustin, J., Luengo~Velasco,
	  I., Melle~Hernandez, A.: Fundamental group of plane curves and related
	  invariants.
	\newblock In: In Contribuciones matemáticas: libro-homenaje al profesor D.
	  Joaquín Arregui Fernández. Editorial Complutense, Madrid, pp. 77--104
	  (2021).
	\newblock
	  \doi{https://docta.ucm.es/rest/api/core/bitstreams/d583f115-5bfc-4a2c-b708-e2eea41a074d/content}
	
	\bibitem{AC}
	Artal-Bartolo, E., Cogolludo~Agustin, J.: Some open questions on arithmetic
	  {Z}ariski pairs.
	\newblock In: Singularities in geometry, topology, foliations and dynamics. A
	  celebration of the 60th birthday of José Seade. Selected papers based on the
	  presentations at the workshop, Mérida, Mexico, December 8--19, 2014., pages
	  31--54. Cham: Birkhäuser, 2017.
	
	\bibitem{artal1}
	Artal-Bartolo, E., Cogolludo-Agustin, J., Martin-Morales, J.: Cremona
	  transformations of weighted projective planes, {Z}ariski pairs, and rational
	  cuspidal curves.
	\newblock In: J.~Fernández~de Bobadilla, T.~László, A.~Stipsicz (eds.)
	  Trends in Mathematics, pp. 117--157. Birkhäuser, Cham (2021).
	\newblock \doi{https://doi.org/10.1007/978-3-030-61958-9_7}
	
	\bibitem{Bartolo_tokunaga_survey}
	Artal-Bartolo, E., Cogolludo-Agustin, J., Tokunaga, H.o.: A survey on {Z}ariski
	  pairs.
	\newblock In Algebraic geometry in East Asia-Hanoi 2005. Proceedings of the 2nd
	  international conference on algebraic geometry in East Asia, Hanoi, Vietnam,
	  October 10-14, 2005, pages 1-100. Tokyo: Mathematical Society of Japan
	  (2008).
	\newblock \doi{10.2969/aspm/05010001}
	
	\bibitem{Dimca}
	Artal-Bartolo, E., Dimca, A.: On fundamental groups of plane curve complements.
	\newblock Ann. Univ. Ferrara Sez. VII Sci. Mat. \textbf{61}, 255--262 (2004).
	\newblock \doi{https://doi.org/10.1007/s11565-015-0231-x}
	
	\bibitem{ben2}
	Bannai, S., Guerville-Ballé, B., Shirane, T.: Zariski pairs of conic-line
	  arrangements with a unique conic.
	\newblock Preprint at https://arxiv.org/abs/2410.04969 (2024)
	
	\bibitem{bannai_tokunaga2019zariski}
	Bannai, S., Tokunaga, H.o.: {Zariski tuples for a smooth cubic and its tangent
	  lines}.
	\newblock Proceedings of the Japan Academy, Series A, Mathematical Sciences
	  \textbf{96}(2), 18--21 (2020).
	\newblock \doi{10.3792/pjaa.96.004}
	
	\bibitem{bannai-tokunaga-yorisaki}
	Bannai, S., Tokunaga, H.o., Yorisaki, E.: The realization spaces of certain
	  conic-line arrangements of degree 7.
	\newblock Preprint at https://arxiv.org/abs/2409.05011 (2024)
	
	\bibitem{Bjorner2005}
	Bjorner, A., Brenti, F.: Combinatorics of {C}oxeter Groups.
	\newblock Springer-Verilag Berlin, Heidelberg (2005).
	\newblock \doi{10.1007/3-540-27596-7}
	
	\bibitem{DEG_isotopy}
	Degtyarev, A.: Isotopic classification of complex plane projective curves of
	  degree~5.
	\newblock Algebra i Analiz \textbf{1}, 78--101 (1989)
	
	\bibitem{DEG}
	Degtyarev, A.: On deformations of singular plane sextics.
	\newblock J. Algebraic Geom. \textbf{17}, 101--135 (2008).
	\newblock \doi{10.1090/S1056-3911-07-00469-9}
	
	\bibitem{Kap}
	Drutu, C., Kapovich, M.: Geometric Group Theory (Colloquium Publications).
	\newblock AMS (2018).
	\newblock \doi{https://www.ams.org/books/coll/063/}
	
	\bibitem{Friedmann_Garber_conicline_fundamental}
	Friedman, M., Garber, D.: On the structure of conjugation-free fundamental
	  groups of conic-line arrangements.
	\newblock J. Homotopy Rel. Struct. \textbf{10}, 685--734 (2015).
	\newblock \doi{10.1007/s40062-014-0081-8}
	
	\bibitem{shustin_appendix}
	Friedman, M., Leyenson, M., Shustin, E.: On ramified covers of the projective
	  plane i: Interpreting {S}egre's theory (with an appendix by {E}ugenii
	  {S}hustin).
	\newblock Internat. J. Math. \textbf{22} (2011).
	\newblock \doi{10.1142/S0129167X11006945}
	
	\bibitem{Ben}
	Guerville-Ballé, B.: Topology and homotopy of lattice isomorphic arrangements.
	\newblock In: Proc. Am. Math. Soc., 148(5), pp. 2193--2200 (2020).
	\newblock \doi{10.1090/proc/14878}
	
	\bibitem{vanKampen33}
	Kampen, E.R.V.: On the fundamental group of an algebraic curve.
	\newblock Amer. J. Math. \textbf{55}(1), 255--260 (1933).
	\newblock \doi{http://www.jstor.org/stable/2371128}
	
	\bibitem{sirocco}
	Marco-Buzunariz, M.{\'A}., Rodr{\'i}guez, M.: Sirocco: A library for certified
	  polynomial root continuation.
	\newblock In: G.M. Greuel, T.~Koch, P.~Paule, A.~Sommese (eds.) Mathematical
	  Software -- ICMS 2016, pp. 191--197. Springer International Publishing, Cham
	  (2016)
	
	\bibitem{Namba2004OnTF}
	Namba, M., Tsuchihashi, H.: On the fundamental groups of {G}alois covering
	  spaces of the projective plane.
	\newblock Geom. Dedicata \textbf{105}, 85--105 (2004).
	\newblock \doi{10.1023/B:GEOM.0000024693.88450.d9}
	
	\bibitem{Nazir_Yoshinaga_line_arrangements}
	Nazir, S., Yoshinaga, M.: On the connectivity of the realization spaces of line
	  arrangements.
	\newblock Ann. Sc. Norm. Super. Pisa Cl. Sci. (5) \textbf{Ser. 5, 11}(4),
	  921--937 (2012).
	\newblock \doi{10.2422/2036-2145.201009_003}
	
	\bibitem{Oka1999FlexCA}
	Oka, M.: Flex curves and their applications.
	\newblock Geom. Dedicata \textbf{75}, 67--100 (1999).
	\newblock \doi{10.1023/A:1005004123844}
	
	\bibitem{Oka2002}
	Oka, M.: A new {A}lexander-equivalent {Z}ariski pair.
	\newblock Acta Math. Vietnam. \textbf{27}, 349--357 (2002)
	
	\bibitem{Randell}
	Randell, R.: Lattice-isotopic arrangements are topologically isomorphic.
	\newblock Proc. Amer. Math. Soc. \textbf{107}(2), 555--559 (1989).
	\newblock \doi{10.1090/S0002-9939-1989-0984812-7}
	
	\bibitem{Rybnikov_line_arrangement}
	Rybnikov, G.: On the fundamental group of the complement of a complex
	  hyperplane arrangement.
	\newblock Funct. Anal. Appl. \textbf{45}, 137--148 (2011).
	\newblock \doi{10.1007/s10688-011-0015-8}
	
	\bibitem{shimada}
	Shimada, I.: A note on {Z}ariski pairs.
	\newblock Compos. Math. \textbf{104} (1995).
	\newblock \doi{10.48550/arXiv.alg-geom/9505007}
	
	\bibitem{takahashi_tokunaga2020explicit}
	Takahashi, A., Tokunaga, H.o.: {An explicit construction for $n$-contact curves
	  to a smooth cubic via divisions and {Z}ariski tuples}.
	\newblock Hokkaido Math. J. \textbf{51}(3), 389--405 (2022).
	\newblock \doi{10.14492/hokmj/2020-391}
	
	\bibitem{tokunaga2014}
	Tokunaga, H.o.: Sections of elliptic surfaces and {Z}ariski pairs for
	  conic-line arrangements via dihedral covers.
	\newblock J. Math. Soc. Japan \textbf{66}(2), 613--640 (2014).
	\newblock \doi{10.2969/jmsj/06620613}
	
	\bibitem{Ye2013ClassificationOM}
	Ye, F.: Classification of moduli spaces of arrangements of 9 projective lines.
	\newblock Pacific J. Math. \textbf{265}, 243--256 (2013).
	\newblock \doi{http://dx.doi.org/10.2140/pjm.2013.265.243}
	
	\bibitem{zariski29}
	Zariski, O.: On the problem of existence of algebraic functions of two
	  variables possessing a given branch curve.
	\newblock Amer. J. Math. \textbf{51}(2), 305--328 (1929).
	\newblock \doi{http://www.jstor.org/stable/2370712}
	
	\end{thebibliography}
\end{document}